%amalgamated version of slsi paper 12.02.2013 

\font\Bigrm = cmr10 scaled \magstep 2
\font\bigrm = cmr10 scaled \magstep 1
\font\Bigmath = cmmi10 scaled \magstep 2
 1
\font\Bigbb = msbm10 scaled \magstep 2
 1
\input amssym

\centerline{\Bigrm SPECTRAL LARGE SIEVE INEQUALITIES}

\smallskip

\centerline{{\Bigrm FOR HECKE CONGRUENCE SUBGROUPS OF
{\Bigmath SL}({\Bigrm 2},$\,${\Bigbb Z}{\bigrm [{\Bigmath i}]})} }

\medskip\centerline{by}\medskip 

\centerline{\bigrm NIGEL WATT}

\vskip 8 mm

\footnote{{}}{\hskip -2 mm
{\bf 2010 Mathematics Subject Classification:}\quad
11F$30^{*}$, 11F37, 11F70, 11F72, 11L05, 11L07, 11M41, 
11N35, 11R42, 22E30, 33C10, 44A15.}

\footline={\hss{\vbox to 2cm{\vfil\hbox{\rm\folio}}}\hss}

\noindent{\bf Abstract:}\ 
We prove, in respect of an arbitrary Hecke congruence subgroup 
$\Gamma =\Gamma_0(q_0)\leq SL(2,{\Bbb Z}[i])$, 
some new upper bounds (or `spectral large sieve inequalities') for sums involving 
Fourier coefficients of $\Gamma$-automorphic cusp forms on $SL(2,{\Bbb C})$. 
The Fourier coefficients in question may arise from the Fourier expansion 
at any given cusp ${\frak c}$ of $\Gamma\,$ (our results are not 
limited to the case ${\frak c}=\infty$). For this reason, our proof 
is reliant upon an  
extension, to arbitrary cusps, of the spectral-Kloosterman sum formula 
for  $\Gamma\backslash SL(2,{\Bbb C})$ obtained by Hristina Lokvenec-Guleska 
in her doctoral thesis (generalising the sum formulae   
of Roelof Bruggeman and Yoichi Motohashi 
for $PSL(2,{\Bbb Z}[i])\backslash PSL(2,{\Bbb C})$ in several respects, 
though not as regards the choice of cusps). A proof of the required 
extension of the sum formula is given in an appendix.  

\medskip 


\vskip 6mm

\noindent{\bf Contents}
 
{\settabs\+\quad &5..\ &1.5.\ 
&Schwartz Spaces, Fourier Integrals, Poisson Summation and the 
Analytic Large Sieve\qquad &\cr  %sample line 
\+&&&&Page\cr 
\+&Outline of results and methods\ .\ .\ .\ .\ .\ .\ .\ .\ .\ .\ .\ .\ .\ 
.\ .\ .\ .\ .\ .\ .\ .\ .\ .\ .\ .\ .\ .\ .\ .\ .\ .\ .\ .\ .\ .\ .\ .\ .\ 
.\ .\ .\ .\ .\ .\ .\ .\ .&&&2\cr 
\smallskip
\+&1.&Definitions and statements of the results\ .\ .\ .\ .\ 
.\ .\ .\ .\ .\ .\ .\ .\ .\ .\ .\ .\ .\ .\ .\ .\ .\ .\ .\ .\ .\ .\ .\ .\ .\ 
.\ .\ .\ .\ .\ .\ .\ .\ .&&2\cr
\+&&1.1&The quotient $\Gamma\backslash G\,$: 
coordinates, measure, cusps and fundamental domains\ \ .\ .\ .\ .\ .\ .\ .\ .\ .\ .&2\cr  
\+&&1.2&$\Gamma$-automorphic functions, Casimir operators and the 
Laplacian for ${\Bbb H}_3$\ .\ .\ .\ .\ .\ .\ .\ .\ .\ .\ .\ .\ .\ .&5\cr  
\+&&1.3&Functions of $K$-type $(\ell,q)\,$\ .\ .\ .\ .\ .\ .\ .\ .\ .\ .\ .\ .\ .\ .\ .\ .\ .\ .\ .\ .\ .\ .\ 
.\ .\ .\ .\ .\ .\ .\ .\ .\ .\ .\ .\ .\ .\ .\ .\ .\ .\ .\ .\ .\ .\ .&8\cr
\+&&1.4&Fourier expansions at cusps; spaces of cusp forms\ \ .\ .\ .\ .\ .\ .\ .\ .\ .\ .\ .\ .\ .\ .\ .\ .\ .\ .\ .\ .\ .\ .\ 
.\ .\ .\ .\ .\ .&9\cr  
\+&&1.5&The Jacquet integral; generalised Kloosterman sums and the Fourier expansion of&\cr 
\+&&&\ Poincar\'{e} series; Fourier coefficients of cusp forms$\,$\ .\ .\ .\ .\ .\ .\ .\ .\ .\ .\ .\ .\ .\ .\ .\ .\ .\ .\ .\ .\ .\ .\ 
.\ .\ .\ .\ .\ .&11\cr  
\+&&1.6&The spaces $H(\nu,p)$ of $K$-finite functions; principal and 
complementary series\ .\ .\ .\ .\ .\ .\ .\ .\ .&13\cr  
\+&&1.7&Decomposing the space $L^2(\Gamma\backslash G)$\ \ .\ .\ .\ .\ .\ .\ .\ .\ .\ .\ .\ .\ .\ .\ .\ .\ .\ .\ .\ .\ .\ .\ 
.\ .\ .\ .\ .\ .\ .\ .\ .\ .\ .\ .\ .\ .\ .\ .\ .\ .&14\cr  
\+&&1.8&Decomposing the subspace ${}^{\rm e}L^2(\Gamma\backslash G)$: 
the Eisenstein series and a Parseval identity\ \ .\ .\ .\ .&17\cr 
\+&&1.9&Results and applications\ \ .\ .\ .\ .\ .\ .\ .\ .\ .\ .\ .\ .\ .\ .\ .\ .\ .\ .\ .\ .\ .\ .\ 
.\ .\ .\ .\ .\ .\ .\ .\ .\ .\ .\ .\ .\ .\ .\ .\ .\ .\ .\ .\ .\ .\ .\ .&19\cr   
\smallskip 
\+&Notation\ .\ .\ .\ .\ .\ .\ .\ .\ .\ .\ .\ .\ .\ .\ .\ .\ .\ .\ .\ .\ .\ .\ 
.\ .\ .\ .\ .\ .\ .\ .\ .\ .\ .\ .\ .\ .\ .\ .\ .\ .\ .\ .\ .\ .\ .\ .\ .\ 
.\ .\ .\ .\ .\ .\ .\ .\ .\ .\ .\ .\ .\ .\ .\ .\ .&&&26\cr  
\smallskip 
\+&2.&Lemmas\ .\ .\ .\ .\ .\ .\ .\ .\ .\ .\ .\ .\ .\ .\ .\ .\ .\ .\ .\ .\ .\ .\ 
.\ .\ .\ .\ .\ .\ .\ .\ .\ .\ .\ .\ .\ .\ .\ .\ .\ .\ .\ .\ .\ .\ .\ .\ .\ 
.\ .\ .\ .\ .\ .\ .\ .\ .\ .\ .\ .\ .\ .&&31\cr 
\smallskip 
\+&3.&The proof of Proposition~2$\,$\ \ .\ .\ .\ .\ .\ .\ .\ .\ .\ .\ .\ .\ .\ .\ .\ .\ .\ .\ .\ .\ .\ .\ 
.\ .\ .\ .\ .\ .\ .\ .\ .\ .\ .\ .\ .\ .\ .\ .\ .\ .\ .\ .\ .\ .\ .\ .\ .&&40\cr 
\smallskip 
\+&4.&Further lemmas\ \ .\ .\ .\ .\ .\ .\ .\ .\ .\ .\ .\ .\ .\ .\ .\ .\ .\ .\ .\ .\ .\ .\ 
.\ .\ .\ .\ .\ .\ .\ .\ .\ .\ .\ .\ .\ .\ .\ .\ .\ .\ .\ .\ .\ .\ .\ .\ .\ 
.\ .\ .\ .\ .\ .\ .\ .&&50\cr 
\smallskip 
\+&5.&The proof of Theorem~1\ \ .\ .\ .\ .\ .\ .\ .\ .\ .\ .\ .\ .\ .\ .\ .\ .\ .\ .\ .\ .\ .\ .\ 
.\ .\ .\ .\ .\ .\ .\ .\ .\ .\ .\ .\ .\ .\ .\ .\ .\ .\ .\ .\ .\ .\ .\ .\ .\ 
.\ .&&58\cr 
\smallskip 
\+&6.&Appendix on the proof of the sum formula$\,$\ .\ .\ .\ .\ .\ .\ .\ .\ .\ .\ .\ .\ .\ .\ .\ .\ .\ .\ .\ .\ .\ .\ 
.\ .\ .\ .\ .\ .\ .\ .\ .\ .\ .\ .\ .\ .&&68\cr 
\smallskip 
\+&&6.1&Generalised Kloosterman sums\ \ .\ .\ .\ .\ .\ .\ .\ .\ .\ .\ .\ .\ .\ .\ .\ .\ .\ .\ .\ .\ .\ .\ 
.\ .\ .\ .\ .\ .\ .\ .\ .\ .\ .\ .\ .\ .\ .\ .\ .\ .\ .&69\cr 
\+&&6.2&Poincar\'{e} series$\,$\ .\ .\ .\ .\ .\ .\ .\ .\ .\ .\ .\ .\ .\ .\ .\ .\ .\ .\ .\ .\ .\ .\ 
.\ .\ .\ .\ .\ .\ .\ .\ .\ .\ .\ .\ .\ .\ .\ .\ .\ .\ .\ .\ .\ .\ .\ .\ .\ 
.\ .\ .\ .\ .\ .&74\cr 
\+&&6.3&The Goodman-Wallach operator ${\bf M}_{\omega}\,$\ .\ .\ .\ .\ .\ .\ .\ .\ .\ .\ .\ .\ .\ .\ .\ .\ .\ .\ .\ .\ .\ .\ 
.\ .\ .\ .\ .\ .\ .\ .\ .\ .\ .\ .\ .\ .\ .&82\cr 
\+&&6.4&The Lebedev transform ${\bf L}^{\omega}_{\ell,q}$ and 
auxilliary test functions\ .\ .\ .\ .\ .\ .\ .\ .\ .\ .\ .\ .\ .\ .\ .\ .\ .\ .\ .\ .\ .\ .\ 
.&84\cr 
\+&&6.5&Poincar\'{e} series revisited\ \ .\ .\ .\ .\ .\ .\ .\ .\ .\ .\ .\ .\ .\ .\ .\ .\ .\ .\ .\ .\ .\ .\ 
.\ .\ .\ .\ .\ .\ .\ .\ .\ .\ .\ .\ .\ .\ .\ .\ .\ .\ .\ .\ .\ .\ .\ .&86\cr 
\+&&6.6&The preliminary spectral summation formula\ \ .\ .\ .\ .\ .\ .\ .\ .\ .\ .\ .\ .\ .\ .\ .\ .\ .\ .\ .\ .\ .\ .\ 
.\ .\ .\ .\ .\ .\ .\ .\ .&108\cr 
\+&&6.7&Completing the proof of the spectral summation formula$\,$\ .\ .\ .\ .\ .\ .\ .\ .\ .\ .\ .\ .\ .\ .\ .\ .\ .\ .\ .\ .\ .\ .\ .&121\cr 
\smallskip 
\+&Acknowledgements$\,$\ \ .\ .\ .\ .\ .\ .\ .\ .\ .\ .\ .\ .\ .\ .\ .\ .\ .\ .\ .\ .\ .\ .\ 
.\ .\ .\ .\ .\ .\ .\ .\ .\ .\ .\ .\ .\ .\ .\ .\ .\ .\ .\ .\ .\ .\ .\ .\ .\ 
.\ .\ .\ .\ .\ .\ .\ .&&&135\cr 
\smallskip 
\+&References$\,$\ \ .\ .\ .\ .\ .\ .\ .\ .\ .\ .\ .\ .\ .\ .\ .\ .\ .\ .\ .\ .\ .\ .\ 
.\ .\ .\ .\ .\ .\ .\ .\ .\ .\ .\ .\ .\ .\ .\ .\ .\ .\ .\ .\ .\ .\ .\ .\ .\ 
.\ .\ .\ .\ .\ .\ .\ .\ .\ .\ .\ .\ .\ .&&&135\cr}

\bigskip 

\centerline{\bf Outline of results and methods.}

\medskip

In 1982, in Theorem~2 of their paper [9], Deshouillers and Iwaniec 
generalised Iwaniec's ground-breaking estimate [20], Theorem~1, to obtain 
similar `spectral large sieve inequalities' for Fourier coefficients of holomorphic cusp forms, 
non-holomorphic cusp forms, or Eisenstein series, automorphic (or `modular') 
with respect to the action of an abitrary Hecke congruence subgroup
of $SL(2,{\Bbb Z})$ upon the upper half complex plane ${\Bbb H}_2$.
This paper concerns corresponding results for Fourier coefficients of 
functions on $SL(2,{\Bbb C})$ that are automorphic with respect 
to some Hecke congruence subgroup $\Gamma$ of $SL(2,{\Bbb Z}[i])$.
Our principal results, obtained in Theorem~1 below,   
are not quite a perfect analogy of the results (1.28) and (1.29) of Theorem~2 of [9], 
and so seem open to further improvement. Another question left open is as to whether it is possible to 
achieve a refinement of our Theorem~1 paralleling  
the significant refinement 
of Deshouillers and Iwaniec's spectral large sieve inequalities
that was obtained by Jutila in Theorem~1.1 of [23]. 

We have modelled our proof of Theorem~1 on the 
proof of Theorem~2 of [9] that is contained in Section~5 of~[9]. 
Just as the proof of Theorem~2 of [9] is dependent on the 
estimates for sums of generalised Kloosterman that are supplied in Proposition~3 of [9], 
so too is our proof of Theorem~1 dependent on the estimates for sums of Kloosterman 
sums that we obtain in our Proposition~2 (stated at the end of Section~1). 
Our proof of Proposition~2 follows the same basic pattern as the 
proof (in [9], Subsection~5.1) of Proposition~3 of [9], but it does have some 
novel features (such as those relating to the `grossencharakter' factor 
$(\omega_1\omega_2/|\omega_1\omega_2|)^m$ which occurs in Equation~(1.9.25), below).

In the work [9] of Deshouillers and Iwaniec a crucial part is played by 
summation formulae of Bruggeman [3] and Kuznetsov [28], [29], 
expressing certain sums involving 
Fourier coefficients of modular forms in terms of sums of Kloosterman sums
(and vice versa). These particular summation formulae apply only to modular forms 
on ${\Bbb H}_2$ (i.e. the homogeneous space $SL(2,{\Bbb R})/SO(2,{\Bbb R})$\/), 
whereas in the present work one needs instead summation formulae for 
Fourier coefficients 
of automorphic functions on $SL(2,{\Bbb C})$. Such formulae were first 
obtained (for functions automorphic with respect to $SL(2,{\Bbb Z}[i])$) 
in Bruggeman and Motohashi's paper [5], and were extended to the 
case of functions automorphic with respect to Hecke congruence subgroups of 
$SL(2,{\Bbb Z}[i])$ in Lokvenec-Guleska's thesis [32]. These authors 
dealt only with Fourier coefficients for the cusp at infinity, but their methods 
can be  adapted to succesfully handle Fourier coefficients at other cusps:  
the relevant summation formula, Theorem~B, is stated at 
the beginning of Subsection~1.9 of this paper; the 
required adaptations of the proofs in [5] and [32] are discussed in 
an appendix (Section~6). 

An introduction to relevant concepts and terminology now follows 
(preparatory to the statement of the principal new results contained in this paper).

\bigskip

\centerline{\bf \S 1. Definitions and statements of the results.}

\medskip 

\centerline{\bf \S 1.1 The quotient $\Gamma\backslash G\,$: coordinates, measure, cusps and fundamental domains.}

\medskip

Let $G=SL(2,{\Bbb C})$ and $K=SU(2)\,$ 
(the maximal compact subgroup of the Lie group $G$\/). 
Let ${\frak O}$ denote the ring ${\Bbb Z}[i]$ of Gaussian integers; we shall use the 
notation $m\sim n$ to signify that $m$ and $n$ are associates 
(in the sense that $m\in n{\frak O}^*\subset{\frak O}$, where 
${\frak O}^*=\{ i,-1,-i,1\}$). 
Suppose moreover 
that $q_0\in{\frak O}-\{ 0\}\,$; and that $\Gamma$ is the Hecke congruence 
subgroup of $SL(2,{\frak O})$ given by
$$\Gamma =\Gamma_0\!\left( q_0\right) 
=\left\{ \pmatrix{a&b \cr c&d}\in SL(2,{\frak O}) : c\in q_0 {\frak O}\right\}\;.\eqno(1.1.1)$$
The group $\Gamma$ is a discrete and cofinite (but not cocompact) subgroup of $G$. 
Both $\Gamma$ and $G$ act by left multiplication on 
the homogeneous space $G/K$. By the Iwasawa decomposition  
each $g\in G$ has a unique factorisation of form 
$$g=n[z]a[r]k[\alpha,\beta]\;,\eqno(1.1.2)$$
where, for $z\in{\Bbb C}$, $r>0$, and $\alpha,\beta\in{\Bbb C}$ such that $|\alpha|^2+|\beta|^2=1$, 
one has:
$$n[z]=\pmatrix{1&z\cr 0&1}\;,\quad 
a[r]=\pmatrix{\sqrt{r}&0\cr 0&1/\sqrt{r}}\quad\hbox{and}\quad 
k[\alpha,\beta]=\pmatrix{\alpha&\beta\cr -\overline{\beta}&\overline{\alpha}}\in K\;.\eqno(1.1.3)$$
The mapping $gK\mapsto (z,r)$ (with $z$ and $r$ are as in (1.1.2)) defines a homeomorphism 
between $G/K$ and the topological space 
${\Bbb H}_3={\Bbb C}\times{\Bbb R}^{+}={\Bbb C}\times\{ r\in{\Bbb R} : r>0\}$. 
The action of $G$ upon $G/K$ (by left multiplication) may be interpreted as a continuous group action of $G$ upon
${\Bbb H}_3$ by putting $g(z,r)=(z',r')$ when $z,z'\in{\Bbb C}$, $r,r'>0$ and  
$gn[z]a[r]K=n[z']a[r']K$: by a calculation, one then has  
$$g(z,r)=\left( 
{(a z + b)\,\overline{(c z + d)} + a\,\overline{c}\,r^2\over 
|c z +d|^2 + |c|^2 r^2}\,,\,{r\over |c z +d|^2 + |c|^2 r^2}\right)\in{\Bbb H}_3
\quad\hbox{for $\,g=\pmatrix{a &b\cr c &d}\in G$, $(z,r)\in{\Bbb H}_3$.}\eqno(1.1.4)$$

By Proposition~7.3.9 of [11], the set
$${\cal F}_{{\Bbb Q}(i)}
=\left\{ (z,r)\in{\Bbb H}_3 : |z|^2+r^2\geq 1\ {\rm and}\ 
|{\rm Re}(z)|,{\rm Im}(z)\in[0,1/2]\right\}\eqno(1.1.5)$$
is a fundamental domain for the action of $SL(2,{\frak O})=\Gamma_0(1)$ upon ${\Bbb H}_3$. 
Consideration of the natural homomorphism from $SL(2,{\frak O})$ 
into $SL(2,{\frak O}/q_0 {\frak O})$ shows that 
$$\left[SL(2,{\frak O}) : \Gamma\right] = \left| q_0\right|^2 \prod_{(\varpi)\ni q_0} 
\left( 1+{1\over |\varpi|^2}\right)\;,\eqno(1.1.6)$$
where the product is taken over prime ideals $(\varpi)=\varpi{\frak O}\subset{\frak O}$
(see Section 2.4 of [21] for the $SL(2,{\Bbb Z})$-analogue of this).  
Therefore there exist $\gamma_1,\ldots ,\gamma_{[SL(2,{\frak O}) : \Gamma]}\in SL(2,{\frak O})$
such that the set
$${\cal F}=\bigcup_{k=1}^{[SL(2,{\frak O}) : \Gamma]}\gamma_k {\cal F}_{{\Bbb Q}(i)}\eqno(1.1.7)$$
is a fundamental domain for the action of $\Gamma$ upon ${\Bbb H}_3$.
Moreover, as noted in Section 2.2 of [11], there do exist fundamental domains 
for the action of 
$\Gamma$ upon ${\Bbb H}_3$ that have a connected interior
(these being the Poincar\'{e} normal polyhedrons ${\cal P}_{Q}(\Gamma)\subset{\Bbb H}_3$
 centred at $Q=(z,r)$\/). Therefore one may assume a choice of 
$\gamma_1,\ldots ,\gamma_{[SL(2,{\frak O}) : \Gamma]}$ 
in (1.1.7) that makes the interior of ${\cal F}$ be connected.

When $k=k[\alpha,\beta]\in K$ one has, for some  
$\theta\in[0,\pi)$ and some real $\varphi$ and $\psi$ 
satisfying $\varphi\pm\psi\in[0,4\pi)$, the factorisation 
$$k=h\bigl[ e^{i\varphi /2}\bigr] v[i\theta] h\bigl[ e^{i\psi /2}\bigr]\;,\eqno(1.1.8)$$
where  
$$h[u]=\pmatrix{u&0\cr 0&u^{-1}}\qquad\hbox{and}\qquad 
v[i\theta]=\pmatrix{\cos(\theta /2)&i\sin(\theta /2)\cr i\sin(\theta /2)&\cos(\theta /2)}\;.\eqno(1.1.9)$$
This factorisation is unique in the cases where $\theta$ satisfies $0<\theta<\pi$, 
so the `Iwasawa coordinates' $z,r,\theta,\varphi,\psi$ 
in (1.1.2) and (1.1.8)    
are a serviceable coordinate system for $G$. 
\par 
We define $C^{\infty}(G)$ to be the space of functions 
$f : G\rightarrow{\Bbb C}$  which are `smooth', in the sense that, for 
each $g_0\in G$, and each $j\in{\Bbb N}$, all $2^j$ of the partial 
derivatives of order $j$ of the function 
$(x,y,r,\theta,\varphi,\psi)\mapsto 
f(g_0 n[x+iy] a[r] h[e^{i\varphi /2}] v[i\theta] h[e^{i\psi /2}])$ 
are defined and continuous on ${\Bbb R}^5$. The space of smooth complex-valued 
functions on $K$ is denoted by $C^{\infty}(K)\,$ (one has $F\in C^{\infty}(K)$ 
if and only if $F$ is the restriction to $K$ of some element of $C^{\infty}(G)$). 
We define $C^0(G)$ to be the space of functions $f : G\rightarrow{\Bbb C}$ that 
are continuous with respect to the topology on $G$ defined in Section~2.1 of [11] 
(this just means that a function $f : G\rightarrow{\Bbb C}$ 
will lie in $C^0(G)$ if and only if $f(g)$ is continuous 
as a function of the Iwasawa coordinates of $g$).    

In terms of the parameterisations introduced in (1.1.3) and (1.1.8)-(1.1.9),  
the groups 
$$K=\left\{ k[\alpha,\beta] : \alpha,\beta\in{\Bbb C}\ {\rm and}\ |\alpha|^2+|\beta|^2=1\right\} , \quad 
A=\{ a[r] : r>0\}\quad\hbox{and}\quad N=\{ n[z] : z\in{\Bbb C}\}$$ 
have left and right Haar measures 
$${\rm d}k=2^{-3}\pi^{-2}\sin(\theta)\,{\rm d}\varphi\,{\rm d}\theta\,{\rm d}\psi ,\quad\  
{\rm d}a=r^{-1}\,{\rm d}r\quad\ \hbox{and}\quad\  
{\rm d}n={\rm d}_{+}z={\rm d}x\,{\rm d}y\;,\eqno(1.1.10)$$
respectively, where $x$ and $y$ are the real and imaginary parts of $z$. 
Note that the compactness of $K$ has here allowed a choice of ${\rm d}k$  
such that 
$$\int_{K}{\rm d}k=2\;.$$
The group $G$ also 
has a left and right Haar measure:
$${\rm d}g=r^{-2}\,{\rm d}n\,{\rm d}a\,{\rm d}k=r^{-3}\,{\rm d}x\,{\rm d}y\,{\rm d}r\,{\rm d}k
\quad\hbox{for}\ \,
g=n[x+iy]a[r]k\ \,\hbox{with}\  \,x,y\in{\Bbb R},\ r>0\ \,\hbox{and}\ \,k\in K\;.\eqno(1.1.11)$$

With respect to the hyperbolic Riemannian metric on ${\Bbb H}_3$, 
$${|{\rm d}z|^2+{\rm d}r^2\over r^2}={{\rm d}x^2+{\rm d}y^2+{\rm d}r^2\over r^2}\;,\eqno(1.1.12)$$
the elements of $G$ act upon ${\Bbb H}_3$ as elements of the group ${\it Iso}^{+}\bigl( {\Bbb H}_3\bigr)$
of orientation preserving isometries:  
one has in effect a homomorphism $g\mapsto g|_{{\Bbb H}_3}$ 
from $G$ into ${\it Iso}^{+}\bigl( {\Bbb H}_3\bigr)$, 
which (see [11], Proposition~1.1.3) is surjective and has kernel $\{ h[1] , h[-1]\}$. 
The hyperbolic metric (1.1.12) 
makes ${\Bbb H}_3$ a model for 
three dimensional hyperbolic space, 
and induces on ${\Bbb H}_3$ a $G$-invariant measure, 
$$r^{-3}\,{\rm d}_{+}z\,{\rm d}r=r^{-3}\,{\rm d}x\,{\rm d}y\,{\rm d}r={\rm d}Q\qquad 
\hbox{($Q=(z,r)=(x+iy,r)$ with $x,y\in{\Bbb R}$ and $r>0$\/),}\eqno(1.1.13)$$
identical to that which is induced 
(via the homeomorphism between $G/K$ and ${\Bbb H}_3$) 
by the Haar measure (1.1.11).
Let $K^{+}\subset K$ be a fundamental domain for $\{ h[1] , h[-1]\}\backslash K$; 
and ${\cal F}'$ any fundamental domain for the action of $\Gamma$ on ${\Bbb H}_3$. 
Then a fundamental domain for $\Gamma\backslash G$ is the set 
$\{ n[z]a[r] : (z,r)\in{\cal F}'\} K^{+}\subset NAK=G$.
One therefore has (using the fundamental domain ${\cal F}$ from (1.1.7)): 
$${\rm vol}(\Gamma\backslash G)=\int_{\Gamma\backslash G} {\rm d}g 
=\int_{\cal F} {\rm d}Q\,\int_{K^{+}} {\rm d}k
=\int_{\cal F} {\rm d}Q={\rm vol}\left({\cal F}\right)\;,\eqno(1.1.14)$$
where, by (1.1.5)-(1.1.7) and Theorem~7.1.1 of~[11], 
$${\rm vol}\left({\cal F}\right)
=\int_{\cal F}r^{-3}\,{\rm d}x\,{\rm d}y\,{\rm d}r
={\rm vol}\left( {\cal F}_{{\Bbb Q}(i)}\right) [SL(2,{\frak O}) : \Gamma]
=2\pi^{-2}\zeta_{{\Bbb Q}(i)}(2)
\left| q_0\right|^2\!\!\!\!\prod_{\scriptstyle {\frak O}
\supset (\varpi)\ni q_0\atop\scriptstyle (\varpi)\ {\rm is\ prime}}
\!\!\left( 1+{1\over |\varpi|^2}\right)\eqno(1.1.15)$$
with $\zeta_{{\Bbb Q}(i)}(s)=\zeta(s)L\bigl( s,\chi_4\bigr)$ being the Dedekind zeta-function for ${\Bbb Q}(i)$, so that 
$$2\pi^{-2}\zeta_{{\Bbb Q}(i)}(2)
={2\over\pi^2}\sum_{0\neq\alpha\in{\frak O}}|\alpha|^{-4} 
={1\over 3}L(2,\chi_4)
={1\over 3}\sum_{n=1}^{\infty}{(-1)^{n-1}\over (2n-1)^2} 
={1\over 3}\left( 1-\sum_{j=1}^{\infty} 2^{-4j}j\zeta(2j+1)\right)$$
(see Page~312 of~[11]). 

The actions of elements of $G$ upon ${\Bbb H}_3$ extend, by continuity, to actions upon 
${\Bbb H}_3\cup{\Bbb P}^{1}({\Bbb C})$, where 
${\Bbb P}^{1}({\Bbb C})={\Bbb C}\cup\{\infty\}$ (the Riemann sphere): 
a projective point ${\frak z}=\bigl[ z_1,z_2\bigr]\in{\Bbb P}^{1}({\Bbb C})$ being mapped 
by the action of 
$$g=\pmatrix{a&b\cr c&d}\in G$$
to the point 
$g{\frak z}=\bigl[ a z_1+b z_2,c z_1+d z_2\bigr]\in{\Bbb P}^{1}({\Bbb C})$ 
(so that $g\infty =\infty$ if and only if $c=0$\/). 
If $Q=(z,r)\in{\Bbb H}_3$, then the subgroup of elements of $\Gamma$ fixing $Q$ 
is finite. In contrast there exist points ${\frak z}\in{\Bbb P}^{1}({\Bbb C})$ for which 
the stabiliser, $\Gamma_{\frak z}=\{ g\in\Gamma : g{\frak z} ={\frak z}\}$, contains a 
free Abelian subgroup of rank 2: such points ${\frak z}$ are called `cusps' of $\Gamma$. 
Since $\Gamma$ is a congruence subgroup (that is, there exists $M\in{\frak O}-\{ 0\}$, 
namely $M=q_0$, such that $\Gamma$ contains the kernel of the natural 
homomorphism from $SL(2,{\frak O})$ into $SL(2,{\frak O}/M{\frak O})$\/), 
the set of all cusps of $\Gamma$ is simply 
${\Bbb Q}(i)\cup\{\infty\}={\Bbb P}^{1}({\Bbb Q}(i))$. 

Each point $Q\in{\Bbb H}_3$, or ${\frak z}\in{\Bbb P}^{1}({\Bbb C})$, has a 
$\Gamma$-orbit: $\Gamma Q=\{\gamma Q : \gamma\in\Gamma\}\subset{\Bbb H}_3$,  
$\Gamma{\frak z}=\{\gamma{\frak z} : \gamma\in\Gamma\}\subset{\Bbb P}^{1}({\Bbb C})$. 
For $Q\in{\Bbb H}_3$ one has $1\leq|{\cal F}\cap\Gamma Q|\ll 1$ and 
$|{\rm Int}({\cal F})\cap\Gamma Q|\leq 1$. 
For a pair of cusps ${\frak a},{\frak b}\in{\Bbb P}^{1}({\Bbb Q}(i))$, the relation 
${\frak a}\sim^{\!\!\!\!\Gamma}{\frak b}$ (`$\Gamma$-equivalence') is deemed to 
hold if and only if $\Gamma{\frak a}=\Gamma{\frak b}$. This is an equivalence relation
under which ${\Bbb P}^{1}({\Bbb Q}(i))$ is partitioned into a finite number 
of distinct equivalence classes, ${\cal P}_1,\ldots ,{\cal P}_{H(\Gamma)}$. 

Let ${\frak c}$ be a cusp for $\Gamma$. Since $G$ acts transitively (even 3-transitively) 
on ${\Bbb P}^{1}({\Bbb C})$, one may choose $g_{{\frak c}}\in G$  such that
$$g_{{\frak c}}\infty ={\frak c}\;.\eqno(1.1.16)$$
One then has 
$$\Gamma_{{\frak c}}=\Gamma\cap g_{{\frak c}}Pg_{{\frak c}}^{-1}\;,\eqno(1.1.17)$$
where
$$P=\{ g\in G : g\infty =\infty\}=\left\{\pmatrix{u&v\cr 0&u^{-1}} : u\in{\Bbb C}^{*},\ v\in{\Bbb C}\right\}\;.\eqno(1.1.18)$$
The maximal free abelian subgroup $\Gamma_{{\frak c}}'\leq\Gamma_{{\frak c}}$ consists of the 
identity element and 
unipotent elements (non-identity elements with trace equal to 2). By (1.1.17), 
$$\Gamma_{{\frak c}}'=\Gamma\cap g_{{\frak c}}Ng_{{\frak c}}^{-1}\;.\eqno(1.1.19)$$
Note that $\Gamma_{\frak c}'$ is a normal subgroup of $\Gamma_{\frak c}\,$: 
for if $\gamma\in\Gamma_{\frak c}'$ and $\eta\in\Gamma_{\frak c}$ then 
${\rm Tr}\bigl(\eta\gamma\eta^{-1}\bigr) ={\rm Tr}(\gamma)=2$.

As is shown in Lemma 4.2 of this paper, the above `scaling matrix' $g_{{\frak c}}\in G$ may be chosen 
so that one has both (1.1.16) and 
$$g_{{\frak c}}^{-1}\Gamma_{{\frak c}}' g_{{\frak c}}
=g_{{\frak c}}^{-1}\Gamma g_{{\frak c}}\cap N 
=B^{+}\;,\eqno(1.1.20)$$
where
$$B^{+}=\left\{ n[\alpha] : \alpha\in{\frak O}\right\}
=SL(2,{\frak O})_{\infty}'\;.\eqno(1.1.21)$$
Such choice of $g_{{\frak c}}$
simplifies Fourier expansions at cusps: see (1.4.1)-(1.4.3) below. 
It is therefore to be assumed throughout this paper that  
one works with a choice of scaling matrices such that  
(1.1.16)-(1.1.21) hold for all cusps ${\frak c}$ of $\Gamma$.

For all cusps ${\frak c}$ of $\Gamma$ one has 
$\left[\Gamma_{\frak c}:\Gamma_{\frak c}'\right]\in\{ 2,4\}$ and, 
by appropriate choice of scaling matrix $g_{\frak c}$ (satisfying (1.1.16) and (1.1.20)-(1.1.21)),  
one may ensure that 
$g_{{\frak c}}^{-1}\Gamma_{\frak c} g_{{\frak c}}|_{{\Bbb P}^{1}({\Bbb C})-\{\infty\}}$ 
has as a fundamental domain the set 
$${\cal R}_{\frak c}
=\cases{\left\{ z\in{\Bbb C} : |{\rm Re}(z)|\leq 1/2,\, |{\rm Im}(z)|\leq 1/2\right\} &if 
$\,\left[\Gamma_{\frak c}:\Gamma_{\frak c}'\right] =2\;$,\cr
&\hbox{\quad}\cr 
\left\{ z\in{\Bbb C} : |{\rm Re}(z)|\leq 1/2,\, 0\leq{\rm Im}(z)\leq 1/2\right\} &if 
$\,\left[\Gamma_{\frak c}:\Gamma_{\frak c}'\right] =4\;$.}\eqno(1.1.22)$$
Any cusp ${\frak c}=u/w\in{\Bbb Q}(i)$ (with $u,w\in{\frak O}$, $w\neq 0$ and 
$(u,w)\sim 1$) has a `width' $|m_{\frak c}|^2$, 
where $m_{\frak c}\sim q_0/\!\left( w^2 , q_0\right)$. 
By defining $m_{\infty}\sim 1$ one ensures that each pair of $\Gamma$-equivalent cusps 
${\frak a},{\frak b}$ has $m_{\frak a}\sim m_{\frak b}$, and equal widths. 

For a suitable set of representatives ${\frak C}(\Gamma)\subset{\Bbb Q}(i)\cup\{\infty\}$
of the $\Gamma$-equivalence classes of cusps, and suitably chosen scaling matrices 
$g_{\frak c}$ (${\frak c}\in{\frak C}(\Gamma)$) satisfying (1.1.16) and (1.1.20)-(1.1.21), 
the sets 
$${\cal E}_{\frak c}=\left\{ g_{\frak c}Q : Q=(z,r)\in{\Bbb H}_3, 
z\in{\cal R}_{\frak c}\ {\rm and}\ r>1/\left| m_{\frak c}\right|\right\} 
\qquad\qquad\hbox{(${\frak c}\in {\frak C}(\Gamma)$)}\eqno(1.1.23)$$
are pairwise disjoint non-compact subsets of ${\Bbb H}_3$ and, for some  
compact hyperbolic polyhedron ${\cal D}\subset {\Bbb H}_3$ (having  
finitely many faces and a connected interior), the union of sets  
$${\cal F}_{*}
={\cal D}\cup\bigcup_{{\frak c}\in{\frak C}(\Gamma)} {\cal E}_{\frak c}\eqno(1.1.24)$$
is a fundamental domain for the action of $\Gamma$ on ${\Bbb H}_3$, 
with a connected interior,  
${\rm Int}({\cal F}_{*})\,$ (it simultaneously being the case that, 
for ${\frak c}\in{\frak C}(\Gamma)$, one has 
${\rm Int}({\cal D})\cap\left(\overline{{\cal E}_{\frak c}}-{\rm Int}\bigl( 
{\cal E}_{\frak c}\bigr)\right) ={\cal E}_{\frak c}'-{\cal E}_{\frak c}^{*}$, where 
the set ${\cal E}_{\frak c}^{*}$ is finite  
and ${\cal E}_{\frak c}'=\left\{ g_{\frak c}Q : 
Q=\left( z , 1/| m_{\frak c}|\right)\in{\Bbb H}_3\, ,\,  
z\in{\rm Int}\bigl( {\cal R}_{\frak c}\bigr) 
\right\}\,$).  
We make use of this type of fundamental domain in the appendix to this paper:  
see the proof of Corollary~6.2.10, and the proof of 
Lemma~6.5.16. 

\bigskip

\centerline{\bf \S 1.2 $\Gamma$-automorphic functions, Casimir operators and the 
Laplacian for ${\Bbb H}_3$.}

\medskip

A function $f : G \rightarrow {\Bbb C}$ is said to be $\Gamma$-automorphic if and only if it satisfies 
$$f(\gamma g)=f(g)\quad\hbox{for $\gamma\in\Gamma$ and $g\in G\,$.}\eqno(1.2.1)$$
Since $\Gamma\ni h[-1]$, the $\Gamma$-automorphic functions $f$ are 
even  (i.e. for $g\in G$ they satisfy $f(h[-1]g)=f(g)$\/).  
The `square integrable' $\Gamma$-automorphic functions $f : G\rightarrow{\Bbb C}$ 
are those satisfying  
$\langle f,f\rangle_{\Gamma\backslash G}<\infty$, where
$$\langle f,h\rangle_{\Gamma\backslash G}=\int_{\Gamma\backslash G}f(g)\overline{h(g)}{\rm d}g
=\int_{\cal F}\int_{K^{+}}f(n[z]a[r]k)\overline{h(n[z]a[r]k)}
\,{\rm d}k\,r^{-3}{\rm d}_{+}z\,{\rm d}r\;.\eqno(1.2.2)$$
The space $L^2(\Gamma\backslash G)$ of all square integrable 
$\Gamma$-automorphic functions is (if one does not 
discriminate between functions that are equal almost everywhere) 
a Hilbert space with 
respect to the inner product in~(1.2.2). 
For $f\in L^2(\Gamma\backslash G)$, the norm 
$\| f\|_{\Gamma\backslash G}$ of $f$ 
is given by $\| f\|_{\Gamma\backslash G}
=\sqrt{\langle f,f\rangle_{\Gamma\backslash G}}$. 
The space of all smooth $\Gamma$-automorphic functions on $G$ is 
$$C^{\infty}(\Gamma\backslash G)
=\left\{ f\in C^{\infty}(G) : f\ {\rm is}\ \hbox{$\Gamma$-automorphic}\right\}\;.\eqno(1.2.3)$$ 
By the above definitions, neither of the spaces 
$L^2(\Gamma\backslash G)$ or $C^{\infty}(\Gamma\backslash G)$
contains the other, 
and all functions contained in 
$L^2(\Gamma\backslash G)\cup C^{\infty}(\Gamma\backslash G)$ 
are measurable (i.e. measurable with respect to the Haar measure 
${\rm d}g$).  
\par 
The $\Gamma$-automorphic functions on ${\Bbb H}_3$ are those complex-valued functions that 
satisfy 
$$f(\gamma Q)=f(Q)\quad\hbox{for $\gamma\in\Gamma$ and $Q\in{\Bbb H}_3\,$;}\eqno(1.2.4)$$
and of these, those that have $\int_{\cal F}\left| f(Q)\right|^2 {\rm d}Q<\infty$ are 
the elements of the Hilbert space $L^2\left(\Gamma\backslash{\Bbb H}_3\right)$. 
The norm on 
$L^2\left(\Gamma\backslash{\Bbb H}_3\right)$ is 
given by $\| f\|_{\Gamma\backslash{\Bbb H}_3}=\sqrt{\langle f , f\rangle_{\Gamma\backslash{\Bbb H}_3}}$, where 
$$\left\langle f_1 , f_2\right\rangle_{\Gamma\backslash{\Bbb H}_3}
=\int_{\cal F} f_1(Q)\overline{f_2(Q)}\,{\rm d}Q\quad
\hbox{for $f_1 , f_2\in L^2\left(\Gamma\backslash{\Bbb H}_3\right)\,$.}\eqno(1.2.5)$$
As explained just prior to (1.2.13), below, 
the space 
$C^{\infty}\left({\Bbb H}_3\right)$ of infinitely differentiable functions 
$f : {\Bbb H}_3\mapsto{\Bbb C}$ may be viewed as 
a certain subpace of `$K$-trivial' functions contained in $C^{\infty}(G)$. 
One may similarly view $L^2\left(\Gamma\backslash{\Bbb H}_3\right)$ 
as the  subspace of $K$-trivial functions in $L^2\left(\Gamma\backslash G\right)$.

The complex Lie algebra of $G$ is ${\frak g}=\frak{sl}(2,{\Bbb C})\otimes_{\Bbb R}{\Bbb C}$, 
where $\frak{sl}(2,{\Bbb C})$ is the real vector space of all complex $2\times 2$ matrices with trace
equal to zero. The elements of ${\frak g}$ may be identified with left-invariant 
first order differential operators on $C^{\infty}(G)$ by setting
$$({\bf X}f)(g)={{\rm d}\over{\rm d}t}\,f\left( g\exp(t{\bf X})\right)\bigl|_{t=0}\quad
\hbox{for ${\bf X}\in\frak{sl}(2,{\Bbb C})$, $f\in C^{\infty}(G)$ and $g\in G$}\eqno(1.2.6)$$
(where ${\rm d}/{\rm d}t$ signifies differentiation of a function of a real variable).
Then the universal enveloping algebra, ${\cal U}({\frak g})\supset{\frak g}$, 
has centre ${\cal Z}({\frak g})={\Bbb C}\left[\Omega_{+},\Omega_{-}\right]$, where, 
in terms of the Iwasawa coordinates, one has 
$$\eqalign{\Omega_{+}
 &=F^{+}_{r,\varphi,\theta}\!\left( {\partial\over\partial z}\,,
\,{\partial\over\partial\overline{z}}\,,
\,{\partial\over\partial r}\,,\,{\partial\over\partial\varphi}\,,
\,{\partial\over\partial\theta}\,,\,{\partial\over\partial\psi}\right) =\cr 
 &=F^{+}_{r,\varphi,\theta}
\!\left( {1\over 2}\left( {\partial\over\partial x} 
-i\,{\partial\over\partial y}\right)\,,
\,{1\over 2}\left( {\partial\over\partial x}+i\,{\partial\over\partial y}\right)\,, 
\,{\partial\over\partial r}\,,\,{\partial\over\partial\varphi}\,,
\,{\partial\over\partial\theta}\,,\,{\partial\over\partial\psi}\right) =\cr 
 &={1\over 2}\,r^2\,{\partial\over\partial z}\,{\partial\over\partial\overline{z}}
+{1\over 2}\,r e^{i\varphi}\cot(\theta)\,{\partial\over\partial z}\,{\partial\over\partial \varphi}
-{1\over 2}\,ir e^{i\varphi}\,{\partial\over\partial z}\,{\partial\over\partial \theta}
-{1\over 2}\,re^{i\varphi}\csc(\theta)\,{\partial\over\partial z}\,{\partial\over\partial \psi}\ +\qquad\quad\cr
 &\qquad +{1\over 8}\,r^2\,{\partial^2\over\partial r^2} 
-{1\over 4}\,ir\,{\partial\over\partial r}\,{\partial\over\partial \varphi} 
-{1\over 8}\,\,{\partial^2\over\partial\varphi^2} -{1\over 8}\,r\,{\partial\over\partial r} 
+{1\over 4}\,i\,{\partial\over\partial \varphi}} 
\eqno(1.2.7)$$
and $\Omega_{-}=F^{-}_{r,\varphi,\theta}\left(\,{\partial /\partial \overline{z}},\,{\partial /\partial z},\,{\partial /\partial r}, 
\,{\partial /\partial \varphi},\,{\partial /\partial \theta},
\,{\partial /\partial \psi}\right)$ 
with, for $\varphi,\theta\in{\Bbb R}$ and $r>0$, each coefficient of 
$F^{-}_{r,\varphi,\theta}\in{\Bbb C}\left[ X_1,\ldots ,X_6\right]$ being 
equal to the complex-conjugate of the corresponding coefficient in the polynomial 
$F^{+}_{r,\varphi,\theta}\left( X_1,\ldots ,X_6\right)$.
A function $f\in C^{\infty}(G)$ is  
said to be a function with character $\Upsilon$ (for ${\cal Z}({\frak g})$\/) if and only if 
$$\Omega_{+}f=\Upsilon\left(\Omega_{+}\right)f\quad\hbox{and}\quad 
\Omega_{-}f=\Upsilon\left(\Omega_{-}\right)f\;.\eqno(1.2.8)$$
   
The complex Lie algebra of $K$ is ${\frak k}=\frak{su}(2)\otimes_{\Bbb R}{\Bbb C}$, 
where  $\frak{su}(2)$ is the set of skew Hermitian elements of $\frak{sl}(2,{\Bbb C})$. 
Both it and its universal enveloping algebra, ${\cal U}({\frak k})\supset{\frak k}$, 
are generated by the elements 
$${\bf H}_2=\pmatrix{i/2&0\cr 0&-i/2},\quad 
{\bf W}_1=\pmatrix{0&1/2\cr -1/2&0},\quad 
{\bf W}_2=\pmatrix{0&i/2\cr i/2&0}\;.\eqno(1.2.9)$$
The elements of ${\cal U}({\frak k})$ 
may be interpreted
(similarly to those of ${\cal U}({\frak g})$) as left-invariant differential operators: 
$$({\bf H}_2 f)(k)={{\rm d}\over{\rm d}t}\,f\left( k\exp(t{\bf H}_2)\right)\bigl|_{t=0}
=\left( {\partial f\over\partial \psi}\right)\!(k)\quad
\hbox{for $f\in C^{\infty}(K)$ and $k\in K$}\eqno(1.2.10)$$
(for example). The centre of ${\cal U}({\frak k})$ is 
${\cal Z}({\frak k})={\Bbb C}\left[\Omega_{\frak k}\right]$, where 
$$\Omega_{\frak k}
={1\over 2}\left( {\bf H}_2^2 +{\bf W}_1^2 +{\bf W}_2^2\right) 
={1\over 2}\,\csc^2(\theta)\left(\,{\partial^2\over\partial\varphi^2} 
+\,{\partial^2\over\partial\psi^2}\right) 
-\csc(\theta)\cot(\theta)\,{\partial\over\partial\varphi} 
\,{\partial\over\partial \psi}
+{1\over 2}\,\,{\partial^2\over\partial\theta^2} 
+{1\over 2}\,\cot(\theta)\,{\partial\over\partial \theta}\;.\eqno(1.2.11)$$ 

In addition to being left-invariant the Casimir operators $\Omega_{\pm}$ 
are also right-invariant:  
$$R_g\Omega_{+}f=\Omega_{+}R_g f ,\quad\hbox{and}\quad 
R_g\Omega_{-}f=\Omega_{-}R_g f\quad
\hbox{for $g\in G$ and $f\in C^{\infty}(G)\,$,}\eqno(1.2.12)$$
where $R_g$ (the right-translation operator) maps $f$ to the function $R_g f\in C^{\infty}(G)$ 
such that $\left( R_g f\right)(h)=f(hg)$ for $h\in G$. This can be proved 
by the methods of [6], Proposition~2.2.4 and Lemma~2.2.2, where an analogous 
result concerning the  
centre of ${\cal U}\left({\frak gl}(n,{\Bbb R})\right)$ is obtained. 
One can show likewise that $\Omega_{\frak k}$ is 
an operator on $C^{\infty}(K)$ that is invariant with respect to  
right-translation by any element of $K$. 

By (1.2.12) whenever 
$f\in C^{\infty}(G/K)=\left\{ \phi\in C^{\infty}(G) : \phi(gk)=\phi(g) 
\ {\rm for\ all}\ k\in K, g\in G\right\}$ 
(the space of `$K$-trivial' elements of $C^{\infty}(G)$\/)  
one will then also have $\Omega_{\pm} f\in C^{\infty}(G/K)$. 
Therefore, and by virtue of the natural bijection
(induced by the homeomorphism between $G/K$ and ${\Bbb H}_3$ 
described below (1.1.3)) from $C^{\infty}(G/K)$ onto  
$C^{\infty}\left({\Bbb H}_3\right)$, one may view   
$\Omega_{+}|_{C^{\infty}(G/K)}$ and $\Omega_{-}|_{C^{\infty}(G/K)}$ as 
being operators from $C^{\infty}\left({\Bbb H}_3\right)$ 
into $C^{\infty}\left({\Bbb H}_3\right)$. 
The hyperbolic Laplacian operator on $C^{\infty}(G/K)$, 
or (equivalently) on $C^{\infty}\left({\Bbb H}_3\right)$,  is then 
$${\bf\Delta} =4\left(\Omega_{+}+\Omega_{-}\right)\bigl|_{C^{\infty}(G/K)}\ 
=r^{2}\left( {\partial^2\over\partial x^2}+{\partial^2\over\partial y^2}
+{\partial^2\over\partial r^2}\right) 
-r\,{\partial\over\partial r}\;,\eqno(1.2.13)$$
where $x,y,r$ signify real-valued coordinates of a point $(x+iy,r)\in{\Bbb H}_3$. 
This Laplacian inherits, from the Casimir operators, 
left-invariance with respect to the actions of elements of $G$, so that 
$$({\bf\Delta} f) \circ (g|_{{\Bbb H}_3}) 
={\bf\Delta}\left( f\circ (g|_{{\Bbb H}_3})\right)\quad
\hbox{for $f\in C^{\infty}\left({\Bbb H}_3\right)$ and $g\in G\,$.}\eqno(1.2.14)$$
Let 
$$C^{\infty}\left(\Gamma\backslash{\Bbb H}_3\right)
=\left\{ f\in C^{\infty}\left({\Bbb H}_3\right) :  f\ \hbox{is $\Gamma$-automorphic}\right\}\;.\eqno(1.2.15)$$
Then by the restriction of (1.2.14) to $g\in\Gamma$, one has  
$${\bf\Delta} f\in C^{\infty}\left(\Gamma\backslash{\Bbb H}_3\right)\quad{\rm if}\quad 
f\in C^{\infty}\left(\Gamma\backslash{\Bbb H}_3\right) .\eqno(1.2.16)$$
Similarly, by the left-invariance of the elements of ${\frak g}$  
(viewed as differential operators on $C^{\infty}(G)$),  
$$\Psi f\in C^{\infty}\left(\Gamma\backslash G\right)\quad{\rm if}\quad 
\Psi\in{\cal U}({\frak g})\quad\hbox{and}\quad f\in C^{\infty}\left(\Gamma\backslash G\right) .\eqno(1.2.17)$$

If $f ,\phi\in L^2\left(\Gamma\backslash{\Bbb H}_3\right)\cap C^{\infty}\left({\Bbb H}_3\right)$ 
are such that ${\bf\Delta} f , {\bf\Delta}\phi\in L^2\left(\Gamma\backslash{\Bbb H}_3\right)$, then  
$$\langle -{\bf\Delta} f , \phi\rangle_{\Gamma\backslash{\Bbb H}_3} 
=\int_{\cal \phi}\left( {\partial f\over\partial x}\overline{{\partial \phi \over\partial x}}
+{\partial f\over\partial y}\overline{{\partial \phi \over\partial y}}
+{\partial f\over\partial r}\overline{{\partial \phi \over\partial r}}\right) 
r^2 {\rm d}Q\;,\eqno(1.2.18)$$
so that  
when operating on functions satisfying the  above constraints $-{\bf\Delta}$ is both symmetric and positive: 
$$\langle -{\bf\Delta} f , \phi\rangle_{\Gamma\backslash{\Bbb H}_3} 
=\langle f , -{\bf\Delta}\phi\rangle_{\Gamma\backslash{\Bbb H}_3}\;,\eqno(1.2.19)$$
and 
$$\langle -{\bf\Delta} f , f\rangle_{\Gamma\backslash{\Bbb H}_3}>0\quad
\hbox{if $f$ is non-constant.}\eqno(1.2.20)$$
These results (1.2.18)-(1.2.20) are contained in Theorem~4.1.7 of~[11], 
where the (effectively) more general case of an arbitrary discrete subgoup $\Gamma\leq PSL(2,{\Bbb C})$ 
is treated, and where it is moreover shown 
that for an appropriately extended domain of definition   
the operator $-{\bf\Delta}$ is essentially self-adjoint.  
As for the Casimir operator $\Omega_{\frak k}\in{\cal Z}({\frak k})$, 
one has by [42], Chapter~2, Equation~(6.3), the result that 
$$\left(\Omega_{\frak k}F , \Phi\right)_K
=\left( F , \Omega_{\frak k}\Phi\right)_K\quad
\hbox{for $F,\Phi\in C^{\infty}(K)\,$,}\eqno(1.2.21)$$
where 
$$\left( F_1 , F_2\right)_K
=\int_K F_1 (k) \overline{F_2 (k)}\,{\rm d}k\;.\eqno(1.2.22)$$

\bigskip

\centerline{\bf \S 1.3 Functions of $K$-type $(\ell,q)$.}

\medskip

The differential operators ${\bf H}_2$ and $\Omega_{\frak k}$, as initially defined, 
have domain $C^{\infty}(K)$.  
One extends ${\bf H}_2$ to the domain $C^{\infty}(G)$ by 
putting $\left({\bf H}_2 f\right) (g)=\left({\bf H}_2 f_{g}\right) (k[1,0])$  for $g\in G-K$, 
where each function $f_{g} : K\rightarrow{\Bbb C}$ is given by $f_{g}(k)=f(gk)$.  
The corresponding extension of any other elements of ${\cal U}({\frak k})$ 
to the domain $C^{\infty}(G)$ (that of $\Omega_{\frak k}$ in particular) is defined similarly.
In what follows ${\bf H}_2$ and $\Omega_{\frak k}$ may denote either the 
extensions just defined, or their restrictions to the domain $C^{\infty}(K)$ (i.e. 
the differential operators ${\bf H}_2,\Omega_{\frak k}$ as initially defined 
in (1.2.10), (1.2.11)): in each instance the reader should infer the 
option suited to the relevant operand. 

Assuming $\ell\geq -1/2$ and $q\in{\Bbb R}$, 
an element $f$ of either $C^{\infty}(K)$ or $C^{\infty}(G)$ is said to be 
`of $K$-type $(\ell ,q)$' if and only if $${\bf H}_2 f=-iqf\quad\hbox{and}\quad
\Omega_{\frak k}f=-{1\over 2}\,\left(\ell^2+\ell\right) f\;.\eqno(1.3.1)$$
If it is the case that $q\in{\Bbb Z}$ then all functions of $K$-Type $(\ell,q)$ are 
necessarily even functions (this follows by virtue of (1.2.10), (1.1.2), (1.1.8) and 
(1.1.9)). To obtain useful examples of such functions, suppose that 
$\nu\in{\Bbb C}$, $p,q,\ell\in{\Bbb Z}$, $\ell\geq |p| , |q|$ and $k=k[\alpha ,\beta]\in K$, 
and let  the function $\varphi_{\ell ,q}(\nu ,p) : G\rightarrow{\Bbb C}$ be defined by
$$\varphi_{\ell ,q}(\nu ,p)\left( n a[r] k\right) =r^{1+\nu}\Phi_{p,q}^{\ell}(k)\qquad
\hbox{($n\in N$, $r>0$, $k\in K$\/),}\eqno(1.3.2)$$
where $\Phi_{p,q}^{\ell}\left( k[\alpha ,\beta]\right)$ is the coefficient of 
$X^{\ell -p}$ in $(\alpha X-\overline{\beta})^{\ell -q} (\beta X+\overline{\alpha})^{\ell +q}\in{\Bbb C}[X]$. 
Then 
$\varphi_{\ell ,q}(\nu ,p)$ lies in the space 
$$C^{\infty}(N\backslash G)=\left\{ f\in C^{\infty}(G) : f(ng)=f(g)\ {\rm for}\ n\in N,\ g\in G\right\}$$  
and is a function of $K$-type $(\ell,q)$ with, moreover, 
character $\Upsilon=\Upsilon_{\nu,p}\,$,  
the unique character for ${\cal Z}({\frak g})$  such that 
$$\Upsilon_{\nu,p}\left(\Omega_{+}\right)
={1\over 8}\,\left( (\nu -p)^2 -1\right)\qquad{\rm and}\qquad 
\Upsilon_{\nu,p}\left(\Omega_{-}\right)
={1\over 8}\,\left( (\nu +p)^2 -1\right)\;.\eqno(1.3.3)$$
For $p,q,\ell\in{\Bbb Z}$ and $\ell\geq |p|,|q|$, the above defined function $\Phi_{p,q}^{\ell}(k)$ 
is an element of $C^{\infty}(K)$ of 
$K$-type $(\ell ,q)$; it is, moreover, an even function of $k$: for it follows directly from 
the definition that 
$$\Phi_{p,q}^{\ell}\left( h\bigl[ e^{i\varphi /2}\bigr] k h\bigl[ e^{i\psi /2}\bigr]\right)  
=e^{-ip\varphi -iq\psi} \Phi_{p,q}^{\ell}(k)\qquad\hbox{($\varphi,\psi\in{\Bbb R}$ and $k\in K$\/).}\eqno(1.3.4)$$
The set 
${\cal L}_{\rm even}(K)=\left\{\Phi_{p,q}^{\ell} : \hbox{$p,q,\ell\in{\Bbb Z}$ and $\ell\geq |p|,|q|$}\right\}$ 
is orthogonal with respect to the inner-product defined in (1.2.22);  
it spans a dense 
subspace of the Hilbert space, $L^2_{\rm even}(K)$, of even functions $f : K\rightarrow{\Bbb C}$ 
such that $\int_K |f|^2 {\rm d}k <\infty$. 
By (1.2.21) and the points since noted (excluding (1.3.3)) one may deduce that 
if $f\neq 0$ is an even function of $K$-type $(\ell, q)$ then $\ell,q\in{\Bbb Z}$,  
$\ell\geq |q|$ and there exist $K$-trivial functions  
$h_{f,p'}\in C^{\infty}(G/K)$ ($p'=-\ell,\ldots ,\ell$)  
such that 
$$f(nak)=\sum_{p'=-\ell}^{\ell} h_{f,p'}(na) \Phi_{p',q}^{\ell}(k)\quad
\hbox{for $n\in N$, $a\in A$ and $k\in K\,$.}\eqno(1.3.5)$$

\medskip 

\centerline{\bf \S 1.4 Fourier expansions at cusps; spaces of cusp forms.}

\medskip

A function $f : G\rightarrow{\Bbb C}$ is said to be of uniform polynomial growth 
along $A$ if and only if there exist real numbers $b\geq 0$ and $r_0\geq 1$ such that 
$$\left| f(n a[r] k)\right|\leq r^b\quad\hbox{for $n\in N$, $k\in K$ and $r\geq r_0\;$.}$$
\par 
A $\Gamma$-automorpic function $f : G\rightarrow{\Bbb C}$ is said to have 
polynomial growth if and only if it is that case that, for all cusps ${\frak c}$ of $\Gamma$, 
the function $(f|{\frak c}) : G\rightarrow{\Bbb C}$ given by 
$(f|{\frak c})(g)=f\left( g_{\frak c}g\right)$ 
is of uniform polynomial growth along~$A$.
\par 
If~${\frak c}$ is a cusp of $\Gamma$, and if $f : G\rightarrow{\Bbb C}$ is $\Gamma$-automorphic, then,  
by (1.1.19)-(1.1.21) and (1.2.1), the above defined function $(f|{\frak c})$ satisfies 
$(f|{\frak c})(n[\alpha]g)=(f|{\frak c})(g)$ for all $\alpha\in{\frak O}$ and all $g\in G$. 
Hence, for any $f\in C^{\infty}(\Gamma\backslash G)$, 
one has the Fourier expansion at ${\frak c}\,$: 
$$f\left( g_{\frak c} g\right) 
=(f|{\frak c})( n[0]g)
=\sum_{\omega\in{\frak O}} \left( F^{\frak c}_{\omega} f\right)(g)\qquad\qquad 
\hbox{($g\in G$),}\eqno(1.4.1)$$
where $( F^{\frak c}_{\omega} f)(g)$, the `Fourier term of order $\omega$ for 
$f$ at ${\frak c}\,$', is given by 
$$\left( F^{\frak c}_{\omega} f\right)(g)
=\int_{B^{+}\backslash N}\left(\psi_{\omega}(n)\right)^{-1} f\left( g_{\frak c} ng\right) {\rm d}n\eqno(1.4.2)$$
with $B^{+}$ as in (1.1.20)-(1.1.21), and 
$$\psi_{\omega}(n[z])={\rm e}\left({\rm Re}(\omega z)\right)\quad 
\hbox{for $\omega\in{\frak O}$, $z\in{\Bbb C}\,$}\eqno(1.4.3)$$
(it being henceforth understood that ${\rm e}(\beta)=e^{2\pi i\beta}$ for $\beta\in{\Bbb C}$). 
For $\ell ,q\in{\Bbb Z}$, with $\ell\geq |q|$, and  any character $\Upsilon$ 
of ${\cal Z}({\frak g})$, let 
$$A_{\Gamma}(\Upsilon;\ell,q)
=\left\{ f\in C^{\infty}(\Gamma\backslash G) : 
f\ \hbox{is a of $K$-type $(\ell ,q)$ with character $\Upsilon$}\right\}\;.\eqno(1.4.4)$$
Each of these spaces $A_{\Gamma}(\Upsilon;\ell,q)$ has the subspaces 
$$A_{\Gamma}^{\rm pol}(\Upsilon;\ell,q)
=\left\{ f\in A_{\Gamma}(\Upsilon;\ell,q) : f\ \hbox{has polynomial growth}\right\}\eqno(1.4.5)$$
and 
$$A_{\Gamma}^{0}(\Upsilon;\ell,q)
=\left\{ f\in A_{\Gamma}^{\rm pol}(\Upsilon;\ell,q) : 
\left( F^{\frak c}_{0} f\right)(g)=0\ \hbox{for all cusps ${\frak c}$ of $\Gamma$ and 
all $g\in G$}\right\}\;.\eqno(1.4.6)$$
The latter of these, $A_{\Gamma}^{0}(\Upsilon;\ell,q)$, is the 
space of all cusp forms of $K$-type $(\ell,q)$ with character $\Upsilon$. 

By (1.4.2)-(1.4.3), the operator $F^{\frak c}_{\omega}$ 
maps each 
$f$ in the space $C^{\infty}(\Gamma\backslash G)$ to an even function $F^{\frak c}_{\omega}f$ in the space
$$C^{\infty}(N\backslash G,\omega)
=\left\{ h\in C^{\infty}(G) : h(ng)=\psi_{\omega}(n)h(g)\ \hbox{for}\ n\in N, g\in G\right\}\eqno(1.4.7)$$
and commutes with the actions (as differential operators upon those spaces)  
of the elements of ${\cal U}({\frak g})$.  
Consequently $F^{\frak c}_{\omega}$ maps functions $f\in A_{\Gamma}(\Upsilon;\ell,q)$ 
to functions $F^{\frak c}_{\omega} f$ lying in the complex vector space 
$$W_{\omega}(\Upsilon;\ell,q)
=\left\{ h\in C^{\infty}(N\backslash G,\omega) : 
\hbox{$h$ is of $K$-type $(\ell,q)$ with character $\Upsilon$}\right\}\;.\eqno(1.4.8)$$
If $f\in C^{\infty}(\Gamma\backslash G)$ has polynomial growth, then by (1.4.2)-(1.4.3) 
the function $F^{\frak c}_{\omega} f\in C^{\infty}(N\backslash G , \omega )$ is of 
uniform polynomial growth along $A$. Therefore if 
$f\in A_{\Gamma}^{\rm pol}(\Upsilon;\ell,q)$, then  
$F^{\frak c}_{\omega} f$ lies in the complex vector space
$$W_{\omega}^{\rm pol}(\Upsilon;\ell,q)
=\left\{ h\in W_{\omega}(\Upsilon;\ell,q) : 
\hbox{$h$ is of uniform polynomial growth along $A$}\right\}\;.\eqno(1.4.9)$$

Given the restriction to $K$-types $(\ell,q)$ with $\ell,q\in{\Bbb Z}$ (and $\ell\geq q$\/) 
it follows by (1.1.8), (1.2.10), (1.3.1) and (1.4.8) that 
all functions in the space $W_{\omega}(\Upsilon;\ell,q)$ are even. 
In [5], Lemma~4.1, Lemma~4.2 and Lemma~4.3, the following is proved. 
If $W_{\omega}(\Upsilon;\ell,q)\neq\{ 0\}$, then there exists  
$(\nu,p)\in{\Bbb C}\times{\Bbb Z}$ with $|p|\leq\ell$ such that 
$\Upsilon=\Upsilon_{\nu,p}=\Upsilon_{-\nu ,-p}$
(the character of ${\cal Z}({\frak g})$ given by (1.3.3)); 
furthermore, for any $(\nu,p)\in{\Bbb C}\times{\Bbb Z}$ 
and all $\ell,q\in{\Bbb Z}$ with $\ell\geq |p|,|q|$, one has 
$${\rm dim}_{\Bbb C}\,W_0^{\rm pol}\!\left(\Upsilon_{\nu,p};\ell,q\right)
={\rm dim}_{\Bbb C}\,W_0\!\left(\Upsilon_{\nu,p};\ell,q\right)=2\eqno(1.4.10)$$
and, when $0\neq\omega\in{\frak O}$, 
$${\rm dim}_{\Bbb C}\,W_{\omega}^{\rm pol}\!\left(\Upsilon_{\nu,p};\ell,q\right)\leq 1\quad 
{\rm and}\quad 
{\rm dim}_{\Bbb C}\,W_{\omega}\!\left(\Upsilon_{\nu,p};\ell,q\right)\leq 2\;,\eqno(1.4.11)$$
with any generator, $h$, of $W_{\omega}^{\rm pol}\left(\Upsilon_{\nu,p};\ell,q\right)$ 
necessarily satisfying 
$$h(n a[r] k)\ll_h r^{\ell+1/2}\exp(-2\pi|\omega|r)\quad\hbox{for}\ n\in N, k\in K\ {\rm and}\ r\geq 1\;.\eqno(1.4.12)$$

Given (1.4.1) and the remark ending with (1.4.8), the above shows too that 
$A_{\Gamma}(\Upsilon;\ell,q)=\{ 0\}$ 
unless it is the case that, for some $\nu\in{\Bbb C}$ and some integer $p\in[-\ell,\ell]$, one has 
$\Upsilon =\Upsilon_{\nu , p}$ (i.e. as in (1.3.3)): one then designates  $(\nu , p)$ as the 
`spectral parameters' of the space $A_{\Gamma}(\Upsilon;\ell,q)$ and its elements. 
Moreover, as is shown in Lemma~5.2.1 of [32], it follows from (1.4.12) and the remark ending 
with (1.4.9) that any cusp form 
$f\in A_{\Gamma}^0\left(\Upsilon_{\nu,p};\ell,q\right)$ must satisfy, for each cusp 
${\frak c}$ of $\Gamma$, 
$$f\left( g_{\frak c}n a[r] k\right)\ll_{f,{\frak c}}r^{\ell +1/2}\exp(-\pi r)\quad 
\hbox{for}\ n\in N, k\in K\ {\rm and}\ r\geq 1\eqno(1.4.13)$$
(note that the implicit constant here does not depend also upon $g_{\frak c}$, since any alternative to $g_{\frak c}$ 
permitted by (1.1.16) and (1.1.20)-(1.1.21) 
must lie in some coset $g_{\frak c}h[\eta]N$ with $\eta^8=1$). 
By (1.4.13), (1.1.5)-(1.1.7), (1.1.14), (1.1.15) and (1.2.2)  all
cusp forms lie not only in $C^{\infty}(\Gamma\backslash G)$, but 
also in $L^2(\Gamma\backslash G)$ (as do all constant functions $f : G\rightarrow{\Bbb C}$\/). 
Another implication of (1.4.13) is that if $f : G\rightarrow{\Bbb C}$ is a non-zero constant function, 
then $f$ is not a cusp form.

By (1.1.8), (1.1.9), (1.2.10), (1.2.11) and (1.3.1), 
the space $A_{\Gamma}(\Upsilon;\ell,q)$ may contain non-zero $K$-trivial functions 
only if  $\ell =q=0$, which (by the first  observation of the previous paragraph) 
requires that $\Upsilon =\Upsilon_{\nu ,p}$ for some $(\nu,p)\in{\Bbb C}\times\{ 0\}$. 
Conversely, by (1.3.5), all even functions $f\in C^{\infty}(G)$ of $K$-type  
$(0,0)$ are $K$-trivial. Suppose now that $f\in A_{\Gamma}^{0}\left(\Upsilon_{\nu ,0};0,0\right)$
(i.e. that $f$ is some $K$-trivial cusp form). Then, by the remarks preceding (1.2.13), 
one may treat $f$ as an element of $L^2\left(\Gamma\backslash{\Bbb H}_3\right)$, and 
hence make (1.2.14), (1.2.16), (1.2.17) and (1.2.18)-(1.2.20) applicable to $f$. By (1.2.13) in particular, 
and (1.4.4), (1.4.5), (1.4.6), (1.2.8) and 
(1.3.3), one has, for $\nu\in{\Bbb C}$ and $f\in A_{\Gamma}^{0}\left(\Upsilon_{\nu ,0};0,0\right)$, 
$$-{\bf\Delta} f =\lambda_{\nu} f\quad\hbox{with}\quad \lambda_{\nu}=1-\nu^2\;.\eqno(1.4.14)$$
By (1.2.20) one must have $\lambda_{\nu}>0$ in (1.4.14) when $f$ is non-constant: 
given that cusp forms (excepting $0$) are non-constant, this shows that
$A_{\Gamma}^{0}\left(\Upsilon_{\nu ,0};0,0\right)\neq\{0\}$ only if $\nu^2<1$.
Consequently two cases may be distinguished: that of the `principal series', in which $\nu^2\leq 0\,$ 
(so that $\lambda_{\nu}\geq 1$\/); 
and that of the `complementary series', where one has $0<\nu^2<1$ (so that $0<\lambda_{\nu}<1$\/). 
If the generalised Selberg conjecture is true, then the complementary series 
is absent when (as is the case here) the relevant discrete group $\Gamma$ is 
a congruence subgroup of $SL(2,{\frak O})$. Though this conjecture remains open,   
the work [26] and [25], Theorem~4.10, of Kim and Shahidi has shown that 
$$A_{\Gamma}^{0}\left(\Upsilon_{\nu ,0};0,0\right)\neq 0\quad\hbox{only if}\quad 
\nu^2\leq (2/9)^2\;,\eqno(1.4.15)$$
so that in (1.4.14) one always has either $\lambda_{\nu}=0$ or $\lambda_{\nu}\geq 77/81$. 
By (1.4.14) and (1.2.19) (and (1.2.5), (1.2.2) and (1.1.14)), 
the spaces 
$A_{\Gamma}^{0}\left(\Upsilon_{\nu ,0};0,0\right), 
A_{\Gamma}^{0}\left(\Upsilon_{\xi ,0};0,0\right)\subset L^2\left(\Gamma\backslash G\right)$ 
are mutually orthogonal if $\xi^2\neq\nu^2$.
Similarly, since $-{\bf\Delta}\varphi_{0,0}(-1,0)=0\varphi_{0,0}(-1,0)$  (whereas  
$\lambda_{\nu}>0$ in (1.4.14) if $f\in A_{\Gamma}^{0}\left(\Upsilon_{\nu ,0};0,0\right) -\{ 0\}$\/) any non-zero $K$-trivial cusp form is orthogonal to the
subspace ${\Bbb C}\varphi_{0,0}(-1,0)\subset L^2\left(\Gamma\backslash G\right)$ 
of constant functions; the same conclusion is true in respect of 
all cusp forms ($K$-trivial or not): for if there is a cusp form $f$ of $K$-type $(\ell,q)$ that is not 
$K$-trivial, then $\ell$ must be a positive integer, and so the conclusion that  
$\left\langle f , \varphi_{0,0}(-1,0)\right\rangle_{\Gamma\backslash G}=0$ follows by  
(1.2.2), (1.1.14) and both (1.3.5)  and the observations on
${\cal L}_{\rm even}(K)$ preceding it.

\bigskip

\goodbreak\centerline{\bf \S 1.5 The Jacquet integral; generalised Kloosterman sums and the Fourier}  
\centerline{\bf expansion of a Poincar\'{e} series; Fourier coefficients of cusp forms.}

\medskip

In the case where $f$ is a cusp form, the Fourier expansion (1.4.1) may be made 
much more explicit if one has, for each non-zero $\omega\in{\frak O}$, an explicitly defined  
non-zero function in the space $W_{\omega}^{\rm pol}\!\left(\Upsilon_{\nu,p};\ell,q\right)$. 
The key to this will be the Jacquet integral, ${\bf J}_{\omega}f : G\mapsto{\Bbb C}$, which   
for $\omega\in{\Bbb C}$ and  
functions $f\in C^{\infty}(G)$ that, for some $\sigma  >0$ and some $r_1 >0$, satisfy 
$$f(n a[r] k)\ll r^{1+\sigma}\quad\hbox{for $n\in N$, $k\in K$ and $0<r\leq r_1\,$,}\eqno(1.5.1)$$
is given by 
$$\left( {\bf J}_{\omega}f\right) (g)
=\int_N \left(\psi_{\omega}(n)\right)^{-1} f\left( k[0,-1] n g\right) {\rm d}n\quad
\hbox{for $g\in G\,$.}\eqno(1.5.2)$$
By this definition, and by~(1.2.6), it follows 
that if $\sigma$ is a positive real number, if $f\in C^{\infty}(G)$ satisfies (1.5.1), 
and if ${\bf X}\in\frak{sl}(2,{\Bbb C})$ is such that  
(1.5.1) remains valid following the substitution of ${\bf X}f$ for $f$, 
then one has 
$$\left( {\bf X}{\bf J}_{\omega} f\right)(g)=\left( {\bf J}_{\omega} {\bf X}f\right)(g)\qquad\qquad  
\hbox{($g\in G$).}\eqno(1.5.3)$$
\par 
The Jacquet integral arises naturally, via the Bruhat decomposition 
$\Gamma=\Gamma_{\infty}\sqcup\left(\Gamma\cap P k[0,-1] N\right)$,  
in connection with the Fourier expansions of certain Poincar\'{e} series   
(see Section~5 of [5] for a concise sketch of the details). 
A digression on Poincar\'{e} series 
and their Fourier expansions
now facilitates the introduction of other important concepts and terminology:
discussion of the Jacquet integral 
as it relates to $W_{\omega}^{\rm pol}\!\left(\Upsilon_{\nu,p};\ell,q\right)$
resumes after that.   

Let ${\frak a}$ be a cusp of $\Gamma$. Then, for 
suitable functions $f$ lying in the space 
$$C^{\infty}(B^{+}\backslash G)=\left\{ \phi\in C^{\infty}(G) : \phi(bg)=\phi(g)
\ {\rm for}\ b\in B^{+},\ g\in G\right\}$$ 
(where $B^{+}$ is as in (1.1.20)-(1.1.21)),  one has a Poincar\'{e} series,  
$P^{\frak a} f\in C^{\infty}(\Gamma\backslash G)$, given by 
$$\left( P^{\frak a} f\right) (g)
={1\over \left[\Gamma_{\frak a} : \Gamma_{\frak a}'\right]}
\sum_{\gamma\in\Gamma_{\frak a}'\backslash\Gamma} 
f\left( g_{\frak a}^{-1}\gamma g\right)\quad\hbox{for $g\in G\,$.}\eqno(1.5.4)$$
In Subsection~6.2 (below) we describe criteria 
with the aid of which we are able to establish,  
ultimately in Subsection~6.5, that certain non-constant functions $f$ 
are `suitable' (in the above sense).

Supposing that $\omega , \omega'\in{\frak O}$, take 
$f=f_{\omega}$ to be a suitable function in the space $C^{\infty}(N\backslash G , 
\omega)\subset C^{\infty}(B^{+}\backslash G)$ 
given by (1.4.7) and (1.4.3). Then via the Bruhat decomposition of $\Gamma_{\frak a}'\backslash\Gamma$  
one can ascertain that,  in the Fourier expansion of $(P^{\frak a} f_{\omega})(g)$ at 
an arbitrary cusp ${\frak a}'$ of $\Gamma$, the Fourier term of order $\omega'$ is  
$$\eqalign{\bigl( F^{{\frak a}'}_{\omega'} P^{\frak a} f_{\omega}\bigr) (g)
 &={1\over \left[\Gamma_{\frak a} : \Gamma_{\frak a}'\right]}
\sum_{\scriptstyle\gamma\in\Gamma_{\frak a}'\backslash\Gamma\ :\ \gamma{\frak a}'={\frak a}\atop
\scriptstyle g_{\frak a}^{-1}\gamma g_{{\frak a}'}\in h\left[ u(\gamma)\right] N} 
f_{\omega}\!\left( g_{\frak a}^{-1}\gamma g_{{\frak a}'} g\right) \delta_{\omega u(\gamma) , \omega'/u(\gamma)}\ +\cr 
 &\quad\ +{1\over\left[\Gamma_{\frak a} : \Gamma_{\frak a}'\right]}
\sum^{\hbox{\quad}}_{c\in\,{}^{\frak a}{\cal C}^{{\frak a}'}} 
S_{{\frak a} , {\frak a}'}\!\left(\omega , \omega' ; c\right) 
\left( {\bf J}_{\omega'}{\bf h}_{1/c} f_{\omega}\right)\!(g)\;,}\eqno(1.5.5)$$
where
$$\delta_{\rho ,\sigma}=\cases{1 &if $\rho =\sigma$,\cr 
0 &otherwise,}\eqno(1.5.6)$$
and 
$${\bf h}_{u}f(g)=f\left( h[u]g\right)\quad\hbox{for $g\in G\,$,}\eqno(1.5.7)$$
and where, with 
$${}^{\frak a}\Gamma^{{\frak a}'}(c)
=\left\{\gamma\in\Gamma : 
g_{\frak a}^{-1}\gamma g_{{\frak a}'}=\pmatrix{*&*\cr c&*}\right\}\quad 
\hbox{for $c\in{\Bbb C}\,$,}\eqno(1.5.8)$$
one has
$${}^{\frak a}{\cal C}^{{\frak a}'}
=\left\{ c\in{\Bbb C}-\{ 0\} : {}^{\frak a}\Gamma^{{\frak a}'}(c)\neq\emptyset\right\}\eqno(1.5.9)$$
and, for $c\in {}^{\frak a}{\cal C}^{{\frak a}'}$, 
$$S_{{\frak a} , {\frak a}'}\left(\omega , \omega' ; c\right) 
=\sum_{\scriptstyle \gamma\in\Gamma_{\frak a}'\backslash 
{}^{\frak a}\Gamma^{{\frak a}'}(c)/\Gamma_{{\frak a}'}'\atop\scriptstyle 
g_{\frak a}^{-1} \gamma g_{{\frak a}'}
=\left( {\scriptstyle\!\!\!\!\!s(\gamma)\quad *\atop\scriptstyle\  c\quad\ d(\gamma)}\right) }
{\rm e}\!\left({\rm Re}\!\left(\omega\,{s(\gamma)\over c}+\omega'\,{d(\gamma)\over c}\right)\right)
 .\eqno(1.5.10)$$
Given the restricted choice of scaling matrices $g_{\frak a}, g_{{\frak a}'}$ 
(as in (1.1.16) and (1.1.20)-(1.1.21)), one can show that 
$${}^{\frak a}{\cal C}^{{\frak a}'}\ni c\quad\hbox{only if}\quad  
0\neq c^2\in{\frak O}\ \,{\rm and}\ \, 
c^2\sim (c')^2\,m_{\frak a} m_{{\frak a}'}\ \,
\hbox{for some}\ \,c'\in{\frak O}\eqno(1.5.11)$$
(where the non-zero Gaussian integers $m_{\frak c}$ are as defined just below (1.1.22)). 
For $c\in {}^{\frak a}{\cal C}^{{\frak a}'}$ and $c'$ as in (1.5.11), the `generalised Kloosterman sum' in (1.5.10) trivially satisfies 
$\left| S_{{\frak a} , {\frak a}'}\left(\omega , \omega' ; c\right)\right|\leq 
\left| m_{\frak a} m_{{\frak a}'}c'\right|^2=|c|^2 \left| m_{\frak a} m_{{\frak a}'}\right|$, 
while work of Bruggeman and Miatello in Proposition~9 and Theorem~10 of [4] 
shows, by means of an exponential sum 
estimate of A. Weil, that one has the generally non-trivial bound: 
$$\left| S_{{\frak a} , {\frak a}'}\left(\omega , \omega' ; c\right)\right|
\leq \sqrt{8}\,\left| m_{\frak a} m_{{\frak a}'}\right|^2 \left|\left( c' , q_0^{\infty}\right)
\left( c' , \omega , \omega '\right) c'\right| \tau\!\left( c'\right) ,\eqno(1.5.12)$$
where 
$\left|\left( c' , q_0^{\infty}\right)\right| =\lim_{n\rightarrow\infty}\left|\left( c' , q_0^n\right)\right|$
and $\tau\left( c'\right)$ is the number of Gaussian integer divisors of $c'$. 
The generalised Kloosterman sums also satisfy some symmetry relations: 
$$S_{{\frak a},{\frak a}'}\!\left(\omega,\omega';c\right) 
=S_{{\frak a}',{\frak a}}\!\left(-\omega',-\omega;-c\right) 
=S_{{\frak a}',{\frak a}}\!\left(-\omega',-\omega;c\right)\eqno(1.5.13)$$
(the first of these equations following from 
the observation that if $\left\{ \gamma_1,\ldots ,\gamma_n\right\}$ is a 
complete set of representatives for the set of double cosets 
$\Gamma_{\frak a}'\backslash{}^{\frak a}\Gamma^{{\frak a}'}(c)/\Gamma_{{\frak a}'}'$, then 
$\left\{ \gamma_1^{-1},\ldots ,\gamma_n^{-1}\right\}$ is a complete set of representatives for 
$\Gamma_{{\frak a}'}'\backslash{}^{{\frak a}'}\Gamma^{\frak a}(-c)/\Gamma_{\frak a}'$; 
the second following since ${}^{{\frak a}'}\Gamma^{\frak a}(-c)=h[-1] {}^{{\frak a}'}\Gamma^{\frak a}(c)$\/).
\par
Returning now to the subject of the Jacquet integral's r\^{o}le in providing an explicit 
generator for the space $W_{\omega}^{\rm pol}\left(\Upsilon_{\nu ,p};\ell,q\right)$, 
suppose that $\nu\in{\Bbb C}$ and $p,q,\ell\in{\Bbb Z}$ with $|p|,|q|\leq\ell$.
Then it can be verified that the function $\varphi_{\ell ,q}(\nu ,p) : G\rightarrow{\Bbb C}$ 
(defined in (1.3.2)) lies in the space $W_0^{\rm pol}\left(\Upsilon_{\nu ,p};\ell,q\right)$;  
and that, given (1.4.10), it is furthermore the case that  
$$W_0\left(\Upsilon_{\nu ,p};\ell,q\right)
={\Bbb C}\,\varphi_{\ell,q}(\nu,p)\oplus
\cases{{\Bbb C}\,\!\displaystyle{\partial\over\partial \nu} 
\,\varphi_{\ell,q}(\nu,0)\Bigr|_{\nu=0} &if $\nu=p=0$, \cr
\quad &\quad \cr 
{\Bbb C}\,\varphi_{\ell,q}(-\nu,-p) &otherwise.}\eqno(1.5.14)$$
In addition, since (1.5.3) implies that the Jacquet integral has, in common with 
the Fourier operators $F^{\frak c}_{\omega}$, the property of commuting 
with the actions (as differentiable operators) of all elements of 
${\cal U}({\frak g})$ 
(commuting, in particular, with $\Omega_{\pm}$ and the extensions of $\Omega_{\frak k}$ and 
${\bf H}_2$\/), 
it therefore follows by (1.3.2), (1.5.12), (1.5.1)-(1.5.2) and the alternative  representation of 
${\bf J}_{\omega}\varphi_{\ell,q}(\nu,p)$ in Equation~(5.10) of [5], that one has 
$$W_{\omega}^{\rm pol}\left(\Upsilon_{\nu ,p};\ell,q\right)
\ni{\bf J}_{\omega}\varphi_{\ell,q}(\nu,p)\quad 
\hbox{for $\omega\in{\Bbb C}$, if ${\rm Re}(\nu)>0\,$.}\eqno(1.5.15)$$

The Jacquet integral $({\bf J}_{\omega}\varphi_{\ell,q}(\nu,p))(g)$ (where $g\in G$)
fails to converge absolutely when ${\rm Re}(\nu)\leq 0$. 
It nevertheless follows from Lemma~5.1 of [5] that when 
$\omega\neq 0$ the holomorphic function 
$\nu\mapsto ({\bf J}_{\omega}\varphi_{\ell,q}(\nu,p))(g)$  
(with domain $\{\nu\in{\Bbb C} : {\rm Re}(\nu)>0\}$\/) 
has an extension via analytic continuation that is entire; 
and that, for each $\nu\in{\Bbb C}$, the resulting function 
${\bf J}_{\omega}\varphi_{\ell,q}(\nu,p) : G\rightarrow{\Bbb C}$ is an 
element of the space $W_{\omega}^{\rm pol}\left(\Upsilon_{\nu ,p};\ell,q\right)$
distinct from $0$. 
Given this extension of (1.5.15), it follows by (1.4.11) that  
$$W_{\omega}^{\rm pol}\left(\Upsilon_{\nu ,p};\ell,q\right)
={\Bbb C}\,{\bf J}_{\omega}\varphi_{\ell,q}(\nu,p)\quad\hbox{for $\omega\neq 0\,$.}\eqno(1.5.16)$$
By this and the remarks encompassing (1.4.9) it moreover follows  that 
if ${\frak c}$ is a cusp of $\Gamma$, 
$f\in A_{\Gamma}^{\rm pol}\left(\Upsilon_{\nu,p};\ell,q\right)$  
and  $0\neq\omega\in{\frak O}$ then  
$$\left(F^{\frak c}_{\omega} f\right) (g) 
=c_f^{\frak c}(\omega)\left({\bf J}_{\omega}\varphi_{\ell,q}(\nu,p)\right)(g)\quad 
\hbox{for $g\in G\,$,}\eqno(1.5.17)$$
where $c_f^{\frak c}(\omega)$  (the Fourier coefficient) is a complex number 
depending only upon $\omega$, $g_{\frak c}$ and $f$.

In the case of $({\bf J}_{0}\varphi_{\ell,q}(\nu,p))(g)$ the extension 
by analytic continuation with respect to $\nu\in{\Bbb C}$ is meromorphic on ${\Bbb C}$, 
and is given by 
$${\bf J}_{0}\varphi_{\ell,q}(\nu,p)
=\pi\,{\Gamma(\ell+1-\nu)\Gamma(|p|+\nu)\over \Gamma(\ell+1+\nu)\Gamma(|p|+1-\nu)} 
\,\varphi_{\ell,q}(-\nu,-p)\eqno(1.5.18)$$
whenever the right-hand side is defined. 

\bigskip

\goodbreak\centerline{\bf \S 1.6 The spaces $H(\nu,p)$ of $K$-finite functions; 
principal and complementary series.}

\medskip

The results in this paper concern sums involving the Fourier coefficients 
(as given by (1.5.17)) of orthogonal systems of cusp forms. A significant aid 
in describing the relevant systems of cusp forms are certain spaces 
$H(\nu,p)\subset C^{\infty}(N\backslash G)$ that have the functions  
$\varphi_{\ell,q}(\nu,p)$ as their generators. 
To motivate the definition of $H(\nu,p)$ (which follows) one may 
first observe that the $(2\ell+1)\times(2\ell+1)$ matrices 
${\bf\Phi}_{\ell}=\left(\Phi_{p,q}^{\ell}\right)$  satisfy
$${\bf\Phi}_{\ell}\left( k_1 k_2\right) ={\bf\Phi}_{\ell}\left( k_1\right){\bf\Phi}_{\ell}\left( k_2\right)\quad 
\hbox{for $k_1,k_2\in K$ and $\ell =0,1,2,\ldots\ \;$}$$
(a slight change in the definition of ${\bf\Phi}_{\ell}$, 
involving  a normalisation of the entries $\Phi_{p,q}^{\ell}$, would in fact transform 
it into a realisation of the unitary representation of degree $2\ell +1$ of $K$).
The functions $\varphi_{\ell ,q}(\nu ,p) : G\mapsto{\Bbb C}$ are therefore $K$-finite (each set of `$K$-translates' 
$\{g\mapsto\varphi_{\ell ,q}(\nu ,p)(gk) : k\in{\Bbb K}\}$ being contained in the span of the finite set 
$\{\varphi_{\ell ,q'}(\nu ,p) : q'\in{\Bbb Z}\ {\rm and}\ |q'|\leq\ell\}$\/). 
Therefore, for each of the pairs 
$(\nu ,p)\in ((i{\Bbb R})\times{\Bbb Z})\cup((-1,1)\times\{ 0\})\,$ (in particular), 
one has the space
$$H(\nu ,p)=\left\{ {\rm finite\ linear\ combinations\ of\ functions\ }\varphi_{\ell ,q}(\nu ,p)
{\rm\ with\ } \ell ,q\in{\Bbb Z}{\rm \ and\ }\ell\geq |q|,|p|\right\}\;,\eqno(1.6.1)$$ 
which contains only $K$-finite functions. Though $H(\nu ,p)$ is invariant under 
the action of $K$ by right translation, 
it is not so under the action of $G$ by right translation. 
Nevertheless, these spaces $H(\nu ,p)$ are ${\frak g}$-invariant,  
and each constitutes an irreducible representation space for ${\frak g}$. 

If one excludes the case in which $\omega =0$ and $\nu=p=0$, then 
the extension by analytic continuation of the Jacquet integral 
$({\bf J}_{\omega}\varphi_{\ell,q}(\nu,p))(g)$ 
(see between (1.5.15) and (1.5.18)) permits a unique further extension to a linear operator, 
$${\bf J}_{\omega}^{\nu,p} : H(\nu ,p)\rightarrow W_{\omega}^{\rm pol}\left(\Upsilon_{\nu,p}\right)
=\bigoplus_{\ell =|p|}^{\infty}\,\bigoplus_{q=-\ell}^{\ell} 
W_{\omega}^{\rm pol}\left(\Upsilon_{\nu,p};\ell,q\right)$$
(the `Jacquet operator'): note the relevance of (1.5.14) and (1.5.18) for the case $\omega =0$. 
The form of the analytic continuation provided by Lemma~5.1 of [5] is such as 
to ensure that the Jacquet operator ${\bf J}_{\omega}^{\nu,p}$ does 
(even when ${\rm Re}(\nu)\leq 0$\/) inherit from the Jacquet integral 
$({\bf J}_{\omega}f)(g)$ 
the property of commuting with the actions of all $\Psi\in{\cal U}({\frak g})$. 

Let $H^{\infty}(\nu ,p)$ be the space of functions $f\in C^{\infty}(N\backslash G)$ 
that satisfy $f\bigl( a[r] h\bigl[e^{it}\bigr] k\bigr) =r^{1+\nu}e^{-2p it} f(k)$ 
for $r>0$, $t\in{\Bbb R}$ and $k\in K$. 
The inner-product appropriate for $H(\nu ,p)$ derives from the duality 
between the pair of spaces $H^{\infty}(\pm\nu ,\pm p)$
given by the bilinear form 
$$\langle f_{+} , f_{-}\rangle_{\natural} =\int_K f_{+}(k) f_{-}(k) {\rm d}k\quad 
\hbox{($f_{\pm}\in H^{\infty}(\pm\nu ,\pm p)$\/).}\eqno(1.6.2)$$
One has $H^{\infty}(\nu ,p)\supset H(\nu ,p)$ and 
$H^{\infty}(-\nu ,-p)\supset H(-\nu ,-p)$; and 
if $(\nu , p)\in (i{\Bbb R})\times{\Bbb Z}\,$ (which is the `unitary principal series' case), 
then $\overline{f}\in H(\overline{\nu} , -p)=H(-\nu , -p)$ when $f\in H(\nu , p)$, 
so that one may define the inner product: 
$$\left( f_1 , f_2\right)_{\rm ps}=\left\langle f_1 , \overline{f_2}\right\rangle_{\natural}\quad
\hbox{for $f_1,f_2\in H(\nu , p)$ and $(\nu , p)\in (i{\Bbb R})\times{\Bbb Z}\;$.}\eqno(1.6.3)$$
If instead $p=0$ and $0<\nu^2<1$ (the `complementary series' case), then 
one has $\overline{f}\in H(\overline{\nu} , 0)=H(\nu , 0)$ when $f\in H(\nu , p)$. 
Then, in order to make applicable the duality of (1.6.2), one passes from 
$\overline{f}\in H(\nu , 0)$ to 
$\pi^{-1}(\Gamma(1-\nu)/\Gamma(\nu)) {\bf J}_0^{\nu,0} \overline{f}\in H(-\nu ,0)$.  
Since ${\bf J}_0^{\nu,0}\overline{f}=\overline{{\bf J}_0^{\nu,0}f}$, this leads to the definition: 
$$\left( f_1 , f_2\right)_{cs}
={\Gamma(1-\nu)\over\pi\Gamma(\nu)}\,\Bigl\langle f_1 , 
\overline{{\bf J}_0^{\nu,0} f_2}\,\Bigr\rangle_{\!\natural}\quad 
\hbox{for $f_1,f_2\in H(\nu , 0)$ with $0<\nu^2<1\;$.}\eqno(1.6.4)$$

Since ${\bf J}_0^{\nu,0}$ is a linear operator that commutes with the actions of 
the elements of ${\cal U}({\frak g})$ 
upon the space $H(\nu , 0)$, and since $\langle f_{+} , f_{-}\rangle_{\natural}$ is invariant under any 
right translation by an element of $G$ (applied simultaneously to 
$f_{\pm}\in H^{\infty}(\pm\nu , 0)$\/), 
completion of $H(\nu , 0)$ with respect to the norm $\| f\|_{\rm cs}=\sqrt{(f,f)_{\rm cs}}$
yields a Hilbert space $H^2(\nu , 0)\supset H^{\infty}(\nu , 0)$ upon which $G$ acts unitarily (as $G$ also does 
in the `principal series' case, where $H^2(\nu , p)\supset H^{\infty}(\nu , p)$ is instead    
the completion of $H(\nu , p)$ with respect to the norm $\| f\|_{\rm ps}=\sqrt{(f,f)_{\rm ps}}\;$\/). 
See [32], Section~2.3, for a fuller discussion of the 
spaces $H^{\infty}(\nu ,p)$ and $H^2(\nu ,p)$. 

Focusing on the generators of $H(\nu , p)$, 
one finds by (1.6.4) and (1.5.18) that, for $p=0$ and $0<\nu^2<1$,  
$$\eqalign{\left(\varphi_{\ell ,q}(\nu , 0), \varphi_{\ell',q'}(\nu ,0)\right)_{\rm cs}
 &=\left( \Phi_{0,q}^{\ell} \,,\, 
{\Gamma(1+\ell' -\nu)\over \Gamma(1+\ell' +\nu)}\,\Phi_{0,q'}^{\ell'} \right)_{\!\!K} =\cr 
 &=\delta_{\ell,\ell'}\delta_{q,q'}\,{\Gamma(1+\ell -\nu)\over \Gamma(1+\ell +\nu)}\,
{1\over (\ell +\textstyle{1\over 2})}\,\pmatrix{2\ell\cr\ell}\pmatrix{2\ell\cr\ell -q}^{-1}\;,}  
\eqno(1.6.5)$$
while, for $(\nu ,p)\in (i{\Bbb R})\times{\Bbb Z}$ with $|p|\leq\ell$, 
$$\left(\varphi_{\ell ,q}(\nu , p), \varphi_{\ell',q'}(\nu ,p)\right)_{\rm ps}
=\bigl( \Phi_{p,q}^{\ell}\,, \Phi_{p,q'}^{\ell'} \bigr)_{\!K}  
=\delta_{\ell,\ell'}\delta_{q,q'}\,{1\over (\ell +\textstyle{1\over 2})}\,
\pmatrix{2\ell\cr\ell -p}\pmatrix{2\ell\cr\ell -q}^{-1}\;.\eqno(1.6.6)$$
Note that by (1.6.1), (1.6.5) and (1.6.6), it is evident that the inner products $(\cdot , \cdot)_{\rm ps}$ and 
$(\cdot , \cdot)_{\rm cs}$ are indeed positive definite on 
the relevant spaces $H(\nu , p)$.

Since the spaces $H(\nu , p)$ and $H^{\infty}(\nu , p)$ 
have been defined only for integer values of $p$ 
(and given their definitions, including (1.3.2), along with 
the points noted in connection with (1.3.4) and (1.3.5)) 
it is implied that these spaces contain only even functions. 
The same is therefore true of the completed spaces $H^2(\nu , p)$. 

\bigskip

\centerline{\bf \S 1.7 Decomposing the space $L^2(\Gamma\backslash G)$.}

\medskip

The utility of the spaces $H(\nu ,p)$ in studying cusp forms derives from 
their r\^{o}le in classifying the irreducible unitary representations of $SL(2,{\Bbb C})$. 
More specifically, it is known that any even non-trivial irreducible unitary 
representation of the Lie group $G=SL(2,{\Bbb C})$ is, for some 
$(\nu,p)\in ((i{\Bbb R})\times{\Bbb Z})\cup((-1,1)\times\{ 0\})$, 
unitarily equivalent to a certain representation 
of $G$ with representation space $H^2(\nu ,p)\,$; this representation,  
being (of necessity) unitary itself, is a strongly continuous homomorphism,  
$$P^{2p,2\nu} : G\rightarrow{U}\!\left( H^2(\nu ,p)\right)\;,\eqno(1.7.1)$$ 
mapping elements $g\in G$ to elements $P^{2p,2\nu}_g$ of the group ${U}\!\left( H^2(\nu ,p)\right)$
of unitary operators on $H^2(\nu ,p)\,$. For all $g\in G$, the defining property of 
the operator $P^{2p,2\nu}_g : H^2(\nu ,p)\rightarrow H^2(\nu ,p)$ is that 
$$\left( P^{2p,2\nu}_g\varphi\right) (h)=\varphi(hg)\quad 
\hbox{for $\varphi\in H^2(\nu ,p)$ and $h\in G\,$.}\eqno(1.7.2)$$
Note that this is the so-called `induced picture' of the 
classification: Theorem~16.2 of [27] gives the `non-compact' picture 
(without the restriction to even representations) and 
Section~7.1 [27] describes the `induced', `compact' and `non-compact' pictures 
in relation to one another. 

The above classification is significant in the current context, since (as already noted) 
each cusp form $f$ lies in $L^2(\Gamma\backslash G)$, and is orthogonal 
to the subspace of constant functions. Hence, on identifying the constant 
functions with elements of ${\Bbb C}$, one has: 
$$L^2(\Gamma\backslash G)
={\Bbb C}\oplus{}^{0}L^2(\Gamma\backslash G)\oplus{}^{\rm e}L^2(\Gamma\backslash G)\;,\eqno(1.7.3)$$
where ${}^{0}L^2(\Gamma\backslash G)$ is the closure of the space spanned by the set of 
all cusp forms (of arbitrary $K$-type $(\ell,q)$, and with any character $\Upsilon_{\nu,p}$\/) and 
${}^{\rm e}L^2(\Gamma\backslash G)$ is the orthogonal complement of 
${\Bbb C}\oplus{}^{0}L^2(\Gamma\backslash G)$ in $L^2(\Gamma\backslash G)$. 
The space ${}^{0}L^2(\Gamma\backslash G)$ is invariant with respect to the 
right-actions of the elements of $G$, and one has 
$${}^{0}L^2(\Gamma\backslash G) =\overline{\bigoplus V}\;,\eqno(1.7.4)$$
where the direct sum is that of countably many 
pairwise orthogonal infinite-dimensional closed subspaces~$V$, each of which is
invariant and irreducible with respect to the right-actions of the elements of $G$. 
Since all functions in the space $L^2(\Gamma\backslash G)$ are 
(by the observation following (1.2.1)) necessarily functions that are even, 
so too are all functions lying in ${}^{0}L^2(\Gamma\backslash G)$, 
or ${}^{\rm e}L^2(\Gamma\backslash G)$, or in any one of 
the above factors $V$. 

Now (1.2.2) and the unimodularity of the Haar measure on $G$ imply that  
$$\int_{\Gamma\backslash G}
f\left( g m\right) \overline{h\left( g m\right)}\,{\rm d}g
=\int_{\Gamma\backslash G}
f\left( g\right) \overline{h\left( g\right)}\,{\rm d}g
=\langle f , h\rangle_{\Gamma\backslash G}\quad 
\hbox{for $f,h\in {}^{0}L^2(\Gamma\backslash G)$ and $m\in G\,$.}\eqno(1.7.5)$$
One has therefore (for each $V$ in (1.7.4)) the non-trivial irreducible unitary 
representation 
$R^V : G\rightarrow U(V)$, which for $m\in G$ maps $m$ to
the right-action $R^V_m$ that has $\left( R^V_m f\right)(g)=f(gm)$ when $f\in V$, $g\in G\,$; 
and, as all functions in the space $V$ are even, this representation $R^V$ 
is also even (i.e. $R^V_{h[-1]}=R^V_{h[1]}$\/). 
By the discussion, around (1.7.1), (1.7.2), concerning the classification of such 
representations of $G$ it follows that, for each $V$ in (1.7.4), there exists  
$\left(\nu_V , p_V\right)\in ((i{\Bbb R})\times{\Bbb Z})\cup ((-1,1)\times\{ 0\})$ 
and a surjective linear isometry 
$\tilde T_V : H^2\left(\nu_V , p_V\right)\rightarrow V$ such that 
$$R^V_g \tilde T_V =  \tilde T_V P_g^{2p_V,2\nu_V}\quad\hbox{for $g\in G\,$.}\eqno(1.7.6)$$
The operator $T_V=\tilde T_V\bigr|_{H\left(\nu_V , p_V\right)}$ has a dense
image $V_K\subset V$ (the $K$-finite subspace of $V$). Hence (and since $T_V$ is an isometry) 
one has by (1.6.1), (1.7.4) and the relevant orthogonality relations,  (1.6.5) or (1.6.6), the decomposition: 
$$V_K
=\bigoplus_{\ell =\left| p_V\right|}^{\infty}\bigoplus_{q=-\ell}^{\ell} 
V_{K,\ell,q}\subset V\subset {}^0L^{2}(\Gamma\backslash G)\;,\eqno(1.7.7)$$
where
$$V_{K,\ell,q}={\Bbb C} T_V \varphi_{\ell,q}\left(\nu_V,p_V\right)\;.\eqno(1.7.8)$$
Though $T_V$ does not inherit from $\tilde T_V$ the property of commuting  (as in (1.7.6))
with the right-actions of all elements of $G$, it does nevertheless  follow from (1.7.6) that 
$${\bf X} T_V=T_V{\bf X}\quad\hbox{for ${\bf X}\in\frak{sl}(2,{\Bbb C})\,$.}\eqno(1.7.9)$$
Therefore $T_V$ commutes with the actions, as differential operators, 
of all elements of ${\cal U}({\frak g})$.

By (1.7.9), the operator $T_V$ shares with the Jacquet operators 
${\bf J}_{\omega}^{\nu,p}$ 
the property of preserving the $K$-type and character $\Upsilon$ of the functions on 
which it operates. 
Therefore, if one supposes now 
that $\ell,q\in{\Bbb Z}$, $\ell\geq q$ and $\ell\geq |p_V|$, then 
it follows by (1.7.7), (1.7.8), (1.5.13) and (1.4.8) 
that $T_V\varphi_{\ell ,q}\left(\nu_V,p_V\right)\in 
A_{\Gamma}\left(\Upsilon_{\nu_V,p_V};
\ell,q\right)\cap{}^{0}L^{2}(\Gamma\backslash G)$. 
This in fact implies that 
$$V_{K,\ell,q}\subseteq 
A_{\Gamma}^0\left(\Upsilon_{\nu_V,p_V};\ell,q\right)\;.\eqno(1.7.10)$$
A proof of (1.7.10) may be given along the following lines. 
Let ${\frak c}$ be any cusp of $\Gamma$. 
Then by (1.4.1)-(1.4.2) (and since 
$T_V\varphi_{\ell ,q}\left(\nu_V,p_V\right)
\in A_{\Gamma}\left(\Upsilon_{\nu_V,p_V};\ell,q\right)\subset C^{\infty}(\Gamma\backslash G)$\/), 
one has the Fourier expansion:  
$$\left( T_V\varphi_{\ell ,q}\left(\nu_V,p_V\right)\right)\left( g_{\frak c} g\right)
=\sum_{\omega\in{\frak O}}
\left( F^{\frak c}_{\omega} T_V\varphi_{\ell,q}\left(\nu_V,p_V\right)\right)
(g)\qquad\hbox{($g\in G$),}\eqno(1.7.11)$$
where, by the point noted in connection with (1.4.8),  
$F^{\frak c}_{\omega} T_V\varphi_{\ell,q}\left(\nu_V,p_V\right)
\in W_{\omega}\left(\Upsilon_{\nu_V,p_V};\ell,q\right)$ for 
$\omega\in{\frak O}$. Since 
$T_V\varphi_{\ell,q}\left(\nu_V,p_V\right)\in{}^{0}L^{2}(\Gamma\backslash G)$, each 
term in this Fourier expansion 
is necessarily square integrable 
over $g_{\frak c}^{-1}{\cal E}_{\frak c} K^{+}$, 
where ${\cal E}_{\frak c}K^{+}\subset G$ is the `cusp sector' defined by 
(1.1.14) and (1.1.22)-(1.1.23)
(to prove this one uses (1.4.7), (1.4.8) and  
the fact that the characters $\psi_{\omega}$ defined in (1.4.3) satisfy 
$\int_{B^{+}\backslash N}\psi_{\omega'}(n)\overline{\psi_{\omega}(n)} {\rm d}n=0$ 
when $\omega$ and $\omega'$ are distinct Gaussian integers). 
In respect of the particular case $\omega =0$, 
it moreover follows  that, since cusp forms span a dense subspace of 
${}^{0}L^{2}(\Gamma\backslash G)$, one must have:  
$$0=\int_{g_{\frak c}^{-1}{\cal E}_{\frak c} K^{+}} 
\left|\left( F^{\frak c}_{0} T_V\varphi_{\ell,q}\left(\nu_V,p_V\right)\right)
(g)\right|^2 {\rm d}g 
={1\over\left[\Gamma_{\frak c} : \Gamma_{\frak c}'\right]} 
\int_{|m_{\frak c}|^{-1}}^{\infty}\int_{K} 
\left|\left( F^{\frak c}_{0} T_V\varphi_{\ell,q}\left(\nu_V,p_V\right)\right)
(a[r] k)\right|^2 {\rm d}k r^{-3}{\rm d}r\;.$$
By combining these observations with the analysis of the space 
$W_0\left(\Upsilon_{\nu,p};\ell,q\right)$ 
worked out (via (1.3.5) and elements of the theory of Bessel functions) 
in the proof of Lemma~4.2 of [5], it may be deduced that 
$$\left( F^{\frak c}_{0} T_V\varphi_{\ell,q}\left(\nu_V,p_V\right)\right) (g)=0\quad\  
\hbox{for $g\in G\,$,}\eqno(1.7.12)$$
and that 
$F^{\frak c}_{\omega} T_V\varphi_{\ell,q}\left(\nu_V,p_V\right)\in 
W_{\omega}^{\rm pol}\left(\Upsilon_{\nu_V,p_V};\ell,q\right)$
for $0\neq\omega\in{\frak O}$. Therefore it follows by (1.5.16)-(1.5.17) that 
the Fourier expansion (1.7.11) has the special form:  
$$\left( T_V\varphi_{\ell,q}\left(\nu_V,p_V\right)\right)\left( g_{\frak c} g\right) 
=\sum_{0\neq\omega\in{\frak O}}
c_V^{\frak c}(\omega)\left( {\bf J}_{\omega} 
\varphi_{\ell,q}\left(\nu_V,p_V\right)\right) (g)\quad\  
\hbox{for $g\in G\,$.}\eqno(1.7.13)$$
The Fourier coefficients $c_V^{\frak c}(\omega)$ in (1.7.13) depend 
only upon $V$, $T_V$, ${\frak c}$, $g_{\frak c}$ and $\omega\,$: for 
both $F^{\frak c}_{\omega}\,T_V$ and ${\bf J}_{\omega}^{\nu_V,p_V}$ commute with 
all $\Psi\in{\cal U}({\frak g})\,$; and 
since $H\left(\nu_V,p_V\right)$ has the definition (1.6.1), and is irreducible with respect to 
the actions of the elements of ${\cal U}({\frak g})$, 
one must have 
$\Psi\varphi_{\left| p_V\right|,0}\left(\nu_V,p_V\right)=\varphi_{\ell,q}\left(\nu_V,p_V\right)$
for some $\Psi\in {\cal U}({\frak g})$.  Finally, as somewhat of a converse 
to (1.5.17) and the results noted between (1.4.7) and (1.4.9), it follows 
by Lemma~5.2.1 of [32] that the existence of the expansion (1.7.13), 
when combined with the fact that $T_V\varphi_{\ell,q}\left(\nu_V,p_V\right)
\in A_{\Gamma}\left(\Upsilon_{\nu_V,p_V};\ell,q\right)$, 
is sufficient to imply that the growth condition (1.4.13) is satisfied when  
$f=T_V\varphi_{\ell,q}\left(\nu_V,p_V\right)$. 
This (given the arbitrary choice of cusp ${\frak c}$) shows that 
$T_V\varphi_{\ell,q}\left(\nu_V,p_V\right)\in 
A_{\Gamma}^{\rm pol}\left(\Upsilon_{\nu_V,p_V};\ell,q\right)$. 
Hence and by (1.7.12) and (1.7.8), one obtains (1.7.10). 
See [11], Theorem~3.3.1, for the corresponding proof 
in the $K$-trivial case.

As indicated by (1.7.10) and (1.7.7), the spectral 
parameters $\left(\nu_V , p_V\right)$ are shared by all elements of $V_K$ (i.e. 
$\Omega_{\pm}f=\Upsilon_{\nu_V,p_V}\left(\Omega_{\pm}\right) f$ for 
$f\in V_K$ and either choice of sign, where the eigenvalues 
$\Upsilon_{\nu_V,p_V}\left(\Omega_{\pm}\right)$ are as given by (1.3.3)). 
Therefore  (and since $V_K$ is dense in $V$) one calls  
$\left(\nu_V , p_V\right)$ the `spectral parameters of $V$'. 
By (1.7.8) and (1.7.10), the 
one-dimensional space ${\Bbb C}T_V\varphi_{\ell,q}\left(\nu_V,p_V\right)$ is 
the subspace $V_{K,\ell,q}\subset V_K$ spanned by all cusp-forms of $K$-type $(\ell,q)$ in $V_K$. 
With the operator $T_V$ being an isometry, one necessarily has: 
$$\left\| T_V \varphi_{\ell,q}\left(\nu_V,p_V\right)\right\|_{\Gamma\backslash G}
=\cases{\left\|\varphi_{\ell,q}\left(\nu_V,p_V\right)\right\|_{\rm ps} 
&if $\left(\nu_V,p_V\right)\in (i{\Bbb R})\times{\Bbb Z}\,$,\cr 
\left\|\varphi_{\ell,q}\left(\nu_V,p_V\right)\right\|_{\rm cs} 
&if $0<\nu_V^2<1$ and $p_V=0\,$,}\eqno(1.7.14)$$
where the norms $\| \cdot \|_{\rm ps}$ and $\| \cdot \|_{\rm cs}$ are as 
described between (1.6.4) and (1.6.6). Therefore $V$, ${\frak c}$ and  $g_{\frak c}$ 
determine both the set 
$\{ e^{i\phi} T_V \varphi_{\ell,q}\left(\nu_V,p_V\right) : 0\leq\phi <2\pi\}$ 
and the set  
$\{ ( e^{i\phi} c_V^{\frak c}(\omega))_{\omega\in{\frak O}-\{ 0\}} : 
0\leq\phi <2\pi\}$. 
It nevertheless later becomes convenient to work instead with modified Fourier coefficients: 
$$C_V^{\frak c}\left(\omega;\nu_V,p_V\right) 
=\left(\pi|\omega|\right)^{\nu_V}\left({\omega /|\omega|}\right)^{-p_V} 
c_V^{\frak c}(\omega).\eqno(1.7.15)$$
By the point noted below (1.4.13), the dependence of $c_V^{\frak c}(\omega)$ upon $g_{\frak c}$ 
is quite simple. Indeed, the Fourier coefficients $c_V^{\frak c}(\omega)$ are, as one would expect, 
essentially determined by the $\Gamma$-equivalence class of the cusp ${\frak c}\,$: 
for if ${\frak a}\sim^{\!\!\!\!\Gamma}{\frak b}$ 
and $\tau\in{\Bbb C}$, $\eta\in\{ u\in{\Bbb C} : u^8=1\}$ are such that  
$h[\eta] n[\tau]\in g_{\frak a}^{-1}\Gamma g_{\frak b}$ 
then 
$$c_V^{\frak b}(\omega)
=\eta^{2 p_V}\,{\rm e}\left({\rm Re}(\tau\omega)\right) 
c_V^{\frak a}\left(\eta^{-2}\omega\right)\quad
\hbox{for $\omega\in{\frak O}-\{ 0\}\,$.}\eqno(1.7.16)$$

\medskip 

\noindent{\bf Remark 1.7.1 (on spectral parameters).}\quad 
Except in those cases where the 
elements of $V$ are functions with character $\Upsilon_{0,0}\,$ 
(so that, by (1.7.10) and (1.3.3), one has $\nu_V=p_V=0$), 
the procedure for assigning spectral parameters to $V$ 
indicated in (1.7.6) will yield exactly two choices 
for the spectral parameters $(\nu_V,p_V)$ of $V$, with the 
two choices in question, 
$(\nu',p')$ and $(\nu'',p'')\,$ (say),  satisfying the relation  
$\nu'+\nu''=p'+p''=0$. Consequently one may substitute 
$(-\nu_V,-p_V)$ for $(\nu_V,p_V)$ in all of the points covered 
between (1.7.7) and (1.7.16) (always provided that the 
appropriate operator on $H(-\nu_V,-p_V)$ is substituted for $T_V$). 
\par 
We assume henceforth that 
each irreducible subspace 
$V\subset{}^0L^2(\Gamma\backslash G)$ is  
assigned a specific choice of spectral parameters $(\nu_V,p_V)$. 
However, since the essential points of what follows in this paper 
are independent of the choices made in the course of assigning 
those spectral parameters, we allow that those choices may be made 
arbitrarily. With regard to this note that, by (1.7.8), (1.6.5), 
(1.6.6), (1.7.13)-(1.7.15) and the functional equation 
$$\left(\pi |\omega|\right)^{-\nu} 
\left( i\omega /|\omega|\right)^p 
\Gamma(\ell +1+\nu) {\bf J}_{\omega}\varphi_{\ell,q}(\nu,p) 
=\left(\pi |\omega|\right)^{\nu} 
\left( i\omega /|\omega|\right)^{-p}  
\Gamma(\ell +1-\nu) {\bf J}_{\omega}\varphi_{\ell,q}(-\nu,-p)\eqno(1.7.17)$$
(which is Equation~(5.29) of~[5]), it follows that the operator 
$T_V : H(\nu_V,p_V)\rightarrow V_K$ and 
its counterpart with domain $H(-\nu_V,-p_V)$ determine a 
real constant $\phi\in[0,2\pi)$ such that, for all 
non-zero $\omega\in{\frak O}$ and all cusps ${\frak c}$ of $\Gamma$, 
one has $C_V^{\frak c}\left(\omega;\nu_V,p_V\right) 
=e^{i\phi}C_V^{\frak c}\left(\omega;-\nu_V,-p_V\right)$. 
Therefore it is in particular the case that each summand of 
the sum over $V$ occurring in the Spectral Sum Formula of   
Theorem~B (below) is unchanged  
if $(-\nu_V,-p_V)$ is substituted for $(\nu_V,p_V)\,$
(it being assumed here that the relevant function $h$ satisfies the 
condition (i) of Theorem~B); the same is true of the summands 
of that sum over $V$ which occurs in the definition of 
$E_{0}^{\frak a}\left( q_0 ,P,K;N,b\right)$ given in 
Theorem~1 (below). 
 
\bigskip

\centerline{\bf \S 1.8 Decomposing the subspace ${}^{\rm e}L^2(\Gamma\backslash G)\,$: the Eisenstein series and a Parseval identity.}

\medskip

The subspace ${}^{\rm e}L^2(\Gamma\backslash G)$ in (1.7.3) 
is generated by integrals of certain Eisenstein series. 
To obtain a set of these series, sufficient for the generation of ${}^{\rm e}L^2(\Gamma\backslash G)$, 
first choose (once and for all) a complete set of representatives ${\frak C}(\Gamma)$
of the $\Gamma$-equivalence classes of cusps, and, for each ${\frak c}\in{\frak C}(\Gamma)$,    
a scaling matrix $g_{\frak c}$ satisfying (1.1.16)-(1.1.21). 
Then, for ${\frak c}\in{\frak C}(\Gamma)$, $\ell,p,q\in{\Bbb Z}$ with $\ell\geq |p|,|q|$    
and $\nu\in{\Bbb C}$ with ${\rm Re}(\nu)>1$, the Eisenstein series
$E_{\ell,q}^{\frak c}(\nu,p) : G\rightarrow{\Bbb C}$ is given by: 
$$E_{\ell,q}^{\frak c}(\nu,p)(g)
={1\over\left[\Gamma_{\frak c} : \Gamma_{\frak c}'\right]} 
\sum_{\gamma\in\Gamma_{\frak c}'\backslash\Gamma}
\varphi_{\ell,q}(\nu,p)\left( g_{\frak c}^{-1}\gamma g\right)\quad
\hbox{for $g\in G\,$,}\eqno(1.8.1)$$
where $\varphi_{\ell,q}(\nu,p)\in C^{\infty}(N\backslash G)$ is as in (1.3.2). 
By (1.1.20)-(1.1.21), the sum in (1.8.1) is well-defined; 
the results on the $K$-trivial case $\ell=p=q=0$ in  
[11], Proposition~3.1.3, Proposition~3.2.1, Proposition~3.2.3 and Corollary~3.1.6,    
imply that, while this sum is divergent (for almost all $g\in G$) when $\nu\leq 1$, 
it does converge uniformly (and absolutely) 
for the pairs $(\nu ,g)\in{\Bbb C}\times N a[r] K$ with ${\rm Re}(\nu)\geq 1+\varepsilon$ and $r\geq\varepsilon$, 
where $\varepsilon$ is an arbitrary positive constant.
The definition (1.8.1) ensures that $E_{\ell,q}^{\frak c}(\nu,p)$ is 
a $\Gamma$-automorphic function: it moreover inherits from $\varphi_{\ell,q}(\nu,p)$  
the property of being a function of $K$-type $(\ell,q)$ with character $\Upsilon_{\nu,p}$. 

By (1.3.2) and (1.3.4), one has 
$$\varphi_{\ell,q}(\nu,p)\left( n h[u] g\right) 
=|u|^{2(1+\nu)}\left({u/|u|}\right)^{-2p}\varphi_{\ell,q}(\nu,p)(g)\quad 
\hbox{for $u\in{\Bbb C}-\{ 0\}$, $n\in N\;$.}\eqno(1.8.2)$$
Consequently, for cusps ${\frak c}$ with 
$g_{\frak c}^{-1}\Gamma_{\frak c} g_{\frak c}\cap h[i] N\neq\emptyset$,    
the sum in (1.8.1) will equal zero whenever $p$ is odd. 
Since  $\left[\Gamma_{\frak c} : \Gamma_{\frak c}'\right]=4$ if 
$g_{\frak c}^{-1}\Gamma_{\frak c} g_{\frak c}\cap h[i] N\neq\emptyset$,
while $\left[\Gamma_{\frak c} : \Gamma_{\frak c}'\right]=2$ otherwise, it 
therefore follows that 
$$E_{\ell,q}^{\frak c}(\nu,p) \neq0\quad\hbox{only if}\quad 
p\in\textstyle{1\over 2}\left[\Gamma_{\frak c} : \Gamma_{\frak c}'\right]{\Bbb Z}\;.\eqno(1.8.3)$$
For $p\in{1\over 2}\left[\Gamma_{\frak c} : \Gamma_{\frak c}'\right]{\Bbb Z}$ 
the definition (1.8.1) coincides with Definition~3.3.2 of [32].

It is almost immaterial exactly which representative ${\frak c}$ and which scaling matrix $g_{\frak c}$, 
are chosen (as above) for use in defining the Eisenstein series: for it follows by 
(1.8.1)-(1.8.3), Lemma~2.1 and Lemma~4.2 that a different choice of ${\frak c}$ or $g_{\frak c}$
(in respect of any one $\Gamma$-equivalence class of cusps) 
will merely replace $E_{\ell,q}^{\frak c}(\nu,p)$ by a function equal to 
$\epsilon^p E_{\ell,q}^{\frak c}(\nu,p)$, for some $\epsilon\in{\frak O}^*$. 

By (1.5.4), and since $\varphi_{\ell,q}(\nu,p)\in C^{\infty}(N\backslash G)$, the Eisenstein series of  
(1.8.1) is a Poincar\'{e} series, $(P^{\frak c}\varphi_{\ell,q}(\nu,p))(g)$, to which one  
might hope the case $\omega =0$ of  (1.5.5)-(1.5.10) would apply;   
it can be shown to follow from the definition (1.3.2) that this 
hope is justified when one has ${\rm Re}(\nu)>1$.
Hence and by (1.4.1), (1.5.18), (1.8.2), (1.8.3), Lemma~4.2 (below) and 
the linearity inherent in the definition~(1.5.2) one finds that,  
if ${\frak a},{\frak b}\in{\frak C}(\Gamma)$, and if  $\ell,p,q\in{\Bbb Z}$ and $\nu\in{\Bbb C}$
are such that $\ell\geq\max\{ |p|,|q|\}$,  $\,{\rm Re}(\nu)>1$ and $E_{\ell,q}^{\frak a}(\nu,p)\neq 0$, 
then, for $g\in G$, 
$$\eqalignno{
E_{\ell,q}^{\frak a}(\nu,p)\left( g_{\frak b} g\right) 
 &=\delta^{\Gamma}_{{\frak a},{\frak b}}\varphi_{\ell,q}(\nu,p)( g) 
+{1\over \left[\Gamma_{\frak a} : \Gamma_{\frak a}'\right]}\,D_{\frak a}^{\frak b}(0;\nu,p)
\,{\pi\Gamma(|p|+\nu)\over\Gamma(\ell+1+\nu)}\,{\Gamma(\ell+1-\nu)\over\Gamma(|p|+1-\nu)}\,
\varphi_{\ell,q}(-\nu,-p)(g)\ +\qquad \cr 
 &\quad\ +{1\over \left[\Gamma_{\frak a} : \Gamma_{\frak a}'\right]}
\sum_{0\neq\psi\in{\frak O}}
D_{\frak a}^{\frak b}\left(\psi;\nu,p\right) 
\left( {\bf J}_{\psi}\varphi_{\ell,q}(\nu,p)\right)\!(g)\;, &(1.8.4)}$$
where
$$\delta^{\Gamma}_{{\frak a},{\frak b}} 
=\cases{1 &if ${\frak a}\sim^{\!\!\!\!\Gamma}{\frak b}$, \cr 
0 &otherwise,}\eqno(1.8.5)$$ 
and 
$$D_{\frak a}^{\frak b}(\psi;\nu,p)
=\sum_{c\in\,{}^{\frak a}{\cal C}^{\frak b}}^{\hbox{\quad }}
S_{{\frak a} , {\frak b}}\!\left( 0 , \psi ; c\right) |c|^{-2(1+\nu)}\left( c/|c|\right)^{2p}\;.\eqno(1.8.6)$$ 
Note that Theorem~3.4.1 of~[11],  
contains the $K$-trivial case of this Fourier expansion.  
The sums $S_{{\frak a} , {\frak b}}\!\left( 0 , \psi ; c\right)$ generalise the 
Ramanujan sum evaluated in [15], Theorem~271; one has in particular 
a better estimate for these sums than that provided by (1.5.12); it can consequently be shown that 
the sum in (1.8.6) is absolutely convergent when either 
$\psi =0$ and ${\rm Re}(\nu)>1$, or $0\neq\psi\in{\frak O}$ and ${\rm Re}(\nu)>0$. 

It follows by (1.3.2) and the uniform convergence of the series in (1.8.1) that 
when ${\frak c}$, $\ell$, $p$, $q$ and $g$ are given, the function 
$\nu\mapsto E_{\ell,q}^{\frak c}(\nu,p)(g)$ is holomorphic for ${\rm Re}(\nu)>1$. 
Furthermore, it is known that this function of $\nu$ has a meromorphic 
continuation to all of ${\Bbb C}$ with (in the particular cases considered in this paper)
a simple pole at $\nu =1$ if and only if $\ell =p=q=0$, and no other poles in 
the closed half plane $\{\nu\in{\Bbb C} : {\rm Re}(\nu)\geq 0\}$. 
Applying this meromorphic continuation for each $g\in G$, one obtains, when 
${\rm Re}(\nu)\geq 0$ and $(\nu ,p)\not\in\{\,(0,0) , (1,0)\}$, 
a function $E_{\ell,q}^{\frak c}(\nu,p) : G\rightarrow{\Bbb C}$ which lies in the 
space $C^{\infty}(\Gamma\backslash G)$ and inherits (from the 
functions $E_{\ell,q}^{\frak c}\left(\nu',p\right)$ with ${\rm Re}\left(\nu'\right)>1$) 
the properties of being of $K$-type $(\ell,q)$ with character $\Upsilon_{\nu,p}$. 
Due to the nature of the first two terms on the right-hand side of the Fourier expansion (1.8.4), 
the function $E_{\ell,q}^{\frak c}(\nu,p)$ does not lie in the space $L^2(\Gamma\backslash G)$ 
(except, possibly, when $(\nu,p)=(0,0)$\/): this can be seen by evaluation of the integral 
$\int_{(z,r)\in {\cal E}_{\frak c}}\int_{k\in K^{+}}\left|\varphi_{\ell,q}\left( g_{\frak c} n[z]a[r]k\right)\right|^2 
r^{-3}{\rm d}_{+}z\,{\rm d}r\,{\rm d}k$, with ${\cal E}_{\frak c}\subset{\Bbb H}_3$ 
as in (1.1.23). Therefore it is only by averaging $E_{\ell,q}^{\frak c}(\nu,p)$ over a range of values 
of $\nu$ that one obtains an element of the space 
${}^{e}L^{2}(\Gamma\backslash G)$ (see Theorem A below). 

For trivial $K$-type (i.e. when $\ell=p=q=0$) the above facts concerning the meromorphic 
continuation of the Eisenstein series are established in 
[11], Theorem~6.1.2 and Theorem~6.1.11; 
this being achieved by means of an elegant general theory  
(valid when one substitutes for $\Gamma$ any discrete subgroup $\Gamma'<SL(2,{\Bbb C})$ 
for which the corresponding findamental domain ${\cal F}'\subset{\Bbb H}_3$ 
is non-compact, yet of finite volume): 
the corresponding facts in respect of Eisenstein series of arbitrary $K$-type 
are contained in Langlands' even more general theory [31].
Unlike the general situation described in Proposition~6.2.2 of~[11],  
there is here no  possibility of a (finite) number of generators of 
the space ${}^{\rm e}L^2(\Gamma\backslash G)$ arising from  
residues of the Eisenstein series $E_{0,0}^{\frak c}(\nu ,0)$ at a poles 
lying in the interval $(0,1]$, for the only pole of $E_{0,0}^{\frak c}(\nu ,0)$
with a positive real part is that at $\nu =1$, and the residue there is a function 
that is constant on $G$ (and therefore orthogonal to  ${}^{\rm e}L^2(\Gamma\backslash G)$\/).

Alternative proof of the above remarks on meromorphic continuation of Eisenstein series 
may be obtained by detailed consideration of 
the particular coefficients $D_{\frak a}^{\frak b}(\psi;\nu,p)$ in~(1.8.4) and~(1.8.6). 
This requires an evaluation of the sum $S_{{\frak a},{\frak b}}(0,\psi;c)$ 
(analogous to the evaluation of the classical Ramanujan 
obtained in Theorem~271 of [15]), 
which turns out to be a not overly complicated affair in the case ${\frak b}=\infty$. 
The result of this calculation (for ${\frak b}=\infty$) leads, via 
(1.8.6), to an expression for $D_{\frak a}^{\infty}(\psi;\nu ,p)$  in terms of 
Hecke zeta-functions 
$$\zeta\left( s , \lambda^{p/2}\chi\right)
={1\over 4}\sum_{0\neq\alpha\in{\frak O}}
{\lambda^{p/2}(\alpha)\chi(\alpha)\over |\alpha|^{2s}}\qquad\quad\   
\hbox{(${\rm Re}(s)>1$),}\eqno(1.8.7)$$
where $\lambda^{m}(\alpha)=(\alpha /|\alpha|)^{4m}$ and either 
$s=\nu$ and $\chi : {\frak O}\rightarrow\{1\}$, or
$s=1+\nu$ and there is a primitive character 
$\tilde\chi : ({\frak O}/d{\frak O})^{*}\rightarrow{\Bbb C}^{*}$ 
such that $\chi(\alpha)=\tilde\chi({\cal A})$ whenever $\alpha\in{\cal A}\in ({\frak O}/d{\frak O})^{*}$
(while $\chi(\alpha)=0$ if $|(\alpha ,d)|>1$). 
It is therefore a corollary of Hecke's work in Section~6 of [17] that 
$D_{\frak a}^{\infty}(\psi;\nu,p)$ can be meromorphically continued into 
all of ${\Bbb C}$; consequently one obtains, 
via (1.8.4) and Lemma~5.1 of~[5], the meromorphic 
continuation of $E_{\ell,q}^{\frak a}(\nu,g)(g)$. 
See Lemma~5.2 of [5] for an explicit determination of 
$D_{\infty}^{\infty}(\psi;\nu,p)$ in the case $q_0=1$ (i.e. for $\Gamma =SL(2,{\frak O})$\/). 
By Lemma~5.1 of [5] (again), and by (1.4.1)-(1.4.2), (1.5.14), (1.5.15) and, 
for ${\rm Re}(\nu)>1$, the equation~(1.8.4), 
the meromorphic continuation of $E_{\ell,q}^{\frak a}(\nu,p)(g)$ implies that of 
$D_{\frak a}^{\frak b}(\psi;\nu,p)$ for all ${\frak b}\in{\frak C}(\Gamma)$ 
(i.e. not only in the special case ${\frak b}=\infty$\/). 

As the remarks of the last three paragraphs might suggest, the Eisenstein series enable one to describe a decomposition 
of the subspace ${}^{e}L^{2}(\Gamma\backslash G)$ in (1.7.3). By combining this decomposition 
with (1.7.4)-(1.7.8) (the decomposition of ${}^{0}L^{2}(\Gamma\backslash G)$\/)
one can obtain a useful decomposition of $L^{2}(\Gamma\backslash G)$. 
One may work instead with restriction of the decomposition (1.7.3) to the 
subspace $L^2(\Gamma\backslash G;\ell,q)$ 
spanned by the elements in $L^2(\Gamma\backslash G)$ of $K$-type $(\ell,q)$. 
Then a key result is the following.

\bigskip

\proclaim Theorem A (a Parseval identity). 
Let $\ell,q\in{\Bbb Z}$ satisfy $\ell\geq |q|$, and suppose  that 
$f_1,f_2\in L^2(\Gamma\backslash G;\ell,q)$ are represented by bounded elements of
$C^{\infty}(\Gamma\backslash G)$. 
Then, when ${\frak c}\in{\frak C}(\Gamma)$, $j\in\{ 1,2\}$ 
and $p\in{\Bbb Z}$ with $|p|\leq\ell$, the inner product 
$\langle f_j , E_{\ell,q}^{\frak c}(it,p)\rangle_{\Gamma\backslash G}
=F_{j,p}^{\frak c}(t)$ (say) is  
defined  for all real $t$; 
the functions so defined, $F_{j,p}^{\frak c} : {\Bbb R}\rightarrow{\Bbb C}$  
(${\frak c}\in{\frak C}(\Gamma)$, $j=1,2$ and $p=-\ell,\ldots,\ell$\/), 
are each square-integrable with respect to the Lebesgue measure on~${\Bbb R}$, and one has: 
$$\eqalign{
\left\langle f_1 , f_2\right\rangle_{\Gamma\backslash G}
 &=\,{1\over{\rm vol}(\Gamma\backslash G)}\, 
\left\langle f_1 , 1\right\rangle_{\Gamma\backslash G}
\left\langle 1 , f_2\right\rangle_{\Gamma\backslash G}\ +\cr
 &\quad\ +\sum_{\scriptstyle V\atop\scriptstyle -\ell\leq p_V\leq\ell} 
{1\over\left\| T_V\varphi_{\ell,q}\left(\nu_V,p_V\right)\right\|_{\Gamma\backslash G}^2}\,
\left\langle f_1 , T_V\varphi_{\ell,q}\left(\nu_V,p_V\right)\right\rangle_{\Gamma\backslash G}
\left\langle T_V\varphi_{\ell,q}\left(\nu_V,p_V\right) , f_2\right\rangle_{\Gamma\backslash G}\ +\cr
 &\quad\ +\sum_{{\frak c}\in{\frak C}(\Gamma)} 
{\left[\Gamma_{\frak c} : \Gamma_{\frak c}'\right]\over 4\pi i}
\sum_{\scriptstyle p\in{1\over 2}\left[\Gamma_{\frak c} : \Gamma_{\frak c}'\right]{\Bbb Z}\atop\scriptstyle |p|\leq\ell}
\ \int\limits_{(0)} {1\over
\left\|\varphi_{\ell,q}\left(\nu,p\right)\right\|_{\rm ps}^2}\,
\left\langle f_1 , E_{\ell,q}^{\frak c}\left(\nu,p\right)\right\rangle_{\Gamma\backslash G}
\left\langle E_{\ell,q}^{\frak c}\left(\nu,p\right) , f_2\right\rangle_{\Gamma\backslash G}
{\rm d}\nu\;,}\eqno(1.8.8)$$
where the `$1$' 
in both $\left\langle f_1 , 1\right\rangle_{\Gamma\backslash G}$ and 
$\left\langle 1 , f_2\right\rangle_{\Gamma\backslash G}$  
denotes the constant function $\varphi_{0,0}(-1,0)$ defined by (1.3.2), 
while $V$, in the second summation on the right-hand side of the equation, 
runs over the pairwise-orthogonal cuspidal subspaces of $L^2(\Gamma\backslash G)$ 
that occur in the direct sum in (1.7.4), and the notation `$(0)$' below the integral sign signifies that 
the integration is along the line $\{\nu\in{\Bbb C} : {\rm Re}(\nu)=0\}$, 
oriented as a contour from $-i\infty$ to $i\infty\,$. 
The sums and integrals in equation (1.8.8) are absolutely convergent. 

\medskip 

\noindent{\bf Proof.}\quad 
These results are a special case of Theorem~8.1 of [32], 
and are (conversely) a slight generalisation of Theorem 8.1 of [5]. 
They are also a special case of the very general Parseval identity  
proved in [31] (see also [16]).  
As is noted in [32],  the restriction of (1.8.8) to pairs of 
$K$-trivial functions $f_1,f_2$ is 
(effectively) a result contained in Theorem~6.3.4 of~[11]\ $\blacksquare$

\bigskip 

The definition (1.8.6) of  the Fourier coefficients $D_{\frak a}^{\frak b}(\omega ;\nu,p)$ 
is not only valid for cusps in the particular set of representatives ${\frak C}(\Gamma)\,$: 
it is in fact the appropriate definition for a arbitrary pair of cusps ${\frak a},{\frak b}$ of $\Gamma$. For if 
${\frak a}$ and ${\frak b}$ are not $\Gamma$-equivalent then one may assume that 
${\frak C}(\Gamma)\supseteq\{ {\frak a} , {\frak b}\}\,$; while if instead 
${\frak a}\sim^{\!\!\!\!\Gamma}{\frak b}$, then one may assume that 
${\frak C}(\Gamma)\ni{\frak b}$, and so make use of the fact that
$E_{\ell,q}^{\frak a}(\nu ,p)\left( g_{\frak b} g\right)={\epsilon}^p E_{\ell,q}^{{\frak b}}(\nu,p)\left( g_{\frak b} g\right)$  
for some $\epsilon\in{\frak O}^{*}$. 
Indeed,  apart from the coefficient $\delta^{\Gamma}_{{\frak a},{\frak b}}$ in (1.8.4) possibly 
being replaced by $\epsilon^p \delta^{\Gamma}_{{\frak a},{\frak b}}$, for some $\epsilon\in{\frak O}^{*}$, 
the whole Fourier expansion (1.8.4)-(1.8.6) is valid for arbitrary cusps ${\frak a},{\frak b}\,$: 
and this does not require that  ${\frak a}={\frak b}$ imply   
$g_{\frak a}= g_{\frak b}$. 
A useful normalisation of $D_{\frak a}^{\frak b}(\psi;\nu,p)$ (analogous to (1.7.15)) is
given by: 
$$B_{\frak a}^{\frak b}(\omega;\nu,p)
=\left(\pi|\omega|\right)^{\nu}\left(\omega/|\omega|\right)^{-p}
D_{\frak a}^{\frak b}(\omega;\nu,p)\;.\eqno(1.8.9)$$

\medskip

\centerline{\bf \S 1.9 Results and applications.}

\medskip

The principal new results of this paper depend on being able to
deduce estimates for a certain mean-value of 
Fourier coefficients of $\Gamma$-automorphic cusp forms and Eisenstein series 
from suitable bounds for sums of the generalised Kloosterman sums defined in (1.5.10).
This is made possible by the following result. 

\bigskip

\proclaim Theorem B (spectral sum formula). 
Let the real numbers $\,\sigma\in(1/2,1)$,  $\,\varrho,\vartheta \in(3,\infty)$, 
and the function 
$\,h : \{\nu\in{\Bbb C} : |{\rm Re}(\nu)|\leq\sigma\}\times{\Bbb Z}\rightarrow{\Bbb C}$,  
satisfy the three conditions\smallskip 
(i)\quad\ $h(\nu,p)=h(-\nu ,-p)\;$;\hfill\smallskip 
(ii)\quad for $p\in{\Bbb Z}$, the function 
$\nu\mapsto h(\nu,p)$ can be holomorphically 
continued into a neighbourhood of the\hfill\break   
$\hbox{\qquad\qquad\ }$strip $\{\nu\in{\Bbb C} : |{\rm Re}(\nu)|\leq\sigma\}$;\smallskip 
(iii)\ $\ \,h(\nu,p)\ll_{h,\varrho,\vartheta }\,(1+|{\rm Im}(\nu)|)^{-\varrho} (1+|p|)^{-\vartheta}$. \medskip 
\noindent Suppose moreover that 
$0\neq q_0\in{\frak O}={\Bbb Z}[i]$, and that 
$\,\Gamma=\Gamma_0\!\left( q_0\right)\leq SL(2,{\frak O})$. 
Then, for all $\omega_1,\omega_2\in{\frak O}-\{ 0\}$, 
all pairs of cusps ${\frak a},{\frak b}$ of 
$\,\Gamma$, and all choices of the associated 
scaling matrices $g_{\frak a}, g_{\frak b}$ that 
satisfy the conditions (1.1.16) and (1.1.20)-(1.1.21), 
one has   
$$\eqalign{ 
 &\sum_V \,\overline{C_V^{\frak a}\left(\omega_1;\nu_V,p_V\right)}\,
C_V^{\frak b}\left(\omega_2;\nu_V,p_V\right) h\left(\nu_V , p_V\right)\ +\cr 
 &\qquad +\sum_{{\frak c}\in{\frak C}(\Gamma)}
{1\over 4\pi i\left[\Gamma_{\frak c} : \Gamma_{\frak c}'\right]} 
\sum_{p\in{1\over 2}\left[\Gamma_{\frak c} : \Gamma_{\frak c}'\right]{\Bbb Z}}
\ \int\limits_{(0)} \overline{B_{\frak c}^{\frak a}\left(\omega_1;\nu,p\right)}\,
B_{\frak c}^{\frak b}\left(\omega_2;\nu,p\right) h(\nu,p)\,{\rm d}\nu =\cr
 &\qquad\qquad\qquad\qquad\qquad\qquad\qquad\qquad\qquad\quad   
=\;{1\over 4\pi^3 i}\,\delta^{{\frak a},{\frak b}}_{\omega_1,\omega_2} 
\sum_{p\in{\Bbb Z}}\ \int\limits_{(0)} h(\nu,p) \left( p^2 -\nu^2\right)\,{\rm d}\nu\ +\cr 
 &\qquad\qquad\qquad\qquad\qquad\qquad\qquad\qquad\qquad\qquad\qquad    
+\sum_{c\in {}^{\frak a}{\cal C}^{\frak b}} 
\,{S_{{\frak a},{\frak b}}\left(\omega_1 , \omega_2 ; c\right)\over |c|^2}\,
({\bf B}h)\!\!\left( {2\pi\sqrt{\omega_1\omega_2}\over c}\right) , 
}\eqno(1.9.1)$$
where  
$$\delta^{{\frak a},{\frak b}}_{\omega_1,\omega_2}
=\sum_{\scriptstyle \gamma\in\Gamma_{\frak a}'\backslash\Gamma\;:\;\gamma{\frak b}={\frak a}\atop\scriptstyle 
g_{\frak a}^{-1} \gamma g_{\frak b}
=\left({\scriptstyle\!\!\!\!\!u(\gamma)\quad\beta(\gamma)\atop\scriptstyle\  
0\quad\ 1/u(\gamma)}\right)} 
{\rm e}\left( {\rm Re}\left( \beta(\gamma)u(\gamma)\omega_1\right)\right) 
\delta_{u(\gamma)\omega_1 , \omega_2/u(\gamma)}\;,\eqno(1.9.2)$$
other notation is as developed in 
(1.1.17)-(1.1.19), (1.5.6), (1.5.8)-(1.5.10), (1.7.4), (1.7.13)-(1.7.15), (1.8.1), (1.8.6) and (1.8.9), 
and   
$$({\bf B}h)(z)
={1\over 4\pi i} \sum_{p\in{\Bbb Z}}\ \int\limits_{(0)} {\cal K}_{\nu,p}(z) h(\nu,p)\left( p^2 -\nu^2\right) 
\,{\rm d}\nu\;,\eqno(1.9.3)$$
with 
$${\cal K}_{\nu,p}(z)={1\over\sin(\pi\nu)}\,\left( {\cal J}_{-\nu,-p}(z)-{\cal J}_{\nu,p}(z)\right)\;,
\eqno(1.9.4)$$
$${\cal J}_{\nu,p}(z)
=|z/2|^{2\nu}\,(z/|z|)^{-2p}\,J_{\nu -p}^{*}(z) J_{\nu +p}^{*}\left( \overline{z}\right)\;,\eqno(1.9.5)$$
and
$$J_{\xi}^{*}(z)=\sum_{m=0}^{\infty} {(-1)^m\,(z/2)^{2m}\over\Gamma(m+1)\Gamma(\xi +m+1)}\;.\eqno(1.9.6)$$
In (1.9.1) the set of representatives ${\frak C}(\Gamma)$ of the $\Gamma$-equivalence classes 
of cusps may be chosen independently of the given pair of cusps ${\frak a},{\frak b}$; and 
nothing more than (1.1.16) and (1.1.20)-(1.1.21) need be assumed in respect of the choice of 
scaling matrices $g_{\frak c}$ for ${\frak c}\in{\frak C}(\Gamma)$ 
(even in the event that  ${\frak C}(\Gamma)\cap\{ {\frak a},{\frak b}\}\neq\emptyset$\/): similarly, 
$g_{\frak a}$ is allowed to differ from $g_{\frak b}$, even when ${\frak a}={\frak b}$. 
All sums and integrals occurring in the equations (1.9.1) and (1.9.3) are 
absolutely convergent; the sum occurring in Equation~(1.9.2) has at most finitely many terms. 

\bigskip 

\noindent{\bf Remark~1.9.1 (on the proof of Theorem~B).}\quad  
Theorem~B is an extension of 
Bruggeman and Motohashi's Spectral-Kloosterman sum formula, Theorem~10.1 of [5],    
which pertains to the case in which one has $\Gamma =SL(2,{\Bbb Z}[i])\,$ (so that  
there is only one $\Gamma$-equivalence class of cusps). It builds also upon the work of 
Lokvenec-Guleska who, in Theorem~11.3.3 of [32], 
succeeded in generalising Bruggeman and Motohashi's method so 
as to obtain a Spectral 
Kloosterman sum formula for Hecke congruence subgroups over 
an arbitrary imaginary quadratic field. Theorem~11.3.3 of~[32] contains   
the case ${\frak a}={\frak b}=\infty$ of Theorem~B. A proof of 
Theorem~B is described in an appendix to this paper;  
in this proof    
the relevant steps of [4] and [32] are adapted so as to deal with any  
choice of the cusps ${\frak a}, {\frak b}$. 

\bigskip

\noindent{\bf Remark~1.9.2 (on a result of Kim and Shahidi).}\quad 
Given the points noted in Subsection~1.7$\,$  
(see, in particular, the case $p_V=\ell =q=0$ of (1.7.10), and 
what is discussed between (1.7.5) and~(1.7.8)), 
it follows 
from the result (1.4.15) of Kim and Shahidi 
that in the first summation in Equation~(1.9.1)  
the spectral parameters 
$\nu_V,p_V$  of the relevant subspaces $V\subset{}^0L^2(\Gamma\backslash G)$ must, 
in each instance, satisfy 
$${\rm either}\quad\left( \nu_V , p_V\right)\in (i{\Bbb R})\times{\Bbb Z}\;, 
\quad{\rm or\ else}\quad 
p_V=0\ \,{\rm and}\ \,\nu_V\in[-2/9,2/9]\;.\eqno(1.9.7)$$

\medskip

\noindent{\bf Remark~1.9.3 (on a Bessel function).}\quad 
By assigning a fixed value to either one of the variables in (1.9.6), one obtains 
a single variable complex function 
(i.e. either $z\mapsto J_{\xi}^{*}(z)$ or $\xi\mapsto J_{\xi}^{*}(z)$\/) that is holomorphic 
on ${\Bbb C}$. When $\xi$ is not an integer, two linearly independent solutions of Bessel's differential equation, 
$x^2 y''+xy'+(x^2-\xi^2)y=0$ ($x>0$), are $y_1=J_{\xi}(x)$ and $y_2=J_{-\xi}(x)$, 
where $$J_{\nu}(z)=(z/2)^{\nu} J_{\nu}^{*}(z)\eqno(1.9.8)$$
(this function $J_{\nu}(z)$ being Bessel's function of order $\nu$).  
When $\xi$ is an integer, the equations $y_1=J_{\xi}(x)$ and  $y_2=J_{-\xi}(x)$ 
(with $J_{\nu}(z)$ as in (1.9.8)) do define solutions of Bessel's differential equation, 
but these solutions are linearly dependent, for it follows from (1.9.8) and (1.9.6) that 
$$J_{-n}(z)=(-1)^n J_n(z)=J_n(-z)\qquad\hbox{($n\in{\Bbb Z}$, $z\in{\Bbb C}$).}\eqno(1.9.9)$$

\medskip 

\noindent{\bf Remark~1.9.4 (on a partial inversion of the operator ${\bf B}$).}\quad 
Let the function $f : {\Bbb C}^{*}\rightarrow{\Bbb C}$ be  
compactly supported and even.  Suppose moreover that 
$f$ is `smooth', in the sense that every partial derivative (of whatever order)  
of the function $(x,y)\mapsto f(x+iy)$ is defined and continuous 
on the set ${\Bbb R}^2-\{ (0,0)\}$. Then, as is shown in Theorem~11.1 of [5], 
one has: 
$$\pi {\bf BK}f = f\eqno(1.9.10)$$
where the operator ${\bf B}$ is that given by (1.9.3)-(1.9.6), while
$$({\bf K}f)(\nu,p)=\int_{{\Bbb C}^{*}} {\cal K}_{\nu,p}(z) f(z) |z|^{-2}{\rm d}_{+}z\quad
\hbox{for $p\in{\Bbb Z}$ and $\nu\in{\Bbb C}$ with $|{\rm Re}(\nu)|<1\,$.}$$
By applying Theorem B with $h={\bf K}f$ one obtains (using (1.9.10)) the corollary 
that the sum 
$$L_{{\frak a},{\frak b}}\!\left(\omega_1,\omega_2;f\right) 
=\sum_{c\in {}^{\frak a}{\cal C}^{\frak b}} 
\,{S_{{\frak a},{\frak b}}\left(\omega_1 , \omega_2 ; c\right)\over |c|^2}\,
f\!\left( {2\pi\sqrt{\omega_1\omega_2}\over c}\right)$$
may be expressed in terms of sums involving Fourier coefficients of 
$\Gamma$-automorphic cusp forms and Eisenstein series. By 
Lemma~11.1 of [5], this inversion of the summation formula (1.9.1)
contains no `diagonal term' (i.e. no counterpart of the term in (1.9.1) with 
coefficient $\delta_{\omega_1,\omega_2}^{{\frak a},{\frak b}}$\/);  
it is in fact only a one-sided (non-surjective) inversion, since, 
as is noted in Section~11 of [5], there are test functions $h$ 
satisfying the conditions (i)-(iii) of Theorem B that do produce a non-zero 
diagonal term on the right-hand side of (1.9.1). 
This inversion of Theorem B is not needed in this paper, but is 
important in [46] and [47].  

\bigskip

\noindent{\bf Remark~1.9.5.}\quad 
Bruggeman and Motohashi showed in Theorem~12.1 of [5] that 
when $0\neq z\in{\Bbb C}$, when $e^{i\theta}=z/|z|$ (so that $\theta\in{\Bbb R}$), 
and when one defines ${\cal K}_{\nu ,p}(z)$ by (1.9.4)-(1.9.6), it then follows that 
$${\cal K}_{\nu ,p}(u)={(-1)^p\over \pi/2}\int\limits_0^{\infty}\!y^{2\nu}
\!\left({ye^{i\theta}+(ye^{i\theta})^{-1}\over\left|
ye^{i\theta}+(ye^{i\theta})^{-1}\right|}\right)^{\!\!2p}\!\!J_{2p}\!\left( |u|\left|
ye^{i\theta}+(ye^{i\theta})^{-1}\right|\right)\!{{\rm d}y\over y}\quad\  
\hbox{($p\in{\Bbb Z}$, $|{\rm Re}(\nu)|<{1\over 4}$).}\eqno(1.9.11)$$ 
In proving Theorem 1 one needs to consider, for a suitable test function $h$, 
the transform ${\bf B}h$ that is defined in (1.9.3). Useful approximations   
to the relevant transformed function $({\bf B}h)(z)$ may be deduced 
with the aid of both the identity (1.9.11) and an addition law for Bessel functions 
(see Lemma 4.5 and the proof of Lemma 4.6, below). 

\bigskip

The principal new result in this paper is Theorem~1, stated next. We prove 
this theorem in Section~5 of this paper, with the help of certain 
bounds for sums of Kloosterman sums. These bounds 
are supplied by Proposition~2 (which we state after 
several remarks following Theorem~1). 

\bigskip

\proclaim Theorem 1. 
Let $\varepsilon >0$, $0\neq q_0\in{\frak O}={\Bbb Z}[i]$, 
$\Gamma =\Gamma_0\!\left( q_0\right)\leq SL(2,{\frak O})$ and 
$K,P,N\geq 1$.  
Suppose further that $b : {\frak O}-\{ 0\}\rightarrow{\Bbb C}$, and that 
$u,w\in{\frak O}$ satisfy $w\neq 0$ and $(u,w)\sim 1$. 
Then, when
${\frak a}$ is a cusp of $\Gamma$ 
with ${\frak a}\sim^{\!\!\!\!\Gamma}u/w$, and when 
$E_{0}^{\frak a}\left( q_0 ,P,K;N,b\right)$, $E_{1}^{\frak a}\left( q_0 ,P,K;N,b\right)$
are the quadratic moments given by  
$$E_{0}^{\frak a}\left( q_0 ,P,K;N,b\right) =\sum_{\scriptstyle V\atop\scriptstyle \left| p_V\right|\leq P ,\ \left|\nu_V\right|\leq K}
\Bigl|\sum_{\scriptstyle\omega\in{\frak O}\atop\scriptstyle N/2<|\omega|^2\leq N}
b(\omega) C_V^{\frak a}\left(\omega;\nu_V,p_V\right)\Bigr|^2\;,\eqno(1.9.12)$$
$$E_{1}^{\frak a}\left( q_0 ,P,K;N,b\right)=\sum_{{\frak c}\in{\frak C}(\Gamma)}
{1\over 4\pi\left[\Gamma_{\frak c} : \Gamma_{\frak c}'\right]}
\sum_{\scriptstyle p\in{1\over 2}\left[\Gamma_{\frak c} : \Gamma_{\frak c}'\right]{\Bbb Z}\atop\scriptstyle 
|p|\leq P}
\int_{-K}^{K}\Bigl|\sum_{\scriptstyle\omega\in{\frak O}\atop\scriptstyle N/2<|\omega|^2\leq N}
b(\omega) B_{\frak c}^{\frak a}(\omega;it,p)\Bigr|^2 {\rm d}t\eqno(1.9.13)$$
(where the terminology used has the same meaning as in Theorem B), one has the upper bounds:    
$$E_{j}^{\frak a}\left( q_0 ,P,K;N,b\right)
\ll \left( P^2+K^2\right)\left(PK+
O_{\varepsilon}\left( {N^{1+\varepsilon}\over (PK)^{1/2}}\,|\mu({\frak a})|^2\right)\right)
\left\|{\bf b}_N\right\|_2^2
\qquad\quad\hbox{($j=0,1$)}\,,\eqno(1.9.14)$$
where $\mu({\frak a})\in\{ 1/\alpha : 0\neq \alpha\in{\frak O}\}$,  
$${1\over\mu({\frak a})}\sim{\left( w , q_0\right) q_0\over \left( w^2 , q_0\right)}
\sim{q_0\over\bigl( \left( w , q_0\right) , q_0/\left( w , q_0\right)\bigr)}\eqno(1.9.15)$$
and 
$$\left\|{\bf b}_N\right\|_2
=\Biggl( \sum_{\scriptstyle\omega\in{\frak O}\atop\scriptstyle 
N/2<|\omega|^2\leq N}\!\!\!\!\!\!|b(\omega)|^2\Biggr)^{1/2}\,.\eqno(1.9.16)$$

\medskip

\noindent{\bf Remark~1.9.6.}\quad   
Since $\left[\Gamma_{\frak c} : \Gamma_{\frak c}'\right]\in\{ 2 , 4\}$ for 
all cusps ${\frak c}$, the factor 
$(4\pi\left[\Gamma_{\frak c} : \Gamma_{\frak c}'\right])^{-1}$ in (1.9.13) may 
be omitted.

\bigskip

\noindent{\bf Remark~1.9.7.}\quad 
One may check that  (1.9.15) makes the ideal $(1/\mu({\frak a})){\frak O}$   
a function of the $\Gamma$-equivalence class  of the cusp ${\frak a}$. The same is therefore true 
of the reciprocal of the norm of this ideal, which is the factor $|\mu({\frak a})|^2$  appearing in (1.9.14). 
Since $\infty\sim^{\!\!\!\!\Gamma}1/q_0$ (for $\Gamma =\Gamma_0(q_0)$), 
one has in particular $1/\mu(\infty)\sim 1/\mu(1/q_0)\sim q_0$. 

\bigskip 

\noindent{\bf Remark~1.9.8 (some comparisons and conjectures).}\quad
In their proof of Corollary~10.1 of~[5], Bruggeman and Motohashi show, 
by an application of their `Spectral-Kloosterman' sum formula  
(Theorem~10.1 of [5]), that if $\Gamma =SL(2,{\Bbb Z}[i])\,$ (so that $q_0\sim 1$) 
then, for $0\neq\omega\in{\frak O}$, $P\geq 1$, $K\geq 1$ and $\varepsilon >0$, 
one has 
$$\eqalign{ 
\sum_V \left| C_V^{\infty}\!\left(\omega :\nu_V , p_V\right)\right|^2 
\exp\!\left( (\nu/K)^2 -(p/P)^2\right) 
 &={1\over 4\pi^2}\,\left( P^2+K^2\right) PK\left( 1 
+O\left( P^2 e^{-\pi^2 P^2}\right)\right)\ + \cr 
 &\quad\ +O_{\varepsilon}\left( |\omega|^{1+\varepsilon}
\left( P^2+K^2\right)^{1+\varepsilon}\right) .}\eqno(1.9.17)$$
We may compare this with the result (5.55) obtained in 
the course of our proof of Theorem~1. Our result (5.55) is certainly 
weaker than (1.9.17) in cases where $|\omega|^2 >(PK)^{1+\varepsilon}$. 
However it follows from (5.55) that when   
$PK >|\omega|^2$ one may substitute 
$O_{\varepsilon}( |\omega|^2 (PK)^{-1/2} (P^2+K^2)^{1+\varepsilon} )$ 
for the final $O_{\varepsilon}$-term in Equation~(1.9.17). 
The proof of this does require certain estimates for 
the relevant instances of the 
sum occurring on the final line of Equation~(5.3); 
estimates sufficient for this purpose 
were already obtained in the proof of Corollary~10.1 of [5], 
by means of the lower bounds 
$\,|\zeta(1+\nu,\lambda^{p/2})|\gg 1/\log(|\nu|+|p|+2)\,$ 
($p\in 2{\Bbb Z}$, $\nu\in i{\Bbb R}$). 
\par 
Based on considerations 
relating to the Eisenstein series $E^{\frak a}_{0,0}(\nu,0)$ and 
Rankin-Selberg function,  
$$L_V^{\frak a}(s) 
=\sum_{0\neq\omega\in{\frak O}}\left| c_V^{\frak a}(\omega)\right|^2 |\omega|^{-2s}\;,$$ 
we conjecture that in general (i.e. not just for $q_0\sim 1$) it is 
the case that if, as in Lemma~2.2 (below), one has 
${\frak a}\sim^{\!\!\!\!\Gamma}u/w$ for some $u,w\in{\frak O}$ satisfying 
$(u,w)\sim 1$, $u\neq 0$ and $w\mid q_0$, then 
$$\sum_{\scriptstyle\omega\in{\frak O}\atop\scriptstyle 
0<|\omega|^2\leq N}\left| C_V^{\frak a}\!\left(\omega :\nu_V , p_V\right)\right|^2 
\sim\rho(\Gamma , {\frak a}) N\qquad\quad{\rm as}\quad N\rightarrow \infty\;,\eqno(1.9.18)$$ 
where 
$$\rho(\Gamma , {\frak a}) 
={2\over {\rm vol}(\Gamma\backslash G)}
\prod_{\scriptstyle\varpi\in ({\frak O}-\{ 0\})/{\frak O}^*\atop 
{\scriptstyle\varpi{\frak O}\ {\rm is\ prime}\atop\scriptstyle 
\varpi\mid (q/w , w/(q/w,w))}}\!\left( 1-|\varpi|^{-2}\right)$$ 
(so that one has, in particular, 
$\rho(\Gamma,\infty)=2/{\rm vol}(\Gamma\backslash G)$). 
We believe that this conjecture might be 
shown to be correct by methods 
analogous to those described in Section~8.2 of [22]. 
\par  
Returning from the conjectural to the proven we note that, 
through an application of the  
Spectral-Kloosterman summation formula obtained in Theorem~11.3.3 of [32],  
Lokvenec-Guleska obtains, in Section~11.5 of~[32], 
new asymptotic estimates for sums involving 
modified Fourier coefficients of cusp forms. In this work of 
Lokvenec-Guleska `${\frak O}$' 
denotes the ring of integers of an arbitrary quadratic number field 
$F={\Bbb Q}(\sqrt{-d})\,$ (with $d\in{\Bbb N}$), and `$\Gamma$' denotes 
a Hecke congruence subgroup of $SL(2,{\frak O})$ associated with 
some (arbitrary) non-zero ideal $I\subseteq{\frak O}$; her definition of 
automorphicity is broader than ours, in that it depends on a character 
$\chi : \Gamma\rightarrow{\Bbb C}$ derived from an arbitrary character for 
$({\frak O}/I)^{*}$;  the relevant 
irreducible cuspidal subspaces (corresponding to the spaces 
denoted by `$V$' in our paper) contain automorphic 
functions $f : G\rightarrow{\Bbb C}$ that 
are even if $\chi(-1)=1$, but odd if $\chi(-1)=-1$. 
We state here only one of 
Lokvenec-Guleska's results, specialised to the particular case 
${\frak O}={\Bbb Z}(i)$, $I=q_0{\frak O}$, $\chi : \Gamma\rightarrow\{ 1\}$,  
in which her `$L^2(\Gamma\backslash G,\chi)$' coincides with  
the space $L^2(\Gamma\backslash G)\,$ that we have defined. 
In respect of this case, it follows by the combination of Part~(i) of 
Theorem~11.5.2 of~[32] and the symmetry used in the deduction of 
Corollary~11.5.4 of [32] that, for $0\neq\omega\in{\frak O}$ and 
$p\in{\Bbb Z}$, one has 
$$\sum_{\scriptstyle V\ :\ p_V=p\atop\scriptstyle 1-\nu_V^2-p_V^2\leq X} 
\left| C_V^{\infty}\left(\omega ;\nu_V,p_V\right)\right|^2 
\sim {1\over 3\pi^3}\,X^{3/2}\qquad\quad{\rm as}\quad X\rightarrow\infty\;.\eqno(1.9.19)$$ 
The second condition of summation in (1.9.19) is motivated by the fact 
that $1-\nu_V^2-p_V^2$ is (by definition) an eigenvalue of  
the operator $-4(\Omega_{+}+\Omega_{-})$. 
\par 
Theorem~11.5.2 of~[32] is 
essentially a corollary (deduced via a Tauberian argument) 
of the asymptotic estimates for smoothly weighted sums 
which Lokvenec-Guleska obtains in Proposition~11.5.1 of~[32]. 
Results very similar to the special case 
$p=0$, $I={\frak O}\,$ ($\,\chi : \Gamma\rightarrow\{ 1\}$) of 
Proposition~11.5.1 of~[32] had 
previously been obtained in the paper [39] of Rhagavan and Sengupta;  
the result (13) of [39] is a somewhat unwieldy special case of the 
spectral-Kloosterman sum formula obtained in Theorem~11.3.3 of [32].

\par 
The combination of the result (1.9.19) and our conjecture (1.9.18) are 
so suggestive as to 
lead us to make the further conjecture that, for 
$\Gamma =\Gamma_0(q_0)\leq SL(2,{\Bbb Z}[i])$ and $p\in{\Bbb Z}$, one 
has 
$$\sum_{\scriptstyle V\,:\,p_V=p\atop\scriptstyle |\nu_V|\leq K} 1 
\sim {{\rm vol}(\Gamma\backslash G)\over 6\pi^2}\,K^3\qquad\quad 
{\rm as}\quad K\rightarrow\infty\;.\eqno(1.9.20)$$ 
This conjecture is at least partly correct: for the case $p=0$ of 
(1.9.20) is a known instance in which Weyl's law holds (see 
[11], Page~308 and 
Section~8.9, for details of this). Moreover, what is hypothesised in (1.9.20)  
is equivalent to something that is (at least) superficially analogous to what 
M\"{u}ller has proved in Theorem~0.1 of [37]. In fact, on the 
basis of the conjecture (1.9.18), 
and the form of the result which we obtain in (5.55) (below), 
we venture to put forward the conjecture that
as $\delta\rightarrow 0+$ and $( p^2+K^2)\delta K\rightarrow\infty$,  
with $p\in{\Bbb Z}$ and $\delta,K\in (0,\infty)$, one has 
$$\sum_{\scriptstyle V\,:\,p_V=p\atop\scriptstyle 
K<|\nu_V|\leq K+\delta K} 1\sim {{\rm vol}(\Gamma\backslash G)\over 2\pi^2}\, 
\left( p^2+K^2\right)\delta K\eqno(1.9.21)$$ 
uniformly for all  $\Gamma=\Gamma_0(q)\leq SL(2,{\Bbb Z}[i])$. 
\par 
It is a significant feature of Theorem~1 
that the bounds (1.9.14) hold for all cusps ${\frak a}$ of $\Gamma$: 
a direct application of this feature occurs in the proof of   
Theorem~4 of [46], and indirect use of it is essential to 
the proof of Theorem~10 of [46]. However, in order to compare 
Theorem~1 with Lokvenec-Guleska's result (1.9.19), we now focus for a moment 
on the case ${\frak a}=\infty$. If one were deprived of all information 
concerning the Fourier coefficients $C_V^{\infty}(\omega ;\nu_V,p_V)$ other 
than what is stated in (1.9.19), then essentially the strongest upper bound for 
$E_0^{\infty}(q_0,P,K;N,b)$ that one could deduce would be the 
result 
$$\limsup_{K\rightarrow\infty} K^{-3} E_0^{\infty}(q_0,P,K;N,b) 
\leq {1\over 6\pi^2}\,(2P+1)(N+o(N)) \left\|{\bf b}_N\right\|_2^2\;\eqno(1.9.22)$$ 
(in which $q_0\in{\frak O}-\{ 0\}$, $P,N\in{\Bbb N}$ and the function 
$b : {\frak O}-\{ 0\}\rightarrow{\Bbb C}$ are presumed given). 
In fact the inequality (1.9.22) is simply what follows from (1.9.12) and (1.9.19) by means 
of the Cauchy-Schwarz inequality,    
$$\biggl|\sum_{\scriptstyle\omega\in{\frak O}\atop\scriptstyle N/2<|\omega|^2\leq N}
b(\omega) C_V^{\frak a}\left(\omega;\nu_V,p_V\right)\biggr|^2
\leq\ \left\|{\bf b}_N\right\|_2^2 
\!\!\sum_{\scriptstyle\omega\in{\frak O}\atop\scriptstyle N/2<|\omega|^2\leq N}
\left| C_V^{\frak a}\left(\omega;\nu_V,p_V\right)\right|^2\;,$$ 
and known estimates (either classical or more recent, as in [18]) for the number of 
lattice points lying in a disc of specified radius. The case $j=0$ of (1.9.14) implies 
that the inequality (1.9.22) remains valid 
if, in place of the factor $N+o(N)$ on the right-hand side of that inequality, 
one substitutes just a factor~$O(1)$; this reveals that 
the extent of `cancellation' occurring within the sums over $\omega$ 
on the right-hand side of Equation~(1.9.12) is, with few exceptions, 
comparable to that which one would expect to find, in the limit as $N\rightarrow\infty$, 
in respect of a sum 
$$R({\cal X};N,b)=\sum_{\scriptstyle\omega\in{\frak O}\atop\scriptstyle 
N/2<|\omega|^2\leq N} b(\omega)\exp\left( i X_{\omega}\right)$$ 
in which $b(\omega)$ is a given complex-valued function,  
while ${\cal X}=(X_{\omega})_{\omega\in{\frak O}-\{ 0\}}$ is a family 
of independent real-valued random variables such that, for 
all $\omega\in{\frak O}-\{ 0\}$, the mean value of $\exp( iX_{\omega})$ is equal to zero. 

\bigskip\par 

\noindent{\bf Remark~1.9.9 (on improving on Theorem~1).}\quad 
The bounds in (1.9.14) are less precise than  
the asymptotic result (Equation~(5.55), below), from which they are deduced. 
Moreover, in view of the 
positivity (exploited in (5.2) and (5.61), below) 
of the terms summed in (1.9.12) and (1.9.13), it is evident 
that Theorem 1 is not optimal when 
$PK=o\left( N^{2/3} |\mu({\frak a})|^{4/3}\right)$
and either $P=o(K)$ or $K=o(P)$. Indeed, by virtue of the positivity of those terms,   
it is a direct corollary of Theorem~1 itself that the 
result (1.9.14) may be improved to:  
$$E_{j}^{\frak a}\left( q_0 ,P,K;N,b\right)
\ll\left(\left( P^2+K^2\right)\left( PK+
O_{\varepsilon}\!\bigl( {\cal Y}^{2/3}\bigr)\right)  
+O_{\varepsilon}\!\bigl( (P+K){\cal Y}\bigr)\right)
\left\|{\bf b}_N\right\|_2^2
\qquad\ \hbox{($j=0,1$)}\;,\eqno(1.9.23)$$
where ${\cal Y}={\cal Y}(N,\varepsilon ; q_0 ,{\frak a})=N^{1+\varepsilon} |\mu({\frak a})|^2$. 
On the basis of the conjecture (1.9.18), the conjecture (1.9.21) and 
the tendency towards cancellation observed in the final paragraph of Remark~1.9.8, 
we are led to make the further conjecture that Theorem~1 would remain 
true if the bound  
$$E_0^{\frak a}\left( q_0 ,P,K;N,b\right)
\ll\left(\left( P^2+K^2\right) PK 
+{N\over{\rm vol}(\Gamma\backslash G)}\right)
\left\|{\bf b}_N\right\|_2^2$$
were substituted in place of the result in the case $j=0$ of (1.9.14). Therefore we  
believe that even the case $j=0$ of (1.9.23)  
falls distinctly short of being best-possible. 
It also seems likely that 
Theorem 1 (or its corollary, (1.9.23)) could be improved upon in other ways:        
one might, for example, be able to prove a `short-spectral-interval' refinement of 
the case $j=0$ of (1.9.14) (i.e. a result analogous to the results 
obtained in Theorem~1.1 of [23], 
Theorem~3.3 of [35] and Lemma~7 of [24]); one might also be able to obtain, 
instead of just an upper bound,  
an asymptotic estimate for the sum $E_{0}^{\frak a}\left( q_0 ,P,K;N,b\right)$. 

We expect that the bound in the case $j=1$ of the corollary (1.9.23)  
(and hence also the bound in the case $j=1$ of (1.9.14)) can  
be significantly sharpened by exploiting the special nature of the 
modified Fourier coefficients $B_{\frak c}^{\frak a}(\omega;\nu,p)$ 
of the Eisenstein series: see
Equation~(10.1) of Theorem~10.1 of [5] for 
the simple explicit form that these coefficients take  
in the case $q_0 =1$. The discussion around (1.8.7) casts 
further light on this interesting possibility. We have not ourselves attempted 
to improve upon our stated bounds for  
$E_1^{\frak a}(q_0,P,K;N,b)$: it turns out that 
Theorem~1, as it stands, is adequate for our needs in [46] and [47]. 
\par 
Given what Lokvenec-Guleska achieved in her thesis [32], we are almost  
certain that methods similar to those of the present paper 
are capable of yielding, when  
$F\in\{ {\Bbb Q}(-d) : d\in{\Bbb N},\,d>1\ {\rm and}\ d\ {\rm is\ squarefree}\}$  
and $\Gamma$ is any Hecke congruence subgroup of $SL(2,{\frak O}_F)\,$ 
(where ${\frak O}_F$ is equal to the ring of integers of $F$),  
both an analogue for $\Gamma\backslash SL(2,{\Bbb C})$ of  
the summation formula  in 
Theorem~B and results in some (useful) sense analogous to those in Theorem~1. 
Although such generalisations of Theorem~B and Theorem~1 would be of 
considerable number-theoretical interest, we have (so far) 
preferred not to do any work in that direction ourselves: 
we found the case $F={\Bbb Q}(\sqrt{-1})$ to be complicated enough.   

\bigskip\par

\noindent{\bf Remark~1.9.10 (on related work of the author).}\quad 
This paper is the first in a planned 
series of three: the other two being [46] and [47]. 
Our work in [46] involves the application of Theorem~1, in combination 
with the result (1.4.15) of Kim and Shahidi 
and the the inversion of Theorem~B mentioned in Remark~1.9.4; the results 
obtained there include new upper bounds for sums of the form 
$$\sigma_{\Gamma}^{\frak a}({\bf b},N;X) 
=\sum_{\scriptstyle V\atop\scriptstyle \nu_V^2>0 ,\ p_V=0} X^{\left|\nu_V\right|} 
\left|\sum_{\scriptstyle n\in{\frak O}\atop\scriptstyle 
N/2<|n|^2\leq N} b_n C_V^{\frak a}\left( n ; \nu_V , 0\right)\right|^2\qquad\qquad 
\hbox{(${\frak a}\in{\Bbb Q}(i)\cup\{\infty\}$)}$$ 
and 
$$S(Q,X,N) 
=\sum_{\scriptstyle q_1\in{\frak O}\atop\scriptstyle Q/2<\left| q_1\right|^2\leq Q}
\sigma_{\Gamma_0\!\left( q_1\right)}^{\infty}({\bf b},N;X) \;,$$ 
where $X$ and $N$ denote real numbers satisfying $X\geq 1$, $N\geq 1$, while the 
coefficients $b_n\,$ ($0\neq n\in{\frak O}$) are complex numbers that 
are allowed to be arbitrary in [46], Theorem~4, Theorem~5, Theorem~6 and Theorem~7, but  
are required to be of a special type in [46], Theorem~8 and Theorem~9. 
Our results in [46] include analogues, for $SL(2,{\Bbb C})$, of 
several of the ($SL(2,{\Bbb R})$ related) results obtained by 
Deshouillers and Iwaniec in [9], 
and also the `$SL(2,{\Bbb C})$ analogue' of Theorem~2 of~[45]. 
\par 
In work only partially written up we have made use of 
the results of [46] in estimating a particular type of 
weighted fourth power moment for 
the family of Hecke zeta functions $\bigl( \zeta(s,\lambda^k)\bigr)_{k\in{\Bbb Z}}\,$ 
given, for ${\rm Re}(s)>1$, by the cases of (1.8.7) in which 
$\chi : {\frak O}\rightarrow\{ 1\}$ and $p$ is even.  
It was shown by Hecke, in [17], that these 
zeta functions have a meromorphic continuation to all of ${\Bbb C}$, with 
no poles except at $s=1$ (and no pole there except when $k=0$). 
By using the results of [46] we have been able to prove that, 
for arbitrary complex coefficients 
$a_n\,$ ($0\neq n\in{\frak O}$), all $\varepsilon >0$ and 
all $D,N\in{\Bbb N}$, one has 
$$\sum_{k=-D}^D\ \int\limits_{-D}^D\left|\zeta\left(\textstyle{1\over 2}+it,\lambda^k\right) 
\right|^4 \Biggl|\sum_{\scriptstyle n\in{\frak O}\atop\scriptstyle 
0<|n|^2\leq N}a_n\lambda^k(n)|n|^{-2it}\Biggr|^2 {\rm d}t 
\ll_{\varepsilon} D^{2+\varepsilon} N\max_{0<|n|^2\leq N}\left| a_n\right|^2\qquad\ 
\hbox{if $\ N^2\leq D$.}$$ 
We are preparing our proof of this result for publication, and hope that it 
may appear in [47]. 

\bigskip 

\proclaim Proposition~2. 
Let $\varepsilon$, $N$, $q_0$ and $b$ satisfy the 
same hypotheses as in Theorem 1, 
and let $\left\|{\bf b}_N\right\|_2$ be given by (1.9.16); 
let $A_1$ and $A_2$ be arbitrary positive absolute constants; 
let $\psi\in{\Bbb R}$ and $M\in{\Bbb N}\cup\{ 0\}$, 
and let ${\frak a}$ be a cusp of $\,\Gamma=\Gamma_0\!\left( q_0\right)\leq SL(2,{\frak O})$; 
let  the `scaling matrix' $g_{\frak a}\in G=SL(2,{\Bbb C})$ be such as to satisfy, 
for ${\frak c}={\frak a}$, each of (1.1.16), (1.1.20) and (1.1.21),  
and let $c$ be an element of the set ${}^{\frak a}{\cal C}^{\frak a}\subset{\Bbb C}^{*}$ 
defined in (1.5.8), (1.5.9) 
(i.e. with ${\frak a}'={\frak a}$ and $g_{{\frak a}'}=g_{\frak a}$ there).  
Then 
$$0\neq c\in {1\over\mu({\frak a})}\,{\frak O}\eqno(1.9.24)$$
(the relationship between $\mu({\frak a})$ and ${\frak a}$ being as described in 
Remark~1.9.7, above), and the sum 
$$U_{\frak a}(\psi,c;M;N,b)=\sum_{m=-M}^{M}
\left|\quad\ \;  
\sum\!\!\!\!\!\!\!\!\!\!\!\sum_{\!\!\!\!\!\!\!\!\!\!\!{\scriptstyle \omega_1,\omega_2\in{\frak O}\atop\scriptstyle 
N/2<\left|\omega_1\right|^2,\left|\omega_2\right|^2\leq N}}
\!\!\!\!\!\overline{b\left(\omega_1\right)}b\left(\omega_2\right) 
\left({\omega_1\omega_2\over\left|\omega_1\omega_2\right|}\right)^{\!\!m}
S_{{\frak a},{\frak a}}\!\left(\omega_1 , \omega_2 ; c\right) 
{\rm e}\!\left(
{\psi\sqrt{\left|\omega_1\omega_2\right|}\over |c|}\right)\right| \eqno(1.9.25)$$
(with $S_{{\frak a},{\frak a}'}(\omega,\omega';c)$ defined as in (1.5.8)-(1.5.10)) 
satisfies both 
$$U_{\frak a}(\psi,c;M;N,b) 
\ll \tau^{3/2}(c) |c|(M+1)N\left\|{\bf b}_N\right\|_2^2\eqno(1.9.26)$$ 
and 
$$U_{\frak a}(\psi,c;M;N,b) 
\ll (1+|\psi|)^{1/2}\left( \left| c\right| (M+1)+N^{1/2}\right)
\left(\left| c\right|+N^{1/2}\right) \left\|{\bf b}_N\right\|_2^2\;,\eqno(1.9.27)$$ 
where $\tau(c)$ is the number of Gaussian integer divisors of $c$. 
\hfill\break $\hbox{\qquad}$If it is moreover the case that  one has 
$$0<|c|^2\leq A_1 N^{1-\varepsilon}\eqno(1.9.28)$$
and 
$$0<|\psi|\leq A_2\eqno(1.9.29)$$
then  
$$U_{\frak a}(\psi,c;M;N,b) 
\ll_{A_1,A_2}\,\left( \left|\psi\right|^{-1/2} +O_{\varepsilon}\!\left( 1\right)\right) 
\left( |c|^{1/2} N^{3/4}+|c|^{3/2} M N^{1/4}\right)
N^{\varepsilon}\left\| {\bf b}_N\right\|_2^2\;.\eqno(1.9.30)$$

\bigskip

Proposition~2 is proved in Section 3: lemmas necessary for its proof 
are collected in Section 2.

\bigskip

\centerline{\bf Notation.}

\medskip

The following index of notation covers only some of the more unusual 
notation used in this paper; it is not comprehensive. 
Some of our other terminology is explained  
in five supplementary paragraphs. 
  
\medskip 

\noindent{\bf Index of notation:} 

\medskip 

{\settabs\+$E_j^{\frak a}(q_0 ,P,K;N,{\bf b})\,$\  
&a complete set of representatives of the $\Gamma$-equivalence classes of cusps\quad\  
&\cr %sample line 
\+{\bf Symbol}&{\bf Description}&{\bf Place defined}\cr 

\smallskip 

\+$|{\cal A}|$  
&{\it (when ${\cal A}$ is a set): the cardinal number of ${\cal A}$ 
}&-- \cr 
\+$|f|$  
&{\it (when $f$ is a complex valued function of $x$): the function 
$x\mapsto |f(x)|$
}&-- \cr 
\+$R^{*}$  
&{\it (when $R$ is a ring with an identity): the group of units of $R$
}&-- \cr 
\+$\overline{U}$
&{\it (when $U$ is a subset of a metric or topological space): 
the closure of $U$
}&-- \cr  
\+$(g\circ f)(x)$  
&{\it equal to $g(f(x))$ 
}&-- \cr 
\+$a\cdot b$  
&{\it $a$ multiplied by $b\,$ ($\,(a\cdot b)(x)=a(x)b(x)$ if $a$ and $b$ are  
functions of $x$) 
}&-- \cr 
\+${\bf v}\cdot{\bf w}$
&{\it equal to $v_1 w_1+\cdots+v_n w_n$, 
the inner product of vectors ${\bf v},{\bf w}\in{\Bbb C}^n$ 
}&--\cr
\+$|z|$ and $\overline{z}$   
&{\it the absolute value (or `modulus') and complex conjugate of $z\in{\Bbb C}$
}&-- \cr 
\+${\frak a}\sim^{\!\!\!\!\Gamma}{\frak b}$
&{\it the relation of $\Gamma$-equivalence (for cusps ${\frak a}, {\frak b}$)
}&above (1.1.16)\cr 
\+$m\mid n$
&{\it (when $m,n\in{\frak O}$): the relation `$n$ is divisible by $m$'  
}&--\cr
\+$m\sim n$
&{\it (when $m,n\in{\frak O}$): the relation `$n$ is an associate of $m$' 
}&above (1.1.1)\cr 
\+$(m_1,\ldots,m_n)$
&{\it a highest common factor (of $m_1,m_2,\ldots ,m_n\in{\frak O}$)
}&--\cr 
\+$(C,q_0^{\infty})$  
&{\it a certain factor of the non-zero Gaussian integer $C$
}&below (6.1.17) \cr 
\+$a\equiv b\,\bmod c{\frak O}$  
&{\it equivalent to the statement that one has 
$a,b,c\in{\frak O}$ and $c\mid (b-a)$
}&-- \cr 
\+$[x]$
&{\it the greatest rational integer less than or equal to $x$
}&--\cr 
\+$\|\beta\|$
&{\it equal to $\min\{ |n-\beta| : n\in{\Bbb N}\}$
}&--\cr 
\+$\|{\bf b}_N\|_2$
&{\it the Euclidean norm of a vector involving coefficients $b_n$ ($0\neq n\in{\frak O}$) 
}&in (1.9.16)\cr 

\smallskip 

\+$\int\limits_{(\alpha)}f(z) {\rm d}z$
&{\it contour integral along the straight line from $\alpha -i\infty$ to 
$\alpha +i\infty$ 
}&--\cr 
\+$\langle f,h\rangle_{\Gamma\backslash G}$ 
&{\it the inner product of $\,f,g\in L^2(\Gamma\backslash G)$ 
}&in (1.2.2)  \cr 
\+$\| f\|_{\Gamma\backslash G}$  
&{\it  a norm on $L^2(\Gamma\backslash G)$, 
equal to $\sqrt{\langle f , f\rangle_{\Gamma\backslash G}}$
}&--  \cr 
\+$(f_1,f_2)_K$ 
&{\it the inner product for square integrable functions $f_1, f_2 : K\rightarrow{\Bbb C}$ 
}&in (1.2.22)  \cr 
\+$\|\Phi\|_K$ 
&{\it a norm, equal to $\sqrt{(\Phi,\Phi)_K}$  
}&below (6.4.2) \cr 
\+$\langle f,F\rangle_{N\backslash G}$  
&{\it  a certain inner product, defined when $f\,\overline{F}\in L^1(N\backslash G)$ 
}&in (6.2.9) \cr 
\+$(f_1,f_2)_{\rm ps}$  
&{\it the inner product for the `principal series' 
}&(1.6.2), (1.6.3) \cr 
\+$\|\varphi\|_{\rm ps}$  
&{\it a norm on the space $H^2(\nu ,p)$, when $(\nu,p)\in (i{\Bbb R})\times{\Bbb Z}$
}&below (1.6.4) \cr 
\+$(f_1,f_2)_{\rm cs}$  
&{\it the inner product for the `complementary series' 
}&(1.6.2), (1.6.4) \cr 
\+$\|\varphi\|_{\rm cs}$  
&{\it a norm on the space $H^2(\nu , p)$, when $0<\nu^2<1$
}&below (1.6.4) \cr 
\+$(f|{\frak c})$
&{\it a `left-translate' of $f : G\rightarrow{\Bbb C}$,  
used for Fourier expansion of $f$ at ${\frak c}$ 
}&above (1.4.2)\cr 
\smallskip 

\+$\hat{F}({\bf y})$, $\hat{f}(w)$
&{\it Fourier transforms 
of $F : {\Bbb R}^n\rightarrow{\Bbb C}$ and 
$f : {\Bbb C}\rightarrow{\Bbb C}$
}&(2.44), (2.46)\cr 
\+$\Gamma$, $\Gamma_0(q_0)$
&{\it the Hecke congruence subgroup of $SL(2,{\frak O})$  
of `level' $q_0$  
}&in (1.1.1)\cr
\+$\Gamma_{\frak c}$, $\Gamma_{\frak c}'$
&{\it `stabiliser' and 
`parabolic stabiliser' subgroups (for the cusp ${\frak c}$) 
}&(1.1.17)-(1.1.21)\cr 
\+${}^{\frak a}\Gamma^{\frak b}(c)$
&{\it a `Bruhat cell'
}&in (1.5.8)\cr  
\+$\Gamma(z)$
&{\it Euler's Gamma function, defined for $z\in {\Bbb C}-\{0,-1,-2,\ldots\ \}$ 
}&--\cr 
\+$\gamma$
&{\it most often used to denote an element of the group $\Gamma$    
}&--\cr
\+${\bf\Delta}$ 
&{\it the hyperbolic Laplacian operator on $L^2(\Gamma\backslash{\Bbb H}_3)$
}&in (1.2.13)\cr 
\+${\bf\Delta}_{{\Bbb R}\times{\Bbb R}}$, ${\bf\Delta}_{\Bbb C}$   
&{\it Euclidean Laplacian operators
}&above (2.48)\cr 
\+$\partial/\partial z$, 
$\partial/\partial\overline{z}$
&{\it Complex partial differentiation operators 
}&see (1.2.7)\cr 
\+$\delta^{\Gamma}_{{\frak a},{\frak b}}$ 
&{\it the `delta symbol' for $\Gamma$-equivalence of the cusps ${\frak a}$ and ${\frak b}$ 
}&in (1.8.5)\cr 
\+$\delta_{\omega,\omega'}^{{\frak a},{\frak b}}$
&{\it the `delta symbol' of the `spectral to Kloosterman' sum formula
}&in (1.9.2)\cr 
\+$\delta_{w,z}$
&{\it the `delta-symbol' for equality of the complex numbers $w$ and $z$
}&in (1.5.6)\cr 
\+$\epsilon(p)$
&{\it a certain function defined on ${\Bbb Z}$; takes values in $\{ -1,1\}$  
}&in (6.4.5) \cr 
\+$\zeta(s)$
&{\it the Riemann zeta-function 
}&-- \cr 
\+$\zeta_{{\Bbb Q}(i)}(s)$
&{\it the Dedekind zeta-function for the algebraic number field ${\Bbb Q}(i)$  
}&below (6.5.59) \cr 
\+$\zeta(s,\lambda^{p/2}\chi)$
&{\it a Hecke zeta function (with gr\"{o}ssencharakter) for $\,{\Bbb Q}(i)$ 
}&in (1.8.7) \cr 
\+$\zeta^{{\frak a},{\frak b}}_{\omega,\omega'}(s)$
&{\it a modified Linnik-Selberg series 
}&in (6.5.61) \cr 
\+$\kappa(\omega_1,\omega_2;c)$
&{\it a linear operator from ${\cal T}^{\ell}_{\sigma}$ into ${\cal T}^{\ell}_{\sigma}$ 
}&in (6.4.12) \cr
\+$\lambda^{p/2}$, $\lambda^k$
&{\it a Hecke gr\"{o}ssencharakter on ${\frak O}-\{ 0\}$ 
}&below (1.8.7) \cr 
\+$\lambda_{\nu}$
&{\it an eigenvalue of the operator $-{\bf\Delta}$ on $L^2(\Gamma\backslash{\Bbb H}_3)$ 
}&in (1.4.14)\cr
\+$\lambda^{*}_{\ell}(\nu,p)$
&{\it a certain function on the set ${\Bbb C}\times\{ p\in{\Bbb Z} : |p|\leq\ell\}$ 
}&in (6.6.4) \cr
\+$\mu({\frak a})$
&{\it $1/|\mu({\frak a})|$ is a useful lower bound for 
the set $\{ |c| : c\in{}^{\frak a}{\cal C}^{\frak a}\}$ 
}&in (1.9.15)\cr 
\+$(\nu_V,p_V)$
&{\it the spectral parameters of the cuspidal irreducible subspace 
$V$
}&above (1.7.6)\cr 
\+$\rho : G\rightarrow(0,\infty)$
&{\it $\rho(g)$ is equal to the Iwasawa coordinate $r$ of the element $g\in G$
}&(6.2.1), (1.1.2) \cr

\smallskip 

\+$\sum\limits_V$
&{\it summation over irreducible subspaces $V\subset{}^0L^2(\Gamma\backslash G)$
}&below (1.8.8)\cr 

\smallskip 

\+$\tau(n)$
&{\it (when $n\in{\frak O}$): the number of Gaussian integer divisors of $n$
}&-- \cr
\+$\tau : G\rightarrow{\Bbb C}$
&{\it an element of $C^{\infty}(G)$ chosen according to certain criteria
}&see (6.5.1) \cr
\+$\tau{\bf M}_{\omega}\varphi_{\ell,q}(\nu,0)$
&{\it the product of the functions $\tau$ and ${\bf M}_{\omega}\varphi_{\ell,q}(\nu,0)$ 
}&in (6.5.3) \cr
\+$(1-\tau){\bf M}_{\omega}\varphi$
&{\it (when $\varphi =\varphi_{\ell,q}(\nu,0)$): 
the product of the functions $1-\tau$ and ${\bf M}_{\omega}\varphi$ 
}&below (6.5.3) \cr
\+$\Upsilon_{\nu,p}$
&{\it a certain character for the centre of ${\frak g}$ 
}&in (1.3.3) \cr
\+$\Phi^{\ell}_{p,q}$
&{\it a certain even and square integrable function defined on $K$ 
}&below (1.3.2) \cr
\+$\varphi_{\ell,q}(\nu,p)$
&{\it a certain function lying in the space $C^{\infty}(N\backslash G)$ 
}&in (1.3.2) \cr
\+$\phi(\alpha)$
&{\it (when $0\neq\alpha\in{\frak O}$): Euler's function, equal to 
$|({\frak O}/\alpha{\frak O})^{*}|$ 
}&below (2.5) \cr
\+$\phi_{\omega'}(\nu,g)$
&{\it an analytic continuation of a Fourier term, 
$\,(F^{\frak b}_{\omega'}P^{\frak a}{\bf M}_{\omega}\varphi_{\ell,q}(\nu,0))(g)$ 
}&in (6.5.72) \cr
\+$\Phi(\nu,g)$
&{\it an analytic continuation of the function 
$\,\nu\mapsto P^{\frak a}{\bf M}_{\omega}\varphi_{\ell,q}(\nu,0))(g_{\frak b}g)$ 
}&below (6.5.72) \cr
\+$\phi_1$, $\phi_2$ 
&{\it in $\S 6.6$ these denote certain pseudo Poincar\'{e} series 
}&in (6.6.1)\cr
\+$\chi_4 : {\Bbb N}\rightarrow\{ -1,1\}$ 
&{\it the real primitive Dirichlet character modulo $4$ 
}&-- \cr 
\+$\chi_{q_0} : G\rightarrow\{ 0,1\}$
&{\it the characteristic function of $\,\Gamma_0(q_0)<G$ 
}&in (6.1.9) \cr
\+$\chi^{{\frak d},{\frak d}'}_{\omega,\omega'}(h)$
&{\it the `delta-term' in the spectral sum formula
}&in (6.7.1) \cr
\+$\psi_{\omega} : N\rightarrow{\Bbb C}$
&{\it a certain character for the group $N$ 
}&in (1.4.3)\cr
\+$\psi(y,x;\phi)$
&{\it a certain elementary function 
}&in (4.23) \cr
\+$\omega(c)$
&{\it (when $0\neq c\in{\frak O}$): the number of prime ideals of $\,{\frak O}$ 
containing $c$
}&-- \cr
\+$\Omega_{\pm}$
&{\it the Casimir operators associated with $G$
}&in/below (1.2.7)\cr 
\+$\Omega_{\frak k}$
&{\it the Casimir operator associated with $K$
}&in (1.2.11) \cr 
\+${\bf x}, {\bf X},\ \ldots\ $
&{\it vectors in ${\Bbb R}^n$ or ${\Bbb C}^n$; sets of coefficients; operators; power set 
}&--\cr
\+${\frak a},{\frak b},{\frak c},\ldots\ $
&{\it cusps of $\,\Gamma$, or (more generally) points 
in ${\Bbb P}^1({\Bbb C})={\Bbb C}\cup\{\infty\}$
}&above (1.1.16)\cr 
\+$a[r]$, $A$
&{\it $A=\{ a[r] : r>0\}<G$
}&see (1.1.3)\cr
\+$A^0_{\Gamma}(\Upsilon;\ell,q)$  
&{\it the space of cusp forms in $C^{\infty}(\Gamma\backslash G)$ 
of $K$-type $(\ell,q)$ with character $\Upsilon$ 
}&in (1.4.6) \cr 
\+$A_M(\phi,\theta)$
&{\it a certain trigonometric sum
}&in (4.24) \cr  
\+${\rm Arg}(z)$ 
&{\it the principal argument of $z\in{\Bbb C}$, satisfying $-\pi<{\rm Arg}(z)\leq\pi$
}&-- \cr
\+$B^{+}$
&{\it the subgroup $\{ n[\alpha] : \alpha\in{\frak O}\}<N<G$
}& in (1.1.21)\cr 
\+$B_{\frak a}^{\frak b}\left(\omega;\nu,p\right)$
&{\it a modified Fourier coefficient of an Eisenstein series
}&in (1.8.9)\cr 
\+${\bf B}h : {\Bbb C}^{*}\rightarrow{\Bbb C}$
&{\it the ${\bf B}$-transform of the function 
$\,h : \{\nu\in{\Bbb C} : |{\rm Re}(\nu)|\leq\sigma\}\times{\Bbb Z}\rightarrow{\Bbb C}$ 
}&(1.9.3)-(1.9.6) \cr  
\+$b(\eta)$, $b(\omega;\ell,q;\eta)$  
&{\it (when $\eta$ is a function in the space ${\cal T}^{\ell}_{\sigma}$): 
a complex constant 
}&in (6.4.6) \cr 
\+${\cal B}$  
&{\it a basis for both $\,\frak{sl}(2,{\Bbb C})$ and ${\frak g}$ 
}&below (6.5.26) \cr 
\+${\cal B}_1$  
&{\it a basis for ${\frak g}$ 
}&above (6.5.27) \cr 
\+${\frak C}(\Gamma)$ &{\it a complete set of representatives of 
the $\Gamma$-equivalence classes of cusps
}&--\cr 
\+${}^{\frak a}{\cal C}^{\frak b}$
&{\it a subset of ${\Bbb C}^{*}$, equal to the domain of the mapping  
$c\mapsto S_{{\frak a},{\frak b}}(\omega_1,\omega_2;c)$  
}&in (1.5.9)\cr
\+$C^{\infty}(G)$  
&{\it the space of all smooth functions $f : G\rightarrow{\Bbb C}$
}&after (1.1.9) \cr 
\+$C^{\infty}(\Gamma\backslash G)$  
&{\it the space of all smooth and $\,\Gamma$-automorphic functions on $G$ 
}&in (1.2.3) \cr 
\+$C^{\infty}(\Gamma\backslash{\Bbb H}_3)$  
&{\it the space of all smooth and $\,\Gamma$-automorphic functions on ${\Bbb H}_3$ 
}&in (1.2.15) \cr 
\+$C^{\infty}(N\backslash G,\omega)$  
&{\it (when $\omega\in{\Bbb C}$): a certain subspace of $C^{\infty}(G)$ 
}&(1.4.7), (1.4.3) \cr 
\+$C^0(G)$  
&{\it the space of all continuous functions $f : G\rightarrow{\Bbb C}$
}&after (1.1.9) \cr 
\+$C^0(B^{+}\backslash G)$  
&{\it a certain subspace of $\,C^0(G)$ 
}&start of $\S 6.2$ \cr 
\+$C^0(N\backslash G,\omega)$  
&{\it (when $\omega\in{\frak O}$): a certain subspace of $\,C^0(B^{+}\backslash G)$ 
}&start of $\S 6.2$ \cr 
\+$C^0(\Gamma\backslash G)$ 
&{\it the space of continuous and $\Gamma$-automorphic functions 
$f : G\rightarrow{\Bbb C}$ 
}&-- \cr 
\+$c_V^{\frak c}\left(\omega\right)$ 
&{\it a Fourier coefficient of a cuspidal irreducible subspace $V$ 
}&in (1.7.13)\cr 
\+$C_V^{\frak c}\left(\omega;\nu_V,p_V\right)$
&{\it a modified Fourier coefficient of a cuspidal subspace
}&in (1.7.15)\cr 
\+${\rm d}a$, ${\rm d}k$, ${\rm d}n$  
&{\it left and right Haar measures on the groups $A$, $N$ and $K$, respectively
}&in (1.1.10) \cr 
\+${\rm d}g$
&{\it a left and right Haar measure on $G$
}&in (1.1.11)\cr 
\+${\rm d}Q$
&{\it a $G$-invariant measure on ${\Bbb H}_3$
}&in (1.1.13)\cr 
\+$D_{\frak a}^{\frak b}\left(\omega;\nu,p\right)$
&{\it a Fourier coefficient of an Eisenstein series
}&as in (1.8.4)\cr
\+${\rm d}_{+}z$ 
&{\it the standard Lebesgue measure on ${\Bbb C}$
}&in (1.1.10)\cr 
\+$E_j^{\frak a}(q_0 ,P,K;N,{\bf b})$ 
&{\it a spectral mean, for cusp forms ($j=0$), or 
Eisenstein series ($j=1$)
}&(1.9.12), (1.9.13)\cr 
\+${\rm e}(x)$
&{\it equal to $\exp(2\pi i x)$, a character for the additive group ${\Bbb R}/{\Bbb Z}$
}&--\cr 
\+${\cal E}_{\frak c}$
&{\it a `cusp sector' in ${\Bbb H}_3$
}&in (1.1.23) \cr 
\+$(E^{\frak c}_{\ell,q}(\nu,p))(g)$
&{\it a $\Gamma$-automorphic Eisenstein series associated with the cusp ${\frak c}$ 
}&in $\S 1.8$ \cr  
\+$E$  
&{\it a subset of $\,{\Bbb C}\times{\Bbb Z}$
containing all pairs $(\nu_V,p_V)$ of spectral parameters 
}&in (6.7.13) \cr 
\+${\cal F}$, ${\cal F}_{*}$ 
&{\it fundamental domains for the action of $\Gamma$ upon ${\Bbb H}_3$ 
}&(1.1.7), (1.1.24) \cr 
\+$(F_m^{\frak c} f)(g)$
&{\it the Fourier term of order $m$ for $f$ at ${\frak c}$
}&(1.4.1), (1.4.2)\cr  
\+$F_{P,K}(\eta,\xi)$
&{\it a certain polynomial function 
}&in (4.23) \cr  
\+$g_{\frak c}$ 
&{\it a scaling matrix for the cusp ${\frak c}$
}&(1.1.16)-(1.1.21)\cr
\+$G$, $g$
&{\it the special linear group $SL(2,{\Bbb C})$, and one of its elements 
}&--\cr 
\+$g(a,d;c)$
&{\it a specific element of the group $G$
}&in (6.1.8) \cr  
\+${\frak g}$
&{\it the complex Lie algebra of $G$, 
equal to $\frak{sl}(2,{\Bbb C})\otimes_{\Bbb R}{\Bbb C}$ 
}&above (1.2.6) \cr  
\+$h[u]$
&{$(h[u])_{u\in{\Bbb C}^*}$ is a family of elements of $G$; $h[u]\in K$ when $|u|=1$
}&in (1.1.9)\cr 
\+${\Bbb H}_3$
&{\it the upper half-space model for three dimensional hyperbolic space 
}&in $\S 1.1$ \cr 
\+${\bf H}_2$  
&{\it an element of ${\frak k}\,$ identified with a certain differential 
operator  
}&(1.2.9), (1.2.10) \cr 
\+${\bf h}_{u}$
&{\it (when $u\in{\Bbb C}^{*}$): a left-translation operator on functions 
$f : G\rightarrow{\Bbb C}$ 
}&(1.5.7), (1.1.9)\cr 
\+$H(\nu,p)$
&{\it an irreducible representation space for ${\frak g}$ 
}&in (1.6.1) \cr 
\+$H^2(\nu,p)$
&{\it a Hilbert space obtained as a certain completion of $H(\nu,p)$ 
}&below (1.6.4) \cr 
\+${\cal H}^{\sigma}_0(\varrho,\vartheta)$  
&{\it the space of functions $h$ satisfying conditions (i)-(iii) of Theorem~B 
}&start of $\S 6.7$ \cr 
\+${\cal H}_{\star}^{\sigma}(\varrho,\vartheta)$   
&{\it a certain subspace of ${\cal H}^{\sigma}_0(\varrho,\vartheta)$ 
}&above (6.7.1) \cr 
\+${\rm Im}(z)$
&{\it the imaginary part of the complex number $z\,$ 
(equal to ${\rm Re}(-iz)$)
}&-- \cr  
\+${\rm Int}(U)$
&{\it (if $U$ is a subset of a metric or topological space): 
the interior of $U$
}&-- \cr  
\+$J_{\nu}(z)$, $J_n(z)$
&{\it a Bessel function 
}&(1.9.8), (1.9.6)\cr  
\+$J_{\nu}^{*}(z)$
&{\it equal to $(z/2)^{-\nu}J_{\nu}(z)$ when $z>0$   
}&in (1.9.6)\cr 
\+${\cal J}_{\mu,k}(z)$, ${\cal K}_{\nu,p}(z)$ 
&{\it functions related to Bessel functions of representations of 
$PSL(2,{\Bbb C})$ 
}&(1.9.5), (1.9.4)\cr 
\+${\cal J}^{*}_{\nu,p}(z)$ 
&{\it a function closely related to ${\cal J}_{\nu,p}(z)$ 
}&in (6.3.12) \cr 
\+$({\bf J}_{\omega}\varphi_{\ell,q}(\nu,p))(g)$  
&{\it a Jacquet integral  
}&in $\S 1.5$\cr 
\+${\bf J}_{\omega}$, ${\bf J}_{\omega}^{\nu,p}$ 
&{\it (when $\omega\in{\Bbb C}$, $\nu\in{\Bbb C}$ and $p\in{\Bbb Z}$): a Jacquet operator 
}&in $\S 1.5$, $\S 1.6$\cr 
\+$K$, $k[\alpha ,\beta]$
&{\it the special unitary group, $SU(2)<G$, and one of its elements  
}&see (1.1.3)\cr 
\+$K^{+}$
&{\it a fundamental domain for $\{ h[1] , h[-1]\}\backslash K$
}&-- \cr 
\+$K$-type $(\ell,q)$  
&{\it classifies elements of $\,C^{\infty}(G)$ or 
$C^{\infty}(K)$ satisfying certain P.D.E.s  
}&see (1.3.1) \cr 
\+${\frak k}$
&{\it the complex Lie algebra of $K$, 
equal to $\frak{su}(2)\otimes_{\Bbb R}{\Bbb C}$ 
}&above (1.2.9) \cr  
\+$\log(x)$
&{\it equal to $\log_{e}(x)$, the natural logarithm
}&--\cr 
\+$L(s,\chi)$
&{\it Dirichlet's $L$-series (with Dirichlet character $\chi$) 
}&-- \cr 
\+$L^2(\Gamma\backslash G)$  
&{\it the Hilbert space of square-integrable $\Gamma$-automorphic functions on $G$ 
}&(1.2.1), (1.2.2) \cr 
\+$L^2(\Gamma\backslash{\Bbb H}_3)$  
&{\it the Hilbert space of square-integrable $\Gamma$-automorphic functions on ${\Bbb H}_3$ 
}&(1.2.4), (1.2.5) \cr 
\+${}^0L^2(\Gamma\backslash G)$  
&{\it the closure of the subspace of 
$L^2(\Gamma\backslash G)$ spanned by cusp forms  
}&below (1.7.3) \cr 
\+${}^{\rm e}L^2(\Gamma\backslash G)$  
&{\it a subspace of $L^2(\Gamma\backslash G)$ spanned by 
mean values of Eisenstein series 
}&below (1.7.3) \cr 
\+$L^p(\Gamma\backslash G)$  
&{\it a certain space of measurable and 
$\Gamma$-automorphic functions on $G$
}&below (6.5.7)\cr 
\+$L^p(\Gamma\backslash G;\ell,q)$  
&{\it a subspace of $L^p(\Gamma\backslash G)$ characterised by the 
$K$-type $(\ell,q)$
}&after (6.5.89) \cr 
\+$L^{\infty}(\Gamma\backslash G)$  
&{\it the space of essentially bounded elements of $L^1(\Gamma\backslash G)$
}&above (6.5.6) \cr 
\+$L^{\infty}(\Gamma\backslash G;\ell,q)$  
&{\it a subspace of $L^{\infty}(\Gamma\backslash G)$ characterised 
by the $K$-type $(\ell,q)$ 
}&in (6.5.6) \cr 
\+$L^p(N\backslash G)$  
&{\it (when $1\leq p<\infty$): a certain space of measurable functions on $G$ 
}&before (6.6.5) \cr 
\+$({\bf L}^{\omega}_{\ell,q} f)(\nu,p)$ 
&{\it the Lebedev transform of $f\in P_{\ell,q}(N\backslash G,\omega)$ 
}&in (6.4.2) \cr 
\+$(\widetilde{\bf L}^{\omega}_{\ell,q}\eta)(g)$ 
&{\it the `inverse Lebedev transform' of $\eta\in{\cal T}^{\ell}_{\sigma}$ 
}&(6.4.4), (6.4.5) \cr 
\+$(\widetilde{\bf L}^{\omega,*}_{\ell,q}\eta)(g)$ 
&{\it first modification of the `inverse Lebedev transform' of $\eta\in{\cal T}^{\ell}_{\sigma}$ 
}&(6.5.1)-(6.5.3) \cr 
\+$(\widetilde{\bf L}^{\omega,\dagger}_{\ell,q}\eta)(g)$ 
&{\it second modification of the `inverse Lebedev transform' of $\eta\in{\cal T}^{\ell}_{\sigma}$ 
}&in (6.5.20) \cr 
\+$m_{\frak c}$  
&{\it a non-zero Gaussian integer: $|m_{\frak c}|^2$ is the `width' of 
the cusp ${\frak c}$
}&below (1.1.22) \cr 
\+${\bf M}_{\omega}$, ${\bf M}^{\nu,p}_{\omega}$
&{\it the Goodman-Wallach operator on the space $H(\nu,p)$ 
}&(6.3.1), (6.3.2) \cr 
\+$N$, $n[z]$ 
&{\it $N=\{ n[z] : z\in{\Bbb C}\}<G$ 
}&see (1.1.3)\cr 
\+${\frak O}$
&{\it equal to ${\Bbb Z}[i]$, the ring of integers of the Gaussian number field 
${\Bbb Q}(i)$
}&--\cr
\+$P$
&{\it the group of those elements of $G$ that are upper triangular matrices  
}&in (1.1.18) \cr  
\+${\Bbb P}^1({\Bbb C})$
&{\it a projective line, identified with the Riemann sphere, ${\Bbb C}\cup\{\infty\}$
}&--\cr
\+${\Bbb P}^1\bigl( {\Bbb Q}(i)\bigr)$
&{\it a projective line, identified with ${\Bbb Q}(i)\cup\{\infty\}$, the set of all cusps
}&above (1.1.16)\cr  
\+$P^{\frak a}f$, $P^{\frak a}f_{\omega}$   
&{\it a Poincar\'{e} series associated with the cusp ${\frak a}$  
}&in (1.5.4) \cr 
\+$P^{{\frak a},*}\widetilde{\bf L}^{\omega}_{\ell,q}\eta$  
&{\it a certain pseudo Poincar\'e series 
}&in (6.5.5) \cr 
\+$P_{\ell,q}(N\backslash G,\omega)$ 
&{\it a subspace of $C^{\infty}(N\backslash G,\omega)$: its elements 
satisfy growth conditions
}&in (6.4.1) \cr 
\+${\cal P}^{\frak a}_{\!\!\!\scriptscriptstyle\hookleftarrow} 
{\bf M}_{\omega}\varphi_{\ell,q}(\nu,0)$  
&{\it an analytic continuation of 
$P^{\frak a}{\bf M}_{\omega}\varphi_{\ell,q}(\nu,0)$
}&(6.5.82)-(6.5.84) \cr 
\+${\cal P}^{\frak a}_{\!\!\!\scriptscriptstyle\hookleftarrow} 
\tau{\bf M}_{\omega}\varphi_{\ell,q}(\nu,0)$  
&{\it an analytic continuation of 
$P^{\frak a}\tau{\bf M}_{\omega}\varphi_{\ell,q}(\nu,0)$
}&see (6.5.89) \cr 
\+$q_0$   
&{\it the `level' of the Hecke congruence subgroup $\Gamma\leq SL(2,{\frak O})$  
}&(1.1.1), (6.5.96) \cr 
\+${\rm Re}(z)$
&{\it the real part of the complex number $z$
}&-- \cr  
\+${\cal R}_{\frak c}$
&{\it a compact subset of $\,{\Bbb C}$ associated with the group $\Gamma_{\frak c}$  
}&in (1.1.22) \cr 
\+${\cal R}(\sigma_1,\sigma_2,t_1)$  
&{\it a certain closed rectangular region in the complex plane 
}&above (6.5.73)\cr 
\+$S(\omega_1,\omega_2;c)$
&{\it a `simple' (or `classical') Kloosterman sum
}&in (2.16)\cr
\+$S_{{\frak a},{\frak b}}\left(\omega_1 , \omega_2 ; c\right)$
&{\it a generalised Kloosterman sum
}&(1.5.8)-(1.5.10)\cr  
\+$T_V\varphi_{\ell,q}(\nu_V,p_V)$
&{\it a $\Gamma$-automorphic cusp form on $G$ 
}&in $\S 1.7$ \cr  
\+${\cal T}^{\ell}_{\sigma}$ 
&{\it a space of `test functions' defined on a 
subset of $\,{\Bbb C}\times{\Bbb Z}$  
}&below (6.4.3) \cr 
\+${\cal U}({\frak g})$
&{\it the universal enveloping algebra of ${\frak g}$   
}&-- \cr  
\+${\cal U}({\frak k})$
&{\it the universal enveloping algebra of ${\frak k}$ 
}&above (1.2.10) \cr  
\+$U_{\frak a}(\psi,c;M;N,b)$
&{\it a sum involving certain sums of Kloosterman sums 
}&in (1.9.25) \cr  
\+${\rm vol}(\Gamma\backslash G)$
&{\it the covolume of $\,\Gamma$ in $G$
}&(1.1.14), (1.1.15) \cr  
\+$V$
&{\it an irreducible cuspidal subspace of\ $\,{}^0L^2(\Gamma\backslash G)$
}&below (1.7.4)\cr 
\+$V_{K,\ell,q}$ 
&{\it a one dimensional subspace of $V$
}&(1.7.5)-(1.7.8)\cr 
\+$W_{\omega}(\Upsilon;\ell,q)$ 
&{\it a subspace of $C^{\infty}(N\backslash G,\omega)$
}&in (1.4.8) \cr 
\+$X^{{\frak d},{\frak d}'}_{\omega,\omega'}(h)$  
&{\it the sum of Kloosterman sums occurring in the spectral sum formula
}&in (6.7.2) \cr 
\+$Y^{{\frak d},{\frak d}'}_{\omega,\omega'}(h)$  
&{\it equal to the `spectral side' of the spectral sum formula 
}&in (6.7.3) \cr 
\+${\cal Z}({\frak g})$, ${\cal Z}({\frak k})$ 
&{\it the centres of ${\cal U}({\frak g})$ and ${\cal U}({\frak k})$, respectively 
}&-- \cr  
}

\medskip

\proclaim Other Algebraic Notation.  
When $U$, $V$ and $W$ are groups, the notation 
$U\leq W$ (resp. $U<W$) is used to indicate 
that $U$ is a subgroup (resp. proper subgroup) of $W$. 
If $U$ and $V$ are subgroups of the group $W$, then 
$W/V$, $U\backslash W$ and $U\backslash W/V$ denote 
the relevant sets of left cosets, right cosets and double cosets (respectively); 
and $\bigl[ W : U\bigr]$ denotes the index of $U$ in $W$, 
so that $\bigl[ W : U\bigr]=|W/U|$.
This notation for `quotients', such as $U\backslash W$ and $U\backslash W/V$, 
may apply in more general contexts. For example, if $U$ is a subgroup of $W$, and if 
$S$ is a subset of the elements of the group $W$ such that $uS\subseteq S$ 
for all $u\in U$, then $S$ can be expressed as a disjoint union of 
certain of the right cosets of $U$ in $W$, and so the notation $U\backslash S$ 
makes sense (as shorthand for the set of right cosets occurring in 
that disjoint union). Similar considerations apply in the 
case of quotients $S/V$ and $U\backslash S/V$, provided that the 
set $S$  is suitably invariant 
(either under left-multiplication by elements of $U\leq W$, or under 
right-multiplication by elements of the group $V\leq W$). 
 
\proclaim Other number-theoretic notation.  
In relations such as $h m^*\equiv\ell\bmod c{\frak O}$, 
or in expressions such as the highest common factor $(h m^* , c)$, 
the rational expression $h m^*/c$, or 
(see (2.16)) the `simple Kloosterman sum' $S(h m^*,\ell;c)$,  
it is to be understood that $m^*$ denotes an
arbitrary element of ${\frak O}$ satisfying 
$m m^*\equiv 1\bmod c{\frak O}$.  
It is therefore implicit in such expressions that one has both $(m,c)\sim1$ and 
$(m^*,c)\sim 1$. 
\hfill\break $\hbox{\qquad}$  
We use the square-brackets notation $[m,n]$ to denote a 
highest common factor of the Gaussian integers $m$ and $n$ (in the event  
that $m$ and $n$ are real it should be clear from the context whether or 
not $[m,n]$ instead denotes a real interval).  
When $\varpi$ is a Gaussian prime and $n$ is a non-negative integer 
the relation 
$\varpi^n \| c$ holds if and only if one has both 
$\varpi^n\mid c$ and $\varpi^{n+1}\!\not\,\mid c$. 
 
\proclaim Summation related conventions.   
When a condition of the form `$m\bmod c{\frak O}$' 
appears below the summation sign, 
it is to be understood that the variable of summation $m$ ranges 
(to the extent permitted 
by any other conditions of summation) 
over some set of representatives of the cosets of 
$c{\frak O}$ in ${\frak O}$. Conditions of summation 
such as $\gamma\in\Gamma_{\frak a}'\backslash\Gamma\,$ (or 
$\gamma\in 
\Gamma_{\frak a}'\backslash{}^{\frak a}\Gamma^{\frak b}(c)/\Gamma_{\frak b}'$)  
indicate that the variable of summation $\gamma$ ranges over 
some set of representatives of the 
relevant cosets (or double cosets); this is an abuse of commonly accepted  
group-theoretic notation, insofar as that in these instances $\gamma$ does not 
itself denote a coset, or double coset. 
\hfill\break $\hbox{\qquad}$
Where there is no indication to the contrary, variables of 
summation range over all values in ${\frak O}$ 
consistent with whatever conditions are attached.

\proclaim Square roots. Our use of the square root sign is mildly ambiguous. 
When the context implies that one has $z>0$, it is then to be understood that 
$\sqrt{z}$ denotes the positive square root of $z\,$ (see (6.6.47) for one such instance);    
otherwise $\sqrt{z}$ denotes an arbitrary solution $w\in{\Bbb C}$ of the 
equation $w^2=z\,$ (this being the case in (1.9.1), (2.10), (6.1.26), 
(6.3.11), (6.5.18) and Lemma~6.6.6, for example).  
 
\proclaim Notation for bounds and asymptotic estimates. 
Where $B\geq 0$, we use the notation $O_{\alpha_1 ,\ldots,\alpha_n}(B)$  
to denote    
a complex-valued variable $\beta$  
satisfying a condition of the form $|\beta|\leq C(\alpha_1,\ldots,\alpha_n) B$, in which  
the `implicit constant' 
$C(\alpha_1,\ldots,\alpha_n)$ is positive and 
depends only on previously declared constants and 
$\alpha_1,\ldots,\alpha_n$.
As alternatives to an expression of the form 
`$\xi =O_{\alpha_1 ,\ldots,\alpha_n}(B)$', we may prefer
to follows Vinogradov in using either `$\xi\ll_{\alpha_1 ,\ldots,\alpha_n} B$', 
or `$B\gg_{\alpha_1 ,\ldots,\alpha_n} \xi$'. 
Where $A\geq 0$ and $B\geq 0$, the notation 
$A\asymp_{\alpha_1 ,\ldots,\alpha_n}B$ may be used to 
signify that one has both $A\ll_{\alpha_1 ,\ldots,\alpha_n}B$ 
and $B\ll_{\alpha_1 ,\ldots,\alpha_n}A$. 
There are a few places where, instead of attaching subscripts (to the $O$, $\ll$, $\gg$ or 
$\asymp$ sign), we have preferred to  
explicitly state the parameters upon which 
the relevant implicit constant may depend. 
\hfill\break $\hbox{\qquad}$
By $f(x)=o(\phi(x))$, we mean that the function 
$\phi$ is positive valued, and that $f(x)/\phi(x)\rightarrow 0$ as $x$ tends to 
a certain limit: in cases where that limit is not specified it should be understood to be the 
limit as $x\rightarrow\infty$, with $x\in{\Bbb R}$.  
By $f(x)\sim g(x)$, we mean that $f(x)-g(x)=o(|g(x)|)$.

\proclaim Measurable sets and functions. 
The term `measurable', 
when applied to a either a subset of $G$ or a function  
$f : G\rightarrow{\Bbb C}$, is to be understood as meaning 
that the set, or function, is measurable with respect to 
the Haar measure ${\rm d}g$ on $G$. 

\bigskip

\goodbreak\centerline{\bf \S 2. Lemmas.}

\medskip 

In this section are collected the lemmas used in the proof of Proposition~2. 
Terminology already defined in Section 1 is used freely (without specific references to those 
definitions). 

\bigskip

\proclaim Lemma 2.1. Suppose that 
${\frak a},{\frak b},{\frak a}',{\frak b}'\in{\Bbb Q}(i)\cup\{\infty\}$, and that 
$\tau_1,\tau_2\in\Gamma=\Gamma_0(q_0)\leq SL(2,{\frak O})$ 
satisfy $\tau_1{\frak a}={\frak a}'$, 
$\tau_2{\frak b}={\frak b}'$ (so that ${\frak a}\sim^{\!\!\!\!\Gamma}{\frak a}'$ 
and ${\frak b}'\sim^{\!\!\!\!\Gamma}{\frak b}$\/). 
Let $g_{\frak a},g_{\frak b},g_{{\frak a}'},g_{{\frak b}'}\in G=SL(2,{\Bbb C})$ 
be chosen so that (1.1.16) and (1.1.20)-(1.1.21) hold for each cusp  
${\frak c}\in\{ {\frak a},{\frak b},{\frak a}',{\frak b}'\}$. 
Put $\rho_1=g_{{\frak a}'}^{-1}\tau_1 g_{\frak a}$ and 
$\rho_2=g_{{\frak b}'}^{-1}\tau_2 g_{\frak b}$. 
Then, for some $\beta_1,\beta_2,\eta_1,\eta_2\in{\Bbb C}$ 
with $\eta_j^2\in{\frak O}^{*}$ ($j=1,2$\/), one has 
$$\rho_j=\pmatrix{\eta_j&0\cr 0&1/\eta_j}
\pmatrix{1&\beta_j\cr 0&1}\qquad 
\hbox{($j=1,2$),}\eqno(2.1)$$ 
$${}^{\frak a}{\cal C}^{\frak b}
=\eta_1\eta_2 {}^{{\frak a}'}\!{\cal C}^{{\frak b}'}\eqno(2.2)$$ 
and, for $c\in {}^{\frak a}{\cal C}^{\frak b}$ and $\omega_1,\omega_2\in{\frak O}$, 
$$S_{{\frak a},{\frak b}}\left(\omega_1,\omega_2;c\right) 
={\rm e}\left({\rm Re}\left(\beta_2\omega_2-\beta_1\omega_1\right)\right) 
 S_{{\frak a}',{\frak b}'}\!\left( \eta_1^{-2}\omega_1\,, 
\,\eta_2^{-2}\omega_2\,;\,\eta_1^{-1}\eta_2^{-1}c\right) .\eqno(2.3)$$

\medskip

\noindent{\bf Proof.}\quad 
The proof is a straightforward adaptation of the proof given on 
Page~239 of [9] for Equation~(1.4) of [9], and is therefore 
omitted here $\blacksquare$ 

\bigskip

\proclaim Lemma 2.2. 
If ${\frak a}$ is a cusp of $\Gamma=\Gamma_0(q_0)\leq SL(2,{\frak O})$ then 
one has 
$${\frak a}\sim^{\!\!\!\!\Gamma}u/w\ {\rm for\ some}\  
u,w\in{\frak O}\ {\rm satisfying}\  (u,w)\sim 1,\ u\neq 0\ {\rm and}\  w\mid q_0\;.\eqno(2.4)$$
Let $u_1,w_1,u_2,w_2\in{\frak O}$ be such that, for $i=1,2$, 
$\left( u_i,w_i\right)\sim 1$ and  $w_i\mid q_0$. 
Then $u_1/w_1\sim^{\!\!\!\!\Gamma}u_2/w_2$ if and only if 
$w_2\sim w_1$ and 
$u_2w_1/w_2\equiv\pm u_1\bmod{\left( w_1,q_0/w_1\right){\frak O}}$. If 
$\,{\frak C}(\Gamma)$ is  
the set of all $\Gamma$-equivalence classes of cusps then  
$$|{\frak C}(\Gamma)|={1\over 8}\sum_{w\mid q_0}\phi\left( \left( w , q_0/w\right)\right) 
+{1\over 8}\!\!\!\!\!\sum_{\scriptstyle w\mid q_0\atop\scriptstyle \left( w , q_0 /w\right)\mid 2}
\!\!\!\!\!\phi\left( \left( w , q_0 /w\right)\right) ,\eqno(2.5)$$
where $\phi$ is Euler's function (i.e. $\phi(\alpha)=\left| ({\frak O}/\alpha{\frak O})^{*}\right|$ 
if $0\neq\alpha\in{\frak O}$\/). 

\medskip

\noindent{\bf Proof.}\quad 
Let ${\frak a}$ be a cusp of $\Gamma=\Gamma_0(q_0)$. Then 
${\frak a}\sim^{\!\!\!\!\Gamma}t/v$ for some $t,v\in{\frak O}$ with $(t,v)\sim 1$. 
Choose $w\sim (v,q_0)$. Then since $(v/w , q_0 /w)\sim 1$ and $(t,w)\sim 1$ it is possible 
to find a pair $\kappa,\delta\in{\frak O}$ with $(\delta , q_0)\sim 1$ that satisfy 
the equation $(q_0 /w)t\kappa +(v/w)\delta =1$. One then has $(\delta , q_0 \kappa)\sim 1$, 
so that, for some $\alpha,\beta\in{\frak O}$, 
$$\Gamma\ni\pmatrix{\alpha &\beta\cr q_0 \kappa &\delta}=\gamma\quad\hbox{(say).}$$
The result (2.4) follows, for 
$\gamma(t/v)=(\alpha t+\beta v)/(q_0 \kappa t+\delta v)=(\alpha t+\beta v)/w$ 
where, since $\gamma\in PSL(2,{\frak O})$, one has $(\alpha t+\beta v,w)\sim 1\,$: 
in the event that $\alpha t+\beta v=0$ one may additionally use 
$0/1\sim^{\!\!\!\!\Gamma}1/1$ (as $\Gamma_0(q_0)\ni n[1]$\/). 
 
Consider now the cusps $u_i/w_i$ ($i=1,2$). Supposing they are $\Gamma$-equivalent,  
there are $\alpha,\beta,\kappa,\delta\in{\frak O}$ such that  
$$\alpha\delta - q_0 \kappa\beta =1\qquad
\hbox{and}\qquad 
{u_2\over w_2}={\alpha u_1+\beta w_1\over q_0 \kappa u_1 + \delta w_1}\,.\eqno(2.6)$$
Since $\left( u_i,w_i\right)\sim 1$ ($i=1,2$), the equations in (2.6) imply 
$w_2\sim q_0 \kappa u_1 + \delta w_1$. Therefore one will have 
$\left( w_2 , q_0\right)\sim\left(\delta w_1, q_0\right)\sim\left( w_1, q_0\right)$, 
which, as $w_i\!\mid\!q_0$ for $i=1,2$ (by hypothesis), implies 
the desired conclusion that $w_2\sim w_1$ when $u_1/w_1\sim^{\!\!\!\!\Gamma} u_2/w_2$.
From the relations $q_0 \kappa u_1+\delta w_1\sim w_2\sim w_1$
(where $w_1\!\mid\!q_0$) and (2.6), one may deduce also that  
$\delta\equiv\varepsilon\bmod\left( q_0 /w_1\right){\frak O}$ 
for the same $\epsilon\in{\frak O}^{*}$ such that $\alpha u_1+\beta w_1=\varepsilon u_2 w_1/w_2$. 
Since the first equation of (2.6) then implies 
$\alpha\equiv\overline{\varepsilon}\bmod\left( q_0 /w_1\right){\frak O}$, one obtains:  
$$u_2 w_1/w_2=\overline{\epsilon}\alpha u_1+\overline{\epsilon}\beta w_1\equiv 
(\overline{\epsilon})^2 u_1\bmod \left( w_1, q_0 /w_1\right){\frak O}.$$
As $(\overline{\epsilon})^2=\pm 1$, it has been shown that 
if $u_2/w_2$ and $u_1/w_1$ are $\Gamma$-equivalent cusps then 
$w_2\sim w_1$ and 
$u_2w_1/w_2\equiv\pm u_1\bmod{\left( w_1,q_0 /w_1\right){\frak O}}$.

In order to establish the converse of what was just found it suffices 
to consider the case of Gaussian integers $u_1$, $u_2$ and $w$ such that 
$$\left( u_1 u_2,w\right)\sim 1,\quad w\!\mid\!q_0 \quad\hbox{and}\quad 
u_2\equiv\pm u_1\bmod(w , q_0 /w){\frak O}\;.\eqno(2.7)$$
If in addition $u_2\equiv u_1\bmod w{\frak O}$ then it is a straightforward deduction that 
$u_2/w\sim^{\!\!\!\!\Gamma}u_1/w$, for 
one has $u_2/w =\gamma\left( u_1/w\right)$  
where $\gamma = n\left[\left( u_2-u_1\right)/w\right]\in\Gamma_0(q_0)$. 

Now let the assumptions be restricted to (2.7) and the congruence 
$u_2\equiv u_1\bmod(q_0 /w){\frak O}$. For $i=1,2$ one can find an element 
$$\sigma_i=\pmatrix{u_i&*\cr w&\tilde{u}_i}\in SL(2,{\frak O})\;.\eqno(2.8)$$
Here $\sigma_i\infty=u_i/w$ ($i=1,2$), so that $u_2/w=\left(\sigma_2\sigma_1^{-1}\right)\left(u_1/w\right)$. 
Therefore, and since the lower left entry of $\sigma_2\sigma^{-1}$ 
is $w\left( \tilde{u}_1-\tilde{u}_2\right)$, one sees  
that if $\tilde{u}_2\equiv\tilde{u}_1\bmod(q_0 /w){\frak O}$ then 
$u_2/w\sim^{\!\!\!\!\Gamma}u_1/w$. Given 
the hypothesis that $u_2\equiv u_1\bmod(q_0 /w){\frak O}$ 
one must have $\tilde{u}_2\equiv\tilde{u}_1\bmod(w , q_0 /w){\frak O}$, for 
(2.8) implies $u_2\tilde{u}_2\equiv u_1\tilde{u}_1\equiv 1\bmod w{\frak O}$.
Consequently $\tilde{u}_2=\tilde{u}_1 + \lambda w + \kappa q_0 /w$ for some 
$\lambda,\kappa\in{\frak O}$. On choosing 
$u_3=u_1$ and $\tilde{u}_3=\tilde{u}_1+\lambda w$, 
the case $i=3$ of (2.8) defines a $\,\sigma_3\in SL(2,{\frak O})$
such that 
$u_2/w=\left(\sigma_2\sigma_3^{-1}\right)\left(u_1/w\right)$ 
and $\sigma_2\sigma_3^{-1}\in\Gamma_0(q_0)$. This shows that (2.7) and 
the congruence $u_2\equiv u_1\bmod(q_0 /w){\frak O}$ are sufficient 
to imply $u_2/w\sim^{\!\!\!\!\Gamma}u_1/w$. Given the conclusion of 
the previous paragraph, it is now proven that sufficient conditions 
for $u_2/w\sim^{\!\!\!\!\Gamma}u_1/w$ to hold are that 
one has (2.7) for the positive choice of sign: alternatively, it would  suffice to 
have (2.7) for the negative choice of sign, since    
$-u_1/w=h[i]\left( u_1/w\right)$ and $h[i]\in\Gamma_0(q_0)$.

By the last two paragraphs, sufficient conditions for the $\Gamma$-equivalence 
of $u_1/w_1$ and $u_2/w_2$ are that 
${w_2\sim w_1}$ and 
$u_2w_1/w_2\equiv\pm u_1\bmod{\left( w_1 , q_0 /w_1\right){\frak O}}\,$;
the necessity of these conditions had already been established, so (as the lemma claims) one does have 
$u_1/w_1\sim^{\!\!\!\!\Gamma}u_2/w_2$ if and only if these conditions hold.

The result (2.5) is a direct consequence of the results preceding it in the lemma,
with the presence of more than one sum over $w$ being explained by the fact that, when 
$( u,w)\sim 1$, one will have $-u\equiv u\bmod(w , q_0 /w){\frak O}$ if and only if 
$(w , q_0 /w)\mid 2\quad\blacksquare$ 

\bigskip 

Proposition~2 concerns the special 
case ${\frak a}'={\frak a}$ of the generalised Kloosterman sum 
$S_{{\frak a},{\frak a}'}\left(\omega,\omega';c\right)$ defined in (1.5.8)-(1.5.10).  
Kloosterman sums of this type are dealt with in the following three lemmas, which  are 
needed for the proof of Proposition~2 (in the next section).

\bigskip 

\proclaim Lemma 2.3. 
Suppose that $u,w\in{\frak O}$ and the cusp ${{\frak a}'}$ 
of $\Gamma =\Gamma_0(q_0)\leq SL(2,{\frak O})$ are such that one has 
$(u,w)\sim 1$, $u\neq 0$, $w\!\mid\!q_0$ and  
$u/w={{\frak a}'}$. Let $v\in{\frak O}$ satisfy 
$$v\sim{q_0 \over\left( q_0 , w^2\right)}\eqno(2.9)$$
(so that $|v|^2$ is the `width' of the cusp ${{\frak a}'}$, as defined below (1.1.22)) 
and put 
$$g_{{\frak a}'}=\pmatrix{u\sqrt{v}&0\cr w\sqrt{v}&(u\sqrt{v})^{-1}} .\eqno(2.10)$$
Then $g_{{\frak a}'}\in G=SL(2,{\Bbb C})$, and $g_{{\frak a}'}$ is such that 
(1.1.16) and (1.1.20)-(1.1.21) hold for ${\frak c}={{\frak a}'}$. 
Suppose moreover that $\,{}^{{\frak a}'}\!{\cal C}^{{\frak a}'}\!$ and 
$S_{{{\frak a}'},{{\frak a}'}}(\omega,\omega';c)$ are as given by (1.5.8)-(1.5.10) 
(i.e. with ${\frak a}={\frak a}'$ and $g_{\frak a}=g_{{\frak a}'}$ there). Then 
$${}^{{\frak a}'}\!{\cal C}^{{\frak a}'}
=\left\{ \gamma vw : 0\neq\gamma\in{\frak O}\ {\rm and}\ 
u\delta^2+\gamma\delta-u\equiv 0\,\bmod(w , q_0 /w){\frak O}\ \,{\rm for\ some}\ 
\delta\in{\frak O}\right\}\eqno(2.11)$$
and, for $c'=\gamma vw\in\,{}^{{\frak a}'}\!{\cal C}^{{\frak a}'}$ and $\omega_1,\omega_2\in{\frak O}$, 
$$S_{{{\frak a}'},{{\frak a}'}}\!\left(\omega_1 ,\omega_2;c'\right)
={\rm e}\!\left( {\rm Re}\!\left( {\omega_2 -\omega_1\over uvw}\right)\right) 
{\sum_{\alpha , \delta\bmod c'{\frak O}}}^{\!\!\!\!\!\!*}\ \  
{\rm e}\!\left( {\rm Re}\!\left( {\omega_1\alpha +\omega_2\delta\over c'}\right)\right) ,\eqno(2.12)$$
where the asterisk signifies that 
the variable of summation, $\delta\bmod c'{\frak O}$, runs over solutions of 
$$\delta (u\delta +\gamma)\equiv u\bmod(w , q_0 /w){\frak O}\quad\hbox{and}\quad 
(\delta,\gamma q_0 /w)\sim 1\sim (u\delta +\gamma,w)\;,\eqno(2.13)$$
while $\alpha\bmod c'{\frak O}$ is determined by:
$$\alpha\delta\equiv 1\bmod(\gamma q_0 /w){\frak O}\quad\hbox{and}\quad 
\left( u_1\alpha -\gamma_1\right)\left( u_1\delta +\gamma_1\right)
\equiv u_1^2\bmod\gamma_1 w{\frak O}\;,\eqno(2.14)$$
where $u_1,\gamma_1\in{\frak O}$ are such that 
$$u_1/\gamma_1=u/\gamma\quad\hbox{and}\quad \left( u_1 ,\gamma_1\right)\sim1\;.\eqno(2.15)$$
The congruences (2.14) imply also that $\delta\bmod c'{\frak O}$ is 
determined by $\alpha\bmod c'{\frak O}$, so that the implied function 
mapping $\delta\bmod c'{\frak O}$ to $\alpha\bmod c'{\frak O}$ 
is a bijection from one subset of ${\frak O}/ c'{\frak O}$ onto another. 

\medskip

\noindent{\bf Proof.}\quad 
This lemma is analogous (even in form) to Lemma~2.5 of [9]: 
the only novelty being that it involves Gaussian integer variables and has 
${\rm e}({\rm Re}(z))$ in place of ${\rm e}(x)$. The 
proof is such a straightforward 
adaptation of the proofs of [9], Lemma~2.4 and Lemma~2.5,  
as to make its inclusion here superfluous\quad $\blacksquare$

\bigskip

The next lemma concerns the more classical type of Kloosterman sum
$S\left(\omega_1,\omega_2 ;c\right)$, given by 
$$S\left(\omega_1,\omega_2 ;c\right)
=\sum_{\scriptstyle\delta\bmod c{\frak O}\atop\scriptstyle (\delta ,c)\sim 1}
{\rm e}\left( {\rm Re}\left( {\omega_1 \delta^{*}+\omega_2\delta\over c}\right)\right)
\qquad\qquad\hbox{($\omega_1,\omega_2\in{\frak O}$ and 
$0\neq c\in{\frak O}$),}\eqno(2.16)$$
where the dependent variable $\delta^{*}\bmod c{\frak O}$ is the solution 
of the congruence $\delta\delta^{*}\equiv 1\bmod c{\frak O}\,$. 
In it are stated some (almost optimal) 
upper bounds for $\left| S\left(\omega_1,\omega_2 ;c\right)\right|$.
These bounds and (2.12) enable one to deduce
the bounds on the generalised Kloosterman sum
$S_{{{\frak a}'},{{\frak a}'}}\left(\omega_1,\omega_2;c'\right)$
that are contained in Lemma 2.5. Lemma 2.5 is directly analogous 
to Lemma~2.6 of [9]: another precedent for this 
type of result may be found in work of Gundlach 
in Section~4 of [14], which includes what is essentially a 
`Weil-Estermann bound' for the analogue of the sum 
$S_{{\frak a},{\frak b}}(\omega_1,\omega_2;c)$ in the 
theory of principal congruence subgroups of 
Hilbert's modular group for any totally real algebraic number field.

\bigskip

\proclaim Lemma 2.4 (a Weil-Estermann bound over ${\Bbb Q}(i)$\/). 
Let $\omega_1,\omega_2\in{\frak O}$, $m\in{\Bbb N}$, and suppose that 
${\varpi}\neq 0$ is a prime element of ${\frak O}$. Then 
$$\left| S\left(\omega_1,\omega_2;{\varpi}^m\right)\right|\leq 
\tau_{\varpi} |{\varpi}|^{\upsilon_{\varpi}} \left|\left(\omega_1,\omega_2,{\varpi}^m\right){\varpi}^m\right|\,,\eqno(2.17)$$
where $\left( \tau_{\varpi},\upsilon_{\varpi}\right)\in{\Bbb R}^2$ is given by 
$$\left( \tau_{\varpi},\upsilon_{\varpi}\right) =\cases{ 
\left( 8\sqrt{2},2\right) &if ${\varpi}\!\mid\!2\,$; \cr 
(2,0) &otherwise.}$$ 
Moreover, for $0\neq c\in{\frak O}$, one has 
$$\left| S\left(\omega_1,\omega_2;c\right)\right|\leq 
2^{7/2} 2^{\omega(c)}\left|\left(\omega_1,\omega_2,c\right) c\right| \;,\eqno(2.18)$$
where $\omega(c)$ is the number of prime ideals of ${\frak O}$ that contain $c$. 

\medskip

\noindent{\bf Proof.}\quad 
The bounds (2.17) and (2.18) are special cases of results obtained by Bruggeman and Miatello,  
in [4], Proposition~9 and Theorem~10\quad $\blacksquare$ 

\bigskip

\proclaim Lemma 2.5. 
Let $q_0$, $\Gamma$, ${{\frak a}'}$, $u$, $v$, $w$, $g_{{\frak a}'}$, 
$\,{}^{{\frak a}'}\!{\cal C}^{{\frak a}'}$ and the generalised Kloosterman sum 
$S_{{{\frak a}'},{{\frak a}'}}(\omega,\omega';c')$
be as described in Lemma 2.3. 
Then  $0\not\in{}^{{\frak a}'}{\cal C}^{{\frak a}'}\subset vw{\frak O}\subseteq{\frak O}$ and 
$$\left| S_{{{\frak a}'},{{\frak a}'}}\left(\omega_1,\omega_2;c'\right)\right| 
\leq \sqrt{8}\,\left|\left(\omega_1,\omega_2,c'\right) c'\right|\tau(c')\qquad 
\hbox{for $\omega_1,\omega_2\in{\frak O}$ and $c'\in\,{}^{{\frak a}'}\!{\cal C}^{{\frak a}'}\,$,}
\eqno(2.19)$$ 
where $\tau(c')$ is the number of Gaussian integer divisors of $c'$. 

\medskip

\goodbreak\noindent{\bf Proof.}\quad 
Suppose that  $c'\in\,{}^{{\frak a}'}\!{\cal C}^{{\frak a}'}$. Then by the hypotheses, and (2.9) and  
(2.11) of Lemma 2.3, one has $u,v,w\in{\frak O}-\{ 0\}$ and 
$c'=\gamma vw$ for some non-zero $\gamma\in{\frak O}$. It therefore only remains to 
prove (2.19). 

Choose $u_1$ and $\gamma_1$ satisfying the conditions (2.15) of Lemma 2.3, and let 
$\omega_1,\omega_2\in{\frak O}$. Then, for $d\in{\frak O}$ with $d\!\mid\!c'$, let 
$K\left(\omega_1,\omega_2;d\right)$ be given by: 
$$K\left(\omega_1,\omega_2;d\right) 
=\qquad\qquad\quad\ \sum\!\!\!\!\!\!\!\!\!\!\!\!\!\!\!\!\!\!\!\!\!\!\!\!\!\!\!\!\sum_{
\!\!\!\!\!\!\!\!\!\!\!\!\!\!{\scriptstyle \alpha,\delta\bmod d{\frak O}\atop
  {\scriptstyle\alpha\delta\equiv 1\bmod(\gamma q_0 /w , d){\frak O}\atop\scriptstyle 
\left( u_1\alpha -\gamma_1\right)
\left( u_1\delta +\gamma_1\right)\equiv u_1^2\bmod\left(\gamma_1 w,d\right){\frak O}
}}}\!\!\!\!\!\!\!\!\!\!\!\!{\rm e}\!\left({\rm Re}\left( {\omega_1\alpha +\omega_2\delta\over d}\right)\right) . 
\eqno(2.20)$$
It may be deduced from (2.12)-(2.15) of Lemma 2.3 that one has 
$$\left| S_{{{\frak a}'},{{\frak a}'}}\left(\omega_1,\omega_2;c'\right)\right| 
=\left| K\left(\omega_1,\omega_2;c'\right)\right|\;.\eqno(2.21)$$
A trivial consequence of (2.20) and (2.21) is that 
$\left| S_{{{\frak a}'},{{\frak a}'}}\left(\omega_1,\omega_2;c'\right)\right|\leq |c'|^2\,$; 
this at least shows that the bound of (2.19)  is true for $c'\in{\frak O}^{*}\,$;
one is therefore to assume henceforth that $c'\in{\frak O}$ is not a unit, or zero. 
Writing 
$c'={\varpi}_1^{e_1}\cdots{\varpi}_r^{e_r}$, where $r,e_1,\ldots ,e_r\in{\Bbb N}$
and ${\varpi}_1,\ldots ,{\varpi}_r$ are non-zero prime elements of ${\frak O}$ 
with ${\varpi}_i\not\sim{\varpi}_j$ for $i\neq j$, it follows by the Chinese Remainder
Theorem that 
$$K\left(\omega_1,\omega_2;c'\right) 
=\prod_{j=1}^r 
K\left(\omega_1\lambda_j,\omega_2\lambda_j;{\varpi}_j^{e_j}\right) ,\eqno(2.22)$$
where, for  $j=1,\ldots ,r$, one has 
$\lambda_j\in{\frak O}$ and 
${\varpi}_j^{-e_j} c'\,\lambda_j\equiv 1\bmod{\varpi}_j^{e_j}{\frak O}$ 
(so that $\left(\lambda_j ,{\varpi}_j\right)\sim 1$\/). 

Let $j\in\{ 1,\ldots ,r\}$ and put $\Omega_k=\omega_k\lambda_j$ ($k=1,2$), 
${\varpi}={\varpi}_j$ and $E=e_j$. Then 
$$\eqalign{K\left(\omega_1\lambda_j,\omega_2\lambda_j;{\varpi}_j^{e_j}\right)
=K\left(\Omega_1,\Omega_2;{\varpi}^{E}\right) 
 &=\qquad\qquad\quad\sum\!\!\!\!\!\!\!\!\!\!\!\!\!\!\!\!\!\!\!\!\!\!\!\!\sum_{
\!\!\!\!\!\!\!\!\!\!\!\!\!\!{\scriptstyle \alpha,\delta\bmod{\varpi}^E{\frak O}\atop
  {\scriptstyle\alpha\delta\equiv 1\bmod{\varpi}^F{\frak O}\atop\scriptstyle 
\left( u_1\alpha -\gamma_1\right)
\left( u_1\delta +\gamma_1\right)\equiv u_1^2\bmod{\varpi}^G{\frak O}
}}}\!\!\!\!\!\!\!\!\!\!\!\!\!\!{\rm e}\!\left({\rm Re}\left( 
{\Omega_1\alpha +\Omega_2\delta\over{\varpi}^E}\right)\right) =\cr 
 &=S_{\varpi}\!\left( E,F,G;\Omega_1,\Omega_2\right)\matrix{\ \cr\ }\quad\hbox{(say),}
}\eqno(2.23)$$
where $F$ and $G$ are the non-negative integers for which one has 
$${\varpi}^F\|\gamma q_0 /w\qquad\hbox{and}\qquad{\varpi}^G\|\gamma_1 w\;.\eqno(2.24)$$
By (2.9), (2.15) and the hypothesis that $(u,w)\sim 1$ one has here 
$[\gamma q_0 /w , \gamma_1 w]\sim\gamma q_0 /(w , q_0 /w)\sim\gamma vw=c'$. 
Therefore (and since ${\varpi}^E={\varpi}_j^{e_j}\|c'$) it must be the case that  
$$\max\{ F , G\} =E\in{\Bbb N}\eqno(2.25)$$
in (2.23), (2.24). It now suffices to show that 
$$\left| S_{\varpi}\!\left( E,F,G;\Omega_1,\Omega_2\right)\right|  
\leq \tau_{\varpi}|{\varpi}|^{\upsilon_{\varpi}}
\left|\left(\Omega_1,\Omega_2,{\varpi}^E\right){\varpi}^E\right|\;,\eqno(2.26)$$
with $\tau_{\varpi}\geq 2$ and $\upsilon_{\varpi}\geq 0$ as defined in Lemma 2.4: 
for, since one has 
$\prod_{k=1}^r \tau_{{\varpi}_k}\left| {\varpi}_k\right|^{\upsilon_{{\varpi}_k}}\leq 
2^{r+7/2}\leq 2^{3/2}\tau(c')$ and $\left(\lambda_k,{\varpi}_k\right)\sim 1$ for $k=1,\ldots ,r$,  
the equations (2.22), (2.23) and the bound (2.26) directly imply (2.19).

If $G>0$, then by (2.24) one has ${\varpi}\!\not\,\mid\!u_1$ 
(since $\left( u_1,\gamma_1\right)\sim 1$, $u_1\!\mid\!u$ and $(u,w)\sim 1$\/) 
and, on choosing $u_1^{*}\in{\frak O}$ so that $u_1 u_1^{*}\equiv 1\bmod{\varpi}^E{\frak O}$, 
one may use the linear change of variables of summation given by  
$\alpha - u_1^{*}\gamma_1=\delta_1$, $\delta+ u_1^{*}\gamma_1=\alpha_1$ 
to deduce that 
$$S_{\varpi}\!\left( E,F,G;\Omega_1,\Omega_2\right)
={\rm e}\left( {\rm Re}\left(\left(\Omega_1 -\Omega_2\right)\gamma_1 u_1^{*}/{\varpi}^E\right)\right) 
S_{\varpi}\!\left( E,G,F;\Omega_2,\Omega_1\right) .\eqno(2.27)$$
By (2.25) and the application of  (2.27) when $G>F$, one may reduce   
to considering the cases in which
$$0\leq G\leq F=E\in{\Bbb N}\eqno(2.28)$$
(it being understood that the values of $F$ and $G$ may 
have been interchanged here, so that (2.24) might no longer be valid). Proving 
(2.26) in only these cases will be sufficient, since the condition (2.25) and 
the upper bound given by (2.26)  remain the same if one interchanges $F$ with $G$, 
or $\Omega_1$ with $\Omega_2$. 

Given (2.28), the congruence conditions on $\alpha$ and $\delta$ in the sum in (2.23) 
simplify to: 
$$\alpha\delta\equiv 1\bmod{\varpi}^E{\frak O}\quad\hbox{and}\quad 
u_1\delta^2 +\gamma_1\delta\equiv u_1\bmod{\varpi}^{G_1}{\frak O}\;,\eqno(2.29)$$
where ${\varpi}^{G_1}\sim{\varpi}^G/\left(\gamma_1,{\varpi}^G\right)$, so that 
$0\leq G_1\leq G$. When $G_1=0$ one therefore has 
$S_{\varpi}\left( E,F,G;\Omega_1,\Omega_2\right)
=S_{\varpi}\left( E,E,G;\Omega_1,\Omega_2\right)
=S\left(\Omega_1,\Omega_2;{\varpi}^E\right)$ (the latter term being one 
of the Kloosterman sums defined by (2.16)), so that by applying the result (2.17) in Lemma 2.4
one obtains exactly the desired bound (2.26). 

Suppose now that $G_1>0$. Then, since $E>0$ and $\left( u_1,\gamma_1\right)\sim 1$, 
either it is the case that the congruences in (2.29) have no simultaneous solutions, 
or else ${\varpi}\!\not\,\mid\!u_1$ 
and the second of those congruences is equivalent to a condition of the 
form $\delta\equiv\mu\pm\sigma\bmod{\varpi}^g{\frak O}$, where 
$\mu,\sigma\in{\frak O}$, $g\in{\Bbb Z}$ and $0<g\leq G_1$.
Therefore either $S_{\varpi}\left( E,F,G;\Omega_1,\Omega_2\right) =0$ 
(in which case (2.26) certainly holds), or 
for some $\mu,\sigma\in{\frak O}$, and some $g\in{\Bbb Z}$ satisfying
$1\leq g\leq G_1\leq G\leq F=E$,
one has: 
$$S_{\varpi}\left( E,F,G;\Omega_1,\Omega_2\right) 
=\sum_{\scriptstyle\delta\bmod{\varpi}^E{\frak O}\atop\scriptstyle 
\delta\equiv\mu\pm\sigma\bmod{\varpi}^g{\frak O}}
{\rm e}\left({\rm Re}\left( {\Omega_1\delta^{*}+\Omega_2\delta\over{\varpi}^E}\right)\right) 
=T_{\varpi}\left( E,g;\Omega_1,\Omega_2;\mu,\sigma\right)\quad\hbox{(say),}\eqno(2.30)$$
where $\delta^{*}$, in the sum over $\delta\bmod\varpi^E {\frak O}$, 
signifies an element of ${\frak O}$ satisfying 
$\delta\delta^{*}\equiv 1\bmod{\varpi}^E{\frak O}$ 
(so that one sums only over $\delta$ with $(\delta ,{\varpi})\sim 1$\/). 

Assuming (2.30) (with $1\leq g\leq E$), take now $H$ to be the non-negative integer with 
$${\varpi}^H\sim\left(\Omega_1,\Omega_2,{\varpi}^E\right)\eqno(2.31)$$ 
and put 
$$E_1=E-H\qquad\hbox{and}\qquad\Omega_k'={\varpi}^{-H}\Omega_k\quad 
\hbox{($k=1,2$),}\eqno(2.32)$$
so that $E_1\in{\Bbb Z}$, $0\leq E_1\leq E$, $\Omega_1',\Omega_2'\in{\frak O}$ and 
$$\left(\Omega_1',\Omega_2',{\varpi}^{E_1}\right)\sim 1\;.\eqno(2.33)$$
From (2.30)-(2.32) it follows trivially that 
$$\eqalign{\left| S_{\varpi}\left( E,F,G;\Omega_1,\Omega_2\right)\right|  
 &=\left|T_{\varpi}\left( E,g;\Omega_1,\Omega_2;\mu,\sigma\right)\right| \leq\cr
 &\leq 2\left| {\varpi}^{E-g}\right|^2
=2|{\varpi}|^{E+H+E_1-2g}
<2|{\varpi}|^{E+H}=2\left|\left(\Omega_1,\Omega_2,{\varpi}^E\right){\varpi}^E\right|\quad
\hbox{if $2g>E_1\,$.}
}\eqno(2.34)$$

By (2.34), one has (2.26) when $2g>E_1$, so (recalling that $g\geq 1$) the only cases  
requiring further consideration are those in which  
$$E_1\geq 2g\geq 2\;.\eqno(2.35)$$
Given (2.31), the summand in the sum appearing in (2.30) is a 
function of $\delta\bmod{\varpi}^{E_1}{\frak O}$. Hence and by (2.32) and (2.35) one 
deduces that 
$$T_{\varpi}\left( E,g;\Omega_1,\Omega_2;\mu,\sigma\right)
=\left| {\varpi}^H\right|^{2}
T_{\varpi}\left( E_1,g;\Omega_1',\Omega_2';\mu,\sigma\right) .\eqno(2.36)$$
Since (2.35) implies $\left[ E_1/2\right]\geq g\geq 1$, one may adapt 
the proofs of Sali\'{e}'s formulae for classical Kloosterman sums to prime-power
moduli ([21], Lemma~4.1 and Lemma~4.2) so as to obtain here: 
$$T_{\varpi}\left( E_1,g;\Omega_1',\Omega_2';\mu,\sigma\right)
=\sum_{\scriptstyle\Delta\bmod{\varpi}^{E_1}{\frak O}\atop{\scriptstyle 
\Delta\equiv\mu\pm\sigma\bmod{\varpi}^g{\frak O}\atop\scriptstyle 
\Omega_2'\Delta\equiv\Omega_1'\Delta^{*}\bmod{\varpi}^{\left[ E_1/2\right]}{\frak O}}}
{\rm e}\left({\rm Re}\left( {\Omega_1'\Delta^{*}+\Omega_2'\Delta\over{\varpi}^{E_1}}\right)\right) ,
\eqno(2.37)$$
where $[x]=\max\{ n\in{\Bbb Z} : n\leq x\}$
and $\Delta^{*}\bmod{\varpi}^{E_1}{\frak O}$ is such that  
$\Delta\Delta^{*}\equiv 1\bmod{\varpi}^{E_1}{\frak O}$. 
Since $\left[ E_1/2\right]\geq 1$, one finds by (2.33) and the last 
condition of summation in (2.37) that 
$T_{\varpi}\left( E_1,g;\Omega_1',\Omega_2';\mu,\sigma\right)\neq 0$
only if ${\varpi}\!\not\,\mid\!\Omega_1'\Omega_2'$. 

If $E_1$ is even, or if ${\varpi}\!\mid\!2$, then one requires now no more than  
a trivial corollary of (2.37): 
$$\left| T_{\varpi}\left( E_1,g;\Omega_1',\Omega_2';\mu,\sigma\right)\right| 
\leq\left|\left\{\Delta\bmod{\varpi}^{E_1}{\frak O} : 
\Omega_2'\Delta^2\equiv\Omega_1'\bmod{\varpi}^{\left[ E_1/2\right]}{\frak O}\right\}\right| 
\leq 2\left| \left( {\varpi}^2 , 2\right) {\varpi}^{E_1-\left[ E_1/2\right]}\right|^2\;.
\eqno(2.38)$$

In cases where $\left( E_1{\varpi},2\right)\sim 1$ one can follow the method of 
proof of [21], Lemma~4.2, one step further, to deduce from (2.37) that 
$$T_{\varpi}\left( E_1,g;\Omega_1',\Omega_2';\mu,\sigma\right)  
=\sum_{\epsilon\in{\cal U}}
{\rm e}\!\left( 2{\rm Re}\left( 
{\Omega_2'\epsilon\Delta_0\over{\varpi}^{E_1}}\right)\right)
\left|{\varpi}^{\left[ E_1/2\right]}\right|^2
\sum_{\beta\bmod\varpi\frak{O}}
{\rm e}\!\left( {\rm Re}\left( {\Omega_2'\epsilon\Delta_0\beta^2\over{\varpi}}\right)\right) ,\eqno(2.39)$$
where ${\cal U}$ is some subset of $\{ 1,-1\}$ and   $\Delta_0$ is (if ${\cal U}\neq\emptyset$\/) 
a Gaussian integer satisfying $\Omega_2'\Delta_0^2\equiv\Omega_1'\bmod{\varpi}^{E_1}$. 
The innermost sum in (2.39) is an analogue of the classical Gauss sum, and 
(by a brief elementary 
manipulation and evaluation of the squared absolute value of the sum) is easily  
seen to have the same absolute value as~${\varpi}$. Therefore, from (2.39) one deduces 
the upper bound 
$$\left| T_{\varpi}\left( E_1,g;\Omega_1',\Omega_2';\mu,\sigma\right)\right|  
\leq 2 |{\varpi}|^{2\left[ E_1/2\right] +1}=2|{\varpi}|^{E_1}\quad\  \,
\hbox{if\quad $\left( E_1{\varpi} , 2\right)\sim 1\,$.}\eqno(2.40)$$

On combining (2.36), (2.38) and (2.40) one obtains: 
$$\left| T_{\varpi}\left( E,g;\Omega_1,\Omega_2;\mu,\sigma\right)\right|
\leq |{\varpi}|^{2H}
\left| T_{\varpi}\left( E_1,g;\Omega_1',\Omega_2';\mu,\sigma\right)\right| 
\leq 2\left|\left( {\varpi}^2 , 2\right)\right|^{5/2} |{\varpi}|^{H+E}\eqno(2.41)$$
in all the cases where (2.35) holds (with $H$ and $E_1$ given by (2.31), (2.32), so that 
$H+E_1=E$\/); it follows by (2.41), (2.30) and (2.31), that one does have 
the desired bound (2.26) when (2.35) holds; and this 
(with the like conclusion having been reached in all other relevant cases) 
establishes the truth of (2.26)  subject to the conditions (2.28);  
as explained  between 
(2.25) and (2.29), this proof is thereby completed\quad 
$\blacksquare$ 

\bigskip

The next result is of a kind long well known in analytic number theory. 
For examples of closely related results already available in the literature 
see, for example, the result obtained by Duke in [10], Theorem~1.1, Part~(i), 
and the result obtained by Coleman in [8], Theorem~6.2.  

\bigskip

\proclaim Lemma 2.6. 
Let $\alpha ,\beta\in{\Bbb R}$, with $\alpha\neq 0$; and let   
$f : (0,\infty)\rightarrow{\Bbb R}$ be a differentiable function such that  
$$f'(x)=\alpha x^{\beta}\qquad\hbox{($x>0$).}$$
Suppose also that $0\neq c\in{\frak O}={\Bbb Z}[i]$, 
$a : {\frak O}-\{ 0\}\rightarrow{\Bbb C}$  and $N\geq 1$, 
and let $s\, :\, {\frak O}\times{\Bbb Z}\times{\Bbb R}\rightarrow{\Bbb C}$ be given by: 
$$s(\delta,m,t)=\sum_{\scriptstyle \omega\in{\frak O}\atop\scriptstyle N/2<|\omega|^2\leq N} 
a(\omega)\left({\omega\over |\omega|}\right)^{\!\!m}
{\rm e}\left({\rm Re}\left({\delta\omega\over c}\right) +tf(|\omega|)\right)\qquad 
\hbox{($\delta\in{\frak O}$, $m\in{\Bbb Z}$, $t\in{\Bbb R}$\/).}$$ 
Then, for $M,T>0$, the sum
$$E_c(N;M,T)=\sum_{\delta\bmod c{\frak O}}\  
\sum_{\scriptstyle m\in{\Bbb Z}\atop\scriptstyle |m|\leq M}\ \int_{-T}^T |s(\delta,m,t)|^2{\rm d}t$$
satisfies
$$E_c(N;M,T)\ll_{\beta}\left( |c|(M+1)+N^{1/2}\right)\left( |c|T+|\alpha|^{-1}N^{-\beta /2}\right)
\sum_{\scriptstyle\omega\in{\frak O}\atop\scriptstyle N/2<|\omega|^2\leq N}|a(\omega)|^2 .$$

\medskip

\noindent{\bf Proof.}\quad 
It is well-known that the real functions 
$$\Lambda(x)=\max\{ 0,\, 1-|x|\}
\qquad\hbox{and}\qquad
{\rm sinc}^2(x)=\cases{ (\pi x)^{-2}\sin^2(\pi x) &if $x\neq 0\,$,\cr 1 &if $x=0\,$,}$$ 
are a pair of mutual Fourier
transforms, so that 
$$\int_{-\infty}^{\infty}{\rm sinc}^2\left( {t\over T_1}\right){\rm e}(tx){\rm d}t
=T_1\Lambda\left( T_1x\right)\quad\ \hbox{for $x\in{\Bbb R}$, $T_1>0\,$.}$$
They have also the following useful properties:  
$$\sum_{m\in{\Bbb Z}}{\rm sinc}^2(\Delta m){\rm e}(mx)
={1\over\Delta}\,\Lambda\left( {\| x\|\over\Delta}\right)\quad\ 
\hbox{$1\geq 2\Delta>0$ and $x\in{\Bbb R}$}$$
(where $\| x\|=\min\{ |x-n| : n\in{\Bbb Z}\}$); 
$${\rm sinc}^2(x)\geq 0\quad\hbox{($x\in{\Bbb R}$);}\quad{\rm and}\quad 
{\rm sinc}^2(x)\geq 4\pi^{-2}\quad\hbox{if $|x|\leq 1/2\,$.}$$
In combination with the orthogonality of the additive characters 
$\psi_{\delta}(\xi)={\rm e}({\rm Re}(\delta\xi/c))$, 
and the non-negativity of $|s(\delta,m,t)|^2$,
the above properties of ${\rm sinc}^2(x)$ and $\Lambda(x)$ enable one 
to see that if $M\geq 1$ and $T>0$, and if one takes $\Delta =1/(2M)$ and $T_1=2T$, then 
$$\eqalign{
E_c(N;M,T) &\leq\left({\pi\over 2}\right)^4
\sum_{\delta\bmod c{\frak O}}\  
\sum_{m\in{\Bbb Z}}{\rm sinc}^2(\Delta m) 
\ \int_{-\infty}^{\infty}{\rm sinc}^2\left( {t\over T_1}\right) |s(\delta,m,t)|^2{\rm d}t =\cr 
 &=\left({\pi\over 2}\right)^4
\quad  
\sum\!\!\!\!\!\!\!\!\!\!\!\sum_{\!\!\!\!\!\!\!\!\!\!\!{\scriptstyle 
\omega_1,\omega_2\in{\frak O}\atop\scriptstyle 
N/2<\left|\omega_1\right|^2,\left|\omega_2\right|^2\leq N}}
\!\!\!a\left(\omega_1\right)\,\overline{a\left(\omega_2\right)}\,
\prod_{j=1}^3 F_j\left(\omega_1,\omega_2\right) ,}$$
where 
$$F_1\left(\omega_1,\omega_2\right)
=\sum_{\delta\bmod c{\frak O}}{\rm e}\left( {\rm Re}\left( 
{\delta\left(\omega_1 -\omega_2\right)\over c}\right)\right) 
=\cases{|c|^2 &if $\omega_1\equiv\omega_2\bmod c{\frak O}\,$, \cr 
0 &otherwise,}$$
$$F_2\left(\omega_1,\omega_2\right)
=2M\Lambda\!\left( 2M\left\| 
{ {\rm Arg}\left(\omega_1\right) -{\rm Arg}\left(\omega_2\right)\over 2\pi}\right\|\right)\quad 
\hbox{and}\quad 
F_3\left(\omega_1,\omega_2\right)
=2T\Lambda\bigl( 2T f\left(\left|\omega_1\right|\right) -2Tf\left(\left|\omega_2\right|\right)\bigr) .$$
Therefore, and since $\Lambda(x)=0$ when $|x|\geq 1$, while 
$0\leq\Lambda(x)\leq 1$ for all $x\in{\Bbb R}$,  
one consequently has: 
$$E_c(N;M,T)\leq {\pi^4\over 4}|c|^2MT\sum_{\omega_1,\omega_2}\left| a\left(\omega_1\right) 
a\left(\omega_2\right)\right| ,$$ where the sum is over pairs 
$(\omega_1,\omega_2)\in{\frak O}\times{\frak O}$ satisfying
$$N/2<|\omega_1|^2,|\omega_2|^2\leq N,\qquad \omega_1\equiv\omega_2\pmod{c},$$
$$\left\| {{\rm Arg}(\omega_1)-{\rm Arg}(\omega_2)\over 2\pi}\right\|
\leq {1\over 2M}\quad\hbox{and}\quad \left| f(|\omega_1|)-f(|\omega_2|)\right| 
\leq {1\over 2T}\;.$$
\par
These conditions imply that $\bigl| |\omega_1|-|\omega_2|\bigr |\leq 2^{|\beta /2|-1}(|\alpha|T)^{-1} N^{-\beta
/2}$. Applying the arithmetic-geometric mean inequality to 
$\left| a\left(\omega_1\right)\right|\left| a\left(\omega_2\right)\right|$, and
appealing to symmetry, one finds that
$$E_c(N;M,T)\leq{\pi^4\over 4}|c|^2MT\sum_{N/2<|\omega|^2\leq N}\left| a\left(\omega\right)\right|^2
V_{\omega} ,\eqno(2.42)$$ where $V_{\omega}$ is the number of elements $\omega'\in{\frak O}$
satisfying
$$\omega'\equiv\omega\bmod c{\frak O},\eqno(2.43)$$
$$\left\|{{\rm Arg}(\omega')-{\rm Arg}(\omega)\over 2\pi}\right\|
\leq{1\over 2M}\quad\hbox{and}\quad \bigl| |\omega'|-|\omega|\bigr|\leq {2^{|\beta /2|}\over
2|\alpha| TN^{\beta /2}} .$$ The latter two conditions may be simplified to the form
$$r\leq |\omega'|\leq R\quad\hbox{and}\quad\phi\leq{\rm Arg}(\omega')\leq\Phi ,$$
where $r,R,\phi,\Phi$ are determined by $\omega ,\alpha ,\beta , T, M, N$, and satisfy
$$0\leq r\leq R\leq N^{1/2},\qquad 0\leq\phi\leq\Phi\leq 2\pi ,$$
$$R-r\ll_{\beta}{N^{-\beta /2}\over |\alpha|T}\quad\hbox{and}\quad
\Phi -\phi\leq{2\pi\over M}.$$ By this, and by (2.43), the complex numbers $(\omega'-\omega)/c$ are
Gaussian integers lying in a simply connected region ${\cal R}\subset{\Bbb C}$ with
$${\rm Area}({\cal R})={\left( R^2-r^2\right)(\Phi -\phi)\over 2|c|^{2}}
\ll_{\beta}|\alpha|^{-1}N^{(1-\beta)/2}T^{-1}M^{-1}|c|^{-2}$$ and
$${\rm Perimeter}({\cal R})={2(R-r)+(\Phi -\phi)(R+r)\over |c|} 
\ll_{\beta} |\alpha|^{-1}N^{-\beta /2}T^{-1}|c|^{-1}+N^{1/2}M^{-1}|c|^{-1}\;.$$
Therefore
$$V_{\omega}\ll 1+{\rm Perimeter}({\cal R})+{\rm Area}({\cal R}) 
\ll_{\beta}\left( 1+N^{1/2}M^{-1}|c|^{-1}\right)
\left( 1+|\alpha|^{-1}N^{-\beta /2}T^{-1}|c|^{-1}\right) ,$$ 
so that the result of the lemma now follows by
(2.42). This completes the proof in cases where $M\geq 1$, and 
implies the validity, when $M\geq 0$, of a bound for $E_c(N,M+1,T)$ 
similar to the bound for $E_c(N,M,T)$ appearing in the lemma
(and differing only in that the implicit constant may be larger by some factor $b\leq 2$\/); 
since one has (trivially) $E_c(N;M+1,T)\geq E_c(N;M,T)$, the remaining cases of the lemma follow\quad $\blacksquare$ 

\bigskip

\proclaim Lemma 2.7 (Poisson summation over ${\Bbb Z}^n$ and over ${\frak O}={\Bbb Z}[i]$\/). 
Let $n\in{\Bbb N}$ and suppose that the function
$F : {\Bbb R}^n\rightarrow{\Bbb C}$ lies in 
the Schwartz space: so that, for all $A\geq 0$ and all integers $j,k\geq 0$,
the function 
$\|{\bf x}\|^A {\partial^{j_1+\cdots j_n}\over\partial x_1^{j_1}\cdots\partial x_n^{j_n}} F({\bf x})$
is continuous and bounded on ${\Bbb R}^n$. Then the Fourier transform
$$\hat F(y_1,\ldots ,y_n)
=\int_{{\Bbb R}^n}F({\bf x}){\rm e}(-{\bf x}\cdot{\bf y}) {\rm d}x_1\cdots{\rm d}x_n\eqno(2.44)$$
is a complex-valued functions defined on ${\Bbb R}^n$, and lies in the Schwartz space. 
One has, moreover: 
$$\sum_{{\bf m}\in{\Bbb Z}^n}F({\bf m})
=\sum_{{\bf m}\in{\Bbb Z}^n}\hat F({\bf m})\;.\eqno(2.45)$$

When $n=2$ and $f : {\Bbb C}\rightarrow{\Bbb C}$ 
is given by $f(x+iy)=F(x,y)$ for $(x,y)\in{\Bbb R}^2$ (with $F({\bf x})$ as above), 
the complex Fourier transform 
$$\hat f(w)=\int_{\Bbb C} f(z){\rm e}\left(-{\rm Re}(wz)\right) {\rm d}_{+}z
=\int_{-\infty}^{\infty}\int_{-\infty}^{\infty}f(x+iy)
{\rm e}\left(-{\rm Re}((x+iy)w)\right) {\rm d}x {\rm d}y\eqno(2.46)$$
is a complex-valued function defined on ${\Bbb C}$, and one has 
$$\sum_{\alpha\in{\frak O}}f(\alpha)
=\sum_{\alpha\in{\frak O}}\hat f(\alpha)\;.\eqno(2.47)$$

\medskip

\noindent{\bf Proof.}\quad
For the results up to and including (2.45) see, for example [30], Chapter~13, 
Section~4 and Section~6.
The results (2.46) and (2.47) amount to a special case of (2.44)-(2.45), for the 
right-hand side of (2.46) is (by definition) equal to $\hat F\left( {\rm Re}(w) , -{\rm Im}(w)\right)$ 
and, since complex conjugation is a permutation of ${\frak O}$, it therefore follows from 
(2.45) that  
$\sum_{\alpha\in{\frak O}}f(\alpha)
=\sum_{\alpha\in{\frak O}}\hat f(\overline{\alpha})
=\sum_{\alpha\in{\frak O}}\hat f(\alpha)$
\ $\blacksquare$

\bigskip

\proclaim Lemma 2.8. Let $f$ and $F$ be as in the case $n=2$ of Lemma 2.7. 
For $z=x+iy$ with $x,y\in{\Bbb R}$ put 
$$({\bf\Delta}_{\Bbb C} f)(z)=({\bf\Delta}_{{\Bbb R}\times{\Bbb R}} F)(x,y)
=\left({\partial^2\over\partial x^2}+{\partial^2\over\partial y^2}\right) f(x+iy)\;.$$
Then ${\bf\Delta}_{{\Bbb R}\times{\Bbb R}} F$ (the Laplacian of $F$) is a member of the Schwartz space.
The functions $f$ and ${\bf\Delta}_{\Bbb C} f$ have 
Fourier transforms $\hat{f},\widehat{{\bf\Delta}_{\Bbb C} f} : {\Bbb C}\rightarrow{\Bbb C}$
(defined as in Lemma 2.7), which are related to one another by:
$$|2\pi w|^2 \hat{f}(w) = -\widehat{{\bf\Delta}_{\Bbb C} f}(w)\quad\hbox{for $w\in{\Bbb C}$.}\eqno(2.48)$$
For all $w\in{\Bbb C}-\{ 0\}$ and $j=0,1,2,\ldots\ $ , one  has
$$\bigl|\hat{f}(w)\bigr| 
=\left( 2\pi |w|\right)^{-2j}\Bigl|\widehat{{\bf\Delta}_{\Bbb C}^j f}(w)\Bigr|
\leq\left( 2\pi |w|\right)^{-2j}
\widehat{\bigl| {\bf\Delta}_{\Bbb C}^j f\bigr|}\,(0) 
=\left( 2\pi |w|\right)^{-2j}
\int_{\Bbb C}\,\Bigl| \bigl( {\bf\Delta}_{\Bbb C}^j f\bigr)(z)\Bigr| 
\,{\rm d}_{+}z\;.\eqno(2.49)$$

\medskip

\noindent{\bf Proof.}\quad
Since $F$ is a member of the Schwartz space, so is $\bigl(\partial^2/\partial x^2\bigr) F(x,y)$.
Therefore two applications of integration by parts suffice to show that, when $u,y\in{\Bbb R}$,
$$\int_{-\infty}^{\infty} {\rm e}(-ux) F_{11}(x,y) {\rm d}x
=(2\pi i u)^2\int_{-\infty}^{\infty} F(x,y){\rm e}(-ux){\rm d}x\;,$$
where
$F_{11}(x,y)=({\partial}/{\partial x})^2 F(x,y)$. It follows that
$\hat{F}_{11}(u,v)=-(2\pi u)^2\hat{F}(u,v)$, for $u,v\in{\Bbb R}$.
Similarly, integrations by parts with respect to $y$ yield
$\hat{F}_{22}(u,v)=-(2\pi v)^2\hat{F}(u,v)$, where
$F_{22}(x,y)=({\partial}/{\partial y})^2 F(x,y)$.
Therefore, for $u,v\in{\Bbb R}$ and $w=u+iv$, one has:
$$|2\pi w|^2 \hat{f}(w) = \bigl( (2\pi u)^2+(2\pi v)^2\bigr)\hat{F}(u,-v)
=-\bigl(\hat{F}_{11}(u,-v)+\hat{F}_{22}(u,-v)\bigr)
=-\widehat{{\bf\Delta}_{{\Bbb R}\times{\Bbb R}} F}(u,-v)
=-\widehat{{\bf\Delta}_{\Bbb C} f}(w)\;,$$
which is (2.48). Since 
${\bf\Delta}_{\Bbb C}\,{\rm e}\bigl({\rm Re}(wz)\bigr)=-|2\pi w|^2{\rm e}\bigl({\rm Re}(wz)\bigr)$,  
the result (2.48) is (in essence) illustrative of the fact that the Laplacian operator 
for functions on ${\Bbb R}^n$ is symmetric on a dense subspace of 
the Lebsgue space $L^2\left( {\Bbb R}^n\right)$. 

From (2.48) it follows by induction that
$$\hat{f}(w)=(-1)^j |2\pi w|^{-2j}\,\widehat{{\bf\Delta}_{\Bbb C}^j f}(w)
\quad\hbox{for $w\in{\Bbb C}-\{ 0\}\,$ and $j=0,1,2,\ldots\ $ .}$$
This implies the first equality of (2.49): the rest of (2.49) follows since 
(by virtue of $\bigl({\bf\Delta}_{{\Bbb R}\times{\Bbb R}}^j F\bigr)(x,y)$ being a member of the Schwartz space,
and as $|{\rm e}(x)|=1$ for $x\in{\Bbb R}$) it is the case that 
the function $g(z)=\bigl({\bf\Delta}_{\Bbb C}^j f\bigr)(z) {\rm e}(-{\rm Re}(wz) )$
satisfies 
$$\left|\,\int_{\Bbb C} g(z)\,{\rm d}_{+}z\right|\leq\int_{\Bbb C}\left| g(z)\right|{\rm d}_{+}z
=\int_{\Bbb C}\,\Bigl|\bigr({\bf\Delta}_{\Bbb C}^j f\bigr)(z)\Bigr|\,{\rm d}_{+}z<+\infty\quad\blacksquare$$

\medskip 

\goodbreak\centerline{\bf \S 3. The proof of Proposition~2.}

\medskip

\noindent{\bf Part I: a proof of (1.9.26) using the Weil-Estermann bound.}\quad 
By Lemma 2.2 one may choose coprime non-zero 
Gaussian integers $u,w$ such that $w\mid q_0$ and 
$u/w\sim^{\!\!\!\!\Gamma}{\frak a}$. 
One then has $\tau{\frak a}=u/w$ for some $\tau\in\Gamma$. 
Choose next a Gaussian integer $v$ satisfying hypothesis (2.9) of Lemma 2.3. 
Then, on putting $\rho =g_{u/w}^{-1}\tau g_{\frak a}$, where $g_{u/w}$ is the scaling 
matrix $g_{{\frak a}'}$ given by (2.10) of Lemma 2.3, it follows by Lemma 2.1 that  
$\rho =h\bigl[\eta] n[\beta]$ for some $\beta\in{\Bbb C}$ 
and some $\eta\in{\Bbb C}$ with $\eta^2\in{\frak O}^{*}$,  
and that this $\beta$  and the unit $\eta^2=\epsilon$ (say) are such that one has both 
$${}^{\frak a}{\cal C}^{\frak a}=\epsilon\,{\cal C}'\;,\quad
\hbox{where}\quad {\cal C}'={}^{u/w}{\cal C}^{u/w}\;,\eqno(3.1)$$ 
and, when $c\in {}^{\frak a}{\cal C}^{\frak a}$ and $\omega_1,\omega_2\in{\frak O}$, 
$$S_{{\frak a},{\frak a}}\left(\omega_1,\omega_2;c\right) 
={\rm e}\left( {\rm Re}\left(\beta\left(\omega_2 -\omega_1\right)\right)\right) 
S_{u/w,u/w}\left({\omega_1\over\epsilon}\,,\,{\omega_2\over\epsilon}\,;{c\over\epsilon}\right) \eqno(3.2)$$
(note that the case of Lemma 2.1 used here is that in which ${\frak b}={\frak a}$, 
${\frak b}'={\frak a}'=u/w$, $g_{\frak b}=g_{\frak a}$, 
$g_{{\frak b}'}=g_{{\frak a}'}=g_{u/w}$ and $\tau_1 =\tau_2 =\tau$:  
one therefore has $\rho_1=\rho_2=\rho$, which, by (2.1),  implies both 
$\eta_1=\eta_2=\eta$, say,  
and $\beta_1=\beta_2=\beta$, say\/).
In view of Lemma 2.5 (where ${\frak a}'=u/w$\/), it follows from (3.1), (2.9) and (3.2) that  
$$0\not\in{}^{\frak a}{\cal C}^{\frak a}\subset \epsilon vw{\frak O} 
=vw{\frak O}={1\over\mu({\frak a})}\,{\frak O}\eqno(3.3)$$
(where $\mu({\frak a})$ is as discussed in Remark~1.9.7, below Theorem 1),  
and that 
$$\left| S_{{\frak a},{\frak a}}\left(\omega_1,\omega_2;c\right)\right| 
\leq\sqrt{8}\left|\left( {\omega_1\over\epsilon}\,,\,{\omega_2\over\epsilon}\,,\,{c\over\epsilon}\right)
{c\over\epsilon}\right| \tau\left( {c\over\epsilon}\right)
=\sqrt{8}\left|\left(\omega_1\omega_2,c\right) c\right|\tau(c)\;,\eqno(3.4)$$
for $c\in {}^{\frak a}{\cal C}^{\frak a}$ and $\omega_1,\omega_2\in{\frak O}$. 
Note that (3.3) is the result (1.9.24) of the proposition. 
Applying (3.4) directly to (1.9.25), one obtains (subject to the hypotheses of the proposition) 
the upper bound: 
$$\left| U_{\frak a}(\psi , c; M; N, b)\right|
\leq \sqrt{8}\,(2M+1)\quad 
\sum\!\!\!\!\!\!\!\!\!\!\!\sum_{\!\!\!\!\!\!\!\!\!\!\!{\scriptstyle \omega_1,\omega_2\in{\frak O}\atop\scriptstyle 
N/2<\left|\omega_1\right|^2,\left|\omega_2\right|^2\leq N}}
\!\!\!\!\!\left| b\left(\omega_1\right) b\left(\omega_2\right)
\left(\omega_1\omega_2,c\right) c\right|\tau(c)\;.$$
This bound for $\left| U_{\frak a}(\psi , c; M; N ,b)\right|$ 
implies the result in (1.9.26): for one has, by the Cauchy-Schwarz inequality,  
$$\quad 
\sum\!\!\!\!\!\!\!\!\!\!\!\sum_{\!\!\!\!\!\!\!\!\!\!\!{\scriptstyle \omega_1,\omega_2\in{\frak O}\atop\scriptstyle 
N/2<\left|\omega_1\right|^2,\left|\omega_2\right|^2\leq N}}
\!\!\!\!\!\left|\left(\omega_1,\omega_2,c\right) b\left(\omega_1\right) b\left(\omega_2\right)\right|
\leq\left( E(c,N)\right)^{1/2} \left\| b_N\right\|_2^2\;,$$
where $\left\| b_N\right\|_2^2$ is as defined in (1.9.16), while 
$$\eqalign{E(c,N) &=\quad 
\sum\!\!\!\!\!\!\!\!\!\!\!\sum_{\!\!\!\!\!\!\!\!\!\!\!{\scriptstyle \omega_1,\omega_2\in{\frak O}\atop\scriptstyle 
N/2<\left|\omega_1\right|^2,\left|\omega_2\right|^2\leq N}}
\!\!\!\!\!\left|\left(\omega_1,\omega_2,c\right)\right|^2 \leq\cr 
 &\leq{1\over 4}\sum_{\delta\mid c} |\delta|^2
\Biggl(\sum_{\scriptstyle\omega\in\delta{\frak O}\atop\scriptstyle 0<|\omega|^2\leq N} 1\Biggr)^{\!\!2}
=\sum_{\delta\mid c} |\delta|^2\left( O\left( {N\over |\delta|^2}\right)\right)^2
\ll N^2\sum_{\delta\mid c} {1\over |\delta|^2}\leq N^2\tau(c)\quad\square}$$

\medskip 

\noindent{\bf Part II: a proof of (1.9.27) by direct use of Lemma 2.6.}\quad  
The bounds (1.9.27) and (1.9.28) remain to be proved. 
Bounds on the absolute values of the relevant Kloosterman sums 
are not sufficient  to achieve this: one may use instead the explicit 
representation given by (2.12) of Lemma 2.3. Bearing in mind (3.1)-(3.3),  
the result (2.12) shows that 
$$S_{{\frak a},{\frak a}}\left(\omega_1,\omega_2;c\right) 
={\rm e}\left( {\rm Re}\left( \epsilon^{-1}\beta'\left(\omega_2 -\omega_1\right)\right)\right) 
{\sum_{\alpha , \delta\bmod c'{\frak O}}}^{\!\!\!\!\!\!*}\quad 
{\rm e}\left( {\rm Re}
\left( {\epsilon^{-1}\omega_1\alpha +\epsilon^{-1}\omega_2\delta\over c'}\right)\right) ,$$
where
$$\eqalignno{
\beta'  &=\epsilon\beta +{1\over uvw}\;,\cr 
c' &=c/\epsilon\in{\cal C}'\subseteq {1\over\mu({\frak a})}\,{\frak O}-\{ 0\}\;,&(3.5)
}$$
and where the asterisk indicates the same conditions of summation as in (2.13)-(2.15) of Lemma 2.3.
By using this to rewrite the right-hand side of (1.9.25), then substituting 
$\epsilon\omega_1$ and $\epsilon\omega_2$ for $\omega_1$ and $\omega_2$ (respectively) 
and noting that $|\epsilon|^{2m}=1^{2m}=1$, 
one obtains: 
$$\eqalignno{
U_{\frak a}(\psi,c;M;N,b) &=\sum_{m=-M}^{M}
\Biggl|\quad\ \;  
\sum\!\!\!\!\!\!\!\!\!\!\!\sum_{\!\!\!\!\!\!\!\!\!\!\!{\scriptstyle \omega_1,\omega_2\in{\frak O}\atop\scriptstyle 
N/2<\left|\omega_1\right|^2,\left|\omega_2\right|^2\leq N}}
\!\!\!\!\!\overline{ b\left(\epsilon\omega_1\right)} b\left(\epsilon\omega_2\right) 
\left({\omega_1\omega_2\over\left|\omega_1\omega_2\right|}\right)^{\!\!m}
{\rm e}\!\left(
{\psi\sqrt{\left|\omega_1\omega_2\right|}\over \left| c'\right|}
+{\rm Re}\left( \beta'\left(\omega_2 -\omega_1\right)\right)\right)\times\qquad\cr
 &\qquad\qquad\qquad\qquad\qquad\qquad\qquad\qquad\qquad\qquad\qquad   
\times{\sum_{\alpha , \delta\bmod c'{\frak O}}}^{\!\!\!\!\!\!*}\quad 
{\rm e}\left( {\rm Re}
\left( {\omega_1\alpha +\omega_2\delta\over c'}\right)\right) 
\Biggr| =\cr 
 &\quad\matrix{\ }\hbox{\quad}\cr 
 &=\sum_{m=-M}^{M}
\left|\quad\ \;  
\sum\!\!\!\!\!\!\!\!\!\!\!\sum_{\!\!\!\!\!\!\!\!\!\!\!{\scriptstyle \omega_1,\omega_2\in{\frak O}\atop\scriptstyle 
N/2<\left|\omega_1\right|^2,\left|\omega_2\right|^2\leq N}}
\!\!\!\!\!\overline{ a\!\left(\omega_1\right)} a\!\left(\omega_2\right) 
\left({\omega_1\omega_2\over\left|\omega_1\omega_2\right|}\right)^{\!\!m}
{\rm e}\!\left(
{\psi\sqrt{\left|\omega_1\omega_2\right|}\over \left| c'\right|}\right)
S^{*}\!\left(\omega_1,\omega_2;c'\right)\right|\;,&(3.6)}$$
where
$$ a(\omega)= b(\epsilon\omega)\,{\rm e}\!\left( {\rm Re}\left(\beta'\omega\right)\right)\qquad 
\hbox{($0\neq\omega\in{\frak O}$)}\eqno(3.7)$$
and 
$$S^{*}\!\left(\omega_1,\omega_2;c'\right) 
={\sum_{\alpha , \delta\bmod c'{\frak O}}}^{\!\!\!\!\!\!*}\quad 
{\rm e}\left( {\rm Re}
\left( {\omega_1\alpha +\omega_2\delta\over c'}\right)\right) \qquad
\hbox{($\omega_1,\omega_2\in{\frak O}$)}\eqno(3.8)$$
with the asterisk superfixed to the summation sign having the same meaning as in Lemma 2.3. 
It will therefore suffice to bound the slightly more general sum: 
$$U_{q_0 ,w,u}^{\circ}\left(\psi , c';M,N\right) 
=\sum_{m=-M}^{M}
\eta_m\quad  
\sum\!\!\!\!\!\!\!\!\!\!\!\sum_{\!\!\!\!\!\!\!\!\!\!\!{\scriptstyle \omega_1,\omega_2\in{\frak O}\atop\scriptstyle 
N/2<\left|\omega_1\right|^2,\left|\omega_2\right|^2\leq N}}
\!\!\!\!\!\overline{ a\!\left(\omega_1\right)} a\!\left(\omega_2\right) 
\left({\omega_1\omega_2\over\left|\omega_1\omega_2\right|}\right)^{\!\!m}
\!{\rm e}\!\left(
{\psi\sqrt{\left|\omega_1\omega_2\right|}\over \left| c'\right|}\right)
S^{*}\!\left(\omega_1,\omega_2;c'\right) ,\eqno(3.9)$$
where the coefficients $\eta_m$ are arbitrary complex numbers satisfying  
$$\left|\eta_m\right| =1\quad\hbox{for}\quad m=-M,-M+1,\ \ldots\ ,M\,.\eqno(3.10)$$
\par
The next step requires the introduction of a `redundant weighting function'  
$f : (0,\infty)\rightarrow [0,1]$ possessing a continuous second derivative, and 
such that 
$$f(x)=\cases{0 &if $x\leq 1/2$ or $x\geq 2\,$;\cr 1 &if $1/\sqrt{2}\leq x\leq 1\,$.}$$
Such a function can be explicitly defined (if need be). Here it will suffice to 
suppose that $f$ is chosen absolutely independently of 
all factors other than the above explicit requirements: 
by this supposition there exists an absolute constant $C_2\in [1,\infty)$ such that  
$\left| f^{(j)}(x)\right|\leq C_2$ for $j=0,1,2$ and $x\in{\Bbb R}$.
For this $f$ one has (see Appendix~A.2 of [19]): 
$$f(x){\rm e}\left({\psi\sqrt{N}\,x\over \left| c'\right|}\right)
=\int_{-\infty}^{\infty} g_{Y}(t)x^{it}{\rm d}t\qquad\hbox{($x>0$),}\eqno(3.11)$$ 
where 
$$Y ={\psi\sqrt{N}\over\left| c'\right|}\in{\Bbb R}\eqno(3.12)$$
and
$g_{Y} : {\Bbb R}\rightarrow{\Bbb C}$ is the Mellin-transform given by 
$$g_{Y}(t)={1\over 2\pi}\int_0^{\infty} f(x){\rm e}\left( Y x\right)
x^{-1-it}{\rm d}x\;.$$
The point of these observations is that, since 
$f\bigl(\sqrt{\left|\omega_1\omega_2\right| /N}\,\bigr)=1$ when 
$N/2<\left|\omega_1\right|^2,\left|\omega_2\right|^2\leq N$, 
the case $x=\sqrt{\left|\omega_1\omega_2\right| /N}$ of (3.11) therefore 
provides a means to effect a `separation of variables' 
in the summand on the right-hand side of equation (3.9)
(the relevant variables being $\omega_1$ and $\omega_2$).
\par
Practical application of (3.11) requires suitable estimates for $\left| g_Y(t)\right|$. 
One has, trivially, 
$$g_Y(t)\ll\int_{0}^{\infty} |f(x)| x^{-1}{\rm d}x\leq\int_{1/2}^2 x^{-1}{\rm d}x\ll 1\;.$$
For less trivial estimates, in cases where $t\neq 0$, it is helpful to note that 
$$g_Y(t)=\int_{1/2}^2 F(x)e^{iu(x)}{\rm d}x\eqno(3.13)$$ 
where 
$$F(x)=(2\pi x)^{-1}f(x)\quad\hbox{and}\quad
u(x)=u_{Y,t}(x)=2\pi Y x-{t}\log x .$$
Here $\inf_{1/2\leq x\leq 2}\left| u''(x)\right| =\inf_{1/2\leq x\leq 2}|t|/x^2\geq |t|/4$, while  
the choice of $f$ ensures that $\sup_{1/2\leq x\leq 2}|F(x)|\ll 1$ and 
$\int_{1/2}^2 \left| F'(x)\right|{\rm d}x\ll 1$; it therefore follows by the `second derivative test'   
(as formulated in Lemma~5.1.3 of~[18])  that one has, uniformly with respect to $Y$,  
$$g_Y(t)\ll |t|^{-1/2}\qquad\hbox{for $0\neq t\in{\Bbb R}$.}$$
A third bound on $g_Y(t)$ is useful in cases where 
$$|t|\geq T_0=8\pi |Y|.\eqno(3.14)$$
In such cases one finds that
$${|t|\over 4}\leq \left| {{\rm d}u\over{\rm d}x}\right|=\left| 2\pi Y-{t\over x}\right|\leq 3|t|\qquad\hbox{($1/2<x<2$),}$$
so that one may rewrite (3.13) as 
$$g_Y(t)=\int_{a}^{b} F(x){{\rm d}x\over{\rm d}u}\, e^{iu}{\rm d}u,$$
with $|b-a|<6|t|$ ($a=u(1/2)$, $b=u(2)$\/).  Then two integrations-by-parts show that 
$$g_Y(t)=-\int_a^b\left( F''(x)\left({{\rm d}x\over{\rm d}u}\right)^3
+3F'(x){{\rm d}x\over{\rm d}u}{{\rm d}^2 x\over{\rm d}u^2} +F(x){{\rm d}^3 x\over{\rm d}u^3}\right)
e^{iu}{\rm d}u$$ 
(boundary terms are absent, since $f$ is supported in the interval $[1/2,2]$, and has a 
continuous second derivative on $(-\infty,\infty)$\/). 
Here it follows by elementary real-variable calculus that when $a<u<b$ (so that $1/2<x<2$\/) one has 
$$F^{(j-1)}(x)\ll 1\qquad\hbox{and}\qquad {{\rm d}^j x\over{\rm d}u^j}\ll |t|^{-j}\qquad
\hbox{($j=1,2,3$).}$$
Consequently one deduces from the last equation involving $g_Y(t)$ 
that if (3.14) holds then 
$$g_Y(t)\ll |b-a||t|^{-3}\ll |t|^{-2}\;.$$
On combining this last bound for $g_Y(t)$ with the others found before, one has: 
$$g_Y(t)\ll\cases{ 1 &if $|t|\leq 1$; \cr |t|^{-1/2} &if $1<|t|\leq 1+T_0$; \cr
|t|^{-2} &if $|t|>1+T_0$.}\eqno(3.15)$$
\par
In light of the final remark in the paragraph containing (3.11), one 
may attach to each summand on the right-hand side of Equation~(3.9) 
an extra factor $f(\sqrt{|\omega_1\omega_2|/N})\,$ (this action 
does not change the value of any summand in (3.9), for it is    
equivalent in effect to multiplying each summand by $1$).   
Hence, by applying the case $x=\sqrt{\left|\omega_1\omega_2\right| /N}$ of (3.11)   
and recalling the definition (3.8) of $S^{*}\!\left(\omega_1,\omega_2;c'\right)$, 
one finds that   
$$U_{q_0 ,w,u}^{\circ}\left(\psi,c';M,N\right) 
=\int_{-\infty}^{\infty}g_Y(t)\sigma_{q_0 ,w,u}\left( c';M,N;t\right) N^{-it/2}{\rm d}t ,$$
where
$$\sigma_{q_0 ,w,u}\left( c';M,N;t\right) 
={\sum_{\alpha , \delta\bmod c'{\frak O}}}^{\!\!\!\!\!\!*}\ \ \, 
\sum_{m=-M}^M\eta_m
A\left( {\alpha\over c'}\,,m,t\right)
\overline{A\!\left( -{\delta\over c'}\,,-m,-t\right) }\;,$$
with
$$A(\theta,n,v)
=\sum_{\scriptstyle \omega\in{\frak O}\atop\scriptstyle 
N/2<\left|\omega\right|^2\leq N}\!\!\!\overline{ a\left(\omega\right)}\,
{\rm e}\left( {\rm Re}\left(\omega\theta\right)\right) 
\left({\omega\over |\omega|}\right)^n
|\omega|^{iv/2}\;.\eqno(3.16)$$ 
From this and the estimates in (3.15) one can
deduce that 
$$U_{q_0 ,w,u}^{\circ}\left(\psi,c';M,N\right) 
\ll\int_1^{2(1+T_0)}T^{-3/2}{\cal E}(T){\rm d}T +\int_{2(1+T_0)}^{\infty}T^{-3}{\cal E}(T)
{\rm d}T\;,\eqno(3.17)$$
where, by (3.10) and the arithmetic-geometric mean inequality, 
$${\cal E}(T)=\int_{-T}^T |\sigma_{q_0 ,w,u}\left( c';M,N;t\right) |{\rm d}t 
\leq{1\over 2}\int_{-T}^T 
{\sum_{\alpha , \delta\bmod c'{\frak O}}}^{\!\!\!\!\!\!*}\ \ \, 
\sum_{m=-M}^M
\left(\left| A\left( {\alpha\over c'}\,,m,t\right)\right|^2 
+\left| A\left( -{\delta\over c'}\,,m,t\right)\right|^2\right) {\rm d}t\;.$$
In this last upper bound the variable $\alpha$ is 
dependent upon the variable of summation $\delta$, 
and the conditions of summation (on both $\delta$ and $\alpha$) 
are those described in (2.13)-(2.15) of Lemma 2.3 (note already having been made of this 
below (3.8)). Since the very last part of Lemma 2.3 implies that the
the function mapping $\delta\bmod c'{\frak O}$ to $\alpha\bmod c'{\frak O}$
is injective (as is the function mapping $\delta\bmod c'{\frak O}$ to 
$-\delta\bmod c'{\frak O}$\/), one is therefore able to deduce from the last bound 
for ${\cal E}(T)$ that 
$${\cal E}(T) 
\leq\int_{-T}^T\,  
\sum_{\lambda\bmod c'{\frak O}}\ 
\sum_{m=-M}^M
\left| A\left( {\lambda\over c'}\,,m,t\right)\right|^2{\rm d}t\;.
$$
Hence (and by (3.16)) an appeal to the case $\alpha =(4\pi)^{-1}$, $\beta =-1$
of Lemma 2.6  yields the bound:   
$${\cal E}(T)
\ll \left( \left| c'\right| (M+1) +N^{1/2}\right)\left( \left| c'\right| T+N^{1/2}\right) 
\!\!\!\sum_{\scriptstyle \omega\in{\frak O}\atop\scriptstyle N/2<|\omega|^2\leq N}
\!\!\!\left| a(\omega)\right|^2\quad\qquad\hbox{($T\geq 1$).}$$ 
It follows by (3.17), the equality in (3.14), (3.12), (3.7) and (3.5) (in which $\epsilon\in{\frak O}^{*}$\/) 
that one has:  
$$\eqalign{
U_{q_0 ,w,u}^{\circ}\left(\psi,c';M,N\right) 
 &\ll \left( \left| c\right| (M+1)+N^{1/2}\right)\left( \left| c\right|  (1+T_0)^{1/2}+N^{1/2}\right)
\left\| {\bf b}_{N}\right\|_2^2 \asymp\cr
 &\asymp \left( \left| c\right| (M+1)+N^{1/2}\right)\left( \left| c\right| 
+\left| c\right|^{1/2}N^{1/4}|\psi|^{1/2}+N^{1/2}\right) \left\| {\bf b}_{N}\right\|_2^2  \leq\cr
 &\leq
{\sqrt{5}\over 2}\,(1+|\psi|)^{1/2}\left( \left| c\right| (M+1)+N^{1/2}\right)
\left(\left| c\right|+N^{1/2}\right) \left\| {\bf b}_N\right\|_2^2\;,}$$
where the $\left\| {\bf b}_N\right\|_2$ notation is as defined in (1.9.16).  
This (given (3.6), (3.9) and (3.10)) completes the proof
of the result (1.9.27) of the proposition\quad$\square$  

\medskip 

\noindent{\bf Part III: a proof of (1.9.30) by use of 
the Cauchy-Schwarz inequality and Lemma 2.6.}\quad  
It now only remains to prove the conditional bound (1.9.30): so 
suppose now that $c$ and $\psi$ satisfy the additional constraints in (1.9.28) and (1.9.29). 
Then, on setting 
$$L=\min\left\{ M\,,\,\left[ {N^{(1-\varepsilon)/2}\over\left| c'\right|}\right]\right\}\in{\Bbb N}\cup\{ 0\}\eqno(3.18)$$ 
(with $c'$ as indicated by (3.5) and (3.1)-(3.3)), 
it follows by subdivision of the outer sum in (3.6) (and positivity of the terms in this sum) 
that, 
for some pair of integers $M_1,M_2$ satisfying both $-M\leq M_1\leq M_2\leq M$ and 
$M_2-M_1=2L$,  one has: 
$$\eqalignno{U_{\frak a}(\psi,c;M;N,b)
 &\ll\left( {M+1\over L+1}\right) 
\sum_{m=M_1}^{M_2}
\Biggl|\quad\ \;  
\sum\!\!\!\!\!\!\!\!\!\!\!\sum_{\!\!\!\!\!\!\!\!\!\!\!{\scriptstyle \omega_1,\omega_2\in{\frak O}\atop
\scriptstyle N/2<\left|\omega_1\right|^2,\left|\omega_2\right|^2\leq N}}
\!\!\!\!\!\overline{ a\!\left(\omega_1\right)}\,a\!\left(\omega_2\right) 
\left({\omega_1\omega_2\over\left|\omega_1\omega_2\right|}\right)^{\!\!m}
{\rm e}\!\left(
{\psi\sqrt{\left|\omega_1\omega_2\right|}\over \left| c'\right|}\right)
S^{*}\!\left(\omega_1,\omega_2;c'\right)\Biggr| =\cr 
 &=\left( {M+1\over L+1}\right) U_{q_0 ,w,u}^{\star}\left(\psi,c';M_1,M_2;N\right)\qquad\quad 
\hbox{(say).}
&(3.19)}$$
After substituting $m+L+M_1$ for $m$, one finds that one has in (3.19):  
$$U_{q_0 ,w,u}^{\star}\left(\psi,c';M_1,M_2;N\right)
=\sum_{m=-L}^{L}
\eta_m'\quad 
\sum\!\!\!\!\!\!\!\!\!\!\!\sum_{\!\!\!\!\!\!\!\!\!\!\!{\scriptstyle \omega_1,\omega_2\in{\frak O}\atop
\scriptstyle N/2<\left|\omega_1\right|^2,\left|\omega_2\right|^2\leq N}}
\!\!\!\!\!\!\!\overline{ a^{-}\!\left(\omega_1\right)}\,a^{+}\!\left(\omega_2\right) 
\left({\omega_1\omega_2\over\left|\omega_1\omega_2\right|}\right)^{\!\!m}
\!{\rm e}\!\left( {\psi\sqrt{\left|\omega_1\omega_2\right|}\over \left| c'\right|}\right)
S^{*}\!\left(\omega_1,\omega_2;c'\right) , 
$$
where 
$$a^{\pm}(\omega)=a(\omega)\left( {\omega\over |\omega|}\right)^{\pm\left( L+M_1\right)}\qquad 
\hbox{($0\neq\omega\in{\frak O}$\/)}\eqno(3.20)$$
and the coefficients $\eta_m'$ are certain complex numbers satisfying 
$$\left|\eta_m'\right| =1\qquad 
\hbox{for $m=-L,-L+1,\ \ldots\ ,L\,$.}\eqno(3.21)$$
Using a change in the order of summation, and a subdivision of the sums over both $m$ and $\omega_2$, 
it may now be deduced that, for each $H\in{\Bbb N}$, there exists some positive integer $h=h(H)\leq H$ 
such that 
$$\eqalign{
\left| U_{q_0 ,w,u}^{\star}\left(\psi,c';M_1,M_2;N\right)\right| 
 &\leq H\Biggl| \sum_{r=0}^1\sum_{\scriptstyle \omega_1\in{\frak O}\atop
\scriptstyle N/2<\left|\omega_1\right|^2\leq N}
\!\!\!\overline{ a^{-}\!\left(\omega_1\right)}\ \times \cr
 &\qquad\times\sum_{\scriptstyle\omega_2\in{\frak O}\atop
\scriptstyle {\cal N}(h)<\left|\omega_2\right|^2\leq {\cal N}(h-1)}
\!\!\!\!\!a^{+}\!\left(\omega_2\right) 
{\cal L}_r\left(\omega_1\omega_2\right) 
{\rm e}\!\left( {\psi\sqrt{\left|\omega_1\omega_2\right|}\over \left| c'\right|}\right)
S^{*}\!\left(\omega_1,\omega_2;c'\right)\Biggr| ,}\eqno(3.22)$$
with 
$$\eqalignno{
{\cal N}(j) &=2^{-j/H} N\in\left[ N/2 , N\right]\qquad\hbox{($j=0,1,\ldots ,H$\/),} &(3.23)\cr 
 &\hbox{\quad} \cr
{\cal L}_r(z) &=\sum_{\scriptstyle -L\leq m\leq L\ \atop\scriptstyle m\equiv r\!\!\!\!\!\pmod{2}}
\eta_m'\left({z\over |z|}\right)^{\!\!m}\qquad 
\hbox{($z\in{\Bbb C}-\{ 0\}$\/)} &(3.24)
}$$
and $a^{\pm}(\omega)$ and $\eta_m'$ as in (3.20)-(3.21).
\par 
For later working it suffices that one takes, in the above,   
$$H=\bigl[ N^{\varepsilon /2}\bigr]\in{\Bbb N}\eqno(3.25)$$
(and, of course, $h=h(H)$\/). Then, by applying the 
Cauchy-Schwarz inequality to bound the sum on the right-hand side of (3.22), one has 
(see (3.20), (3.7) and (1.9.16)): 
$$\left|U_{q_0 ,w,u}^{\star}\left(\psi,c';M_1,M_2;N\right)\right|^2
\leq 2 H^2 {\cal U}_H\left\| {\bf a}^{-}_N\right\|_2^2 = 2H^2 {\cal U}_H\left\| {\bf b}_N\right\|_2^2
\ll N^{\varepsilon}{\cal U}_H\left\| {\bf b}_N\right\|_2^2\;,\eqno(3.26)$$
where
$${\cal U}_H=\sum_{r=0}^1\sum_{\scriptstyle \omega\in{\frak O}\atop\scriptstyle N/2<|\omega|^2\leq N}
\left|\sum_{\scriptstyle \omega'\in{\frak O}\atop\scriptstyle {\cal N}(h)<|\omega'|^2\leq{\cal N}(h-1)}
\!\!\! a^{+}\!\left(\omega'\right) {\cal L}_r\left(\omega\omega'\right) 
{\rm e}\!\left({\psi\sqrt{|\omega\omega'|}\over\left| c'\right|}\right)
S^{*}\!\left(\omega,\omega';c'\right)\right|^2\;.$$
\par
Choose now an infinitely differentiable function $G : {\Bbb R}\rightarrow [0,1]$ 
such that 
$$G(x)=\cases{1 &if $1/2\leq x\leq 1$,\cr 0 &if $x\leq 1/4$ or $x\geq 2$,}$$
and define $g$ to be the complex function such that 
$$g(z)=G\left( {|z|^2\over N}\right)\qquad\hbox{for $z\in{\Bbb C}\,$.}$$
Given that $G$ is chosen independently of all 
factors other than the above explicit requirements (as may, and shall, be assumed), 
one has $G^{(j)}(x)\ll_j 1$ for $j=0,1,2,\ldots\ $ and $x\in{\Bbb R}$. 
Moreover, since $g(\omega)=1$ if $N/2<|\omega|^2\leq N$, and is otherwise positive or zero, 
it follows that 
$${\cal U}_H\leq
\sum_{r=0}^1\sum_{\omega\in{\frak O}} g(\omega)
\Biggl|\sum_{\scriptstyle \omega'\in{\frak O}\atop\scriptstyle {\cal N}(h)<|\omega'|^2\leq{\cal N}(h-1)}
\!\!\!\!\!\!a^{+}\!\left(\omega'\right) S^{*}\!\left(\omega,\omega';c'\right) {\cal L}_r\left(\omega\omega'\right) 
{\rm e}\!\left({\psi\sqrt{|\omega\omega'|}\over\left| c'\right|}\right)\Biggr|^2\;.$$ 
This, via recall of the definitions, (3.8) and (3.24),  of $S^{*}\!\left(\omega,\omega';c'\right)$
and ${\cal L}_r(z)$,  gives: 
$$\eqalign{{\cal U}_H &
\leq\quad{\sum_{\delta_1\bmod c'{\frak O}}}^{\!\!\!\!\!\!\!\!*}\quad 
\ {\sum_{\delta_2\bmod c'{\frak O}}}^{\!\!\!\!\!\!\!\!*}
\qquad\     
\sum\!\!\!\!\sum_{\!\!\!\!\!\!\!\!\!\!\!\!\!\!\!\!\!{\scriptstyle -L\leq m_1,m_2\leq L\atop
\scriptstyle m_1\equiv m_2\!\!\!\!\!\pmod{2}}}
\eta_{m_1}'\,\overline{\eta_{m_2}'}\qquad\  
\sum\!\!\!\!\!\!\!\!\!\!\!\!\!\!\!\!\!\sum_{\!\!\!\!\!\!\!\!\!\!\!\!\!\!{\scriptstyle \omega_1',\omega_2'\in{\frak O}\atop
\scriptstyle {\cal N}(h)<\left|\omega_1'\right|^2,\left|\omega_2'\right|^2\leq{\cal N}(h-1)}}
\!\!\!\!\!a^{+}\!\left(\omega_1'\right)\,\overline{a^{+}\!\left(\omega_2'\right)}\,
\ \times\cr
 &\qquad\quad\times \left({\omega_1'\over |\omega_1'|}\right)^{\!\! m_1}
\left({\omega_2'\over |\omega_2'|}\right)^{\!\! -m_2}
{\rm e}\!\left( {\rm Re}\left( {\delta_1\omega_1'-\delta_2\omega_2'\over c'}\right)\right) 
W\!\left(\delta_1,\delta_2;m_1,m_2;\omega_1',\omega_2'\right) ,}\eqno(3.27)$$
with the asterisks modifying the sums $\bmod\ c'{\frak O}$ in the same way as 
in (2.12)-(2.13) of Lemma 2.3, while   
$$W\!\left(\delta_1,\delta_2;m_1,m_2;{\omega_1'},\omega_2'\right) 
=\sum_{\omega\in{\frak O}}g(\omega)\left({\omega\over |\omega|}\right)^{m_1-m_2} 
{\rm e}\!\left( {\rm Re}\!\left({\left(\alpha_1-\alpha_2\right)\over c'}\,\omega\right) 
+{\psi\left(\sqrt{|\omega_1'|} -\sqrt{|\omega_2'|}\right)\over\left| c'\right|}\,
|\omega|^{1/2}\right)$$
where $\alpha_i\bmod c'{\frak O}$ is the (unique) solution 
of the case $\delta\equiv\delta_i\bmod c'{\frak O}$ of (2.14)-(2.15) of Lemma 2.3. 
\par
The summand of the last sum above may be expressed as $F\bigl({\rm Re}(\omega),{\rm Im}(\omega) \bigr)$  
where, since $g$ was chosen so that $g(x+iy)=G\left(\left( x^2+y^2\right) /N\right)$ for $x,y\in{\Bbb R}$ 
(with the real function $G$ having derivatives of all orders, and support that is a compact subset of $(0,\infty)$\/), 
one is able to deduce that the function $F : {\Bbb R}^2\rightarrow {\Bbb C}$ lies in the Schwartz space. 
It therefore follows by (2.46)-(2.47) of Lemma 2.7 (Poisson's summation formula) that 
$$W\!\left(\delta_1,\delta_2;m_1,m_2;{\omega_1'},\omega_2'\right) 
=\sum_{\nu\in{\frak O}}\hat f(\nu) ,\eqno(3.28)$$
where 
$$f(z)=g(z)\left( {z\over |z|}\right)^{m_1-m_2}
{\rm e}\!\left( {\psi\left(\sqrt{|\omega_1'|} -\sqrt{|\omega_2'|}\right)\sqrt{|z|}\over
\left| c'\right|}\right)
{\rm e}\!\left( {\rm Re}\!\left({\left(\alpha_1-\alpha_2\right)\over c'}\,z\right) \right) \qquad
\hbox{($z\in{\Bbb C}$\/),}$$
so that in cases relevant to (3.27) one has (see (2.46)): 
$$\hat f(\nu) =\hat\varphi\left(\nu +{\alpha_2-\alpha_1\over c'}\right)\qquad\hbox{($\nu\in{\frak O}$)}\eqno(3.29)$$
with $\varphi : {\Bbb C}\rightarrow{\Bbb C}$ given by 
$$\varphi(z)=g(z)\left( {z\over |z|}\right)^{\!\!2d} {\rm e}\!\left( D\sqrt{|z|}\right)\qquad
\hbox{($z\in{\Bbb C}$),}\eqno(3.30)$$
where 
$$d= {m_1-m_2\over 2}\in{\Bbb Z}\qquad\hbox{and}\qquad 
D={\psi\left(\sqrt{|\omega_1'|} -\sqrt{|\omega_2'|}\right)\over\left| c'\right|}\in{\Bbb R}\;.\eqno(3.31)$$
Note here that, for 
the choices of $m_1$, $m_2$, $\omega_1'$ and $\omega_2'$ 
permitted in the sum on the right-hand side of (3.27),
one has (using (3.18) and 
(1.9.29), (3.23) and (3.25)): 
$$|d|\leq L\leq{N^{(1-\varepsilon)/2}\over\left| c'\right|}\;.\eqno(3.32)$$
and
$$|D|\leq A_2\,{\left( {\cal N}(h-1)-{\cal N}(h)\right)\over 4\left| c'\right| (N/2)^{3/4}}
\asymp {\left( 2^{1/H}-1\right) N^{1/4}\over\left| c'\right|}\asymp {N^{(1/4)-(\varepsilon /2)}\over\left| c'\right|}\;. 
\eqno(3.33)$$
\par
By the result (2.49) of Lemma 2.8, estimates for $\hat\varphi (w)$ ($w\in{\Bbb C}$) 
follow from bounds on the functions ${\bf\Delta}_{\Bbb C}^j\varphi\,$ ($j=0,1,2,\ldots\ $), where 
${\bf\Delta}_{\Bbb C}$ is Laplace's operator. Obtaining these bounds requires some preparation. 
Firstly, on recalling that $g(z)=G\left( |z|^2 /N\right)$, one may 
rewrite (3.30) to obtain: 
$$\varphi(z)=\Phi\!\left( |z|^2\right) z^d\left(\overline{z}\right)^{-d}\qquad\hbox{($z\in{\Bbb C}-\{ 0\}$),}\eqno(3.34)$$
where 
$$\Phi(x)=G(x/N)\,{\rm e}\!\left( D x^{1/4}\right)\qquad\hbox{($x>0$).}$$
It is here moreover the case that, since $G^{(j)}(x)\ll_j 1$ and 
$\left({\partial/\partial x}\right)^j {\rm e}\left( D x^{1/4}\right)\ll_j\left( 1+D x^{1/4}\right)^j x^{-j}$ 
for $x>0$ and $j\in{\Bbb N}\cup\{ 0\}$, and since the support of $G$ is contained in the interval $[1/4,2]$, 
the formulae of Leibniz for the 
higher order derivatives of a product enable one to deduce that 
$$x^j \Phi^{(j)}(x)\ll_j \left( 1+N^{1/4} |D|\right)^j\qquad\hbox{($x>0$ and $j=0,1,2,\ldots\ $).}\eqno(3.35)$$
One may next use the decomposition 
$${\bf\Delta}_{\Bbb C} =4\,{\partial\over\partial z}\,{\partial\over\partial\overline{z}}\;,\eqno(3.36)$$
where the linear operators 
${\partial /\partial z}$ and ${\partial /\partial\overline{z}}$ 
are defined as is indicated in (1.2.7). 
The decomposition (3.36) is 
useful for the matter in hand: for, when 
$q(z)$ is a complex function holomorphic on ${\Bbb C}-\{ 0\}$ (say), one has 
$${\partial\over\partial z}\,q(z)=q'(z),\quad 
{\partial\over\partial z}\,q\left(\overline{z}\right)=0,\quad 
{\partial\over\partial \overline{z}}\,q(z)=0,\quad\hbox{and}\quad
{\partial\over\partial \overline{z}}\,q\left(\overline{z}\right)=q'\left(\overline{z}\right)\qquad\quad
\hbox{($z\neq 0$).}$$
By the last observation, and the chain and product rules, one finds that 
$${\partial\over\partial z}\,\Phi^{(j)}\left( |z|^2\right) ={\overline z}\,\Phi^{(j+1)}\left( |z|^2\right)\quad 
\hbox{and}\quad 
{\partial\over\partial \overline{z}}\,\Phi^{(j)}\left( |z|^2\right) =z\Phi^{(j+1)}\left( |z|^2\right)\qquad 
\hbox{($z\neq 0$ and $j=0,1,2,\ldots\ $).}$$
By combining the facts just noted with Leibniz's formulae for derivatives of a product one
may deduce from (3.34) and (3.36) that, for $0\neq z\in{\Bbb C}$ and $j=0,1,2,\ldots\ $, 
$$\eqalign{ 
{\bf\Delta}_{\Bbb C}^j\varphi(z) 
 &=\left( 4\,{\partial\over\partial z}\,{\partial\over\partial\overline{z}}\right)^j 
\Phi\left( |z|^2\right) z^d\left(\overline{z}\right)^{-d} =\cr 
 &=\left( 4\,{\partial\over\partial z}\right)^j
\sum_{r=0}^j\pmatrix{j\cr r}z^{d+r}\,\Phi^{(r)}\!\left( |z|^2\right)
(-d)(-d-1)\cdots(-d-(j-r-1))\left(\overline{z}\right)^{-d-(j-r)} =\cr 
 &=4^j\sum_{r=0}^j {j!\over r!}\pmatrix{-d\cr j-r}\left(\overline{z}\right)^{-d-(j-r)}
\left( {\partial\over\partial z}\right)^j z^{d+r}\,\Phi^{(r)}\!\left( |z|^2\right) ,
}$$
where $\pmatrix{m\cr n}$ is the coefficient of $x^n$ in the binomial expansion of 
$(1+x)^m$. 
Similarly, 
$$\left( {\partial\over\partial z}\right)^j z^{d+r}\,\Phi^{(r)}\!\left( |z|^2\right) 
=\sum_{s=0}^j {j!\over s!}\pmatrix{d+r\cr j-s}\left(\overline{z}\right)^s
\Phi^{(r+s)}\left( |z|^2\right) z^{d+r-(j-s)}\;,$$
so that one obtains, 
for $z\neq 0\,$ and $j=0,1,2,\ldots\,$,  
$${\bf\Delta}_{\Bbb C}^j\varphi(z) 
=\left( {4\over |z|^2}\right)^j \left( {z\over\overline{z}}\right)^d
\sum_{r=0}^j \sum_{s=0}^j {(j!)^2\over r!s!}\pmatrix{-d\cr j-r}\pmatrix{d+r\cr j-s}
|z|^{2(r+s)}\Phi^{(r+s)}\left( |z|^2\right) .$$
One may therefore deduce from the bounds in (3.35) that,  
for $0\neq z\in{\Bbb C}$ and $j=0,1,2,\ldots\ $, one has 
$${\bf\Delta}_{\Bbb C}^j\varphi(z) 
=O_j\!\left( |z|^{-2j}\right) \sum_{r=0}^j\sum_{s=0}^j
O_{j,r,s}\!\left( (1+|d|)^{(j-r)+(j-s)}\right) 
O_{r+s}\!\left( \left( 1+N^{1/4} |D|\right)^{r+s}\right)  
\ll_j \left( {1+V\over |z|^2}\right)^{\!j} ,\eqno(3.37)$$
where (given (3.32) and (3.33)) 
$$V=\max\left\{  d^2\,,\, N^{1/2}D^2\right\}\ll {N^{1-\varepsilon}\over\left| c'\right|^2}\;.\eqno(3.38)$$
Since one has here $N^{1-\varepsilon}/\left| c'\right|^2\gg 1$ (by (3.5) and (1.9.28)), 
and since $\varphi(z)=0$ unless $N/4\leq |z|^2\leq 2N$ (by (3.30) and the choice of $g$), 
one finds, by means of the result (2.49) of Lemma 2.8, that the bounds in (3.37) and (3.38) imply 
that one has  
$$\hat\varphi (w)  
\ll_j N |w|^{-2j}\left( {N^{1-\varepsilon}/\left| c'\right|^2\over N}\right)^j
=N^{1-j\varepsilon}\left| c' w\right|^{-2j}\qquad\hbox{($0\neq w\in{\Bbb C}$ and $j=0,1,2,\ldots\ $).}
\eqno(3.39)$$

The aim now is to apply (3.39) in order to estimate the sum 
$\sum_{\nu\in{\frak O}}\hat f(\nu)$ on the right-hand side of (3.28). 
Let $\nu_0$ be the unique Gaussian integer such that 
$\nu_0+\left(\alpha_2 -\alpha_1\right) /c'\in{\Bbb Q}(i)$ 
has both its real and imaginary parts lying in the interval $(-1/2,1/2]$. 
Then, for $\nu\in{\frak O}$ with $\nu\neq\nu_0$, one has 
$$\left|\nu+{\alpha_2 -\alpha_1\over c'}\right|
\geq\left| {\nu}-{\nu}_0\right|
-\left|\nu_0+{\alpha_2 -\alpha_1\over c'}\right|
\geq\left| {\nu}-{\nu}_0\right| -{1\over\sqrt{2}}
\gg\left| {\nu}-{\nu}_0\right|=\left|\nu'\right|\quad\hbox{(say).}$$
From this, (3.29) and the case $j=2+[2/\varepsilon]$ of (3.39)
it follows that
$$\eqalignno{
\sum_{{\nu}\in{\frak O}}\hat f(\nu)
 &=\hat\varphi\left( {\nu}_0+{\alpha_2 -\alpha_1\over c'}\right)
+\sum_{0\neq\nu'\in{\frak O}} 
O_{\varepsilon}\left( N^{1-j\varepsilon}\left| c'\right|^{-2j}\left| \nu'\right|^{-2j}\right) =\cr
 &=\hat\varphi\left( {\nu}_0+{\alpha_2 -\alpha_1\over c'}\right)
+O_{\varepsilon}\left( N^{-1}\left| c'\right|^{-2j}\right) , &(3.40)}$$
where, by (3.3) and (3.5) (in which $vw\in{\frak O}$ and 
$\epsilon\in{\frak O}^{*}$\/), one has 
$0\neq c'\in{\frak O}$, so that $\left| c'\right|\geq 1$. If it is here the case that 
$\nu_0+\left(\alpha_2 -\alpha_1\right)\!/c'\neq 0$, then one has 
$$\left|\nu_0+{\alpha_2 -\alpha_1\over c'}\right| 
={\left| c'\nu_0 +\alpha_2 -\alpha_1\right|\over\left| c'\right|}
\geq {1\over\left| c'\right|}\,,$$
so that, by (3.39), 
$$\hat\varphi\left( \nu_0+{\alpha_2 -\alpha_1\over c'}\right)
\ll_j N^{1-j\varepsilon}\qquad\hbox{($j=0,1,2,\ldots\ $).}$$
Hence (and by (3.28) and (3.40)) one obtains: 
$$W\!\left(\delta_1,\delta_2;m_1,m_2;{\omega_1'},\omega_2'\right) 
=O_{\varepsilon}\left( N^{-1}\right) 
+\cases{\hat\varphi (0) &if  $\left(\alpha_2 -\alpha_1\right) /c'\in{\frak O}\,$; \cr
0\matrix{\hbox{\ }\cr\hbox{\ }} &otherwise.}\eqno(3.41)$$
It should be noted here that $\left(\alpha_2 -\alpha_1\right) /c'\in{\frak O}$
if and only if $\alpha_1\equiv\alpha_2\bmod c'{\frak O}$. 
Moreover, for $i=1,2$, the relationship between $\alpha_i\bmod c'{\frak O}$ and 
$\delta_i\bmod c'{\frak O}$ is (as noted below (3.27)) the same as 
that existing between  the variables $\alpha\bmod c'{\frak O}$ and 
$\delta\bmod c'{\frak O}$ in (2.13)-(2.15) of Lemma 2.3, so that by virtue of 
the final point noted in Lemma 2.3 one has 
$\alpha_1\equiv\alpha_2\bmod c'{\frak O}$ if and only if 
$\delta_1\equiv\delta_2\bmod c'{\frak O}$. In addition to this, it is easily seen 
that $\hat\varphi (0)=0$ unless a further independent condition (on $d$\/) is satisfied: 
for, by (3.30), one has 
$$\hat\varphi (0)=\int\limits_{-\infty}^{\infty}\int\limits_{-\infty}^{\infty} g( z )
\left( {z\over |z|}\right)^{2d}{\rm e}\left( D\sqrt{| z |}\right) {\rm d}x{\rm d}y$$
(where $z$ is a dependent variable satisfying $z=x+iy$\/), and so, given that 
$g(z)=g\left( |z|^2 /N\right) =g(|z|)$ and that $d\in{\Bbb Z}$, 
a change to polar coordinates shows that  
$$\hat\varphi (0)
=\int\limits_0^{2\pi}\int\limits_0^{\infty} g(r)e^{2id\theta}
\,{\rm e}\!\left( D\sqrt{r}\right)
 r{\rm d}r \,{\rm d}\theta
=\int\limits_0^{2\pi}e^{2id\theta}{\rm d}\theta
\int\limits_0^{\infty} g(r)\,{\rm e}\!\left( D\sqrt{r}\right) r{\rm d}r
=E(2d)\int\limits_0^{\infty} r g(r)\,{\rm e}\!\left( D\sqrt{r}\right) {\rm d}r\;,$$
where $E(n)=2\pi$ if $n=0$, while 
$E(n)=0$ if $0\neq n\in{\Bbb Z}$.
Since the relevant values of $d$ and $D$ are those given by (3.31), it  follows from 
(3.41) and the points just noted that one has: 
$$W\!\left(\delta_1,\delta_2;m_1,m_2;{\omega_1'},\omega_2'\right) 
=O_{\varepsilon}\left( N^{-1}\right)
+\cases{2\pi\tilde{g}\left( \omega_1',\omega_2'\right) &if $m_1=m_2$ 
and $\delta_1\equiv \delta_2\bmod c'{\frak O}$, \cr 0 &otherwise ,}\eqno(3.42)$$ 
where
$$\eqalignno{
\tilde{g}\left( \omega_1',\omega_2'\right)
 &=\int\limits_0^{\infty}rg(r){\rm e}\left(
{\psi\left(\sqrt{|\omega_1'|}-\sqrt{|\omega_2'|}\right)\sqrt{r}\over \left| c'\right|}\right) {\rm d}r =\cr 
 &=2\!\!\!\!\!\int\limits_0^{\quad (2N)^{1/4}}\!\!t^3 G\left( {t^4\over N}\right) {\rm e}\left(
{\psi\left(\sqrt{|\omega_1'|}-\sqrt{|\omega_2'|}\right)t\over \left| c'\right|}\right) {\rm d}t\;.
&(3.43)}$$
\par
By (3.27), (3.21) and (3.42), one obtains the upper bound
$${\cal U}_H\leq {\cal U}_H'+2\pi {\cal U}_H'',\eqno(3.44)$$
where, by (3.20), (3.23) and (3.7),  
$$\eqalignno{{\cal U}_H'
 &={\sum_{\delta_1\bmod c'{\frak O}}}^{\!\!\!\!\!\!\!*}\quad 
{\sum_{\delta_2\bmod c'{\frak O}}}^{\!\!\!\!\!\!\!*}\qquad\,    
\sum\!\!\!\!\sum_{\!\!\!\!\!\!\!\!\!\!\!\!\!\!\!\!\!{\scriptstyle -L\leq m_1,m_2\leq L\atop
\scriptstyle m_1\equiv m_2\!\!\!\!\!\pmod{2}}}
\qquad\quad\, 
\sum\!\!\!\!\!\!\!\!\!\!\!\!\!\!\!\sum_{\!\!\!\!\!\!\!\!\!\!\!\!\!\!\!\!{\scriptstyle 
\omega_1',\omega_2'\in{\frak O}\atop\scriptstyle 
{\cal N}(h)<|\omega_1'|^2,|\omega_2'|^2\leq{\cal N}(h-1)}} 
\!\!\!\!\!\!\!\!\left| a^{+}\!\left(\omega_1'\right) 
a^{+}\!\left(\omega_2'\right)\right| 
O_{\varepsilon}\left( N^{-1}\right) \ll\cr
 &\ll \left| c'\right|^4 (L+1)^2 
\quad  
\sum\!\!\!\!\!\!\!\sum_{\!\!\!\!\!\!\!\!\!\!\!\!\!\!\!\!\!{\scriptstyle 
\omega_1',\omega_2'\in{\frak O}\atop\scriptstyle 
N/2<|\omega_1'|^2,|\omega_2'|^2\leq N}} 
\!\!\!\!\!\!\!\!\left(\left| b\!\left(\omega_1'\right)\right|^2 + 
\left| b\!\left(\omega_2'\right)\right|^2\right) O_{\varepsilon}\!\left( N^{-1}\right) 
\ll_{\varepsilon} \left| c'\right|^4 \left( L+1\right)^2 \left\| {\bf b}_N\right\|_2^2&(3.45)
 }$$
(with $\left\| {\bf b}_N\right\|_2$ as in (1.9.16)), while by (3.43) one has  
$$\eqalignno{
{\cal U}_H''
 &={\sum_{\delta\bmod c'{\frak O}}}^{\!\!\!\!\!\!*}\quad 
\!\!\sum_{m=-L}^L\qquad\quad
\sum\!\!\!\!\!\!\!\!\!\!\!\!\!\!\!\sum_{\!\!\!\!\!\!\!\!\!\!\!\!\!\!\!\!{\scriptstyle 
\omega_1',\omega_2'\in{\frak O}\atop\scriptstyle 
{\cal N}(h)<|\omega_1'|^2,|\omega_2'|^2\leq{\cal N}(h-1)}} 
\!\!\!\!\!\!\!\!a^{+}\!\left(\omega_1'\right)\!\left({\omega_1'\over |\omega_1'|}\right)^{\!\!m}
\overline{a^{+}\!\left(\omega_2'\right)\!\left({\omega_2'\over |\omega_2'|}\right)^{\!\! m}}
{\rm e}\!\left( {\rm Re}\!\left(\! {\delta\left( \omega_1'-\omega_2'\right) \over c'}\!\right)\!\right) 
\tilde{g}\left( \omega_1',\omega_2'\right) =\cr 
 &=2{\sum_{\delta\bmod c'{\frak O}}}^{\!\!\!\!\!\!*}\quad 
\!\!\sum_{m=-L}^L\int_0^{\quad (2N)^{1/4}} t^3 G\left( {t^4\over N}\right) 
\left| s(\delta , m, t )\right|^2 {\rm d}t  \ll\cr 
 &\ll N^{3/4}{\sum_{\delta\bmod c'{\frak O}}}^{\!\!\!\!\!\!*}\quad 
\!\!\sum_{m=-L}^L\int_{-T}^T \left| s(\delta , m, t )\right|^2 {\rm d}t\;, &(3.46)
}$$ 
with $T=(2N)^{1/4}$ and 
$$s(\delta ,m,t)
=\sum_{\scriptstyle\omega\in{\frak O}\atop\scriptstyle {\cal N}(h)<|\omega|^2\leq{\cal N}(h-1)} 
a^{+}(\omega)\left( {\omega\over |\omega|}\right)^m 
{\rm e}\!\left( {\rm Re}\left( {\delta\omega\over c'}\right) 
+{\psi\sqrt{|\omega|}\ t\over\left| c'\right|}\right) .\eqno(3.47)$$
Given (3.23), and given that the hypothesis (1.9.29) implies $|\psi|>0$,
it follows by (3.46) and (3.47) that one may bound ${\cal U}_H''$ by applying the case 
$\alpha =\psi /|2c'|$, $\beta =-1/2$ of Lemma 2.6. 
By (3.23), (3.20) and (3.7), this application of Lemma 2.6 shows that 
$$\eqalignno{
{\cal U}_H''
 &\ll N^{3/4}\left( \left| c'\right| (L+1) +N^{1/2}\right) 
\left( \left| c'\right| T +|\psi|^{-1}\left| c'\right| N^{1/4}\right)
 \left\| {\bf b}_N\right\|_2^2 \ll\cr 
&\ll \left( \left| c'\right| (L+1) +N^{1/2}\right) |\psi|^{-1} \left| c'\right| N \left\| {\bf b}_N\right\|_2^2 
&(3.48)}$$
(the latter bound following since one has $N^{1/4}\asymp T$ and, by the hypothesis (1.9.29), 
$|\psi|^{-1}\gg 1$\/). 
By (3.5), the hypothesis (1.9.28) and (3.18), one has 
$$\left| c'\right|\leq\left| c'\right| (L+1)\ll N^{(1-\varepsilon)/2}\leq N^{1/2}\;,$$ 
so that from (3.44), (3.45) and (3.48) one may deduce: 
$$\eqalign{ 
{\cal U}_H &\leq\biggl( 
O_{\varepsilon}\left( \left| c'\right|^4 (L+1)^2\right) 
+O\left( \left(\left| c'\right| (L+1) +N^{1/2}\right)\left|\psi\right|^{-1}\left| c'\right| N\right)\biggr)
 \left\| {\bf b}_N\right\|_2^2 \ll\cr 
 &\ll\left( \left|\psi\right|^{-1} +O_{\varepsilon}\left( N^{-3\varepsilon /2}\right)\right) 
\left| c'\right| N^{3/2}\left\| {\bf b}_N\right\|_2^2\;.
}$$
By using this estimate for ${\cal U}_H$ in (3.26) and noting that $N^{\varepsilon}\geq 1$, 
one finds that 
$$U_{q_0 ,w,u}^{\star}\left(\psi,c';M_1,M_2;N\right)
\ll \left( \left|\psi\right|^{-1/2} +O_{\varepsilon}\left( 1\right)\right) 
\left| c'\right|^{1/2} N^{(3+2\varepsilon)/4}\left\| {\bf b}_N\right\|_2^2\;.\eqno(3.49)$$

Finally, since $\left| c'\right| =|c|$ (by (3.5)), and since the definition (3.18) implies that  
$${M+1\over L+1}\leq 1+{M\over L+1}\ll 1+\left| c'\right| M N^{-(1-\varepsilon)/2}\;,$$ 
one finds that the bounds (3.19) and (3.49) combine to give (1.9.30)\quad$\blacksquare$

\bigskip\par

\centerline{\bf \S 4. Further lemmas.}

\medskip

These lemmas are needed for the proof of Theorem 1, in the next section.

\bigskip

\proclaim Lemma 4.1. Let $0\neq q_0\in{\frak O}={\Bbb Z}[i]$ and 
$\tau\in\Gamma =\Gamma_0(q_0)\leq SL(2,{\frak O})$. 
Suppose that ${\frak a}$ and ${\frak a}'$ are cusps of $\Gamma$ such that 
$\tau{\frak a}={\frak a}'$. Let $g_{\frak a},g_{{\frak a}'}\in SL(2,{\Bbb C})$, with 
$g_{{\frak a}'}$ such that (1.1.16) and (1.1.20)-(1.1.21) hold for ${\frak c}={\frak a}'$. 
Then (1.1.16) and (1.1.20)-(1.1.21) hold for ${\frak c}={\frak a}$ if and only if one has 
$$\tau g_{\frak a}=g_{{\frak a}'}\,h\!\left[\eta\right] n[\beta]
\qquad\hbox{for some $\beta,\eta\in{\Bbb C}$ with $\eta^2\in{\frak O}^*$.}\eqno(4.1)$$

\medskip

\noindent{\bf Proof.}\quad 
The stated condition is necessary, since if (1.1.16) and (1.1.20)-(1.1.21) 
hold for ${\frak c}={\frak a}'$, and for ${\frak c}={\frak a}$, then by the result (2.1) of Lemma 2.1 
there must exist some   
$\eta\in{\Bbb C}$ with $\eta^2\in{\frak O}^*$ and some $\beta\in{\Bbb C}$ such that 
$g_{{\frak a}'}^{-1}\tau g_{\frak a}=h\!\left[\eta\right] n[\beta]$.
The stated condition is also sufficient, since if (1.1.16) and (1.1.20)-(1.1.21) hold for 
${\frak c}={\frak a}'$, then (4.1) implies both that 
$$g_{\frak a}\infty =\tau^{-1}g_{{\frak a}'}\,h\!\left[\eta\right] n[\beta]\infty
=\tau^{-1}g_{{\frak a}'}\infty =\tau^{-1}{\frak a}'={\frak a}$$
and that 
$$\eqalign{g_{\frak a}^{-1}\Gamma_{\frak a}' g_{\frak a}
 =g_{\frak a}^{-1}\tau^{-1}\Gamma_{{\frak a}'}'\tau g_{\frak a}
 &=\left(\tau g_{\frak a}\right)^{-1}\Gamma_{{\frak a}'}' \tau g_{\frak a} =\cr
 &=n[-\beta]\,h\!\left[ 1/\eta\right] g_{{\frak a}'}^{-1}\,\Gamma_{{\frak a}'}'\,g_{{\frak a}'}
\,h\!\left[\eta\right] n[\beta] =\cr
 &=n[-\beta]\,h\!\left[ 1/\eta\right] B^{+}
\,h\!\left[\eta\right] n[\beta] =\cr
 &=\left\{ n[-\beta]\,h\!\left[ 1/\eta\right] n[\alpha] 
\,h\!\left[\eta\right] n[\beta] : \alpha\in{\frak O}\right\} =\cr
 &=\left\{ n\!\left[\alpha /\eta^2\right] : \alpha\in{\frak O}\right\} =\cr
 &=\left\{ n\left[\alpha'\right] : \alpha'\in{\frak O}\right\} =B^{+}\quad\blacksquare
}$$

\bigskip

\proclaim Lemma 4.2. Let $0\neq q_0\in{\frak O}={\Bbb Z}[i]$; let  
$\Gamma =\Gamma_0(q_0)\leq SL(2,{\frak O})$; and let $B=B^{+}\cup h[-1]B^{+}$, 
where $B^{+}$ is as in (1.1.21). 
Suppose also that 
${\frak a}$ is a cusp of $\Gamma$, and let $\mu({\frak a})$ be as 
described in Theorem 1. 
Then there exists $g_{\frak a}\in SL(2,{\Bbb C})$ such that (1.1.16) and 
(1.1.20)-(1.1.21) hold for ${\frak c}={\frak a}$. 
Moreover, for each such $g_{\frak a}$ one has either 
$$q_0 \mu({\frak a})\mid 2\quad\hbox{and}\quad 
g_{\frak a}^{-1}\Gamma_{\frak a}\,g_{\frak a}
=B\cup h[i]\,n\!\left[\beta_{\frak a}\right] B\  
\ \hbox{for some $\ \beta_{\frak a}\in{\Bbb C}$,}\eqno(4.2)$$
or else 
$$q_0 \mu({\frak a})\!\not\;\mid 2\quad\hbox{and}\quad 
g_{\frak a}^{-1}\Gamma_{\frak a} g_{\frak a}=B\;.\eqno(4.3)$$

\medskip

\noindent{\bf Proof.}\quad 
Note firstly that, if ${\frak a}=u/w$ with non-zero $u,w\in{\frak O}$ such that $(u,w)\sim 1$ and 
$w\mid q_0$, then there exists a choice of 
$g_{\frak a}\in\left\{ g\in SL(2,{\Bbb C}) : g\infty ={\frak a}\right\}$ 
such that one of the two statements (4.2), (4.3) is true. 
Indeed, supposing $u,w\in{\frak O}$ to have the properties just listed, one may 
choose $\tilde u,\tilde w\in{\frak O}$ such that 
$$SL(2,{\frak O})\ni\pmatrix{u &-\tilde w\cr w &\tilde u}
=\tilde g_{u/w}\quad\hbox{(say).}\eqno(4.4)$$
Then $\tilde g_{u/w}\infty =u/w$. 
Moreover, since  $SL(2,{\frak O})$ is a group, it follows from (4.4) and 
(1.1.17)-(1.1.19) that 
$$\eqalign{ 
\tilde g_{u/w}^{-1}\Gamma_{u/w}\tilde g_{u/w} 
 &=\tilde g_{u/w}^{-1}\Gamma\tilde g_{u/w}\cap P =\cr 
 &=\left\{ p\in P\cap SL(2,{\frak O}) : \tilde g_{u/w} p\tilde g_{u/w}^{-1}\in\Gamma 
=\Gamma_0\!\left( q_0\right)\right\} =\cr
 &=\left\{\pmatrix{\alpha &z\cr 0&\alpha^{-1}} : \alpha\in{\frak O}^* , 
z\in{\frak O}\ {\rm and}\ \alpha \tilde u w-zw^2 
-\alpha^{-1}\tilde u w\in q_0 {\frak O}\right\} =\cr 
 &=\left\{\pm\pmatrix{1&z\cr 0&1} : z\in{\frak O},\,  
wz\in(q_0 /w){\frak O}\right\}\cup 
\left\{\pm\pmatrix{i&z\cr 0&-i} : z\in{\frak O},\,  
wz-2i\tilde u\in (q_0 /w){\frak O}\right\}\;.
}$$
Given that $\left( w , \tilde u\right)\sim 1$, the congruence  
$wz\equiv 2i\tilde u\bmod (q_0 /w){\frak O}$ is soluble if and only if 
$(w , q_0 /w)\mid 2$. Hence, and by (1.9.15), it follows from the above that 
one has either 
$$q_0 \mu(u/w)\mid 2\quad\hbox{and}\quad 
\tilde g_{u/w}^{-1}\Gamma_{\frak u/w} \tilde g_{u/w}
=\left\{\,h[\pm 1]\,n\!\left[ q'\zeta\right] : \zeta\in{\frak O}\right\}\cup 
\left\{ h[\pm i]\,n\!\left[-iz_0 +q'\zeta'\right] : \zeta'\in{\frak O}\right\}$$
or 
$$q_0 \mu(u/w)\!\!\not\ \mid 2\quad\hbox{and}\quad 
\tilde g_{u/w}^{-1}\Gamma_{u/w}\tilde g_{u/w}=
\left\{\,h[\pm 1]\,n\!\left[ q'\zeta\right] : \zeta\in{\frak O}\right\}\;,$$
where (in either case) $q'$ is any Gaussian integer such that 
$q'\sim (q_0 /w)/(w , q_0 /w)$, and where (in the former case) $z_0$ denotes an arbitrary 
Gaussian integer $z_0$ satisfying 
$${w\over\left( w , q_0 /w\right)}\,z_0 
\equiv {2i\tilde u\over\left( w , q_0 /w\right)}\bmod q'{\frak O}\;.$$
Therefore, 
and since $h[1/u] h[\alpha] n[\beta] h[u] =h[\alpha]\,n\!\left[\beta /u^2\right]$, 
$h[\pm i]=h[\pm 1]\,h[i]$ and 
$n[\beta' +\zeta']=n[\beta'] n[\zeta']$ 
(for $\alpha ,u\in{\Bbb C}^*$ and $\beta,\beta',\zeta'\in{\Bbb C}$\/), 
one may ensure that 
one of the two statements (4.2), (4.3) holds for ${\frak a}=u/w$
by putting $g_{u/w}=\tilde g_{u/w}\,h\!\left[\sqrt{q'}\right]$: 
note that for this choice of $g_{u/w}$ one also obtains the 
case ${\frak c}=u/w$ of (1.1.16), since one has  
$g_{u/w}\infty =\tilde g_{u/w}\,h\!\left[\sqrt{q'}\right]\infty =\tilde g_{u/w}\infty =u/w$.
 
Since the claim made at the beginning of this proof has now been verified, it 
now only remains to be shown that what was claimed there does imply the lemma. In accordance with 
this aim, let it now be supposed that ${\frak a}$ is some cusp of $\Gamma$. 
By the result (2.4) of Lemma 2.2, one has  $\tau{\frak a}=u/w$ 
for some $\tau\in\Gamma$, and some 
pair of non-zero Gaussian integers $u,w$ such that $(u,w)\sim 1$ and $w\mid q_0$; 
moreover, by virtue of the claim verified in the first part of this proof, 
it follows that one may assign to any such pair $u,w$  some  
$g_{u/w}\in SL(2,{\Bbb C})$ such that one has $g_{u/w}\infty =u/w$ and either
$$q_0 \mu({u/w})\mid 2\quad\hbox{and}\quad 
g_{u/w}^{-1}\Gamma_{u/w} g_{u/w}=B\cup h[i]\,n\!\left[\beta_{u/w}\right] B\  
\ \hbox{for some $\,\beta_{u/w}\in{\Bbb C}$,}\eqno(4.5)$$
or 
$$q_0 \mu({u/w})\!\not\,\mid 2\quad\hbox{and}\quad 
g_{u/w}^{-1}\Gamma_{u/w} g_{u/w}=B\;.\eqno(4.6)$$
Since $B=B^{+}\cup h[-1] B^{+}$, where $B^{+}\leq N\leq P$ and 
$N$ contains $n[\beta_{u/w}]$, but not $h[-1]$, $h[i]$ or $h[-i]$, 
it follows by (1.1.17) and (1.1.19) that in both 
the cases (4.5), (4.6) just mentioned one has 
$$g_{u/w}^{-1}\Gamma_{u/w}' g_{u/w}
=g_{u/w}^{-1}\Gamma_{u/w}\,g_{u/w}\cap N
=B^{+}\cap N=B^{+}\;,$$ 
so that (1.1.16) and (1.1.20)-(1.1.21) hold for ${\frak c}=u/w$. 
Therefore, and since one has $\tau{\frak a}=u/w$ (with $\tau\in\Gamma$\/), 
it follows by Lemma 4.1 that the elements $g_{\frak a}\in SL(2,{\Bbb C})$ 
such that (1.1.16) and (1.1.20)-(1.1.21) hold for ${\frak c}={\frak a}$ are 
precisely those elements $g_{\frak a}\in SL(2,{\Bbb C})$ which, for some 
$\beta,\eta\in{\Bbb C}$ with $\eta^2\in{\frak O}^*$, satisfy: 
$$\tau g_{\frak a}
=g_{u/w}\,h\!\left[\eta\right] n[\beta]\;.\eqno(4.7)$$
Since ${\frak O}^*\ni 1^2$, ${\Bbb C}\ni 0$ 
and $\tau^{-1}g_{u/w}\,h[1]n[0]=\tau^{-1}g_{u/w}\in SL(2,{\Bbb C})$, 
the observation just made 
completes the proof that there exists some $g_{\frak a}\in SL(2,{\Bbb C})$ such that 
(1.1.16) and (1.1.20)-(1.1.21) hold for ${\frak c}={\frak a}$. 

For the final result of the lemma one may note
that if (4.5) and (4.7) hold (with $\eta^2\in{\frak O}^*$\/) then, 
since $\tau{\frak a}=u/w$ (where $\tau\in\Gamma$), one will have 
$$\eqalignno{
g_{\frak a}^{-1}\Gamma_{\frak a}\,g_{\frak a} 
=g_{\frak a}^{-1}\tau^{-1}\Gamma_{u/w}\,\tau g_{\frak a}
 &=\left(\tau g_{\frak a}\right)^{-1}\Gamma_{u/w}\,\tau g_{\frak a} =\cr 
 &=n[-\beta]\,h\!\left[ 1/\eta\right] g_{u/w}^{-1} \Gamma_{u/w}\,g_{u/w}\, 
h\!\left[\eta\right] n[\beta] =\cr
 &=n[-\beta]\,h\!\left[ 1/\eta\right] \left( 
B\cup h[i]\,n\!\left[\beta_{u/w}\right] B\right) h\!\left[\eta\right] n[\beta] =\cr
 &=B\cup h[i]\,n\!\left[\eta^{-2}\beta_{u/w}+2\beta\right] B\;.&(4.8)}$$
If instead (4.6) and (4.7) hold  
(with $\beta,\eta\in{\Bbb C}$ and $\eta^2\in{\frak O}^*$)   
then one will have (similarly):   
$$g_{\frak a}^{-1}\Gamma_{\frak a}\,g_{\frak a} 
=n[-\beta]\,h\!\left[ 1/\eta\right] B h\!\left[\eta\right] n[\beta] 
=B\;.$$
In both the above two cases, the equation  
$\tau{\frak a}=u/w$ implies the relation $\mu( {\frak a})\sim\mu(u/w)\quad\blacksquare$ 

\medskip

\goodbreak\noindent{\bf Remark~4.3.}\quad 
By (4.8) it is evident that when $q_0 \mu( {\frak a} )\mid 2$ one may choose 
$g_{\frak a}$ such that (4.2) holds with $\beta_{\frak a}=0$. There will not, however, be any 
use made of this observation in what follows.

\bigskip

\proclaim Lemma 4.4 (a Fourier transform). For real $y$ put 
$$G_n(y)=\int_{-\infty}^{\infty} x^n \exp\left( 2ixy -x^2\right) {\rm d}x\qquad 
\hbox{($n=0,1,2,\ldots\ $),}$$
and set $G_{-1}(y)$ equal to zero. Then, when $y\in{\Bbb R}$, one has: 
$$G_{2m}(y)=\int_{-\infty}^{\infty} x^{2m} e^{-x^2}\cos(2xy) {\rm d}x
=2\int_{0}^{\infty} x^{2m} e^{-x^2}\cos(2xy) {\rm d}x\qquad 
\hbox{($m=0,1,2,\ldots\ $),}\eqno(4.9)$$
$$G_{2m+1}(y)=i\int_{-\infty}^{\infty} x^{2m+1} e^{-x^2}\sin(2xy) {\rm d}x
=2i\int_{0}^{\infty} x^{2m+1} e^{-x^2}\sin(2xy) {\rm d}x\qquad 
\hbox{($m=0,1,2,\ldots\ $),}\eqno(4.10)$$
$$G_0(y)=\sqrt{\pi}\,e^{-y^2}\eqno(4.11)$$ 
and 
$$G_{n+1}(y)=iyG_n(y)+{n\over 2}\,G_{n-1}(y)\qquad 
\hbox{($n=0,1,2,\ldots\ $).}\eqno(4.12)$$

\medskip

\noindent{\bf Proof.}\quad 
Absolute convergence of the above integrals (for $n,m\geq 0$) may be established by means of the 
upper bound $x^k\exp\left( -x^2\right)\ll_k \exp(-|x|)$ (valid for $x\in{\Bbb R}$, when $k\geq 0$). 
The results (4.9) and (4.10) follow from Euler's identity (applied to $\exp(2ixy)$\/), 
coupled with the fact that if $f : {\Bbb R}\rightarrow{\Bbb C}$ is an 
odd integrable function then the integral $\int_{-\infty}^{\infty}f(x){\rm d}x$ will, if 
it is absolutely convergent,  be equal to zero. Indications as to a proof of (4.11) are 
given in Exercise~10.22 of [1]. The definition of $G_n(y)$ implies that
$-2G_{n+1}(y)=\int_{-\infty}^{\infty} x^n\exp(2ixy)\,{\rm d}\!\left(\exp\left( -x^2\right)\right)$:
integrating by parts, one obtains (4.12)\quad$\blacksquare$ 
\bigskip

\proclaim Lemma 4.5 (Bessel functions of integer order and the Neumannn-Graf addition formula). 
Let $n\in{\Bbb Z}$, and take $J_n : {\Bbb C}\rightarrow{\Bbb C}$ to be the 
entire function which, for $z\neq 0$, satisfies $J_n(z)=(z/2)^{n} J_n^{*}(z)$, 
with $J_n^{*}(z)$ as defined in (1.9.6). Then, for all $z\in{\Bbb C}$, one has: 
$$J_n(z)
={1\over 2\pi}\int_{-\pi}^{\pi} \exp(-ni\Theta +iz\sin(\Theta)) {\rm d}\Theta
={1\over 2\pi}\int_{0}^{2\pi} \cos(z\sin(\Phi)-n\Phi) {\rm d}\Phi\;,\eqno(4.13)$$
$$\left| J_n(z)\right|
\leq\min\left\{ \exp\left(\left| {\rm Im}(z)\right|\right)\,,\,
{|z/2|^{|n|}\over |n|!}\,\exp\left( |z|\right)\right\}\eqno(4.14)$$
and, if  $|{\rm Re}(z)|\gg n^2+1$, then 
$$\left| J_n(z)\right|
\ll\exp\left(\left| {\rm Im}(z)\right|\right) \left| {\rm Re}(z)\right|^{-1/2}\;.\eqno(4.15)$$
If $n/2=p\in{\Bbb Z}$, $0\neq u\in{\Bbb C}$ and $e^{i\theta}=u/|u|$ then, for 
$y>0$ such that $y^2\neq e^{-2i\theta}$, one has 
$$(-1)^p  \left({y e^{i\theta}+(y e^{i\theta})^{-1}\over\left|
y e^{i\theta}+(y e^{i\theta})^{-1}\right|}\right)^{\!2p}\!\!J_{2p}\left( |u|\left|
y e^{i\theta}+(y e^{i\theta})^{-1}\right|\right)
=\sum_{m=-\infty}^{\infty} (-1)^m J_{m+p}(y|u|) J_{m-p}\!\left( {|u|\over y}\right) e^{2im\theta}\;,\eqno(4.16)$$
and the sum on the right-hand side of this equation is continuous, as a function of $y$, on 
the interval $(0,\infty)$. 

\medskip

\noindent{\bf Proof.}\quad 
The first equality in (4.13) is a result established in Section~17.23 of [48]; 
the subsequent equality is an elementary deduction (utilising the identity 
$2\cos(\alpha)=e^{i\alpha}+e^{-i\alpha}$, 
a substitution $\Phi =-\phi$, and the periodicity of the relevant integrands). 
The first bound implicit in (4.14) is a trivial corollary of (4.13), for one has 
$|\exp(-ni\Theta +iz\sin(\Theta))|=\exp(-{\rm Im}(z)\sin(\Theta))\leq\exp(|{\rm Im}(z)|)$ 
when $\Theta$ is real. 
Since (1.9.9), (1.9.8) and (1.9.6) imply the inequality 
$\bigl| J_n(z)\bigr|\leq |z/2|^{|n|}(|n|!)^{-1}J_0(i|z|)$, 
the second bound implicit in (4.14) follows from the first. 
\par 
The upper bound (4.15) is obtained by applying the first derivative test,  
Lemma~5.1.2 of [18], in order to estimate the first integral appearing in 
(4.13): after partitioning the range of integration $[-\pi ,\pi]$ 
into intervals within which the requisite 
monotonicity conditions are fulfilled, and having chosen some $\delta >0$, one proceeds to  
refine the partition into a disjoint union of subintervals on which one has 
$|\cos(\Theta)|> (\delta +|n|)/|{\rm Re}(z)|$, 
and subintervals on which $|\cos(\Theta)|\leq (\delta +|n|)/|{\rm Re}(z)|$; 
the refined partition need consist of no more than eight subintervals; 
those on which $\cos(\Theta)$ is bounded away from zero 
will (by the first derivative test) contribute no more than $(4/\pi)\delta^{-1}\exp(|{\rm Im}(z)|)$ 
to the absolute value of the first integral in (4.13); the remaining subintervals contribute (by virtue of 
the implied bound on their length) no more than $(\delta +|n|)/|{\rm Re}(z)|$; 
hence one obtains the bound 
$\bigl| J_n(z)\bigr|\leq\bigl( 4|{\rm Re}(\pi z) |^{-1/2}+|n||{\rm Re}(\pi z) |^{-1}\bigr)\exp(|{\rm Im}(z)|)$ 
by choosing $\delta =2|{\rm Re}(z/\pi)|^{1/2}$; when $|{\rm Re}(z)|\gg n^2+1$ 
this gives (4.15). 
\par
The identity (4.16) is an application of Graf's generalisation, 
in Section~2 of [13], of Neumann's addition formula. Graf's result, as 
presented in Section~11.3 of [44], is that 
$$J_{\nu}(\varpi)\left( {Z-z e^{-i\phi}\over Z-z e^{i\phi}}\right)^{\nu /2} 
=\sum_{k=-\infty}^{\infty} J_{\nu +k}(Z) J_k(z) e^{ik\phi}\;,\eqno(4.17)$$
where 
$$\varpi =\sqrt{ Z^2+ z^2 -2Z z\cos(\phi)}=\sqrt{\left( Z-z e^{-i\phi}\right)\left( Z-z e^{i\phi}\right)}\eqno(4.18)$$ 
and $\nu$, $Z$, $z$ and $\phi$ may take any complex values satisfying 
$$|z/Z|<\exp(-|{\rm Im}(\phi)|)\;.\eqno(4.19)$$
Suppose now that $n/2=p\in{\Bbb Z}$. Then, according to a remark 
on Page~361 of [44], the case $\nu =2p$ of (4.17)-(4.18) 
holds independently of (4.19). Since that remark is made without proof, 
it is worthwhile to verify it here. In doing so one is to assume that 
$\nu =2p\in 2{\Bbb Z}$, and that $Z$ and $\phi$ are given complex numbers 
with $Z\neq 0$. Then, by (1.9.9), (1.9.8) and (1.9.6), the complex function 
$w\mapsto J_{\nu}(w)$ is entire, even, and has a zero of order $|\nu|=2|p|$ 
at $w=0$. Consequently, given (4.18), the left-hand side of (4.17) is a function of  
$z$ that is meromorphic on ${\Bbb C}$, and has no singularities other than 
(possibly) a removable singularity at one of the points $z=Ze^{\pm i\phi}$. 
Moreover, since $z\mapsto J_k(z)$ is an entire function when $k\in{\Bbb Z}$, 
and since the bound (4.14) shows the sum over $k$ in (4.17) to be 
uniformly convergent for the values of $z$ lying in any given compact subset of ${\Bbb C}$, 
it follows that the right-hand side of (4.17) is an entire function of the complex variable $z$. 
Since both sides of (4.17) are (apart from at most one removable singularity) entire functions 
of $z$, and since they are equal when $z$ lies in the open disc where 
(4.19) is satisfied,  it follows that (4.17) must hold for all $z\in{\Bbb C}$ 
(subject to suitable definition of the left-hand side at any removable singularity). 
This justifies the use of (4.17)-(4.18) whenever $\nu$ is an even integer, and 
$Z$, $z$ and $\phi$ are complex numbers with $Z\neq 0$. 
\par
To deduce (4.16) one applies  the identity (4.17)-(4.18)
with $\nu =n=2p$, $Z=|u| y$, $z=-|u|/y$ and $\phi =2\theta$ 
(where, by hypothesis, $p\in{\Bbb Z}$, $0\neq u\in{\Bbb C}$, $e^{i\theta}=u/|u|$ and $y>0$\/). 
Then, in (4.17) and (4.18), one has: 
$$\varpi^2 =|u|^2\left( y^2+y^{-2}+2\cos(2\theta)\right) 
=|u|^2\left|\,y e^{i\theta}+( y e^{i\theta})^{-1}\right|^2\;,$$
$$\left( {Z-z e^{-i\phi}\over Z-z e^{i\phi}}\right)^{\nu /2}
=\left({y+y^{-1} e^{-2i\theta}\over 
y+y^{-1} e^{2i\theta}}\right)^{\!p}
=\left({y e^{i\theta}+y^{-1} e^{-i\theta}\over 
y e^{-i\theta}+y^{-1} e^{i\theta}}\right)^{\!p} e^{-2i\theta p}$$
and (using (1.9.9))
$$\eqalign{
\sum_{k=-\infty}^{\infty} J_{\nu +k}(Z) J_k(z) e^{ik\phi}
 &=\sum_{k=-\infty}^{\infty} J_{k+2p}(y|u|) J_k(-|u|/y) e^{2im\theta} =\cr 
 &=\sum_{m=-\infty}^{\infty}J_{m+p}(y|u|) (-1)^{m-p} J_{m-p}(|u|/y) 
e^{2i(m-p)\theta}\;;}$$
since $J_{2p}(\varpi)$ is an even function of $\varpi$, and since 
$(\alpha /\overline{\alpha})^p=(\alpha /|\alpha|)^{2p}$, one therefore obtains the 
desired result (4.16) on multiplying by $(-1)^p e^{2ip\theta}$ both sides of the 
specified case of equation (4.17). Just as the right-hand side of (4.17) was found 
to be entire, as a function of $z$, so it can also be established (by similar reasoning) that, 
for arbitrary fixed $p\in{\Bbb Z}$, $u\in{\Bbb C}-\{ 0\}$ and $\theta \in{\Bbb R}$, 
the right-hand side of (4.16)  will be a function of $y$ that is holomorphic on ${\Bbb C}-\{ 0\}$, 
and so will (in particular) be continuous for $y>0$ 
\quad $\blacksquare$

\bigskip

\proclaim Lemma 4.6 (the ${\bf B}$-transform of a certain test function). 
Let $\sigma\in(1/2,1)$ and $K,P\geq 1\,$; let   
$h$ be the function on 
${\cal S}^{\star}_{\sigma}=\left\{\nu\in{\Bbb C} : |{\rm Re}(\nu)|\leq\sigma\right\}\times{\Bbb Z}$ 
given by 
$$h(\nu,p)=\exp\left( (\nu/K)^2-(p/P)^2\right)\quad\hbox{for $(\nu,p)\in{\cal S}^{\star}_{\sigma}\,$.}\eqno(4.20)$$
Then there exist real numbers $\varrho,\vartheta >3$ such that 
the conditions (i)-(iii) of Theorem~B are satisfied. \hfill\break 
$\hbox{\qquad}$Suppose moreover that ${\Bbb C}\ni u=|u| e^{i\theta}\neq 0$,   
that one has $1\leq\Delta\in{\Bbb R}$, $M\in{\Bbb N}$ and  
$$\Delta\leq {M\over 1+|u|}\leq 2\Delta\;,\eqno(4.21)$$
and that ${\bf B}h$ is the transform of $h$ defined by (1.9.3)-(1.9.4) of Theorem B. 
Then
$$({\bf B}h)(u)
={1\over 4\pi^{3}}\int\limits_{-\pi}^{\pi}\int\limits_{-\infty}^{\infty}
\int\limits_{-\infty}^{\infty}\!\!F_{P,K}(\eta ,\xi) e^{-\xi^2-\eta^2}
A_M(\phi,\theta) \cos( |u|\psi(\xi /K,\eta /P;\phi)){\rm d}\eta {\rm d}\xi {\rm d}\phi
+E_M(P,K;u)\;,\eqno(4.22)$$
where
$$F_{P,K}(\eta ,\xi) =
\left( {1\over 2}-\eta^2\right)P^2 +\left( {1\over 2} -\xi^2\right) K^2\;,\qquad 
\psi(y,x;\phi)=e^y\sin(\phi -x)+e^{-y}\sin(\phi +x)\;,\eqno(4.23)$$ 
$$A_M(\phi,\theta) =\sum_{m=-M}^M (-1)^m \cos(2m\phi) e^{2im\theta} \eqno(4.24)$$
and 
$$E_M(P,K;u)\ll_j \left( P^2+K^2\right) (1+|u|)\Delta^{1-2j}\qquad 
\hbox{($j\in{\Bbb N}$\/).}\eqno(4.25)$$
At the same time, one has also  
$$({\bf B}h)(u)
={|u|^2\over 8\pi^{3}}\int\limits_{-\pi}^{\pi}\int\limits_{-\infty}^{\infty}
\int\limits_{-\infty}^{\infty}\!\!G_{P,K}(\eta ,\xi) e^{-\xi^2-\eta^2}
A_M(\phi,\theta) \cos( |u|\psi(\xi /K,\eta /P;\phi)){\rm d}\eta {\rm d}\xi {\rm d}\phi 
+E_M(P,K;u)\;,\eqno(4.26)$$
where 
$$G_{P,K}(\eta ,\xi) = \cosh(2\xi /K)-\cos( 2\eta /P) =2\left( \sinh^2(\xi /K)+\sin^2(\eta /P)\right)
\eqno(4.27)$$
(while $\psi(x,y;\phi)$, $A_M(\phi,\theta)$ and $E_M(P,K;u)$ remain as in 
(4.22)-(4.25)\/).

\medskip

\noindent{\bf Proof.}\quad 
By (4.20) it is evident that $h(-\nu ,p)=h(\nu ,-p)=h(\nu ,p)$, so that condition (i) of Theorem B is satisfied. 
Since the function $\nu\mapsto \exp\left( \nu^2 /K^2\right)$ is entire, the definition (4.20) also 
ensures that $h(\nu,p)$ satisfies condition (ii) of Theorem~B. One can moreover
check that $h(\nu ,p)$ satisfies condition (iii) of Theorem~B for abitrary $\varrho,\vartheta >3$ 
(a short calculation using the  inequality $\exp\left( -x^2\right)\leq 1/\left( 1+x^4 /2\right)$
gives, in particular, the case $\varrho =\vartheta =4$ of the condition (iii), with an implicit constant 
not greater than $O\left( K^4 P^4\right)=O_h(1)$\/).
\par 
Suppose now that $0\neq u\in{\Bbb C}$, and that $\theta\in{\Bbb R}$ satisfies 
$e^{i\theta}=u/|u|$. Since the definitions (1.9.4)-(1.9.6) of ${\cal K}_{\nu ,p}(z)$ 
are equivalent to the equations (6.21) and (7.21) of [5], 
it follows by Bruggeman and Motohashi's identity 
(1.9.11) and the definition (1.9.3) of $({\bf B}h)(z)$ that, for 
the chosen test function $h$ (i.e. that in (4.20)), one has 
$$({\bf B}h)(u)={2\over\pi}\sum_{p\in{\Bbb Z}}{(-1)^p\over e^{(p/P)^2}}
\int\limits_0^{\infty}  \left({y e^{i\theta}+(y e^{i\theta})^{-1}\over\left|
y e^{i\theta}+(y e^{i\theta})^{-1}\right|}\right)^{2p} J_{2p}\left( |u|\left|
y e^{i\theta}+(y e^{i\theta})^{-1}\right|\right)f_p(y)\,{{\rm d}y\over y}\;,\eqno(4.28)$$ 
where, by an appeal to Lemma 4.4, 
$$\eqalign{
f_p(y) ={1\over 4\pi i}\int\limits_{(0)}y^{2\nu}
e^{(\nu /K)^2}\left( p^2-\nu^2\right){\rm d}\nu 
 &={K\over 4\pi}\left( p^2 G_0(K\log y)+K^2 G_2(K\log y)\right) =\cr 
&={K\over 4\pi}\left( p^2 +K^2 \left( {1\over 2}-(K\log y)^2\right) \right) G_0(K\log y) =\cr 
 &={K\over 4\sqrt{\pi}}\left( p^2+{K^2\over 2}-K^4\log^2 y\right) e^{-K^2\log^2y}}\eqno(4.29)$$
(since $\left| y^{2\nu}\right|=1$ for $y>0$ and ${\rm Re}(\nu)=0$, the change in the 
order of integration required to attain (4.28) is justified by the absolute convergence 
of the integrals in (4.29) and (1.9.11): see Theorem~10.40 of [1] and note 
that, by the bound (4.15) of Lemma~4.5, one has 
$J_{2p}\left( |u| \bigl| y e^{i\theta}+(y e^{i\theta})^{-1}\bigr|\right)\ll_{p,|u|} ( y+y^{-1})^{-1/2}$ 
for~$y>0$). 
\par
By the result (4.16) of Lemma 4.5 (i.e. the Neumann-Graf addition law), the integrand 
in (4.28) is equal to 
$$(-1)^p y^{-1} f_p(y)\sum_{m=-\infty}^{\infty}(-1)^m J_{m+p}(y|u|)J_{m-p}(y^{-1}|u|)e^{2im\theta}\;,$$ 
and so, by substituting $-p$ for $p$, and $e^{\xi /K}$ for $y$,
one may rewrite (4.28) as: 
$$({\bf B}h)(u)={2\over\pi K}\sum_{p\in{\Bbb Z}} e^{-(p/P)^2} 
\int\limits_{-\infty}^{\infty}f_p\left( e^{\xi /K}\right)
\sum_{m\in{\Bbb Z}}J_{m-p}\left( e^{\xi /K}|u|\right)
J_{m+p}\left( e^{-\xi /K}|u|\right) \left({iu\over |u|}\right)^{\!\!2m} {\rm d}\xi\eqno(4.30)$$
(\/$f_p(y)$ being, by (4.29), an even function of $p$\/). Then, by using the latter of the two 
integral representations of $J_n(z)$ appearing in the equations (4.13) of Lemma 4.5,  
one obtains (via some consideration of the periodicity of the relevant integrands) the result that  
$$\eqalign{J_{m-p}\left( e^{\xi /K}|u|\right) & J_{m+p}\left( e^{-\xi /K}|u|\right) =\cr 
 &={1\over 8\pi^2}\qquad\int\!\!\!\!\!\!\!\!\!\!\!\!\int\limits_{\!\!\!\!\!\!\!\!\!\!\!\!\scriptstyle
-2\pi\leq\Psi-\Phi\leq 2\pi\atop\scriptstyle 0\leq\Psi+\Phi\leq
4\pi}^{\matrix{\ }}\!\!\!\!\!\!\!\!\cos\!\left( e^{\xi /K}|u|\sin\Phi -(m-p)\Phi\right)\!\cos\!\left(e^{-\xi
/K}|u|\sin\Psi -(m+p)\Psi\right)\! {\rm d}\Phi{\rm d}\Psi\;.}$$
On substituting 
$\delta =(\Psi - \Phi)/2$ and $\phi=(\Psi+\Phi)/2$, the integrand becomes a function 
of $\phi$ with period $2\pi$. Hence one may take 
$[-\pi ,\pi]$ as the range of integration for both $\delta$ and $\phi$.
Then, on using the identity $\cos(x)\cos(y)=\bigl(\cos(x+y)+\cos(x-y)\bigr) /2$ to express the 
integral as a sum of two integrals, one finds (by interchange of $\delta$ and $\phi$, and 
the identity $\cos(-x)=\cos(x)$) that both integrals are equal. The outcome of
these manipulations  is that one obtains: 
$$J_{m-p}\left( e^{\xi /K}|u|\right) J_{m+p}\left( e^{-\xi /K}|u|\right)
={1\over 4\pi^2}\int\limits_{-\pi}^{\pi}\int\limits_{-\pi}^{\pi}
\cos( |u|\psi(\xi /K,\delta;\phi) -2m\phi-2p\delta){\rm d}\delta{\rm d}\phi\;,\eqno(4.31)$$
where $\psi(y,x;\phi)$ is as given by (4.23). The inner range of integration here 
may be changed to $[-\pi /2,3\pi /2]$ (by periodicity of the integrand as a function of $\delta$\/). 
Then, on writing the integral concerned as the sum of integrals 
over the subintervals $[-\pi /2,\pi /2]$ and $[\pi /2,3\pi /2]$ (respectively), 
and substituting $\delta +\pi$ for $\delta$ in the latter integral, one finds 
that elementary trigonometric identities suffice to express the sum of the 
integrals over the two subintervals as the single integral 
$$\int\limits_{-\pi /2}^{\pi /2}
2\cos( |u|\psi(\xi /K,\delta;\phi))\cos(2m\phi+2p\delta){\rm d}\delta\;.$$ 
Since $\cos(2m\phi+2p\delta)=\cos(2m\phi)\cos(2p\delta)-\sin(2m\phi)\sin(2p\delta)$, 
one may therefore reformulate (4.31) as:  
$$J_{m-p}\left( e^{\xi /K}|u|\right) J_{m+p}\left( e^{-\xi /K}|u|\right) 
={1\over\pi}\int\limits_{-\pi /2}^{\pi /2}\left(\cos(2p\delta)c_{2m}(|u|,\xi /K,\delta) 
-\sin(2p\delta)s_{2m}(|u|,\xi /K,\delta)\right) {\rm d}\delta\;,\eqno(4.32)$$
where
$$c_{2m}(r,y,x)={1\over 2\pi}\int\limits_{-\pi}^{\pi}\cos( r\psi(y,x;\phi))\cos(2m\phi){\rm
d}\phi\;,\qquad
s_{2m}(r,y,x)={1\over 2\pi}\int\limits_{-\pi}^{\pi}\cos( r\psi(y,x;\phi))\sin(2m\phi){\rm
d}\phi\;.$$
\par
It is trivial that $\left| c_{2m}(r,y,x)\right|$ and $\left| s_{2m}(r,y,x)\right|$ are each 
bounded above by $1$. Moreover, 
since the relevant functions have period $2\pi$, two integrations by parts 
suffice in order to show that 
$$c_{2m}(r,y,x)
=-{1\over 2\pi}\int\limits_{-\pi}^{\pi}
\left( {\partial^2\over\partial\phi^2}\,\cos( r\psi(y,x;\phi))\right) {\cos(2m\phi)\over 4 m^2}\,
{\rm d}\phi\;,$$
with a similar result holding in respect of $s_{2m}(r,y,x)$. By these results, along 
with the observations that, for $x,y,\phi\in{\Bbb R}$, the function $\psi$ given by (4.23) is real-valued, and satisfies  
$\bigl| (\partial/\partial\phi)^j\,\psi(y,x;\phi)\bigr|\leq e^y+e^{-y}$ for $j=1,2$, 
one has:  
$$\left| c_{2m}(r,y,x)\right|\,,\ \left| s_{2m}(r,y,x)\right| 
\ll m^{-2}(1+r)r\cosh(2y)\qquad\hbox{($0\neq m\in{\Bbb Z}$, $x,y\in{\Bbb R}$, $r>0$).}$$
Therefore, and by (4.32) and (4.29), the
integral and innermost sum in (4.30) are uniformly convergent with respect to $p$.  
This justifies a change in the order of summation and integration, 
so that one may sum firstly over $p$, in (4.30), before going on to sum over $m$ and integrate 
with respect to $\xi$. Terms involving the factor $\sin(2p\delta)$ (from (4.32)) cancel, 
since all other relevant factors are even functions of $p$. Hence (and by (4.29)) one obtains: 
$$({\bf B}h)(u)=\pi^{-3/2}\int\limits_{-\infty}^{\infty}e^{-\xi^2}\sum_{m\in{\Bbb Z}}
\left({iu\over |u|}\right)^{2m}T_m(|u|,\xi ;P,K){\rm d}\xi\;,\eqno(4.33)$$ where
$$T_m(r,\xi ;P,K)={1\over 2\pi}\int\limits_{-\pi /2}^{\pi /2}
c_{2m}(r,\xi /K,\delta)\left(\sum_{p\in{\Bbb Z}}\cos(2p\delta)\left( p^2+\left({1\over
2}-\xi^2\right) K^2\right) e^{-(p/P)^2}\right){\rm d}\delta .\eqno(4.34)$$ 
\par 
One now repeats a calculation from the proof of Corollary~10.1 of [5], by 
applying Poisson summation (i.e. the case $n=1$ of the equation (2.45) of Lemma 2.7)
to the sum over $p$ in the last equation;   
then, after evaluating the relevant Fourier integrals (very similar to the integral appearing in (4.29))  
by means of Lemma 4.4, one finds that the sum over $p$ in (4.34) equals  
$$\sqrt{\pi}P\sum_{v\in{\Bbb Z}}\left(\left( {1\over 2}-P^2(\pi v+\delta)^2\right)P^2
+\left( {1\over 2} -\xi^2\right) K^2\right) e^{-P^2(\pi v+\delta)^2}\;.$$
Therefore, 
and given that $c_{2m}(r,\xi /K,\delta) =c_{2m}(r,\xi /K,v\pi +\delta)$ for $v\in{\Bbb Z}$ 
(as follows by (4.23) and the definition of $c_{2m}(r,y,x)$\/), one finds that  
$$T_m(r,\xi ;P,K)={P\over 2\sqrt{\pi}}\sum_{v\in{\Bbb Z}}\ \int\limits_{-\pi /2}^{\pi /2}
g(\pi v +\delta){\rm d}\delta = {P\over 2\sqrt{\pi}}\int\limits_{-\infty}^{\infty} g(x){\rm d}x,$$
where $$g(x)=\left(\left( {1\over 2}-P^2 x^2\right)P^2 +\left( {1\over 2} -\xi^2\right) K^2\right)
e^{-P^2 x^2}c_{2m}(r,\xi /K,x)\;.$$ 
Recalling now the definition of $c_{2m}(r,y,x)$, 
and substituting $\eta /P$ for $x$ (in the above),  one obtains: 
$$T_m(r,\xi ;P,K)={1\over 4\pi^{3/2}}\int\limits_{-\infty}^{\infty}F_{P,K}(\eta ,\xi)
 e^{-\eta^2}
\int\limits_{-\pi}^{\pi}\cos( r\psi(\xi /K,\eta /P;\phi))\cos(2m\phi){\rm d}\phi{\rm d}\eta\;,
\eqno(4.35)$$
where $F_{P,K}(\eta ,\xi)$ is as given by (4.23).
\par
For $j\in{\Bbb N}$, one can show (using $2j$ integrations-by-parts) that 
$$\int\limits_{-\pi}^{\pi}\cos( r\psi(y,x;\phi))\cos(2m\phi){\rm d}\phi
\ll_j m^{-2j}\left( 1+r^{2j}\cosh(2jy)\right)\qquad 
\hbox{($0\neq m\in{\Bbb Z}$, $x,y\in{\Bbb R}$).}$$
By these bounds, and (4.35) and (4.23), it follows that if $j\in{\Bbb N}$ then 
$$T_m(r,\xi ;P,K)
\ll_j m^{-2j}\left( P^2 +\left(1+\xi^2\right)K^2\right) 
\left( 1+r^{2j}\cosh(2j\xi /K)\right)\qquad 
\hbox{($0\neq m\in{\Bbb Z}$,  $\xi\in{\Bbb R}$, $r>0$),}$$
so that, since $K\geq 1$, one obtains: 
$$\sum_{\scriptstyle m\in{\Bbb Z}\atop\scriptstyle |m|>M}
\left|T_m(r,\xi ;P,K)\right|
\ll_j (M+1)^{1-2j}\left( P^2 +\left(1+\xi^2\right)K^2\right) 
\left( 1+r^{2j}\cosh(2j\xi)\right)\qquad 
\hbox{($M\geq 0$, $\xi\in{\Bbb R}$, $r>0$).}$$

Suppose now that $\Delta\geq 1$, and that $M\in{\Bbb N}$ 
satisfies the conditions in (4.21). Then it follows, by (4.33) and the last bound obtained, that 
for $j\in{\Bbb N}$ one has 
$$\eqalign{
({\bf B}h)(u)-{1\over\pi^{3/2}}\int\limits_{-\infty}^{\infty}e^{-\xi^2} \sum_{m=-M}^M  
\left({iu\over |u|}\right)^{2m}T_m(|u|,\xi ;P,K){\rm d}\xi
&=O_j\left( M^{1-2j}\left( P^2+K^2\right)\left( 1+|u|^{2j}\right)\right) \ll_j\cr 
&\ll_j\left( P^2+K^2\right) ( 1+|u|)\Delta^{1-2j}\;.}$$
This, together with (4.35) (and justifiable changes in the order of summation and integration),
yields the first main estimate of the lemma, as expressed in (4.22)-(4.25).
\par
It now only remains to show how (4.26)-(4.27) follow.   
Let $I(P',K')$ denote the triple integral got by substituting $F_{P',K'}(\eta ,\xi)$ in 
place of the factor $F_{P,K}(\eta ,\xi)$  in the integrand in equation (4.22). 
Then, by (4.23), the integral in (4.22) may be expressed as the sum $I(P,0)+I(0,K)$.
To transform $I(P,0)$ one integrates by parts twice (with respect to $\eta$), using: 
$$\int (1/2-\eta^2)e^{-\eta^2}{\rm d}\eta
=(1/2)\eta e^{-\eta^2}+C_1\quad\hbox{and}\quad
\int \eta e^{-\eta^2}{\rm d}\eta =-(1/2)e^{-\eta^2}+C_2\;.$$
After transforming $I(0,K)$ similarly (through integrations by parts with respect to $\xi$), 
one observes a certain obvious cancellation between terms of the 
expression obtained for $I(P,0)$ and terms of the expression obtained for $I(0,K)$: it then
requires only a few more steps
(involving use of elementary trigonometric identities) to
arrive at the result expressed in (4.26)-(4.27)\quad$\blacksquare$

\medskip\par
\noindent{\bf Remark~4.7.}\quad 
The definition (4.24) expresses the function $A_M(\phi,\theta)$ in the form 
most convenient for use in this paper. It may nevertheless be worth noting 
that, by (4.24), one in fact has: 
$$A_M(\phi,\theta)
=\sum_{m=-M}^M (-1)^m\cos(2m\phi)\cos(2m\theta)
={(-1)^M\over 2}\left( 
{\cos((2M+1)(\phi+\theta))\over\cos(\phi+\theta)}
+{\cos((2M+1)(\phi-\theta))\over\cos(\phi-\theta)}\right) ,$$
for $M\in{\Bbb N}\cup\{ 0\}$ and $\phi,theta\in{\Bbb R}$ 
such that $\cos(\phi+\theta)\cos(\phi -\theta)\neq 0$. 

\bigskip

\proclaim Lemma 4.8. 
Let $-\pi /2<x<\pi /2$; and let $y\in{\Bbb R}$. Suppose that 
$$J=\int\limits_{-\pi}^{\pi}\left|\psi(y,x;\phi)\right|^{-1/2}{\rm d}\phi\;,$$
where $\psi(y,x;\phi)$ is as defined in (4.23). 
Then $J\ll (\cos(x))^{-1/2}$.

\medskip 

\noindent{\bf Proof.}\quad 
Using $\sin(A\pm B)=\sin(A)\cos(B)\pm\cos(A)\sin(B)$, one  finds that
$$\psi(y,x;\phi)=2\left(\cosh(y)\cos(x)\sin(\phi)-\sinh(y)\sin(x)\cos(\phi)\right) 
=2\left( A\cos(\phi) +B\sin(\phi)\right) ,$$
where $A=-\sinh(y)\sin(x)$ and $B=\cosh(y)\cos(x)$. Therefore, on putting
$$Z=Z(x,y)=\sqrt{A^2+B^2}=\sqrt{\cosh^2 (y)-\sin^2 (x)},$$
one has $A/Z=\cos(\alpha)$ and $B/Z=\sin(\alpha)$ for some $\alpha =\alpha(x,y)$ (obviously
independent of $\phi$), and so
$$\psi(y,x;\phi)=2Z(x,y)\cos\left(\phi -\alpha(x,y)\right)\qquad\hbox{($\phi\in{\Bbb R}$).}$$
Using this, one obtains: 
$$J\ll Z(x,y)^{-1/2}\int\limits_{-\pi}^{\pi}\left|\cos\left( \phi -\alpha(x,y)\right)\right|^{-1/2}
{\rm d}\phi 
=4Z(x,y)^{-1/2}\int\limits_{0}^{\pi /2}\left(\sin(\theta)\right)^{-1/2}{\rm d}\theta 
\ll Z(x,y)^{-1/2}\;.$$ 
The lemma follows, since $Z^2=\cosh^2 (y)-\sin^2 (x)\geq 1-\sin^2 (x) 
=\cos^2 (x)\ \,\blacksquare$

\bigskip

\centerline{\bf \S 5. The proof of Theorem 1.}

\medskip

Let the hypotheses of the theorem hold; 
let ${\frak a}$ be a cusp of $\Gamma$ with 
${\frak a}\sim^{\!\!\!\!\Gamma} u/w$; and let 
$E_j^{\frak a}(q_0 ,P,K;N,b)$ ($j=0,1$) and $\mu({\frak a})$ be as 
given by (1.9.12), (1.9.13) and (1.9.15). Suppose also that $1/2<\sigma <1$, and that   
$$h(\nu ,p)=\exp\left( \left( {\nu\over K}\right)^2 -\left( {p\over P}\right)^2\right)\qquad\quad  
\hbox{for $p\in{\Bbb Z}$ and $\nu\in{\Bbb C}$ with $|{\rm Re}(\nu)|\leq\sigma\,$}\eqno(5.1)$$
(so that $h(\nu ,p)$ is the function $h : {\cal S}^{\star}_{\sigma}\rightarrow{\Bbb C}$ defined in 
(4.20) of Lemma 4.6). 
Then  $h(\nu ,p)\geq\exp(-2)$ for all pairs 
$(\nu ,p)\in ((i{\Bbb R})\cup [-\sigma ,\sigma])\times{\Bbb Z}$ 
such that $|\nu|\leq K$ and $|p|\leq P$. 
Moreover, each irreducible subspace $V\subset{}^0L^2(\Gamma\backslash G)$   
that is the index of a summand on the right-hand side of Equation~(1.9.12) 
is also the index of a summand on the left-hand side of Equation~(1.9.1), 
and so it may be inferred from Remark~1.9.2, below Theorem~B, 
that each pair of spectral parameters 
$\nu_V, p_V$ associated with a summand on the right-hand side of Equation~(1.9.12) 
is such that the condition~(1.9.7) is satisfied.    
Therefore, with $h(\nu ,p)$ as in (5.1) (where, by hypothesis, $\sigma>1/2>2/9$),  
it follows from the definitions (1.9.12) and (1.9.13) that one has 
$$0\leq
E_j^{\frak a}\!\left( q_0 ,P,K;N,b\right)
\leq\sum_{k=0}^1\,E_k^{\frak a}\!\left( q_0 ,P,K;N,b\right) 
\leq e^2{\cal E}^{\frak a}\!\left( q_0 ,P,K;N,b\right)\quad\ \hbox{for $j=0,1\,$,}\eqno(5.2)$$
where 
$$\eqalign{\!\!\!{\cal E}^{\frak a}\!\left( q_0 ,P,K;N,b\right)  
 &=\sum_V\,\Biggl|\sum_{\scriptstyle\omega\in{\frak O}\atop\scriptstyle N/2<|\omega|^2\leq N}
b(\omega) C_V^{\frak a}\left(\omega;\nu_V,p_V\right)\Biggr|^2 
\!\exp\!\left(\!{\nu_V^2\over K^2} - {p_V^2\over P^2}\!\right)\ +\cr 
 &\quad\ +\sum_{{\frak c}\in{\frak C}(\Gamma)}
{1\over 4\pi i\left[\Gamma_{\frak c} : \Gamma_{\frak c}'\right]} 
\sum_{p\in{1\over 2}\left[\Gamma_{\frak c} : \Gamma_{\frak c}'\right]{\Bbb Z}}
\ \,\int\limits_{(0)}\,\Biggl|\sum_{\scriptstyle\omega\in{\frak O}\atop\scriptstyle N/2<|\omega|^2\leq N}
\!\!\!b(\omega) B_{\frak c}^{\frak a}\!\left(\omega;\nu,p\right)\Biggr|^2 
\!\exp\!\left(\!{\nu^2\over K^2} - {p^2\over P^2}\!\right)\!{\rm d}\nu\;. 
}\eqno(5.3)$$

After using the identity $|S|^2 =\overline{S}\,S$ to expand the squared absolute values 
in (5.3), one may change the order of summation so as to rewrite the right-hand side of 
equation (5.3) in the form 
$$\sum\!\!\!\!\!\!\!\!\!\!\sum_{\!\!\!\!\!\!\!\!\!\!\!\!\!{\scriptstyle\omega_1,\omega_2\in{\frak O}\atop\scriptstyle 
N/2<\left|\omega_1\right|^2,\left|\omega_2\right|^2\leq N}}  
\!\!\!\overline{b\left(\omega_1\right) }\,b\left(\omega_2\right) 
H_{{\frak a},{\frak a}}\left(\omega_1,\omega_2\right) ,$$
where $H_{{\frak a},{\frak b}}\left(\omega_1,\omega_2\right)$ denotes the 
left-hand side of the special case of equation (1.9.1) in which $h$ is the function 
defined in (5.1). 
By Lemma~4.6 it follows that, when $\sigma$ and 
the function $h$ are as we suppose in~(5.1), 
there exists some pair of real numbers  
$\varrho,\vartheta >3$ 
such that the conditions (i)-(iii) of Theorem~B 
are satisfied (the proof of Lemma~4.6 shows, specifically, 
that one may take $\varrho =\vartheta =4$ here).  
Hence one may apply the spectral sum formula (Theorem B) to each 
term  $H_{{\frak a},{\frak a}}\left(\omega_1,\omega_2\right)$ in the above sum, 
and so obtain: 
$${\cal E}^{\frak a}\!\left( q_0 ,P,K;N,b\right) 
={\cal D}^{\frak a}\!\left( q_0 ,P,K;N,b\right) 
+{\cal L}^{\frak a}\!\left( q_0 ,P,K;N,b\right)\;,\eqno(5.4)$$
where 
$${\cal D}^{\frak a}\!\left( q_0 ,P,K;N,b\right) 
={1\over 4\pi^3 i}\Biggl( 
\ \sum_{p\in{\Bbb Z}}\ \int\limits_{(0)} h(\nu,p) \left( p^2 -\nu^2\right) {\rm d}\nu\Biggr)\quad\   
\sum\!\!\!\!\!\!\!\!\!\!\sum_{\!\!\!\!\!\!\!\!\!\!\!\!\!{\scriptstyle\omega_1,\omega_2\in{\frak O}\atop\scriptstyle 
N/2<\left|\omega_1\right|^2,\left|\omega_2\right|^2\leq N}}  
\!\!\!\delta^{{\frak a},{\frak a}}_{\omega_1,\omega_2}\,\overline{b\left(\omega_1\right) }
\,b\left(\omega_2\right)\eqno(5.5)$$
and 
$${\cal L}^{\frak a}\!\left( q_0 ,P,K;N,b\right) 
=\qquad\sum\!\!\!\!\!\!\!\!\!\!\sum_{\!\!\!\!\!\!\!\!\!\!\!\!\!{\scriptstyle\omega_1,\omega_2\in{\frak O}\atop\scriptstyle 
N/2<\left|\omega_1\right|^2,\left|\omega_2\right|^2\leq N}}\!\!\!\overline{b\left(\omega_1\right) }
\,b\left(\omega_2\right) 
\,\sum_{c\in {}^{\frak a}{\cal C}^{\frak a}} 
\,{S_{{\frak a},{\frak a}}\left(\omega_1 , \omega_2 ; c\right)\over |c|^2}\,
({\bf B}h)\!\!\left( {2\pi\sqrt{\omega_1\omega_2}\over c}\right) ,\eqno(5.6)$$
with $h$ as in (5.1), and with $\delta^{{\frak a},{\frak a}}_{\omega_1,\omega_2}$ 
and $({\bf B}h)(z)$ as given by 
(1.9.2) and  (1.9.3)-(1.9.6). 
\par
By (1.9.2) one has, in equation (5.5), 
$$\delta^{{\frak a},{\frak a}}_{\omega_1,\omega_2}
=\sum_{\scriptstyle \gamma\in\Gamma_{\frak a}'\backslash\Gamma_{\frak a}\atop\scriptstyle 
g_{\frak a}^{-1} \gamma g_{\frak a}
=\left({\scriptstyle\!\!\!\!\!u(\gamma)\quad\beta(\gamma)\atop\scriptstyle\  0\quad\ 1/u(\gamma)}\right)} 
{\rm e}\left( {\rm Re}\left( \beta(\gamma)u(\gamma)\omega_1\right)\right) 
\delta_{u(\gamma)\omega_1 , \omega_2/u(\gamma)}\;.$$
For a more explicit representation of the last sum note firstly that 
the scaling matrix $g_{\frak a}$ is, by hypothesis, such that 
that (1.1.16) and (1.1.20)-(1.1.21) hold for ${\frak c}={\frak a}$. 
Therefore $\Gamma_{\frak a}'=g_{\frak a} B^{+} g_{\frak a}^{-1}$, 
where $B^{+}=\{ n[\alpha] : \alpha\in{\frak O}\}$. 
This, when combined with the results of Lemma 4.2 concerning 
$g_{\frak a}^{-1}\Gamma_{\frak a} g_{\frak a}$, 
shows that there exists a complex number $\beta_{\frak a}$ such that 
the set ${\cal T}_{\frak a}$ given by
$${\cal T}_{\frak a}
=\cases{ g_{\frak a}\bigl\{ h[1] , h[-1] , 
h[i]\,n\!\left[\beta_{\frak a}\right] , h[-i]\,n\!\left[\beta_{\frak a}\right]\bigr\} g_{\frak a}^{-1}
 &if $q_0 \mu({\frak a})\mid 2$, \cr 
 \quad &\hbox{\quad} \cr
g_{\frak a}\left\{ h[1] , h[-1]\right\} g_{\frak a}^{-1} &otherwise,}\eqno(5.7)$$
is a complete set of coset representatives of $\Gamma_{\frak a}'$ in $\Gamma_{\frak a}\,$: 
as noted below (1.1.19), the group $\Gamma_{\frak a}'$ is a normal subgroup of $\Gamma_{\frak a}$ 
(one has, in particular, 
$h[i]\,n\!\left[\beta_{\frak a}\right] B^{+}= B^{+} h[i]\,n\!\left[\beta_{\frak a}\right]$\/).
Consequently one finds that if $q_0 \mu({\frak a})\mid 2$ then
$$\delta^{{\frak a},{\frak a}}_{\omega_1,\omega_2}=2\delta_{\omega_1,\omega_2}
+2{\rm e}\left( {\rm Re}\left( i^2\beta_{\frak a}\omega_1\right)\right) \delta_{i\omega_1,\omega_2 /i}
=2\bigl(\delta_{\omega_1,\omega_2}
+{\rm e}\left( -{\rm Re}\left( \beta_{\frak a}\omega_1\right)\right) 
\delta_{-\omega_1,\omega_2}\bigr)\;,$$ 
while if instead $q_0 \mu({\frak a})\!\not\,\mid 2$ then one has simply 
$\delta^{{\frak a},{\frak a}}_{\omega_1,\omega_2}=2\delta_{\omega_1,\omega_2}$. 
Therefore, and since (as (5.7) shows)  
$$\left[\Gamma_{\frak a} : \Gamma_{\frak a}'\right]
=\cases{4 &if $q_0 \mu({\frak a})\mid 2\,$,\cr 2 &otherwise,}\eqno(5.8)$$
it follows that one has: 
$$\qquad\qquad\quad       
\sum\!\!\!\!\!\!\!\!\!\!\sum_{\!\!\!\!\!\!\!\!\!\!\!\!\!{\scriptstyle\omega_1,\omega_2\in{\frak O}\atop\scriptstyle 
N/2<\left|\omega_1\right|^2,\left|\omega_2\right|^2\leq N}}  
\!\!\!\delta^{{\frak a},{\frak a}}_{\omega_1,\omega_2}\,\overline{b\left(\omega_1\right) }
\,b\left(\omega_2\right)
=\left[\Gamma_{\frak a} : \Gamma_{\frak a}'\right]
\!\!\!\sum_{\scriptstyle\omega\in{\frak O}\atop\scriptstyle N/2<|\omega|^2\leq N} 
\left| b^{\frak a}(\omega)\right|^2  
=\left[\Gamma_{\frak a} : \Gamma_{\frak a}'\right]\,\left\| {\bf b}_N^{\frak a}\right\|_2^2\;,
\eqno(5.9)$$
where, for $0\neq\omega\in{\frak O}$, 
$$b^{\frak a}(\omega) 
=\cases{{1\over 2}\,b(\omega){\rm e}\left( {\rm Re}\left( {1\over 2}\beta_{\frak a}\omega\right)\right)
+{1\over 2}\,b(-\omega){\rm e}\left( -{\rm Re}\left( {1\over 2}\beta_{\frak a}\omega\right)\right)  
 &if $q_0 \mu({\frak a})\mid 2\,$; \cr 
\quad &\hbox{\quad} \cr 
b(\omega) &otherwise.}\eqno(5.10)
$$

Regarding the other factors on the right-hand side of (5.5), one finds that, for 
$h$ as in (5.1) and $p\in{\Bbb Z}$, 
$$\eqalignno{
{1\over 4\pi^3 i}\int\limits_{(0)} h(\nu,p) \left( p^2 -\nu^2\right) {\rm d}\nu
 &={1\over 4\pi^3}\,\exp\left( -(p/P)^2\right)\int_{-\infty}^{\infty} 
\left( p^2+t^2\right)\exp\left( -(t/K)^2\right) {\rm d}t =\cr 
 &={K\over 4\pi^3}\,\exp\left( -(p/P)^2\right)\int_{-\infty}^{\infty} 
\left( p^2+K^2 x^2\right)\exp\left( -x^2\right) {\rm d}x =\cr
 &={K\over 4\pi^3}\,\exp\left( -(p/P)^2\right) \left( p^2 G_0(0)+K^2 G_2(0)\right) ,&(5.11)
}$$
where $G_0(0)$ and $G_2(0)$ are the constants given by (4.11)-(4.12) of Lemma 4.4. 
Moreover, by Poisson sumation over ${\Bbb Z}$ (i.e. the case $n=1$ of results 
(2.44)-(2.45) of Lemma 2.7) , one finds that if $k$ is a non-negative integer then 
$$\sum_{p=-\infty}^{\infty} p^{2k}\exp\left( -(p/P)^2\right) 
=\sum_{v=-\infty}^{\infty}\int_{-\infty}^{\infty} y^{2k}\exp\left( -(y/P)^2\right) 
{\rm e}(-vy){\rm d}y
=P^{1+2k}\sum_{v=-\infty}^{\infty} G_{2k}(-\pi Pv)\;,\eqno(5.12)$$ 
where $G_0(Y), G_1(Y), G_2(Y),  \ldots\ $ are given by (4.11)-(4.12) of Lemma 4.4, so that 
one has (in particular): 
$$G_2(Y)=iYG_1(Y)+{1\over 2}\,G_0(Y)
=\left( {1\over 2} - Y^2\right) G_0(Y)
=\left( {1\over 2} - Y^2\right) \sqrt{\pi}\,\exp\left( -Y^2\right)\qquad\hbox{($Y\in{\Bbb R}$).}
\eqno(5.13)$$
By (5.11), (5.12) (for $k=0,2$) and (5.13), one obtains 
$$\eqalignno{
{1\over 4\pi^3 i}\sum_{p=-\infty}^{\infty}
\int\limits_{(0)} h(\nu,p) \left( p^2 -\nu^2\right) {\rm d}\nu
 &={KP\over 4\pi^2}\sum_{v=-\infty}^{\infty}\left( P^2\left( {1\over 2} -(\pi Pv)^2\right) 
+{K^2\over 2}\right) \exp\left( -(\pi Pv)^2\right) =\qquad\qquad\cr
 &={1\over 8\pi^2}\,KP\left( K^2+P^2\right)
\left( 1+O\left( P^2 e^{-\pi^2 P^2}\right)\right)\quad\ \hbox{for $K,P\geq 1\,$.} &(5.14)
}$$ 
By applying the results of (5.9) and (5.14) one may now conclude that 
$${\cal D}^{\frak a}\!\left( q_0 ,P,K;N,b\right) 
={\left[\Gamma_{\frak a} : \Gamma_{\frak a}'\right]\over 8\pi^2 }
\left( P K^3+P^3 K\right)\left( 1+O\left( P^2 e^{-\pi^2 P^2}\right)\right)
\left\| {\bf b}_N^{\frak a}\right\|_2^2\;,\eqno(5.15)
$$
where the function $b^{\frak a} : {\frak O}-\{ 0\}\rightarrow{\Bbb C}$ is 
given by (5.10), and the meaning of $\left\| {\bf b}_N^{\frak a}\right\|_2$ 
is consistent with (1.9.16). 
 
In order to complete this proof it is necessary to obtain a 
suitable upper bound for the absolute value of the term 
${\cal L}^{\frak a}(q_0 ,P,K;N,b)$ in (5.4). As a first step towards this 
one may observe that, by (5.6), one has: 
$${\cal L}^{\frak a}\!\left( q_0 ,P,K;N,b\right) 
=\sum_{c\in {}^{\frak a}{\cal C}^{\frak a}}
{\Phi(c;P,K;N)\over |c|^2}\;,\eqno(5.16)$$
where, for $0\neq c\in{}^{\frak a}{\cal C}^{\frak a}\,$,
$$\Phi(c;P,K;N)
=\qquad\sum\!\!\!\!\!\!\!\!\!\!\sum_{\!\!\!\!\!\!\!\!\!\!\!\!\!{\scriptstyle\omega_1,\omega_2\in{\frak O}\atop\scriptstyle 
N/2<\left|\omega_1\right|^2,\left|\omega_2\right|^2\leq N}}
\!\!\!({\bf B}h)\!\!\left( {2\pi\sqrt{\omega_1\omega_2}\over c}\right)
S_{{\frak a},{\frak a}}\!\left(\omega_1 , \omega_2 ; c\right) 
\,\overline{b\left(\omega_1\right) }
\,b\left(\omega_2\right) .\eqno(5.17)$$

No attempt will be made to extract a saving from the 
averaging over $c$ that is apparent in the above: 
the approach taken will instead be to 
bound, individually, the absolute value of each  
term in the sum on the right-hand side of (5.16). 
Therefore suppose now that $c\in{}^{\frak a}{\cal C}^{\frak a}\,$. 
Then, as is recorded in Proposition~2, the number $c$  
satisfies (1.9.24), and is consequently (see (1.9.15), and Remark~1.9.7 below it) 
a non-zero Gaussian integer. Moreover, given Lemma 4.2, 
it follows trivially from the results (2.3), (2.4) and (2.12) of 
Lemma~2.1, Lemma~2.2 and Lemma~2.3 that  
$$\left| S_{{\frak a},{\frak a}}\!\left(\omega_1,\omega_2;c\right)\right| 
\leq\sum_{\delta\bmod c{\frak O}}1 
=|c|^2\qquad\hbox{($\omega_1,\omega_2\in{\frak O}$\/).}\eqno(5.18)$$

Given (5.1), one may prepare for the estimation of $\Phi(c;p,K;N)$ by applying Lemma 4.6
to the factor $({\bf B}h)( 2\pi\sqrt{\omega_1\omega_2}/c)$ on the 
right-hand side of (5.17); let it be supposed that 
Lemma 4.6 is applied for 
$$\Delta =\left(\left( |c|^2 +N\right) PK\right)^{\zeta}\;,\eqno(5.19)$$
where $\zeta$ is an arbitrary real number satisfying $0<\zeta\leq 1/5$ 
(one then has, by the hypotheses of the theorem, $\Delta\geq (NPK)^{\zeta}\geq 1$\/). 
The conditions of summation in (5.17) mean that one need only 
consider $({\bf B}h)(u)$ in cases where $|u|$ lies in the interval $[r,R]$ with endpoints 
$R=2\pi\sqrt{N}/|c|>0$ and $r=R/\sqrt{2}$. 
Since $\Delta\geq 1$ and $R/r<2$, there exists an $M\in{\Bbb N}$ such that 
(4.21) holds for all $u\in{\Bbb C}$ such that $|u|\in [r,R]$. This $M$ (considered fixed henceforth) 
will certainly satisfy 
$$1<M\asymp \left( 1+|c|^{-1} N^{1/2}\right)\Delta\;.\eqno(5.20)$$
After first rewriting the summand on the right-hand side of (5.17) by means of 
the result (4.22)-(4.25) of Lemma~4.6 (applied for  
$u=|u|e^{i\theta}=2\pi\sqrt{\omega_1\omega_2}/c$\/), and then 
making a permissible change to the order of summation and integration, 
one obtains: 
$$\eqalign{\Phi(c;P,K;N)
 &={1\over 4\pi^3}\int_{-\pi}^{\pi}\int_{-\infty}^{\infty}\int_{-\infty}^{\infty}
F_{P,K}(\eta,\xi) e^{-\xi^2-\eta^2}
\widetilde\Phi\!\left( c,M;N;\psi\left( {\xi\over K}\,,\,{\eta\over P}\,;\phi\right) ,\phi\right) 
{\rm d}\eta\,{\rm d}\xi\,{\rm d}\phi\ +\cr 
 &\quad\ +\Phi_M^{*}(c;P,K;N)\;,\matrix{\hbox{\ }\cr\hbox{\ }}
}\eqno(5.21)$$
where, with $\theta_z={\rm Arg}(z)$ (for $z\in{\Bbb C}$\/) and $F_{P,K}(\eta,\xi)$, $\psi(y,x;\phi)$ and $A_M(\phi,\theta)$ as in 
(4.23) and (4.24), one has   
$$\eqalign{ &\!\!\!\!\widetilde\Phi(c,M;N;\psi,\phi) =\cr 
 &\qquad\quad =\qquad\sum\!\!\!\!\!\!\!\!\!\!\sum_{\!\!\!\!\!\!\!\!\!\!\!\!\!{\scriptstyle\omega_1,\omega_2\in{\frak O}\atop\scriptstyle 
N/2<\left|\omega_1\right|^2,\left|\omega_2\right|^2\leq N}}
\!\!\!A_M\!\!\left(\phi,{\theta_{\omega_1}+\theta_{\omega_2}\over 2}-\theta_c\right) 
\cos\!\left( {2\pi\psi\sqrt{\left|\omega_1\omega_2\right|}\over |c|}\right) 
S_{{\frak a},{\frak a}}\left(\omega_1 , \omega_2 ; c\right) 
\,\overline{b\left(\omega_1\right) }
\,b\left(\omega_2\right)}\eqno(5.22)$$
and, for $j=1,2,\ldots\ $, 
$$\eqalignno{
 \Phi_M^{*}(c;P,K;N)
 &=\qquad\sum\!\!\!\!\!\!\!\!\!\!\sum_{\!\!\!\!\!\!\!\!\!\!\!\!\!{\scriptstyle\omega_1,\omega_2\in{\frak O}\atop\scriptstyle 
N/2<\left|\omega_1\right|^2,\left|\omega_2\right|^2\leq N}}
\!\!\!\overline{b\left(\omega_1\right) }
\,b\left(\omega_2\right) 
S_{{\frak a},{\frak a}}\left(\omega_1 , \omega_2 ; c\right)
E_M\!\!\left( P,K;{2\pi\sqrt{\omega_1\omega_2}\over c}\right) \ll_j 
\hbox{\qquad\qquad\qquad\qquad}\cr 
 &\ll_j\qquad\sum\!\!\!\!\!\!\!\!\!\!\sum_{\!\!\!\!\!\!\!\!\!\!\!\!\!{\scriptstyle\omega_1,\omega_2\in{\frak O}\atop\scriptstyle 
N/2<\left|\omega_1\right|^2,\left|\omega_2\right|^2\leq N}}
\!\!\!\left| b\left(\omega_1\right) b\left(\omega_2\right) 
S_{{\frak a},{\frak a}}\left(\omega_1 , \omega_2 ; c\right)\right| 
\left( P^2+K^2\right)\left( 1+{N^{1/2}\over |c|}\right)\Delta^{1-2j}\;. &(5.23)}$$
For later reference, note here that by (5.18), (5.19) and the Cauchy-Schwarz inequality, 
the case $j=[1/\zeta]+2$ of (5.23) implies the  bounds: 
$$\eqalignno{
\Phi_M^{*}(c;P,K;N)
 &=O_{\zeta}\!\!\left( |c|^2\left( P^2+K^2\right)\left( 1+|c|^{-1}N^{1/2}\right) 
\left(\left( |c|^2+N\right) PK\right)^{-2}
\Biggl(
\sum_{\scriptstyle\omega\in{\frak O}\atop\scriptstyle N/2<|\omega|^2\leq N} |b(\omega)|
\Biggr)^{\!\!2}\right) =\cr 
 &=O_{\zeta}\!\left(\left( P^{-2}+K^{-2}\right)|c|\left( |c|+N^{1/2}\right)
\left( |c|^2+N\right)^{-2} O(N)\left\|{\bf b}_N\right\|_2^2\right) \ll_{\zeta}\cr
 &\ll_{\zeta}\left( P^{-2}+K^{-2}\right)|c| N\left( |c|+N^{1/2}\right)^{\!-3}\left\|{\bf b}_N\right\|_2^2\;.
 &(5.24)}$$

Since $A_M(\phi ,\theta)$ is real-valued for $\phi ,\theta\in{\Bbb R}$, and 
since, by (1.5.10) and (1.5.13), 
$$\overline{S_{{\frak a},{\frak a}}\left(\omega_1,\omega_2;c\right)}
=S_{{\frak a},{\frak a}}\left(-\omega_1,-\omega_2;c\right)
=S_{{\frak a},{\frak a}}\left(\omega_2,\omega_1;c\right)\qquad
\hbox{($\omega_1,\omega_2\in{\frak O}$\/),}$$
it follows by Euler's identity, $e^{it}=\cos(t)+i\sin(t)$, that 
one may reformulate (5.22) as: 
$$\widetilde\Phi(c,M;N;\psi,\phi) 
={\rm Re}\left( \Phi^{\circ}(c,M;N;\psi,\phi) \right) ,$$
where 
$$\eqalign{ &\Phi^{\circ}(c,M;N;\psi,\phi) =\cr 
 &\qquad\qquad =\qquad\sum\!\!\!\!\!\!\!\!\!\!\sum_{\!\!\!\!\!\!\!\!\!\!\!\!\!{\scriptstyle\omega_1,\omega_2\in{\frak O}\atop\scriptstyle 
N/2<\left|\omega_1\right|^2,\left|\omega_2\right|^2\leq N}}
\!\!\!A_M\!\!\left(\phi,{\theta_{\omega_1}+\theta_{\omega_2}\over 2}-\theta_c\right) 
{\rm e}\!\left( {\psi\sqrt{\left|\omega_1\omega_2\right|}\over |c|}\right) 
S_{{\frak a},{\frak a}}\left(\omega_1 , \omega_2 ; c\right) 
\,\overline{b\left(\omega_1\right) }
\,b\left(\omega_2\right) ,
}\eqno(5.25)$$
Hence, and by (5.21), 
$$\eqalignno{\left|\Phi(c;P,K;N)\right| 
 &\leq{1\over 4\pi^3}\int_{-\pi}^{\pi}\int_{-\infty}^{\infty}\int_{-\infty}^{\infty}
e^{-\xi^2-\eta^2}\left| F_{P,K}(\eta,\xi) 
\,\Phi^{\circ}\!\!\left( c,M;N;\psi\!\left( {\xi\over K}\,,\,{\eta\over P}\,;\phi\right) ,\phi\right) \right| 
{\rm d}\eta\,{\rm d}\xi\,{\rm d}\phi\ +\cr 
 &\quad\ +\left|\Phi_M^{*}(c;P,K;N)\right|\matrix{\hbox{\ }\cr\hbox{\ }} &(5.26)}$$
(one can in fact show that by omitting both pairs of absolute value parentheses from the 
upper bound given in (5.26) one obtains, instead of that upper bound, 
an expression that is identically equal to $\Phi(c;P,K;N)$\/).

As an alternative to the use of (4.22) in (5.17), one may 
choose instead to apply 
in (5.17) the other result (4.26) of Lemma 4.6; then, 
via steps similar to those that produced (5.28) 
(and with the same 
choice of $\Delta$ and $M$ as before), one obtains:
$$\eqalignno{\left|\Phi(c;P,K;N)\right| 
 &\leq{1\over 2\pi}\int_{-\pi}^{\pi}\int_{-\infty}^{\infty}\int_{-\infty}^{\infty}
e^{-\xi^2-\eta^2} G_{P,K}(\eta,\xi) 
\,\left|\Phi^{\bullet}\!\!\left( c,M;N;\psi\!\left( {\xi\over K}\,,\,{\eta\over P}\,;\phi\right) ,\phi\right) \right| 
{\rm d}\eta\,{\rm d}\xi\,{\rm d}\phi\ +\cr 
 &\quad\ +\left|\Phi_M^{*}(c;P,K;N)\right|\matrix{\hbox{\ }\cr\hbox{\ }} &(5.27)}$$
where $\psi(y,x;\phi)$ and  $G_{P,K}(\eta,\xi)$ are as in (4.23)   
and (4.27), while 
$$\eqalign{ &\Phi^{\bullet}(c,M;N;\psi,\phi) =\cr 
 &\qquad\qquad =\qquad\sum\!\!\!\!\!\!\!\!\!\!\sum_{\!\!\!\!\!\!\!\!\!\!\!\!\!{\scriptstyle\omega_1,\omega_2\in{\frak O}\atop\scriptstyle 
N/2<\left|\omega_1\right|^2,\left|\omega_2\right|^2\leq N}}
\!\!\!\!A_M\!\!\left(\phi,{\theta_{\omega_1}+\theta_{\omega_2}\over 2}-\theta_c\right) 
{\rm e}\!\left( {\psi\sqrt{\left|\omega_1\omega_2\right|}\over |c|}\right) 
S_{{\frak a},{\frak a}}\left(\omega_1 , \omega_2 ; c\right) 
\,\overline{b^{\times}\!\left(\omega_1\right) }
\,b^{\times}\!\left(\omega_2\right) ,
}\eqno(5.28)$$
with $A_M(\phi,\theta)$ as in (4.24) and 
$$b^{\times}(\omega) =\left| {\omega\over c}\right| b(\omega)\qquad 
\hbox{($0\neq\omega\in{\frak O}$\/),}\eqno(5.29)$$
and $\Phi_M^{*}(c;P,K;N)$ is the same term seen in both (5.21) and (5.24). 
Now, since $\exp\bigl( i\theta_z\bigr) =z/|z|$ for $0\neq z\in{\Bbb C}$, 
it is a trivial consequence of (5.25) and (4.24) that 
$$\left|\Phi^{\circ}(c,M;N;\psi,\phi)\right| 
\leq U_{\frak a}(\psi ,c;M;N,b)\qquad\hbox{($\psi,\phi\in{\Bbb R}$),}\eqno(5.30)$$
where $U_{\frak a}(\psi ,c;M;N,b)$  is as given by equation (1.9.25) of Proposition~2. 
Similarly, by (5.28) and (4.24), 
$$\left|\Phi^{\bullet}(c,M;N;\psi,\phi)\right| 
\leq U_{\frak a}\left(\psi ,c;M;N,b^{\times}\right)\qquad\hbox{($\psi,\phi\in{\Bbb R}$),}\eqno(5.31)$$
where $U_{\frak a}\left(\psi ,c;M;N,b^{\times}\right)$ is given by 
(1.9.25) with $b^{\times}$ substituted for $b$ throughout. 

Given  (5.20), and given that $\Delta\geq 1$, it 
follows from (5.30) and the result (1.9.27) of Proposition~2 that  
$$\eqalignno{
\Phi^{\circ}(c,M;N;\psi,\phi)
 &\ll (1+|\psi|)^{1/2}\left( |c|M+N^{1/2}\right)\left( |c|+N^{1/2}\right)\left\|{\bf b}_N\right\|_2^2 \ll\cr
 &\ll (1+|\psi|)^{1/2}\Delta\left( N^{1/2}+|c|\right)^2\left\|{\bf b}_N\right\|_2^2\qquad\qquad\qquad  
\hbox{for  $\psi,\phi\in{\Bbb R}\,$.}&(5.32)
}$$
Moreover, since (5.29) and (1.9.16) imply 
$$\left\|{\bf b}_N^{\times}\right\|_2^2\leq |c|^{-2} N\left\|{\bf b}_N\right\|_2^2\;,\eqno(5.33)$$
it likewise follows by (5.31) and (1.9.27) (with $b^{\times}$ substituted for $b$\/) 
that 
$$\eqalignno{
\Phi^{\bullet}(c,M;N;\psi,\phi) 
 &\ll (1+|\psi|)^{1/2}\Delta\left( |c|+N^{1/2}\right)^2 |c|^{-2} N\left\|{\bf b}_N\right\|_2^2 \ll\cr 
 &\ll (1+|\psi|)^{1/2}\Delta\left( N^{1/2}+|c|^{-1}N\right)^2\left\|{\bf b}_N\right\|_2^2\qquad\quad   
\hbox{for $\psi,\phi\in{\Bbb R}\,$.} 
&(5.34) }$$

By (4.23) and (4.27), 
$$|\psi(y,x;\phi)|\leq 2\cosh(y)\leq 2\exp(|y|)\qquad\hbox{($\phi,x,y\in{\Bbb R}$\/)}\eqno(5.35)$$
and 
$$0\leq G_{P,K}(\eta,\xi)
\leq 2\left( (\xi /K)^2\cosh^2(\xi /K)+(\eta /P)^2\right) 
\ll e^{2|\xi|/K}(\xi /K)^2+(\eta /P)^2\;,\eqno(5.36)$$
so that (given that $K\geq 1$) one has, for $-\pi\leq\phi\leq\pi$, 
$$\eqalignno{ 
 &\int_{-\infty}^{\infty}
\int_{-\infty}^{\infty}
e^{-\xi^2 -\eta^2}G_{P,K}(\eta ,\xi)
\left( 1+\left|\psi\left( {\xi\over K} , {\eta\over P} ; \phi\right)\right|\right)^{1/2} 
{\rm d}\eta\,{\rm d}\xi \ll\matrix{\hbox{\quad}\cr\hbox{\quad}\cr\hbox{\quad}}
\hbox{\qquad\qquad\qquad\qquad\qquad\qquad\qquad\qquad\quad\ }\cr
 &\qquad\qquad\qquad\qquad\qquad\qquad\qquad\qquad\qquad\ll 
\int_{-\infty}^{\infty}
\int_{-\infty}^{\infty}
e^{-\xi^2 -\eta^2}
\left( e^{2|\xi|/K}(\xi /K)^2+(\eta /P)^2\right) e^{|\xi|/(2K)}
{\rm d}\eta\,{\rm d}\xi \leq\cr
 &\qquad\qquad\qquad\qquad\qquad\qquad\qquad\qquad\qquad\leq 2\left( K^{-2}+P^{-2}\right) 
\int_{-\infty}^{\infty}
\int_{-\infty}^{\infty}
\left(\xi^2 +\eta^2\right) e^{(5/2)\xi-\xi^2 -\eta^2}
{\rm d}\eta\,{\rm d}\xi \ll\cr 
 &\qquad\qquad\qquad\qquad\qquad\qquad\qquad\qquad\qquad\ll K^{-2}+P^{-2}\;. 
&(5.37)}$$
Therefore, by applying the estimates (5.24) and (5.34) 
for terms on the right-hand side of (5.27), one finds (given (5.19), and since $N\geq 1$\/) that 
$$\eqalignno{ 
\Phi(c;P,K;N) 
 &\ll \left( P^{-2}+K^{-2}\right)\left(\Delta \left( N +|c|^{-2}N^2\right)
+O_{\zeta}\left( |c|^{-2}N\right)\right)\left\|{\bf b}_N\right\|_2^2 \ll\cr 
 &\ll\left( P^{-2}+K^{-2}\right) (PK)^{\zeta}
\left( |c|^2 +O_{\zeta}(N)\right)^{1+\zeta} 
|c|^{-2} N\left\|{\bf b}_N\right\|_2^2 \ll_{\zeta}\cr 
 &\ll_{\zeta} (PK)^{\zeta}\left( P^{-2}+K^{-2}\right) 
\left( N |c|^{2\zeta}+N^{2+\zeta}|c|^{-2}\right) \left\|{\bf b}_N\right\|_2^2 \;. 
 &(5.38)
}$$

For a further alternative upper bound on $\Phi(c;P,K;N)$ 
(useful when $|c|$ is large), one may apply (5.27), and then (5.31) 
combined with the result (1.9.26) of Proposition~2  
(with $b^{\times}$ substituted for $b$ there):   
given (5.33) and (5.37), and given the upper bound   
$\tau(c)\ll_{\zeta} |c|^{2\zeta}$ (which is an elementary corollary of the fundamental 
theorem of arithmetic for the Gaussian integers), 
one may in this way obtain the estimate 
$$\eqalign{
\Phi(c;P,K;N)
 &\ll \left( P^{-2}+K^{-2}\right) \tau^{3/2}(c) |c| M N \left\| {\bf b}_N^{\times}\right\|_2^2 
+\left| I_M^{*}(c;P,K;N)\right| =\cr 
 &=O_{\zeta}\!\left( 
\left( P^{-2}+K^{-2}\right) |c|^{2\zeta -1} M N^2 \left\|{\bf b}_N\right\|_2^2\right) 
+\left| I_M^{*}(c;P,K;N)\right| .
}$$
Moreover, by (5.19), (5.20) and (5.24), and since $|c|,N,\Delta\geq 1$, 
it follows from this 
last bound that one has 
$$\eqalignno{
\Phi(c;P,K;N) 
 &=O_{\zeta}\left(  
\left( P^{-2}+K^{-2}\right)\left( |c|^{2\zeta -1} 
\bigl( 1+|c|^{-1}N^{1/2}\bigr)\Delta N^2
+|c|^{-2}N\right) \left\|{\bf b}_N\right\|_2^2\right) \ll_{\zeta}\cr 
 &\ll_{\zeta} 
(PK)^{\zeta}\left( P^{-2}+K^{-2}\right) |c|^{4\zeta -1} 
N^2 \bigl( 1+|c|^{-1}N^{1/2}\bigr)^{1+2\zeta} 
\left\|{\bf b}_N\right\|_2^2\;. 
&(5.39)}$$

Given that $0<\zeta\leq 1/5$, one may deduce from the bounds (5.38) and (5.39) that 
$$\Phi(c;P,K;N)
\ll_{\zeta} (PK)^{\zeta}\left( P^{-2}+K^{-2}\right) |c|^{-\zeta} N^{1+5\zeta}
\left\|{\bf b}_N\right\|_2^2\quad\ \hbox{if\  $\,|c|^2>N^{1-2\zeta}\,$}\eqno(5.40)$$
(i.e. the bound (5.39) implies this for $|c|>N$, while the bound (5.38) implies this 
for $N\geq |c|>N^{(1/2)-\zeta}$\/). 
Now the sum $\sum_{0\neq\alpha\in{\frak O}}|\alpha|^{-2\sigma}$ is convergent when  
$\sigma >1$, and so, 
given the definition (1.9.15) and Remark~1.9.7,  
it follows by the result (1.9.24) of Proposition~2 that 
$$\sum_{c\in{}^{\frak a}{\cal C}^{\frak a}}
|c|^{-2\sigma}
\leq\sum_{0\neq\alpha\in{\frak O}}\left|{\alpha\over\mu({\frak a})}\right|^{-2\sigma}
\ll_{\sigma} |\mu({\frak a})|^{2\sigma}\leq |\mu({\frak a})|^{2}\qquad\quad\hbox{($\sigma >1$).}\eqno(5.41)$$
Therefore application of the bound (5.40) shows that 
$$\eqalignno{
\sum_{\scriptstyle c\in{}^{\frak a}{\cal C}^{\frak a}\atop\scriptstyle 
|c|^2 > N^{1-2\zeta}}
{|\Phi(c;P,K;N)|\over |c|^2}
 &\ll_{\zeta} (PK)^{\zeta}\left( P^{-2}+K^{-2}\right)
N^{1+5\zeta}\left\|{\bf b}_N\right\|_2^2
\sum_{c\in{}^{\frak a}{\cal C}^{\frak a}}
{1\over |c|^{2+\zeta}} \ll_{\zeta}\cr 
 &\ll_{\zeta} (PK)^{\zeta}\left( P^{-2}+K^{-2}\right)
|\mu({\frak a})|^2 N^{1+5\zeta}\left\|{\bf b}_N\right\|_2^2\;.\matrix{\hbox{\ }\cr\hbox{\ }}
 &(5.42)}$$

Suppose now that 
$$0<|c|^2\leq N^{1-2\zeta}\;.\eqno(5.43)$$
Since $\zeta >0$ and $N\geq 1$, it is certainly implied by 
(5.43) that $c$ and $N$ satisfy the case $A_1=1$, $\varepsilon =\zeta$ of 
the condition (1.9.28) of Proposition~2. 
Indeed, given (5.19), (5.20), (5.30), (5.31), (5.33) and (5.43), 
the result (1.9.28)-(1.9.30) of Proposition~2 implies that if 
$\phi ,\psi\in{\Bbb R}$ and 
$$0<|\psi|\leq 2e\quad\hbox{(say)}\eqno(5.44)$$
then one has: 
$$\eqalignno{ 
\Phi^{\circ}(c,M;N;\psi,\phi) 
 &=O_{\zeta}\left( |\psi|^{-1/2}\right) 
\left( |c|^{1/2}N^{3/4}+M |c|^{3/2} N^{1/4}\right) N^{\zeta}\left\|{\bf b}_N\right\|_2^2 =\cr 
 &=O_{\zeta}\left( |\psi|^{-1/2}\right) 
\left( |c|^{1/2}N^{3/4}+\Delta\left( |c|^{3/2} N^{1/4}+|c|^{1/2}N^{3/4}\right)\right) 
N^{\zeta}\left\|{\bf b}_N\right\|_2^2 \ll\cr 
 &\ll O_{\zeta}\left( |\psi|^{-1/2}\right) 
\Delta |c|^{1/2}N^{(3/4)+\zeta}\left\|{\bf b}_N\right\|_2^2 \ll_{\zeta}\cr 
 &\ll_{\zeta} |\psi|^{-1/2} (PK)^{\zeta}
N^{(3/4)+2\zeta}|c|^{1/2}\left\|{\bf b}_N\right\|_2^2
 &(5.45)}$$
and, similarly, 
$$\eqalignno{ 
\Phi^{\bullet}(c,M;N;\psi,\phi) 
 &=O_{\zeta}\left( |\psi|^{-1/2} (PK)^{\zeta}
N^{(3/4)+2\zeta}|c|^{1/2}\left\|{\bf b}_N^{\times}\right\|_2^2\right) \ll_{\zeta}\cr 
 &\ll_{\zeta} |\psi|^{-1/2} (PK)^{\zeta}
N^{(7/4)+2\zeta}|c|^{-3/2}\left\|{\bf b}_N\right\|_2^2 \;.
 &(5.46)}$$

In order to facilitate the application of the last bound one may note that,  
by (5.27) and (4.27), it is certainly the case that one has 
$$|\Phi(c;P,K;N)|
\leq {1\over 2\pi}\Biggl(\sum_{j=-2}^2 I_j\Biggr)
+\Phi_M^{*}(c;P,K;N)\;,\eqno(5.47)$$
where 
$$I_0
=\int_{-K}^K\int_{-P}^P e^{-\xi^2 -\eta^2} G_{P,K}(\eta,\xi) 
\,\iota\!\!\left( {\xi\over K} , {\eta\over P}\right) {\rm d}\eta {\rm d}\xi\;,$$
$$I_{\pm 1}
=\int_{K}^{\infty}\int_{-\infty}^{\infty} e^{-\xi^2 -\eta^2} G_{P,K}(\eta,\xi) 
\,\iota\!\!\left( \pm{\xi\over K} , {\eta\over P}\right) {\rm d}\eta {\rm d}\xi\;,$$
$$I_{\pm 2}
=\int_{-\infty}^{\infty}\int_{P}^{\infty} e^{-\xi^2 -\eta^2} G_{P,K}(\eta,\xi) 
\,\iota\!\!\left( {\xi\over K} , \pm{\eta\over P}\right) {\rm d}\eta {\rm d}\xi $$
and 
$$\iota(y,x)
=\int_{-\pi}^{\pi}
\left|\Phi^{\bullet}\bigl( c,M;N;\psi(y,x;\phi),\phi\bigr)\right| {\rm d}\phi\;.$$

By (4.23) and (5.35), the condition (5.44) will hold for $\psi =\psi(y,x;\phi)$ 
if $y\in[-1,1]$, $x,\phi\in(-\pi,\pi)$ and $\tan(\phi)\neq \tanh(y)\tan(x)$. 
Moreover, since $0<1<\pi/2$, it follows by Lemma 4.8 that 
$$\int_{-\pi}^{\pi}|\psi(y,x;\phi)|^{-1/2} {\rm d}\phi 
\ll {1\over\sqrt{\cos(1)}}\ll 1\qquad\hbox{for $x\in[-1,1]$ and $y\in{\Bbb R}\,$.}$$
Therefore it follows  from the bound (5.46) for 
$\left|\Phi^{\bullet}\bigl( c,M;N;\psi(y,x;\phi),\phi\bigr)\right|$ that  one has: 
$$\iota(y,x)
\ll_{\zeta}  (PK)^{\zeta}
N^{(7/4)+2\zeta}|c|^{-3/2}\left\|{\bf b}_N\right\|_2^2 \qquad 
\hbox{for $x,y\in[-1,1]\,$.}$$ 
This, together with the bound obtained in (5.37),  enables one to conclude that 
$$I_0
\ll_{\zeta} (PK)^{\zeta} 
\left( P^{-2}+K^{-2}\right) N^{(7/4)+2\zeta} |c|^{-3/2} 
\left\|{\bf b}_N\right\|_2^2\;.\eqno(5.48)$$

In estimating $I_{\pm 1}$ and $I_{\pm 2}$ 
one may use the bound 
$$\iota(y,x)
\ll\exp(|y|/2) (PKN)^{\zeta}
N^2 |c|^{-2}\left\|{\bf b}_N\right\|_2^2 \qquad 
\hbox{($x,y\in{\Bbb R}$\/),}\eqno(5.49)$$
which follows, given (5.43), from (5.34), (5.35) and (5.19). 
Indeed, since one has 
$$\eqalign{
\int_K^{\infty}\int_{-\infty}^{\infty}
e^{-\xi^2-\eta^2}
\left( e^{2\xi /K}\!\left( {\xi\over K}\right)^{\!\!2} +\left( {\eta\over P}\right)^{\!2}\right) 
e^{\xi /(2K)} {\rm d}\eta\,{\rm d}\xi  &\ll \left( {1\over P^{2}}+{1\over K^{2}}\right) 
\int_K^{\infty} e^{(5/2)\xi -\xi^2}\xi^2 {\rm d}\xi \leq\cr
 &\leq \left( {1\over P^{2}}+{1\over K^{2}}\right) \exp\left( -{K^2\over 2}\right) 
\int_{-\infty}^{\infty} e^{\left( 5\xi -\xi^2\right) /2}\xi^2 {\rm d}\xi \;,
}$$
it follows from (5.49) and (5.36) that 
$$I_{\pm 1}
\ll \exp\left( -K^2 /2\right) (PK)^{\zeta} \left( P^{-2}+K^{-2}\right) 
 N^{2+\zeta} |c|^{-2} \left\|{\bf b}_N\right\|_2^2 \;.$$
Similar reasoning shows that the estimates (5.49) and (5.36) imply also that 
$$I_{\pm 2} 
\ll \exp\left( -P^2 /2\right) (PK)^{\zeta} \left( P^{-2}+K^{-2}\right) 
 N^{2+\zeta} |c|^{-2} \left\|{\bf b}_N\right\|_2^2 \;.$$

By (5.24), (5.47), (5.48) and the bounds just obtained for the integrals 
$I_{\pm 1}$ and $I_{\pm 2}$, it follows that, 
subject to the condition (5.43) holding, one has: 
$$\eqalignno{ 
\Phi(c;P,K;N) &\ll O_{\zeta}\left(\left( P^{-2}+K^{-2}\right) 
\bigl( (PK)^{\zeta} N^{(7/4)+2\zeta} |c|^{-3/2} 
+N^{-1/2} |c|\bigr) \left\|{\bf b}_N\right\|_2^2\right)\ +\cr
&\quad\ +\bigl( e^{-P^2 /2}+e^{-K^2 /2}\bigr) (PK)^{\zeta} \left( P^{-2}+K^{-2}\right) 
N^{2+\zeta} |c|^{-2} \left\|{\bf b}_N\right\|_2^2 \ll 
\hbox{\qquad\qquad\qquad\qquad\qquad\qquad} \cr 
 &\ll (PK)^{\zeta}\bigl( P^{-2}+K^{-2}\bigr) 
\!\left( O_{\zeta}\bigl( N^{(7/4)+2\zeta} |c|^{-3/2}\bigr) 
+\bigl( e^{-P^2 /2}+e^{-K^2 /2}\bigr) N^{2+\zeta} |c|^{-2}\right) 
\!\left\|{\bf b}_N\right\|_2^2 .
 &(5.50)}$$
Note that (5.50) was derived principally from (5.27) and the bounds for 
$\left| \Phi^{\bullet}(c,M;N;\psi,\phi)\right|$ given by (5.46) and (5.34). 
By  using instead (5.26) and 
the bounds (5.45) and (5.32)  
found for $\left|\Phi^{\circ}(c,M;N;\psi,\phi)\right|$, 
one similarly obtains (as an alternative to (5.50)) the estimate 
$$\Phi(c;P,K;N) \ll (PK)^{\zeta}\bigl( P^2+K^2\bigr) 
\left( O_{\zeta}\bigl( N^{(3/4)+2\zeta} |c|^{1/2}\bigr) 
+\bigl( e^{-P^2 /2}+e^{-K^2 /2}\bigr) N^{1+\zeta}\right) 
\left\|{\bf b}_N\right\|_2^2\;,\eqno(5.51)$$
subject to the condition (5.43) holding. 
\par 
Since $0<\zeta\leq 1/5<3/2$, it follows by (5.41) that one has:  
$$\sum_{\scriptstyle c\in{}^{\frak a}{\cal C}^{\frak a}\atop\scriptstyle |c|>X} 
|c|^{-7/2}< X^{\zeta -(3/2)}\sum_{c\in{}^{\frak a}{\cal C}^{\frak a}} |c|^{-2-\zeta} 
\ll_{\zeta} X^{\zeta -(3/2)} |\mu({\frak a})|^2\qquad 
\hbox{for $X>0\,$.}$$
Hence, and with the aid of the case $\sigma =2$ of (5.41), one may deduce from (5.50) that 
$$\eqalignno{ 
\sum_{\scriptstyle c\in{}^{\frak a}{\cal C}^{\frak a}\atop\scriptstyle 
N^{(1/2)-\zeta}\geq |c|>\textstyle{N^{(1/2)-\zeta}\over PK}} 
 &\!\!\!\!\!{|\Phi(c;P,K;N)|\over |c|^2} \ll 
\hbox{\qquad\qquad\qquad\qquad\qquad\qquad\qquad\qquad\qquad\qquad\qquad\qquad\qquad\qquad\qquad} \cr 
 &\qquad\qquad\ll (PK)^{\zeta}\!\left( {1\over P^2}+{1\over K^2}\right) 
\!\Biggl( O_{\zeta}\!\biggl( 
N^{(7/4)+2\zeta}\left( {N^{(1/2)-\zeta}\over PK}\right)^{\!\!\zeta -(3/2)}
|\mu({\frak a})|^2\biggr)\ +\cr
&\qquad\qquad\qquad\qquad\qquad\qquad\qquad\quad\ +\bigl( e^{-P^2 /2}+e^{-K^2 /2}\bigr) N^{2+\zeta} |\mu({\frak a})|^4\Biggr) 
\!\left\|{\bf b}_N\right\|_2^2 \ll_{\zeta}\cr 
 &\qquad\qquad\ll_{\zeta} 
(PK)^{\zeta}\!\left( {1\over P^2}+{1\over K^2}\right) 
\!\Biggl( (PK)^{3/2} N^{1+4\zeta} |\mu({\frak a})|^2\ +\cr 
&\qquad\qquad\qquad\qquad\qquad\qquad\qquad\quad\     
\,+\bigl( e^{-P^2 /2}+e^{-K^2 /2}\bigr) N^{2+\zeta} |\mu({\frak a})|^4 
\Biggr)\!\left\|{\bf b}_N\right\|_2^2 .  
 &(5.52)}$$
Moreover, since the bound (5.51) shows that if $|c|\leq N^{(1/2)-\zeta}/(PK)$ then
$$|c|^{\zeta}\Phi(c;P,K;N) \ll (PK)^{\zeta}\bigl( P^2+K^2\bigr) 
\left( O_{\zeta}\bigl( (PK)^{-1/2}\bigr) 
+e^{-P^2 /2}+e^{-K^2 /2}\right) N^{1+2\zeta}
\left\|{\bf b}_N\right\|_2^2\;,$$
one may therefore use the case $\sigma =1+(\zeta /2)$ of (5.41) to deduce from (5.51) that  
$$\eqalignno{
\sum_{\scriptstyle c\in{}^{\frak a}{\cal C}^{\frak a}\atop\scriptstyle 
|c|\leq N^{(1/2)-\zeta}/(PK)} 
\!\!&{|\Phi(c;P,K;N)|\over |c|^2} \ll_{\zeta}\cr 
&\qquad\quad\ll_{\zeta} 
(PK)^{\zeta}\bigl( P^2+K^2\bigr) 
\left( (PK)^{-1/2}+e^{-P^2 /2}+e^{-K^2 /2}\right) 
N^{1+2\zeta}|\mu({\frak a})|^2 
\left\|{\bf b}_N\right\|_2^2 =\qquad\quad\cr 
&\qquad\quad = (PK)^{\zeta}\left( {1\over P^2}+{1\over K^2}\right) 
\Bigl( (PK)^{3/2}N^{1+2\zeta}|\mu({\frak a})|^2\ +\cr 
 &\qquad\qquad\qquad\qquad\qquad\qquad\qquad      
+\bigl( e^{-P^2 /2}+e^{-K^2 /2}\bigr) (PK)^2 N^{1+2\zeta}|\mu({\frak a})|^2\Bigr) 
\left\|{\bf b}_N\right\|_2^2\;. 
 &(5.53)}$$

Observe now that, given the result (1.9.24) of Proposition~2, 
the sum over $c$ on the left-hand side of (5.53) is in fact an empty sum unless 
$|\mu({\frak a})| N^{(1/2)-\zeta} /(PK)\geq 1$. 
One may therefore replace the upper bound on the right-hand side of 
(5.53) by  
$$O_{\zeta}\left(  
  (PK)^{\zeta}\!\left( {1\over P^2}+{1\over K^2}\right) 
\!\Bigl( (PK)^{3/2}N^{1+2\zeta}|\mu({\frak a})|^2  
+\bigl( e^{-P^2 /2}+e^{-K^2 /2}\bigr) N^2 |\mu({\frak a})|^4\Bigr)
\!\left\|{\bf b}_N\right\|_2^2\right) .$$
Given this modification of (5.53), together with the complementary bounds found in (5.52) and (5.42),
it now follows by the triangle inequality that the sum 
${\cal L}^{\frak a}(q_0 ,P,K;N,b)$ defined in (5.16) must satisfy 
$$\eqalign{
 &{\cal L}^{\frak a}\!\left( q_0 ,P,K;N,b\right) \ll_{\zeta}\cr
 &\qquad\qquad\ \ll_{\zeta}
 (PKN)^{5\zeta}\!\left( P^{-2}+K^{-2}\right) 
\!\Bigl( (PK)^{3/2}N |\mu({\frak a})|^2  
+\bigl( e^{-P^2 /2}+e^{-K^2 /2}\bigr) N^2 |\mu({\frak a})|^4\Bigr)
\!\left\|{\bf b}_N\right\|_2^2 .}\eqno(5.54)$$
 
Since (5.54) has been established for an arbitrary $\zeta\in (0,1/5]$, it follows 
now by (5.4), (5.15) and (5.54) that, for $0<\eta\leq 1$,  
$$\eqalign{
\!\!{\cal E}^{\frak a}\!\left( q_0 ,P,K;N,b\right)   
 &={\left[\Gamma_{\frak a} : \Gamma_{\frak a}'\right]\over 8\pi^2}
\left( P^2 +K^2\right) PK\left( 1+O\bigl( P^2 e^{-\pi^2 P^2}\bigr)\right)
\left\| {\bf b}_N^{\frak a}\right\|_2^2\ +\cr
 &\quad\,+O_{\eta}\Biggl( 
 (PKN)^{\eta}\!\left( P^{2}+K^{2}\right) 
\!\biggl( {N |\mu({\frak a})|^2 \over (PK)^{1/2}} 
+\bigl( e^{-P^2 /2}+e^{-K^2 /2}\bigr) {N^2 |\mu({\frak a})|^4 \over (PK)^{2}}\biggr)
\!\left\|{\bf b}_N\right\|_2^2\Biggr) . 
}\eqno(5.55)$$
Note here that, by (5.10), (1.9.16) and the arithmetic-geometric mean inequality, one has 
$$\left\|{\bf b}_N^{\frak a}\right\|_2^2\leq\left\|{\bf b}_N\right\|_2^2\;.$$
Therefore if one supposes now that $j\in\{ 0,1\}$, then 
it may be deduced from (5.2), (5.8) and (5.55) that, for $P_1,K_1\geq 1$ and $0<\eta\leq 1$, 
one has: 
$$\eqalign{
 &E_j^{\frak a}\!\left( q_0 , P_1 , K_1 ;N,b\right) \ll\cr  
 &\qquad\quad\ll \left( P_1^2+K_1^2\right)
\!\Biggl(\!P_1 K_1 
+ O_{\eta}\bigl( (P_1 K_1 N)^{\eta}\bigr) 
\!\biggl( {N |\mu({\frak a})|^2 \over (P_1 K_1)^{1/2}} 
+\bigl( e^{-P_1}+e^{-K_1}\bigr) {N^2 |\mu({\frak a})|^4 \over (P_1 K_1)^{2}}\biggr)\!\Biggr) 
\!\left\|{\bf b}_N\right\|_2^2\;.
}\eqno(5.56)$$
By the definition (1.9.15) and Remark~1.9.7, one has 
$1/\mu({\frak a})\in{\frak O}$, so that $0<|\mu({\frak a})|\leq 1$. 
Therefore, and since $P,K\geq 1$, it follows by 
(5.56) for $\eta =1/2$, $P_1=P$ and $K_1=K$ that 
one has (in particular): 
$$E_j^{\frak a}\!\left( q_0 ,P,K;N,b\right)  
\ll\left( P^2 +K^2\right) 
\left( PK + N^{3/2} +(PK)^{-3/2} N^{5/2}\right) 
\left\|{\bf b}_N\right\|_2^2\;.\eqno(5.57)$$

The aim now is to show that (5.57) and (5.56) (for any given $\eta\in(0,1]$\/) 
imply the bound 
$$E_j^{\frak a}\!\left( q_0 , P , K ;N,b\right)
\ll\left( P^2 +K^2\right) 
\Biggl( PK 
+O_{\eta}\!\left( {N^{1+4\sqrt{\eta}}\,|\mu({\frak a})|^2\over (PK)^{1/2}}\right)\Biggr) 
\left\|{\bf b}_N\right\|_2^2\;.\eqno(5.58)$$
This (since (5.56) holds for all $\eta\in(0,1]$)  
will be enough to prove the theorem: for, in cases where $0<\varepsilon\leq 4$, 
the bound (1.9.14) follows 
immediately from (5.58) for $\eta =(\varepsilon /4)^2$;  
while, in cases where $\varepsilon >4$, the bound (1.9.14) is 
(given that $N\geq 1$\/) a trivial corollary of  (1.9.14) for  $\varepsilon =4$.

Now since (5.57) implies the desired result (5.58) if  
$PK>N^{3/2}$,  one may henceforth suppose that 
$$PK\leq N^{3/2}\;.\eqno(5.59)$$
Independently of the conclusion just reached, one may observe also that if 
one has 
$$P_1\geq P,\quad K_1\geq K\quad {\rm and}\quad P_1 K_1\leq N^2\;,\eqno(5.60)$$
then, by (1.9.12)-(1.9.13) and (5.56), it follows that 
$$\eqalignno{
E_j^{\frak a}\!\left( q_0 ,P,K;N,b\right)  
 &\leq E_j^{\frak a}\!\left( q_0 ,P_1 ,K_1 ;N,b\right) \ll 
\hbox{\qquad\qquad\qquad\qquad\qquad\qquad\qquad\qquad\qquad\qquad\qquad\qquad\qquad\qquad\ } \cr 
 &\ll \left( P_1^2+K_1^2\right)
\!\Biggl(\!P_1 K_1 
+ O_{\eta}\!\left( N^{3\eta}\right) 
\!\biggl( {N |\mu({\frak a})|^2 \over (PK)^{1/2}} 
+\bigl( e^{-P_1}+e^{-K_1}\bigr) {N^2 |\mu({\frak a})|^4 \over (PK)^{2}}\biggr)\!\Biggr) 
\!\left\|{\bf b}_N\right\|_2^2 . &(5.61) 
}$$

Given (5.59), the condition (5.60) is (in particular) satisfied when 
$P_1=P$ and $K_1 =K$, so that (5.61) holds in this case. 
Hence, on noting that $e^{-P}+e^{-K}\leq 2/e$, one obtains the bound 
$$E_j^{\frak a}\!\left( q_0 ,P,K;N,b\right)  
\ll\left( P^2+K^2\right)
\!\Biggl(\!PK 
+ O_{\eta}\!\left( N^{3\eta}\right) 
\!\biggl( {N |\mu({\frak a})|^2 \over (PK)^{1/2}} 
+{N^2 |\mu({\frak a})|^4 \over (PK)^{2}}\biggr)\!\Biggr) 
\!\left\|{\bf b}_N\right\|_2^2 ,\eqno(5.62)$$
which (since $PK\geq 1$, $N\geq 1$ and $\sqrt{\eta}\geq\eta >0$\/) 
implies the desired result (5.58) if $N |\mu({\frak a})|^2 < N^{\sqrt{\eta}}$. 
Therefore it is only cases where  one has 
$$N |\mu({\frak a})|^2\geq N^{\sqrt{\eta}}$$
that require any further consideration. In these cases 
$N^{\eta}\leq \left( N |\mu({\frak a})|^2\right)^{\sqrt{\eta}}$, so 
that, by (5.62), one obtains: 
$$E_j^{\frak a}\!\left( q_0 ,P,K;N,b\right)  
\ll\left( P^2+K^2\right)
\!\Biggl(\!PK 
+ O_{\eta}\!\biggl( {\left( N |\mu({\frak a})|^2\right)^{1+3\sqrt{\eta}} \over (PK)^{1/2}} 
+{\left( N |\mu({\frak a})|^2\right)^{2+3\sqrt{\eta}}\over (PK)^{2}}\biggr)\!\Biggr) 
\!\left\|{\bf b}_N\right\|_2^2\;.$$  
Since $N\geq 1\geq |\mu({\frak a})|$, the bound just noted certainly implies 
the desired result (5.58) if $(PK)^{3/2}>N |\mu({\frak a})|^2$, 
and so one may henceforth suppose that 
$$PK\leq {N |\mu({\frak a})|^2\over (PK)^{1/2}}\;.\eqno(5.63)$$

Take now
$$P_1=\max\left\{ P , \left( N |\mu({\frak a})|^2\right)^{\eta /2}\right\}\quad\hbox{and}\quad 
K_1=\max\left\{ K , \left( N |\mu({\frak a})|^2\right)^{\eta /2}\right\}\;.$$
Since $N,P,K\geq 1$ and $\left( N |\mu({\frak a})|^2\right)^{\eta /2}\leq N^{\eta /2}\leq N^{1/2}$, 
it follows by (5.59) that (5.60) holds. Therefore (5.61) also holds. 
Moreover, since 
$$e^{-P_1}+e^{-K_1}\leq 2\exp\left( -\left( N |\mu({\frak a})|^2\right)^{\eta /2}\right) 
\ll_{\eta} \left( N |\mu({\frak a})|^2\right)^{-1}\quad\hbox{and}\quad 
PK\geq 1\;,$$
it follows from (5.61) that 
$$E_j^{\frak a}\!\left( q_0 ,P,K;N,b\right)  
\ll\left( P_1^2+K_1^2\right) 
\left( P_1 K_1 + 
O_{\eta}\!\left( {N^{1+3\eta} |\mu({\frak a})|^2 \over (PK)^{1/2}} \right) 
\right) \left\|{\bf b}_N\right\|_2^2\;.$$
Now observe that this implies the desired result (5.58): 
for one has, in the above bound,  
$$P_1^2+K_1^2\leq \left( N |\mu({\frak a})|^2\right)^{\eta}\left( P^2+K^2\right) 
\leq N^{\eta}\left( P^2+K^2\right)$$ 
and, by (5.63), 
$$P_1 K_1\leq\left( N |\mu({\frak a})|^2\right)^{\eta}PK 
\leq {\left( N |\mu({\frak a})|^2\right)^{1+\eta}\over (PK)^{1/2}}
\leq {N^{1+\eta} |\mu({\frak a})|^2\over (PK)^{1/2}}\;,$$
where $N^{\eta}\leq N^{\sqrt{\eta}}$ (given that $0<\eta\leq 1$ and $N\geq 1$).
Since no other cases remain to be considered, it has therefore now 
been shown that, subject to the hypotheses of the theorem, the bound 
(5.58) holds for $j\in\{ 0,1\}$ and $0<\eta\leq 1$. 
For the reasons mentioned below (5.58), this completes proof of the theorem\quad$\blacksquare$

\bigskip

\centerline{\bf \S 6. Appendix on the proof of the sum formula.}

\medskip 

In this appendix we give a description of the proof of Theorem~B. 
The proof we shall describe is obtained through an adaptation 
of the work [5] of Bruggeman and Motohashi; it also owes 
much to Lokvenec-Guleska's thesis [32], in which a significant 
generalisation of the Bruggeman-Motohashi 
sum formula is obtained. 
\par 
It is to be assumed throughout this appendix that $q_0$ is a given non-zero 
Gaussian integer, and that $\Gamma =\Gamma_0(q_0)\,$ (the Hecke congruence subgroup 
of $SL(2,{\Bbb Z}[i])$ defined in Equation~(1.1.1)). Notations 
already introduced in Section~1 of the paper remain in use: we 
shall define additional terminology as the need arises, and our useage of any such  
additional terminology shall not be limited to the subsection of the appendix 
in which the relevant definition is stated (in particular,  
a full understanding of Subsection~6.5 and Subsection~6.6 
requires some familiarity with  
terminology defined in the first four subsections of this appendix). 
 
\bigskip 

\goodbreak\centerline{\bf\S 6.1 Generalised Kloosterman sums.}

\medskip 

In this subsection we aim to justify what is stated 
in and between (1.5.11) and (1.5.12), concerning the generalised 
Kloosterman sums, and associated sets ${}^{\frak a}{\cal C}^{\frak b}$,   
defined in (1.5.8)-(1.5.10). For our proof of 
the upper bound (1.5.12) we shall need the 
analogue, due to Bruggeman and Miatello [4], of 
the `Weil-Estermann' bound for classical Kloosterman sums. We need also the 
following lemma, the proof of which is modelled  
very closely on Motohashi's work, in Section~15 of [36],  
on Hecke congruence sugroups of $SL(2,{\Bbb Z})$. 

\bigskip 

\proclaim Lemma 6.1.1. 
Let $u_1,w_1,u_2,w_2\in{\frak O}$ satisfy    
$$w_j\mid q_0\qquad{\rm and}\qquad\left( u_j , w_j\right)\sim 1\qquad\qquad 
\hbox{($j=1,2$).}\eqno(6.1.1)$$ 
Let ${\frak a}'$ and ${\frak b}'$ be the cusps of  
$\Gamma=\Gamma_0(q_0)\leq SL(2,{\frak O})$ given by   
${\frak a}'=u_1/w_1$ and ${\frak b}'=u_2/w_2$. 
Suppose moreoever that 
the scaling matrices, $g_{u_1/w_1}=g_{{\frak a}'}$ and $g_{u_2/w_2}=g_{{\frak b}'}$ 
are chosen similarly to $g_{u/w}$ in the proof of Lemma~4.2, so that one has: 
$$g_{u_j/w_j}=\varpi_{u_j/w_j}\tau_{v_j}\in SL(2,{\Bbb C})\qquad\qquad 
\hbox{($j=1,2$),}\eqno(6.1.2)$$ 
where 
$$\varpi_{u/w}=\pmatrix{u &-\tilde w\cr w &\tilde u}\in SL(2,{\frak O})\qquad\qquad 
\hbox{($u,w\in{\frak O}$, $w\neq 0$ and $(u,w)\sim 1$),}\eqno(6.1.3)$$
with $\tilde u,\tilde w$ denoting an arbitrary pair of 
Gaussian integers such that $u\tilde u+w\tilde w=1$, while 
$$\tau_v=\pmatrix{\sqrt{v} &0\cr 0 &1/\sqrt{v}}\in SL(2,{\Bbb C})\qquad\qquad 
\hbox{($v\in{\Bbb C}^*$),}\eqno(6.1.4)$$
and $v_1,v_2$ are an arbitrary pair of Gaussian integers such that 
$$v_j\sim {q_0 /w_j\over\left( q_0 /w_j\,,\,w_j\right)} 
\qquad\qquad\hbox{($j=1,2$).}\eqno(6.1.5)$$
Then $g_{{\frak a}'},g_{{\frak b}'}\in SL(2,{\Bbb C})$ are such that 
the conditions (1.1.16) and (1.1.20)-(1.1.21) are satisfied when either 
${\frak c}={\frak a}'$ or ${\frak c}={\frak b}'$. The corresponding set 
${}^{{\frak a}'}{\cal C}^{{\frak b}'}$ (for which see (1.5.8)-(1.5.9)) satisfies 
$${}^{{\frak a}'}{\cal C}^{{\frak b}'}\subseteq \sqrt{v_1 v_2}\,{\frak O} -\{ 0\}\eqno(6.1.6)$$
(with $v_1,v_2$ as in (6.1.5)) and, for 
all $m,n\in{\frak O}$ and all $C\in{\frak O}-\{ 0\}$ such that 
$\sqrt{v_1v_2}\,C\in {}^{{\frak a}'}{\cal C}^{{\frak b}'}$,   
one has    
$$S_{{\frak a}',{\frak b}'}\left( m , n ; C\sqrt{v_1 v_2}\right) 
=\qquad\sum\quad\sum_{
\!\!\!\!\!\!\!\!\!\!\!\!\!\!\!\!\!\!\!\!\!\!\!\!\!\!\!\!\!\!\!\!{\scriptstyle 
A\bmod v_1 C{\frak O},\  
D\bmod v_2 C{\frak O}\atop\scriptstyle AD\equiv 1\bmod C{\frak O}}}
\chi_{q_0}\!\left(\varpi_{{\frak a}'} g(A,D;C)\varpi_{{\frak b}'}^{-1}\right) 
{\rm e}\left( {\rm Re}\left( {mA\over v_1 C}+{nD\over v_2 C}\right)\right) ,\eqno(6.1.7)$$ 
where the matrices $\varpi_{{\frak a}'}=\varpi_{u_1/w_1}$ and 
$\varpi_{{\frak b}'}=\varpi_{u_2/w_2}$ are as in (6.1.2)-(6.1.3), while 
$$g(a,d;c)=\pmatrix{1 &a/c\cr 0 &1}
\pmatrix{0 &-1/c\cr c &0}
\pmatrix{1 &d/c\cr 0 &1} 
=\pmatrix{a &*\cr c &d}\in SL(2,{\Bbb C})\qquad\quad 
\hbox{($a,d\in{\Bbb C}$, $c\in{\Bbb C}^*$),}\eqno(6.1.8)$$
and where, for $g\in SL(2,{\Bbb C})$, 
$$\chi_{q_0}(g)=\cases{1 &if $g\in\Gamma_0(q_0)$,\cr 0 &otherwise.}\eqno(6.1.9)$$
In particular, the terms of the sum on the right-hand side of (6.1.7) are 
well defined, so that, for $a,d\in{\Bbb C}$, and all relevant choices 
of $C$, one has 
$$\chi_{q_0}\!\left(\varpi_{{\frak a}'} g(a+sv_1C,d+tv_2C;C)
\varpi_{{\frak b}'}^{-1}\right) 
=\chi_{q_0}\!\left(\varpi_{{\frak a}'} g(a,d;C)\varpi_{{\frak b}'}^{-1}\right) 
\quad\hbox{when $s,t\in{\frak O}$.}\eqno(6.1.10)$$

\medskip 

\noindent{\bf Proof.}\quad 
The scaling matrices $g_{{\frak a}'},g_{{\frak b}'}$ are chosen similarly to 
the scaling matrix $g_{u/w}$ which features in the proof of Lemma 4.2. 
That proof shows $g_{u/w}$ to be such that 
the conditions (1.1.16) and (1.1.20)-(1.1.21) are satisfied when ${\frak c}=u/w$. 
One may therefore 
infer that $g_{{\frak a}'}$ and $g_{{\frak b}'}$ are such that 
the same is true both for ${\frak c}={\frak a}'$, and for ${\frak c}={\frak b}'$. 
\par
By (6.1.2) and (6.1.3), one has 
$g_{{\frak a}'}^{-1}\Gamma g_{{\frak b}}
\subseteq\tau_{v_1}^{-1}\varpi_{u_1/w_1}^{-1}SL(2,{\frak O})\varpi_{u_2/w_2}\tau_{v_2}
\subseteq\tau_{v_1}^{-1}SL(2,{\frak O})\tau_{v_2}$. Given this, and the definitions 
(6.1.4), (1.5.8) and (1.5.9), a very short calculation suffices to  
show that (6.1.6) holds. 
\par 
We now have only to prove the results in (6.1.7)-(6.1.10). 
We may suppose that $m$ and $n$ are Gaussian  integers, and that 
$\sqrt{v_1 v_2}\,C=c\in{}^{{\frak a}'}{\cal C}^{{\frak b}'}$.  
By (6.1.6), we have $0\neq C\in{\frak O}$. 
Given that the relevant cases of (1.1.20)-(1.1.21) hold, it follows 
by the definition (1.5.10) that 
$$S_{{\frak a}',{\frak b}'}\left( m,n;c\right) 
=\sum_{\scriptstyle g\in B^{+}\backslash 
g_{{\frak a}'}^{-1}\,{}^{{\frak a}'}\!\Gamma^{{\frak b}'}\!\!(c)g_{{\frak b}'}/B^{+}\atop\scriptstyle 
g=\left( {\scriptstyle\!\!\!\!\!\alpha(g)\quad *\atop\scriptstyle\  c\quad\ \delta(g)}\right) }
{\rm e}\!\left({\rm Re}\!\left( m\,{\alpha(g)\over c}+n\,{\delta(g)\over c}\right)\right) ,\eqno(6.1.11)$$ 
where the summation is such that $g$ runs over the elements of a complete set of 
representatives, $\left\{ g^{(r)} : r\in R\right\}$  (say), for the set of double cosets  
$B^{+}\backslash 
g_{{\frak a}'}^{-1}\,{}^{{\frak a}'}\!\Gamma^{{\frak b}'}\!\!(c)g_{{\frak b}'}/B^{+} 
=B^{+}\backslash \tau_{v_1}^{-1}\varpi_{{\frak a}'}^{-1}
\,{}^{{\frak a}'}\!\Gamma^{{\frak b}'}\!\!(c)
\varpi_{{\frak b}'}\tau_{v_2}/B^{+}$.  
The corresponding set $\left\{\tau_{v_1} g^{(r)}\tau_{v_2}^{-1} : r\in R\right\}$ is, 
given (6.1.4), a complete set of representatives for the set of double cosets 
$B^{+}_1\backslash \varpi_{{\frak a}'}^{-1}
\,{}^{{\frak a}'}\!\Gamma^{{\frak b}'}\!\!(c)
\varpi_{{\frak b}'}/B^{+}_2$, where, for $j=1,2$, 
$$B^{+}_j=\tau_{v_j} B^{+}\tau_{v_j}^{-1} 
=\left\{ h\left[\sqrt{v_j}\right] n[\xi] h\left[ 1/\sqrt{v_j}\right] : \xi\in{\frak O}\right\} 
=\left\{ n\left[ v_j\xi\right] : \xi\in{\frak O}\right\} .\eqno(6.1.12)$$ 
Hence, upon rewriting the summation in (6.1.11) in terms of 
$$\tilde g=\tau_{v_1} g\tau_{v_2}^{-1} 
=\pmatrix{\alpha(g)\sqrt{v_1/v_2} &*\cr c/\sqrt{v_1 v_2} &\delta(g)\sqrt{v_2/v_1}} 
=\pmatrix{A(\tilde g) &*\cr C &D(\tilde g)}\quad\hbox{(say),}$$
one finds that 
$$S_{{\frak a}',{\frak b}'}\left( m,n;C\sqrt{v_1 v_2}\right) 
=\sum_{\scriptstyle\tilde g\in B^{+}_1\backslash 
\varpi_{{\frak a}'}^{-1}\,{}^{{\frak a}'}\!\Gamma^{{\frak b}'}\!\!\left( C\sqrt{v_1 v_2}\right)
\,\varpi_{{\frak b}'}/B^{+}_2\atop\scriptstyle 
\tilde g=\left( {\scriptstyle\!\!\!\!\!A(\tilde g)\quad *\atop\scriptstyle\  C\quad\ D(\tilde g)}\right) }
{\rm e}\!\left({\rm Re}\!\left( m\,{A(\tilde g)\over v_1 C}+n\,{D(\tilde g)\over v_2 C}\right)\right) ,
\eqno(6.1.13)$$ 
where the set $\varpi_{{\frak a}'}^{-1}\,{}^{{\frak a}'}\!\Gamma^{{\frak b}'}\!\!\left( C\sqrt{v_1 v_2}\right)
\,\varpi_{{\frak b}'}$ is (of course) invariant under both multiplication on the left by 
elements of $B^{+}_1$, and multiplication on the right by elements of $B^{+}_2$, and is  
moreover a subset of $SL(2,{\frak O})$ (for, by (1.5.8), and (6.1.1) and (6.1.3),  we have 
${}^{{\frak a}'}\!\Gamma^{{\frak b}'}\!\!\left( C\sqrt{v_1 v_2}\right)
\subset\Gamma\leq SL(2,{\frak O})$ and 
$\varpi_{{\frak a}'},\varpi_{{\frak b}'}\in SL(2,{\frak O})$). 
By the last observation, we have, in the above, $A(\tilde g),D(\tilde g)\in{\frak O}$, 
and, given that $C\neq 0$, may express $\tilde g$ as $g(A(\tilde g), D(\tilde g);C)$, 
where the notation $g(a,d;c)$ is defined in (6.1.8). Furthermore, by (6.1.8) it is evident 
that, for $z_1,z_2\in{\Bbb C}$, one has 
$$n\left[z_1\right] g(A,D;C) n\left[ z_2\right] 
=n\!\left[z_1+{A\over C}\right] g(0,0;C) n\!\left[{D\over C}+z_2\right]
=g\left( A + C z_1 , D + C z_2;C\right)\;,\eqno(6.1.14)$$
and so the 
result stated in (6.1.7)-(6.1.9) is merely a more 
explicit formulation of (6.1.13) and (6.1.12). 
\par 
Finally, given (6.1.12) and (6.1.14), it suffices for proof of the result 
(6.1.10) that we can show $\varpi_{{\frak a}'}B^{+}_1\varpi_{{\frak a}'}^{-1}$ 
and $\varpi_{{\frak b}'}B^{+}_2\varpi_{{\frak b}'}^{-1}$ to be subgroups 
of the group $\Gamma =\Gamma_0(q_0)$. This is easily achieved.  
Indeed, by the first equality of (6.1.12), 
and by (6.1.2) and (1.1.20), one has 
$$\varpi_{u_j/w_j} B^{+}_j\varpi_{u_j/w_j}^{-1} 
=\varpi_{u_j/w_j} \tau_{v_j} B^{+}\tau_{v_j}^{-1}\varpi_{u_j/w_j}^{-1} 
=g_{u_j/w_j} B^{+} g_{u_j/w_j}^{-1}
=\Gamma_{u_j/w_j}'<\Gamma\qquad\qquad\hbox{($j=1,2$),}$$
which (given that $u_1/w_1={\frak a}'$ and $u_2/w_2={\frak b}'$) is 
just what was required\quad $\blacksquare$ 

\bigskip 

\proclaim Corollary 6.1.2. Subject to the hypotheses of the above lemma, one has 
the trivial bound: 
$$\left| S_{{\frak a}',{\frak b}'}\left( m , n ; c\right)\right| 
\leq |c|^2 \left| v_1 v_2\right|\qquad\qquad 
\hbox{($m,n\in{\frak O}$ and $c\in{}^{{\frak a}'}{\cal C}^{{\frak b}'}$).}\eqno(6.1.15)$$

\medskip 

\noindent{\bf Proof.}\quad 
In the sum on the right-hand side of equation (6.1.7) there are 
$\left| ({\frak O}/C{\frak O})^*\right|\leq |C|^2$ choices for 
$A\bmod C{\frak O}$; and when that residue class is given, 
the congruence condition $AD\equiv 1\bmod C{\frak O}$ determines 
$D\bmod C{\frak O}$, so that there remain just $|v_1|^2$ choices for $A\bmod v_1C{\frak O}$ and 
$|v_2|^2$ choices for $D\bmod v_2 C{\frak O}$. Hence the sum appearing in (6.1.7) contains 
at most $|Cv_1v_2|^2$ non-zero terms. The bound (6.1.15) 
therefore follows
from (6.1.6) and the case $C=c/\sqrt{v_1 v_2}$ of (6.1.7)-(6.1.9), by 
way of the triangle inequality\ $\blacksquare$  

\bigskip 

The next lemma is the Gaussian integer analogue of Equation~(15.14) of [36]. 
Before stating the lemma, it is helpful to clarify that it 
assumes a more special choice of all the relevant scaling matrices than 
was the case in Lemma 6.1.1. To be precise, although it is still 
to be assumed that $g_{u/w}=\varpi_{u/w}\tau_{v}$, with 
$\varpi_{u/w}$ and $\tau_v$ as in (6.1.3)-(6.1.5), it is now 
also to be supposed that $\varpi_{u/w}$ is chosen in such a way that 
the relevant Gaussian integer $\tilde u$ 
in (6.1.3) satisfies 
$$u\tilde u\equiv 1\bmod q_0 {\frak O}\;.\eqno(6.1.16)$$ 
Since this is only possible when $u$ is coprime to $q_0$, 
the choice of cusps is also more restricted than is 
the case in Lemma~6.1.1. 

\bigskip 

\proclaim Lemma 6.1.3. Let the hypotheses of Lemma 6.1.1 be satisfied. 
Suppose moreover that $(u_1 u_2 , q_0 )\sim 1$, and that 
the scaling matrices $g_{u_1/w_1}=g_{{\frak a}'}$ and 
$g_{u_2/w_2}=g_{{\frak b}'}$ satisfy the 
additional constraint imposed in (6.1.16). 
Let $C$ be a non-zero Gaussian integer such that 
$C\sqrt{v_1 v_2}\in{}^{{\frak a}'}{\cal C}^{{\frak b}'}$.  
Then, for some $C_{q_0},C_{q_0}',\tilde C_{q_0}\in{\frak O}$, one has 
$$C_{q_0} C_{q_0}'=C\;,\qquad\quad C_{q_0}'\sim\left( C , q_0^{\infty}\right)\;,\quad\qquad 
C_{q_0} \tilde C_{q_0}\equiv 1\bmod \left[\left[v_1 , v_2\right] C_{q_0}' , q_0\right] 
{\frak O}\eqno(6.1.17)$$
(where $[r,s]$ denotes an arbitrary least common multiple of $r$ and $s$, 
and where the notation `$q_0^{\infty}$' signifies 
that one has $C_{q_0}'\sim (C, q_0^N)$ for all sufficiently large $N\in{\Bbb N}$) and also: 
$${}^{{\frak a}''}{\cal C}^{{\frak b}''}\ni C_{q_0}'\sqrt{v_1 v_2}\eqno(6.1.18)$$
and, for all $m,n\in{\frak O}$,   
$$S_{{\frak a}',{\frak b}'}\left( m , n ; C\sqrt{v_1 v_2}\right) 
=S_{{\frak a}'',{\frak b}''}\!\left( \tilde C_{q_0} m , \tilde C_{q_0} n ; C_{q_0}'\sqrt{v_1 v_2}\right) 
S\!\left( \left( C_{q_0}' v_1\right)^{\!*} m , \left( C_{q_0}' v_2\right)^{\!*} n ; C_{q_0}\right)\;,\eqno(6.1.19)$$
where 
$${\frak a}''=\tilde C_{q_0} u_1/w_1\;,\qquad\qquad 
{\frak b}''=\tilde C_{q_0} u_2/w_2\eqno(6.1.20)$$ 
and the associated scaling matrices are 
$g_{{\frak a}''}=\varpi_{\tilde C_{q_0} u_1/w_1}\tau_{v_1}$ and 
$g_{{\frak b}''}=\varpi_{\tilde C_{q_0} u_2/w_2}\tau_{v_2}$ 
(with $\tau_v$, $v_1$ and $v_2$ as in (6.1.4)-(6.1.5), and with 
$\varpi_{u/w}$ as indicated in (6.1.3) and (6.1.16)), 
while $S(a,b;c)$ denotes the `simple Kloosterman sum' defined 
in (2.16), and the $*$-notation has the significance explained 
in our subsection on notation 
(under the sub-heading `Other Number-Theoretic Notation'). 

\medskip 

\noindent{\bf Proof.}\quad 
Since the proof involves calculations very similar 
to those of Section~15 of [36], we only sketch the main points. 

Note firstly 
that, since ${\frak O}$ is a unique factorisation domain, 
we can certainly find a Gaussian integer $C_{q_0}'$ satisfying 
the condition $C_{q_0}'\sim (C,q_0^{\infty})$ in (6.1.17). 
Then $C_{q_0}'\mid C$, and the Gaussian integer $C_{q_0}=C/C_{q_0}'$ is 
coprime to each of $C_{q_0}'$, $q_0$, $v_1$ and $v_2$ 
(the latter two being factors of $q_0$, by virtue of (6.1.5) and (6.1.1)). 
Hence there exist  
$C_{q_0},C_{q_0}',\tilde C_{q_0}\in{\frak O}$ such that 
all the conditions in (6.1.17) are satisfied.  

Let $C_{q_0},C_{q_0}',\tilde C_{q_0}\in{\frak O}$ satisfy the conditions in (6.1.17), 
and let the cusps ${\frak a}'', {\frak b}''$ be given by (6.1.20). 
Since $(u_1 u_2 , q_0)\sim 1$ (by hypothesis), $(\tilde C_{q_0} , q_0)\sim 1$ (by (6.1.17)) 
and $w_1,w_2\mid q_0$ (by (6.1.1)), we may choose the scaling matrices 
$g_{{\frak a}''}$ and $g_{{\frak b}''}$ to be as described below (6.1.20). 
Then, given the hypotheses of the lemma concerning $g_{{\frak a}'}$ 
and $g_{{\frak b}'}$, it may be shown by a calculation that, 
with regard to the result (6.1.7) of Lemma 6.1.1,  one has: 
$$\chi_{q_0}\!\left(\varpi_{{\frak a}'} g(A,D;C)\varpi_{{\frak b}'}^{-1}\right) 
=\chi_{q_0}\!\left(\varpi_{{\frak a}''} g\left( A,D;C_{q_0}'\right)\varpi_{{\frak b}''}^{-1}\right)\quad\ 
\hbox{for $\ A,D\in{\frak O}\ $ with $\ AD\equiv 1\bmod C{\frak O}$}\eqno(6.1.21)$$ 
(where $\varpi_{u/w}$, $g(a,d;c)$ and $\chi_{q_0}$ are given by (6.1.3), (6.1.8) and (6.1.9)). 
We remark that the full force of the congruence in (6.1.17) is not 
required in the above: it would suffice there to have 
just $C_{q_0} \tilde C_{q_0}\equiv 1\bmod q_0 {\frak O}$. 
\par 
Given that ${}^{{\frak a}'}{\cal C}^{{\frak b}'}\ni C\sqrt{v_1 v_2}$, 
and given the assumptions made concerning scaling matrices, 
it may be deduced from (6.1.21) that  
${}^{{\frak a}''}{\cal C}^{{\frak b}''}\ni C_{q_0}'\sqrt{v_1 v_2}$, so that one   
may apply Lemma~6.1.1 with ${\frak a}''$, ${\frak b}''$ and $C_{q_0}'$ 
substituted for ${\frak a}'$, ${\frak b}'$ and $C$, respectively (and with 
no change in the values of $v_1$ and $v_2$). By making these 
substitutions in (6.1.10), one finds that the right-hand side of 
the equation in (6.1.21) is a function of the residue class of 
$A\bmod v_1 C_{q_0}'{\frak O}$ and the residue class of $D\bmod v_2 C_{q_0}'{\frak O}$. 
Therefore, and since $(v_1 C_{q_0}' , C_{q_0})\sim 1$ and 
$(v_2 C_{q_0}' , C_{q_0})\sim 1$, 
one may apply the Chinese Remainder Theorem to deduce from  
(6.1.7) and (6.1.21) that 
$$S_{{\frak a}',{\frak b}'}\left( m , n ; C\sqrt{v_1 v_2}\right) 
=XY\;,$$
where   
$$X=\qquad\sum\quad\sum_{
\!\!\!\!\!\!\!\!\!\!\!\!\!\!\!\!\!\!\!\!\!\!\!\!\!\!\!\!\!\!\!\!{\scriptstyle 
A\bmod v_1 C_{q_0}'{\frak O},\  
D\bmod v_2 C_{q_0}'{\frak O}\atop\scriptstyle A D\equiv 1\bmod C_{q_0}'{\frak O}}}
\chi_{q_0}\!\left(\varpi_{{\frak a}''} g\left( A,D;C_{q_0}'\right)
\varpi_{{\frak b}''}^{-1}\right) 
{\rm e}\!\left(\!{\rm Re}\!\left( {\tilde C_{q_0} m A\over v_1 C_{q_0}'}
+{\tilde C_{q_0} n D\over v_2 C_{q_0}'}\right)\!\right)$$ 
and  
$$Y=\qquad\sum\quad\sum_{
\!\!\!\!\!\!\!\!\!\!\!\!\!\!\!\!\!\!\!\!\!\!\!\!\!\!\!\!\!\!\!\!{\scriptstyle 
A\bmod C_{q_0} {\frak O} ,\  
D\bmod C_{q_0} {\frak O}\atop\scriptstyle A D\equiv 1\bmod C_{q_0} {\frak O}}}
{\rm e}\left( {\rm Re}\left( {\left( v_1 C_{q_0}'\right)^{*} m A\over C_{q_0}}
+{\left( v_2 C_{q_0}'\right)^{*} n D\over C_{q_0}}\right)\right)
=S\!\left( \left( C_{q_0}' v_1\right)^{\!*} m , 
\left( C_{q_0}' v_2\right)^{\!*} n ; C_{q_0}\right)\;.$$
The result (6.1.19) follows:  
for, by Lemma~6.1.1 (applied with 
${\frak a}''$ and ${\frak b}''$ substituted for 
${\frak a}'$ and ${\frak b}'$, respectively), one has 
$X=S_{{\frak a}'',{\frak b}''}\bigl( \tilde C_{q_0} m , 
\tilde C_{q_0} n ; C_{q_0}'\sqrt{v_1 v_2}\bigr)$\ $\,\blacksquare$ 

\bigskip 

\proclaim Corollary 6.1.4. Let the combined hypotheses of Lemma~6.1.1 and 
Lemma~6.1.3 be satisfied. Then, for all $m,n\in{\frak O}$, one has 
$$\left| S_{{\frak a}',{\frak b}'}\left( m , n ; C\sqrt{v_1v_2}\right)\right| 
\leq 2^{3/2}\tau(C)\left|(m,n,C)C\left( C , q_0^{\infty}\right) v_1^2 v_2^2\right|\;,
\eqno(6.1.22)$$
where $\tau(k)$ equals the number of Gaussian integer divisors of $k$, 
and `$q_0^{\infty}$' has the same meaning as in (6.1.17). 

\medskip 

\noindent{\bf Proof.}\quad 
By (6.1.20), one has ${\frak a}''=u_1'/w_1$ and ${\frak b}''=u_2'/w_2$, 
where $u_j'=\tilde C_{q_0} u_j\in{\frak O}$ ($j=1,2$), so that, by (6.1.1) and (6.1.17), 
$(u_j',w_j)\sim 1$ for $j=1,2$. Hence (and given the result in (6.1.18)) 
Corollary~6.1.2 may be applied with 
${\frak a}''$ and ${\frak b}''$ substituted for ${\frak a}'$ and 
${\frak b}'$ (respectively), and with $v_1$ and $v_2$ unchanged. 
Consequently one has the upper bound 
$$\left| S_{{\frak a}'',{\frak b}''}\left( \tilde C_{q_0} m , 
\tilde C_{q_0} n ; C_{q_0}'\sqrt{v_1v_2}\right)\right| 
\leq\left| C_{q_0}'\sqrt{v_1v_2}\right|^2 \left| v_1v_2\right| 
=\left| C_{q_0}' v_1v_2\right|^2\;.\eqno(6.1.23)$$
By the result (2.18) of Lemma 2.4 (the `Weil-Estermann' bound obtained by Bruggeman and 
Miatello in [4]), we have also the upper bound 
$$\left| 
S\!\left( \left( C_{q_0}' v_1\right)^{\!*}\!m , \left( C_{q_0}' v_2\right)^{\!*}\!n ; C_{q_0}\right)
\right| 
\leq 2^{3/2}\tau\!\left( C_{q_0}\right)\left| \left( m , n , C_{q_0}\right) C_{q_0}\right|\eqno(6.1.24)$$ 
for the second factor on the right-hand side of equation (6.1.19) 
(note that we have used here the fact that, in equation (6.1.19), 
both $(C_{q_0}' v_1)^*$ and $(C_{q_0}' v_2)^*$ are, by definition, coprime to $C_{q_0}$). 
By (6.1.24), (6.1.23), (6.1.19) and (6.1.17), 
the bound (6.1.22) follows\ $\blacksquare$ 

\bigskip 

The next (and final) lemma in this subsection contains the results stated 
in (1.5.11) and (1.5.12). 
In those results, and in the lemma, we make use of the notation `$m_{\frak c}$' 
introduced below equation (1.1.22). Hence, for each cusp ${\frak c}$ 
of the group $\Gamma =\Gamma_0(q_0)\leq SL(2,{\frak O})$, we take 
$m_{\frak c}$ to signify an arbitrary Gaussian integer satisfying 
$$m_{\frak c}\sim
\cases{ q_0 /(w^2 , q_0) &if ${\frak c}=u/w$ with $u,w\in{\frak O}$, $w\neq 0$ and $(u,w)\sim 1$; \cr 
1 &if ${\frak c}=\infty$.}\eqno(6.1.25)$$
 (see (1.5.8)-(1.5.10) for the definitions of 
the relevant generalised Kloosterman sums 
$S_{{\frak a},{\frak b}}(m,n;c)$ and the associated sets  
${}^{\frak a}{\cal C}^{\frak b}$). It follows from (6.1.25) 
that the ideal $m_{\frak c}{\frak O}$ depends only on  
the $\Gamma$-equivalence class of the cusp ${\frak c}$ (we skip the 
easy proof of this). For the notations `$\tau(k)$' and `$q_0^{\infty}$', see below 
(6.1.17) or (6.1.22). 

\proclaim Lemma 6.1.5. Let $m,n\in{\frak O}$. Let ${\frak a}$ and ${\frak b}$ 
be cusps of $\Gamma =\Gamma_0(q_0)$, and let the associated scaling matrices 
$g_{\frak a},g_{\frak b}\in SL(2,{\Bbb C})$ be such that the conditions 
(1.1.16) and (1.1.20)-(1.1.21) are satisfied for ${\frak c}\in\{ {\frak a} , {\frak b}\}$. 
Then one has, for some $\epsilon\in{\frak O}^*$, both the relation 
$${}^{\frak a}{\cal C}^{\frak b}\subseteq\sqrt{\epsilon m_{\frak a} m_{\frak b}}\ {\frak O} 
-\{ 0\}\eqno(6.1.26)$$ 
and, for all $c=\sqrt{\epsilon m_{\frak a} m_{\frak b}}\ C\in{}^{\frak a}{\cal C}^{\frak b}$, 
the upper bounds   
$$\left| S_{{\frak a},{\frak b}}(m,n;c)\right| 
\leq 2^{3/2}\tau(C)\left| (m,n,C)\,C \left( C , q_0^{\infty}\right) m_{\frak a}^2 m_{\frak b}^2\right|  
\eqno(6.1.27)$$
and 
$$\left| S_{{\frak a},{\frak b}}(m,n;c)\right| 
\leq \left| C m_{\frak a} m_{\frak b}\right|^2\;.\eqno(6.1.28)$$

\medskip 

\noindent{\bf Proof.}\quad 
By the result (2.4) of Lemma 2.2, there exist $u_1,w_1,u_2,w_2\in{\frak O}$ satisfying 
(6.1.1), and such that 
$${\frak a}\sim^{\!\!\!\!\Gamma} u_1/w_1\qquad\qquad{\rm and}\qquad\quad 
{\frak b}\sim^{\!\!\!\!\Gamma} u_2/w_2\;.\eqno(6.1.29)$$
It may at the same time be assumed (in the above) that the hypothesis 
$(u_1 u_2 , q_0)\sim 1$ of Lemma~6.1.3 is satisfied: note that this does not 
entail any loss of generality, for one has 
$u_j/w_j\sim (u_j+k w_j)/w_j$, for $j=1,2$ and all $k\in{\frak O}$, and 
the coprimality conditions in (6.1.1) imply that both of the sets   
$\{ u_1+kw_2 : k\in{\frak O}\}$ and $\{ u_2+kw_2 : k\in{\frak O}\}$ 
do contain elements coprime to $q_0$. 

We now write $u_1/w_1={\frak a}'$ and $u_2/w_2={\frak b}'$, and take 
the associated scaling matrices, 
$g_{{\frak a}'}=g_{u_1/w_1}$ and $g_{{\frak b}'}=g_{u_2/w_2}$,  
to be of the form indicated in (6.1.2)-(6.1.5) of Lemma~6.1.1, so that 
$g_{{\frak a}'}=\varpi_{u_1/w_1}\tau_{v_1}$ and 
$g_{{\frak b}'}=\varpi_{u_2/w_2}\tau_{v_2}$
where $v_1$, $v_2$, $\tau_{v_1}$ and $\tau_{v_2}$ are as in (6.1.5) and 
(6.1.4), while $\varpi_{u_1/w_1}$ and $\varpi_{u_2/w_2}$ are as indicated by 
(6.1.3). Since $(u_1 u_2 , q_0)\sim 1$, and since $w_1,w_2\mid q_0$, we may suppose  
moreover that, for $j=1,2$, one has   
$$\varpi_{u_j/w_j}=\pmatrix{u_j &-\tilde w_j\cr w_j &\tilde u_j}
\in SL(2,{\frak O})\;,$$ 
with $\tilde u_j\in{\frak O}$ such that $u_j \tilde u_j\equiv 1\bmod q_0 {\frak O}$ 
(i.e. one can find $r_j,s_j\in{\frak O}$ 
such that $u_j r_j + q_0 s_j=1$, and so may take 
$\tilde u_j=r_j$ and $\tilde w_j=(q_0 /w_j)s_j$ in the above). 
Note that (as is stated in Lemma~6.1.1) such a choice of 
the scaling matrices $g_{{\frak a}'}$ and $g_{{\frak b}'}$ ensures that 
the conditions (1.1.16) and (1.1.20)-(1.1.21) are satisfied 
when ${\frak c}\in\{ {\frak a}' , {\frak b}'\}$. Therefore, given 
the relations in (6.1.29), it follows by Lemma~2.1, (2.2) and (2.3), that, 
for some pair of units $\epsilon_1,\epsilon_2\in{\frak O}^*$, 
one has both 
$${}^{\frak a}{\cal C}^{\frak b}=\sqrt{\epsilon_1\epsilon_2}
\ {}^{{\frak a}'}{\cal C}^{{\frak b}'}\eqno(6.1.30)$$ 
and, for all $c\in{}^{\frak a}{\cal C}^{\frak b}$, 
$$\left| S_{{\frak a},{\frak b}}(m,n;c)\right| 
=\left| S_{{\frak a}',{\frak b}'}\!\left( 
\overline{\epsilon_1}\,m , \overline{\epsilon_2}\,n ; c/\sqrt{\epsilon_1\epsilon_2}\right) 
\right|\eqno(6.1.31)$$
(note that, by (1.5.13), the choice of square roots in the above 
is immaterial). 

Now, by (6.1.30) and the result (6.1.6) of Lemma 6.1.1, we have 
$${}^{\frak a}{\cal C}^{\frak b}\subseteq\sqrt{\epsilon_1\epsilon_2 v_1 v_2}
\ {\frak O}-\{ 0\}\;,\eqno(6.1.32)$$ 
where the Gaussian integers $v_1,v_2$ satisfy (6.1.5), and so (see 
(6.1.25)) are such that $v_1\sim m_{{\frak a}'}$ and 
$v_2\sim m_{{\frak b}'}$. Given the relations in (6.1.29), we have moreover 
$m_{\frak a}\sim m_{{\frak a}'}\sim v_1$ and $m_{\frak b}\sim m_{{\frak b}'}\sim v_2$, 
and so, by (6.1.32), we obtain the result (6.1.26) with 
$\epsilon =(v_1/m_{\frak a})(v_2/m_{\frak b})\epsilon_1\epsilon_2\in{\frak O}^*$. 
Moreover, with $\epsilon$ as just indicated, it follows from 
(6.1.30) and (6.1.32) that, if $c,C\in{\Bbb C}^*$ are 
such that $c\in{}^{\frak a}{\cal C}^{\frak b}$ 
and $\sqrt{\epsilon m_{\frak a} m_{\frak b}}\,C=c$, then one has 
$0\neq C\in{\frak O}$ and $\sqrt{v_1 v_2}\,C=c/\sqrt{\epsilon_1\epsilon_2}
\in{}^{{\frak a}'}{\cal C}^{{\frak b}'}$, and so obtains, 
by (6.1.31) and Corollary~6.1.4,  
$$\left| S_{{\frak a},{\frak b}}(m,n;c)\right|
=\left| S_{{\frak a}',{\frak b}'}\!\left( 
\overline{\epsilon_1}\,m , \overline{\epsilon_2}\,n ; \sqrt{v_1 v_2}\ C\right) 
\right|
\leq 2^{3/2}\tau(C)\left|\left(\overline{\epsilon_1}\,m,\overline{\epsilon_2}\,n,C\right)
\,C\left( C , q_0^{\infty}\right) v_1^2 v_2^2\right|\;.$$
Since we have here $\epsilon_1,\epsilon_2\in{\frak O}^*$, 
$v_1\sim m_{\frak a}$ and $v_2\sim m_{\frak b}$, the result 
(6.1.27) follows immediately. 
The proof of (6.1.28) is similar, differing only in that one   
uses Corollary~6.1.2 
in place of the `Weil-Estermann' bound, (6.1.22) 
\quad$\blacksquare$

\bigskip 

\centerline{\bf\S 6.2 Poincar\'{e} series.} 

\medskip 

\noindent{\bf Convergence, continuity and a lemma on inner products.}\quad  
For $\omega\in{\frak O}$, we define  
$C^0(N\backslash G, \omega)$ to be the space of continuous 
functions $h : G\rightarrow{\Bbb C}$ with the property that, 
for $n\in N$ and $g\in G$,  
one has $h(ng)=\psi_{\omega}(n)h(g)\,$ 
(where $\psi_{\omega}$ is the character of $N$ defined 
in (1.4.3)). We remark that for each $\omega\in{\frak O}$ one has the relation  
$C^0(N\backslash G ,\omega)
\subset C^0(B^{+}\backslash G)$, where 
$$C^0\!\left( B^{+}\backslash G\right) 
=\left\{\phi\in C^0(G)\,:\,\phi(bg)=\phi(g)\ \,{\rm for}
\ \,b\in B^{+},\,g\in G\right\} .$$
\par 
Assume now (and for the remainder of this subsection)  
that ${\frak a}$ is some cusp 
of $\Gamma$, that ${\omega}\in{\frak O}$, and that 
one has $f=f_{\omega}\in C^0(N\backslash G ,\omega)$. 
Then the 
conditions (1.1.16) and (1.1.20)-(1.1.21), for ${\frak c}={\frak a}$,  
suffice to ensure that if, for all $g\in G$, the series on the right-hand 
side of the equation (1.5.4) converges absolutely 
then (1.5.4) defines a function 
$P^{\frak a} f : G\rightarrow{\Bbb C}$ which is $\Gamma$-automorphic. 

In determining sufficient conditions for the 
absolute convergence of the Poincar\'e series $P^{\frak a}f$ 
we shall generalise the approach taken in Section~7.1 of [32], where only the 
case ${\frak a}=\infty$ is considered. As a matter of notational 
convenience, we begin by defining  
$$\rho(na[r]k)=r\qquad\quad\hbox{for $n\in N$, $r>0$ and $k\in K$,}\eqno(6.2.1)$$
so that, for each $\nu\in{\Bbb C}$, one has 
$$\left( \rho(g)\right)^{1+\nu}=\varphi_{0,0}(\nu ,0)(g)\qquad\qquad\hbox{($g\in G$),}
\eqno(6.2.2)$$
where $\varphi_{\ell ,q}(\nu ,p) : G\rightarrow {\Bbb C}$ is the function defined 
in (1.3.2). The following lemma will prove useful. 

\bigskip 

\proclaim Lemma 6.2.1. Let ${\frak a}$ and ${\frak b}$ be cusps of $\Gamma$,  
and let $g_{\frak a},g_{\frak b}\in G=SL(2,{\Bbb C})$ be such that the conditions 
(1.1.16) and (1.1.20)-(1.1.21) are satisfied for ${\frak c}\in\{ {\frak a} , {\frak b}\}$. 
Then the family of sets 
$( {}^{\frak a}\Gamma^{\frak b}(c) )_{c\in {}^{\frak a}{\cal C}^{\frak b}\cup\{ 0\}}$,  
defined by (1.5.8) and (1.5.9), is a partitioning 
of the set of elements of $\Gamma$.
For $c\in {}^{\frak a}{\cal C}^{\frak b}\cup\{ 0\}$, 
$\gamma\in {}^{\frak a}\Gamma^{\frak b}(c)$ and $g\in G$, one has 
$$\rho\left( g_{\frak a}^{-1}\gamma g_{\frak b} g\right) =\rho(g)\qquad\quad\hbox{if $\,c=0$;}
\eqno(6.2.3)$$
$$\rho\left(g_{\frak a}^{-1}\gamma g_{\frak b} g\right) \rho(g)\leq |c|^{-2}\leq\left| m_{\frak a} m_{\frak b}\right|^{-1}\qquad\quad 
\hbox{if $\,0\neq c\in{}^{\frak a}{\cal C}^{\frak b}$,}\eqno(6.2.4)$$
where $|m_{\frak c}|^2\in{\Bbb N}$ is the `width' of the cusp ${\frak c}$ 
(as defined below (1.1.22)). 

\medskip 

\noindent{\bf Proof.}\quad 
Given that $\{ g_{\frak a},g_{\frak b}\}\cup\Gamma\subset G=SL(2,{\Bbb C})$, 
the assertion concerning the partitioning of $\Gamma$ is an immediate  
corollary of the definitions in (1.5.8) and (1.5.9). 

By (1.5.8), one has $\gamma\in {}^{\frak a}\Gamma^{\frak b}(0)$ if and 
only if $\gamma\in\Gamma$ and $g_{\frak a}^{-1}\gamma g_{\frak b}\infty =\infty$. 
Since the latter relation is equivalent to the relation $\gamma{\frak b}={\frak a}$, 
it therefore follows by Equation (2.1) of Lemma~2.1 that if  
$\gamma\in {}^{\frak a}\Gamma^{\frak b}(0)$ then 
$g_{\frak a}^{-1}\gamma g_{\frak b}=h[\eta]n[\beta]$, for some  
$\eta,\beta\in{\Bbb C}$ with $\eta^2\in{\frak O}^{*}$. 
The result (6.2.3) follows: for, by (1.1.4), (1.1.3), 
(1.1.9) and (6.2.1), one has 
$\rho( h[u]n[\beta]g)=|u|^2\rho(n[\beta]g)=|u|^2\rho(g)$ 
for $u\in{\Bbb C}^{*}$, $\beta\in{\Bbb C}$, $g\in G$. 

Supposing now that $g\in G$, and that 
$\gamma\in{}^{\frak a}\Gamma^{\frak b}(c)$ for some non-zero $c\in{\Bbb C}$, 
it follows by (1.1.4) and the definition (1.5.8) that 
$\rho( g_{\frak a}^{-1}\gamma g_{\frak b} g)
=\rho(g)/(|cz+d|^2+|c|^2\rho^2(g))$, where $z\in{\Bbb C}$ is the Iwasawa 
`$z$-coordinate' of $g$, while $-d/c=g_{\frak b}^{-1}\gamma^{-1}g_{\frak a}\infty\in{\Bbb C}$.  
This shows, since $|cz+d|^2\geq 0$ 
(while $\rho(g)$ and $|c|$ are strictly positive), 
that one has 
$\rho( g_{\frak a}^{-1}\gamma g_{\frak b} g)\leq 1/(|c|^2\rho(g))\in(0,\infty)$, 
which implies the first inequality in (6.2.4). What remains of (6.2.4) is 
immediate from (6.1.25) and the result (6.1.26) of Lemma 6.1.5\quad$\blacksquare$ 

\bigskip 

We now suppose that, for some $\sigma_0 >1$ and some $R_0>0$ (both of which may 
depend on $f_{\omega}$), one has 
$$f_{\omega}(g)\ll_{f_{\omega}}\,\left(\rho(g)\right)^{1+\sigma_0}\qquad\qquad 
\hbox{for all $g\in G$ such that $\rho(g)\leq R_0$.}\eqno(6.2.5)$$
Note that, since $|f_{\omega}|\in C^0(N\backslash G,0)$, 
it is implied by the hypothesis 
(6.2.5) that if $R>0$ then 
$$f_{\omega}(g)\ll_{f_{\omega},\sigma_0,R_0,R}\ \left(\rho(g)\right)^{1+\sigma_0}\qquad\qquad 
\hbox{for all $g\in G$ such that $\rho(g)\leq R$.}\eqno(6.2.6)$$

\bigskip 

\proclaim Lemma 6.2.2. Let ${\frak a}$ be a cusp of $\Gamma$, let 
$g_{\frak a}\in G$ be such that (1.1.16) and (1.1.20)-(1.1.21) hold for 
${\frak c}={\frak a}$, let 
$\omega\in{\frak O}$, $\sigma_0>1$ and $R_0>0$, and let 
$f_{\omega}\in C^0(N\backslash G,\omega)$ be such that the condition (6.2.5) is satisfied. 
Then, when $J\subseteq G$ is compact, the series 
$\sum_{\Gamma_{\frak a}'\backslash\Gamma} |f_{\omega}(g_{\frak a}^{-1}\gamma g)|$ 
is uniformly convergent for all $g\in J$. 

\medskip 

\noindent{\bf Proof.}\quad 
Let $J$ be a compact subset of the set of elements of $G$. 
By (6.1.25) and 
the case ${\frak b}=\infty$, $g_{\frak b}=h[1]$ of Lemma 6.2.1, 
it follows that, for $g\in G$ and $\gamma\in\Gamma$, one has 
$\rho( g_{\frak a}^{-1}\gamma g)\leq\exp(|\log\rho(g)|)$. 
Since it is moreover the case that the compactness of $J$ implies the existence of   
$\max\{\;\!\exp(|\log\rho(g)|) : g\in J\}$, and since the 
hypothesis (6.2.5) implies that we have (6.2.6) when (in particular) 
$R$ is equal to this maximum, we therefore find that 
$$f_{\omega}\left( g_{\frak a}^{-1}\gamma g\right) 
\ll_{f_{\omega},\sigma_0,R_0,J}\ (\rho(g))^{1+\sigma_0}
=\varphi_{0,0}\left(\sigma_0,0\right)(g)\qquad\qquad 
\hbox{for $g\in J$, $\gamma\in\Gamma$}\eqno(6.2.7)$$
(the final equality having been noted in (6.2.2)). 
It follows that the series 
$\sum_{\Gamma_{\frak a}'\backslash\Gamma} |f_{\omega}(g_{\frak a}^{-1}\gamma g)|$ 
converges uniformly, for $g\in J$, if and only if the same is true 
of the series $E_{0,0}^{\frak a}(\sigma_0,0)(g)$ 
defined in (1.8.1). The latter series is an instance of 
the Eisenstein series extensively discussed in 
Chapter~3 of~[11], and would there be designated by 
the notation `$E_{g_{\frak a}^{-1}}(gj,\sigma_0)$' (in which 
`$j$' signifies the point $(0,1)\in{\Bbb H}_3$). 
Given (6.2.2), it is shown by the case `$\eta =\infty$, $B=I$' of 
Proposition~3.2.3 of~[11], and by [11], Definition~3.1.2, Proposition~3.1.3 and 
Corollary~3.1.6 (which imply that 
the `abscissa of convergence' of 
the cofinite subgroup $\Gamma\leq SL(2,{\Bbb C})$ equals $1$),     
that if $\alpha >0$ and $\beta >1$ then the series 
$$\sum_{\gamma\in\Gamma_{\frak a}'\backslash\Gamma} 
\left( {\rho\left( g_{\frak a}^{-1}\gamma g\right)\over 
\rho(g)}\right)^{1+\nu}
=\sum_{\gamma\in\Gamma_{\frak a}'\backslash\Gamma} 
\left(\rho(g)\right)^{-1-\nu} 
\varphi_{0,0}\left(\nu,0\right)\left( g_{\frak a}^{-1}\gamma g\right)$$
converges uniformly for all $(g,\nu)\in G\times{\Bbb C}$ 
satisfying both $\rho(g)\geq\alpha$ and ${\rm Re}(\nu)\geq\beta$; in particular, 
since $\sigma_0>1$, and since 
the compactness of $J$ implies that one has 
$\min\{\rho(g) : g\in J\} >0$, 
this series converges uniformly for all $(g,\nu)\in J\times \{\sigma_0\}$. 
It follows that the series 
$E_{0,0}^{\frak a}(\sigma_0,0)(g)
=\sum_{\gamma\in\Gamma_{\frak a}'\backslash\Gamma} 
[\Gamma_{\frak a}:\Gamma_{\frak a}']^{-1}
\varphi_{0,0}(\sigma_0,0)( g_{\frak a}^{-1}\gamma g)$ 
is (likewise) uniformly convergent for $g\in J$: 
for one has  $\min\{ (\rho(g))^{-1-\sigma_0} : g\in J\}>0$
(by virtue of the continuity of the function 
$\rho : G\rightarrow (0,\infty)$ and the compactness 
of $J$), and the factor $[\Gamma_{\frak a}:\Gamma_{\frak a}']^{-1}\in (0,1/2]$ 
is independent of $g$. 
This, with (6.2.7), completes the proof of the lemma\quad$\blacksquare$ 

\bigskip 

\proclaim Corollary~6.2.3. If the hypotheses of the above lemma are satisfied, 
and if $f=f_{\omega}$, then
the equation (1.5.4) defines a $\Gamma$-automorphic function 
$P^{\frak a} f : G\rightarrow{\Bbb C}$ which is continuous on $G$. 

\medskip 

\noindent{\bf Proof.}\quad 
By the hypothesis that $f_{\omega}\in C^0(N\backslash G,\omega)$,  
each term $f_{\omega}(g_{\frak a}^{-1}\gamma g)$ in the the series in (1.5.4) 
is a continuous function of $g$. Therefore, given that 
$G$ is locally compact, 
it follows by the uniform convergence established in the lemma 
that the series in (1.5.4) is convergent for all $g\in G$, and has a sum that is  
a continuous function of $g$. The absolute convergence of the series in (1.5.4) implies 
that all rearrangements of that series have the same sum, and so the $\Gamma$-automorphicity 
of $P f_{\omega}$ may be seen to follow by observing that, when $\tau_0\in\Gamma$ is given,  
the mapping $\Gamma_{\frak a}'\gamma\mapsto\Gamma_{\frak a}'\gamma\tau_0$ 
is a permutation on the set of right cosets of $\Gamma_{\frak a}'$ in $\Gamma$  
(this implying that if    
${\cal X}$ is 
a complete set of right coset representatives of $\Gamma_{\frak a}'$ in $\Gamma$ then 
so too is ${\cal X}\tau$, if $\tau\in\Gamma$)\quad$\blacksquare$ 

\bigskip 

The next lemma is of a well-known type: it is, in particular, a  
minor extension of Lemma~7.3.1 of [32] (which dealt only with 
Poincar\'e series associated with the cusp $\infty$). 
We shall later obtain an extension of it (Lemma~6.6.2) through which   
the key part that certain Poincar\'{e} series have to play 
in the proof of the sum formula (Theorem B) is mediated. 
Before stating the lemma we clarify that henceforth 
$C^0(\Gamma\backslash G)$ denotes the space of 
those functions $F : G\rightarrow{\Bbb C}$ that are both 
continuous and $\Gamma$-automorphic, while $L^1(\Gamma\backslash G)$ 
will denote the space of those measurable and $\Gamma$-automorphic functions 
$f : G\rightarrow{\Bbb C}$ that satisfy $\int_{\Gamma\backslash G} |f|\,{\rm d}g<\infty$. 

\bigskip 

\proclaim Lemma 6.2.4. Let the hypotheses of Lemma~6.2.2 be satisfied. 
Suppose, moreover, that $\phi\in C^0(\Gamma\backslash G)$ is such 
that $(P^{\frak a}|f_{\omega}|)\cdot\phi\in L^1(\Gamma\backslash G)$. 
Then 
$$\left[\Gamma_{\frak a} : \Gamma_{\frak a}'\right] 
\left\langle P^{\frak a} f_{\omega} , \phi\right\rangle_{\Gamma\backslash G} 
=\left\langle f_{\omega} , F^{\frak a}_{\omega}\phi\right\rangle_{N\backslash G}\;,
\eqno(6.2.8)$$ 
where, for $f,F\in C^0(N\backslash G , \omega)$,  
$$\langle f , F\rangle_{N\backslash G} 
=\int_{N\backslash G} f(g)\,\overline{F(g)}\,{\rm d}\dot{g} 
=\int_A\int_K {f(ak)\,\overline{F(ak)}\,{\rm d}k\,{\rm d}a\over (\rho(a))^2}  
=\int\limits_0^{\infty} \int_K 
{f\left( a[r] k\right)\,\overline{F\left( a[r] k\right)}\,{\rm d}k\,{\rm d}r\over 
r^3}\;.\eqno(6.2.9)$$

\medskip 

\noindent{\bf Proof.}\quad 
Let $\{\gamma_j : j\in{\Bbb N}\}$ be a complete set of representatives 
of the right cosets of $\Gamma_{\frak a}'$ in $\Gamma$. 
Then, by the triangle inequality and (1.5.4), 
it follows that   
$$\left| \sum_{j=1}^J f_{\omega}\left( g_{\frak a}^{-1}\gamma_j g\right) 
\overline{\phi(g)}\right| 
\leq\left[\Gamma_{\frak a} : \Gamma_{\frak a}'\right] 
\left( P^{\frak a}\left| f_{\omega}\right|\right) (g) 
\left|\phi(g)\right| \qquad\quad
\hbox{for all $J\in{\Bbb N}$ and all $g\in G$.}$$
This, combined with Lemma~6.2.2 and the hypothesis that 
$(P^{\frak a}|f_{\omega}|)\cdot\phi\in L^1(\Gamma\backslash G)$, 
enables the application of Lebesgue's `dominated convergence' theorem, so 
that one obtains:  
$$\int_{\Gamma\backslash G} 
\left[\Gamma_{\frak a} : \Gamma_{\frak a}'\right] 
\left( P^{\frak a} f_{\omega}\right) (g) 
\,\overline{\phi(g)}\,{\rm d}g 
=\lim_{J\rightarrow\infty}\int_{\Gamma\backslash G} 
\sum_{j=1}^J f_{\omega}\left( g_{\frak a}^{-1}\gamma_j g\right) 
\overline{\phi(g)}\,{\rm d}g 
=\sum_{\gamma\in\Gamma_{\frak a}'\backslash\Gamma} 
\int_{\Gamma\backslash G}  f_{\omega}\left( g_{\frak a}^{-1}\gamma g\right) 
\overline{\phi(g)}\,{\rm d}g\;.$$
Hence, given that $\phi$ is $\Gamma$-automorphic, 
one finds (by the usual `unfolding' method) that   
$$\left[\Gamma_{\frak a} : \Gamma_{\frak a}'\right] 
\int_{\Gamma\backslash G} 
\left( P^{\frak a} f_{\omega}\right) (g) 
\,\overline{\phi(g)}\,{\rm d}g 
=\sum_{\gamma\in\Gamma_{\frak a}'\backslash\Gamma} 
\int_{\gamma{\cal F}_{\Gamma\backslash G}}  
f_{\omega}\left( g_{\frak a}^{-1}g\right) 
\overline{\phi(g)}\,{\rm d}g\;,\eqno(6.2.10)$$
where ${\cal F}_{\Gamma\backslash G}$ denotes a (measurable) fundamental domain 
for $\Gamma\backslash G$, while 
$\gamma{\cal F}_{\Gamma\backslash G}=\{ \gamma g : g\in {\cal F}_{\Gamma\backslash G}\}$. 

One may repeat the above steps with $|f_{\omega}|$ and $|\phi|$ substituted 
for $f_{\omega}$ and $\phi$ (respectively). Consequently, and by virtue of 
the countable additivity of the relevant integrals, it follows that 
the measurable function $g\mapsto |f_{\omega}(g_{\frak a}^{-1} g)\,\overline{\phi(g)}\,|\in[0,\infty)$ 
is integrable over $\Gamma_{\frak a}'\backslash G$ (with respect to the 
measure induced by ${\rm d}g$); the same is therefore true of 
the measurable function $g\mapsto f_{\omega}(g_{\frak a}^{-1} g)\,\overline{\phi(g)}$. 
This shows that the integrals on the right-hand side of Equation~(6.2.10) are countably additive, so that 
(6.2.10) implies the equality 
$$\left[\Gamma_{\frak a} : \Gamma_{\frak a}'\right] 
\int_{\Gamma\backslash G} 
\left( P^{\frak a} f_{\omega}\right) (g) 
\,\overline{\phi(g)}\,{\rm d}g 
=\int_{\Gamma_{\frak a}'\backslash G}  
f_{\omega}\left( g_{\frak a}^{-1}g\right) 
\overline{\phi(g)}\,{\rm d}g\;.\eqno(6.2.11)$$
\par
If ${\cal F}_{\frak a}$ is any fundamental domain for $\Gamma_{\frak a}'\backslash G$, 
then (by virtue of the hypothesis that (1.1.20) holds for ${\frak c}={\frak a}$) 
the set $g_{\frak a}^{-1}{\cal F}_{\frak a} 
=\{ g_{\frak a}^{-1} g : g\in {\cal F}_{\frak a}\}\subset G$ is 
a fundamental domain for $B^{+}\backslash G$, where $B^{+}<N$ is given by (1.1.21). 
Therefore it follows that, by virtue of 
the left invariance of the Haar measure ${\rm d}g$, by the 
hypothesis that $f_{\omega}\in C^0(N\backslash G , \omega)$, 
and by the Iwasawa decomposition $G=NAK$, one has: 
$$\eqalign{
\int_{\Gamma_{\frak a}'\backslash G}  
f_{\omega}\left( g_{\frak a}^{-1}g\right) 
\overline{\phi(g)}\,{\rm d}g
 &=\int_{B^{+}\backslash G}  
f_{\omega}\left( g\right) 
\overline{\phi\left( g_{\frak a} g\right)}\,{\rm d}g =\cr 
 &=\int_{B^{+}\backslash N}\int_A\int_K 
\psi_{\omega}(n) f_{\omega}(ak)\,\overline{\phi\left( g_{\frak a} nak\right)}\,
(\rho(a))^{-2}{\rm d}k\,{\rm d}a\,{\rm d}n\;,}\eqno(6.2.12)$$
where $\rho : G\rightarrow (0,\infty)$ is given by (6.2.1), and the Haar measures 
${\rm d}n$, ${\rm d}a$ and ${\rm d}k$ are as in (1.1.10). 
By Fubini's theorem, one may change the order of integration in 
the final iterated integral in (6.2.12), so as to give priority 
to the integration with respect to the variable $n$; given the 
definitions of $F^{\frak c}_{\omega} f$ and $\psi_{\omega}(n)$ in (1.4.2) and (1.4.3), 
one thereby obtains from (6.2.11) and (6.2.12)  
the results stated in (6.2.8)-(6.2.9)\quad$\blacksquare$  

\bigskip 

\noindent{\bf Fourier expansions at cusps.}\quad  
Suppose now that ${\frak a}$, ${\frak b}$, $g_{\frak a}$, $g_{\frak b}$, 
$\omega$, $\sigma_0$, $R_0$ and $f_{\omega}$ satisfy the combined hypotheses 
of Lemma~6.2.1 and Lemma~6.2.2. Then it follows by Corollary~6.2.3 that, 
for each $\omega'\in{\frak O}$, one may define the Fourier term of 
order $\omega'$ for $P^{\frak a} f_{\omega}$ at ${\frak b}$ to 
be the function $F^{\frak b}_{\omega'} P^{\frak a} f_{\omega} : G\rightarrow{\Bbb C}$ 
given by 
$$\left( F^{\frak b}_{\omega'} P^{\frak a} f_{\omega}\right) (g) 
=\int_{B^{+}\backslash N}\left(\psi_{\omega'}(n)\right)^{-1} 
\left( P^{\frak a} f_{\omega}\right)\left( g_{\frak b} n g\right) {\rm d}n\qquad\qquad 
\hbox{($g\in G$).}\eqno(6.2.13)$$ 
Note that this accords with the definition given in (1.4.2), though there 
we dealt only with the Fourier expansions of functions lying in the space 
$C^{\infty}(\Gamma\backslash G)$; our present hypotheses do not even imply that 
$P^{\frak a}f_{\omega}$ is differentiable on $G$, nor do they imply that 
the Fourier series $\sum_{\omega'\in{\frak O}} (F^{\frak b}_{\omega'}P^{\frak a}f_{\omega})(g)$ 
is convergent for all $g\in G$ (see Section~13.41 of~[43] for a relevant example). 
Therefore our present hypotheses do not suffice to ensure that 
the Fourier expansion (1.4.1) is valid, for all $g\in G$, when one substitutes 
for ${\frak c}$ and $f$ (there) the cusp ${\frak b}$ and function $P^{\frak a}f_{\omega}$, 
respectively. 
\par 
We shall address the question of the representation of Poincar\'{e} series 
by their Fourier expansions (at cusps) in an ad hoc manner, 
and only as the need arises: see Remark~6.2.6 and the proof of the result 
(6.5.76) of Lemma~6.5.14, below. 
\par 
Regardless of whether or not it is the case that  
$\sum_{\omega'\in{\frak O}} F^{\frak b}_{\omega'}P^{\frak a}f_{\omega} 
=((P^{\frak a}f_{\omega}) | {\frak b})$, it does 
follow from the definitions (6.2.13) and (1.5.4) that, for $\omega'\in{\frak O}$, one has 
$$\eqalignno{
\left( F^{\frak b}_{\omega'} P^{\frak a} f_{\omega}\right) (g) 
 &={1\over\left[\Gamma_{\frak a} : \Gamma_{\frak a}'\right]} 
\int\limits_{B^{+}\backslash N}
\sum_{\gamma\in\Gamma_{\frak a}'\backslash\Gamma} 
\left(\psi_{\omega'}(n)\right)^{-1} 
f_{\omega}\left( g_{\frak a}^{-1} \gamma g_{\frak b} n g\right) {\rm d}n =\cr 
 &={1\over\left[\Gamma_{\frak a} : \Gamma_{\frak a}'\right]} 
\sum_{\gamma\in\Gamma_{\frak a}'\backslash\Gamma} 
\ \int\limits_{B^{+}\backslash N}
\left(\psi_{\omega'}(n)\right)^{-1} 
f_{\omega}\left( g_{\frak a}^{-1} \gamma g_{\frak b} n g\right) {\rm d}n\;. &(6.2.14)
}$$
Note that $B^{+}\backslash N$ 
has a compact fundamental domain, namely the set 
$\{ n[z] : -1/2\leq {\rm Re}(z), {\rm Im}(z)\leq 1/2\}$, 
so that the uniform convergence established in Lemma~6.2.2 
justifies the term by term integration by which the final 
equality in (6.2.14) is obtained. 

\bigskip 

\proclaim Lemma 6.2.5. 
Let the combined hypotheses 
of Lemma~6.2.1 and Lemma~6.2.2 be satisfied, and let $\omega'\in{\frak O}$.  
Then (1.5.4) and (6.2.13) define a function 
$F^{\frak b}_{\omega'} P^{\frak a} f_{\omega}$ which lies in the 
space $C^{0}(N\backslash G,\omega')$.  
For $g\in G$ the formula for  
$\left( F^{\frak b}_{\omega'} P^{\frak a} f_{\omega}\right) (g)$ 
implied by the case ${\frak a}'={\frak b}$ of (1.5.5)-(1.5.10) holds, 
and the sums and integrals occurring 
(explicitly or implicitly) in the equation (1.5.5) are absolutely convergent. 

\medskip

\noindent{\bf Proof.}\quad 
By Corollary~6.2.3 and the case ${\frak c}={\frak b}$ of 
(1.1.20), the function $g\mapsto (P^{\frak a} f_{\omega})(g_{\frak b} g)$ 
is continuous on $G$, and satisfies  $(P^{\frak a} f_{\omega})(g_{\frak b} b g) 
=(P^{\frak a} f_{\omega})(g_{\frak b} g)$, for $g\in G$ and $b\in B^{+}$. 
It follows that, for all $r_1,r_2\in(0,\infty)$, the 
function $g\mapsto (P^{\frak a} f_{\omega})(g_{\frak b} g)$ is  
uniformly continuous on the set  
$\{ n a[r] k : n\in N, r_1\leq r\leq r_2, k\in K\}$.   
This fact, combined with the continuity of $(\psi_{\omega'}(n))^{-1}$ and 
the fact that $B^{+}\backslash N$ is compact (and of finite measure 
with respect to ${\rm d}n$), is sufficient to   
establish both the existence and 
continuity (as a function of $g$)   
of the integral on the right-hand side of Equation~(6.2.13): we may  
conclude that (1.5.4) and (6.2.13) define a continuous function  
$F^{\frak b}_{\omega'} P^{\frak a} f_{\omega} : G\rightarrow{\Bbb C}$. 
Moreover, since the measure ${\rm d}n$ on $N$ induces a Haar measure 
on the group $B^{+}\backslash N$  
it is an immediate consequence of 
(6.2.13) and (1.4.3) that 
$(F^{\frak b}_{\omega'} P^{\frak a} f_{\omega})(ng)=\psi_{\omega'}(n)
(F^{\frak b}_{\omega'} P^{\frak a} f_{\omega})(g)$   
for all $g\in G$ and all $n\in N$; this  completes the proof that one has 
$F^{\frak b}_{\omega'} P^{\frak a} f_{\omega}\in C^0(N\backslash G,\omega')$. 

Let $g\in G$. Our approach to the proof of the formula (1.5.5)  
is modelled on Section~7.2 of [32] (in which the 
case ${\frak a}={\frak b}=\infty$ is treated); this differs from  
the approach taken on Pages~39 and~40 of [5], where 
the case $\Gamma=SL(2,{\frak O})$ of (1.5.5) is 
obtained via the Poisson summation formula. We observe 
firstly that, by (6.2.14) and the partitioning of 
$\Gamma$ noted in Lemma~6.2.1, and by (1.1.20) (for ${\frak c}={\frak b}$) and (6.2.5), 
one has 
$$\eqalign{\left[\Gamma_{\frak a} : \Gamma_{\frak a}'\right] 
\left( F^{\frak b}_{\omega'} P^{\frak a} f_{\omega}\right) (g) 
 &=\sum_{\gamma\in\Gamma_{\frak a}'\backslash{}^{\frak a}\Gamma^{\frak b}(0)} 
\ \int\limits_{B^{+}\backslash N}\left(\psi_{\omega'}(n)\right)^{-1} 
f_{\omega}\left( g_{\frak a}^{-1}\gamma g_{\frak b} ng\right) {\rm d}n\ + \cr
 &\quad +\sum_{c\in{}^{\frak a}{\cal C}^{\frak b}} 
\,\sum_{\gamma\in\Gamma_{\frak a}'\backslash{}^{\frak a}\Gamma^{\frak b}(c)/\Gamma_{\frak b}'} 
\ \int\limits_{N}\left(\psi_{\omega'}(n)\right)^{-1} 
f_{\omega}\left( g_{\frak a}^{-1}\gamma g_{\frak b} ng\right) {\rm d}n 
}\eqno(6.2.15)$$
(we have used here the fact that if $0\neq c\in{}^{\frak a}{\cal C}^{\frak b}$, 
and if ${\cal Z}$ is a complete set of representatives for 
the set of double cosets 
$\{ \Gamma_{\frak a}'\gamma\Gamma_{\frak b}' : \gamma\in{}^{\frak a}\Gamma^{\frak b}(c)\}$, 
then the family of sets $(\gamma\Gamma_{\frak b}')_{\gamma\in{\cal Z}}$ is  a 
partitioning of a complete set of representatives for the 
set of right cosets $\{ \Gamma_{\frak a}'\gamma : \gamma\in{}^{\frak a}\Gamma^{\frak b}(c)\}$). 

With regard to the first sum on the right-hand side of (6.2.15), we recall, 
from the proof of Lemma 6.2.1, that 
${}^{\frak a}\Gamma^{\frak b}(0)=\{\gamma\in\Gamma : \gamma{\frak b}={\frak a}\}$, 
and that for each $\gamma\in {}^{\frak a}\Gamma^{\frak b}(0)$ there is some  
$\eta=\eta(\gamma)\in {\Bbb C}$, with $\eta^2\in{\frak O}^{*}$ and 
some $\beta=\beta(\gamma)\in{\Bbb C}$ such that 
$h[\eta]n[\beta]=g_{\frak a}^{-1}\gamma g_{\frak b}$, so that  
$g_{\frak a}^{-1}\gamma g_{\frak b}n[z]=n[\eta^2 z]g_{\frak a}^{-1}\gamma g_{\frak b}$ for 
$z\in{\Bbb C}$. Hence, by the hypothesis that 
$f_{\omega}\in C^0(N\backslash G,\omega)$, and by (1.4.3), 
we find that  
$$\eqalignno{
\sum_{\gamma\in\Gamma_{\frak a}'\backslash{}^{\frak a}\Gamma^{\frak b}(0)} 
\ \int\limits_{B^{+}\backslash N}\!\!\left(\psi_{\omega'}(n)\right)^{-1} 
f_{\omega}\!\left( g_{\frak a}^{-1}\gamma g_{\frak b} ng\right) {\rm d}n
 &=\sum_{\gamma\in\Gamma\,:\,\gamma{\frak b}={\frak a}}
\ \,\int\limits_{B^{+}\backslash N} 
\!\!\left(\psi_{\omega'}(n)\right)^{-1} 
\psi_{(\eta(\gamma))^2\omega}(n)\,{\rm d}n
\,f_{\omega}\!\left( g_{\frak a}^{-1}\gamma g_{\frak b} g\right) =\cr 
 &=\sum_{\gamma\in\Gamma\,:\,\gamma{\frak b}={\frak a}}
\delta_{\omega', (\eta(\gamma))^2\omega} 
\,f_{\omega}\left( g_{\frak a}^{-1}\gamma g_{\frak b} g\right) . &(6.2.16)}$$
As for the second sum on the right-hand side of (6.2.15), it follows 
by (1.5.8) that if $0\neq c\in{}^{\frak a}{\cal C}^{\frak b}$ then, 
for each $\gamma\in{}^{\frak a}\Gamma^{\frak b}(c)$, there 
is some $(s,d)=(s(\gamma),d(\gamma))\in{\Bbb C}^2$ such that 
$$g_{\frak a}^{-1}\gamma g_{\frak b}=n[s/c]h[1/c]k[0,-1]n[d/c]\eqno(6.2.17)$$
(this being essentially the same device introduced in 
the equation (6.1.8) of Lemma~6.1.1). Hence, supposing now that 
$c\in{}^{\frak a}{\cal C}^{\frak b}$ is given, 
it follows (by the hypothesis that $f_{\omega}\in C^0(N\backslash G,\omega)$, 
and the fact that ${\rm d}n$ is a Haar measure for $N$) that if  
$\gamma\in{}^{\frak a}\Gamma^{\frak b}(c)$ then
$$\eqalign{
\int\limits_{N}\left(\psi_{\omega'}(n)\right)^{-1} 
f_{\omega}\left( g_{\frak a}^{-1}\gamma g_{\frak b} ng\right) {\rm d}n 
 &=\int\limits_{N}\left(\psi_{\omega'}(n)\right)^{-1} 
f_{\omega}\left(  n[s(\gamma)/c]h[1/c]k[0,-1]n[d(\gamma)/c]ng\right) {\rm d}n =\cr 
 &=\psi_{\omega}\left( n[s(\gamma)/c]\right)\psi_{\omega'}\left( n[d(\gamma)/c]\right) 
\int\limits_{N}\left(\psi_{\omega'}(n)\right)^{-1} 
f_{\omega}\left(  h[1/c]k[0,-1]ng\right) {\rm d}n\;.
}$$
By (1.5.2) and (1.5.7), the last integral in the above may be expressed 
as the Jacquet integral $({\bf J}_{\omega'}{\bf h}_{1/c} f_{\omega})(g)$; 
since this integral is independent of the `$\gamma$' in the above equations 
(`$c$' being fixed there), and since   
$$\sum_{\gamma\in\Gamma_{\frak a}'\backslash{}^{\frak a}\Gamma^{\frak b}(c)/\Gamma_{\frak b}'} 
\psi_{\omega}\left( n[s(\gamma)/c]\right)\psi_{\omega'}\left( n[d(\gamma)/c]\right) 
=S_{{\frak a},{\frak b}}\left(\omega ,\omega';c\right)$$
(see (1.4.3), (1.5.10) and (6.2.17)), it therefore follows that, 
for $\,c\in{}^{\frak a}{\cal C}^{\frak b}$, 
$$\sum_{\gamma\in\Gamma_{\frak a}'\backslash{}^{\frak a}\Gamma^{\frak b}(c)/\Gamma_{\frak b}'} 
\ \int\limits_{N}\left(\psi_{\omega'}(n)\right)^{-1} 
f_{\omega}\left( g_{\frak a}^{-1}\gamma g_{\frak b} ng\right) {\rm d}n 
=S_{{\frak a},{\frak b}}\left(\omega ,\omega';c\right) 
\left( {\bf J}_{\omega'}{\bf h}_{1/c} f_{\omega}\right)\!(g)\;.$$
By combining this last result with (6.2.15) and (6.2.16), we arrive 
at the case ${\frak a}'={\frak b}$ of the 
formula stated in (1.5.5). By consideration of the steps 
used to obtain this formula, it may be seen that 
the absolute convergence of all the sums and integrals on the 
right-hand side of (1.5.5) is a direct consequence of 
Lemma~6.2.2 (this becomes clear if one lets 
$|f_{\omega}|$ be substituted for $f_{\omega}$, and then
considers the case $\omega' =\omega =0$)\quad$\blacksquare$ 

\bigskip 

\noindent{\bf Remark~6.2.6.}\quad 
The Fourier expansion (1.8.4)-(1.8.6) of the Eisenstein series 
$E_{\ell,q}^{\frak a}(\nu ,p)$ may be 
shown to follow, via (1.4.1)-(1.4.3), from the case $\omega =0$, 
$f=f_{0}=\varphi_{\ell,k}(\nu ,p)$ of Lemma~6.2.5: 
note in particular that, by the definitions 
in (1.3.2) and (6.2.1), the hypothesis that ${\rm Re}(\nu)>1$ suffices to  
ensure that the condition (6.2.5) is satisfied 
when $f_{\omega}=\varphi_{\ell ,q}(\nu ,p)$ and $R_0=1$ (say). 
The use of the Fourier expansion~(1.4.1) in this context may be justified through 
an appeal  to Theorem~67 of~[2]: we skip the relevant details,  
which are similar to what occurs in the last two paragraphs of the 
proof of Lemma~6.5.14, below.
\par 
By making use of the meromorphic continuation of the functions 
$\nu\mapsto E_{\ell,q}^{\frak a}(\nu ,p)$ and  
$\nu\mapsto D_{\frak a}^{\frak b}(\psi;\nu,p)$, which is discussed in 
Subsection 1.8, one may dispense with the condition ${\rm Re}(\nu)>1$;  
the Fourier expansion in (1.8.4) is thereby obtained for 
all $(\nu,p)\in({\Bbb C}\times{\Bbb Z})-\{ (0,0),(1,0)\}$ 
such that ${\rm Re}(\nu)\geq 0$ and $|p|\leq\ell$. 

\bigskip 

For the proof of the spectral sum formula, Theorem~B, 
we require certain `cusp sector estimates' 
for the Eisenstein series $E_{\ell,q}^{\frak a}(\nu ,p)$    
and other Poincar\'e series; we deduce the required estimates 
(those of Lemma~6.2.8 and Lemma~6.2.9, below) 
with the help of the following lemma, which is 
taken from Section~5.2 of [32]. 

\bigskip 

\proclaim Lemma 6.2.7. Let $\ell,p,q\in{\Bbb Z}$ satisfy $\ell\geq\max\{ |p|,|q|\}$; 
let $\nu\in{\Bbb C}$; let ${\frak b}$ be a cusp of $\Gamma$; and let 
$g_{\frak b}\in G$ be such that (1.1.16) and (1.1.20)-(1.1.21) hold for 
${\frak c}={\frak b}$. Suppose moreover that $f : G\rightarrow{\Bbb C}$ satisfies 
$$f\left( g_{\frak b} g\right) 
=\sum_{0\neq\omega\in{\frak O}} c(\omega)
\left( {\bf J}_{\omega}\varphi_{\ell ,q}(\nu ,p)\right)\!(g)\qquad\qquad 
\hbox{($g\in G$),}\eqno(6.2.18)$$
where, for $0\neq\omega\in{\frak O}$, the coefficient $c(\omega)\in{\Bbb C}$ is 
independent of $g$. 
Then there exists some $r_0\in [1,\infty)$ such that  
$$\left| f\left( g_{\frak b} g\right)\right| \leq e^{-\pi\rho(g)}\qquad\quad 
\hbox{for all $g\in G$ with $\rho(g)\geq r_0$.}\eqno(6.2.19)$$

\medskip 

\noindent{\bf Proof.}\quad 
This lemma is a slight variation on the special case of 
Part~(i) of Lemma~5.2.1 of [32] in which one 
has, in the notation of [32], $p\in{\Bbb Z}\,$ (rather than $2p\in{\Bbb Z}-2{\Bbb Z}$), $\Lambda_{\kappa}'=(1/2){\frak O}$ 
and $\chi : \Gamma\rightarrow\{ 1\}$. Note, in particular, that 
it is clear from the proof given in [32] that the 
upper bounds stated in the equations~(5.9) and~(5.10) of [32] are not optimal, and 
may indeed be sharpened by any factor of the form 
$\exp(-(1-\varepsilon)2\pi\omega_0 r)$, where $\varepsilon$ denotes 
an arbitrarily small positive absolute constant, while  
$\omega_0=\min\{ |\omega| : 0\neq\omega\in\Lambda_{\kappa}'\}\,$ (so that 
$\omega_0=1/2$ when $\Lambda_{\kappa}'=(1/2){\frak O}$).  
Hence, by assuming greater lower bounds for $r=\rho(g)$ than 
those implicit in the equation~(5.9) of [32], one may  
sharpen the upper bound stated there by any factor of 
the form $O(r^{\ell+1/2})\quad\blacksquare$ 

\bigskip 

\proclaim Lemma 6.2.8. 
Let $\ell,q\in{\Bbb Z}$ satisfy $\ell\geq |q|$,  
let $(\nu ,p)\in{\Bbb C}\times{\Bbb Z}-\{ (0,0),(1,0)\}$ 
be such that ${\rm Re}(\nu)\geq 0$ and 
$|p|\leq\ell$; let ${\frak a} , {\frak b}$ be cusps of $\Gamma$, and let 
$g_{\frak a},g_{\frak b}\in G$ be such that (1.1.16) and (1.1.20)-(1.1.21) hold for 
${\frak c}\in\{ {\frak a},{\frak b}\}$. 
Then, for some $r_0=r_0(\Gamma ,\ell ,\nu)\in [1,\infty)$, and some 
$\epsilon =\epsilon(g_{\frak a},g_{\frak b})\in{\frak O}^{*}$ satisfying 
$\epsilon =1$ if ${\frak a}={\frak b}$ and $g_{\frak a}=g_{\frak b}$, 
one has: 
$$\epsilon^p E_{\ell ,q}^{\frak a}(\nu ,p)\left( g_{\frak b} g\right) 
=\delta^{\Gamma}_{{\frak a},{\frak b}}\varphi_{\ell,q}(\nu ,p) (g) 
+O_{\Gamma,\ell ,\nu}\left( \left(\rho(g)\right)^{1-{\rm Re}(\nu)}\right) 
\qquad\  
\hbox{for all $g\in G$ with $\rho(g)\geq r_0$.}\eqno(6.2.20)$$

\medskip 

\noindent{\bf Proof.}\quad 
We may suppose that either 
${\frak a}={\frak b}$ and 
$g_{\frak a}=g_{\frak b}$, or else 
the cusps ${\frak a}$ and ${\frak b}$ are not 
$\Gamma$-equivalent (by virtue of the point noted 
two lines below (1.8.3), these two cases of the lemma imply every 
other case).
Then, in light of Remark~6.2.6, we have the Fourier expansion of $E_{\ell,q}^{\frak a}(\nu,p)$ 
shown in equation (1.8.4). 
It follows that by putting, for $g\in G$,  
$$f(g)=E_{\ell,q}^{\frak a}(\nu,p)(g) 
-\delta^{\Gamma}_{{\frak a},{\frak b}}\varphi_{\ell,q}(\nu,p)\left( g_{\frak b}^{-1} g\right)
-{D_{\frak a}^{\frak b}(0;\nu,p)\over \left[\Gamma_{\frak a} : \Gamma_{\frak a}'\right]}
\,{\pi\Gamma(|p|+\nu)\over\Gamma(\ell+1+\nu)}\,{\Gamma(\ell+1-\nu)\over\Gamma(|p|+1-\nu)}\,
\varphi_{\ell,q}(-\nu,-p)\left( g_{\frak b}^{-1} g\right)$$
we define a function $f : G\rightarrow {\Bbb C}$ satisfying, 
for a certain choice of coefficients $c(\omega)\,$ ($\omega\in{\frak O}-\{ 0\}$), 
the hypothesis (6.2.18) of Lemma 6.2.7; the result (6.2.19) of that lemma   
therefore implies that there exists some $r_0\in [1,\infty)$ such that, 
for all $g\in G$ with $\rho(g)\geq r_0$, one has the equation  
$$E_{\ell,q}^{\frak a}(\nu,p)\left( g_{\frak b} g\right) 
=\delta^{\Gamma}_{{\frak a},{\frak b}}\varphi_{\ell,q}(\nu,p)(g)  
+O_{\Gamma,g_{\frak a},g_{\frak b},\ell,q,p,\nu}\left(\left(\rho(g)\right)^{1-{\rm Re}(\nu)}\right) 
+O\left( e^{-\pi\rho(g)}\right)\eqno(6.2.21)$$
(in which the first of the two $O$-terms  
represents an estimate for the third term on the right-hand 
side of the equation defining $f(g)$, and is justified by virtue of the 
definitions (1.3.2), (1.5.10), (1.8.6), (6.2.1) and the 
case $\psi =0$ of the meromorphic continuation of 
$D_{\frak a}^{\frak b}(\psi ;\nu,p)$ discussed in Subsection 1.8). 

Since $r_0\geq 1$, the final $O$-term in (6.2.21) may 
be omitted: for one has $e^{-\pi\rho}=O_{\nu}(\rho^{1-{\rm Re}(\nu)})$ 
when $\rho\geq 1$. Moreover, since $p$ and $q$ must lie in the finite set 
$\{ -\ell,-\ell+1,\ldots ,\ell\}$, the implicit constant 
associated with the first $O$-term of (6.2.21) may be chosen so 
as to depend only upon $\Gamma$, $g_{\frak a}$, $g_{\frak b}$, $\ell$ and $\nu$; 
indeed, by the point noted two lines below equation (1.8.3), and by the fact 
that the number of $\Gamma$-equivalence classes of cusps is finite, 
it follows that this implicit constant need depend only upon $\Gamma$, $\ell$ and 
$\nu$; for the same reasons, a suitable choice of $r_0$ (in the above)  
may be determined from just $\Gamma$, $\ell$ and $\nu\quad\blacksquare$

\bigskip 

\noindent{\bf Square integrable Poincar\'e series.}\quad  
We assume (as in the preceding discussion) that ${\frak a}$ is a 
cusp of $\Gamma$, that $\omega\in{\frak O}$, and that, for some 
$\sigma_0 >1$, and some $R_0 >0$, the function 
$f_{\omega}\in C^0(N\backslash G,\omega)$ satisfies the 
condition (6.2.5). In order to establish the square-integrability 
(over $\Gamma\backslash G$) of the Poincar\'e series 
$P^{\frak a} f_{\omega}$, we require also that, for some 
$\sigma_{\infty} >0$, and some $R_{\infty}\in (0,\infty)$, the function 
$f_{\omega}$ satisfies: 
$$f_{\omega}(g)\ll_{f_{\omega}}\, \left(\rho(g)\right)^{1-\sigma_{\infty}}\qquad\quad 
\hbox{for all $g\in G$ such that $\rho(g)\geq R_{\infty}$.}\eqno(6.2.22)$$ 
This implies (just as (6.2.5) implies (6.2.6)) that if $R>0$ then 
$$f_{\omega}(g)\ll_{f_{\omega},\sigma_{\infty},R_{\infty},R}
\ \left(\rho(g)\right)^{1-\sigma_{\infty}}\qquad\quad 
\hbox{for all $g\in G$ such that $\rho(g)\geq R$.}\eqno(6.2.23)$$ 

\bigskip 

\proclaim Lemma 6.2.9. Let $0\neq q_0\in{\frak O}$ and 
$\Gamma =\Gamma_0(q_0)\leq SL(2,{\frak O})$; let ${\frak a}$ and ${\frak b}$ 
be cusps of $\Gamma$; let $g_{\frak a},g_{\frak b}\in G$ be such that 
(1.1.16) and (1.1.20)-(1.1.21) hold for ${\frak c}\in\{ {\frak a} ,{\frak b}\}$; 
let $m_{\frak b}\in {\frak O}-\{ 0\}$ be as described below (1.1.22); 
let $\omega\in{\frak O}$, $\sigma_0 >1$, $R_0>0$, $\sigma_{\infty}>0$ and 
$R_{\infty}\in (0,\infty)$, and let $f=f_{\omega}\in C^0(N\backslash G,\omega)$ 
be such that both of the conditions (6.2.5) and (6.2.22) are satisfied. Then, for all 
$g\in G$ such that $\rho(g) > 1/|m_{\frak b}|$, one has: 
$$\left( P^{\frak a} f_{\omega}\right) \left( g_{\frak b} g\right) 
=\delta^{\Gamma}_{{\frak a},{\frak b}}\,O_{\Gamma,f,\sigma_{\infty},R_{\infty}}\left( 
\left(\rho(g)\right)^{1-\sigma_{\!\infty}}\right) 
+O_{\Gamma,f,\sigma_{0},R_{0}}\left( 
\left(\rho(g)\right)^{1-\sigma_{0}}\right) .\eqno(6.2.24)$$

\medskip 

\noindent{\bf Proof.}\quad 
It will suffice to prove that (6.2.24) holds if 
if ${\frak a}={\frak b}$ 
and $g_{\frak a}=g_{\frak b}$, or if  the cusps ${\frak a}$ and ${\frak b}$ 
are not $\Gamma$-equivalent: for in every other case one has 
${\frak a}=\tau {\frak b}$ for some $\tau\in\Gamma$, which 
(by (1.1.4), (6.2.1) and the result (2.1) of Lemma~2.1,  
and by virtue the fact that $P^{\frak a} f_{\omega}$ is $\Gamma$-automorphic) 
implies that, for $g\in G$, one has 
$(P^{\frak a} f_{\omega})(g_{\frak b} g)=(P^{\frak a} f_{\omega})(\tau g_{\frak b} g) 
=(P^{\frak a} f_{\omega})( g_{\frak a} \tilde g)$, where 
$\tilde g =g_{\frak a}^{-1}\tau g_{\frak b} g$ is such that $\rho(\tilde g)=\rho(g)$. 

In cases where ${\frak a}$ and ${\frak b}$ are not $\Gamma$-equivalent 
the set ${}^{\frak a}\Gamma^{\frak b}(0)$ (defined in (1.5.8)) is empty,  
and so in these cases it follows by the result (6.2.4) of Lemma~6.2.1 
that if $g\in G$, and if $\rho(g)>1/|m_{\frak b}|$, then 
$$\rho\left( g_{\frak a}\gamma g_{\frak b} g\right) <1\qquad\quad 
\hbox{for all $\,\gamma\in\Gamma$.}\eqno(6.2.25)$$
Therefore one may, in such cases, apply the definition (1.5.4) and hypothesis 
(6.2.5) (with its corollary (6.2.6)) so as to obtain,  
for $g\in G$ such that  
$\rho(g)>1/|m_{\frak b}|$,  
$$\eqalignno{ 
\left( P^{\frak a} f_{\omega}\right)\left( g_{\frak b} g\right) 
 &={1\over\left[\Gamma_{\frak a} :\Gamma_{\frak a}'\right]} 
\sum_{\gamma\in\Gamma_{\frak a}'\backslash\Gamma} 
f_{\omega}\left( g_{\frak a}^{-1}\gamma g_{\frak b} g\right) = 
\qquad\qquad\qquad\qquad\qquad\qquad\qquad\qquad\qquad\qquad\qquad\ \cr 
 &={1\over\left[\Gamma_{\frak a} :\Gamma_{\frak a}'\right]} 
\sum_{\gamma\in\Gamma_{\frak a}'\backslash\Gamma} 
O_{f_{\omega},\sigma_0,R_0}\left(
\left(\rho\left( g_{\frak a}^{-1}\gamma g_{\frak b} g\right)\right)^{1+\sigma_0}\right)  
\ll_{f_{\omega},\sigma_0,R_0} E_{0,0}^{\frak a}\left(\sigma_0,0\right)\left( g_{\frak b} g\right) 
&(6.2.26)}$$
(the final upper bound here following by (6.2.2) and (1.8.1)). 
By (6.2.26) and the result (6.2.20) of Lemma~6.2.8, 
we obtain proof of those cases of the lemma in which ${\frak a}$ and ${\frak b}$ 
are not $\Gamma$-equivalent. 
\par 
Given the conclusion just reached, and point noted at the beginning of this proof, 
we may assume henceforth that ${\frak a}={\frak b}$ and $g_{\frak a}=g_{\frak b}$. 
Consequently we have only to obtain a suitable bound for 
$(P^{\frak a} f_{\omega})(g_{\frak a} g)$, and (in doing so) may assume that 
$g\in G$ is such that $\rho(g)>1/|m_{\frak a}|$. 
\par 
By the definition (1.5.4) and Lemma~6.2.1, and by (1.5.8) and (1.1.16), 
we find that 
$$\left( P^{\frak a} f_{\omega}\right)\left( g_{\frak a} g\right) 
=S_0(g)+S_1(g)\;,\eqno(6.2.27)$$
where 
$$S_0(g) 
={1\over\left[\Gamma_{\frak a} :\Gamma_{\frak a}'\right]} 
\sum_{\gamma\in\Gamma_{\frak a}'\backslash {}^{\frak a}\Gamma^{\frak a}(0)} 
f_{\omega}\left( g_{\frak a}^{-1}\gamma g_{\frak a} g\right)
={1\over\left[\Gamma_{\frak a} :\Gamma_{\frak a}'\right]} 
\sum_{\gamma\in\Gamma_{\frak a}'\backslash\Gamma_{\frak a}} 
f_{\omega}\left( g_{\frak a}^{-1}\gamma g_{\frak a} g\right)$$ 
and 
$$S_1(g) 
={1\over\left[\Gamma_{\frak a} :\Gamma_{\frak a}'\right]} 
\sum_{0\neq c\in{}^{\frak a}{\cal C}^{\frak a}}
\sum_{\gamma\in\Gamma_{\frak a}'\backslash{}^{\frak a}\Gamma^{\frak a}(c)} 
f_{\omega}\left( g_{\frak a}^{-1}\gamma g_{\frak a} g\right) .$$ 
As a consequence of both the result (6.2.4) of Lemma~6.2.1 and the hypothesis 
(6.2.5) (with corollary (6.2.6)), it follows (similarly to how (6.2.25) and (6.2.26) 
were obtained) that one has 
$$\eqalign{ 
S_1(g) &\ll_{f_{\omega},\sigma_0,R_0} 
\,{1\over\left[\Gamma_{\frak a} :\Gamma_{\frak a}'\right]} 
\sum_{0\neq c\in{}^{\frak a}{\cal C}^{\frak a}}
\sum_{\gamma\in\Gamma_{\frak a}'\backslash{}^{\frak a}\Gamma^{\frak a}(c)} 
\left(\rho\left( g_{\frak a}^{-1}\gamma g_{\frak a} g\right)\right)^{1+\sigma_0} = \cr 
 &\qquad =E_{0,0}^{\frak a}\left(\sigma_0,0\right)\left( g_{\frak a} g\right) 
- {1\over\left[\Gamma_{\frak a} :\Gamma_{\frak a}'\right]}  
\sum_{\gamma\in\Gamma_{\frak a}'\backslash\Gamma_{\frak a}} 
\varphi_{0,0}\left(\sigma_0,0\right)\left( g_{\frak a}^{-1}\gamma g_{\frak a} g\right) . 
}$$
By Lemma~4.2 and (1.8.2) (or (6.2.1), (6.2.2) and (1.1.4)), we have here 
$\varphi_{0,0}(\sigma_0,0)(g_{\frak a}^{-1}\gamma g_{\frak a} g) 
=\varphi_{0,0}(\sigma_0,0)(g)$ when $\gamma\in\Gamma_{\frak a}$, and therefore may 
deduce (using Lemma~6.2.8) that  
$$S_1(g) 
=O_{f_{\omega},\sigma_0,R_0}\left( 
E_{0,0}^{\frak a}\left(\sigma_0,0\right)\left( g_{\frak a} g\right) 
-\varphi_{0,0}\left(\sigma_0,0\right) (g)\right) 
= O_{f_{\omega},\sigma_0,R_0,\Gamma}\left( 
\left(\rho(g)\right)^{1-\sigma_0}\right) .\eqno(6.2.28)$$ 
\par 
We now consider the sum $S_0(g)$, defined below (6.2.27). By the result (6.2.3) of Lemma~6.2.1, 
we have $\rho(g_{\frak a}^{-1}\gamma g_{\frak a} g)=\rho(g)$ for all 
$\gamma\in\Gamma_{\frak a}={}^{\frak a}\Gamma^{\frak a}(0)$. Therefore, given that 
$\rho(g)>1/|m_{\frak b}|\geq 1/|q_0|$ (the last inequality following by (6.1.25)), 
it is a consequence of the hypothesis (6.2.22) (with corollary (6.2.23)) 
that we have 
$$S_0(g) 
={1\over\left[\Gamma_{\frak a} :\Gamma_{\frak a}'\right]} 
\sum_{\gamma\in\Gamma_{\frak a}'\backslash\Gamma_{\frak a}} 
O_{f_{\omega},\sigma_{\infty},R_{\infty},\Gamma} 
\left( 
\left(\rho\left( g_{\frak a}^{-1}\gamma g_{\frak a} g\right)\right)^{1-\sigma_{\infty}} 
\right) 
\ll_{f_{\omega},\sigma_{\infty},R_{\infty},\Gamma} 
\,\left(\rho(g)\right)^{1-\sigma_{\infty}}\;.$$ 
By this last bound, that in (6.2.28), and the equation (6.2.27), we obtain the 
case ${\frak a}={\frak b}$, $g_{\frak a}=g_{\frak b}$ of the lemma, and so 
complete its proof\quad$\blacksquare$ 

\bigskip 

\proclaim Corollary~6.2.10. Let those of the hypotheses of the above lemma 
that concern $q_0$, $\Gamma$, ${\frak a}$, $g_{\frak a}$, $\omega$, 
$\sigma_0$, $R_0$, $\sigma_{\infty}$, $R_{\infty}$ and $f_{\omega}$ 
be satisfied. Then one has $P^{\frak a} f_{\omega}\in L^2(\Gamma\backslash G)$; 
if, moreover, $\sigma_{\infty}\geq 1$ then $P^{\frak a} f_{\omega}$ is bounded on $G$. 

\medskip 

\noindent{\bf Proof.}\quad 
It may be proved, similarly to Proposition~2.2.4 of~[11], 
that the inner product $\langle f , g\rangle_{\Gamma\backslash G}$ defined 
in (1.2.2) is independent of the choice of fundamental domain for the 
action of $\Gamma$ upon ${\Bbb H}_3$. Hence, by choosing to 
replace ${\cal F}$ (in (1.2.2)) 
by a fundamental domain ${\cal F}_{*}$ fitting the description 
given in (1.1.22)-(1.1.24), we find that 
$$\left\| P^{\frak a} f_{\omega}\right\|_{\Gamma\backslash G}^2 
=\langle P^{\frak a} f_{\omega} , P^{\frak a} f_{\omega}\rangle_{\Gamma\backslash G} 
=\int_{{\cal F}_{*}}\int_{K^{+}}  
\left| \left( P^{\frak a} f_{\omega}\right) (n[z] a[r] k)\right|^2 {\rm d}k\, 
r^{-3} d_{+}z\,{\rm d}r$$
(provided that the last integral exists). 
Therefore, and since the union formed in (1.1.24) is (see Lemma~2.2) a union 
of a finite number of sets, we have 
$P^{\frak a} f_{\omega}\in L^2(\Gamma\backslash G)$
if, when ${\frak C}(\Gamma)$, ${\cal D}$ and the family of sets 
$({\cal E}_{\frak c})_{{\frak c}\in{\frak C}(\Gamma)}$ are as 
described in the paragraph containing (1.1.23)-(1.1.24), one has: 
$$\int_{\cal X}\int_{K}  
\left| \left( P^{\frak a} f_{\omega}\right) (n[z] a[r] k)\right|^2 {\rm d}k\, 
r^{-3} d_{+}z\,{\rm d}r < \infty\qquad\quad 
\hbox{for $\,{\cal X}\in\{ {\cal D}\}\cup\left\{ {\cal E}_{\frak c} : 
{\frak c}\in{\frak C}(\Gamma)\right\}$.}\eqno(6.2.29)$$
\par 
We consider firstly the case ${\cal X}={\cal D}$ (some compact hyperbolic polygon, 
contained in ${\Bbb H}_3$). By the Iwasawa decomposition of $G$, one has 
$$\int_{\cal D}\int_{K}  
\left| \left( P^{\frak a} f_{\omega}\right) (n[z] a[r] k)\right|^2 {\rm d}k\, 
r^{-3} d_{+}z\,{\rm d}r
=\int_{\tilde{\cal D}}\Phi(g)\,{\rm d}g\;,\eqno(6.2.30)$$
where  
$\tilde{\cal D}=\{ n[z] a[r] : (z,r)\in{\cal D}\} K\subset G$ and  
$\Phi : \tilde{\cal D}\rightarrow [0,\infty)$ is given by 
$\Phi(g)=|(P^{\frak a}f_{\omega})(g)|^2\,$ (for $g\in\tilde{\cal D}$). 
The mapping $(z,r)\mapsto n[z]a[r]$ is a homeomorphism from the 
compact subset ${\cal D}$ of ${\Bbb H}_3$ onto 
the subset $\{ n[z]a[r] : (z,r)\in{\cal D}\}$  of $G$, and so the latter 
set is a compact subset of $G$; since the subgroup 
$K=SU(2)$ of $G$ is also compact, one can deduce that 
the set $\tilde{\cal D}=\{ n[z] a[r] : (z,r)\in{\cal D}\} K$ 
has the Bolzano-Weierstrass property and is, therefore, a compact subset of $G$.  
Corollary~6.2.3 implies that 
$\Phi$ is continuous on $\tilde{\cal D}$. 
By the compactness of $\tilde{\cal D}$ and the continuity of $\Phi$, 
it follows that we have $\int_{\tilde{\cal D}}\Phi(g)\,{\rm d}g<\infty$ 
in the equation (6.2.30). This proves the case ${\cal X}={\cal D}$ of (6.2.29). 
\par 
Suppose now that ${\frak c}\in{\frak C}(\Gamma)$. Then, by (1.1.23), the 
definition (above (1.1.4)) of the action of $G$ on ${\Bbb H}_3$, and the fact that 
the measure ${\rm d}g$ defined in (1.1.11) is a Haar measure, one has: 
$$\int_{{\cal E}_{\frak c}}\int_{K}  
\left| \left( P^{\frak a} f_{\omega}\right) (n[z] a[r] k)\right|^2 {\rm d}k\, 
r^{-3} d_{+}z\,{\rm d}r   
=\int_{{\cal R}_{\frak c}}\int_{K} 
\ \int\limits_{1/|m_{\frak c}|}^{\infty}
\left| \left( P^{\frak a} f_{\omega}\right)\left( g_{\frak c}n[z]a[r]k\right)\right|^2 
r^{-3}\,{\rm d}r\,{\rm d}k\,d_{+}z\eqno(6.2.31)$$
(this equation being justified by the Tonelli-Hobson Test for integrabililty,  
Theorem~15.8 of [1], provided that the iterated integral on the right-hand side exists). 
Now, by (1.1.10) and (1.1.22), one has 
$$\int_{{\cal R}_{\frak c}}\int_{K} 
\ \int\limits_{1/|m_{\frak c}|}^{\infty}
\left( r^{1-\sigma}\right)^2 r^{-3}\,{\rm d}r\,{\rm d}k\,d_{+}z 
={4\over\left[\Gamma_{\frak c} : \Gamma_{\frak c}'\right]} 
\ \int\limits_{1/|m_{\frak c}|}^{\infty}
r^{-(1+2\sigma)}\,{\rm d}r <\infty\qquad\ \hbox{for $\,\sigma >0$.}$$ 
Therefore, given the definition (6.2.1) and Corollary~6.2.3, and 
given that we have both $\sigma_0>0$ and $\sigma_{\infty}>0$,  
it follows from (6.2.31),  by an application of the result (6.2.24) of Lemma~6.2.9,  
that the case ${\cal X}={\cal E}_{\frak c}$ of (6.2.29) 
holds. This completes the proof that the function $P^{\frak a} f_{\omega}$ 
lies in the space $L^2(\Gamma\backslash G)$. 
\par 
To obtain the final result of the corollary we first note that  
$P^{\frak a}f_{\omega}$ is a continuous function on $G$, and so is  
bounded on the compact set $\tilde{\cal D}=\{ n[z] a[r] : (z,r)\in{\cal D}\} K\subset G$. 
Moreover, if $\sigma_{\infty}\geq 1$ then it is implied by Lemma~6.2.9 that, for 
${\frak c}\in{\frak C}(\Gamma)$, the function $P^{\frak a}f_{\omega}$ is 
bounded on the set $\{ n[z]a[r] : (z,r)\in{\cal E}_{\frak c}\} K\subset G$. 
Given these facts, and given the definition of ${\cal F}_{*}$ in (1.1.24), 
we may conclude that in cases where $\sigma_{\infty}\geq 1$ the function  
$P^{\frak a} f_{\omega}$ is bounded on the set 
$\tilde{\cal F}_{*}=\{ n[z]a[r] : (z,r)\in{\cal F}_{*}\} K^{+}\subset G$;  
since this set $\tilde{\cal F}_{*}$ is a fundamental domain for $\Gamma\backslash G$, 
and since the function $P^{\frak a} f_{\omega} : G\rightarrow{\Bbb C}$ is $\Gamma$-automorphic, 
it therefore follows that 
$P^{\frak a} f_{\omega}$ is bounded on $G$ if $\sigma_{\infty}\geq 1$\quad$\blacksquare$ 

\bigskip 

\centerline{\bf\S 6.3 The Goodman-Wallach operator ${\bf M}_{\omega}$.}

\medskip 

Let $\nu\in{\Bbb C}$ and $p\in{\Bbb Z}$. 
Bruggeman and Motohashi, in Section~6 of [5],  
employ a method of Goodman and Wallach [12] in constructing 
a certain family of linear operators 
$({\bf M}_{\omega}^{\nu,p})_{\omega\in{\Bbb C}}$, 
with common domain $H(\nu ,p)$. 
In this subsection we detail the salient properties of 
these operators; we omit the relevant proofs, which may be found in our 
sources, [5] and [32]. 
\par 
Let $\omega\in{\Bbb C}$. Then, 
for $\varphi\in H(\nu ,p)$, Bruggeman and Motohashi 
define  ${\bf M}_{\omega}\varphi={\bf M}_{\omega}^{\nu ,p}\varphi$
by putting 
$$\left( {\bf M}_{\omega}\varphi\right) (g) 
=\sum_{m=0}^{\infty}\sum_{n=0}^{\infty} 
a_{\omega}(\nu,p;m,n) 
\left( {\partial\over\partial z}\right)^{\!m} 
\!\left( {\partial\over\partial\overline{z}}\right)^{\!n} 
\!\varphi\!\left( k[0,-1]n[z](k[0,-1])^{-1}g\right)\!\Big\vert_{z=0}\qquad\  
\hbox{($g\in G$),}\eqno(6.3.1)$$ 
where 
$$a_{\omega}(\nu ,p;m,n) 
={(\pi i\omega)^m \left(\pi i\,\overline{\omega}\right)^n\over 
(m!)(n!)\Gamma(\nu +1-p+m)\Gamma(\nu+1+p+n)}\;.\eqno(6.3.2)$$
The double series 
in (6.3.1) is absolutely convergent: this follows  
(as observed in Section~6 of [5]) by virtue of 
the function  
$(x,y)\mapsto  \varphi( k[0,-1]n[x+iy]k[0,-1]^{-1}g)$ being real analytic.  
\par 
The operator ${\bf M}_{\omega}={\bf M}_{\omega}^{\nu ,p}$ commutes with 
all elements of ${\cal U}({\frak g})$, and satisfies  
$$\left( {\bf M}_{\omega}\varphi\right) (ng) 
=\psi_{\omega}(n)\left( {\bf M}_{\omega}\varphi\right) (g)\qquad\qquad 
\hbox{for $\varphi\in H(\nu ,p)$, $g\in G$ and $n\in N$}\eqno(6.3.3)$$
(see Section~6 of [5], where the correspondence between 
(6.3.2) and (6.3.3) is established). 
Therefore one has 
${\bf M}_{\omega}\varphi_{\ell ,q}(\nu,p)\in W_{\omega}(\Upsilon_{\nu,p};\ell,q)$  
for $\ell ,q\in{\Bbb Z}$ with $\ell\geq |p|$ and $|q|\leq\ell$, so that 
$${\bf M}_{\omega} : H(\nu , p)\rightarrow 
\bigoplus_{\ell =|p|}^{\infty}\,\bigoplus_{q=-\ell}^{\ell} 
W_{\omega}\left(\Upsilon_{\nu ,p};\ell,q\right) .\eqno(6.3.4)$$
We assume, for the remainder of this subsection, that 
$\ell$ and $q$ are integers satisfying $\ell\geq\max\{ |p|,|q|\}$. 
Since $(h[u])^{-1}k[0,-1]n[z](k[0,-1])^{-1}h[u]
=k[0,-1]n[u^2 z](k[0,-1])^{-1}$, it follows directly from (6.3.1) that  
for $u\in{\Bbb C}^{*}$ one has ${\bf h}_u {\bf M}_{\omega}{\bf h}_u^{-1} 
={\bf M}_{u^2\omega}\,$ (with ${\bf h}_u$ as defined in (1.5.7)). By this and (1.8.2), one finds that 
$${\bf h}_u {\bf M}_{\omega}\varphi_{\ell,q}(\nu,p)
=|u|^{2(1+\nu)} (u/|u|)^{-2p} {\bf M}_{u^2\omega}\varphi_{\ell,q}(\nu,p)\qquad\qquad 
\hbox{($u\in{\Bbb C}^{*}$).}\eqno(6.3.5)$$
\par 
In Lemma 6.1 of [5] it is shown that if $\omega\neq 0$ then one may  
express $({\bf M}_{\omega}\varphi_{\ell,q}(\nu,p))(n a[r] k)$, 
for $n\in N$, $r>0$ and $k\in K$, as a sum of finitely many terms of the form 
$c \psi_{\omega}(n) r^{a+1}I_{\nu +a-b}(2\pi |\omega|r)\Phi_{m,q}^{\ell}(k)$, 
where $c$ is a constant and $a$, $b$ and $m$ are integers satisfying 
$a,b\geq 0$ and $|m|\leq\ell$, while 
$I_{\mu}(z)$ is the `modified' 
Bessel function $i^{-\mu}J_{\mu}(iz)$ (with $J_{\mu}(w)$ as in (1.9.6)-(1.9.8))  
and  $\Phi_{m,q}^{\ell}(k)$ is as in (1.3.2). 
Since $I_{\mu}(y)\sim (2\pi y)^{-1/2}e^y$ as $y\rightarrow +\infty$ (with $\mu\in{\Bbb C}$ a constant, and 
$y\in{\Bbb R}$), it follows from Lemma~6.1 of [5] that if $\omega\neq 0$ then the function 
$r\mapsto |({\bf M}_{\omega}\varphi_{\ell,q}(\nu,p))(n a[r]k)|$ is exponentially 
increasing as $r\rightarrow +\infty$. This implies (given (1.5.16), (1.4.9) and (1.4.11)) that  
for $\omega\neq 0$ the Jacquet integral 
${\bf J}_{\omega}\varphi_{\ell,q}(\nu,p)$ 
and `Goodman-Wallach transform' ${\bf M}_{\omega}\varphi_{\ell,q}(\nu,p)$ 
are two linearly independent functions  
that together span the space $W_{\omega}(\Upsilon_{\nu,p};\ell,q)$. 
The proof of Lemma~6.1 of [5] shows moreover that if $\omega\neq 0$ and $(\nu,p)\neq (0,0)$ then another 
basis for $W_{\omega}(\Upsilon_{\nu,p};\ell,q)$ is 
$\{ {\bf M}_{\omega}\varphi_{\ell,q}(\nu,p) , 
{\bf M}_{\omega}\varphi_{\ell,q}(-\nu,-p)\}$. 
One has, in particular, the following relation of linear dependence 
$$\eqalign{ 
(\pi |\omega|)^{-\nu} \left({-i\omega\over |\omega|}\right)^p 
\Gamma(\ell +1+\nu) 
 &{\bf J}_{\omega}\varphi_{\ell,q}(\nu,p) =\cr 
 &=\sum_{(\mu,\varpi)=\pm (\nu ,p)} 
{\pi^{2}(\pi |\omega|)^{\mu}\over\sin(-\pi\mu)}
\,\left( {i\omega\over |\omega|}\right)^{-\varpi}\Gamma(\ell +1+\mu) 
{\bf M}_{\omega}\varphi_{\ell,q}(\mu,\varpi)\;,}\eqno(6.3.6)$$ 
which is Equation~(6.15) of Lemma~6.1 of [5] (and which is, moreover,  
valid whenever $\omega\neq 0$, provided that in cases 
where $\nu$ is an integer one defines both sides of the equation via 
the relevant analytic continuation). 
\par 
The next lemma supplies a collection of useful estimates for Goodman-Wallach transforms 
on subsets of $G$ of the form $\{ a[r] k : 0<r\leq r_1, k\in K\}$, 
with $r_1$ small and positive; by (6.3.3), each such estimate 
implies an equally strong estimate on the larger set  
$\{ n a[r] k : n\in N, 0<r\leq r_1, k\in K\}$. 

\bigskip 

\proclaim Lemma 6.3.1 (Bruggeman and Motohashi). Let $\omega\neq 0$ and $r_1\in (0,\infty)$. 
Then one has 
$$\eqalignno{
\left( {\bf M}_{\omega}\varphi_{\ell,q}(\nu,p)\right) (a[r]k) 
 &={r^{1+\nu}\over\Gamma(\nu +1-p)\Gamma(\nu +1+p)}\,\Phi_{p,q}^{\ell}(k) 
+O\left( r^{2+{\rm Re}(\nu)}\right) ,&(6.3.7)\cr 
 & \cr 
\left( {\bf M}_{\omega}\varphi_{\ell,q}(0,p)\right) (a[r]k) 
 &={1\over |p|!}\pmatrix{\ell\cr\ell -|p|} 
(\pi |\omega|)^{|p|}\left( {-i\omega\over |\omega|}\right)^p 
r^{1+|p|}\Phi_{0,q}^{\ell}(k) +O\left( r^{2+|p|}\right) , &(6.3.8)}$$
uniformly for $r\in (0,r_1]$ and $k\in K$ (the implicit constant 
in (6.3.7) depends only on $\omega$, $\ell$, $\nu$ and $r_1$; that in 
(6.3.8) depends only on $\omega$, $\ell$ and $r_1$).  
If, moreover, $\sigma_0\in (0,\infty)$ then one has also 
$$\left( {\bf M}_{\omega}\varphi_{\ell,q}(\nu,p)\right) (a[r]k) 
\ll r^{1+{\rm Re}(\nu)}\left( 1+|{\rm Im}(\nu)|\right)^{-2{\rm Re}(\nu)-1} 
e^{\pi |{\rm Im}(\nu)|}\;,\eqno(6.3.9)$$
uniformly for $0<r\leq r_1$, $k\in K$ and $|{\rm Re}(\nu)|\leq\sigma_0$ 
(the implicit constant depends only on $\omega$, $\ell$, $r_1$ and $\sigma_0$). 

\medskip 

\noindent{\bf Proof.}\quad 
The results (6.3.7)-(6.3.8) follow from 
the equations~(6.13) and~(6.14) of Lemma~6.1 of [5],   
the power-series expansion 
$I_{\mu}(z)=\sum_{m\geq 0}(m!\Gamma(\mu +1+m))^{-1}(z/2)^{\mu +2m}$ 
and identities $I_{-n}(z)=I_n(z)$ ($n\in{\Bbb Z}$). 
The calculation~(7.15) of [5] makes implicit use of 
(6.3.8) and the case $\nu =1$, $p=0$ of (6.3.7). Equation~(6.3.7) and the bound~(6.3.9) are the results~(4.53) 
and~(4.55) of [32]. 
The result (6.3.9) is a corollary of  
the upper bound  for $|I_{\mu}(z)|$ stated in Relation~(1.32) of [32]\quad$\blacksquare$ 

\bigskip 

The next lemma is used in the proof of Lemma~6.5.10, where 
has a part to play in the calculation of the Fourier expansion of 
the Poincar\'{e} series $P^{\frak a}{\bf M}_{\omega}\varphi_{\ell,q}(\nu ,0)$. 

\bigskip 

\proclaim Lemma 6.3.2 (Bruggeman and Motohashi). Let $0\neq\omega_1\in{\Bbb C}$ and 
${\rm Re}(\nu)>0$. Then 
$${\bf J}_0 {\bf M}_{\omega_1}\varphi_{\ell,q}(\nu,p) 
=(-1)^p\,{\sin(\pi\nu)\over\left(\nu^2 -p^2\right)}\,
{\Gamma(\ell +1-\nu)\over\Gamma(\ell +1+\nu)}\,\varphi_{\ell,q}(-\nu,-p)\eqno(6.3.10)$$ 
and, for $\omega_2\in{\Bbb C}$ with $\omega_2\neq 0$, one has 
$${\bf J}_{\omega_2} {\bf M}_{\omega_1}\varphi_{\ell,q}(\nu,p) 
={\cal J}_{\nu,p}^{*}\left( 2\pi\sqrt{\omega_1\omega_2}\right) 
{\bf J}_{\omega_2}\varphi_{\ell,q}(\nu,p)\;,\eqno(6.3.11)$$
where, in terms of the notation defined in (1.9.5)-(1.9.6) (Theorem B), 
$${\cal J}_{\nu,p}^{*}(z)
=|z/2|^{-2\nu} (z/|z|)^{2p} {\cal J}_{\nu,p}(z) 
=J_{\nu -p}^{*}(z) J_{\nu +p}^{*}\left(\overline{z}\right)\;.\eqno(6.3.12)$$ 

\medskip 

\noindent{\bf Proof.}\quad 
This is Lemma~6.2 of [5]\quad$\blacksquare$ 

\bigskip 

\centerline{\bf\S 6.4 The Lebedev transform ${\bf L}_{\ell,q}^{\omega}$  
and auxilliary test functions.} 

\medskip

Let $0\neq\omega\in{\Bbb C}$, let $\ell,q\in{\Bbb Z}$ satisfy $\ell\geq |q|$,  
and let $\rho : G\rightarrow{\Bbb C}$ be the function defined in (6.2.1). 
Following Bruggeman and Motohashi [5], we define 
$$P_{\ell,q}(N\backslash G,\omega) 
=\bigcup_{\varepsilon\in(0,1]} 
\left\{ f\in C^{\infty}(N\backslash G,\omega) : 
f\ {\rm is\ of}\ K\!\!{\rm -type}\ (\ell,q)\ 
\,{\rm and}\ 
\sup_{g\in G}\,{|f(g)| e^{\varepsilon\left|\log\rho(g)\right|}\over \rho(g)} 
<\infty\right\}\eqno(6.4.1)$$
(this being equivalent to Equation~(7.3) of [5]). Let 
$f_{\omega}\in P_{\ell,q}(N\backslash G,\omega)$. 
Then $f_{\omega}\in C^{\infty}(N\backslash G,\omega)$, and there exist 
$\sigma_0,R_0,\sigma_{\infty},R_{\infty}\in(0,\infty)$ such that the conditions 
(6.2.5) and (6.2.22) are satisfied;   
such a choice of $\sigma_0$, $R_0$, $\sigma_{\infty}$ and $R_{\infty}$ 
is assumed in what remains of this paragraph.  
In Equation~(7.4) of [5] Bruggeman and Motohashi define their 
`Lebedev transform' 
${\bf L}_{\ell,q}^{\omega} f_{\omega} : \{ \nu\in{\Bbb C} : |{\rm Re}(\nu)|<\sigma_0\}\times 
\{ p\in{\Bbb Z} : |p|\leq\ell\}\rightarrow{\Bbb C}$ by 
$$\left({\bf L}_{\ell,q}^{\omega} f_{\omega}\right) (\nu,p) 
={1\over \pi^2\left\|\Phi_{p,q}^{\ell}\right\|_K}
\left\langle f_{\omega}\,,\, 
(\pi |\omega|)^{\overline{\nu}}\left( {-i\omega\over |\omega|}\right)^{\!\!p} 
\Gamma\left(\ell +1-\overline{\nu}\right) {\bf J}_{\omega}\varphi_{\ell,q}\left( 
-\overline{\nu} , p\right)\right\rangle_{N\backslash G}\;,\eqno(6.4.2)$$
where $\|\Phi\|_K=\sqrt{(\Phi,\Phi)_K}$ is the norm associated with the 
inner product defined in Equation~(1.2.22), while 
the inner product $\langle f , F\rangle_{N\backslash G}$ is that 
defined in (6.2.9). Given (6.2.5), (6.2.22) and the expansion of 
$({\bf J}_{\omega}\varphi_{\ell,q}(\nu,p))(g)$ obtained in Lemma~5.1 of [5],  
well known estimates for the relevant Bessel functions enable one to show  
that (6.4.2)  defines, for each $p\in\{-\ell,1-\ell,\ldots ,\ell\}$, a function 
$\nu\mapsto({\bf L}_{\ell,q}^{\omega} f_{\omega})(\nu,p)$ which is 
holomorphic in the open strip $|{\rm Re}(\nu)|<\sigma_0$. 
By the functional equation (1.7.17), one has the identity 
$( {\bf L}_{\ell,q}^{\omega} f_{\omega})(\nu,p)
=( {\bf L}_{\ell,q}^{\omega} f_{\omega})(-\nu,-p)$. 
\par 
We now describe Bruggeman and Motohashi's one-sided inversion, in Theorem~7.1 of [5],  
of their Lebedev transform operator ${\bf L}_{\ell,q}^{\omega}$. 
Let $\sigma >1$, let $S_{\sigma}=\{\nu\in{\Bbb C} : |{\rm Re}(\nu)|\leq\sigma\}$,  
and define the space ${\cal T}_{\sigma}^{\ell}$ of `test functions'   
to be the linear space of those of the functions 
$$\eta : S_{\sigma}\times\{ p\in{\Bbb Z} : |p|\leq\ell\}\rightarrow{\Bbb C}\eqno(6.4.3)$$
that satisfy all of the following three conditions:   

\medskip

(T1)\quad $\eta(\nu,p)=\eta(-\nu,-p)$; 

\medskip 

(T2)\quad for $p\in\{-\ell,1-\ell,\ldots ,\ell\}$, the function 
$\nu\mapsto\eta(\nu,p)$ can be holomorphically continued into a\break 
$\hbox{\qquad\qquad\quad}$neighbourhood of the strip $S_{\sigma}$; 

\medskip 

(T3)\quad for all $A>0$, one has $\,\eta(\nu ,p)\ll_{\eta ,A} (1+|{\rm Im}(\nu)|)^{-A} e^{-(\pi/2)|{\rm Im}(\nu)|}$. 

\medskip 

\noindent{It} is shown by Theorem~7.1 of [5] that, when $\eta\in{\cal T}_{\sigma}^{\ell}$, 
one may define $\widetilde{\bf L}_{\ell,q}^{\omega}\eta : G\rightarrow{\Bbb C}$ 
(the `inverse Lebedev transform') by putting, for 
$g\in G$, 
$$\bigl(\widetilde{\bf L}_{\ell,q}^{\omega}\eta\bigr) (g) 
={1\over 2\pi^3 i}\sum_{|p|\leq\ell} 
{(-i\omega /|\omega|)^p\over\left\|\Phi_{p,q}^{\ell}\right\|_K} 
\int\limits_{(0)} \eta(\nu,p) (\pi|\omega|)^{-\nu} \Gamma(\ell +1+\nu) 
\left( {\bf J}_{\omega}\varphi_{\ell,q}(\nu,p)\right)\!(g) 
\,\nu^{\epsilon(p)} \sin(\pi\nu) 
{\rm d}\nu\;,\eqno(6.4.4)$$
where 
$$\epsilon(p)=\cases{1 &if $p=0$;\cr -1 &otherwise.}\eqno(6.4.5)$$
One has, in particular, the following theorem. 

\bigskip 

\proclaim Theorem 6.4.1 (Bruggeman and Motohashi). Let 
$\rho : G\rightarrow (0,\infty)$ be given by (6.2.1),  
let $0\neq\omega\in{\Bbb C}$, and let ${\bf M}_{\omega}={\bf M}_{\omega}^{1,0}$ 
be the Goodman-Wallach operator on the space $H(1,0)$ (defined as in Subsection~6.3). 
Suppose moreover that $\ell,q\in{\Bbb Z}$ satisfy $\ell\geq |q|$, and that one has  
$\sigma >1$ and $\eta\in{\cal T}_{\sigma}^{\ell}$. Put  
$$b(\eta)=b(\omega;\ell,q;\eta) 
=\cases{-2\pi |\omega| \ell (\ell!)\left\|\Phi_{1,q}^{\ell}
\right\|_K^{-1} \eta(0,1) &if $\ell\geq 1$;\cr 0 &otherwise.}\eqno(6.4.6)$$
Then 
$$\widetilde{\bf L}_{\ell,q}^{\omega}\eta\in P_{\ell,q}(N\backslash G,\omega)\;,\eqno(6.4.7)$$
$$\bigl(\widetilde{\bf L}_{\ell,q}^{\omega}\eta\bigr)(g) 
\ll_{\omega ,\eta ,A}\ (\rho(g))^{-A}\qquad\quad\    
\hbox{($A\in[-2,\infty)$, $g\in G$)}\eqno(6.4.8)$$
and, when $g\in G$,  
$$\bigl( \widetilde{\bf L}_{\ell,q}^{\omega}\eta\bigr)(g) 
=b(\eta)\left( {\bf M}_{\omega}\varphi_{\ell,q}(1,0)\right)(g) 
+O_{\omega ,\eta}\!\left( 
(\rho(g))^{\min\{ 1+\sigma , 3\}}\right)\qquad\  
\hbox{if $\,\rho(g)\leq 1$.}\eqno(6.4.9)$$
Moreover, for all $(\nu,p)\in{\Bbb C}\times\{ -\ell,1-\ell,\ldots ,\ell\}$ 
such that $|{\rm Re}(\nu)|<1$, one has 
$$\bigl( {\bf L}_{\ell,q}^{\omega}\widetilde{\bf L}_{\ell,q}^{\omega}\eta\bigr)(\nu,p) 
=-{2\over\pi}\,\Gamma(\ell +1-\nu)\Gamma(\ell +1+\nu)\,{\sin(\pi\nu)\over\pi\nu}\, 
{\nu^{1+\epsilon(p)}\over\left(\nu^2 -p^2\right)}\,\eta(\nu,p)\;,\eqno(6.4.10)$$ 
where $\epsilon(p)$ is given by Equation~(6.4.5). 

\medskip 

\noindent{\bf Proof.}\quad 
The results (6.4.7) and (6.4.10) are contained in Theorem~7.1 of [5]  
(see also the extension of that theorem obtained in Theorem~9.1.4 of [32]). 
\par 
The result (6.4.9) is a corollary of the combination of 
Equation~(7.14) of [5] and the calculation implicit in [5],~(7.15): note that the 
the conditions (T2) and (T3) on $\eta\in{\cal T}_{\sigma}^{\ell}$, and the 
bound (6.3.9) of Lemma~6.3.1, enable one to substitute 
$\min\{\sigma,2\}$ for the $\alpha$ in Equation~(7.14) of [5], and so   
to deduce that the first sum on the right-hand side of that equation is  
not greater than 
$O_{\omega,\eta}\bigl( (\rho(g))^{1+\min\{\sigma,2\}}\bigr)$ when $\rho(g)\leq 1$. 
The definition of the constant $b(\eta)$ given in (6.4.6) represents a 
slight correction of the corresponding definition below [5], (7.15). 
We have used (6.3.3) and the equations (6.3.7) and (6.3.8) of 
Lemma~6.3.1 in performing our own check upon the calculation (7.15) of [5]. 
\par 
The case $\rho(g)\leq 1$ of (6.4.8) is 
a corollary of (6.4.9), (6.4.6) and the case $\nu =1$, $p=0$ of the 
estimate (6.3.7) of Lemma~6.3.1.
The case $\rho(g)\geq 1$ of (6.4.8) is a consequence of 
the equations~(5.26)-(5.27) of Lemma~5.1 of [5],  
and is noted within the proof of Theorem~7.1 of [5]: see  the discussion 
around the bound (7.12) of [5] for the Bessel function $K_{\xi}(u)\quad\blacksquare$ 

\bigskip 

In [5], Lemma~7.1 and Lemma~7.2, Bruggeman and Motohashi investigate the
inverse Lebedev transform $\widetilde{\bf L}_{\ell,q}^{\omega}$ further. 
We reproduce those two lemmas here (without the proofs), 
since the results they contain are needed for the 
proof of the spectral Kloosterman sum formula (Theorem B). Before stating these lemmas 
we clarify that henceforth $L^2(N\backslash G)$ denotes the space of those 
measurable functions $f : G\rightarrow{\Bbb C}$ that satisfy $f(ng)=f(g)$, for all 
$n\in N$, $g\in G$, and are such that $\int_{N\backslash G} |f(g)|^2\,{\rm d}\dot g<\infty$,  
where the measure ${\rm d}\dot g$ is that which occurs in (6.2.9); 
the inner product defined 
in (6.2.9) makes $L^2(N\backslash G)$ a Hilbert space. 

\bigskip 

\proclaim Lemma 6.4.2 (Bruggeman and Motohashi). 
Let $\ell\in{\Bbb N}\cup\{0\}$ and $\sigma\in (1,\infty)$,  
and let $\eta,\theta\in{\cal T}_{\sigma}^{\ell}$. 
Then $|\widetilde{\bf L}_{\ell,q}^{\omega}\eta|\in L^2(N\backslash G)$ and, 
for $0\neq\omega\in{\Bbb C}$ and 
$q=-\ell,-\ell+1,\ldots ,\ell$, one has 
$$\left\langle\widetilde{\bf L}_{\ell,q}^{\omega}\eta\,,
\,\widetilde{\bf L}_{\ell,q}^{\omega}\theta\right\rangle_{N\backslash G} 
=-{1\over\pi i}\sum_{|p|\leq\ell}\ \int\limits_{(0)} 
\eta(\nu,p)\,\overline{\theta(\nu,p)}\, 
\,\Gamma(\ell+1-\nu) \Gamma(\ell+1+\nu)\,  
{\sin^2(\pi\nu)\over\pi^2\nu^2}\,{\nu^{2+2\epsilon(p)}\over\left(\nu^2 -p^2\right)}\,
{\rm d}\nu\;,\eqno(6.4.11)$$
where $\langle f , F\rangle_{N\backslash G}$ and $\epsilon(p)$ 
are as defined in (6.2.9) and~(6.4.5).  

\bigskip 

\goodbreak\proclaim Lemma 6.4.3 (Bruggeman and Motohashi). Let $\ell\in{\Bbb N}\cup\{0\}$, 
$\sigma\in (1,2)$ and $\omega_1,\omega_2,c\in{\Bbb C}-\{ 0\}$,  
let $S_{\sigma}=\{ \nu\in{\Bbb C} : |{\rm Re}(\nu)|\leq\sigma\}$,  
and let $\kappa(\omega_1,\omega_2;c)$ be that mapping from 
${\cal T}_{\sigma}^{\ell}$ into the space of all functions 
$f : S_{\sigma}\times\{ p\in{\Bbb Z} : |p|\leq\ell\}\rightarrow {\Bbb C}$
which is given by 
$$\left(\kappa\left(\omega_1,\omega_2;c\right)\eta\right)(\nu,p) 
={\cal K}_{\nu,p}\left( 2\pi\sqrt{\omega_1\omega_2}/c\right) 
\eta(\nu ,p)\qquad\quad\hbox{($\eta\in{\cal T}_{\sigma}^{\ell}$, 
$\nu\in S_{\sigma}$, $p\in{\Bbb Z}$ and $|p|\leq\ell$),}\eqno(6.4.12)$$
where ${\cal K}_{\nu,p}(u)\in{\Bbb C}$ is defined 
by the equations (1.9.4)-(1.9.6) of Theorem~B. Then $\kappa(\omega_1,\omega_2;c)$ 
is a linear operator from ${\cal T}_{\sigma}^{\ell}$ into ${\cal T}_{\sigma}^{\ell}$,  
and, for $\eta\in{\cal T}_{\sigma}^{\ell}$, $q=-\ell,1-\ell,\ldots ,\ell$ and 
$g\in G$,  
one has (with ${\bf h}_u$ as in (1.5.7)): 
$$\bigl({\bf J}_{\omega_2}{\bf h}_{1/c}\widetilde{\bf L}_{\ell,q}^{\omega_1}\eta\bigr)(g) 
=\pi^2 |c|^{-2}\bigl(\widetilde{\bf L}_{\ell,q}^{\omega_2}
\kappa(\omega_1,\omega_2;c)\eta\bigr)(g)\;.\eqno(6.4.13)$$

\medskip 

\noindent{\bf Remark~6.4.4.}\quad 
The first result of Lemma~6.4.3 (i.e. that 
$\kappa(\omega_1,\omega_2;c)\eta\in {\cal T}_{\sigma}^{\ell}$ for all 
$\eta\in{\cal T}_{\sigma}^{\ell}$) has been taken from Lemma~9.1.8 of [32]: its proof 
in respect of cases where one has $1<\sigma<3/2$ 
is indicated there (see, in particular, the upper bound on $|{\cal K}_{\nu,p}(z)|$ 
obtained in Lemma~9.1.7 of [32]), and a trivial extension of this proof  
yields the required result in the remaining cases, where $3/2\leq\sigma<1$.  
The proof of (6.4.13) is, in part, an application of 
Lemma~6.3.2 and the identity (6.3.6). 

\bigskip 

\centerline{\bf \S 6.5 Poincar\'{e} series revisited.} 

\medskip 

It is to be assumed throughout this subsection 
that $\omega$, $\sigma$, $\eta$, 
${\frak a}$, ${\frak b}$, $g_{\frak a}$, $g_{\frak b}$ and 
the $K$-type $(\ell,q)$ are given, with  
$0\neq\omega\in{\frak O}$, $\sigma\in (1,2)$, $\ell,q\in{\Bbb Z}$, $|q|\leq\ell$ and    
$\eta\in{\cal T}_{\sigma}^{\ell}$ (the space defined in, and below, 
(6.4.3)), and with ${\frak a},{\frak b}\in{\Bbb Q}(i)\cup\{\infty\}$ and  
$g_{\frak a},g_{\frak b}\in G$ such that (1.1.16) and (1.1.20)-(1.1.21) hold for 
${\frak c}\in\{ {\frak a} , {\frak b}\}$. We suppose also 
that $\rho : G\rightarrow{\Bbb C}$ is the function given by (6.2.1). 
\par 
Were it guaranteed to be absolutely convergent, for all $g\in G$, 
the Poincar\'e series $(P^{\frak a}\widetilde{\bf L}^{\omega}_{\ell,q}\eta)(g)$ 
might, by itself, serve as the principal `fulcrum' in a proof of    
the spectral sum formula, Theorem~B: for the related identities (6.4.10), 
(6.4.11) and (6.4.13) are key results for the proof of Theorem~B that we 
are going to describe. However, in view of (6.4.6), the estimate  
(6.4.9) and the case $\nu =1$, $p=0$ of (6.3.7), it follows by [11], Corollary~3.1.6 and Proposition~3.2.1~(2),   
that the series $(P^{\frak a}\widetilde{\bf L}^{\omega}_{\ell,q}\eta)(g)$
will be absolutely convergent, for all $g\in G$,  only if one has 
either $\ell =0$, or else $\ell\geq 1$ and $\eta(0,1)=0$.  
This inconvenient fact led Bruggeman and Motohashi   
to construct, in the equations~(9.7) and~(9.8) of [5], a suitable `substitute' 
for $P^{\infty}\widetilde{\bf L}^{\omega}_{\ell,q}\eta$. The construction of this 
`substitute' is a critical step in their proof of Theorem~10.1 of [5]
(the spectral sum formula for $PSL(2,{\frak O})\backslash 
PSL(2,{\Bbb C})$). In this subsection we adapt the method of Bruggeman and Motohashi 
in defining the corresponding `substitute' $P^{{\frak a},*}\widetilde{\bf L}^{\omega}_{\ell,q}\eta$ 
for the Poincar\'e series 
$P^{\frak a}\widetilde{\bf L}^{\omega}_{\ell,q}\eta$. 
\par 
We remark that the 
choice of Poincar\'{e} series in Section~9 of~[5] is 
acknowledged by Bruggeman and Motohashi to have been influenced 
by the widely applicable method which Miatello and Wallach 
developed in [33] and [34]. Consequently our work here has been indirectly 
influenced by those papers of Miatello of Wallach. In particular, we might have 
obtained our Lemma~6.5.16 by an application of Theorem~2.5 of~[33] and the 
result (1.4.15) of Kim and Shahidi:  
Lokvenec-Guleska has shown how to do this when ${\frak a}=\infty\,$ 
(see Section~9.2 of~[32]). We have instead chosen to  
rely on bounds for generalised Kloosterman sums 
(i.e. the results of Subsection~6.1). This is essentially 
what Bruggeman and Motohashi chose to do in Section~9 of~[5]; 
it enables us to 
give proofs which are (on the whole) more self-contained than 
would be the case if we had opted to derive 
our results from the results of Miatello and Wallach. 
\par 
We choose, once and for all, a function $\tau\in C^{\infty}(G)$ such 
that, for $n\in N$, $r>0$ and $k\in K$, one has 
$$[0,1]\ni\tau(n a[r] k)=\tau( a[r])=\cases{1 &if $r\leq 1$;\cr 0 &if $r\geq 2$.}
\eqno(6.5.1)$$ 
Then, motivated by the estimate (6.4.9), we put 
$$\widetilde{\bf L}^{\omega,*}_{\ell,q}\eta 
=\widetilde{\bf L}^{\omega}_{\ell,q}\eta  - b(\omega;\ell,q;\eta) 
\,\tau{\bf M}_{\omega}\varphi_{\ell,q}(1,0)\;,\eqno(6.5.2)$$
where $b(\omega;\ell,q;\eta)$ is the constant in (6.4.6), and where 
$$\left( \tau{\bf M}_{\omega}\varphi_{\ell,q}(\nu,0)\right)(g) 
=\tau(g)\left( {\bf M}_{\omega}\varphi_{\ell,q}(\nu,0)\right)(g)\qquad\quad 
\hbox{for $\nu\in{\Bbb C}$, $g\in G$}\eqno(6.5.3)$$
(we shall use the similar notation $(1-\tau){\bf M}_{\omega}\varphi_{\ell,q}(\nu,0)$ 
to denote the function $g\mapsto (1-\tau(g))({\bf M}_{\omega}\varphi_{\ell,q}(\nu,0))(g)$). 
Note that, by (6.5.2) and (6.4.6), one has 
$\widetilde{\bf L}^{\omega,*}_{\ell,q}\eta =\widetilde{\bf L}^{\omega}_{\ell,q}\eta$ 
if $\ell =0$, or if $\ell\geq 1$ and $\eta(0,1)=0$. Both the transform 
$(\widetilde{\bf L}^{\omega,*}_{\ell,q}\eta)(g)$ and the absolute value of this transform 
may be used in constructing 
Poincar\'e series: for, by (6.5.2), the result (6.4.7) of Theorem~6.4.1, 
our choice of $\tau\in C^{\infty}(G)$ (satisfying (6.5.1)), the observation 
preceding (6.3.4) and the definitions (1.4.7) and (1.4.3), it follows that 
$$\widetilde{\bf L}^{\omega,*}_{\ell,q}\eta\in C^{\infty}(N\backslash G,\omega)\qquad\ 
{\rm and}\qquad\ \bigl|\widetilde{\bf L}^{\omega,*}_{\ell,q}\eta\bigr|\in 
C^0(N\backslash G,0)\;.\eqno(6.5.4)$$
\par 
Using what is effectively the same construction as occurs in 
the equations~(9.7)-(9.9) of [5], we now define the function 
$P^{{\frak a},*}\widetilde{\bf L}^{\omega}_{\ell,q}\eta : G\rightarrow{\Bbb C}$ 
by putting, for $g\in G$, 
$$\bigl( P^{{\frak a},*}\widetilde{\bf L}^{\omega}_{\ell,q}\eta\bigr)(g) 
=\lim_{\nu\rightarrow 1+} 
\left( P^{\frak a}\Bigl(\widetilde{\bf L}^{\omega,*}_{\ell,q}\eta 
+b(\omega;\ell,q;\eta)\,\tau{\bf M}_{\omega}\varphi_{\ell,q}(\nu,0)\Bigr)\right)(g)\;,\eqno(6.5.5)$$ 
where (as we shall henceforth suppose) the 
function $\tau{\bf M}_{\omega}\varphi_{\ell,q}(\nu,0)$ is given by~(6.5.3). 
We have, therefore, to show that the limit in Equation~(6.5.5) exists. 
This will be achieved by means of the analytic continuation, with 
respect to the complex variable $\nu$, of 
the Poincar\'e series $(P^{\frak a}\tau{\bf M}_{\omega}\varphi_{\ell,q}(\nu ,0))(g)$. 
We also aim to show that the `pseudo Poincar\'e series' 
$P^{{\frak a},*}\widetilde{\bf L}^{\omega}_{\ell,q}\eta$ given by (6.5.5) has the specific 
properties that will enable us to make use of it in proving Theorem~B. 
\par 
We find some new terminology convenient in stating the lemmas which follow. 
The space of all of the 
measurable $\Gamma$-automorphic functions 
$f : G\rightarrow{\Bbb C}$ that are essentially bounded is denoted by 
$L^{\infty}(\Gamma\backslash G)$. We define also: 
$$L^{\infty}(\Gamma\backslash G;\ell,q) 
=\left\{ f\in L^{\infty}(\Gamma\backslash G)\,: f\ {\rm is\ of}\ 
K{-\rm type}\ (\ell,q)\,\right\}\;.\eqno(6.5.6)$$ 
Since ${\rm vol}(\Gamma\backslash G)<\infty$,  one has 
$$L^{p_1}(\Gamma\backslash G)\supseteq L^{p_2}(\Gamma\backslash G)\qquad\qquad  
\hbox{($1\leq p_1\leq p_2\leq\infty$),} 
\eqno(6.5.7)$$  
where $L^p(\Gamma\backslash G)$ denotes the space of those of the 
measurable and $\Gamma$-automorphic functions $f : G\rightarrow{\Bbb C}$ 
which are such that $\int_{\Gamma\backslash G}|f|^p\,{\rm d}g<\infty$.

\bigskip 

\proclaim Lemma 6.5.1. Let ${\rm Re}(\nu)>1$. 
Then $P^{\frak a}{\bf M}_{\omega}\varphi_{\ell,q}(\nu,0)\in C^0(\Gamma\backslash G)$,    
$P^{\frak a}|\tau{\bf M}_{\omega}\varphi_{\ell,q}(\nu,0)|\in 
C^0(\Gamma\backslash G)$ and 
$P^{\frak a}|\tau{\bf M}_{\omega}\varphi_{\ell,q}(\nu,0)|\in L^{\infty}(\Gamma\backslash G)$. 

\medskip 

\noindent{\bf Proof.}\ 
Let $F={\bf M}_{\omega}\varphi_{\ell,q}(\nu ,0)$. By the observation preceding (6.3.4), 
$F\in W_{\omega}(\Upsilon_{\nu,0};\ell,q)\subset C^0(N\backslash G,\omega)$, and 
so, given the choice of $\tau$ in and above (6.5.1), we have also 
$|\tau F|\in C^0(N\backslash G,0)$. The estimate~(6.3.7) of Lemma~6.3.1 shows, moreover, 
that the case $\sigma_0={\rm Re}(\nu)$, $R_0=1$ (say) of the condition (6.2.5) is 
satisfied when $f_{\omega}=F$; and obviously the same is true 
if one substitutes $|\tau F|$ for $F$ (since $|\tau(g)|\leq 1$ for all $g\in G$). 
Therefore, given that we have ${\rm Re}(\nu)>1$, the hypotheses of Lemma~6.2.2 
are satisfied when one takes there: $\sigma_0={\rm Re}(\nu)$, $R_0=1$ and 
either $f_{\omega}=F$, or $\omega =0$ and $f_0=|\tau F|$. Hence the 
first two results of the lemma follow by virtue of Corollary~6.2.3. 
\par 
We now have only to prove that 
$P^{\frak a}|\tau F|\in L^{\infty}(\Gamma\backslash G)$. 
Given the results of the lemma that have already been proved, it suffices  
to show that the function 
$P^{\frak a}|\tau F| : G\rightarrow [0,\infty)$ 
is bounded. One may deduce this from Corollary~6.2.10: 
for (6.5.1) trivially implies that 
$(\tau F)(g)=\tau(g) F(g)=0 F(g)=0$ for all $g\in G$ such that $\rho(g)\geq 2$, 
and so the condition (6.2.22) is satisfied when 
$f_{\omega}=|\tau F|$, $R_{\infty}=2$ and $\sigma_{\infty}=1$ (say)
\quad$\blacksquare$ 

\bigskip 

\proclaim Lemma 6.5.2. For $\nu\in{\Bbb C}$, the function 
$P^{\frak a}(1-\tau){\bf M}_{\omega}\varphi_{\ell,q}(\nu,0) : G\rightarrow{\Bbb C}$ 
is defined, lies in $\ C^{\infty}(\Gamma\backslash G)$, and is of 
$K$-type $(\ell,q)$. 
The mapping 
$(\nu,g)\mapsto (P^{\frak a}(1-\tau){\bf M}_{\omega}\varphi_{\ell,q}(\nu,0))(g)$ 
is a continuous function on ${\Bbb C}\times G$; and, for each $g\in G$, 
the complex function 
$\nu\mapsto (P^{\frak a}(1-\tau){\bf M}_{\omega}\varphi_{\ell,q}(\nu,0))(g)$ is entire. 
One has, moreover, 
$$\left( P^{\frak a}(1-\tau){\bf M}_{\omega}\varphi_{\ell,q}(\nu,0)\right)\left( 
g_{\frak b} g\right) 
={1\over\left[\Gamma_{\frak a} : \Gamma_{\frak a}'\right]} 
\sum_{\scriptstyle\gamma\in\Gamma_{\frak a}'\backslash\Gamma\atop\scriptstyle 
\gamma{\frak b}={\frak a}} 
\left( {\bf M}_{\omega}\varphi_{\ell,q}(\nu,0)\right)\left( 
g_{\frak a}^{-1}\gamma g_{\frak b} g\right)\quad\ 
\hbox{if $\,\rho(g)\geq 2$.}\eqno(6.5.8)$$

\medskip 

\noindent{\bf Proof.}\quad 
Let ${\nu}\in{\Bbb C}$. Put $f_{\omega}=(1-\tau){\bf M}_{\omega}\varphi_{\ell,q}(\nu,0)$, where 
the meaning of the term `$(1-\tau){\bf M}_{\omega}\varphi_{\ell,q}(\nu,0)$' is that 
indicated below (6.5.3).
Then, given our choice of $\tau$ (described in, and above, (6.5.1)), 
it follows by the observation preceding (6.3.4) that we have:   
$$f_{\omega}\in C^{\infty}(N\backslash G,\omega)\subset C^{\infty}(G)\;,\qquad 
\hbox{$f_{\omega}$ is of $K$-type $(\ell,q)$}\eqno(6.5.9)$$
and 
$$f_{\omega}(g)=\cases{\left( {\bf M}_{\omega}\varphi_{\ell,q}(\nu,0)\right)(g) &if $\rho(g)\geq 2$;\cr 
0 &if $\rho(g)\leq 1$.}\eqno(6.5.10)$$ 
\par 
By (1.5.4), we have also  
$$\left( P^{\frak a}f_{\omega}\right)\left( g_{\frak b}g\right) 
={1\over\left[\Gamma_{\frak a}:\Gamma_{\frak a}'\right]} 
\sum_{\gamma\in\Gamma_{\frak a}'\backslash\Gamma} f_{\omega}\left( 
g_{\frak a}^{-1}\gamma g_{\frak b} g\right)\;,\eqno(6.5.11)$$
for any $g\in G$ such that the sum on the right-hand side of (6.5.11) converges. 
Suppose that $Z_1>0$ and $R_1>R_0>0$. Put 
$U=\{ (z,r)\in{\Bbb H}_3 : |z|\leq Z_1\ {\rm and}\ R_0\leq r\leq R_1\}$ and 
$\tilde U=\{ n[z]a[r] : (z,r)\in~\!\!U\}K$. 
Then, when $g\in\tilde U$, it follows by (6.5.10) and Lemma~6.2.1 that 
each non-zero term occurring in the sum in (6.5.11) corresponds to a  
$\gamma$ contained in the set $(\Gamma_{\frak a}'\backslash I_0)\cup(\Gamma_{\frak a}'\backslash I_1(g))$, 
where $I_0={}^{\frak a}\Gamma^{\frak b}(0)=\{\gamma\in\Gamma : 
\gamma{\frak b}={\frak a}\}\,$  
(so that $\Gamma_{\frak a}'\backslash I_0$ is a finite set, empty 
unless ${\frak a}\sim^{\!\!\!\!\Gamma}{\frak b}$), while 
$I_1(g)=\{\gamma\in\Gamma : \gamma{\frak b}\neq{\frak a}\ {\rm and}\ 
1<\rho(g_{\frak a}^{-1}\gamma g_{\frak b}g)\leq 1/R_0\}$. 
\par 
If, in particular, $R_0=2$ and $g\in\tilde U$, then the equality in (6.5.8) holds:  
for in this case the set $I_1(g)$ is empty, and it 
moreover follows by Lemma~6.2.1 and (6.5.10) that 
for all $\gamma\in I_0={}^{\frak a}\Gamma^{\frak b}(0)$ one has 
$\rho(g_{\frak a}^{-1}\gamma g_{\frak b}g)=\rho(g)\geq R_0=2$,  
so that the first case of (6.5.10) applies to 
$f_{\omega}(g_{\frak a}^{-1}\gamma g_{\frak b}g)$. Since $G=NAK$, and 
since the assumptions concerning  
$Z_1$, $R_0$ and $R_1$ are just that $Z_1>0$ and $R_1>R_0=2$, the above therefore  
completes the proof of (6.5.8). 
\par 
We now revert to considering the more general case of 
any $Z_1>0$, and any $R_1$ and $R_0$ 
with $R_1>R_0>0$. However, since the results that remain to be proved are 
independent of the cusp ${\frak b}$, it will be convenient to be 
more specific in another respect, by assuming 
henceforth that  ${\frak b}=\infty$, and that $g_{\frak b}=g_{\infty}=h[1]$.  
\par 
Returning  
to the matters discussed in the paragraph containing Equation (6.5.11), 
we note that each $\gamma\in\Gamma_{\frak a}'\backslash I_1(g)$ represents a coset 
$\Gamma_{\frak a}'\gamma_{*}$ (say), where $\gamma_{*}\in I_1(g)$ may 
(by (1.1.20)-(1.1.21) for ${\frak c}={\frak a}$) be chosen so that one has:  
$$g_{\frak a}^{-1}\gamma_{*}g_{\infty}g\in\tilde V 
=\left\{ n[z] a[r] : (z,r)\in V\right\}K\;,$$
where 
$V=\{ (z,r)\in{\Bbb H}_3 : |{\rm Re}(z)|,|{\rm Im}(z)|\leq 1/2\ {\rm and}\ 
1\leq r\leq 1/R_0\}$. 
Hence, when $g\in\tilde U$, the summation in (6.5.11) is effectively 
over all $\gamma\in \Gamma_{\frak a}'\backslash\Gamma$ that are either 
contained in the finite set $\Gamma_{\frak a}'\backslash I_0$, or else represent 
cosets of the form $\Gamma_{\frak a}'\gamma_{*}$ with $\gamma_{*}\in\Gamma$ 
such that $(\gamma_{*}g_{\infty}\tilde U)\cap(g_{\frak a}\tilde V)\neq\emptyset$. 
Given the definition of the action of $G$ upon ${\Bbb H}_3$, 
and given that $g_{\infty}=h[1]\,$ (the identity element of $G$), 
the latter 
condition is equivalent to the condition that
$(\gamma_{*}U)\cap(g_{\frak a}V)\neq\emptyset$, and so it is 
satisfied only if one has 
$$\left(\gamma_{*}W\right)\cap W\neq\emptyset\;,\eqno(6.5.12)$$
where $W=U\cup(g_{\frak a}V)$. Since this set $W$ is 
a compact subset of ${\Bbb H}_3$, and since (by Theorem~2.1.2 of [11]) 
the discrete group $\Gamma<SL(2,{\Bbb C})$ is discontinuous, it follows that 
the condition (6.5.12) is satisfied for at most finitely many 
choices of $\gamma_{*}\in\Gamma$. The condition (6.5.12) is, furthermore, 
independent of the variable $g\in\tilde U$ (as is the finite set 
$\Gamma_{\frak a}'\backslash I_0$).  Therefore (and since 
we have $g_{\infty}g=h[1]g=g$ for all $g\in G$)  
there exists a finite set $\{\gamma_1,\ldots ,\gamma_J\}\subset\Gamma$ 
such that 
$$\left( P^{\frak a}f_{\omega}\right)(g) 
={1\over\left[\Gamma_{\frak a}:\Gamma_{\frak a}'\right]} 
\sum_{j=1}^J f_{\omega}\left( 
g_{\frak a}^{-1}\gamma_j g\right)\qquad\quad 
\hbox{for all $g\in\tilde U$.}\eqno(6.5.13)$$
\par 
Since $G=NAK$, and since $Z_1>0$ and $R_1>R_0>0$ are arbitrary, the 
fact that we obtain (6.5.13) is enough to prove that, for each $g\in G$, 
the Poincar\'e series $(P^{\frak a} f_{\omega})(g)$ is absolutely convergent 
(by virtue of it being a series that contains only finitely many non-zero terms). 
The mapping $g\mapsto (P^{\frak a} f_{\omega})(g)$ is therefore (given the 
definition (1.5.4) and the first 
part of (6.5.9)) a well-defined $\Gamma$-automorphic function on~$G$. 
Similarly, since it follows by (6.5.13), (6.5.9), (1.3.1) and the left-invariance of all  
elements of ${\cal U}({\frak g})$ and ${\cal U}({\frak k})$ that 
the restriction of $P^{\frak a}f_{\omega}$ to the open set $\ddot U={\rm Int}(\tilde U)$ is 
a smooth function (i.e. an element of $C^{\infty}(\ddot U)\,$) 
satisfying $({\bf H}_2 P^{\frak a} f_{\omega})(g)=-iq (P^{\frak a} f_{\omega})(g)$ 
and $(\Omega_{\frak k} P^{\frak a} f_{\omega})(g)
=-{1\over 2}\,\left(\ell^2+\ell\right)(P^{\frak a} f_{\omega})(g)$ 
for $g\in{\ddot U}$, we may deduce that $P^{\frak a}f_{\omega}$ is in fact 
a smooth function on $G$ of $K$-type $(\ell,q)$. 
\par 
By (6.5.13) one moreover obtains, as a corollary of the 
expansion of $({\bf M}_{\omega}\varphi_{\ell,q}(\nu,p))(g)$ in terms of Bessel functions 
that is given by the equations~(6.13) and~(6.14) of Lemma~6.1 of [5],  
the result that, for each $g\in G$, the mapping 
$\mu\mapsto (P^{\frak a} (1-\tau){\bf M}_{\omega}\varphi_{\ell,q}(\mu,0))(g)$ 
is an entire complex function. The same expansion of 
$({\bf M}_{\omega}\varphi_{\ell,q}(\nu,p))(g)$ in terms of Bessel function 
also enables one to show, in particular, that 
the mapping $(\mu ,g)\mapsto({\bf M}_{\omega}\varphi_{\ell,q}(\mu,0))(g)$ 
is a continuous function on ${\Bbb C}\times G$. Hence and by (6.5.13) one finds that,   
since $G$ is a topological group, and since $(1-\tau)\in C^{\infty}(G)$,  the mapping 
$(\mu,g)\mapsto (P^{\frak a} (1-\tau){\bf M}_{\omega}\varphi_{\ell,q}(\mu,0))(g)$ is a continuous function on ${\Bbb C}\times G$
\quad$\blacksquare$ 

\bigskip 

\proclaim Lemma 6.5.3. Let $0\neq\omega'\in{\Bbb C}$,  
and let $p\in{\Bbb Z}$ satisfy 
$|p|\leq\ell$. Then, for all $g\in G$, each of 
the two mappings $\nu\mapsto ({\bf J}_{\omega'}\varphi_{\ell,q}(\nu,p))(g)$ and 
$\nu\mapsto({\bf M}_{\omega'}\varphi_{\ell,q}(\nu,p))(g)$ is an 
entire function of the complex variable $\nu$. Suppose moreover that 
$r_0\in(0,\infty)$, $\nu\in{\Bbb C}$, $|{\rm Re}(\nu)|\leq\sigma_1<\infty$
and $\omega'\in{\frak O}$ (so that, in particular, $|\omega'|\geq 1$), 
and that $n\in N$, $r>0$, $k\in K$ and   
$g=n a[r]k$ (so that $\rho(g)=r$). Put  
$$f_{\omega'}(\nu,p;g)=\left(\pi |\omega'|\right)^{-\nu}(i\omega'/|\omega'|)^p
\Gamma(\ell +1+\nu) \bigl( {\bf J}_{\omega'}\varphi_{\ell,q}(\nu,p)\bigr) (g)\;.$$ 
Then 
$$f_{\omega'}(\nu,p;g) 
\ll_{\ell,\sigma_1,r_0}\,\left|\omega'\right|^{-1} 
\left( 1+|{\rm Im}(\nu)|\right)^{\ell-|p|}\left|\omega' r\right|^{\ell +1} 
e^{-2\pi |\omega'| r}\qquad\hbox{if $\,|\omega'|r\geq r_0$.}\eqno(6.5.14)$$
Moreover, for all $\varepsilon\in (0,1/4]$ and all $d\in{\Bbb N}$ 
such that $d/2>\sigma_1 +\ell$, one has  
$$f_{\omega'}(\nu,p;g)
={r\over e^{(\pi /2)|{\rm Im}(\nu)|}}\cdot 
\cases{O_{\ell,\sigma_1,r_0,\varepsilon} 
\!\left( (1+|{\rm Im}(\nu)|)^{|{\rm Re}(\nu)|-1/2+\ell} 
\left|\omega' r\right|^{-|{\rm Re}(\nu)|-\varepsilon}\right) 
 &if $\,|\omega'|r\leq r_0$,\cr 
O_{d,\sigma_1,r_0}\!\left(
(1+|{\rm Im}(\nu)|)^{{\rm Re}(\nu)-1/2+\ell +d} 
\left|\omega' r\right|^{-{\rm Re}(\nu)+\ell +|p|-d}\right) 
 &if $\,|\omega'| r\geq r_0$.}\eqno(6.5.15)$$ 
The case $|\omega'|r\geq r_0$ of (6.5.15) remains valid if,  
in place of the $O$-term 
appearing there, one substitutes the term   
$O_{d,\sigma_1,r_0}(\min\{ (1+|{\rm Im}(\nu)|)^{-{\rm Re}(\nu)-1/2+\ell +d} 
\left|\omega' r\right|^{{\rm Re}(\nu)+\ell +|p|-d} , 
(1+|{\rm Im}(\nu)|)^{-1/2+\ell +d} 
\left|\omega' r\right|^{\ell +|p|-d}\})$. 

\medskip 

\noindent{\bf Proof.}\quad 
The first assertion, concerning the mappings 
$\nu\mapsto ({\bf J}_{\omega'}\varphi_{\ell,q}(\nu,p))(g)$, 
$\,\nu\mapsto({\bf M}_{\omega'}\varphi_{\ell,q}(\nu,p))(g)$, is a corollary of 
the relevant expansions in terms of Bessel functions $K_{\mu}(2\pi r)$ and $I_{\mu}(2\pi r)$ 
obtained in [5], Lemma~5.1 and Lemma~6.1. 
\par 
The result (6.5.15) is proved similarly to the equation~(4.28) of Lemma~4.1.3 of [32], 
and (see our Remark~6.5.4, following this proof) coincides with that 
result in the respect of cases with $|\omega'|r\leq r_0$ and 
${\rm Re}(\nu)\leq 0$.
Its proof (which we omit) involves the application of 
the equations~(5.26)-(5.27) of Lemma~5.1 of [5],  and the estimates    
$$K_{\mu}(2\pi R) 
\ll_{\varepsilon,\sigma_2,r_0}
\,e^{-(\pi /2)|{\rm Im}(\mu)|} 
(1+|{\rm Im}(\mu)|)^{|{\rm Re}(\mu)|-1/2} 
R^{-|{\rm Re}(\mu)|-\varepsilon}\qquad\quad   
\hbox{($R\in (0 , r_0]$, $|{\rm Re}(\mu)|\leq\sigma_2<\infty$)}$$
and 
$$K_{\mu}(2\pi R) 
\ll_{d,\sigma_2}
\,e^{-(\pi /2)|{\rm Im}(\mu)|} 
(1+|{\rm Im}(\mu)|)^{{\rm Re}(\mu)+d} 
R^{-{\rm Re}(\mu)-d}\qquad\quad 
\hbox{($R\in(0,\infty)$, $|{\rm Re}(\mu)|\leq\sigma_2<d/2$),}$$
where $K_{\mu}(z)$ is the modified Bessel function defined in 
the equations~10.27.4-10.27.5 of [38]; the first of these estimates  is 
[32],~(1.33); the second is [32],~(1.37) (and is also equivalent 
to [5],~(7.12)). 
\par 
The result stated immediately below (6.5.15) is a consequence of (6.5.15) and 
the functional equation $f_{\omega'}(\nu,p;g)=f_{\omega'}(-\nu,-p;g)$, which is (1.7.17). 
\par 
For the proof of (6.5.14), we use (in place of the above estimates for $K_{\mu}(2\pi R)$) 
the integral representation 
$$K_{\mu}(2\pi R)
={1\over 2}\int_0^{\infty} e^{-\pi R( t+t^{-1})} t^{-\mu -1}{\rm d}t\qquad\quad 
\hbox{($R>0$),}$$ 
which is derived from Equation~10.32.10 of [38]. This integral representation 
of $K_{\mu}(2\pi R)$ implies that  
$$e^{2\pi R}\left| K_{\mu}(2\pi R)\right| 
\leq\int_1^{\infty} e^{-\pi R t( 1-t^{-1})^2} t^{M-1}{\rm d}t\qquad\quad 
\hbox{($R>0$ and $M\geq |{\rm Re}(\mu)|+1$).}\eqno(6.5.16)$$
By setting $\tau =\cosh^{-1}(1+1/2\pi R)\in(0,\infty)$, and then 
applying the inequalities 
$$\left( 1-t^{-1}\right)^2\geq\cases{1/(\pi R e^{\tau}) &if $t\geq e^{\tau}$, \cr 
0 &if $1\leq t<e^{\tau}$,}$$ 
we deduce from (6.5.16) that 
$$\left| K_{\mu}(2\pi R)\right| 
\leq 2^M\left( M^{-1}+\Gamma(M)\right) \exp\left( {M\over 2\pi R}-2\pi R\right)\qquad\quad 
\hbox{for $|{\rm Re}(\mu)|\leq M-1$ and $R>0$.}$$ 
By combining this last result with the case $\omega =1$ of 
the equations~(5.26)-(5.27) of Lemma~5.1 of [5] we obtain, when 
$n_1\in N$, $k_1\in K$, $|{\rm Re}(\nu)|\leq\sigma_1$ and $R\geq r_0$, 
the upper bound: 
$$\pi^{-\nu} i^p\,\Gamma(\ell+1+\nu)
\bigl( {\bf J}_1\varphi_{\ell,q}(\nu,p)\bigr)\!\left( n_1 a[R] k_1\right)  
\ll_{\ell ,\sigma_1, r_0}\,\left( 1+|{\rm Im}(\nu)\right)^{\ell -|p|} 
R^{\ell +1} e^{-2\pi R}\;.\eqno(6.5.17)$$
It is a property of the Jacquet operator that 
$|u|^4 {\bf J}_{u^2\xi}={\bf h}_u {\bf J}_{\xi}{\bf h}_u$ for $u\in{\Bbb C}^{*}$, 
and so, by Equation (1.8.2) and the linearity of ${\bf J}_{\xi}$, 
one finds  that 
$$\bigl( {\bf J}_{\omega'}\varphi_{\ell,q}(\nu,p)\bigr)(g) 
=\left|\omega'\right|^{\nu -1} \left(\omega' /\left|\omega'\right|\right)^{-p} 
\bigl( {\bf J}_1\varphi_{\ell,q}(\nu,0)\bigr)
\bigl( h\bigl[\sqrt{\omega'}\,\bigr] g\bigr)\qquad\quad 
\hbox{for $\nu\in{\Bbb C}$, $g\in G$.}\eqno(6.5.18)$$
Since, moreover, ${\bf h}_u\rho =|u|^2\rho$ for $u\in{\Bbb C}^{*}$, 
the combination of (6.5.17) and (6.5.18) yields (6.5.14)\quad$\blacksquare$ 

\bigskip 

\noindent{\bf Remark~6.5.4.}\quad 
It follows by 
Stirling's asymptotic formula for $\log\Gamma(z)$ that one has   
$$|\Gamma(\mu +1)|\gg_{\sigma_2}
\,\left( 1+|{\rm Im}(\mu)|\right)^{{\rm Re}(\mu)+1/2} e^{-(\pi /2)|{\rm Im}(\mu)|} 
\qquad\qquad\hbox{($\mu\in{\Bbb C}$, $|{\rm Re}(\mu)|\leq\sigma_2<\infty$).}
\eqno(6.5.19)$$

\bigskip 

Before stating the next lemma, we find it convenient to first  
define one relevant piece of new terminology. For all 
$\lambda,\kappa\in{\Bbb Z}$ with $\lambda\geq |\kappa|$, and for all 
$\theta\in{\cal T}^{\lambda}_{\sigma}$, we put 
$$\widetilde{\bf L}^{\omega,\dagger}_{\lambda,\kappa}\theta 
=\widetilde{\bf L}^{\omega}_{\lambda,\kappa}\theta  
-b(\omega ;\lambda,\kappa;\theta ) {\bf M}_{\omega}\varphi_{\lambda,\kappa}(1,0)\;,
\eqno(6.5.20)$$
where the constant $b(\omega;\lambda,\kappa;\theta)$ is defined as in (6.4.6). 
It is then an immediate corollary of the relation (6.4.7) of Theorem~6.4.1 and 
the observation preceding (6.3.4) that, subject to the same hypotheses 
under which $\widetilde{\bf L}^{\omega,\dagger}_{\lambda,\kappa}\theta$ was just 
defined, one has: 
$$\widetilde{\bf L}^{\omega,\dagger}_{\lambda,\kappa}\theta 
\in C^{\infty}(N\backslash G,\omega)\quad\ {\rm and}\quad\ 
\widetilde{\bf L}^{\omega,\dagger}_{\lambda,\kappa}\theta\ 
\,\hbox{is of $K$-type $(\lambda,\kappa)$.}\eqno(6.5.21)$$ 
For later reference we note here also that, given 
our choice of $\tau\in C^{\infty}(G)$ (satisfying the conditions in (6.5.1)), 
it follows by the 
estimate (6.4.9) of Theorem~6.4.1,  
and the definitions (6.5.2) and (6.5.20), that if $\lambda,\kappa\in{\Bbb Z}$, 
$\lambda\geq |\kappa|$ and $\theta\in{\cal T}^{\lambda}_{\sigma}$ then  
$$\bigl(\widetilde{\bf L}^{\omega,*}_{\lambda,\kappa}\theta\bigr)(g) 
=\bigl(\widetilde{\bf L}^{\omega,\dagger}_{\lambda,\kappa}\theta\bigr)(g) 
\ll_{\omega,\theta} \left(\rho(g)\right)^{1+\sigma}\qquad\ 
\hbox{for all $g\in G$ such that $\rho(g)\leq 1$.}\eqno(6.5.22)$$

\bigskip

\proclaim Lemma 6.5.5. Put  
$$I(\ell,q)
=\left\{ (\lambda ,\kappa)\in{\Bbb Z}\times{\Bbb Z} \,:\, \lambda\geq |\kappa|\ {\rm and}\ 
\max\{ |\lambda -\ell| , |\kappa -q|\}\leq 1\right\}\;.\eqno(6.5.23)$$
Then, for each ${\bf X}\in{\frak g}$, there exists a family 
$\bigl( [{\bf X}]_{\ell,q}^{\lambda,\kappa}\bigr)_{(\lambda,\kappa)\in I(\ell,q)}$ 
of linear operators on the space ${\cal T}^{\ell}_{\sigma}$  
such that one has both 
$$[{\bf X}]_{\ell,q}^{\lambda,\kappa} : {\cal T}^{\ell}_{\sigma}\rightarrow 
{\cal T}^{\lambda}_{\sigma}\;,\qquad\quad\hbox{for each $\,(\lambda,\kappa)\in I(\ell,q)$,}\eqno(6.5.24)$$ 
and 
$$\sum_{(\lambda,\kappa)\in I(\ell,q)} 
\widetilde{\bf L}^{\omega,\dagger}_{\lambda,\kappa}[{\bf X}]_{\ell,q}^{\lambda,\kappa}\theta 
={\bf X}\widetilde{\bf L}^{\omega,\dagger}_{\ell,q}\theta\;,\qquad\quad   
\hbox{for all $\,\theta\in{\cal T}^{\ell}_{\sigma}$.}\eqno(6.5.25)$$ 

\medskip 

\noindent{\bf Proof.}\quad 
As in (3.11) of [5], we put  
$${\bf H}_1=\pmatrix{1/2&0\cr 0&-1/2}\;,\qquad 
{\bf V}_1 =\pmatrix{0&1/2\cr 1/2&0}\;,\qquad\ \,
{\bf V}_2 =\pmatrix{0&i/2\cr -i/2&0}\eqno(6.5.26)$$
(having already defined ${\bf H}_2$, ${\bf W}_1$ and ${\bf W}_2$ in (1.2.9)). 
The set 
${\cal B}=\{ {\bf H}_1,{\bf H}_2,{\bf V}_1,{\bf V}_2,{\bf W}_1,{\bf W}_2\}$ 
is a basis 
of the real Lie algebra $\frak{sl}(2,{\Bbb C})$ of $G$, and so is also a 
${\Bbb C}$-basis of the complex Lie algebra ${\frak g}$. 
Another ${\Bbb C}$-basis for ${\frak g}$ is the set 
${\cal B}_1=\{ {\bf H}_1,{\bf H}_2,{\bf F}^{+},{\bf F}^{-},{\bf E}^{+},{\bf E}^{-}\}$, where 
$${\bf F}^{\pm}={\bf V}_1\pm i{\bf V}_2\qquad\ {\rm and}\qquad\ 
{\bf E}^{\pm}={\bf W}_1\pm i{\bf W_2}\eqno(6.5.27)$$  
(with the factor `$i$' here signifying complexification). Therefore 
we may confine ourselves, in this proof, to a discussion of 
the cases in which one has ${\bf X}\in{\cal B}_1$: for, by linearity, these 
special cases of the lemma imply the general case.
\par 
Let ${\bf X}\in{\cal B}_1$. Then, for $(\nu,p)\in{\Bbb C}\times{\Bbb Z}$ with $|p|\leq\ell$, 
one has 
$${\bf X}{\bf J}_{\omega}\varphi_{\ell,q}(\nu,p) 
={\bf J}_{\omega}{\bf X}\varphi_{\ell,q}(\nu,p)\;,\eqno(6.5.28)$$
where, by virtue of the ${\frak g}$-invariance of the space 
$H(\nu,p)$ defined in Equation (1.6.1), 
there exist certain 
complex constants $c^{\bf X}_{\ell,q}(\lambda,\kappa;\nu,p)\,$ 
($\lambda,\kappa\in{\Bbb Z}$), all but finitely many of which are equal to zero, 
such that 
$${\bf X}\varphi_{\ell,q}(\nu,p) 
=\sum_{\lambda =|p|}^{\infty}\sum_{\kappa=-\lambda}^{\lambda} 
c^{\bf X}_{\ell,q}(\lambda,\kappa;\nu,p)\varphi_{\lambda,\kappa}(\nu,p)\;.\eqno(6.5.29)$$ 
In his preliminary notes preparatory to work reported on in~[5] 
Bruggeman has computed, for each ${\bf X}\in{\cal B}_1$, 
the coefficients in (6.5.29). The results of his computations 
(kindly made available to us by personal communication) 
are crucial for this proof; they show, in particular, that 
$c^{\bf X}_{\ell,q}(\lambda,\kappa;\nu,p)$ is a polynomial function 
of the complex variable $\nu$, and is identically zero when 
$(\lambda,\kappa)\not\in I(\ell,q)$. 
\par 
We confine ourselves, in what follows, 
to a discussion of the single case in which ${\bf X}={\bf F}^{+}$. 
This is justifiable, since one can deal similarly with the cases 
in which ${\bf X}\in\{ {\bf H}_1,{\bf H}_2,{\bf F}^{-},{\bf E}^{+},{\bf E}^{-}\}$. 
Bruggeman found that 
$$c^{{\bf F}^{+}}_{\ell,q}(\lambda,\kappa;\nu,p)
=0\qquad\quad\hbox{unless $(\lambda,\kappa)\in I(\ell,q)$ and $\kappa =q+1$,}\eqno(6.5.30)$$ 
and that, when $(\lambda,q+1)\in I(\ell,q)$ and $(\nu,p)\in{\Bbb C}\times{\Bbb Z}$ is 
such that $|p|\leq\ell$, one has: 
$$c^{{\bf F}^{+}}_{\ell,q}(\lambda,q+1;\nu,p)
=\cases{(\ell +1)^{-1}(2\ell+1)^{-1}(\nu +\ell+1)\left( (\ell+1)^2 -p^2\right) 
&if $\lambda =\ell+1$;\cr 
-\ell^{-1}(\ell+1)^{-1}(\ell-q)\nu p &if $\lambda =\ell$;\cr 
\ell^{-1}(2\ell+1)^{-1}(\ell-q)(\ell-q-1)(\ell-\nu) &if $\lambda=\ell-1\geq |p|$;\cr 
0 &otherwise.}\eqno(6.5.31)$$
Note that in (6.5.31) the case $\lambda =\ell$ will arise only if $\ell\geq |q+1|$, 
and that this inequality implies $\ell\neq 0$ (since it is assumed that we have $\ell\geq |q|$). 
\par 
Suppose now that $\theta\in{\cal T}^{\ell}_{\sigma}$. 
Motivated by the definition of the transform $\widetilde{\bf L}_{\ell,q}^{\omega}$ 
in Equation (6.4.4) we observe now that it follows by (6.5.28) and (6.5.29), for ${\bf X}={\bf F}^{+}$, 
and by (6.5.30), (6.5.31), (6.5.23),  
the substitution $\nu =it$ and 
the linearity of the Jacquet operator ${\bf J}_{\omega}$, that if $g\in G$ and $p$ is an integer 
satisfying $|p|\leq\ell$ then one has: 
$$\eqalign{
\int\limits_{(0)} &\theta(\nu,p)(\pi|\omega|)^{-\nu}\Gamma(\ell+1+\nu) 
\left( {\bf F}^{+} {\bf J}_{\omega}\varphi_{\ell,q}(\nu,p)\right)\!(g)\,\nu^{\epsilon(p)} 
\sin(\pi\nu) {\rm d}\nu =\cr 
 &\quad=-\int_{-\infty}^{\infty}\theta(it,p)(\pi|\omega|)^{-it}\,\Gamma(\ell+1+it)\ \times\cr 
 &\quad\qquad\qquad\ \times\biggl(\,\sum_{\lambda =\max\{\ell -1,|p|,|q+1|\}}^{\ell+1} 
c^{{\bf F}^{+}}_{\ell,q}(\lambda,q+1;it,p)
\bigl( {\bf J}_{\omega}\varphi_{\lambda,q+1}(it,p)\bigr) (g)\biggr) 
(it)^{\epsilon(p)} 
\sinh(\pi t) {\rm d}t\;.  
}\eqno(6.5.32)$$
Given the condition (T2) below (6.4.3), the expansion of 
$({\bf J}_{\omega}\varphi_{\lambda,\kappa}(\nu,p))(g)$ obtained in 
the equations~(5.26)-(5.27) of Lemma~5.1 of [5], the definition~(6.4.5) and 
the equation~(6.5.31), one may check that 
(once any inessential discontinuities at $t=0$ or $t=i$ are removed) 
the last integrand above 
is a continuous function 
$(t,g)\mapsto f(t,g)$ 
from the set $\{ t\in{\Bbb C} : |{\rm Im}(t)|\leq\sigma\}\times G$ into ${\Bbb C}$. This 
integrand $f(t,g)$ is, in particular, continuous on ${\Bbb R}\times G$; 
by  the case $\omega'=\omega$, 
$\,\sigma_1=1/2$, $\,\varepsilon =1/4$, $d\rightarrow\infty$  of the bound (6.5.15) of Lemma~6.5.3, 
the equation (6.5.31) and the condition (T3) below (6.4.3), it moreover satisfies  
$$f(t,g)\ll_{\theta,\ell,|\omega|,R,A}\ (1+|t|)^{\ell +(1/2)+\epsilon(p)-A}\qquad  
\hbox{for $\,A,R\in[1,\infty)$, $t\in{\Bbb R}$ and $g\in G$ such that $\rho(g)\leq R$.}$$
Since we have here $\epsilon(p)\leq 1$ (by (6.4.5)), it may be 
deduced from the special case $A=\ell +3$ 
of the above bound for the integrand $f(t,g)$ that integral on the right-hand side 
of Equation~(6.5.32) converges uniformly for all $g$ lying in any given 
compact subset of $G$. 
\par 
Similarly to the above, it may be shown that that the corresponding 
`undifferentiated' integral 
$$-\int_{-\infty}^{\infty}\theta(it,p)(\pi|\omega|)^{-it}\,\Gamma(\ell+1+it) 
\left( {\bf J}_{\omega}\varphi_{\ell,q}(it,p)\right)\!(g) (it)^{\epsilon(p)} 
\sinh(\pi t) {\rm d}t$$
is absolutely convergent (for $|p|\leq\ell$), and that its integrand $F(t,g)$ (say) 
is continuous on ${\Bbb R}\times G$. 
\par 
Given the observations of the preceding two paragraphs, and bearing in 
mind the definition of the differential operator 
${\bf F}^{+}$ (via the cases ${\bf X}={\bf V}_1$ and ${\bf X}={\bf V}_2$ of 
(1.2.6), and the first part of (6.5.27)), 
it follows by the second proposition of 
Section~1.88 of [43] that one may, by `differentiating inside the integral', 
deduce from the equation (6.4.4) (with $\theta$ substituted for $\eta$) 
that one has  
$$\eqalignno{
\bigl({\bf F}^{+}\widetilde{\bf L}_{\ell,q}^{\omega}\theta\bigr) (g) 
 &={1\over 2\pi^3 i}\sum_{|p|\leq\ell} 
{(-i\omega /|\omega|)^p\over\left\|\Phi_{p,q}^{\ell}\right\|_K} 
\int\limits_{(0)} \theta(\nu,p) (\pi|\omega|)^{-\nu}\,\Gamma(\ell +1+\nu) 
\left( {\bf F}^{+} {\bf J}_{\omega}\varphi_{\ell,q}(\nu,p)\right)\!(g) 
\,\nu^{\epsilon(p)} \sin(\pi\nu) 
{\rm d}\nu =\cr 
 &={1\over 2\pi^3 i}\sum_{\lambda =\max\{\ell -1,|q+1|\}}^{\ell+1} 
\sum_{|p|\leq\min\{\ell,\lambda\}} 
{(-i\omega /|\omega|)^p\over\left\|\Phi_{p,q}^{\ell}\right\|_K}\ \times \cr 
 &\qquad\times\int\limits_{(0)}^{\matrix{\ }} \theta(\nu,p) (\pi|\omega|)^{-\nu}
\,\Gamma(\ell +1+\nu) 
c^{{\bf F}^{+}}_{\ell,q}(\lambda,q+1;\nu,p)
\left( {\bf J}_{\omega}\varphi_{\lambda,q+1}(\nu,p)\right)\!(g)
\,\nu^{\epsilon(p)} \sin(\pi\nu) 
{\rm d}\nu =\cr 
 &=\sum_{(\lambda,\kappa)\in I(\ell,q)}^{\matrix{\ }} 
\Bigl(\widetilde{\bf L}^{\omega}_{\lambda,\kappa} 
\bigl[ {\bf F}^{+}\bigr]^{\lambda,\kappa}_{\ell,q}\theta\Bigr)(g)\;, &(6.5.33)}$$ 
where, for $(\lambda,\kappa)\in I(\ell,q)$, the function $\bigl[ {\bf F}^{+}\bigr]^{\lambda,\kappa}_{\ell,q}\theta : 
\{\nu\in{\Bbb C} : |{\rm Re}(\nu)|\leq\sigma\}\times\{ p\in{\Bbb Z} : |p|\leq\lambda\} 
\rightarrow{\Bbb C}$ satisfies 
$$\Bigl( \bigl[ {\bf F}^{+}\bigr]^{\lambda,\kappa}_{\ell,q}\theta\Bigr)(\nu,p) 
=\cases{\displaystyle{\|\Phi^{\lambda}_{p,\kappa}\|_{K}\over \|\Phi^{\ell}_{p,q}\|_{K}}\, 
{\Gamma(\ell +1+\nu)\over\Gamma(\lambda +1+\nu)}
\,c^{{\bf F}^{+}}_{\ell,q}(\lambda,\kappa;\nu,p)\theta(\nu,p) 
&if $\kappa=q+1$ and $|p|\leq\ell$,\cr 
\ & \ \cr
0 &otherwise,}\eqno(6.5.34)$$
while the summand in (6.5.33) is given by the 
case $\eta=\bigl[ {\bf F}^{+}\bigr]^{\lambda,\kappa}_{\ell,q}\theta$, 
$(\ell,q)=(\lambda,\kappa)$ of the equation (6.4.4). 
\par 
With regard to the definition (6.5.34) of the function $\widetilde{\bf L}^{\omega}_{\lambda,\kappa} 
\bigl[ {\bf F}^{+}\bigr]^{\lambda,\kappa}_{\ell,q}\theta$, 
we observe that if $(\lambda,\kappa)\in I(\ell,q)$ then 
the function 
$p\mapsto \|\Phi^{\lambda}_{p,\kappa}\|_{K}/\|\Phi^{\ell}_{p,q}\|_{K}$ is 
(by~(1.6.6)) an even function from $\{-\lambda,1-\lambda,\ldots ,\lambda\}$ 
into $(0,\infty)$. 
Moreover, when $(\lambda,q+1)\in I(\ell,q)$ and $(\nu,p)\in{\Bbb C}\times{\Bbb Z}$ is 
such that $|p|\leq\min\{\lambda,\ell\}$, it follows by 
Equation (6.5.31) that  
$${\Gamma(\ell +1+\nu)\over\Gamma(\lambda +1+\nu)}
\,c^{{\bf F}^{+}}_{\ell,q}(\lambda,q+1;\nu,p)
=\cases{(\ell +1)^{-1}(2\ell+1)^{-1}\left( (\ell+1)^2 -p^2\right) 
&if $\lambda =\ell+1$;\cr 
-\ell^{-1}(\ell+1)^{-1}(\ell-q)\nu p &if $\lambda =\ell$;\cr 
\ell^{-1}(2\ell+1)^{-1}(\ell-q)(\ell-q-1)\left(\ell^2-\nu^2\right) &if $\lambda=\ell-1$.}$$
Hence (and since $|p|$ is an even function of $p$) one may check that,  
when $(\lambda,\kappa)\in I(\ell,q)$, the conditions (T1)-(T3) below (6.4.3) 
continue to hold if one substitutes there  
$\bigl(\bigl[ {\bf F}^{+}\bigr]^{\lambda,\kappa}_{\ell,q}\theta\bigr)(\nu,p)$ for 
$\eta(\nu,p)$. Therefore we obtain the case ${\bf X}={\bf F}^{+}$ of the 
result (6.5.24) stated in the lemma. 
\par 
Given the result obtained in (6.5.33), and given the definition of 
$\widetilde{\bf L}^{\omega,\dagger}_{\lambda,\kappa}\theta$ in (6.5.20), 
the case ${\bf X}={\bf F}^{+}$ of the result (6.5.25) will follow 
if we are able to show that 
$$\sum_{(\lambda,\kappa)\in I(\ell,q)} b\!\left(\omega;\lambda,\kappa; 
\bigl[{\bf F}^{+}\bigr]^{\lambda,\kappa}_{\ell,q}\theta\right)  
{\bf M}_{\omega}\varphi_{\lambda,\kappa}(1,0) 
=b(\omega;\ell,q;\theta) {\bf F}^{+} {\bf M}_{\omega}\varphi_{\ell,q}(1,0)\;. 
\eqno(6.5.35)$$ 
Since the Goodman-Wallach operator ${\bf M}_{\omega}$ 
commutes with all elements of the Lie algebra ${\frak g}$, and since we 
have the case ${\bf X}={\bf F}^{+}$ of the equation (6.5.29), in which 
the coefficients $c^{{\bf F}^{+}}_{\ell,q}(\lambda,\kappa;\nu,p)$ 
($\lambda\geq |p|$, $|\kappa|\leq\lambda$) are 
given by (6.5.30)-(6.5.31), it follows by virtue of the linearity of the operator 
${\bf M}_{\omega}$ (in combination with the observation preceding (6.3.4), 
and the definition of `$K$-type' in the equations (1.3.1)) that the identity 
(6.5.35) is valid if and only if 
$$b\!\left(\omega;\lambda,\kappa; 
\bigl[{\bf F}^{+}\bigr]^{\lambda,\kappa}_{\ell,q}\theta\right)  
=b(\omega;\ell,q;\theta) c^{{\bf F}^{+}}_{\ell,q}(\lambda,\kappa;1,0)\eqno(6.5.36)$$
for all $(\lambda ,\kappa)\in I(\ell,q)$. 
\par 
Let $(\lambda ,\kappa)\in I(\ell,q)$. 
By the definition (6.4.6) of $b(\omega;\ell,q;\theta)$, and by (6.5.34) and (6.5.30), 
one finds that both sides of the equality sign in (6.5.36) equal zero if either 
$\ell =0$ or $\kappa\neq q+1$. If, on the other hand, we have $\ell\geq 1$ and 
$\kappa =q+1$ (so that $(\lambda ,q+1)=(\lambda,\kappa)\in I(\ell,q)$) then it follows from 
(6.4.6) and (6.5.34) that the equation (6.5.36) holds if and only if it is the case that  
$$\lambda c_{\ell,q}^{{\bf F}^{+}}(\lambda,q+1;0,1)
=\ell c_{\ell,q}^{{\bf F}^{+}}(\lambda,q+1;1,0)\eqno(6.5.37)$$
(note that the factor $\lambda$ 
on the left-hand side of this equation makes it unnecessary to distinguish 
the case $\lambda =0$). Since the equation (6.5.31) allows one to show 
that if $\ell\geq 1$ then the condition (6.5.37) is satisfied when 
$(\lambda,q+1)\in I(q,\ell)\,$ (with only 
the case $\lambda =0$ requiring any thought at all), we are therefore able 
to conclude that the condition (6.5.36) is satisfied 
for all $(\lambda ,\kappa)\in I(\ell,q)$; 
it follows that the identity (6.5.35) is valid, and so our proof of the case 
${\bf X}={\bf F}^{+}$ of the lemma is complete. Similar proofs exist for the cases 
in which 
${\bf X}\in\{ {\bf H}_1,{\bf H}_2,{\bf F}^{-},{\bf E}^{+},{\bf E}^{-}\}
={\cal B}_1-\{ {\bf F}^{+}\}$. Therefore, and since ${\cal B}_1$ is a ${\Bbb C}$-basis of 
${\frak g}$, the results stated in the lemma are valid 
in general (i.e. for all ${\bf X}\in{\frak g}$)\quad$\blacksquare$ 

\bigskip 

\noindent{\bf Remark~6.5.6.}\quad 
The above proof is not as we would ideally like it:  
for the validation of the identity (6.5.35) 
is reliant on detailed (and quite mindless) calculations. Despite giving 
it some thought, we were unable  
to discover a general principle that might `explain' the identity (6.5.35). 

\bigskip 

\proclaim Lemma 6.5.7. Let $\ell_1,q_1\in{\Bbb Z}$ be such that $\ell_1\geq |q_1|$,  
let $I(q_1,\ell_1)$ be defined as in Equation~(6.5.23) of Lemma~6.5.5, let 
$\theta\in{\cal T}^{\ell_1}_{\sigma}$, and let ${\bf X}\in{\frak g}$. 
Then ${\bf X}\widetilde{\bf L}^{\omega,\dagger}_{\ell_1,q_1}\theta\in 
C^{\infty}(N\backslash G,\omega)$. Moreover, the Poincar\'e series 
$(P^{\frak a}\widetilde{\bf L}^{\omega,\dagger}_{\ell_1,q_1}\theta)(g)$ 
and $(P^{\frak a}{\bf X}\widetilde{\bf L}^{\omega,\dagger}_{\ell_1,q_1}\theta)(g)$ 
are well-defined and absolutely convergent for all $g\in G$, and one has:  
$${\bf X}P^{\frak a}\widetilde{\bf L}^{\omega,\dagger}_{\ell_1,q_1}\theta 
=P^{\frak a}{\bf X}\widetilde{\bf L}^{\omega,\dagger}_{\ell_1,q_1}\theta\eqno(6.5.38)$$ 
and, for some $\theta^{\bf X}\in\prod_{(\lambda,\kappa)\in I(\ell_1,q_1)} 
{\cal T}^{\lambda}_{\sigma}$, 
$${\bf X}\widetilde{\bf L}^{\omega,\dagger}_{\ell_1,q_1}\theta 
=\sum_{(\lambda,\kappa)\in I(\ell_1,q_1)}^{\displaystyle ^{\quad}} \widetilde{\bf L}^{\omega,\dagger}_{\lambda,\kappa}
\theta^{\bf X}_{(\lambda,\kappa)}\qquad\ {\rm and}\qquad\ 
P^{\frak a}{\bf X}\widetilde{\bf L}^{\omega,\dagger}_{\ell_1,q_1}\theta 
=\sum_{(\lambda,\kappa)\in I(\ell_1,q_1)} 
P^{\frak a}\widetilde{\bf L}^{\omega,\dagger}_{\lambda,\kappa}\theta^{\bf X}_{(\lambda,\kappa)}\;.
\eqno(6.5.39)$$

\medskip 

\noindent{\bf Proof.}\quad 
The result that ${\bf X}\widetilde{\bf L}^{\omega,\dagger}_{\ell_1,q_1}\theta\in 
C^{\infty}(N\backslash G,\omega)$ is an immediate corollary of the first part of 
(6.5.21) and the left-invariance of the differential operator ${\bf X}$. 
By Lemma~6.5.5, we moreover have the identity 
$${\bf X}\widetilde{\bf L}^{\omega,\dagger}_{\ell_1,q_1}\theta 
=\sum_{(\lambda,\kappa)\in I(\ell_1,q_1)} \widetilde{\bf L}^{\omega,\dagger}_{\lambda,\kappa}
\theta^{\bf X}_{(\lambda,\kappa)}\;,\eqno(6.5.40)$$
where $\theta^{\bf X}_{(\lambda,\kappa)} 
=[X]^{\lambda,\kappa}_{\ell_1,q_1}\theta$ (with the operator $[X]^{\lambda,\kappa}_{\ell_1,q_1}$ 
as described in that lemma) so that, for $(\lambda,\kappa)\in I(\ell_1,q_1)$, 
$$\theta^{\bf X}_{(\lambda,\kappa)} 
=[X]^{\lambda,\kappa}_{\ell_1,q_1}\theta\in{\cal T}^{\lambda}_{\sigma}\;,\eqno(6.5.41)$$
and, by (6.5.21) (again), 
$$\widetilde{\bf L}^{\omega,\dagger}_{\lambda,\kappa}
\theta^{\bf X}_{(\lambda,\kappa)}\in C^{\infty}(N\backslash G,\omega)\;.\eqno(6.5.42)$$
By (6.5.41) and (6.5.22), we have the upper bounds 
$$\bigl( \widetilde{\bf L}^{\omega,\dagger}_{\lambda,\kappa}
\theta^{\bf X}_{(\lambda,\kappa)}\bigr)(g)\ll_{\omega,\theta,{\bf X}} 
\,\left(\rho(g)\right)^{1+\sigma} 
\qquad\ \hbox{for $(\lambda,\kappa)\in I\left(\ell_1,q_1\right)$ 
and all $g\in G$ such that $\,\rho(g)\leq 1$.}\eqno(6.5.43)$$ 
Therefore, and by virtue of the identity (6.5.40), it follows that one has also: 
$$\bigl( {\bf X}\widetilde{\bf L}^{\omega,\dagger}_{\ell_1,q_1}
\theta\bigr)(g)\ll_{\omega,\theta,{\bf X}} 
\,\left(\rho(g)\right)^{1+\sigma}\qquad\ \hbox{for $\,g\in G$ such that 
$\,\rho(g)\leq 1$.}\eqno(6.5.44)$$
\par 
Given it was shown that ${\bf X}\widetilde{\bf L}^{\omega,\dagger}_{\ell_1,q_1}\theta\in 
C^{\infty}(N\backslash G,\omega)$, and given the results of (6.5.21) and (6.5.22), 
and the results (6.5.42), (6.5.43) and (6.5.44) just obtained, it follows by Lemma~6.2.2 that 
the Poincar\'e series $(P^{\frak a}\widetilde{\bf L}^{\omega,\dagger}_{\ell_1,q_1}\theta)(g)$, 
$(P^{\frak a} {\bf X}\widetilde{\bf L}^{\omega,\dagger}_{\ell_1,q_1}\theta)(g)$ and 
$(P^{\frak a} \widetilde{\bf L}^{\omega,\dagger}_{\lambda,\kappa}
\theta^{\bf X}_{(\lambda,\kappa)})(g)\,$ (for all $(\lambda,\kappa)\in I(\ell_1,q_1)$) 
are each well-defined and absolutely convergent for all $g\in G$. 
Indeed, by Lemma~6.2.2 and Corollary~6.2.3, the Poincar\'e series 
$(P^{\frak a} {\bf X}\widetilde{\bf L}^{\omega,\dagger}_{\ell_1,q_1}\theta)(g)$ is 
uniformly convergent in each compact subset of $G$, and the sum of this series 
is a function $P^{\frak a} {\bf X}\widetilde{\bf L}^{\omega,\dagger}_{\ell_1,q_1}\theta : 
G\rightarrow{\Bbb C}$ which is continuous on $G$. 
Consequently the proposition of Section~1.72 of [43] (relating to 
`term by term' differentiation of a series) implies that 
if ${\bf X}\in\frak{sl}(2,{\Bbb C})$ and $g\in G$ then 
$( {\bf X}P^{\frak a}\widetilde{\bf L}^{\omega,\dagger}_{\ell_1,q_1}\theta)(g)$ 
exists, and one has: 
$$\eqalign{ 
\bigl( P^{\frak a}{\bf X}\widetilde{\bf L}^{\omega,\dagger}_{\ell_1,q_1}\theta\bigr)(g) 
 &={1\over\left[\Gamma_{\frak a} : \Gamma_{\frak a}'\right]}
\sum_{\gamma\in\Gamma_{\frak a}'\backslash\Gamma} 
\bigl( {\bf X}\widetilde{\bf L}^{\omega,\dagger}_{\ell_1,q_1}\theta\bigr) 
\!\left( g_{\frak a}^{-1}\gamma g\right) = \cr 
 &={1\over\left[\Gamma_{\frak a} : \Gamma_{\frak a}'\right]}
\sum_{\gamma\in\Gamma_{\frak a}'\backslash\Gamma} 
{{\rm d}\over{\rm d}t}
\,\bigl(\widetilde{\bf L}^{\omega,\dagger}_{\ell_1,q_1}\theta\bigr) 
\!\left( g_{\frak a}^{-1}\gamma g\exp(t{\bf X})\right)\Bigr|_{t=0} = \cr 
 &={{\rm d}\over{\rm d}t}
\Biggl( {1\over\left[\Gamma_{\frak a} : \Gamma_{\frak a}'\right]}
\sum_{\gamma\in\Gamma_{\frak a}'\backslash\Gamma} 
\bigl(\widetilde{\bf L}^{\omega,\dagger}_{\ell_1,q_1}\theta\bigr) 
\!\left( g_{\frak a}^{-1}\gamma g\exp(t{\bf X})\right)\Biggr)\Bigr|_{t=0} = \cr 
 &={{\rm d}\over{\rm d}t}
\Bigl( \bigl( P^{\frak a}\widetilde{\bf L}^{\omega,\dagger}_{\ell_1,q_1}\theta\bigr) 
\!\left( g\exp(t{\bf X})\right)\Bigr)\Bigr|_{t=0} 
=\bigl( {\bf X}P^{\frak a}\widetilde{\bf L}^{\omega,\dagger}_{\ell_1,q_1}\theta\bigr)(g)\;.}$$ 
This proves that (6.5.38) holds if ${\bf X}\in \frak{sl}(2,{\Bbb C})$. 
In the remaining cases, where one has ${\bf X}\in{\frak g}-\frak{sl}(2,{\Bbb C})$, 
the differential operator ${\bf X}$ may be expressed as a linear combination 
(with complex coefficients) of the six elements of the set ${\cal B}\subset \frak{sl}(2,{\Bbb C})$ 
that is defined just below (6.5.26). 
Therefore, given the linearity inherent in the definition (1.5.4) of Poincar\'e series, 
one may deduce these remaining cases 
of the result in (6.5.38) from those cases of the result that have already been proved. 
Similarly, given the simple fact (already noted) of the convergence of the relevant 
Poincar\'e series, one may deduce 
that (6.5.40) implies both of the identities stated in (6.5.39) (with, moreover, 
$\theta^{\bf X}=(\theta^{\bf X}_{(\lambda,\kappa)})_{(\lambda,\kappa)\in I(\ell_1,q_1)}\in 
\prod_{(\lambda,\kappa)\in I(\ell_1,q_1)} {\cal T}^{\lambda}_{\sigma}$, 
by virtue of (6.5.41))\quad$\blacksquare$ 

\bigskip 
 
\proclaim Lemma 6.5.8. Let $\widetilde{\bf L}^{\omega,*}_{\ell,q}\eta$ be as defined 
in (6.5.2). Then one has 
$P^{\frak a}\widetilde{\bf L}^{\omega,*}_{\ell,q}\eta\in 
C^{\infty}(G)\cap L^{\infty}(\Gamma\backslash G;\ell,q)$ and 
$P^{\frak a}|\widetilde{\bf L}^{\omega,*}_{\ell,q}\eta|\in 
C^0(G)\cap L^{\infty}(\Gamma\backslash G)$. 

\medskip 

\noindent{\bf Proof.}\quad 
By the relevant definitions in (6.5.2) and (6.5.20), we have 
$$\widetilde{\bf L}^{\omega,*}_{\ell,q}\eta 
=\widetilde{\bf L}^{\omega,\dagger}_{\ell,q}\eta 
+b(\omega;\ell,q;\eta)(1-\tau){\bf M}_{\omega}\varphi_{\ell,q}(1,0)\;.\eqno(6.5.45)$$ 
One may moreover observe, given the relations in (6.5.4), (6.5.21) and (6.5.22), 
that it follows by Lemma~6.2.2 and Corollary~6.2.3 that 
the Poincar\'e series 
$(P^{\frak a}\widetilde{\bf L}^{\omega,*}_{\ell,q}\eta)(g)$, 
$(P^{\frak a}|\widetilde{\bf L}^{\omega,*}_{\ell,q}\eta|)(g)$ and 
$(P^{\frak a}\widetilde{\bf L}^{\omega,\dagger}_{\ell,q}\eta)(g)$ 
are each convergent for all $g\in G$, and have sums that are 
continuous $\Gamma$-automorphic functions on $G$: 
$$P^{\frak a}\widetilde{\bf L}^{\omega,*}_{\ell,q}\eta, 
P^{\frak a}|\widetilde{\bf L}^{\omega,*}_{\ell,q}\eta|,  
P^{\frak a}\widetilde{\bf L}^{\omega,\dagger}_{\ell,q}\eta\in C^0(\Gamma\backslash G)\;.
\eqno(6.5.46)$$ 
By (6.5.45), (1.5.4) and the convergence just noted, we may deduce that 
$$P^{\frak a}\widetilde{\bf L}^{\omega,*}_{\ell,q}\eta 
=P^{\frak a}\widetilde{\bf L}^{\omega,\dagger}_{\ell,q}\eta 
+b(\omega;\ell,q;\eta)P^{\frak a}(1-\tau){\bf M}_{\omega}\varphi_{\ell,q}(1,0)\;.\eqno(6.5.47)$$ 
\par 
We show next that 
$$P^{\frak a}\widetilde{\bf L}^{\omega,\dagger}_{\ell,q}\eta\in C^{\infty}(G)\qquad\ 
{\rm and}\qquad\ 
P^{\frak a}\widetilde{\bf L}^{\omega,\dagger}_{\ell,q}\eta\!\quad 
\hbox{is of $K$-type $(\ell,q)$.}\eqno(6.5.48)$$ 
The proof is by induction; before coming to it  
we firstly introduce some relevant terminology. 
For each non-negative integer $j$, 
we define ${\cal P}(j)$ to be the proposition that 
$$\hbox{if $\,\ell_1,q_1\in{\Bbb Z}$, $\,\ell_1\geq |q_1|\,$ 
and $\,\theta\in{\cal T}^{\ell_1}_{\sigma}$}\  
\,{\rm then}\    
\,P^{\frak a}\widetilde{\bf L}^{\omega,\dagger}_{\ell_1,q_1}\theta\in C^j(G)\;,\eqno(6.5.49)$$
where $C^j(G)$ denotes the space of all functions $f : G\rightarrow{\Bbb C}$ 
that satisfy,  
for each $({\bf X}_1,\ldots ,{\bf X}_j)\in{\frak g}^j$, 
the condition that the $j$-th order derivative 
${\bf X}_j{\bf X}_{j-1}\cdots {\bf X}_1 f$ be defined and continuous on $G$ 
(so that $C^0(G)$ is, as previously, the space 
of all continuous functions $f : G\rightarrow{\Bbb C}$). 
Note that it is trivially implicit in the definition of 
of $C^j(G)$, just given, that $C^j(G)\supseteq C^{j+m}(G)$ for all 
non-negative integers $j$ and $m$; less trivially, one has 
$$\bigcap_{j=0}^{\infty} C^j(G) =C^{\infty}(G)\;.\eqno(6.5.50)$$ 
\par 
The equality (6.5.50) may be proved by determining  
local coordinates $x_1(g),\ldots ,x_6(g)\in{\Bbb R}$ for $G$, 
such that, for $j=1,\ldots ,6$, one has an operator identity of the form 
$$\partial /\partial x_j = c_{1,j}(g){\bf H}_1+c_{2,j}(g){\bf H}_2+c_{3,j}(g){\bf V}_1 
+c_{4,j}(g){\bf V}_2+c_{5,j}(g){\bf W}_1+c_{6,j}(g){\bf W}_2\;,$$
in which 
${\bf H}_1,{\bf H}_2,{\bf V}_1,{\bf V}_2,{\bf W}_1,{\bf W}_2$ are the elements 
of the basis ${\cal B}$ of $\frak{sl}(2,{\Bbb C})$ 
defined below (6.5.26), and 
the coeffients $c_{1,j}\ldots ,c_{6,j}$ are (verifiably) 
functions in the space $C^{\infty}(G)$. 
We omit the (not very enlightening) 
details of this proof, and merely 
note that it requires, in addition to the Iwasawa coordinates, 
some alternative system of coordinates: this is due to singularities 
which occur in the relevant coefficients when the Iwasawa coordinate 
$\theta$ is an integer multiple of $\pi$. 
\par 
As a starting point for our proof by induction of (6.5.48), we observe that 
since our only assumption concerning $\eta$ is that $\eta\in{\cal T}^{\ell}_{\sigma}$, 
and since we assume nothing more of $\ell$ and $q$ than that 
$\ell,q\in{\Bbb Z}$ and $\ell\geq |q|$, the fact 
(recorded in (6.5.46)) of our having shown that 
$P^{\frak a}\widetilde{\bf L}^{\omega,\dagger}_{\ell,q}\eta\in C^0(\Gamma\backslash G)$ 
is enough for us to infer that 
$${\cal P}(0)\ \,\hbox{is true.}\eqno(6.5.51)$$
\par 
Suppose now that $J$ is a non-negative integer such that 
$${\cal P}(J)\ \,\hbox{is true.}\qquad\eqno(6.5.52)$$
Let $\ell_1,q_1\in{\Bbb Z}$ satisfy $\ell_1\geq |q_1|$,  
let $\theta\in{\cal T}^{\ell_1}_{\sigma}$, and let 
${\bf X}_1,\ldots ,{\bf X}_{J+1}\in{\frak g}$. Then, by the results 
(6.5.38)-(6.5.39) of Lemma~6.5.7, we have 
$${\bf X_1}P^{\frak a}\widetilde{\bf L}^{\omega,\dagger}_{\ell_1,q_1}\theta 
=\sum_{(\lambda,\kappa)\in I(\ell_1,q_1)} 
P^{\frak a}\widetilde{\bf L}^{\omega,\dagger}_{\lambda,\kappa}
\theta^{{\bf X}_1}_{(\lambda,\kappa)}\;,\eqno(6.5.53)$$
for some $\theta^{{\bf X}_1}\in\prod_{(\lambda,\kappa)\in I(\ell_1,q_1)} 
{\cal T}^{\lambda}_{\sigma}$. 
Since it therefore follows by our induction hypothesis (6.5.52) 
(i.e. by the case $j=J$ of what is stated in (6.5.49)) that every 
summand on the right-hand side of Equation~(6.5.53) is a function 
in the space $C^J(G)$, we may deduce that the first order derivative 
$f={\bf X}_1 P^{\frak a}\widetilde{\bf L}^{\omega,\dagger}_{\ell_1,q_1}\theta$
lies in the same space; it follows that the $(J+1)$-st order derivative 
${\bf X}_{J+1}{\bf X}_J\cdots{\bf X}_1 
P^{\frak a}\widetilde{\bf L}^{\omega,\dagger}_{\ell_1,q_1}\theta$ (which is 
a $J$-th order derivative of $f$) 
is both defined and continuous on $G$. 
Since our only assumption concerning ${\bf X}_1,{\bf X}_2,\ldots ,{\bf X}_{J+1}$ 
is that these are elements of ${\frak g}$, the conclusion just reached 
is proof that $P^{\frak a}\widetilde{\bf L}^{\omega,\dagger}_{\ell_1,q_1}\theta$ 
lies in the space $C^{J+1}(G)$. 
Therefore, given that our assumptions concerning $\ell_1,q_1$ and $\theta$ 
match the conditions stated in (6.5.49), it may (similarly) be inferred that the 
proposition stated in (6.5.49) is true for $j=J+1$. 
\par 
The above discussion (subsequent to (6.5.52)) has shown that for each 
non-negative integer $J$ the proposition ${\cal P}(J)$ implies 
its successor ${\cal P}(J+1)$. This, combined with the result (6.5.51), 
implies (by induction) that the proposition ${\cal P}(j)$ stated 
in (6.5.49) is true for all non-negative integers $j$,  
and so we may conclude that, subject to the conditions 
on $\ell_1,q_1$ and $\theta$ in (6.5.49) being satisfied, one has 
$P^{\frak a}\widetilde{\bf L}^{\omega,\dagger}_{\ell_1,q_1}\theta\in 
\bigcap_{j=0}^{\infty} C^j(G)$.  
Hence, given (6.5.50), we obtain the first part of 
what is asserted in (6.5.48). 
\par 
We begin our verification of the second part of (6.5.48) 
with the observation that, with ${\bf H}_2\in\frak{sl}(2,{\Bbb C})$ 
defined as in (1.2.9),  
it follows by the result (6.5.38) of Lemma~6.5.7, 
the second part of (6.5.21), and the definitions (1.3.1) and (1.5.4), that 
$${\bf H}_2 P^{\frak a}\widetilde{\bf L}^{\omega,\dagger}_{\ell,q}\eta 
=P^{\frak a}{\bf H}_2\widetilde{\bf L}^{\omega,\dagger}_{\ell,q}\eta 
=P^{\frak a}(-iq)\widetilde{\bf L}^{\omega,\dagger}_{\ell,q}\eta   
=-iq P^{\frak a}\widetilde{\bf L}^{\omega,\dagger}_{\ell,q}\eta\;.\eqno(6.5.54)$$
Moreover, given any choice of ${\bf X},{\bf Y}\in{\frak g}$, 
the results of Lemma~6.5.7 imply that 
$${\bf X}P^{\frak a}\widetilde{\bf L}^{\omega,\dagger}_{\ell,q}\eta 
=\sum_{(\lambda,\kappa)\in I(\ell,q)} 
P^{\frak a}\widetilde{\bf L}^{\omega,\dagger}_{\lambda,\kappa}
\eta^{\bf X}_{(\lambda,\kappa)}\;,$$
where $\eta^{\bf X}\in\prod_{(\lambda,\kappa)\in I(\ell,q)} 
{\cal T}^{\lambda}_{\sigma}$. Two further applications of 
Lemma~6.5.7 then show that one also has 
$$\sum_{(\lambda,\kappa)\in I(\ell,q)} 
{\bf Y} P^{\frak a}\widetilde{\bf L}^{\omega,\dagger}_{\lambda,\kappa}
\eta^{\bf X}_{(\lambda,\kappa)} 
=\sum_{(\lambda,\kappa)\in I(\ell,q)} 
P^{\frak a}{\bf Y}\widetilde{\bf L}^{\omega,\dagger}_{\lambda,\kappa}
\eta^{\bf X}_{(\lambda,\kappa)}   
= P^{\frak a}{\bf Y}\sum_{(\lambda,\kappa)\in I(\ell,q)} 
\widetilde{\bf L}^{\omega,\dagger}_{\lambda,\kappa}
\eta^{\bf X}_{(\lambda,\kappa)} 
= P^{\frak a}{\bf Y}{\bf X}\widetilde{\bf L}^{\omega,\dagger}_{\ell,q}\eta$$
(the penultimate equality here being a consequence of the linearity of 
the differential operator ${\bf Y}$, and the linearity inherent in 
the definition (1.5.4) of Poincar\'e series); this proves that 
${\bf Y}{\bf X}P^{\frak a}\widetilde{\bf L}^{\omega,\dagger}_{\ell,q}\eta 
=P^{\frak a}{\bf Y}{\bf X}\widetilde{\bf L}^{\omega,\dagger}_{\ell,q}\eta$ 
for all ${\bf X},{\bf Y}\in{\frak g}$, and so, given the first equality of (1.2.11), 
we obtain (similarly to (6.5.54)) the result that 
$$\Omega_{\frak k}P^{\frak a}\widetilde{\bf L}^{\omega,\dagger}_{\ell,q}\eta  
=P^{\frak a}\Omega_{\frak k}\widetilde{\bf L}^{\omega,\dagger}_{\ell,q}\eta 
=-{1\over 2}\left(\ell^2 +\ell\right) 
P^{\frak a}\widetilde{\bf L}^{\omega,\dagger}_{\ell,q}\eta \;.\eqno(6.5.55)$$
\par 
By (6.5.54) and (6.5.55), the Poincar\'e series 
$P^{\frak a}\widetilde{\bf L}^{\omega,\dagger}_{\ell,q}\eta$ is a function 
of $K$-type $(\ell,q)$. This completes the verification of 
what is asserted in (6.5.48). Hence, and by Lemma~6.5.2, 
each of the two Poincar\'e series 
$P^{\frak a}\widetilde{\bf L}^{\omega,\dagger}_{\ell,q}\eta$ and 
$P^{\frak a}(1-\tau){\bf M}_{\omega}\varphi_{\ell,q}(1,0)$ 
is a function of $K$-type $(\ell,q)$ lying in the space $C^{\infty}(G)$. 
Therefore we may deduce from the identity (6.5.47) that 
$$P^{\frak a}\widetilde{\bf L}^{\omega,*}_{\ell,q}\eta\in C^{\infty}(G)\qquad\ 
{\rm and}\qquad\ 
P^{\frak a}\widetilde{\bf L}^{\omega,*}_{\ell,q}\eta\quad 
\hbox{is of $K$-type $(\ell,q)$.}\eqno(6.5.56)$$
\par 
Given (6.5.56), the results recorded in (6.5.46), and the 
definitions of the two spaces $L^{\infty}(\Gamma\backslash G;\ell,q)$ and 
$L^{\infty}(\Gamma\backslash G)$ (in, and above, (6.5.6)), 
it suffices for completion of the proof of the lemma that 
we show that the Poincar\'e series 
$P^{\frak a}\widetilde{\bf L}^{\omega,*}_{\ell,q}\eta$ and 
$P^{\frak a}|\widetilde{\bf L}^{\omega,*}_{\ell,q}\eta|$ 
are bounded functions on $G$. To this end we note that, by 
(6.5.2), (6.5.1) and the result (6.4.8) of Theorem~6.4.1, one has: 
$$\bigl(\widetilde{\bf L}^{\omega,*}_{\ell,q}\eta\bigr)(g) 
=\bigl(\widetilde{\bf L}^{\omega}_{\ell,q}\eta\bigr)(g)  
\ll_{\omega,\eta,A} \left(\rho(g)\right)^{-A}\qquad\ 
\hbox{for all $g\in G$ such that $\rho(g)\geq 2$, and all $A\in [0,\infty)$.}$$
We also have (6.5.22), and so it certainly follows  
that both of the conditions (6.2.5) and (6.2.22) are satisfied 
when one has (there): $R_0=1$, $\sigma_0=\sigma$, $R_{\infty}=2$, 
$\sigma_{\infty}=10^{10}$ (say), and either 
$f_{\omega}=\widetilde{\bf L}^{\omega,*}_{\ell,q}\eta$, or else 
$\omega =0$ and $f_0=|\widetilde{\bf L}^{\omega,*}_{\ell,q}\eta|$. 
Since $\sigma\in(1,2)$ (by hypothesis), since $10^{10}\geq 1$, and 
since the relations stated in (6.5.4) hold, it therefore follows 
by Corollary~6.2.10 that each of the Poincar\'e series 
$P^{\frak a}\widetilde{\bf L}^{\omega,*}_{\ell,q}\eta$ and 
$P^{\frak a}|\widetilde{\bf L}^{\omega,*}_{\ell,q}\eta|$ 
is indeed a bounded function on $G\quad\blacksquare$ 

\bigskip 

Our next lemma 
generalises the case ${\frak O}={\Bbb Z}[i]$ of Theorem~7.6.9 of~[11] 
and is a minor addition  
to the theory of the Linnik-Selberg series 
$$Z^{{\frak a},{\frak b}}_{\omega,\omega'}(s) 
=\sum_{c\in{}^{\frak a}{\cal C}^{\frak b}} 
{S_{{\frak a},{\frak b}}\left(\omega,\omega';c\right)\over |c|^{4s}}\qquad\qquad 
\hbox{( $\omega'\in{\frak O}$, ${\rm Re}(s)>3/4$).}\eqno(6.5.57)$$ 
In Theorem~2.16 of [7] Cogdell, Li, Piatetski-Shapiro and Sarnak have shown 
that if  
${\frak a}$ and ${\frak b}$ are $\Gamma$-equivalent cusps then the Linnik-Selberg 
series $Z^{{\frak a},{\frak b}}_{\omega,\omega'}(s)$ can be 
meromorphically continued into all of ${\Bbb C}$. 
In the case of the analogous series involving generalised Kloosterman sums 
associated with any group $\Gamma'<SL(2,{\Bbb R})$ that is a congruence subgroup 
with respect to $SL(2,{\Bbb Z})$, the corresponding 
meromorphic continuation was obtained by Selberg, in Section~3 of [41].  

\bigskip 

\proclaim Lemma 6.5.9. Let $\omega,\omega'\in{\frak O}$, and let $\omega$ be non-zero. 
Suppose that $\sigma^{*}>3/4$. Then the series 
$$\sum_{c\in{}^{\frak a}{\cal C}^{\frak b}} 
{\left| S_{{\frak a},{\frak b}}\left(\omega,\omega';c\right)\right|\over |c|^{4s}}
\eqno(6.5.58)$$
is uniformly convergent in the half-plane where ${\rm Re}(s)\geq \sigma^{*}$. 
For $\sigma'\geq\sigma^{*}$, one has the upper bound 
$$\sum_{c\in{}^{\frak a}{\cal C}^{\frak b}} 
{\left| S_{{\frak a},{\frak b}}\left(\omega,\omega';c\right)\right|\over |c|^{4\sigma'}} 
\leq 2^{-1/2}|\omega|\left| m_{\frak a} m_{\frak b}\right|^{2-2\sigma^{*}} 
\biggl(\,\sum_{\alpha\mid q_0} {1\over |\alpha|}\,\biggr)^{\!\!2}  
\,\zeta_{{\Bbb Q}(i)}^2\!\left( 2\sigma^{*} - 1/2\right) 
\;,\eqno(6.5.59)$$ 
where $\zeta_{{\Bbb Q}(i)}(s)=(1/4)\sum_{0\neq\alpha\in{\frak O}} |\alpha|^{-2s}\,$  
and $m_{\frak c}$ and $q_0$ 
are as in (6.1.25) 

\medskip 

\noindent{\bf Proof.}\quad 
Since $|c|^{4s}$ has absolute value $|c|^{4{\rm Re}(s)}$ it will 
certainly suffice to prove that the series is uniformly convergent 
in the real interval $[\sigma^{*} ,\infty)$. 
Moreover, by (6.1.25) and the result (6.1.26) of Lemma~6.1.5, the series (6.5.58) 
is a Dirichlet series of the special type dealt with in Chapter~9 of [43]: 
in particular, the condition $c\in{}^{\frak a}{\cal C}^{\frak b}$ 
implies that $|c|^4$ is a positive integer. Hence 
(given that, for $n\in{\Bbb N}$, the real function $s\mapsto n^{-s}$ is 
monotonic decreasing), it suffices for proof of the lemma 
that we show that the series in (6.5.58) is convergent for all $s>3/4$. 
A quite standard calculation (which we omit) shows that the convergence of 
this series, for $s>3/4$, may be deduced from the results 
(6.1.26) and (6.1.27) of Lemma~6.1.5: this same calculation
also yields the upper bound in (6.5.59)\quad$\blacksquare$ 

\bigskip 

\proclaim Lemma 6.5.10. Let ${\rm Re}(\nu)>1$, let $0\neq\omega\in{\frak O}$ and 
$\omega'\in{\frak O}$,  
and let $f_{\omega}={\bf M}_{\omega}\varphi_{\ell,q}(\nu,0)$. 
Then, for each $g\in G$, the integral on the right-hand side of Equation~(6.2.13) 
exists, and the equation (6.2.13) defines a function 
$F^{\frak b}_{\omega'} P^{\frak a} f_{\omega}\in C^0(N\backslash G,\omega)$ satisfying, 
for $g\in G$, 
$$\eqalign{
\bigl(F^{\frak b}_{\omega'} P^{\frak a} f_{\omega}\bigr)(g) 
 &=\,\zeta^{{\frak a},{\frak b}}_{\omega,\omega'}\!\!\left({1+\nu\over 2}\right) 
\bigl({\bf J}_{\omega'}\varphi_{\ell,q}(\nu,0)\bigr)(g)\ + \cr 
 &\quad +{1^{\textstyle{\ \atop\ }}\over \left[\Gamma_{\frak a} : \Gamma_{\frak a}'\right]}
\sum_{\scriptstyle\gamma\in\Gamma_{\frak a}'\backslash\Gamma 
\ :\ \gamma{\frak b}={\frak a}\atop
{\scriptstyle g_{\frak a}^{-1}\gamma g_{\frak b}\in h\left[ u(\gamma)\right] N  
\atop\scriptstyle\hbox{\qquad}}} 
\delta_{\omega u(\gamma) , \omega'/u(\gamma)} 
\bigl({\bf M}_{\omega}\varphi_{\ell,q}(\nu,0)\bigr)
\left( g_{\frak a}^{-1}\gamma g_{\frak b} g\right) ,}\eqno(6.5.60)$$ 
where 
$$\zeta^{{\frak a},{\frak b}}_{\omega,\omega'}(s) 
={1\over\left[\Gamma_{\frak a} : \Gamma_{\frak a}'\right]}
\sum_{c\in\,{}^{\frak a}{\cal C}^{\frak b}}^{\hbox{\quad }}
{\cal J}^{*}_{2s-1,0}\!\left({2\pi\sqrt{\omega\omega'}\over c}\right) 
{S_{{\frak a} , {\frak b}}\!\left(\omega , \omega' ; c\right)\over |c|^{4s}} 
\eqno(6.5.61)$$ 
(with ${\cal J}^{*}_{\nu,p}(z)$  as given by Equation~(6.3.12) of Lemma~6.3.2). 
Regarding the case $\omega'=0$, one has moreover:  
$$\eqalignno{
F^{\frak b}_0 P^{\frak a} f_{\omega}  
 &={1\over\left[\Gamma_{\frak a} : \Gamma_{\frak a}'\right]} 
\,{\sin(\pi\nu)\over\nu^2}\,{\Gamma(\ell+1-\nu)\over\Gamma(\ell+1+\nu)} 
\,Z^{{\frak a},{\frak b}}_{\omega,0}\!\left({1+\nu\over 2}\right)
\varphi_{\ell,q}(-\nu,0) = &(6.5.62)\cr 
 &=\pi\,{\Gamma(\nu)\Gamma(\ell +1-\nu)\over\Gamma(1-\nu)\Gamma(\ell +1+\nu)} 
\,\zeta^{{\frak a},{\frak b}}_{\omega, 0}\!\!\left({1+\nu\over 2}\right) 
\varphi_{\ell,q}(-\nu,0) = &(6.5.63)\cr 
 &=\zeta^{{\frak a},{\frak b}}_{\omega, 0}\!\!\left({1+\nu\over 2}\right) 
{\bf J}_0\varphi_{\ell,q}(\nu,0) &(6.5.64)
}$$
(with $Z^{{\frak a},{\frak b}}_{\omega,0}(s)$ as in (6.5.57)). 
When $s=(1+\nu)/2$ the sum on the right-hand side of Equation~(6.5.61) 
is absolutely convergent. 

\medskip 

\noindent{\bf Proof.}\quad 
By Lemma~6.5.1 we have $P^{\frak a} f_{\omega}\in C^0(\Gamma\backslash G)$; 
the integral in Equation~(6.2.13) therefore exists. 
Moreover, by the observation preceding (6.3.4), we have 
$f_{\omega}\in C^{\infty}(N\backslash G,\omega)$, and by the estimate 
(6.3.7) of Lemma~6.3.1 the case $\sigma_0 ={\rm Re}(\nu)$, 
$R_0=1$ of the condition (6.2.5) is satisfied. 
Therefore it follows by Lemma~6.2.5 that the equations (1.5.4) and (6.2.13) 
define a function 
$F^{\frak b}_{\omega'} P^{\frak a} f_{\omega}\in C^0(N\backslash G,\omega)$, 
and that a valid formula for this function is given by 
the case ${\frak a}'={\frak b}$ of Equation~(1.5.5) 
(in which equation, moreover, all the relevant summations are absolutely convergent). 
\par 
Assuming now that $g\in G$, the formula (1.5.5) may be applied: 
we thereby obtain (by firstly recalling that 
$f_{\omega}={\bf M}_{\omega}\varphi_{\ell,q}(\nu,0)$, and then making 
appropriate use of Equation~(6.3.5) and the results (6.3.10)-(6.3.12) of Lemma~6.3.2)
both the result (6.5.62) and the case $\omega'\neq 0$ of the 
result stated in (6.5.60)-(6.5.61). 
The results (6.5.63) and (6.5.64) follow by virtue of 
the definitions (6.3.12), (1.9.6) (which imply  
${\cal J}^{*}_{\nu,0}(0)=(J^{*}_{\nu}(0))^2=(\Gamma(\nu +1))^{-2}\,$) 
and the evaluation of ${\bf J}_0\varphi_{\ell,q}(\nu,p)$ in (1.5.18). 
Given that $\omega\neq 0$, 
the case $\omega' =0$ of the result in (6.5.60)-(6.5.61) follows 
from (6.5.64). 
The absolute convergence, at $s=(1+\nu)/2$, of the sum in (6.5.61) 
is a corollary of the point noted, in parenthesis, at the end of the last 
paragraph\quad$\blacksquare$ 

\bigskip

In stating the remaining lemmas of this subsection we shall use 
(without comment) the same terminology as is used in the statement 
of Lemma~6.5.10. In particular, ${\cal J}^{*}_{\nu,p}(z)$ denotes 
the function defined by (6.3.12) and (1.9.6) (or (1.9.8)),   
and $\zeta^{{\frak a},{\frak b}}_{\omega,\omega'}(s)$ denotes the 
function defined by (6.5.61) and (1.5.8)-(1.5.10). 

\bigskip 

\proclaim Lemma 6.5.11. Let $z\in{\Bbb C}$. Then the complex function 
$\nu\mapsto {\cal J}^{*}_{\nu,0}(z)$ is entire, and one has 
$$\left|  {\cal J}^{*}_{\mu -1/2,0}(z)\right| 
\leq {\Gamma^2\!\left(\,{\rm Re}(\mu)\right) e^{2|{\rm Im}(z)|}\over 
\Gamma^2\!\left(\,{\rm Re}(\mu +1/2)\right) |\Gamma(\mu)|^2}\qquad\ 
\hbox{for ${\rm Re}(\mu)>0$.}\eqno(6.5.65)$$

\medskip 

\noindent{\bf Proof.}\quad 
For all $\xi\in{\Bbb C}$, the series in Equation~(1.9.6) is 
uniformly convergent in some (or indeed any) compact neighbourhood 
of $\xi$. It may therefore be deduced that the equation (1.9.6) 
defines a complex function $\nu\mapsto J^{*}_{\nu}(z)$ which is entire; 
given that one may substitute $\overline{z}\,$ for $z$ in the foregoing, 
it therefore follows that the complex function 
$\nu\mapsto J^{*}_{\nu}(z)J^{*}_{\nu}\left(\overline{z}\right)={\cal J}^{*}_{\nu,0}(z)$ 
is entire. 
\par 
Suppose now that ${\rm Re}(\mu)>0$. Then, by (1.9.8) and Poisson's formula, 
Equation~10.9.4 of [38], one has 
$$J^{*}_{\mu -1/2}(z) 
=(z/2)^{1/2-\mu} J_{\mu -1/2}(z) 
={2\over\pi^{1/2}\Gamma(\mu)}\int_0^1 \left( 1-t^2\right)^{\mu -1}\cos(zt) {\rm d}t\;.$$
Hence, by using the bound 
$\max_{t\in [0,1]}|\cos(zt)|\leq\exp(|{\rm Im}(z)|)$, the 
substitution $t^2 =u$, and Euler's evaluation of his Beta-function, $B(m,n)$, 
one finds that 
$|J^{*}_{\mu -1/2}(z)|\leq \exp(|{\rm Im}(z)|) |\Gamma(\mu)|^{-1} 
\Gamma({\rm Re}(\mu))/\Gamma(1/2+{\rm Re}(\mu))$. 
One may substitute $\overline{z}\,$ for $z$ in the last inequality, and so 
obtain (by (6.3.12)) the result stated in (6.5.65)
\ $\blacksquare$ 

\bigskip 

\proclaim Lemma 6.5.12. Let $0\neq\omega\in{\frak O}$ and $\omega'\in{\frak O}$. 
Then the complex function $s\mapsto\zeta^{{\frak a},{\frak b}}_{\omega,\omega'}(s)$ 
is holomorphic on the half-pane where ${\rm Re}(s)>3/4$, and if 
${\cal D}$ is a non-empty compact subset of that half-plane then one has 
$$\max_{s\in{\cal D}}\bigl|\,\zeta^{{\frak a},{\frak b}}_{\omega,\omega'}(s)\bigr| 
=O_{\Gamma,{\cal D}}\left( |\omega| e^{4\pi\sqrt{|\omega\omega'|}}\right) .\eqno(6.5.66)$$

\medskip 

\noindent{\bf Proof.}\quad 
Let ${\cal D}$ be any non-empty compact subset of the half-plane 
$H=\{ s\in{\Bbb C} : {\rm Re}(s)>3/4\}$.  
Then there is some $\sigma^{*}=\sigma^{*}({\cal D})\in (3/4,\infty)$ and some  
$r_1=r_1({\cal D})\in [\sigma^{*} ,\infty)$ such that one has: 
$$|s|\leq r_1<\infty\qquad{\rm and}\qquad{\rm Re}(s)\geq\sigma^{*} >3/4\;,\qquad 
\hbox{for all $s\in{\cal D}$.}\eqno(6.5.67)$$ 
\par 
Suppose moreover that $c\in{}^{\frak a}{\cal C}^{\frak b}$.  
Then, by (6.1.25) and the result (6.1.26) of Lemma~6.1.5, it follows that 
$|c|\geq\sqrt{|m_{\frak a} m_{\frak b}|}\geq 1$. Hence, and since 
(6.5.67) implies that $\min_{s\in{\cal D}} {\rm Re}(2s)>3/2>1/2$, 
the bound (6.5.65) of Lemma~6.5.11 implies that for all $s\in{\cal D}$ one has 
$$\eqalignno{
\left| {\cal J}^{*}_{2s-1,0}\!\left(\!{2\pi\sqrt{\omega\omega'}\over c}\right) 
\!{S_{{\frak a} , {\frak b}}\!\left(\omega , \omega' ; c\right)\over |c|^{4s}}\right| 
 &\leq {\Gamma^2\!\left(\,{\rm Re}(2s-1/2)\right)\over 
\Gamma^2\!\left(\,{\rm Re}(2s)\right) |\Gamma(2s-1/2)|^2}
\,\exp\!\left(\!2\left|{\rm Im}\!\left(\!{2\pi\sqrt{\omega\omega'}\over c}\right)\right|\right) 
\!{\left| 
S_{{\frak a} , {\frak b}}\!\left(\omega , \omega' ; c\right)\right| 
\over |c|^{4{\rm Re}(s)}} 
\ll_{r_1} \cr 
 &\ll_{r_1}\,\exp\!\left( 4\pi\sqrt{\left|\omega\omega'\right|}\right) 
{\left| 
S_{{\frak a} , {\frak b}}\!\left(\omega , \omega' ; c\right) 
\right|\over |c|^{4{\rm Re}(s)}} 
&(6.5.68)}$$ 
(the last bound following by virtue of the fact that 
both $\Gamma(w)$ and $1/\Gamma(w)$ are holomorphic on the compact set  
$\{ w\in{\Bbb C} : |w|\leq 2 r_1\ {\rm and}\ {\rm Re}(w)\geq 1\}$). 
The bound (6.5.68), in combination with the 
first part of Lemma~6.5.9, is enough to imply the uniform convergence, 
for all $s\in{\cal D}$, of the series on the right-hand side 
of Equation~(6.5.61). It therefore follows (given the first part of 
Lemma~6.5.11) that the equation (6.5.61) defines a complex function 
$s\mapsto \zeta^{{\frak a},{\frak b}}_{\omega,\omega'}(s)$ which is 
holomorphic on the given compact set ${\cal D}\subset H$. 
From this we may infer, by reason of $H=\{ s\in{\Bbb C} : {\rm Re}(s)>3/4\}$ 
being contained in the union of 
its compact subsets, that the 
equation (6.5.61) defines a holomorphic function 
$\zeta^{{\frak a},{\frak b}}_{\omega,\omega'} : H\rightarrow{\Bbb C}$. 
\par 
Finally, by (6.5.61), (6.5.68), (6.5.67), the 
upper bound (6.5.59) of Lemma~6.5.9, and the result in (6.1.25), 
we find that, for all $s\in{\cal D}$, one has 
$$\zeta^{{\frak a},{\frak b}}_{\omega,\omega'}(s) 
\ll_{r_1} \,\exp\!\left( 4\pi\sqrt{\left|\omega\omega'\right|}\right) 
|\omega|\left| m_{\frak a} m_{\frak b}\right|^{1/2} 
\biggl(\,\sum_{\alpha\mid q_0} {1\over |\alpha|}\,\biggr)^{\!\!2}  
\,\zeta_{{\Bbb Q}(i)}^2\!\left( 2\sigma^{*} - 1/2\right) ,$$
where $| m_{\frak a} m_{\frak b}|\leq |q_0|^2$. The result (6.5.66) follows 
from this: for $2\sigma^{*}-1/2>1$, and we have also that  
$r_1$ and $\sigma^{*}$ need depend only on ${\cal D}$, while (since $\Gamma =\Gamma_0(q_0)$) 
the group $\Gamma$ determines the ideal 
$q_0 {\frak O}\subset{\frak O}$\quad$\blacksquare$ 

\bigskip 

\proclaim Lemma 6.5.13. Let $0\neq\omega\in{\frak O}$. 
Suppose that ${\frak a}\sim^{\!\!\!\!\Gamma}{\frak b}$, and let $\gamma\in\Gamma$ 
be such that $\gamma{\frak b}={\frak a}$. Then there exists a unique pair 
$(u,z)=(u(\gamma),z(\gamma))\in{\Bbb C}^{*}\times{\Bbb C}$ such that 
$g_{\frak a}^{-1}\gamma g_{\frak b}=h[u]n[z]$; the relevant $u\in{\Bbb C}^{*}$ 
is a square root of some unit $\epsilon=\epsilon(\gamma)\in{\frak O}^{*}$,  
and one has: 
$$\bigl( {\bf M}_{\omega}\varphi_{\ell,q}(\nu,0)\bigr)  
(g_{\frak a}^{-1}\gamma g_{\frak b} g) 
=\psi_{\epsilon\omega}\!\left( n[z]\right) 
\bigl({\bf M}_{\epsilon\omega}\varphi_{\ell,q}(\nu,0)\bigr)(g)\qquad\quad 
\hbox{($\nu\in{\Bbb C}$, $g\in G$),}\eqno(6.5.69)$$
$$\sum_{\omega'\in{\frak O}} \delta_{\omega u , \omega'/u} 
=\sum_{\omega'\in{\frak O}} \delta_{\epsilon\omega , \omega'} =1\eqno(6.5.70)$$
and 
$$\Gamma_{\frak a}\gamma =\left\{ \gamma_1\in\Gamma : \gamma_1{\frak b}={\frak a}\right\}\;. 
\eqno(6.5.71)$$ 

\medskip 

\noindent{\bf Proof.}\quad  
The results preceding that in (6.5.69) are contained in Lemma~4.1;  
when combined with (6.3.5) (for $p=0$) and (6.3.3), those results imply 
the equality in (6.5.69). The result (6.5.70) is self-evident 
(given that $u^2=\epsilon\in{\frak O}^{*}$). 
To verify the equality in (6.5.71) one need only 
note the equivalence (for $\gamma_1\in\Gamma$) of the four binary relations 
$\gamma_1{\frak b}={\frak a}$, $\gamma_1{\gamma}^{-1}{\frak a}={\frak a}$, 
$\gamma_1\gamma^{-1}\in\Gamma_{\frak a}$ and $\gamma_1\in\Gamma_{\frak a}\gamma$
\quad$\blacksquare$ 

\bigskip 

By Lemma~6.5.12 and Lemma~6.5.3 we may define, for ${\rm Re}(\nu)>1/2$, 
$g\in G$ and $\omega'\in{\frak O}$, the term 
$$\eqalign{
\phi_{\omega'}(\nu,g) 
 &=\,\zeta^{{\frak a},{\frak b}}_{\omega,\omega'}\!\!\left({1+\nu\over 2}\right) 
\bigl({\bf J}_{\omega'}\varphi_{\ell,q}(\nu,0)\bigr)(g)\  + \cr 
 &\quad\ +{1^{\textstyle{\ \atop\ }}\over\left[\Gamma_{\frak a} : \Gamma_{\frak a}'\right]} 
\sum_{\scriptstyle\gamma\in\Gamma_{\frak a}'\backslash\Gamma 
\ :\ \gamma{\frak b}={\frak a}\atop
{\scriptstyle g_{\frak a}^{-1}\gamma g_{\frak b}\in h\left[ u(\gamma)\right] N  
\atop\scriptstyle\hbox{\qquad}}} 
\delta_{\omega u(\gamma) , \omega'/u(\gamma)} 
\,\bigl({\bf M}_{\omega}\varphi_{\ell,q}(\nu,0)\bigr)
\bigl( g_{\frak a}^{-1}\gamma g_{\frak b} g\bigr)\;.}\eqno(6.5.72)$$
When $(\nu,g)\in{\Bbb C}\times G$ is such that the series 
$\sum_{\omega'\in{\frak O}}\phi_{\omega'}(\nu,g)$ is convergent, 
we denote the sum of that series by $\Phi(\nu,g)$. 
Both $\phi_{\omega'}(\nu,g)$ and $\Phi(\nu,g)$ are (of course) 
dependent on $\Gamma$, $g_{\frak a}$, $g_{\frak b}$, 
$\omega$ and the $K$-type $(\ell,q)$, but, since these other parameters 
are (for the purposes of the current discussion) effectively constants,  
our main concern in the next lemma is with the dependence of 
$\phi_{\omega'}(\nu,g)$ and $\Phi(\nu,g)$ 
on the pair $(\nu,g)\in{\Bbb C}\times G$. 

\bigskip 

\proclaim Lemma 6.5.14. Let $0\neq\omega\in{\frak O}$. 
Then, when $r_0,\sigma_2,t_1\in(0,\infty)$, $\sigma_1\in(1/2,\sigma_2)$,  
$G(r_0)=\{ g\in G : \rho(g)\geq r_0\}$ and 
${\cal R}={\cal R}(\sigma_1,\sigma_2,t_1)
=\{\nu\in{\Bbb C} : \sigma_1\leq{\rm Re}(\nu)\leq\sigma_2\ {\rm and}\ 
|{\rm Im}(\nu)|\leq t_1\}$, the series 
$\sum_{\omega'\in{\frak O}}|\phi_{\omega'}(\nu,g)|$ is uniformly 
convergent for all $(\nu,g)\in{\cal R}\times G(r_0)$, and its sum 
$\Phi(\nu,g)$ satisfies 
$$\Phi(\nu,g) 
-{1^{\textstyle{\ \atop\ }}\over\left[\Gamma_{\frak a} : \Gamma_{\frak a}'\right]} 
\sum_{\scriptstyle\gamma\in\Gamma_{\frak a}'\backslash\Gamma\atop\scriptstyle 
\gamma{\frak b}={\frak a}} 
\bigl({\bf M}_{\omega}\varphi_{\ell,q}(\nu,0)\bigr)
\left( g_{\frak a}^{-1}\gamma g_{\frak b} g\right)
	\ll_{\Gamma,\omega,\ell,{\cal R},r_0}\ 
\left(\rho(g)\right)^{1-{\rm Re}(\nu)}\;,\eqno(6.5.73)$$ 
for $(\nu,g)\in{\cal R}\times G(r_0)$. In particular, the series 
$\sum_{\omega'\in{\frak O}}\phi_{\omega'}(\nu,g)$ is 
absolutely convergent for all $(\nu,g)\in{\Bbb C}\times G$ such that 
${\rm Re}(\nu)>1/2$. Properties of its sum $\Phi(\nu,g)$ are:  
$$\eqalignno{ 
 &\hbox{when $g\in G$, the function $\nu\mapsto\Phi(\nu,g)$ is 
holomorphic for ${\rm Re}(\nu)>1/2$;} &(6.5.74)\cr 
 &\hbox{the function $(\nu,g)\mapsto\Phi(\nu,g)$ is continuous 
on $\{\nu\in{\Bbb C} : {\rm Re}(\nu)>1/2\}\times G$;} &(6.5.75)\cr 
 &\hbox{when $\,{\rm Re}(\nu)>1$, one has} 
\ \ \Phi(\nu,g)=\bigl( P^{\frak a}{\bf M}_{\omega}\varphi_{\ell,q}(\nu,0)\bigr)
\!\left( g_{\frak b} g\right)\ 
\,\hbox{for all $g\in G$.} &(6.5.76)}$$ 

\medskip 

\noindent{\bf Proof.}\quad 
By Lemma~4.2, we have $[\Gamma_{\frak a} : \Gamma_{\frak a}']\in\{ 2,4\}$, so 
that (by (6.5.71)) the sum over `$\gamma$' in Equation~(6.5.72) is always finite. 
That sum is, moreover, empty unless $\omega'\sim\omega$: for it follows from Lemma~6.5.13 
that $u(\gamma)$, in (6.5.72), is always a square root of a unit of ${\frak O}$. 
Therefore, assuming now that $\sigma_1$, $\sigma_2$, $r_0$, $t_1$, ${\cal R}$ and $G(r_0)$, 
are as specified in the lemma, and that 
$$\phi^{*}_{\omega'}(\nu,g)=\,\zeta^{{\frak a},{\frak b}}_{\omega,\omega'}\!\!\left({1+\nu\over 2}\right) 
\bigl({\bf J}_{\omega'}\varphi_{\ell,q}(\nu,0)\bigr) (g)\qquad\quad 
\hbox{(for ${\rm Re}(\nu)>1/2$, $g\in G$ and $\omega'\in{\frak O}$),}$$
it will suffice for proof of the first result of the 
lemma that we show the series 
$\sum_{0\neq\omega'\in{\frak O}}|\phi^{*}_{\omega'}(\nu,g)|$ 
to be uniformly convergent for $(\nu,g)\in{\cal R}\times G(r_0)$. 
\par 
Since the conditions ${\rm Re}(\nu)\geq\sigma_1 >1/2$ imply 
${\rm Re}((\nu+1)/2))>3/4$, it follows by Lemma~6.5.12 that 
$$\zeta^{{\frak a},{\frak b}}_{\omega,\omega'}\!\!\left({1+\nu\over 2}\right) 
\ll_{\Gamma,{\cal R}}\ |\omega| e^{4\pi\sqrt{|\omega\omega'|}}\qquad\ 
\hbox{for $\,\nu\in{\cal R}$, $\omega'\in{\frak O}$.}\eqno(6.5.77)$$
By the estimate (6.5.14) of Lemma~6.5.3, we have also the upper bound 
$$\bigl({\bf J}_{\omega'}\varphi_{\ell,q}(\nu,0)\bigr) (g) 
\ll_{\ell,\sigma_2,r_0} 
|\Gamma(\nu +1)|^{-1} \left|\omega'\right|^{\ell+{\rm Re}(\nu)} 
\left(\rho(g)\right)^{\ell+1} e^{-2\pi |\omega'|\rho(g)}\;,$$
for $0\neq\omega'\in{\frak O}$ and $(\nu,g)\in{\cal R}\times G(r_0)$: 
we rely here on the fact the conditions imply the lower bound  
$|\omega'|\geq 1$, and so ensure that $|\omega'|\rho(g)\geq r_0$ when $g\in G(r_0)$. 
Given the definition of the compact set ${\cal R}\subset{\Bbb C}$, 
it follows from this last bound that, since $1/\Gamma(w)$ is an entire complex function,  
one has:  
$$\bigl({\bf J}_{\omega'}\varphi_{\ell,q}(\nu,0)\bigr) (g) 
\ll_{\ell,{\cal R},r_0} \left|\omega'\right|^{\ell+\sigma_2} 
\left(\rho(g)\right)^{\ell+1} e^{-2\pi |\omega'|\rho(g)}\qquad 
\hbox{for $0\neq\omega'\in{\frak O}$, $\,(\nu,g)\in{\cal R}\times G(r_0)$.}\eqno(6.5.78)$$
Combining the bounds in (6.5.77) and (6.5.78) we find that,  
if $0\neq\omega'\in{\frak O}$ and $\,(\nu,g)\in{\cal R}\times G(r_0)$, then 
$$\eqalignno{
\phi^{*}_{\omega'}(\nu,g) 
 &\ll_{\Gamma,\ell,{\cal R},r_0}  
\,|\omega| \left|\omega'\right|^{\ell+\sigma_2}  
\left(\rho(g)\right)^{\ell+1} e^{-\pi(2\rho(g)|\omega'| 
-4\sqrt{|\omega\omega'|}\,)}  \leq \cr
 &\leq\left|\omega/\omega'\right| e^{4\pi |\omega|/r_0} r_0^{-\sigma_2} 
\left(\rho(g)\left|\omega'\right|\right)^{\ell+1+\sigma_2} 
e^{-\pi\rho(g)|\omega'|} \ll_{\ell,\omega,\sigma_2,r_0} \cr 
 &\ll_{\ell,\omega,\sigma_2,r_0} 
e^{-(\pi/2)\rho(g)|\omega'|}\leq e^{-(\pi/2)r_0|\omega'|} 
\leq{48\over\pi^3 r_0^3 |\omega'|^3}\;. &(6.5.79)}$$
From this follows (by virtue of the convergence of the 
series $\sum_{0\neq\omega'\in{\frak O}}|\omega'|^{-3}$) the 
uniform convergence, for 
$(\nu,g)\in{\cal R}\times G(r_0)$, of the series 
$\sum_{0\neq\omega'\in{\frak O}}|\phi^{*}_{\omega'}(\nu,g)|$. 
By earlier remarks, the first result of the lemma follows. 
\par 
We prove next the estimate (6.5.73). By (6.5.79) it is seen 
that, when $(\nu,g)\in{\cal R}\times G(r_0)$, one has: 
$$\eqalign{ 
\sum_{0\neq\omega'\in{\frak O}}\phi^{*}_{\omega'}(\nu,g) 
 &\ll_{\Gamma,\ell,\omega,{\cal R},r_0}  
\sum_{0\neq\omega'\in{\frak O}} e^{-(\pi/2)\rho(g)|\omega'|} \leq \cr 
 &\qquad\leq e^{-(\pi/3)\rho(g)}\sum_{0\neq\omega'\in{\frak O}} 
\biggl( {4!\over r_0^4\left|\omega'\right|^4}\biggr)^{\!\!\pi/6} 
\ll_{r_0}  e^{-(\pi/3)\rho(g)}\leq e^{-\rho(g)}\;.}$$
By (6.5.77), (1.5.18) and (1.3.2)  one has, moreover,   
$$\phi^{*}_0(\nu,g) 
=\,\zeta^{{\frak a},{\frak b}}_{\omega,0}\!\!\left({1+\nu\over 2}\right) 
\varphi_{\ell,q}(-\nu,0)(g)
{\pi\over\nu}\prod_{m=1}^{\ell}\left( {m-\nu\over m+\nu}\right) 
\ll_{\Gamma,\ell,{\cal R}}\ |\omega|\left(\rho(g)\right)^{1-{\rm Re}(\nu)}\;,$$ 
when $(\nu,g)\in{\cal R}\times G(r_0)$. We therefore find that 
$$\sum_{\omega'\in{\frak O}}\phi^{*}_{\omega'}(\nu,g) 
\ll_{\Gamma,\ell,\omega,{\cal R},r_0}   
\,\left(\rho(g)\right)^{1-{\rm Re}(\nu)} +e^{-\rho(g)} 
\ll_{\sigma_2,r_0}\,\left(\rho(g)\right)^{1-{\rm Re}(\nu)}\qquad\ 
\hbox{for $(\nu,g)\in{\cal R}\times G(r_0)$,}$$
and from this the result (6.5.73) follows: for, by (6.5.72) and result  
(6.5.70) of Lemma~6.5.13, it follows 
(certainly for all $(\nu,g)\in{\cal R}\times G(r_0)\,$) that one has 
$$\Phi(\nu,g)-\sum_{\omega'\in{\frak O}}\phi^{*}_{\omega'}(\nu,g)  
=\sum_{\omega'\in{\frak O}}\left(\phi_{\omega'}(\nu,g) 
-\phi^{*}_{\omega'}(\nu,g)\right) 
={1^{\textstyle{\ \atop\ }}\over\left[\Gamma_{\frak a} : \Gamma_{\frak a}'\right]} 
\sum_{\scriptstyle\gamma\in\Gamma_{\frak a}'\backslash\Gamma\atop\scriptstyle 
\gamma{\frak b}={\frak a}} 
\bigl({\bf M}_{\omega}\varphi_{\ell,q}(\nu,0)\bigr)
\left( g_{\frak a}^{-1}\gamma g_{\frak b} g\right)\;.$$ 
\par 
The first two results of the lemma (stated in and above (6.5.73)) have now been proved,  
and so, in proving the other results of the lemma, we may freely apply 
those first two results (in their full generality). 
In particular, the third result of the lemma 
(absolute convergence of the series $\sum_{\omega'\in{\frak O}}\phi_{\omega'}(\nu,g)$ 
when $(\nu,g)\in{\Bbb C}\times G$ 
and ${\rm Re}(\nu)>1/2$) is an immediate corollary of the 
case $\sigma_1={\rm Re}(\nu)$, $\sigma_2={\rm Re}(\nu+1)$, $t_1=|{\rm Im}(\nu)|+1$, 
$r_0=\rho(g)$ of 
the uniform convergence noted in the first result of the lemma. 
\par 
Assume henceforth that $(\nu_0,g_0)\in{\Bbb C}\times G$, and that ${\rm Re}(\nu_0)>1/2$. 
Given the definition (6.5.72), it follows by (1.5.18), Lemma~6.5.3 and Lemma~6.5.12 
that each term of the series $\sum_{\omega'\in{\frak O}}\phi_{\omega'}(\nu,g_0)$ 
is holomorphic (as a function of the complex variable $\nu$) in the half-plane 
where ${\rm Re}(\nu)>1/2$. Therefore (as an application of the first result of the lemma) 
we may deduce from the uniform convergence of the series 
$\sum_{\omega'\in{\frak O}}\phi_{\omega'}(\nu,g)$
for 
$(\nu,g)\in{\cal R}({\rm Re}(\nu_0/2)+1/4,{\rm Re}(2\nu_0),|{\rm Im}(\nu_0)|+1) 
\times\{ g_0\}\subset 
{\cal R}({\rm Re}(\nu_0/2)+1/4,{\rm Re}(2\nu_0),|{\rm Im}(\nu_0)|+1) 
\times G(\rho(g_0))$ 
that the function $\nu\mapsto\Phi(\nu,g_0)
=\sum_{\omega'\in{\frak O}}\phi_{\omega'}(\nu,g_0)$ 
is holomorphic in a neighbourhood of the point $\nu_0\in{\Bbb C}$. 
This proves the result (6.5.74): for we have assumed nothing more than that 
${\rm Re}(\nu_0)>1/2$ and $g_0\in G$. 
\par 
Considering now the series $\sum_{\omega'\in{\frak O}}\phi_{\omega'}(\nu_0,g)$, 
we observe that, given the definition (6.5.72), it follows by (1.5.15) and 
the observation preceding (6.3.4) that each term of this series 
is a continuous (as a function of $g$) on all of $G$. Therefore we may deduce 
from the uniform convergence of the series 
$\sum_{\omega'\in{\frak O}}\phi_{\omega'}(\nu,g)$
for $(\nu,g)\in\{\nu_0\}\times G(\rho(g_0)/2)\subset 
{\cal R}({\rm Re}(\nu_0),{\rm Re}(\nu_0)+1,|{\rm Im}(\nu_0)|+1)\times G(\rho(g_0)/2)$ 
(which follows by the first result of the lemma) 
that the function $g\mapsto\Phi(\nu_0,g)
=\sum_{\omega'\in{\frak O}}\phi_{\omega'}(\nu_0,g)$ 
is continuous at the point $g_0\in G$. 
We may infer from this that, when ${\rm Re}(\nu)>1/2$, the function 
$g\mapsto\Phi(\nu,g)$ is continuous on $G$. This advances us 
one step towards a proof of the result (6.5.75). 
\par 
A locally uniform upper bound for $\Phi(\nu,g)$ suffices 
to complete the proof of (6.5.75). We begin our proof of such a 
bound by putting: $\sigma_1={\rm Re}(\nu_0/3)+1/3$, $\sigma_2={\rm Re}(3\nu_0)$, 
$t_1=|{\rm Im}(\nu_0)|+2$, $r_0=\rho(g_0)/2$, $r_1=2\rho(g_0)$, 
${\cal R}={\cal R}(\sigma_1,\sigma_2,t_1)
=\{\nu\in{\Bbb C} : \sigma_1\leq{\rm Re}(\nu)\leq\sigma_2\ {\rm and}\ 
|{\rm Im}(\nu)|\leq t_1\}$ and 
$G(r_0,r_1)=\{ g\in G : r_0\leq\rho(g)\leq r_1\}$. 
By the result (6.5.73) (proved earlier), the results of Lemma~6.5.13 
(excluding (6.5.70)) and the estimate (6.3.9) of Lemma~6.3.1, 
we find that if $(\nu,g)\in{\cal R}\times G(r_0,r_1)$ then, 
given the definition of the set ${\cal R}={\cal R}(\sigma_1,\sigma_2,t_1)\subset{\Bbb C}$, one has    
$$\eqalign{ 
|\Phi(\nu,g)| &\leq\max_{\epsilon\in{\frak O}^{*}} 
\left|\bigl( {\bf M}_{\epsilon\omega}\varphi_{\ell,q}(\nu,0)\bigr) (g)\right| 
+O_{\Gamma,\omega,\ell,{\cal R},r_0}\left(\left(\rho(g)\right)^{1-{\rm Re}(\nu)}\right) =\cr 
 &=O_{\omega,\ell,r_1,\sigma_2}\left( \left(\rho(g)\right)^{1+{\rm Re}(\nu)} 
\left( 1+|{\rm Im}(\nu)|\right)^{-2{\rm Re}(\nu)-1} 
e^{\pi|{\rm Im}(\nu)|}\right) 
+O_{\Gamma,\omega,\ell,{\cal R},r_0}\left(\left(\rho(g)\right)^{1-{\rm Re}(\nu)}\right) , 
}$$
$$0<\left(\rho(g)\right)^{1+{\rm Re}(\nu)} 
\left( 1+|{\rm Im}(\nu)|\right)^{-2{\rm Re}(\nu)-1} 
e^{\pi|{\rm Im}(\nu)|}\leq\exp\left(\left( 1+\sigma_2\right) r_1 +\pi t_1\right)$$
and 
$$0<\left(\rho(g)\right)^{1-{\rm Re}(\nu)}\leq r_1\exp\left(\sigma_2/r_0\right) .$$ 
One sees from this that the function 
$(\nu,g)\mapsto\Phi(\nu,g)$ is bounded on ${\cal R}(\sigma_1,\sigma_2,t_1)\times G(r_0,r_1)$. 
It therefore follows, 
by an application of Cauchy's integral formula 
for the derivative of a holomorphic function, that the function 
$(\nu,g)\mapsto (\partial/\partial\nu)\Phi(\nu,g)$ is bounded on 
the set ${\cal R}'\times G(r_0,r_1)$, where 
${\cal R}'={\cal R}(2\sigma_1-1/2,2\sigma_2 /3,t_1-1)$ $=\{\nu\in{\Bbb C} : {\rm Re}(2\nu_0/3)+1/6\leq{\rm Re}(\nu)\leq {\rm Re}(2\nu_0)\ 
{\rm and}\ |{\rm Im}(\nu)|\leq |{\rm Im}(\nu_0)|+1\}$; consequently there exists 
some $M\in(0,\infty)$ such that,  
for all $\nu\in{\cal R}'$ and all $g\in G(r_0,r_1)$, one has  
$|\Phi(\nu,g)-\Phi(\nu_0,g)|\leq M |\nu-\nu_0|\,$ 
(it should be noted here that, since ${\rm Re}(\nu_0)>1/2$, 
the definitions of $G(r_0,r_1)$ and ${\cal R}'$ ensure 
that the set ${\cal R}'\times G(r_0,r_1)$ contains 
a neighbourhood of the point $(\nu_0,g_0)\in{\Bbb C}\times G$).  
Hence, by the triangle inequality,  
$$\left|\Phi(\nu,g)-\Phi\left(\nu_0,g_0\right)\right| 
\leq M\left|\nu -\nu_0\right|  
+\left|\Phi\left(\nu_0,g\right)-\Phi\left(\nu_0,g_0\right)\right|\qquad\ 
\hbox{for $(\nu,g)\in{\cal R}'\times G(r_0,r_1)$.}\eqno(6.5.80)$$
\par 
Suppose now that $\varepsilon >0$. Since we assume that ${\rm Re}(\nu_0)>1/2$, 
we know (by work done above) that  
the function $g\mapsto\Phi(\nu_0,g)$ is continuous on $G$. 
Hence there exists a neightbourhood $U_0$ of $g_0$ such that 
$|\Phi(\nu_0,g)-\Phi(\nu_0,g_0)|<\varepsilon/2$ for all $g\in U_0$; 
we may assume, moreover, that $U_0\subset G(r_0,r_1)\,$ 
(if necessary we replace $U_0$ by ${\rm Int}(G(r_0,r_1))\cap U_0$).  
By this observation, combined with (6.5.80), we find that 
$|\Phi(\nu,g)-\Phi(\nu_0,g_0)|<\varepsilon$ for all 
$(\nu,g)$ lying in the neighbourhood 
$\{\nu\in{\rm Int}({\cal R}') : |\nu-\nu_0|<\varepsilon /(2M)\}\times U_0$ 
of the point $(\nu_0,g_0)$; since our only assumptions concerning 
$\nu_0\in{\Bbb C}$, $g_0\in G$ and $\varepsilon\in{\Bbb R}$ are that 
we have ${\rm Re}(\nu_0)>1/2$ and $\varepsilon>0$, this completes the proof of the 
result (6.5.75). 
\par 
We now aim to complete the proof of the lemma, by proving the result (6.5.76). 
Accordingly, it is to be 
assumed henceforth that we have ${\rm Re}(\nu)>1$ and $g\in G$. 
Since ${\rm Re}(\nu)>1$, it follows by the result (6.5.60) of Lemma~6.5.10, 
combined with the definitions (6.5.72), (6.2.13), (1.4.3), (1.1.21), (1.1.10) 
and (1.1.3), that for $(n_1,n_2)\in{\Bbb Z}\times{\Bbb Z}$ and $\omega'=-n_1+in_2$ 
we have:  
$$\phi_{-n_1+in_2}(\nu,g) 
=\bigl( F^{\frak b}_{\omega'} P^{\frak a} {\bf M}_{\omega}\varphi_{\ell,q}(\nu,0)\bigr) (g) 
=\int_{-\infty}^{\infty}\int_{-\infty}^{\infty} f\left( x_1 , x_2\right) 
{\rm e}\left( n_1 x_1 + n_2 x_2\right) {\rm d}x_1\,{\rm d}x_2\;,$$
where the function $f : {\Bbb R}\times{\Bbb R}\rightarrow{\Bbb C}$ is defined by: 
$$f\left( x_1 , x_2\right) 
=\cases{\bigl( P^{\frak a}{\bf M}_{\omega}\varphi_{\ell,q}(\nu,0)\bigr) 
(g_{\frak b} n[x_1+ix_2] g) &if $\,-1/2\leq x_1,x_2<1/2$; \cr 
0 &otherwise.}\eqno(6.5.81)$$ 
We now seek to apply the two-variable case of Bochner's theorem,  
Theorem~67 of [2], on Poisson summation in several variables. 
By the first result of the lemma (proved in the paragraph containing (6.5.79)), 
we know that the series 
$\sum_{\omega'\in{\frak O}}\phi_{\omega'}(\nu,g)$ is absolutely convergent, 
and so `Assumption (c)' of Bochner's theorem is satisfied. 
The function $f$ defined in (6.5.81) may be shown also 
to satisfy the other hypotheses of the case $k=2$ of Bochner's theorem 
(i.e. his `Assumptions (a) and (b)'). Indeed, with regard to 
Assumption~(a) of Theorem~67 of [2] 
we need only note that, since ${\rm Re}(\nu)>1$, it follows by 
Lemma~6.5.1 that $P^{\frak a}{\bf M}_{\omega}\varphi_{\ell,q}(\nu,0)$ 
is a continuous function on the topological group $G$;  with regard to 
Assumption (b) of Theorem~67 of [2] it is enough to observe 
(as a consequence of (6.5.81)) that if $-1/2\leq x,y<1/2$ then  
$f(x+m,y+n)=0$ for all $(m,n)\in{\Bbb Z}\times{\Bbb Z}-\{ (0,0)\}$. 
\par 
Bochner's theorem justifies the application 
of the two variable form of Poisson's summation formula; 
in light of the point just observed in 
connection with his `Assumption (b)', we therefore find that 
$$\sum_{n_1=-\infty}^{\infty}\ \sum_{n_2=-\infty}^{\infty}
\phi_{-n_1+in_2}(\nu,g)  
=\,\sum_{m_1=-\infty}^{\infty}\ \sum_{m_2=-\infty}^{\infty} f\left( m_1 , m_2\right) 
=f(0,0)\;.$$
The sum over $n_1$ and $n_2$ here is $\Phi(\nu,g)$, and so 
(by (6.5.81) for $x_1=x_2=0$) the result (6.5.76) follows\quad$\blacksquare$ 

\bigskip 

The results (6.5.76) and (6.5.74) of Lemma~6.5.14 show 
that the sum 
$\Phi(\nu,g_{\frak b}^{-1}g) 
=\sum_{\omega'\in{\frak O}}\phi_{\omega'}(\nu,g_{\frak b}^{-1}g)$ 
provides a means of extending the domain of the function 
$(\nu,g)\mapsto (P^{\frak a}{\bf M}_{\omega}\varphi_{\ell,q}(\nu,0))(g)$ 
that is `natural' (in that the extended function is determined by 
a process of analytic continuation). 
In referring to the extended function one might avoid the use of any new terminology. 
However, in the interest of clarity, we choose to let   
${\cal P}^{\frak a}_{\!\!\!\scriptscriptstyle\hookleftarrow}{\bf M}_{\omega}\varphi_{\ell,q}(\nu,0)$ denote, when ${\rm Re}(\nu)>1/2$,  
the  function on $G$ satisfying 
$$\bigl( {\cal P}^{\frak a}_{\!\!\!\scriptscriptstyle\hookleftarrow}
{\bf M}_{\omega}\varphi_{\ell,q}(\nu,0)\bigr) (g) 
=\Phi\left(\nu,g_{\frak b}^{-1}g\right) 
=\sum_{\omega'\in{\frak O}}\phi_{\omega'}\left(\nu,g_{\frak b}^{-1}g\right)\qquad\quad 
\hbox{($g\in G$),}\eqno(6.5.82)$$ 
where the terms of the sum over $\omega'\in{\frak O}$ are as indicated in  
Equation~(6.5.72). Note that each function 
${\cal P}^{\frak a}_{\!\!\!\scriptscriptstyle\hookleftarrow}
{\bf M}_{\omega}\varphi_{\ell,q}(\nu,0)$ 
defined in this way is independent of the particular choice of cusp ${\frak b}$
and scaling matrix $g_{\frak b}$: for it follows by (6.5.82) and the results (6.5.76) 
and (6.5.74) of 
Lemma~6.5.14 that 
$$\bigl( {\cal P}^{\frak a}_{\!\!\!\scriptscriptstyle\hookleftarrow}
{\bf M}_{\omega}\varphi_{\ell,q}(\nu,0)\bigr) (g)  
=\bigl( P^{\frak a}{\bf M}_{\omega}\varphi_{\ell,q}(\nu,0)\bigr) (g) \qquad\ 
\hbox{for $(\nu,g)\in{\Bbb C}\times G$ such that ${\rm Re}(\nu)>1$,}\eqno(6.5.83)$$
and that, for each $g\in G$, the function 
$\nu\mapsto\bigl( {\cal P}^{\frak a}_{\!\!\!\scriptscriptstyle\hookleftarrow}
{\bf M}_{\omega}\varphi_{\ell,q}(\nu,0)\bigr) (g)$ 
is holomorphic for ${\rm Re}(\nu)>1/2$, and so is   
completely determined (as a function with domain $\{ \nu\in{\Bbb C} : {\rm Re}(\nu)>1/2\}\,$) 
by the data in (6.5.83).
\par 
The new terminology just introduced aids in the clarification of 
certain steps in subsequent proofs. Although ${\frak b}$, $g_{\frak b}$ and 
the $K$-type $(\ell,q)$ are fixed (for the purposes of the present discussion), 
we nevertheless take (6.5.82) to imply corresponding definitions of 
${\cal P}^{\frak a}_{\!\!\!\scriptscriptstyle\hookleftarrow}
{\bf M}_{\omega}\varphi_{\lambda,\kappa}(\nu,0)$ for all $\kappa,\lambda\in{\Bbb Z}$ 
such that $\lambda\geq|\kappa|$, and we allow substitution of 
other pairs ${\frak c}$, $g_{\frak c}$ in place of the pair ${\frak b}$, $g_{\frak b}$ 
(so long as the conditions (1.1.16) and (1.1.20) are satisfied); this obviously 
has to be accompanied by matching substitutions on the right-hand side 
of Equation (6.5.72), where the term $\phi_{\omega'}(\nu,g)$ is defined; in order 
to make this quite clear, we let $\phi^{\frak c}_{\omega'}(\lambda,\kappa;\nu,g)$ 
denote the 
term $\phi_{\omega'}(\nu,g)$ which Equation~(6.5.72) would define were ${\frak c}$, 
$g_{\frak c}$ and 
the $K$-type $(\lambda,\kappa)$ substituted for ${\frak b}$, $g_{\frak b}$ and the 
$K$-type $(\ell,q)$, respectively, and we also put 
$\Phi^{\frak c}(\lambda,\kappa;\nu,g)=\sum_{\omega'\neq 0}
\phi^{\frak c}_{\omega'}(\lambda,\kappa;\nu,g)$. 
\par 
By the case ${\frak b}=\infty$, $g_{\frak b}=h[1]$ of 
(6.5.82), (6.5.72), Lemma~6.5.13 and Lemma~6.5.14, it follows that,  
when $0\neq\omega\in{\frak O}$, $\lambda,\kappa\in{\Bbb Z}$, $\lambda\geq|\kappa|$, 
$g\in G$ and ${\rm Re}(\nu)>1/2$, we have: 
$$\eqalign{
\bigl({\cal P}^{\frak a}_{\!\!\!\scriptscriptstyle\hookleftarrow}
{\bf M}_{\omega}\varphi_{\lambda,\kappa}(\nu,0)\bigr)(g) 
 &={1^{\textstyle{\ \atop\ }}\over\left[\Gamma_{\frak a} : \Gamma_{\frak a}'\right]} 
\sum_{\scriptstyle\gamma\in\Gamma_{\frak a}'\backslash\Gamma\atop\scriptstyle 
\gamma\infty={\frak a}} 
\bigl({\bf M}_{\omega}\varphi_{\lambda,\kappa}(\nu,0)\bigr)
\bigl( g_{\frak a}^{-1}\gamma g\bigr)\  + \cr 
 &\quad +\sum_{\omega'\in{\frak O}} 
\zeta^{{\frak a},\infty}_{\omega,\omega'}\!\!\left({1+\nu\over 2}\right) 
\bigl({\bf J}_{\omega'}\varphi_{\lambda,\kappa}(\nu,0)\bigr) (g)\;.}\eqno(6.5.84)$$
We remark that the sum over $\gamma\in\Gamma_{\frak a}'\backslash\Gamma$ in 
Equation~(6.5.84) is evidently empty unless the cusp ${\frak a}$ is 
$\Gamma$-equivalent to the cusp $\infty$; by 
Lemma~4.1, Lemma~4.2, the equations~(6.3.3) and~(6.3.5) and the definition (1.4.3), 
it follows moreover that if ${\frak a}\sim^{\!\!\!\!\Gamma}\infty$ then there exists some 
$\epsilon\in{\frak O}^{*}$ and some $\beta\in{\Bbb C}$ 
(both depending only upon $\Gamma$ and $g_{\frak a}$) such that 
$${1^{\textstyle{\ \atop\ }}\over\left[\Gamma_{\frak a} : \Gamma_{\frak a}'\right]} 
\sum_{\scriptstyle\gamma\in\Gamma_{\frak a}'\backslash\Gamma\atop\scriptstyle 
\gamma\infty={\frak a}} 
\bigl({\bf M}_{\omega}\varphi_{\lambda,\kappa}(\nu,0)\bigr)
\left( g_{\frak a}^{-1}\gamma g\right) 
={\rm e}\!\left({\rm Re}(\beta\epsilon\omega)\right) 
\sum_{\alpha=\pm 1}
\,{1\over 2}\bigl({\bf M}_{\alpha\epsilon\omega}\varphi_{\lambda,\kappa}(\nu,0)\bigr)(g)\;,$$ 
when the conditions attached to (6.5.84) are satisfied. 

\bigskip 

\proclaim Lemma 6.5.15. The equation (6.5.82) defines, for ${\rm Re}(\nu)>1/2$, 
a function 
${\cal P}^{\frak a}_{\!\!\!\scriptscriptstyle\hookleftarrow}
{\bf M}_{\omega}\varphi_{\ell,q}(\nu,0) : G\rightarrow{\Bbb C}$ 
which lies in the space $C^{\infty}(\Gamma\backslash G)$ and is 
of $K$-type $(\ell,q)$. 

\medskip 

\noindent{\bf Proof.}\quad 
To avoid both ambiguity and unnecessary repetition we make it our rule 
that, when (in the course of this proof) there is any application 
made of some previous result (or definition) of this subsection, it is to  
be understood that the application in question is made in respect of the 
special case ${\frak b}=\infty$, $g_{\frak b}=h[1]$. 
\par 
Let $\nu_0\in{\Bbb C}$ be such that ${\rm Re}(\nu_0)>1/2$. Then, by 
Lemma~6.5.14, the equation (6.5.82) defines a 
complex valued function 
$g\mapsto({\cal P}^{\frak a}_{\!\!\!\scriptscriptstyle\hookleftarrow}
{\bf M}_{\omega}\varphi_{\ell,q}(\nu_0,0))(g)$ with domain $G$. 
By Lemma~6.5.1, and the results (6.5.76) and (6.5.74) of 
Lemma~6.5.14, it may moreover be deduced that 
this function ${\cal P}^{\frak a}_{\!\!\!\scriptscriptstyle\hookleftarrow}
{\bf M}_{\omega}\varphi_{\ell,q}(\nu_0,0) : G\rightarrow{\Bbb C}$ 
is $\Gamma$-automorphic (i.e. we obtain this not only when 
${\rm Re}(\nu_0)>1$, but also when $1/2<{\rm Re}(\nu_0)\leq 1$). 
\par 
Suppose that ${\bf X}\in{\cal B}_1=\{ {\bf H}_1, 
{\bf H}_2, {\bf F}^{+}, {\bf F}^{-}, {\bf E}^{+}, {\bf E}^{-} \}\,$ 
(the basis of ${\frak g}$ utilised in the 
proof of Lemma~6.5.5). Since the Jacquet operator ${\bf J}_{\omega}$ and 
Goodman-Wallach operator ${\bf M}_{\omega}$ are both linear operators 
that commute with all elements of ${\cal U}({\frak g})$, and since 
all elements of ${\cal U}({\frak g})$ act as left-invariant 
differential operators on the space $C^{\infty}(G)$, it therefore 
follows from the definition (6.5.72) and the equation (6.5.29) 
(together with the remarks following it) that, when $\omega'\in{\frak O}$, we have: 
$${\bf X}\phi_{\omega'}\left(\nu_0,g\right) 
={\bf X}\phi^{\infty}_{\omega'}\left(\ell,q;\nu_0,g\right) 
=\sum_{(\lambda,\kappa)\in I(\ell,q)} 
c_{\ell,q}^{\bf X}\left(\lambda,\kappa;\nu_0,0\right) 
\phi^{\infty}_{\omega'}\left(\lambda,\kappa;\nu_0,g\right)\qquad\quad 
\hbox{($g\in G$),}\eqno(6.5.85)$$
where $I(\ell,q)$ is the finite subset of ${\Bbb Z}\times{\Bbb Z}$ 
defined in (6.5.23), while  each coefficient 
$c_{\ell,q}^{\bf X}(\lambda,\kappa;\nu_0,0)$ is a 
certain polynomial function of $\nu_0$ (and is of course 
independent of the variable $g$). 
By Lemma~6.5.14, it is moreover the case that, for each pair 
$(\lambda,\kappa)\in I(\ell,q)$, the series 
$\sum_{\omega'\in{\frak O}}\phi^{\infty}_{\omega'}(\lambda,\kappa;\nu_0,g)$ 
converges uniformly on any given compact subset of $G$; 
by this and the equation (6.5.85), we may deduce that 
the series 
$\sum_{\omega'\in{\frak O}}{\bf X}\phi^{\infty}_{\omega'}(\ell,q;\nu_0,g)$ 
is (likewise) uniformly convergent on any given compact subset of $G$. 
Moreover, by (6.5.85), (6.5.72) (with the $K$-type $(\lambda,\kappa)$ substituted 
for the $K$-type $(\ell,q)$), the relation (1.5.15), and the remark preceding (6.3.4), 
one sees that each term ${\bf X}\phi^{\infty}_{\omega'}(\ell,q;\nu_0,g)$ 
in the latter series is a continuous function of the variable $g$. 
It therefore follows by the proposition in Section~1.72 of [43] that for 
$g\in G$ the derivative 
${\bf X}\Phi^{\infty}(\ell,q;\nu_0,g)
={\bf X}(\sum_{\omega'\in{\frak O}}\phi^{\infty}_{\omega'}(\ell,q;\nu_0,g))$ 
exists, and is equal to 
$\sum_{\omega'\in{\frak O}}{\bf X}\phi^{\infty}_{\omega'}(\ell,q;\nu_0,g)$:  
note that, although we omit to give a detailed justification of the reasoning here, 
the details omitted are similar to details provided in our proof of  
the result (6.5.38) of Lemma~6.5.7. 
\par 
Now, by the definition (6.5.82), it follows  that 
${\bf X}(\sum_{\omega'\in{\frak O}}
\phi^{\infty}_{\omega'}(\ell,q;\nu_0,g)) 
=({\bf X}{\cal P}^{\frak a}_{\!\!\!\scriptscriptstyle\hookleftarrow}
{\bf M}_{\omega}\varphi_{\ell,q}(\nu_0,0))(g)\,$ for all $g\in G$. 
On the other hand, by (6.5.85), (6.5.23) and (6.5.82) 
(the last applied, this time, with the $K$-type $(\lambda,\kappa)$ substituted for 
$(\ell,q)$), we have, for $g\in G$, 
$$\eqalign{ 
\sum_{\omega'\in{\frak O}} {\bf X}\phi^{\infty}_{\omega'}\left(\ell,q;\nu_0,g\right) 
 &=\sum_{(\lambda,\kappa)\in I(\ell,q)} 
c_{\ell,q}^{\bf X}\left(\lambda,\kappa;\nu_0,0\right) 
\sum_{\omega'\in{\frak O}} 
\phi^{\infty}_{\omega'}\left(\lambda,\kappa;\nu_0,g\right) =\cr 
 &=\sum_{(\lambda,\kappa)\in I(\ell,q)} 
c_{\ell,q}^{\bf X}\left(\lambda,\kappa;\nu_0,0\right)  
\bigl({\cal P}^{\frak a}_{\!\!\!\scriptscriptstyle\hookleftarrow}
{\bf M}_{\omega}\varphi_{\lambda,\kappa}(\nu_0,0)\bigr)(g) \;,
}$$
where, by (6.5.82) and the result (6.5.75) of Lemma~6.5.14, the 
relevant functions 
$g\mapsto ({\cal P}^{\frak a}_{\!\!\!\scriptscriptstyle\hookleftarrow}
{\bf M}_{\omega}\varphi_{\lambda,\kappa}(\nu_0,0))(g)$
are each continuous on $G$. The convergence of the relevant series 
here is not at issue (see the previous paragraph). 
Given the identity 
${\bf X}(\sum_{\omega'\in{\frak O}}\phi^{\infty}_{\omega'}(\ell,q;\nu_0,g))  
=\sum_{\omega'\in{\frak O}}{\bf X}\phi^{\infty}_{\omega'}(\ell,q;\nu_0,g)\,$ 
obtained in the previous paragraph, and the points just noted 
in this paragraph, we find that 
the function $g\mapsto({\bf X}{\cal P}^{\frak a}_{\!\!\!\scriptscriptstyle\hookleftarrow}
{\bf M}_{\omega}\varphi_{\ell,q}(\nu_0,0))(g)\,$ is defined 
and continuous on $G$, and that one has:  
$${\bf X}{\cal P}^{\frak a}_{\!\!\!\scriptscriptstyle\hookleftarrow}
{\bf M}_{\omega}\varphi_{\ell,q}(\nu_0,0) 
=\sum_{(\lambda,\kappa)\in I(\ell,q)} 
c_{\ell,q}^{\bf X}\left(\lambda,\kappa;\nu_0,0\right)  
{\cal P}^{\frak a}_{\!\!\!\scriptscriptstyle\hookleftarrow}
{\bf M}_{\omega}\varphi_{\lambda,\kappa}(\nu_0,0)\;.\eqno(6.5.86)$$
\par 
Since the definition (6.5.23) implies that the set $I(\ell,q)$ is finite, 
and since ${\cal B}_1$ is a ${\Bbb C}$-basis of ${\frak g}$, 
it may be shown, through the iterative application of Equation~(6.5.86) 
(with the choice of ${\bf X}\in{\frak g}$ and 
$K$-type~$(\ell,q)$ varying at each iteration), 
that one has ${\cal P}^{\frak a}_{\!\!\!\scriptscriptstyle\hookleftarrow}
{\bf M}_{\omega}\varphi_{\ell,q}(\nu_0,0)\in\bigcap_{j=0}^{\infty} C^j(G)$,  
where each space $C^j(G)$ is as defined as below (6.5.49), in the proof of Lemma~6.5.8. 
Therefore, given the equality in (6.5.50), it follows that we have  
${\cal P}^{\frak a}_{\!\!\!\scriptscriptstyle\hookleftarrow}
{\bf M}_{\omega}\varphi_{\ell,q}(\nu_0,0)\in C^{\infty}(G)$;  
since we showed earlier that the function 
${\cal P}^{\frak a}_{\!\!\!\scriptscriptstyle\hookleftarrow}
{\bf M}_{\omega}\varphi_{\ell,q}(\nu_0,0) : G\rightarrow{\Bbb C}$ 
is $\Gamma$-automorphic, this completes our proof of 
the result that ${\cal P}^{\frak a}_{\!\!\!\scriptscriptstyle\hookleftarrow}
{\bf M}_{\omega}\varphi_{\ell,q}(\nu_0,0)\in C^{\infty}(\Gamma\backslash G)$.  
To complete the proof of the lemma we now have only to show 
that the function ${\cal P}^{\frak a}_{\!\!\!\scriptscriptstyle\hookleftarrow}
{\bf M}_{\omega}\varphi_{\ell,q}(\nu_0,0) : G\rightarrow{\Bbb C}$ 
is of $K$-type~$(\ell,q)$. 
\par 
By (6.5.86), the definition (6.5.82), 
the result (6.5.76) of Lemma~6.5.14 (applied with 
the $K$-type $(\lambda,\kappa)$ substituted for the $K$-type $(\ell,q)$) and 
the equation (6.5.29), we may infer that, for all ${\rm Re}(\nu)>1$, 
and all ${\bf Y}\in{\frak g}$, one has: 
$$\eqalign{ 
{\bf Y}{\cal P}^{\frak a}_{\!\!\!\scriptscriptstyle\hookleftarrow}
{\bf M}_{\omega}\varphi_{\ell,q}(\nu,0) 
 &=\sum_{(\lambda,\kappa)\in I(\ell,q)} 
c_{\ell,q}^{\bf Y}\left(\lambda,\kappa;\nu,0\right)  
P^{\frak a}{\bf M}_{\omega}\varphi_{\lambda,\kappa}(\nu,0) = \cr 
 &=P^{\frak a}\biggl( 
\,\sum_{(\lambda,\kappa)\in I(\ell,q)} 
c_{\ell,q}^{\bf Y}\left(\lambda,\kappa;\nu,0\right)  
{\bf M}_{\omega}\varphi_{\lambda,\kappa}(\nu,0)
\biggr) = \cr 
 &=P^{\frak a}{\bf M}_{\omega}\biggl( 
\,\sum_{(\lambda,\kappa)\in I(\ell,q)} 
c_{\ell,q}^{\bf Y}\left(\lambda,\kappa;\nu,0\right)  
\varphi_{\lambda,\kappa}(\nu,0)
\biggr)
=P^{\frak a} {\bf M}_{\omega}{\bf Y}\varphi_{\ell,q}(\nu,0)\;.
}$$
Since we have also ${\bf Y}P^{\frak a}{\bf M}_{\omega}\varphi_{\ell,q}(\nu,0) 
={\bf Y}{\cal P}^{\frak a}_{\!\!\!\scriptscriptstyle\hookleftarrow}
{\bf M}_{\omega}\varphi_{\ell,q}(\nu,0)\,$ (for ${\rm Re}(\nu)>1$ and 
${\bf Y}\in{\frak g}$), it follows from the above equations 
(similarly to how the results in (6.5.54) and (6.5.55) were obtained) 
that, when ${\rm Re}(\nu)>1$, one has 
${\bf H}_2 P^{\frak a}{\bf M}_{\omega}\varphi_{\ell,q}(\nu,0) 
=-iq P^{\frak a}{\bf M}_{\omega}\varphi_{\ell,q}(\nu,0)$ 
and $\Omega_{\frak k}P^{\frak a}{\bf M}_{\omega}\varphi_{\ell,q}(\nu,0) 
=-{1\over 2}(\ell^2+\ell) P^{\frak a}{\bf M}_{\omega}\varphi_{\ell,q}(\nu,0)$. 
By this, combined with the result (6.5.76) of Lemma~7.5.12,  
and the definition (6.5.82), it follows that 
$$\eqalignno{
{\bf H}_2 {\cal P}^{\frak a}_{\!\!\!\scriptscriptstyle\hookleftarrow} 
{\bf M}_{\omega}\varphi_{\ell,q}(\nu,0) 
 &=-iq {\cal P}^{\frak a}_{\!\!\!\scriptscriptstyle\hookleftarrow} 
{\bf M}_{\omega}\varphi_{\ell,q}(\nu,0)\qquad\quad  
\hbox{(${\rm Re}(\nu)>1$);} &(6.5.87) \cr  
\Omega_{\frak k}{\cal P}^{\frak a}_{\!\!\!\scriptscriptstyle\hookleftarrow} 
{\bf M}_{\omega}\varphi_{\ell,q}(\nu,0) 
 &=-{1\over 2}\left(\ell^2+\ell\right) 
{\cal P}^{\frak a}_{\!\!\!\scriptscriptstyle\hookleftarrow} 
{\bf M}_{\omega}\varphi_{\ell,q}(\nu,0)\qquad\quad
\hbox{(${\rm Re}(\nu)>1$).} &(6.5.88)}$$ 
Moreover, by iteration of the equation (6.5.86) (in the manner previously 
discussed) and application of the result (6.5.74) of Lemma~6.5.14, 
it may be shown that, when $g\in G$, $J\in{\Bbb N}$ and 
${\bf X}_1,\ldots ,{\bf X}_J\in{\frak g}$, the function 
$\nu\mapsto ({\bf X}_J {\bf X}_{J-1}\cdots {\bf X}_1 
{\cal P}^{\frak a}_{\!\!\!\scriptscriptstyle\hookleftarrow} 
{\bf M}_{\omega}\varphi_{\ell,q}(\nu,0))(g)$ is holomorphic 
on the half-plane where ${\rm Re}(\nu)>1/2$; and so one finds, in particular, 
that both sides of each equation in (6.5.87) and (6.5.88) are 
functions of $\nu$ that are holomorphic when ${\rm Re}(\nu)>1/2$. 
The equations in (6.5.87) and (6.5.88) must therefore be 
valid for all $\nu\in{\Bbb C}$ such that ${\rm Re}(\nu)>1/2$, so we 
have the required proof that 
${\cal P}^{\frak a}_{\!\!\!\scriptscriptstyle\hookleftarrow}
{\bf M}_{\omega}\varphi_{\ell,q}(\nu_0,0)$ is of $K$-type $(\ell,q)\quad\blacksquare$  

\bigskip 

For each $\nu\in{\Bbb C}$ such that ${\rm Re}(\nu)>1/2$, we define now the function 
${\cal P}^{\frak a}_{\scriptscriptstyle\hookleftarrow}\tau{\bf M}_{\omega}\varphi_{\ell,q}(\nu,0) 
: G\rightarrow{\Bbb C}$ by setting 
$${\cal P}^{\frak a}_{\!\!\!\scriptscriptstyle\hookleftarrow}\tau{\bf M}_{\omega}\varphi_{\ell,q}(\nu,0) 
={\cal P}^{\frak a}_{\!\!\!\scriptscriptstyle\hookleftarrow}{\bf M}_{\omega}\varphi_{\ell,q}(\nu,0) 
-P^{\frak a}(1-\tau){\bf M}_{\omega}\varphi_{\ell,q}(\nu,0)\eqno(6.5.89)$$
(it being a corollary of Lemma~6.5.2, the definition (6.5.82) and Lemma~6.5.14
that, when ${\rm Re}(\nu)>1/2$, the expression on the right-hand side of Equation~(6.5.89) denotes a 
well-defined complex-valued function on~$G$). 
In the next lemma we establish certain useful properties of the 
functions defined in (6.5.89): the lemma implies (amongst other things) 
the existence of the limit in the definition~(6.5.5). Note that 
in (6.5.94) we use `$L^p(\Gamma\backslash G;\ell,q)\,$' (for $p=2/\varepsilon$) to denote  
the space $\{ f\in L^p(\Gamma\backslash G) : f\ {\rm is\ of}\ K{\rm -type}\ (\ell,q)\}$.  

\bigskip 

\proclaim Lemma 6.5.16. Let $0\neq\omega\in{\frak O}$ and 
${\cal N}=\{\nu\in{\Bbb C} : {\rm Re}(\nu)>1/2\}$;  
for $\nu\in{\cal N}$, let 
${\cal P}^{\frak a}_{\!\!\!\scriptscriptstyle 
\hookleftarrow}\tau{\bf M}_{\omega}\varphi_{\ell,q}(\nu,0) : G\rightarrow{\Bbb C}$ 
be as defined in (6.5.89); for $(\nu,g)\in{\cal N}\times G$, put 
$\Phi_{\tau}^{\infty}(\nu,g)=
({\cal P}^{\frak a}_{\!\!\!\scriptscriptstyle 
\hookleftarrow}\tau{\bf M}_{\omega}\varphi_{\ell,q}(\nu,0))(g)$. 
Then the following are valid statements: 
$$\eqalignno{ 
 &\hbox{when $\,{\rm Re}(\nu)>1$, one has} 
\ \ {\cal P}^{\frak a}_{\!\!\!\scriptscriptstyle 
\hookleftarrow}\tau{\bf M}_{\omega}\varphi_{\ell,q}(\nu,0) 
=P^{\frak a}\tau{\bf M}_{\omega}\varphi_{\ell,q}(\nu,0)\;; &(6.5.90)\cr 
 &\hbox{when $g\in G$, the function $\nu\mapsto\Phi^{\infty}_{\tau}(\nu,g)$ is 
holomorphic on ${\cal N}$;} &(6.5.91)\cr 
 &\hbox{the function $(\nu,g)\mapsto\Phi^{\infty}_{\tau}(\nu,g)$ is continuous 
on ${\cal N}\times G$;} &(6.5.92) \cr 
 &\hbox{when $\,\nu\in{\cal N}$, the function 
${\cal P}^{\frak a}_{\!\!\!\scriptscriptstyle 
\hookleftarrow}\tau{\bf M}_{\omega}\varphi_{\ell,q}(\nu,0)$ 
lies in $C^{\infty}(\Gamma\backslash G)$ and is of $K$-type $(\ell,q)$;} &(6.5.93) \cr 
 &\hbox{when $|\nu -1|<\varepsilon <1/2$, one has 
$\,{\cal P}^{\frak a}_{\!\!\!\scriptscriptstyle 
\hookleftarrow}\tau{\bf M}_{\omega}\varphi_{\ell,q}(\nu,0)\in 
L^{2/\varepsilon}(\Gamma\backslash G;\ell,q)$;} &(6.5.94) \cr 
 &\hbox{when $\,{\rm Re}(\nu)\geq 1$,  one has 
$\,{\cal P}^{\frak a}_{\!\!\!\scriptscriptstyle 
\hookleftarrow}\tau{\bf M}_{\omega}\varphi_{\ell,q}(\nu,0)\in 
L^{\infty}(\Gamma\backslash G;\ell,q)$.} &(6.5.95)}$$ 
Moreover, when $t_1\in(0,\infty)$, $1/2<\sigma_1<\sigma_2<\infty$ and 
$${\cal R}={\cal R}\left(\sigma_1,\sigma_2,t_1\right) =\left\{\nu\in{\Bbb C} \,:\,  
\sigma_1\leq {\rm Re}(\nu)\leq\sigma_2\ {\rm and}\ |{\rm Im}(\nu)|\leq t_1\right\} ,$$  
one then has 
$$\Phi^{\infty}_{\tau}\left(\nu , g_{\frak b} g\right) 
\ll_{\Gamma,\omega,\ell,{\cal R}}\ \left(\rho(g)\right)^{1-{\rm Re}(\nu)}\qquad\ 
\hbox{for $\ (\nu,g)\in {\cal R}\times\{ g\in G : \rho(g)\geq 1/|q_0|\}$,}\eqno(6.5.96)$$
where $q_0$ denotes the `level' of the Hecke congruence 
subgroup $\Gamma\leq SL(2,{\frak O})$.  

\medskip 

\noindent{\bf Proof.}\quad 
Let $g_0\in G$; and let $\nu_0\in{\Bbb C}$ be such that ${\rm Re}(\nu_0)>1$. 
It is then a corollary of what was found in the proofs of 
Lemma~6.5.1 and Lemma~6.5.2 that the Poincar\'e series 
$(P^{\frak a}{\bf M}_{\omega}\varphi_{\ell,q}(\nu_0,0))(g_0)$  and 
$(P^{\frak a}(1-\tau){\bf M}_{\omega}\varphi_{\ell,q}(\nu_0,0))(g_0)$ 
are absolutely convergent. Hence, given the linearity (implicit in (1.5.4)) 
of the operator $P^{\frak a}$,  
it follows that the Poincar\'e series 
$(P^{\frak a}\tau{\bf M}_{\omega}\varphi_{\ell,q}(\nu_0,0))(g_0) 
=P^{\frak a}\bigl({\bf M}_{\omega}\varphi_{\ell,q}(\nu_0,0) 
-(1-\tau){\bf M}_{\omega}\varphi_{\ell,q}(\nu_0,0)\bigr)(g_0)$ 
is absolutely convergent, and that 
$(P^{\frak a}\tau{\bf M}_{\omega}\varphi_{\ell,q}(\nu_0,0))(g_0) 
=(P^{\frak a}{\bf M}_{\omega}\varphi_{\ell,q}(\nu_0,0))(g_0)  
-(P^{\frak a}(1-\tau){\bf M}_{\omega}\varphi_{\ell,q}(\nu_0,0))(g_0)$. 
By this, combined with the result 
(6.5.76) of Lemma~6.5.14, and the definitions (6.5.82), (6.5.89), 
we infer what is stated in (6.5.90). 
Moreover, given those definitions, (6.5.82) and (6.5.89), 
the results in (6.5.91), (6.5.92) and (6.5.93) are an immediate 
corollary of the combined results of Lemma~6.5.2, Lemma~6.5.14 
(the results (6.5.74), (6.5.75) in particular) and Lemma~6.5.15. 
The remainder of this proof is therefore devoted to 
the demonstration of what is asserted in (6.5.94)-(6.5.96). 
\par 
Henceforth let ${\cal R}(a,b,c)$ denote 
(when $a,b,c\in{\Bbb R}$) the set $\{\nu\in{\Bbb C} : a\leq{\rm Re}(\nu)\leq b\ {\rm and}\ 
|{\rm Im}(\nu)|\leq c\}$. Suppose that $t_1\in (0,\infty)$, that 
$1/2<\sigma_1 <\sigma_2 <\infty$, and that ${\cal R}={\cal R}(\sigma_1,\sigma_2,t_1)$. 
Then, by the definitions (6.5.89), (6.5.82) and the results 
(6.5.73) and (6.5.8) of Lemma~6.5.14 and Lemma~6.5.2, one has 
$$\Phi^{\infty}_{\tau}\left(\nu , g_{\frak b} g\right) 
\ll_{\Gamma,\omega,\ell,{\cal R}}\ \left(\rho(g)\right)^{1-{\rm Re}(\nu)}\qquad\ 
\hbox{for $\ (\nu,g)\in {\cal R}\times\{ g\in G : \rho(g)\geq 2\}$.}\eqno(6.5.97)$$ 
\par 
Observe now that $\Phi^{\infty}_{\tau}(\nu,g_{\frak b}n[\alpha]g)= 
\Phi^{\infty}_{\tau}(\nu,g_{\frak b}g)$ for $(\nu,g)\in{\cal N}\times G$ and 
all $\alpha\in{\frak O}\,$ (this follows, given our assumption that 
$g_{\frak b}^{-1}\Gamma_{\frak b}' g_{\frak b}  
=\{ n[\alpha] : \alpha\in{\frak O}\}$, by virtue of the result 
(6.5.93) proved earlier). Therefore, given the compactness 
of the subset 
${\cal R}\times\{ n[z]a[r]k : {\rm Re}(z),{\rm Im}(z)\in[-1/2,1/2], 
r\in[1/|q_0|,2]\ {\rm and}\ k\in K\}$ of ${\cal N}\times G$, 
and since we have (see (6.5.92)) already established  
the continuity of the function 
$(\nu,g)\mapsto\Phi^{\infty}_{\tau}(\nu, g)$, 
it follows that the function 
$(\nu,g)\mapsto\Phi^{\infty}_{\tau}(\nu, g_{\frak b}g)$ is 
bounded on the set 
${\cal R}\times\{ g\in G : 1/|q_0|\leq\rho(g)\leq 2\}$. 
Hence, by considering all the factors upon which that  
upper bound may depend,  we find that 
$$\Phi^{\infty}_{\tau}\left(\nu , g_{\frak b} g\right) 
\ll_{\Gamma,\omega,\ell,{\cal R}}\ \left(\rho(g)\right)^{1-{\rm Re}(\nu)}\qquad\ 
\hbox{for $\ (\nu,g)\in {\cal R}\times\{ g\in G : 1/|q_0|\leq\rho(g)\leq 2\}$.}
\eqno(6.5.98)$$
This last result merits some further explanation, since it would 
appear (at first sight) that the relevant implicit constant might have to depend 
on $g_{\frak a}$, $g_{\frak b}$ and $\tau$. In fact the question of 
the dependence of this implicit constant on $\tau$ does not arise, since 
our choice of $\tau$ is fixed (i.e. it is clear from (6.5.1), and the 
discussion preceding it, that we may define our fixed choice of 
the function $\tau$ in purely absolute terms). 
\par 
In considering whether or not the implicit constant  in (6.5.98) has 
to depend on the scaling matrix $g_{\frak b}$, we note firstly that, if 
${\frak b}'\sim^{\!\!\!\!\Gamma}{\frak b}$, and if $g_{\frak b},g_{{\frak b}'}\in G$ 
are chosen so that (1.1.16) and (1.1.20) are satisfied for ${\frak c}\in\{ {\frak b}, 
{\frak b}'\}$, then, since $g\mapsto\Phi^{\infty}_{\tau}(\nu,g)$ is 
$\Gamma$-automorphic (when $\nu\in{\cal N}$), it may be deduced from Lemma~4.1 
that when $M$ and $M'$ denote the global maxima attained on the set 
${\cal R}\times\{ g\in G : 1/|q_0|\leq\rho(g)\leq 2\}$ by the functions 
$(\nu,g)\mapsto |\Phi^{\infty}_{\tau}(\nu,g_{\frak b} g)|$ and 
$(\nu,g)\mapsto |\Phi^{\infty}_{\tau}(\nu,g_{\frak b}' g)|\,$ 
(respectively) one will have $M'=M$. Therefore the implicit constant 
in (6.5.98) depends on $g_{\frak b}$ only to the extent that 
it depends on the $\Gamma$-equivalence class of the cusp ${\frak b}$. 
A similar phenomenon may be observed in respect of the 
dependence of the same constant upon the scaling matrix $g_{\frak a}$: 
for it follows by Lemma~4.1 and the equations (1.5.4), (6.3.3) and (6.3.5) that if 
${\frak a}'\sim^{\!\!\!\!\Gamma}{\frak a}$, and if 
$g_{\frak a},g_{{\frak a}'}\in G$ 
are chosen so that (1.1.16) and (1.1.20) are satisfied for ${\frak c}\in\{ {\frak a}, 
{\frak a}'\}$, then there exists some $\beta\in{\Bbb C}$ and some  
$\epsilon\in{\frak O}^{*}$ such that one has 
$P^{\frak a}\tau{\bf M}_{\omega}\varphi_{\ell,q}(\nu,0)
={\rm e}(-{\rm Re}(\omega\beta)) 
P^{{\frak a}'}\tau{\bf M}_{\epsilon\omega}\varphi_{\ell,q}(\nu,0)$ 
for ${\rm Re}(\nu)>1$; and this, by (6.5.90), implies that one must in fact have 
${\cal P}^{\frak a}_{\!\!\!\scriptscriptstyle 
\hookleftarrow}\tau{\bf M}_{\omega}\varphi_{\ell,q}(\nu,0) 
={\rm e}(-{\rm Re}(\omega\beta)) 
{\cal P}^{{\frak a}'}_{\!\!\!\scriptscriptstyle 
\hookleftarrow}\tau{\bf M}_{\epsilon\omega}\varphi_{\ell,q}(\nu,0)$ 
for all $\nu\in{\cal N}$. 
Therefore, by using a relation of the form 
$$c_0\left( \Gamma{\frak a} , \Gamma{\frak b}\right) 
\leq\max_{\left( {\frak a}' , {\frak b}'\right)\in 
{\frak C}(\Gamma)\times{\frak C}(\Gamma)} 
c_0\left( \Gamma{\frak a}' , \Gamma{\frak b}'\right)$$ 
(in which ${\frak C}(\Gamma)$ denotes any complete set of 
representatives of the $\Gamma$-equivalence classes of cusps, so that by Lemma~2.2 
one has $|{\frak C}(\Gamma)|<\infty$), we are able to  
obtain (6.5.98) with an implicit 
constant $c_0^{*}$ independent of $g_{\frak a}$, $g_{\frak b}$, ${\frak a}$, ${\frak b}$ 
and the $\Gamma$-equivalence classes $\Gamma{\frak a}$, $\Gamma{\frak b}$ 
(except inasmuch as $c_0^{*}$ may depend upon the group $\Gamma$). 
\par 
By (6.5.98) and (6.5.97) we obtain the result stated in (6.5.96). In order to complete 
the proof of the lemma, we show next that the bound (6.5.96) implies 
the results stated in (6.5.94) and (6.5.95). 
\par 
We embark firstly upon the deduction of (6.5.95). 
By reasoning similar to (but simpler than) that which precedes our  
statement of the bound (6.5.98), we deduce from the results in (6.5.93) that,  
when $\nu\in{\cal N}$, 
the function 
$g\mapsto\Phi^{\infty}_{\tau}(\nu,g) 
=({\cal P}^{\frak a}_{\!\!\!\scriptscriptstyle 
\hookleftarrow}\tau{\bf M}_{\omega}\varphi_{\ell,q}(\nu,0))(g)$ is 
$\Gamma$-automorphic (on $G$), and is bounded on each 
compact set ${\cal D}$ contained in $G$.  
Hence, since 
there exists (in respect of the action of $\Gamma$ on the upper half-space ${\Bbb H}_3$)  
a fundamental domain  
${\cal F}_{*}$ fitting the description given in (1.1.22)-(1.1.24), 
and since, for such a domain ${\cal F}_{*}$, the 
set $\bigcup_{(z,r)\in{\cal F}_{*}}n[z]a[r]K$ will contain a fundamental 
domain for the action of $\Gamma$ on the group $G$, 
it follows by (1.1.24), (1.1.23) and the case ${\cal R}={\cal R}(1,{\rm Re}(\nu)+1,|{\rm Im}(\nu)|+1)$ of the 
result (6.5.96) that, when ${\rm Re}(\nu)\geq 1$, the $\Gamma$-automorphic 
function $g\mapsto 
({\cal P}^{\frak a}_{\!\!\!\scriptscriptstyle 
\hookleftarrow}\tau{\bf M}_{\omega}\varphi_{\ell,q}(\nu,0))(g)$ 
is bounded on $G\,$ (we use here the fact that, by (6.1.25), one has 
$1/|m_{\frak c}|\geq 1/|q_0|$ for all cusps 
${\frak c}\in{\Bbb Q}(i)\cup\{\infty\}$). 
The result (6.5.95) therefore follows (given the content of 
the result (6.5.93) proved earlier). 
\par 
We now have only to show that (6.5.94) holds. 
In light of the result (6.5.93) obtained earlier, it will 
suffice that we show that 
$$\int_{\Gamma\backslash G} 
\left| 
\bigl({\cal P}^{\frak a}_{\!\!\!\scriptscriptstyle 
\hookleftarrow}\tau{\bf M}_{\omega}\varphi_{\ell,q}(\nu,0)\bigr)(g)
\right|^{2/\varepsilon} {\rm d}g <\infty\qquad\ 
\hbox{when $\ |\nu -1|<\varepsilon <1/2$.}\eqno(6.5.99)$$
\par 
Let $0<\varepsilon <1/2$. By the case 
${\cal R}={\cal R}(1-\varepsilon , 1+\varepsilon , \epsilon)$ of the result (6.5.96), 
it follows that, for $|\nu -1|<\varepsilon$ and $g\in G$ such that 
$\rho(g)\geq 1/|q_0|$, one has 
$$\bigl({\cal P}^{\frak a}_{\!\!\!\scriptscriptstyle 
\hookleftarrow}\tau{\bf M}_{\omega}\varphi_{\ell,q}(\nu,0)\bigr)
\bigl( g_{\frak b} g\bigr) 
\ll_{\Gamma,\omega,\ell,\varepsilon}\ 
\left(\rho(g)\right)^{1-{\rm Re}(\nu)} 
\leq\left| q_0\right|^{2\varepsilon} 
\left(\rho(g)\right)^{|1-\nu|}\;.$$
Hence when $|\nu -1|<\varepsilon$ we have, in particular, 
$$\left| 
\bigl({\cal P}^{\frak a}_{\!\!\!\scriptscriptstyle 
\hookleftarrow}\tau{\bf M}_{\omega}\varphi_{\ell,q}(\nu,0)\bigr)
\bigl( g_{\frak b} g\bigr)
\right|^{2/\varepsilon} 
\ll_{\Gamma,\omega,\ell,\varepsilon}\ 
\left(\rho(g)\right)^{2|\nu -1|/\varepsilon}\quad\ 
\hbox{for all $\,g\in G\,$ such that $\,\rho(g)\geq 1/|m_{\frak b}|$.}\eqno(6.5.100)$$
Since the hypothesis that $|\nu -1|<\varepsilon$ 
implies that, in (6.5.100), the final exponent $2|\nu -1|/\varepsilon$ 
is strictly less than $2$, it may therefore be deduced, 
by reasoning similar to that seen in the proof 
of Corollary~6.2.10, that the bound (6.5.100) 
is sufficient to justify what is asserted in (6.5.99)\quad$\blacksquare$ 

\bigskip 

\proclaim Lemma 6.5.17. When $g\in G$ the limit on the right-hand side of 
Equation~(6.5.5) exists. The function $P^{{\frak a},*}\widetilde{\bf L}^{\omega}_{\ell,q}\eta 
: G\rightarrow{\Bbb C}$ defined by Equation~(6.5.5) is equal to the 
function $P^{\frak a}\widetilde{\bf L}^{\omega,*}_{\ell,q}\eta 
+b(\omega;\ell,q;\eta) {\cal P}^{\frak a}_{\!\!\!\scriptscriptstyle\hookleftarrow}\tau 
{\bf M}_{\omega}\varphi_{\ell,q}(1,0)$, and lies in 
$C^{\infty}(G)\cap L^{\infty}(\Gamma\backslash G;\ell,q)$. 

\medskip 

\noindent{\bf Proof.}\quad 
By the definition (1.5.4), the results of Lemma~6.5.1 and Lemma~6.5.8 concerning 
$P^{\frak a}|\tau{\bf M}_{\omega}\varphi_{\ell,q}(\nu,0)|$ and 
$P^{\frak a}|\widetilde{\bf L}^{\omega,*}_{\ell,q}\eta|$, and the result (6.5.90) 
of Lemma~6.5.16, it follows that when ${\rm Re}(\nu)>1$ one has 
$$\eqalign{
P^{\frak a}\bigl(\widetilde{\bf L}^{\omega,*}_{\ell,q}\eta 
+b(\omega;\ell,q;\eta)\tau{\bf M}_{\omega}\varphi_{\ell,q}(\nu,0)\bigr) 
 &=P^{\frak a}\widetilde{\bf L}^{\omega,*}_{\ell,q}\eta 
+b(\omega;\ell,q;\eta) P^{\frak a}\tau{\bf M}_{\omega}\varphi_{\ell,q}(\nu,0) = \cr 
 &=P^{\frak a}\widetilde{\bf L}^{\omega,*}_{\ell,q}\eta 
+b(\omega;\ell,q;\eta){\cal P}^{\frak a}_{\!\!\!\scriptscriptstyle\hookleftarrow}\tau 
{\bf M}_{\omega}\varphi_{\ell,q}(\nu,0)\;.}$$
By this observation, and the result (6.5.91) of Lemma~6.5.16, we find that 
when $g\in G$ the limit on the right-hand side of Equation~(6.5.5) exists, and 
is equal to  
$(P^{\frak a}\widetilde{\bf L}^{\omega,*}_{\ell,q}\eta)(g)  
+b(\omega;\ell,q;\eta) ({\cal P}^{\frak a}_{\!\!\!\scriptscriptstyle\hookleftarrow}\tau 
{\bf M}_{\omega}\varphi_{\ell,q}(1,0))(g)$. 
This completes the proof of the first two assertions of the lemma;   
since it has, in particular, been shown that   
$P^{{\frak a},*}\widetilde{\bf L}^{\omega}_{\ell,q}\eta 
=P^{\frak a}\widetilde{\bf L}^{\omega,*}_{\ell,q}\eta 
+b(\omega;\ell,q;\eta) {\cal P}^{\frak a}_{\!\!\!\scriptscriptstyle\hookleftarrow}\tau 
{\bf M}_{\omega}\varphi_{\ell,q}(1,0)$, 
the result that $P^{{\frak a},*}\widetilde{\bf L}^{\omega}_{\ell,q}\eta \in 
C^{\infty}(G)\cap L^{\infty}(\Gamma\backslash G;\ell,q)$ 
follows by virtue of Lemma~6.5.8 and the results (6.5.93) and (6.5.95) of Lemma~6.5.16\quad$\blacksquare$ 

\bigskip 

\goodbreak\centerline{\bf \S 6.6 The preliminary spectral summation formula.} 

\medskip 

Throughout this subsection we suppose that ${\frak a}$, ${\frak b}$, 
$g_{\frak a}$, $g_{\frak b}$, $\sigma$ and the $K$-type $(\ell,q)$ are 
as stated at the beginning of Subsection~6.5: in particular, we suppose that 
$1<\sigma <2$. We assume, moreover, 
that $\omega_1$ and $\omega_2$ are non-zero Gaussian integers, and 
that $\eta$ and $\theta$ are functions that lie in 
${\cal T}^{\ell}_{\sigma}\,$ (the space defined in and below (6.4.3)). 
Subject to these hypotheses, we seek to establish the following result. 

\bigskip 

\proclaim Proposition 6.6.1 (preliminary sum formula). Let 
$$\phi_1=P^{{\frak a},*}\widetilde{\bf L}^{\omega_1}_{\ell,q}\eta\qquad\quad 
{\rm and}\qquad\quad 
\phi_2=P^{{\frak b},*}\widetilde{\bf L}^{\omega_2}_{\ell,q}\theta\;.\eqno(6.6.1)$$
Then 
$$\eqalign{ 
 &\sum_{V}
\,\overline{C_V^{\frak a}\left(\omega_1;\nu_V,p_V\right)}
\,C_V^{\frak b}\left(\omega_2;\nu_V,p_V\right) 
h_{\ell}\left(\nu_V , p_V\right)\,+ \cr 
 &\qquad +\sum_{{\frak c}\in{\frak C}(\Gamma)}
{1\over 4\pi i\left[\Gamma_{\frak c} : \Gamma_{\frak c}'\right]} 
\sum_{p\in{1\over 2}\left[\Gamma_{\frak c} : 
\Gamma_{\frak c}'\right]{\Bbb Z}} 
\ \int\limits_{(0)} \overline{B_{\frak c}^{\frak a}\left(\omega_1;\nu,p\right)}\,
B_{\frak c}^{\frak b}\left(\omega_2;\nu,p\right) h_{\ell}(\nu,p)\,{\rm d}\nu = \cr 
 &\qquad\qquad       
=\;{\left[\Gamma_{\frak a}:\Gamma_{\frak a}'\right]
\left[\Gamma_{\frak b}:\Gamma_{\frak b}'\right]\over 4\pi^2}
\,\left\langle\phi_1 , 
\phi_2\right\rangle_{\Gamma\backslash G} = \cr 
 &\qquad\qquad  
=\;{1\over 4\pi^3 i}\,\delta^{{\frak a},{\frak b}}_{\omega_1,\omega_2}\,
\sum_{p\in{\Bbb Z}}\ \int\limits_{(0)}^{\hbox{\ }}h_{\ell}(\nu,p) 
\left( p^2 -\nu^2\right){\rm d}\nu
+\,\sum_{c\in {}^{\frak a}{\cal C}^{\frak b}} 
\,{S_{{\frak a},{\frak b}}\left(\omega_1 , \omega_2 ; c\right)\over |c|^2}\,
\left({\bf B}h_{\ell}\right)
\!\left( {2\pi\sqrt{\omega_1\omega_2}\over c}\right) , 
}\eqno(6.6.2)$$
where, for $(\nu,p)\in{\Bbb C}\times{\Bbb Z}$ such that 
$|{\rm Re}(\nu)|\leq\sigma$, one has  
$$h_{\ell}(\nu,p) 
=\cases{\lambda^{*}_{\ell}(\nu,p)\,\overline{\theta\left( -\overline{\nu}\,, p\right)}
\,\eta(\nu,p) 
 &if $|p|\leq\ell$, \cr 0 &otherwise,}\eqno(6.6.3)$$ 
with 
$$\lambda^{*}_{\ell}(\nu,p) 
=\Gamma(\ell +1+\nu)\Gamma(\ell +1-\nu)\, 
{\sin^2(\pi\nu)\over (\pi\nu)^2}\,{\nu^{2+2\epsilon(p)}\over\left(\nu^2-p^2\right)^2} 
={1\over \Gamma(\ell +1+\nu)\Gamma(\ell +1-\nu)}
\prod_{\scriptstyle 0<m\leq\ell\atop\scriptstyle m\neq |p|} 
\!\left( \nu^2 -m^2\right)^2 \eqno(6.6.4)$$ 
($\epsilon(p)$ being as in (6.4.5)), while the term  
$\delta^{{\frak a},{\frak b}}_{\omega_1,\omega_2}$ and   
${\bf B}$-transform are as defined in (1.9.2)-(1.9.6),  
and all other nonstandard notation has the meaning assigned to it 
in Subsection~1.1, Subsection~1.5, Subsection~1.7, Subsection~1.8 and (1.2.2). 
The sums and integrals occurring 
in Equation~(6.6.2) are absolutely convergent.   

\bigskip 

 We assume henceforth (in this subsection) that $\phi_1$ and $\phi_2$ are 
as stated in (6.6.1) (the relevant terminology having been defined in the equation (6.5.5)). 
By Lemma~6.5.17 and the observation (6.5.7), we have 
$L^2(\Gamma\backslash G)\supseteq L^{\infty}(\Gamma\backslash G)\ni 
\phi_j$ for $j=1,2$. Therefore it follows by the Cauchy-Schwarz 
inequality of Section~12.41 of [43] that the inner product 
$\langle\phi_1 , \phi_2\rangle_{\Gamma\backslash G}$ exists. 
\par 
Our proof of the above proposition occupies the 
remainder of this subsection. It is modelled on Bruggeman and Motohashi's 
original proof (in Section~9 of [5]) of the special case 
$\Gamma =SL(2,{\frak O})$, ${\frak a}={\frak b}=\infty$ of (6.6.2).  
In particular, we `compute' the value of the inner product 
$\langle\phi_1 , \phi_2\rangle_{\Gamma\backslash G}$ in two ways; 
one computation yielding the `geometric description' of 
$\langle\phi_1 , \phi_2\rangle_{\Gamma\backslash G}$ implied by the final 
equality in (6.6.2); the other supplying the `spectral description' (of the same 
quantity) implied by the first equality of (6.6.2). 
We remark that these `geometric' and `spectral' descriptions 
depend on the $K$-type $(\ell,q)$ only insofar as they depend on the parameter $\ell$: 
we shall in fact only ever need to apply the `preliminary sum formula' (6.6.2)   
in respect of cases where the $K$-type is $(\ell,0)$. 
\par 
The next lemma (derived from Lemma~6.2.4) is of key importance in 
enabling both the computation of the geometric description of 
$\langle\phi_1 , \phi_2\rangle_{\Gamma\backslash G}\,$ and the 
computation of the corresponding spectral description. 
Before stating the lemma we clarify that henceforth $L^p(N\backslash G)$ 
will denote (when $1\leq p<\infty$) the space of those measurable functions 
$f : G\rightarrow{\Bbb C}$ that satisfy $f(ng)=f(g)$, for all $n\in N$, $g\in G$, 
and are such that 
$\int_{N\backslash G} |f(g)|^p\,{\rm d}\dot g<\infty\,$ 
(where the measure ${\rm d}\dot g$ is that which occurs in (6.2.9)).  

\bigskip 

\proclaim Lemma 6.6.2. Let $0\neq\omega\in{\frak O}$, let 
$\alpha , \beta\in(0,\infty)$, and let 
$\phi\in C^0(G)\cap L^{1+\alpha}(\Gamma\backslash G)$.  
Suppose moreover that one has $|F^{\frak a}_{\omega}\phi|\in L^{1+\beta}(N\backslash G)$. 
Then 
$$\left[\Gamma_{\frak a} : \Gamma_{\frak a}'\right] 
\bigl\langle P^{{\frak a},*} \widetilde{\bf L}^{\omega}_{\ell,q}\eta\,, 
\phi\bigr\rangle_{\Gamma\backslash G} 
=\bigl\langle \widetilde{\bf L}^{\omega}_{\ell,q}\eta\,, 
F^{\frak a}_{\omega}\phi\bigr\rangle_{N\backslash G}\;.
\eqno(6.6.5)$$ 

\medskip 

\noindent{\bf Proof.}\quad 
It will suffice to show that one has both 
$$\left[\Gamma_{\frak a} : \Gamma_{\frak a}'\right] 
\bigl\langle P^{\frak a} \widetilde{\bf L}^{\omega,*}_{\ell,q}\eta\,, 
\phi\bigr\rangle_{\Gamma\backslash G} 
=\bigl\langle \widetilde{\bf L}^{\omega,*}_{\ell,q}\eta\,, 
F^{\frak a}_{\omega}\phi\bigr\rangle_{N\backslash G}\eqno(6.6.6)$$
and 
$$\left[\Gamma_{\frak a} : \Gamma_{\frak a}'\right] 
\bigl\langle 
{\cal P}^{\frak a}_{\!\!\!\scriptscriptstyle\hookleftarrow}\tau 
{\bf M}_{\omega}\varphi_{\ell,q}(1,0)\,, 
\phi\bigr\rangle_{\Gamma\backslash G} 
=\bigl\langle \tau{\bf M}_{\omega}\varphi_{\ell,q}(1,0)\,, 
F^{\frak a}_{\omega}\phi\bigr\rangle_{N\backslash G}\;.\eqno(6.6.7)$$ 
Indeed, the equalities asserted in (6.6.6) and (6.6.7) imply the equality   
$$\left[\Gamma_{\frak a} : \Gamma_{\frak a}'\right] 
\bigl\langle P^{\frak a} \widetilde{\bf L}^{\omega,*}_{\ell,q}\eta 
+b(\omega;\ell,q;\eta) {\cal P}^{\frak a}_{\!\!\!\scriptscriptstyle\hookleftarrow}\tau 
{\bf M}_{\omega}\varphi_{\ell,q}(1,0)\,, 
\phi\bigr\rangle_{\Gamma\backslash G} 
=\bigl\langle \widetilde{\bf L}^{\omega,*}_{\ell,q}\eta 
+b(\omega;\ell,q;\eta) \tau{\bf M}_{\omega}\varphi_{\ell,q}(1,0)\,, 
F^{\frak a}_{\omega}\phi\bigr\rangle_{N\backslash G}\;,$$
which, by Lemma~6.5.17 and the definition (6.5.2), is equivalent 
to Equation (6.6.5). 
\par 
The equality in (6.6.6) is just the case $f_{\omega}=
\widetilde{\bf L}^{\omega,*}_{\ell,q}\eta$ of the result (6.2.8) of Lemma~6.2.4. 
It therefore suffices for proof of (6.6.6) that we verify that 
the hypotheses of Lemma~6.2.4 are satisfied when 
$f_{\omega}=
\widetilde{\bf L}^{\omega,*}_{\ell,q}\eta$ and $\phi$ is as we suppose (in this proof). 
This (since the relation 
$\phi\in C^0(\Gamma\backslash G)$ is implicit in our current hypotheses) merely 
entails our showing that $(P^{\frak a}|\widetilde{\bf L}^{\omega,*}_{\ell,q}\eta|)\cdot\phi
\in L^1(\Gamma\backslash G)$, that 
$\widetilde{\bf L}^{\omega,*}_{\ell,q}\eta\in C^0(N\backslash G,\omega)$, and that, for 
some $\sigma_0>1$ and some $R_0>0$, one has 
$$\bigl(\widetilde{\bf L}^{\omega,*}_{\ell,q}\eta\bigr)(g) 
\ll_{\omega,\eta,\sigma_0,R_0}\,\left(\rho(g)\right)^{1+\sigma_0}\qquad\ 
\hbox{for all $g\in G$ such that $\rho(g)\leq R_0$.}\eqno(6.6.8)$$
Since we assume (throughout this subsection) that 
$\eta\in{\cal T}^{\ell}_{\sigma}$, and that $\sigma\in (1,2)$, 
it follows by the observation recorded in (6.5.22) that, if one puts 
$R_0=1$ and $\sigma_0=\sigma >1$, then one does have $R_0>0$, $\sigma_0 >1$ and 
the desired growth estimate (6.6.8). By~(6.5.4), we have also that 
$\widetilde{\bf L}^{\omega,*}_{\ell,q}\eta\in C^0(N\backslash G,\omega)$.  
Therefore the equality (6.6.6) follows if 
$(P^{\frak a}|\widetilde{\bf L}^{\omega,*}_{\ell,q}\eta|)\cdot\phi
\in L^1(\Gamma\backslash G)$. By H\"{o}lder's 
inequality (as formulated in Section~12.42 of [43]), 
we find that the last condition on  $(P^{\frak a}|\widetilde{\bf L}^{\omega,*}_{\ell,q}\eta|)\cdot\phi$   
is satisfied: for we 
have $L^{1+\alpha}(\Gamma\backslash G)\ni\phi$, by hypothesis, and it follows 
by Lemma~6.5.8 and the observation (6.5.7) that 
$L^{1+1/\alpha}(\Gamma\backslash G)\supseteq L^{\infty}(\Gamma\backslash G)
\ni P^{\frak a}|\widetilde{\bf L}^{\omega,*}_{\ell,q}\eta|$.
This completes the proof of (6.6.6). 
\par 
In what follows we take $\Phi_1(\nu)$ and $\Phi_2(\nu)$ to denote the 
inner products $\bigl\langle 
{\cal P}^{\frak a}_{\!\!\!\scriptscriptstyle\hookleftarrow}\tau 
{\bf M}_{\omega}\varphi_{\ell,q}(\nu,0)\,, 
\phi\bigr\rangle_{\Gamma\backslash G}$ and  
$\bigl\langle \tau{\bf M}_{\omega}\varphi_{\ell,q}(\nu,0)\,, 
F^{\frak a}_{\omega}\phi\bigr\rangle_{N\backslash G}$, respectively 
(so that each of  $\Phi_1(\nu)$, $\Phi_2(\nu)$ is defined only for those 
$\nu\in{\Bbb C}$ such that the relevant inner product exists). 
We show next that 
$$\left[\Gamma_{\frak a} : \Gamma_{\frak a}'\right]\Phi_1(\nu)=\Phi_2(\nu)\qquad\ 
\hbox{for all $\nu\in{\Bbb C}$ such that ${\rm Re}(\nu)>1$.}\eqno(6.6.9)$$ 

By the result (6.5.90) of Lemma~6.5.16, 
we have $\Phi_1(\nu)=\langle 
P^{\frak a}\tau {\bf M}_{\omega}\varphi_{\ell,q}(\nu,0)\,, 
\phi\rangle_{\Gamma\backslash G}$ when ${\rm Re}(\nu)>1$. This 
suggests that we might prove (6.6.9) 
by another application of Lemma~6.2.4, similar 
to that (above) by which the equality (6.6.6) was obtained (though   
with $f_{\omega}=\tau{\bf M}_{\omega}\varphi_{\ell,q}(\nu,0)$ on this occasion). 
As previously, we note that our hypotheses imply that $\phi\in C^0(\Gamma\backslash G)$.   
Furthermore, by (6.5.1), (6.5.3) and the result (6.3.7) of Lemma~6.3.1, we have 
$(\tau{\bf M}_{\omega}\varphi_{\ell,q}(\nu,0))(g) 
\ll_{\omega,\ell,\nu}\,(\rho(g))^{1+{\rm Re}(\nu)}$ 
for $\nu\in{\Bbb C}$, $g\in G$, and, 
given our choice of $\tau\in C^{\infty}(G)\,$ (as in (6.5.1)) and  
the observation preceding (6.3.4),  we have also 
$\tau{\bf M}_{\omega}\varphi_{\ell,q}(\nu,0)\in C^{\infty}(N\backslash G,\omega)$ 
for all $\nu\in{\Bbb C}$. Therefore, if it is the case that 
$(P^{\frak a}|\tau{\bf M}_{\omega}\varphi_{\ell,q}(\nu,0)|)\cdot\phi
\in L^1(\Gamma\backslash G)$ when ${\rm Re}(\nu)>1$, then Lemma~6.2.4 applies, 
giving (6.6.9). By H\"{o}lder's inequality and the hypothesis that 
$\phi\in L^{1+\alpha}(\Gamma\backslash G)$, we find 
that $(P^{\frak a}|\tau{\bf M}_{\omega}\varphi_{\ell,q}(\nu,0)|)\cdot\phi
\in L^1(\Gamma\backslash G)$ when  
$P^{\frak a}|\tau{\bf M}_{\omega}\varphi_{\ell,q}(\nu,0)| 
\in L^{1+1/\alpha}(\Gamma\backslash G)$; since Lemma~6.5.1 and the observation 
(6.5.7) imply that one has 
$L^{1+1/\alpha}(\Gamma\backslash G)\supseteq 
L^{\infty}(\Gamma\backslash G)
\ni P^{\frak a}|\tau{\bf M}_{\omega}\varphi_{\ell,q}(\nu,0)|$ 
when  ${\rm Re}(\nu)>1$, this 
completes the proof of (6.6.9). 
\par 
We complete the proof of the lemma by showing that 
$$\Phi_j(1)=\lim_{\nu\rightarrow 1+}\Phi_j(\nu)\qquad\ \hbox{for $j=1,2$.}\eqno(6.6.10)$$
This suffices for completion of the proof, since the combination of 
(6.6.9) and (6.6.10) implies the equation 
$\left[\Gamma_{\frak a} : \Gamma_{\frak a}'\right]\Phi_1(1)=\Phi_2(1)$, which 
is (6.6.7). The approach that we take to our proof of (6.6.10) is to 
establish the stronger result that each of the functions 
$\Phi_1(\nu),\Phi_2(\nu)$ is continuous in some neighbourhood of 
the point $\nu =1$. Before we proceed, note that by the relations in (6.5.7) there is no  
loss of generality in assuming henceforth that $0<\alpha < 1/3\,$ (this  
helps to simplify some of the calculations below). 
\par 
We prove first the case $j=1$ of (6.6.10). In doing so we assume 
(as we may) that the sets ${\frak C}(\Gamma)$ and~${\cal D}$, and the family 
$({\cal E}_{\frak c})_{{\frak c}\in{\frak C}(\Gamma)}$, 
are each as described in the final paragraph of Subsection~1.1;  
consequently ${\frak C}(\Gamma)$ is (by Lemma~2.2) a finite 
set of cusps, and the set ${\cal F}_{*}$ given by Equation~(1.1.24) 
is a fundamental domain for the action of $\Gamma$ on ${\Bbb H}_3$.  
By reasoning similar to that already seen 
in the first paragraph of the proof of Corollary~6.2.10, 
we find that the case $j=1$ of (6.6.10) certainly follows if, when 
${\cal X}\in\{{\cal D}\}\cup\{ {\cal E}_{\frak c} : {\frak c}\in{\frak C}(\Gamma)\}$, 
the complex function 
$$\nu\mapsto{1\over 2}\int_{\cal X}\int_K  
\bigl({\cal P}^{\frak a}_{\!\!\!\scriptscriptstyle\hookleftarrow}\tau 
{\bf M}_{\omega}\varphi_{\ell,q}(\nu,0)\bigr)(n[z]a[r]k) 
\,\overline{\phi(n[z]a[r]k)}\,{\rm d}k\,r^{-3} {\rm d}_{+}z{\rm d}r\eqno(6.6.11)$$ 
is continuous in some neighbourhood of the point $\nu =1$. 
\par 
Note that, since $\phi\in C^0(\Gamma\backslash G)$,  
the result (6.5.92) of Lemma~6.5.16 implies that the integrand in (6.6.11) 
is continuous, as a function of $(\nu,g)$, on  
the set $\{\nu\in{\Bbb C} : {\rm Re}(\nu)>1/2\}\times G$. Therefore,   
by a proof similar in principle to that in 
Section~1.52 of [43] (on `The continuity theorem'), it follows that if the 
set $\tilde{\cal X}=\{ n[z]a[r]k : k\in K ,\,(z,r)\in{\cal X}\}$ 
is a compact measurable subset of $G=SL(2,{\Bbb C})$ then the function 
(6.6.11) is continuous in the open half-plane where ${\rm Re}(\nu)>1/2$. 
Hence, in the particular case where 
${\cal X}={\cal D}$ (a compact hyperbolic polyhedron, with finitely 
many faces), it follows by virtue of the compactness 
of $K=SU(2)$ that the function (6.6.11) is certainly 
continuous in the neighbourhood  
$\{\nu\in{\Bbb C} : |\nu -1|<1/2\}$ of the point $\nu =1$. 
\par 
Suppose that, for some ${\frak c}\in{\frak C}(\Gamma)$, 
we have ${\cal X}={\cal E}_{\frak c}\,$ (so that, in (6.6.11), `${\cal X}$' denotes 
some non-compact `cusp-sector' within the hyperbolic 
upper half-space ${\Bbb H}_3$). Then, by (1.1.23), 
the integral in (6.6.11) is 
$$\int_{\tilde{\cal E}_{\frak c}}  
\bigl({\cal P}^{\frak a}_{\!\!\!\scriptscriptstyle\hookleftarrow}\tau 
{\bf M}_{\omega}\varphi_{\ell,q}(\nu,0)\bigr)(g) 
\,\overline{\phi(g)}\,{\rm d}g 
=\int\limits_{1/|m_{\frak c}|}^{\infty}\int_{{\cal R}_{\frak c}}\int_K 
\bigl({\cal P}^{\frak a}_{\!\!\!\scriptscriptstyle\hookleftarrow}\tau 
{\bf M}_{\omega}\varphi_{\ell,q}(\nu,0)\bigr)\!\left( g_{\frak c} n[z]a[r]k\right) 
\overline{\phi\left( g_{\frak c} n[z]a[r]k\right)}
\,{\rm d}k\,{\rm d}_{+}z\,{{\rm d}r\over r^3}\;,$$
where 
$\tilde{\cal E}_{\frak c}=\{ n[z]a[r]k : k\in K ,\,(z,r)\in{\cal E}_{\frak c}\}$ 
and ${\cal R}_{\frak c}$ is the bounded rectangular region of the complex plane  
defined in (1.1.22). Since $G$ is a topological group it follows,  
by (1.1.3) and what has been noted in the preceding paragraph, that 
the latter integrand (above) is continuous, as a function 
of $(\nu,z,r,k)$, on the set 
$\{\nu\in{\Bbb C} : {\rm Re}(\nu)>1/2\}\times {\Bbb C}\times (0,\infty)\times K$. 
Therefore, in order to justify a similar conclusion 
to that reached at the end of the previous paragraph 
(in respect of the case ${\cal X}={\cal D}$), it is enough  
that we establish a certain uniformity of convergence 
of the above integral over the set 
$(1/|m_{\frak c}|,\infty)\times{\cal R}_{\frak c}\times K$: 
indeed, since each set in the family  
$([1/|m_{\frak c}|,r]\times{\cal R}_{\frak c}\times K)_{r\in{\Bbb N}}$ is compact, 
it suffices that we find some $\delta >0$ such that   
$$\lim_{r_0\rightarrow +\infty}\,\sup_{|\nu -1|\leq\delta}  
\,\biggl|\ \int\limits_{r_0}^{\infty} 
\int_{{\cal R}_{\frak c}}\int_K 
\bigl({\cal P}^{\frak a}_{\!\!\!\scriptscriptstyle\hookleftarrow}\tau 
{\bf M}_{\omega}\varphi_{\ell,q}(\nu,0)\bigr)\!\left( g_{\frak c} n[z]a[r]k\right) 
\overline{\phi\left( g_{\frak c} n[z]a[r]k\right)}
\,{\rm d}k\,{\rm d}_{+}z\,{{\rm d}r\over r^3}\,\biggr| = 0\;.\eqno(6.6.12)$$
\par 
A $\delta >0$ such that (6.6.12) holds may be determined  
by applying the results (6.5.94) and (6.5.96) of Lemma~6.5.16: for, 
by the H\"{o}lder inequality of Section~12.42 of [43], 
the hypothesis that $\phi\in L^{1+\alpha}(\Gamma\backslash G)$, 
the case $\epsilon =2/(1+1/\alpha)\in (0,1/2)$ of (6.5.94), and the 
case $t_1=2/(1+1/\alpha)$, $\sigma_1=1-t_1$, $\sigma_2=1+t_1$ of (6.5.96), 
it follows that, when $|\nu -1|<2/(1+1/\alpha)$, one has 
$({\cal P}^{\frak a}_{\!\!\!\scriptscriptstyle\hookleftarrow}\tau 
{\bf M}_{\omega}\varphi_{\ell,q}(\nu,0))\cdot\overline{\phi}\in L^1(\Gamma\backslash G)$ 
and, for $r_0\geq 1$, 
$$\eqalign{ 
\biggl|\ \int\limits_{r_0}^{\infty} 
\int_{{\cal R}_{\frak c}}\int_K  
 &\bigl({\cal P}^{\frak a}_{\!\!\!\scriptscriptstyle\hookleftarrow}\tau 
{\bf M}_{\omega}\varphi_{\ell,q}(\nu,0)\bigr)\!\left( g_{\frak c} n[z]a[r]k\right) 
\overline{\phi\left( g_{\frak c} n[z]a[r]k\right)}
{\rm d}k\,{\rm d}_{+}z\,{{\rm d}r\over r^3}\,\biggr|^{1+1/\alpha} \leq \cr 
 &\qquad\leq 
\biggl( 2\int_{\Gamma\backslash G} |\phi(g)|^{1+\alpha} {\rm d}g\biggr)^{\!\!1/\alpha} 
\int\limits_{r_0}^{\infty} 
\int_{{\cal R}_{\frak c}}\int_K  
\left|\bigl({\cal P}^{\frak a}_{\!\!\!\scriptscriptstyle\hookleftarrow}\tau 
{\bf M}_{\omega}\varphi_{\ell,q}(\nu,0)\bigr)\!\left( g_{\frak c} n[z]a[r]k\right) 
\right|^{1+1/\alpha} 
\,{\rm d}k\,{\rm d}_{+}z\,{{\rm d}r\over r^3} = \cr 
 &\qquad =O_{\alpha ,\phi}(1) 
\cdot\int\limits_{r_0}^{\infty} 
\int_{{\cal R}_{\frak c}}\int_K  
O_{\Gamma,\omega,\ell,\alpha}\left( r^{(1-{\rm Re}(\nu))(1+1/\alpha)-3}\right) 
{\rm d}k\,{\rm d}_{+}z\,{\rm d}r 
\ll_{\Gamma,\omega,\ell,\alpha,\phi}\ 
\mu^{-1} r_0^{-\mu}\;,\cr 
}$$
where we have $\mu=2-(1+1/\alpha)(1-{\rm Re}(\nu))\geq 2-(1+1/\alpha)|1-\nu|>0$. 
Hence we find that (6.6.12) holds for any $\delta\in (0,2/(1+1/\alpha))$, 
and so may deduce that, when 
${\cal X}\in\{ {\cal E}_{\frak c} : {\frak c}\in {\frak C}(\Gamma)\}$, 
the function (6.6.11) is certainly continuous on the neighbourhood 
$\{ \nu\in{\Bbb C} : |\nu -1|\leq 1/(1+1/\alpha)\}$ of the point $\nu =1$. 
Given the similar result obtained earlier in respect of 
the case ${\cal X}={\cal D}$, our proof of the case $j=1$ of (6.6.10) is 
complete. 
\par 
We now have only to prove the case $j=2$ of (6.6.10): the lemma will then follow. 
Given the definition (6.2.9) of the inner product $\langle f , F\rangle_{N\backslash G}$, 
it follows by (6.5.3) and (6.5.1) that 
we have 
$$\Phi_2(\nu) 
= \bigl\langle \tau{\bf M}_{\omega}\varphi_{\ell,q}(\nu,0)\,, 
F^{\frak a}_{\omega}\phi\bigr\rangle_{N\backslash G} 
= \int\limits_0^2\int_K 
\bigl({\bf M}_{\omega}\varphi_{\ell,q}(\nu,0)\bigr)(a[r]k)
\,\overline{\bigl( F^{\frak a}_{\omega}\phi\bigr)(a[r]k)} 
\,{\rm d}k\,{\tau\!\left( a[r]\right) {\rm d}r\over r^3}\;,\eqno(6.6.13)$$
for all $\nu\in{\Bbb C}$ such that the integral to the right 
of the second equality sign in (6.6.13) exists. 
As noted in the final paragraph of the proof of Lemma~6.5.2, 
it is a corollary of Lemma~6.1 of [5] that the function 
$(\nu,g)\mapsto ({\bf M}_{\omega}\varphi_{\ell,q}(\nu,0))(g)$ is 
continuous on ${\Bbb C}\times G$. We have also 
$\tau\in C^{\infty}(G)$ (by choice), and, since our 
hypothesis that $\phi\in C^0(G)\cap L^{1+\alpha}(\Gamma\backslash G)$ 
implies that $\phi\in C^0(\Gamma\backslash G)$, 
it moreover follows from the definition (1.4.2) that 
$F^{\frak a}_{\omega}\phi$ lies in the space 
$C^0(N\backslash G,\omega)$. Therefore the integrand 
which appears in (6.6.13) is continuous, as a function 
of $(\nu,r,k)$, on the set ${\Bbb C}\times(0,2]\times K$. 
Consequently we may now complete the proof of the case 
$j=2$ of (6.6.10) by finding some $\delta_2>0$ such that 
$$\lim_{r_1\rightarrow 0+}\,\sup_{|\nu -1|\leq\delta_2}  
\,\biggl|\ \int\limits_0^{r_1} 
\int_K 
\bigl({\bf M}_{\omega}\varphi_{\ell,q}(\nu,0)\bigr)(a[r]k) 
\,\overline{\bigl( F^{\frak a}_{\omega}\phi\bigr)(a[r]k)} 
\,{\rm d}k\,{\tau\!\left( a[r]\right) {\rm d}r\over r^3}\,\biggr| = 0\;.\eqno(6.6.14)$$
\par 
By H\"{o}lder's inequality, the hypothesis that 
$|F^{\frak a}_{\omega}\phi|\in 
L^{1+\beta}(N\backslash G)$, the equations in (6.5.1), and the case 
$r_1=1$, $\sigma_0=2$ of the estimate (6.3.9) of Lemma~6.3.1, 
it follows that, when $|\nu -1|<2/(\beta +1)$, one has 
$(\tau{\bf M}_{\omega}\varphi_{\ell,q}(\nu,0))\cdot 
\,\overline{(F^{\frak a}_{\omega}\phi)}\in L^1(N\backslash G)$ and, for 
$0<r_1\leq 1$, 
$$\eqalign{ 
\,\biggl|\ \int\limits_0^{r_1} 
\int_K 
 &\bigl({\bf M}_{\omega}\varphi_{\ell,q}(\nu,0)\bigr)(a[r]k) 
\,\overline{\bigl( F^{\frak a}_{\omega}\phi\bigr)(a[r]k)} 
\,{\rm d}k\,{\tau\!\left( a[r]\right) {\rm d}r\over r^3}\,\biggr|^{1+1/\beta}  \leq \cr 
 &\qquad\quad\leq \biggl(\int_{N\backslash G} 
\left|\bigl(F^{\frak a}_{\omega}\phi\bigr)(g)\right|^{1+\beta}\,{\rm d}\dot{g}
\biggr)^{\!1/\beta}\,\int\limits_0^{r_1} 
\int_K 
\left|\bigl({\bf M}_{\omega}\varphi_{\ell,q}(\nu,0)\bigr)(a[r]k)\right|^{1+1/\beta} 
\,{\rm d}k\,{{\rm d}r\over r^3} = \cr 
 &\qquad\quad = O_{{\frak a},\omega,\phi,\beta}(1)  
\cdot\int\limits_0^{r_1} 
\int_K 
O_{\ell,\omega}\left( r^{(1+{\rm Re}(\nu))(1+1/\beta)-3}\right) 
\,{\rm d}k\,{\rm d}r 
\ll_{{\frak a},\ell,\omega,\phi,\beta}\ \lambda^{-1} r_1^{\lambda}\;, 
}$$
where $\lambda =(1+1/\beta)(1+{\rm Re}(\nu))-2>(1+1/\beta)(2-2/(\beta +1))-2=0$. 
This implies that (6.6.14) holds for 
any $\delta_2\in (0 , 2/(\beta +1))$. 
We are therefore able to conclude that the function $\Phi_2(\nu)$ is certainly 
defined and continuous on the neighbourhood 
$\{\nu\in{\Bbb C} : |\nu -1|\leq 1/(\beta +1)\}$ of the point $\nu =1$; 
since the case $j=2$ of (6.6.10) follows, this 
completes our proof of the lemma\quad$\blacksquare$

\medskip 

\noindent{\bf Remark~6.6.3.}\quad 
Each of the two functions $\Phi_1(\nu)$, $\Phi_2(\nu)$  
considered in the above proof is in fact holomorphic at all points $\nu$ satisfying 
a condition of the form ${\rm Re}(\nu)>a_j\,$ (where, in each case, $a_j$ is less than $1$).  
This may (for example) be shown by adapting the proofs of the 
propositions in Section~2.83 and Section~2.84 of [43]. 

\bigskip 

In applying Lemma~6.6.2 to obtain the geometric description of 
$\langle\phi_1,\phi_2\rangle_{\Gamma\backslash G}$ we require the 
assistance of the following six supplementary lemmas. 

\bigskip 

\proclaim Lemma 6.6.4. Let $\delta^{{\frak a},{\frak b}}_{\omega_1,\omega_2}\in{\Bbb C}$ be 
given by the equation (1.9.2). Then 
$$F^{\frak b}_{\omega_2} P^{{\frak a},*}\widetilde{\bf L}^{\omega_1}_{\ell,q}\eta 
={1\over\left[\Gamma_{\frak a}:\Gamma_{\frak a}'\right]}
\,\delta^{{\frak a},{\frak b}}_{\omega_1,\omega_2} 
\widetilde{\bf L}^{\omega_2}_{\ell,q}\eta 
+{\pi^2\over\left[\Gamma_{\frak a}:\Gamma_{\frak a}'\right]} 
\sum_{c\in{}^{\frak a}{\cal C}^{\frak b}} 
{S_{{\frak a},{\frak b}}\left(\omega_1,\omega_2;c\right)\over |c|^2}
\,\widetilde{\bf L}^{\omega_2}_{\ell,q}\kappa\!\left(\omega_1,\omega_2;c\right)\eta\;, 
\eqno(6.6.15)$$
where, for $c\in{}^{\frak a}{\cal C}^{\frak b}$, the linear operator  
$\kappa(\omega_1,\omega_2;c)$ from ${\cal T}^{\ell}_{\sigma}$ into 
${\cal T}^{\ell}_{\sigma}$ is defined as in Lemma~6.4.3. 

\medskip 

\noindent{\bf Proof.}\quad  
By Lemma~6.5.17 and the definition~(6.5.89), 
$$P^{{\frak a},*}\widetilde{\bf L}^{\omega_1}_{\ell,q}\eta 
=P^{\frak a}\widetilde{\bf L}^{\omega_1,*}_{\ell,q}\eta 
+b(\omega_1;\ell,q;\eta) {\cal P}^{\frak a}_{\!\!\!\scriptscriptstyle\hookleftarrow} 
{\bf M}_{\omega_1}\varphi_{\ell,q}(1,0)
-b(\omega_1;\ell,q;\eta) P^{\frak a}(1-\tau){\bf M}_{\omega_1}\varphi_{\ell,q}(1,0)\;.$$ 
Since it is moreover shown by Lemma~6.5.2, Lemma~6.5.8 and Lemma~6.5.15 
that each of the functions 
$P^{\frak a}(1-\tau){\bf M}_{\omega_1}\varphi_{\ell,q}(1,0)$, 
${\cal P}^{\frak a}_{\!\!\!\scriptscriptstyle\hookleftarrow} 
{\bf M}_{\omega_1}\varphi_{\ell,q}(1,0)$ and 
$P^{\frak a}\widetilde{\bf L}^{\omega_1,*}_{\ell,q}\eta$ is  
both continuous and $\Gamma$-automorphic on $G$, we may therefore 
deduce that 
$$F^{\frak b}_{\omega_2} P^{{\frak a},*}\widetilde{\bf L}^{\omega_1}_{\ell,q}\eta 
=F^{\frak b}_{\omega_2} P^{\frak a}\widetilde{\bf L}^{\omega_1,*}_{\ell,q}\eta 
-b(\eta) 
F^{\frak b}_{\omega_2} P^{\frak a}(1-\tau){\bf M}_{\omega_1}\varphi_{\ell,q}(1,0) 
+b(\eta) 
F^{\frak b}_{\omega_2} {\cal P}^{\frak a}_{\!\!\!\scriptscriptstyle\hookleftarrow} 
{\bf M}_{\omega_1}\varphi_{\ell,q}(1,0)\;, \eqno(6.6.16)$$
where $b(\eta)=b(\omega_1;\ell,q;\eta)$. 
\par 
Let $g\in G$. By (1.4.2) and (6.5.82), 
$$\left( F^{\frak b}_{\omega_2} {\cal P}^{\frak a}_{\!\!\!\scriptscriptstyle\hookleftarrow} 
{\bf M}_{\omega_1}\varphi_{\ell,q}(1,0)\right)(g) 
=\int_{B^{+}\backslash N}\left(\psi_{\omega_2}(n)\right)^{-1} 
\biggl(\ \sum_{\omega'\in{\frak O}}\phi_{\omega'}(1,ng)\biggr){\rm d}n\;,\eqno(6.6.17)$$
where the terms of the sum over $\omega'\in{\frak O}$ are given by the case 
$\omega =\omega_1$ of Equation~(6.5.72). For each $\omega'\in{\frak O}$, the 
mapping $n\mapsto(\psi_{\omega_2}(n))^{-1}\phi_{\omega'}(1,ng)$ 
is a continuous function on $N$;  
since one has both $|\psi_{\omega_2}(n)|=1$ and 
$\rho(ng)=\rho(g)$ for $n\in N$ and $g\in G$, it therefore follows by the 
uniformity of convergence established in Lemma~6.5.14 that one may 
integrate term by term in (6.6.17), so as to obtain the result  
$$\left( F^{\frak b}_{\omega_2} {\cal P}^{\frak a}_{\!\!\!\scriptscriptstyle\hookleftarrow} 
{\bf M}_{\omega_1}\varphi_{\ell,q}(1,0)\right)(g) 
=\sum_{\omega'\in{\frak O}} \int_{B^{+}\backslash N}\left(\psi_{\omega_2}(n)\right)^{-1} 
\phi_{\omega'}(1,ng)\,{\rm d}n\;,\eqno(6.6.18)$$ 
where, as follows by (6.5.72) (with $\omega_1$ substituted for $\omega$), 
the equalities (1.5.16) and (1.4.7)-(1.4.9), 
the identities $n[z]n[w]=n[w]n[z]$ and $h[u]n[z]=n[u^2z]h[u]$, 
the equations in (6.3.3) and (1.4.3), and the definition (1.5.6) of $\delta_{\alpha,\beta}$, 
one has: 
$$\eqalign{
\phi_{\omega'}(1,ng) 
 &=\,{1^{\textstyle{\ \atop\ }}\over\left[\Gamma_{\frak a} : \Gamma_{\frak a}'\right]} 
\sum_{\scriptstyle\gamma\in\Gamma_{\frak a}'\backslash\Gamma 
\ :\ \gamma{\frak b}={\frak a}\atop
{\scriptstyle g_{\frak a}^{-1}\gamma g_{\frak b}\in h\left[ u(\gamma)\right] N  
\atop\scriptstyle\hbox{\qquad}}} 
\delta_{\omega_1 u(\gamma) , \omega'/u(\gamma)} 
\,\bigl({\bf M}_{\omega_1}\varphi_{\ell,q}(1,0)\bigr) 
\left( g_{\frak a}^{-1}\gamma g_{\frak b} g\right)\psi_{\omega'}(n)\  + \cr 
 &\quad\ +\zeta^{{\frak a},{\frak b}}_{\omega_1,\omega'}(1) 
\,\bigl({\bf J}_{\omega'}\varphi_{\ell,q}(1,0)\bigr) (g)\,\psi_{\omega'}(n)\;.}\eqno(6.6.19)$$
Since the family  
$(\psi_{\omega})_{\omega\in{\frak O}}$ is an orthonormal system on 
$B^{+}\backslash N$, it follows from (6.6.18) and (6.6.19) that 
$$\left( F^{\frak b}_{\omega_2} {\cal P}^{\frak a}_{\!\!\!\scriptscriptstyle\hookleftarrow} 
{\bf M}_{\omega_1}\varphi_{\ell,q}(1,0)\right)(g) 
=\phi_{\omega_2}(1,n[0]g)\;.\eqno(6.6.20)$$ 
\par 
By the observations (6.5.9) and (6.5.10) noted within the proof of Lemma~6.5.2, and 
by the first part of (6.5.4), and the case $\theta =\eta$ of (6.5.22), 
it follows (since we assume $\sigma\in(1,2)$) that Lemma~6.2.5 
implies,  for ${\frak a}'={\frak b}$, $\omega'=\omega_2$, $\omega =\omega_1$ and 
$f_{\omega}=f_{\omega_1}\in\{ (1-\tau){\bf M}_{\omega_1}\varphi_{\ell,q}(1,0) , 
\widetilde{\bf L}^{\omega_1,*}_{\ell,q}\eta\}$, the applicability of the 
the formula for $(F^{{\frak a}'}_{\omega'} P^{\frak a} f_{\omega})(g)$ 
stated in (1.5.5)-(1.5.10). Hence, given the definition~(6.5.2), we find that 
$$\eqalign{ 
 &\left[\Gamma_{\frak a} : \Gamma_{\frak a}'\right]\Bigl(
\bigl(F^{\frak b}_{\omega_2} P^{\frak a}\widetilde{\bf L}^{\omega_1,*}_{\ell,q}\eta\bigr)(g)  
-b\left(\omega_1;\ell,q;\eta\right) 
\left(F^{\frak b}_{\omega_2} P^{\frak a}(1-\tau){\bf M}_{\omega_1}\varphi_{\ell,q}(1,0)\right) (g)\Bigr) = \cr 
 &\quad\ = 
\sum_{\scriptstyle\gamma\in\Gamma_{\frak a}'\backslash\Gamma\ :\ \gamma{\frak b}={\frak a}\atop
\scriptstyle g_{\frak a}^{-1}\gamma g_{\frak b}\in h\left[ u(\gamma)\right] N} 
\!\!\!\!\!\!\!\delta_{\omega_1 u(\gamma) , \omega_2/u(\gamma)} 
\,\bigl(\widetilde{\bf L}^{\omega_1,\dagger}_{\ell,q}\eta\bigr)
(g_{\frak a}^{-1}\gamma g_{\frak b} g)  
\,+\sum_{c\in\,{}^{\frak a}{\cal C}^{\frak b}}^{\hbox{\quad }}
S_{{\frak a} , {\frak b}}\!\left(\omega_1 , \omega_2 ; c\right) 
\bigl({\bf J}_{\omega_2}{\bf h}_{1/c}\widetilde{\bf L}^{\omega_1,\dagger}_{\ell,q}\eta\bigr) (g)\;,}
\eqno(6.6.21)$$ 
where the function $\widetilde{\bf L}^{\omega_1,\dagger}_{\ell,q}\eta : G\rightarrow{\Bbb C}$ has 
the definition indicated by Equation~(6.5.20). 
\par 
Regarding now what was observed in (6.6.19) and (6.6.20), we 
note that, by what was found below (6.5.68) concerning the convergence 
of the sum in (6.5.61), and by the result~(6.3.11) of Lemma~6.3.2, 
the equation~(6.3.5) and 
the linearity of the operator ${\bf J}_{\omega_2}$, it follows that  
$$\eqalign{ 
\left[\Gamma_{\frak a} : \Gamma_{\frak a}'\right] 
\zeta^{{\frak a},{\frak b}}_{\omega_1,\omega_2}(1) 
\,{\bf J}_{\omega_2}\varphi_{\ell,q}(1,0) 
 &=\sum_{c\in\,{}^{\frak a}{\cal C}^{\frak b}} 
{S_{{\frak a} , {\frak b}}\!\left(\omega_1 , \omega_2 ; c\right)\over |c|^4} 
\,{\cal J}^{*}_{1,0}\!\left( 2\pi\sqrt{c^{-2}\omega_1\omega_2}\right) 
{\bf J}_{\omega_2}\varphi_{\ell,q}(1,0) = \cr 
 &=\sum_{c\in\,{}^{\frak a}{\cal C}^{\frak b}}^{\hbox{\ }^{\hbox{\ }}} 
S_{{\frak a} , {\frak b}}\!\left(\omega_1 , \omega_2 ; c\right)
{\bf J}_{\omega_2}{\bf h}_{1/c}{\bf M}_{\omega_1}\varphi_{\ell,q}(1,0)\;.}$$
Therefore, given that the definition (6.5.20) 
implies the identity 
$\widetilde{\bf L}^{\omega_1,\dagger}_{\ell,q}\eta 
+b(\omega_1;\ell,q;\eta){\bf M}_{\omega_1}\varphi_{\ell,q}(1,0) 
=\widetilde{\bf L}^{\omega_1}_{\ell,q}\eta$, it follows by (6.6.16), (6.6.19), 
(6.6.20) and (6.6.21) that 
$$\eqalign{ 
\left[\Gamma_{\frak a} : \Gamma_{\frak a}'\right] 
\bigl(F^{\frak b}_{\omega_2} P^{{\frak a},*} 
\widetilde{\bf L}^{\omega_1}_{\ell,q}\eta\bigr)(g) 
 &=\ \sum_{\scriptstyle\gamma\in\Gamma_{\frak a}'\backslash\Gamma\ :\ \gamma{\frak b}={\frak a}\atop
\scriptstyle g_{\frak a}^{-1}\gamma g_{\frak b}\in h\left[ u(\gamma)\right] N} 
\!\!\!\!\delta_{\omega_1 u(\gamma) , \omega_2/u(\gamma)} 
\,\bigl(\widetilde{\bf L}^{\omega_1}_{\ell,q}\eta\bigr)
(g_{\frak a}^{-1}\gamma g_{\frak b} g)\  + \cr 
 &\quad\ +\sum_{c\in\,{}^{\frak a}{\cal C}^{\frak b}}^{\hbox{\ }^{\hbox{\ }}}
S_{{\frak a} , {\frak b}}\!\left(\omega_1 , \omega_2 ; c\right) 
\bigl({\bf J}_{\omega_2}{\bf h}_{1/c}\widetilde{\bf L}^{\omega_1}_{\ell,q}\eta\bigr) (g)\;, 
}\eqno(6.6.22)$$
where $\widetilde{\bf L}^{\omega_1}_{\ell,q}\eta : G\rightarrow{\Bbb C}$ is given 
by the case $\omega =\omega_1$ of (6.4.4)-(6.4.5). 
\par 
By Lemma~2.1, the conditions imposed in (6.6.22) 
on the variable of summation $\gamma$ imply that 
one has there $u^2(\gamma)=(u(\gamma))^2\in{\frak O}^{*}$. 
Consequently, given that the relevant summand is non-zero 
only when $\omega_1 u(\gamma)=\omega_2/u(\gamma)$, the sum over 
$\gamma$ in (6.6.22) is effectively empty unless one has $\omega_1\sim\omega_2$;  
furthermore, if $\omega_1\sim\omega_2$, then that summation 
is effectively restricted to $\gamma\in\Gamma_{\frak a}'\backslash\Gamma$ such that 
$g_{\frak a}^{-1}\gamma g_{\frak b} =h[u]n$ for some (necessarily unique) 
pair $(u,n)=(u(\gamma),n[z(\gamma)])\in{\Bbb C}^{*}\times N$ such that the 
number $\epsilon =\epsilon(\gamma)=u^2\in{\Bbb C}^{*}$ satisfies 
$\epsilon =\omega_2/\omega_1\in{\frak O}^{*}$. 
Since it moreover follows by (6.4.4), the case $n=n[0]$ of (1.8.2), and 
the identity ${\bf h}_{u}{\bf J}_{\omega}{\bf h}_{u}=|u|^4 {\bf J}_{u^2\omega}$ 
that ${\bf h}_{\pm\sqrt{\epsilon}}\,\widetilde{\bf L}^{\omega}_{\ell,q}\eta 
=\widetilde{\bf L}^{\epsilon\omega}_{\ell,q}\eta$ for 
$0\neq\omega\in{\Bbb C}$ and any $\epsilon\in{\Bbb C}$ of unit modulus, 
and since one has 
$\widetilde{\bf L}^{\omega_2}_{\ell,q}\eta\in C^{\infty}(N\backslash G,\omega_2)\,$  
(directly by (6.4.4), or by (6.4.7)), we therefore find that 
$$\eqalign{
\sum_{\scriptstyle\gamma\in\Gamma_{\frak a}'\backslash\Gamma\ :\ \gamma{\frak b}={\frak a}\atop
\scriptstyle g_{\frak a}^{-1}\gamma g_{\frak b}\in h\left[ u(\gamma)\right] N} 
\!\!\!\delta_{\omega_1 u(\gamma) , \omega_2/u(\gamma)} 
\,\bigl(\widetilde{\bf L}^{\omega_1}_{\ell,q}\eta\bigr)
(g_{\frak a}^{-1}\gamma g_{\frak b} g) 
 &=\ \sum_{\scriptstyle\gamma\in\Gamma_{\frak a}'\backslash\Gamma\ :\ \gamma{\frak b}={\frak a}\atop
\scriptstyle g_{\frak a}^{-1}\gamma g_{\frak b}=h\left[ u(\gamma)\right] 
n\left[z(\gamma)\right]} 
\!\!\!\!\!\!\!\!\!\delta_{\omega_1 u(\gamma) , \omega_2/u(\gamma)} 
\,\bigl({\bf h}_{u(\gamma)}\widetilde{\bf L}^{\omega_1}_{\ell,q}\eta\bigr)
(n\!\left[z(\gamma)\right] g) = \cr 
 &=\ \sum_{\scriptstyle\gamma\in\Gamma_{\frak a}'\backslash\Gamma\ :\ \gamma{\frak b}={\frak a}\atop
\scriptstyle g_{\frak a}^{-1}\gamma g_{\frak b}=h\left[ u(\gamma)\right] 
n\left[z(\gamma)\right]}^{\matrix{\ }}  
\!\!\!\!\!\!\!\!\!\delta_{\omega_1 u(\gamma) , \omega_2/u(\gamma)} 
\,\bigl(\widetilde{\bf L}^{\omega_2}_{\ell,q}\eta\bigr)
(n\!\left[z(\gamma)\right] g) = \cr 
 &=\ \sum_{\scriptstyle\gamma\in\Gamma_{\frak a}'\backslash\Gamma\ :\ \gamma{\frak b}={\frak a}\atop
\scriptstyle g_{\frak a}^{-1}\gamma g_{\frak b}=h\left[ u(\gamma)\right] 
n\left[z(\gamma)\right]}^{\matrix{\ }}  
\!\!\!\!\!\!\!\!\!\delta_{\omega_1 u(\gamma) , \omega_2/u(\gamma)} 
\,\psi_{\omega_2}\!\left( n\!\left[ z(\gamma)\right]\right) 
\bigl(\widetilde{\bf L}^{\omega_2}_{\ell,q}\eta\bigr)(g) = \cr 
&=\ \delta^{{\frak a},{\frak b}}_{\omega_1,\omega_2}\cdot  
\bigl(\widetilde{\bf L}^{\omega_2}_{\ell,q}\eta\bigr)(g)\;.\matrix{\ \cr\ }}$$
By this, the equation (6.6.22), and the result (6.4.13) of Lemma~6.4.3, 
we obtain the result in (6.6.15)\quad$\blacksquare$ 

\bigskip 

\proclaim Lemma 6.6.5. Let $0\neq\omega\in{\Bbb C}$. Then, for each 
$\beta\in(0,\infty)$, one has $L^{1+\beta}(N\backslash G)\ni 
|\widetilde{\bf L}^{\omega}_{\ell,q}\eta|$. 

\medskip 

\noindent{\bf Proof.}\quad  
In view of the definition (implicit in (6.2.9)) of the measure ${\rm d}\dot{g}$ on $N\backslash G$, 
this lemma is a straightforward corollary of the results (6.4.7) and (6.4.8) 
of Theorem~6.4.1\quad$\blacksquare$ 

\bigskip 

\proclaim Lemma 6.6.6 (Bruggeman and Motohashi). Let $c\in{\Bbb C}^{*}$,  
and let $u=2\pi\sqrt{\omega_1\omega_2}/c$. Then, for $g\in G$ and 
$0<\alpha\leq\sigma$, one has 
$$\eqalignno{ 
 &\bigl(\widetilde{\bf L}^{\omega_2}_{\ell,q}\kappa\!\left(\omega_1,\omega_2;c\right) 
\eta\bigr)(g) = \cr 
 &\qquad\quad = {(-1)\over\pi^3 i}\sum_{|p|\leq\ell} 
\!{\left( -i\omega_2/\left|\omega_2\right|\right)^p\over 
\left\|\Phi^{\ell}_{p,q}\right\|_K} 
\int\limits_{(\alpha)}\!{\cal J}_{\nu,p}(u)\,\eta(\nu,p) 
\left(\pi\left|\omega_2\right|\right)^{-\nu} 
\Gamma(\ell+1+\nu)\,\bigl({\bf J}_{\omega_2}\varphi_{\ell,q}(\nu,p)\bigr) (g) 
\,\nu^{\epsilon(p)}{\rm d}\nu\ +\qquad \cr 
 &\qquad\qquad\ +{\ell !\over\pi^2}\sum_{0<|p|\leq\ell} 
{\left( -i\omega_2/\left|\omega_2\right|\right)^p\over 
\left\|\Phi^{\ell}_{p,q}\right\|_K} 
\,{\cal J}_{0,p}(u)\,\eta(0,p) 
\,\bigl({\bf J}_{\omega_2}\varphi_{\ell,q}(0,p)\bigr) (g)
&(6.6.23)}$$
where ${\cal J}_{\nu,p} : {\Bbb C}^{*}\rightarrow{\Bbb C}$ and 
$\epsilon : {\Bbb Z}\rightarrow\{ -1 , 1\}$ are given by (1.9.5)-(1.9.6) and (6.4.5). 

\medskip 

\noindent{\bf Proof.}\quad  
This identity is implied by the equations~(7.22) and~(7.24) of [5]\quad$\blacksquare$ 

\bigskip 

\proclaim Lemma 6.6.7. Let $p\in{\Bbb Z}$, and let $u\in{\Bbb C}^{*}$. Then 
$$0\leq (-1)^p {\cal J}_{0,p}(u)\leq 
(|p|!)^{-2} |u/2|^{2|p|}\exp\left( |u|^2 /2\right)\;.\eqno(6.6.24)$$ 
Suppose moreover that $\delta\in(0,2]$, and that $\nu\in{\Bbb C}$ is such that 
$\min\{ |m-\nu| : m\in{\Bbb Z}\}\geq\delta$. Then one has also 
$$\left|{\cal J}_{\nu,p}(u)\right| 
\leq {4\delta^{-2} |u/2|^{2{\rm Re}(\nu)}\exp\left( |u|^2 /2\right)\over 
|\Gamma(\nu -p+1)\Gamma(\nu +p+1)|}\;.\eqno(6.6.25)$$ 

\medskip 

\noindent{\bf Proof.}\quad  
By (1.9.5) and (1.9.9), one has 
$(-1)^p {\cal J}_{0,p}(u)=|u/2|^{2|p|} 
J^{*}_{|p|}(u) J^{*}_{|p|}\left(\overline{u}\right)$. 
It is moreover implied by the power series representation (1.9.6) 
of $J_{\xi}^{*}(z)$ that 
$J^{*}_{|p|}(u) J^{*}_{|p|}\left(\overline{u}\right) 
=|J^{*}_{|p|}(u)|^2\leq ((|p|!)^{-1}\exp(|u/2|^2))^2$, and so the result 
(6.6.24) follows. 
By (1.9.6) (again), we have also 
$$\left| J^{*}_{\xi}(z)\right| 
\leq\sum_{m=0}^{\infty} {|z/2|^{2m}\over (m!)(1/2)\delta |\Gamma(\xi +1)|} 
={2\delta^{-1}\exp\left( |z/2|^2\right)\over |\Gamma(\xi +1)|}\qquad\quad 
\hbox{($z\in{\Bbb C}^{*}$)}$$ 
whenever $\xi\in{\Bbb C}$ and 
$\min\{ |m'-\xi| : m'\in{\Bbb Z}\}\geq\delta >0$. We therefore find (subject to 
$\delta$ and $\nu$ satisfying the hypotheses stated above (6.6.25)) that  
$|J^{*}_{\nu -p}(u) J^{*}_{\nu +p}\left(\overline{u}\right) |
\leq 4\delta^{-1}|\Gamma(\nu -p+1)\Gamma(\nu +p+1)|^{-1}\exp(|u|^2 /2)$, and 
so (given the definition (1.9.5)) we obtain the result (6.6.25)\quad$\blacksquare$ 

\bigskip 

\proclaim Lemma 6.6.8. Let $\alpha\in(1/2,1)$, let $\varepsilon\in(0,1/2]$, 
and let $j\in{\Bbb N}$. Suppose moreover that $c\in{\Bbb C}^{*}$, and 
that $|c|\geq c_0>0$. Then, when $g\in G$ and $r=\rho(g)$, one has: 
$$\bigl(\widetilde{\bf L}^{\omega_2}_{\ell,q}
\kappa\!\left(\omega_1,\omega_2;c\right)\eta\bigr)(g) 
=\cases{ 
O_{\eta,\alpha,c_0,\omega_1,\omega_2,\epsilon}\left( |c|^{-2\alpha} 
r^{(1-\alpha)(1-\varepsilon)}\right)  
&if $\,r\leq 1$; \cr 
O_{\eta,\alpha,c_0,\omega_1,\omega_2,j}\left( |c|^{-2\alpha} r^{-j}\right) 
&if $\,r\geq 1$.}\eqno(6.6.26)$$

\medskip 

\noindent{\bf Proof.}\quad  
Let $g\in G$, and let $r=\rho(g)$. If $r\leq 1$ then, by the 
result (6.6.23) of Lemma~6.6.6, 
the inequalities (6.6.24) and (6.6.25) of Lemma~6.6.7, the bound (6.5.15) of 
Lemma~6.5.3$\,$ (for $\omega'=\omega_2$, $\sigma_1=1$, $r_0=|\omega_2|$ and 
$d=2\ell +3\,$ (say), and 
with $(1-\alpha)\varepsilon$ substituted for $\varepsilon$)  
and the lower bound for $|\Gamma(\mu +1)|$ in (6.5.19), one finds 
(since $\epsilon(p)\leq 1$ for $p\in{\Bbb Z}$) that 
$$\eqalign{
\bigl(\widetilde{\bf L}^{\omega_2}_{\ell,q}
\kappa\!\left(\omega_1,\omega_2;c\right)\eta\bigr)(g) 
\ll_{\ell,\omega_1,\omega_2,c_0,\alpha,\varepsilon}\ \,  
 &|u|^{2\alpha} r^{(1-\alpha)(1-\varepsilon)}\sum_{|p|\leq\ell} 
\ \int\limits_{-\infty}^{\infty} |\eta(\alpha +it,p)| e^{(\pi /2)|t|} 
(1+|t|)^{\ell -\alpha -1/2} {\rm d}t\ + \cr 
 &\ +|u|^2 r^{1-(1-\alpha)\varepsilon}\sum_{0<|p|\leq\ell}^{\ }|\eta(0,p)|\;,}$$ 
where $u=2\pi\sqrt{\omega_1\omega_2}/c$. 
By this, our hypotheses concerning $\alpha$ and $c$, and the conditions (T2) and (T3) 
stated below (6.4.3), the case $\rho(g)=r\leq 1$ of (6.6.26) follows. 
The other case of the result~(6.6.26) may be proved similarly: for, 
when $\rho(g)=r\geq 1$, it follows by Lemma~6.6.6, Lemma~6.6.7 and Lemma~6.5.3 
(with the result (6.5.15) of the latter being applied for 
$\omega'=\omega_2$, $\sigma_1=(\alpha +1)/2$, $r_0=|\omega_2|$ and $d=2\ell +1+j$) 
that one has 
$$\eqalign{
\bigl(\widetilde{\bf L}^{\omega_2}_{\ell,q}
\kappa\!\left(\omega_1,\omega_2;c\right)\eta\bigr)(g) 
\ll_{\ell,\omega_1,\omega_2,c_0,\alpha,j}\ \,  
 &|u|^{2\alpha} r^{-j-\alpha}\sum_{|p|\leq\ell} 
\ \int\limits_{-\infty}^{\infty} |\eta(\alpha +it,p)| e^{(\pi /2)|t|} 
(1+|t|)^{3\ell +j-\alpha +1/2} {\rm d}t\ + \cr 
 &\ +|u|^2 r^{-j}\sum_{0<|p|\leq\ell}^{\ }|\eta(0,p)|\;,}$$ 
where, by hypothesis, $1/2<\alpha <1<\sigma$  and 
$|u|^2=|2\pi\sqrt{\omega_1\omega_2}/c|^2\leq 4\pi^2|c_0|^{-2}|\omega_1\omega_2|\quad\blacksquare$ 

\bigskip 

\proclaim Lemma 6.6.9. For each real $\beta >3$, one has 
$L^{1+\beta}(N\backslash G)\ni |F^{\frak b}_{\omega_2} P^{{\frak a},*} 
\widetilde{\bf L}^{\omega_1}_{\ell,q}\eta|$. 

\medskip 

\noindent{\bf Proof.}\quad  
Put $\delta = [\Gamma_{\frak a} : \Gamma_{\frak a}']^{-1}\delta^{{\frak a},{\frak b}}_{\omega_1 , 
\omega_2}$ (a complex constant), $f=\widetilde{\bf L}^{\omega_2}_{\ell,q}\eta\,$ and  
$F=F^{\frak b}_{\omega_2} P^{{\frak a},*} 
\widetilde{\bf L}^{\omega_1}_{\ell,q}\eta - \delta f$, so that by Lemma~6.6.4 one has 
$$F(g)={\pi^2\over\left[\Gamma_{\frak a}:\Gamma_{\frak a}'\right]} 
\sum_{c\in{}^{\frak a}{\cal C}^{\frak b}} 
{S_{{\frak a},{\frak b}}\left(\omega_1,\omega_2;c\right)\over |c|^2}
\ \bigl( \widetilde{\bf L}^{\omega_2}_{\ell,q}\kappa\left(\omega_1,\omega_2;c\right)
\eta\bigr) (g)\qquad\qquad\hbox{($g\in G$).}\eqno(6.6.27)$$
By Theorem~6.4.1 and Lemma~6.5.17, we have   
$\widetilde{\bf L}^{\omega_2}_{\ell,q}\eta\in C^{\infty}(N\backslash G,\omega_2)$ and 
$P^{{\frak a},*}\widetilde{\bf L}^{\omega_1}_{\ell,q}\eta\in C^{\infty}(\Gamma\backslash G)$.  
We consequently 
have $\{ f , F+\delta f\}\subset C^{\infty}(N\backslash G,\omega_2)$, and so 
have also $|f|,|F|,|F+\delta f|\in C^0(N\backslash G,0)\,$ 
(the functions $|f|$, $|F|$ and $|F+\delta f|$ are, in particular, 
measurable). It therefore follows by the 
H\"{o}lder inequality 
$|F+\delta f|^{1+\beta}\leq (1+|\delta|^{1+1/\beta})^{\beta}(|F|^{1+\beta}+|f|^{1+\beta})$ 
and Lemma~6.6.5 that, for each $\beta\in(0,\infty)$ such that  
$$\int_{N\backslash G} |F(g)|^{1+\beta} {\rm d}\dot{g} <\infty\;,\eqno(6.6.28)$$ 
one has   
$\bigl| F^{\frak b}_{\omega_2} P^{{\frak a},*} 
\widetilde{\bf L}^{\omega_1}_{\ell,q}\eta\bigr| 
=|F+\delta f|=\bigl| F+\delta\widetilde{\bf L}^{\omega_2}_{\ell,q}\eta\bigr|\in L^{1+\beta}(N\backslash G)$. 
Hence this proof will be complete once we have  
verified that the condition (6.6.28) is satisfied for all real $\beta >3$. 
\par 
Suppose now that $\beta$ is a positive real number. Then, by (6.6.27), the bounds 
(6.6.26) of Lemma~6.6.8, the bound (6.5.59) of Lemma~6.5.9, the equation~(6.1.25) and 
the result (6.1.26) of Lemma~6.1.5, 
it follows that when $g\in G$, $r=\rho(g)$ and $1/2<\alpha<1$ one has 
$$F(g)
= O_{\Gamma,|\omega_1|}\!\left(\zeta_{{\Bbb Q}(i)}^2\Bigl(\alpha + {1\over 2}\Bigr) 
\right)\cdot 
\cases{O_{\eta,\alpha,\omega_1,\omega_2}\!\left( r^{(1-\alpha)(1-(\alpha -1/2))}\right)  
 &if $r\leq 1$; \cr 
O_{\eta,\alpha,\omega_1,\omega_2}\!\left( r^{-1}\right) 
 &if $r\geq 1$.}\eqno(6.6.29)$$
Given the definition of the measure ${\rm d}\dot{g}$ implicit in (6.2.9), 
it follows from (6.6.29) that when $\alpha\in (1/2,\infty)\,$ 
(so that $(1-\alpha)(1-(\alpha -1/2))=(1/2)-(3/2)(\alpha -1/2)+(\alpha -1/2)^2 
>(1/2)-3(\alpha -1/2)=2-3\alpha$) one has:  
$$\int_{N\backslash G} |F(g)|^{1+\beta} {\rm d}\dot{g} 
\ll_{\Gamma,\eta,\alpha,\beta,\omega_1,\omega_2} 
\,\int\limits_0^1 {r^{(2-3\alpha)(1+\beta)} {\rm d}r\over r^3} 
+\int\limits_1^{\infty} {r^{-(1+\beta)}{\rm d}r\over r^3}  
=\int\limits_0^1 \!r^{(2-3\alpha)(1+\beta)-3}{\rm d}r 
+{1\over\beta +3}\;.\hbox{\ }\eqno(6.6.30)$$
Therefore the condition (6.6.28) is satisfied (for the given 
choice of $\beta\in (0,\infty)$) if there is some $\alpha >1/2$ 
such that $(2-3\alpha)(1+\beta)>2$. The latter is the case if and only if 
one has $(2-3\alpha)(1+\beta)>2$ when $\alpha =1/2\quad\blacksquare$ 

\bigskip 

\noindent{\bf Part~I of the proof of Proposition~6.6.1: 
the geometric description of 
$\langle\phi_1 , \phi_2\rangle_{\Gamma\backslash G}$.}\quad  
By (6.6.1) and Lemma~6.6.2 (with ${\frak b}$, $\theta$ and $\omega_2$ substituted for 
${\frak a}$, $\eta$ and $\omega$, respectively), one obtains 
$$\eqalign{ 
\left[\Gamma_{\frak b} : \Gamma_{\frak b}'\right]
\langle\phi_1 , \phi_2\rangle_{\Gamma\backslash G} 
 &=\overline{\left[\Gamma_{\frak b} : \Gamma_{\frak b}'\right]
\langle P^{{\frak b},*}\widetilde{\bf L}^{\omega_2}_{\ell,q}\theta , \phi_1\rangle_{\Gamma\backslash G}}\,=\cr 
 &=\overline{\langle\widetilde{\bf L}^{\omega_2}_{\ell,q}\theta , 
F^{\frak b}_{\omega_2}\phi_1\rangle_{N\backslash G}} 
=\langle F^{\frak b}_{\omega_2} P^{{\frak a},*}\widetilde{\bf L}^{\omega_1}_{\ell,q}\eta\,, 
\widetilde{\bf L}^{\omega_2}_{\ell,q}\theta\rangle_{N\backslash G}\;. 
}\eqno(6.6.31)$$ 
This application of Lemma~6.6.2 is justified: for, given (6.6.1), 
it follows by Lemma~6.5.17 and (6.5.7) that one has 
$\phi_1\in C^0(G)\cap L^2(\Gamma\backslash G)$,  
while by Lemma~6.6.9 one has (for example)  
$\,| F^{\frak b}_{\omega_2}\phi_1 |\in L^5(N\backslash G)$. 
\par 
In preparation for an application of Lemma~6.6.4 we note  
that, by  Lemma~6.4.2, the inner product 
$\langle\widetilde{\bf L}^{\omega_2}_{\ell,q}\eta \,, \widetilde{\bf L}^{\omega_2}_{\ell,q}\theta 
\rangle_{N\backslash G}$ exists (as an integral with respect to the measure 
${\rm d}\dot g$ on $N\backslash G$) and  one has, 
for any $\delta\in{\Bbb C}$,   
$$\eqalign{ 
\langle F^{\frak b}_{\omega_2} P^{{\frak a},*}\widetilde{\bf L}^{\omega_1}_{\ell,q}\eta 
 &\,, \widetilde{\bf L}^{\omega_2}_{\ell,q}\theta\rangle_{N\backslash G} 
-{\delta\over\pi i}\sum_{p\in{\Bbb Z}}\ \int\limits_{(0)} h_{\ell}(\nu,p) 
\left( p^2-\nu^2\right) {\rm d}\nu = \cr 
 &=\langle F^{\frak b}_{\omega_2} P^{{\frak a},*}\widetilde{\bf L}^{\omega_1}_{\ell,q}\eta \,, 
\widetilde{\bf L}^{\omega_2}_{\ell,q}\theta\rangle_{N\backslash G} 
-\delta\,\langle\widetilde{\bf L}^{\omega_2}_{\ell,q}\eta \,, 
\widetilde{\bf L}^{\omega_2}_{\ell,q}\theta\rangle_{N\backslash G} 
=\langle F^{\frak b}_{\omega_2} P^{{\frak a},*}\widetilde{\bf L}^{\omega_1}_{\ell,q}\eta 
-\delta\widetilde{\bf L}^{\omega_2}_{\ell,q}\eta \,, 
\widetilde{\bf L}^{\omega_2}_{\ell,q}\theta\rangle_{N\backslash G}\;,}$$ 
where $h_{\ell} : \{\nu\in{\Bbb C} : 
|{\rm Re}(\nu)|\leq\sigma\}\times{\Bbb Z}\rightarrow 
{\Bbb C}$ is the function given by (6.6.3)-(6.6.4). 
In the above we may put $\delta = [\Gamma_{\frak a} : \Gamma_{\frak a}']^{-1} 
\delta^{{\frak a},{\frak b}}_{\omega_1,\omega_2}\,$ (with 
$\delta^{{\frak a},{\frak b}}_{\omega_1,\omega_2}\in{\Bbb C}$ as defined in (1.9.2)). 
Hence, given the result (6.6.15) of Lemma~6.6.4, we are able to deduce 
that 
$$\eqalign{ 
\langle F^{\frak b}_{\omega_2} P^{{\frak a},*}\widetilde{\bf L}^{\omega_1}_{\ell,q}\eta   
 &\,, \widetilde{\bf L}^{\omega_2}_{\ell,q}\theta\rangle_{N\backslash G} = \cr 
 &={1\over\pi i\left[\Gamma_{\frak a} : \Gamma_{\frak a}'\right]}  
\,\delta^{{\frak a},{\frak b}}_{\omega_1,\omega_2}\sum_{p\in{\Bbb Z}}\ \int\limits_{(0)} h_{\ell}(\nu,p) 
\left( p^2-\nu^2\right) {\rm d}\nu  
+\int_{N\backslash G} F(g)\,\overline{\bigl(\widetilde{\bf L}^{\omega_2}_{\ell,q}\theta\bigr) 
(g)}\,{\rm d}\dot{g}\;,}\eqno(6.6.32)$$
where $F=F^{\frak b}_{\omega_2} P^{{\frak a},*}\widetilde{\bf L}^{\omega_1}_{\ell,q}\eta 
-[\Gamma_{\frak a} : \Gamma_{\frak a}']^{-1} 
\delta^{{\frak a},{\frak b}}_{\omega_1,\omega_2}\widetilde{\bf L}^{\omega_2}_{\ell,q}\eta$ 
is the very same function as occurs (within the proof of Lemma~6.6.9) 
on the left-hand side of the identity in (6.6.27). 
\par 
In order to progress beyond (6.6.32) 
we must first show that the  
sum over $c\in{}^{\frak a}{\cal C}^{\frak b}$ occurring on the right-hand side of 
the identity in (6.6.27) may be integrated over $N\backslash G$ term by term. 
We shall achieve this through an application of 
Lebesgue's theorem on `dominated convergence', Theorem~1.34 of~[40]. 
\par 
By the relation (6.1.26) of Lemma~6.1.5, the set ${}^{\frak a}{\cal C}^{\frak b}$ 
is a countable subset of ${\Bbb C}^{*}$. Let $c_{*}$ be any 
one-to-one function with domain ${\Bbb N}$ and range ${}^{\frak a}{\cal C}^{\frak b}$,  
and let $(F_M)_{M\in{\Bbb N}}$ be the sequence of functions on $G$ given by: 
$$F_M(g)={\pi^2\over\left[\Gamma_{\frak a} : \Gamma_{\frak a}'\right]} 
\sum_{m=1}^M {S_{{\frak a},{\frak b}}\!\left(\omega_1,\omega_2;c_{*}(m)\right)\over 
\left| c_{*}(m)\right|^2}
\,\bigl(\widetilde{\bf L}^{\omega_2}_{\ell,q}\kappa\!\left(\omega_1,\omega_2;c_{*}(m)\right) 
\eta\bigr)(g)\qquad\quad\hbox{($M\in{\Bbb N}$, $g\in G$).}\eqno(6.6.33)$$
Given how (within the proof of Lemma~6.6.9) 
we obtained the estimate (6.6.29), it may be inferred that for all $g\in G$   
the sum over $c\in{}^{\frak a}{\cal C}^{\frak b}$ in (6.6.27) is 
absolutely convergent, and that for any $M\in{\Bbb N}$ one may 
substitute $F_M(g)$ in place of $F(g)$ in (6.6.29). 
Therefore, and since (6.6.29) implies (6.6.30), we 
may deduce that for each $(\alpha,\beta)\in{\Bbb R}^2$ such that 
$1+\beta>2/(2-3\alpha)>4$ there 
exists some function $D_{\alpha}\in L^{1+\beta}(N\backslash G)$ 
satisfying the condition 
$$D_{\alpha}(g)\geq\left| F_M(g)\right|\qquad\ 
\hbox{for all $\,M\in{\Bbb N}$, $g\in G$.}$$
In particular (by the case $\alpha =4/7$, $\beta =7$ of the above), 
there must exist some function $D\in L^8(N\backslash G)$ 
such that $D : G\rightarrow [0,\infty)$ and 
$$\bigl| F_M(g)
\,\overline{\bigl(\widetilde{\bf L}^{\omega_2}_{\ell,q}\theta\bigr)(g)}\,\bigr| 
\leq \bigl| \bigl(\widetilde{\bf L}^{\omega_2}_{\ell,q}\theta\bigr)(g)\bigr|\,D(g)\qquad\ 
\hbox{for all $\,M\in{\Bbb N}$, $g\in G$.}\eqno(6.6.34)$$
Since Lemma~6.6.5 implies that 
$|\widetilde{\bf L}^{\omega_2}_{\ell,q}\theta |\in L^{8/7}(N\backslash G)$, 
it follows by the H\"{o}lder inequality of Section~12.42 of [43] that the above 
function $D$ will, since it lies in $L^8(N\backslash G)$, be such that   
$$|\widetilde{\bf L}^{\omega_2}_{\ell,q}\theta|\cdot D\in L^1(N\backslash G)\;.\eqno(6.6.35)$$  
\par 
By Lemma~6.4.3, and by Theorem~6.4.1 and Lemma~6.4.2 (applied 
with $\kappa(\omega_1,\omega_2;c)\eta$ substituted for $\eta$), 
$$\bigl(\widetilde{\bf L}^{\omega_2}_{\ell,q}\kappa\left(\omega_1,\omega_2;c\right)\eta\bigr) 
\cdot\,\overline{\bigl(\widetilde{\bf L}^{\omega_2}_{\ell,q}\theta\bigr)} 
\in L^1(N\backslash G)\qquad\ 
\hbox{for $\,c\in{}^{\frak a}{\cal C}^{\frak b}\subset{\Bbb C}^{*}$.}\eqno(6.6.36)$$ 
Since the sum over $c\in{}^{\frak a}{\cal C}^{\frak b}$ in (6.6.27) 
is absolutely convergent for all $g\in G$, it follows by the definition (6.6.33) 
and our hypotheses concerning the function $c_{*}$ that, for each $g\in G$,  
$$F(g)\,\overline{\bigl(\widetilde{\bf L}^{\omega_2}_{\ell,q}\theta\bigr)(g)} 
=\lim_{M\rightarrow\infty} F_M(g)
\,\overline{\bigl(\widetilde{\bf L}^{\omega_2}_{\ell,q}\theta\bigr)(g)}  
=\lim_{M\rightarrow\infty} \left( F_M
\cdot\,\overline{\bigl(\widetilde{\bf L}^{\omega_2}_{\ell,q}\theta\bigr)}\,\right)\!(g)\qquad\ 
\hbox{(say),}$$
and that the relations in (6.6.36) imply the relation   
$L^1(N\backslash G)\supseteq\{  
F_M\cdot\,\overline{(\widetilde{\bf L}^{\omega_2}_{\ell,q}\theta)} 
\,: M\in{\Bbb N}\}$.
Therefore, given that we have also (6.6.34) and (6.6.35), 
it follows by Lebesgue's theorem on dominated convergence   
that 
$$\lim_{M\rightarrow\infty}
\int_{N\backslash G} F_M(g)\,\overline{\bigl(\widetilde{\bf L}^{\omega_2}_{\ell,q}\theta\bigr) 
(g)}\,{\rm d}\dot{g} 
=\int_{N\backslash G} F(g)\,\overline{\bigl(\widetilde{\bf L}^{\omega_2}_{\ell,q}\theta\bigr) 
(g)}\,{\rm d}\dot{g}\;.\eqno(6.6.37)$$ 
\par 
Equation~(6.6.37) enables us to justify 
term by term integration (over $N\backslash G$) of the sum over 
$c\in{}^{\frak a}{\cal C}^{\frak b}$ seen in (6.6.27). 
Indeed, by the definitions (6.6.33) and (6.2.9) and the relations in (6.6.36), 
and by Lemma~6.4.3 and Lemma~6.4.2, we find that for $M\in{\Bbb N}$ one has 
$$\eqalign{ 
{\left[\Gamma_{\frak a} : \Gamma_{\frak a}'\right]\over 4\pi^2} 
 &\int_{N\backslash G} F_M(g)
\,\overline{\bigl(\widetilde{\bf L}^{\omega_2}_{\ell,q}\theta\bigr) 
(g)}\,{\rm d}\dot{g} = \cr 
 &=\biggl\langle {1\over 4} 
\sum_{m=1}^M {S_{{\frak a},{\frak b}}\!\left(\omega_1,\omega_2;c_{*}(m)\right)\over 
\left| c_{*}(m)\right|^2}
\,\widetilde{\bf L}^{\omega_2}_{\ell,q}\kappa\!\left(\omega_1,\omega_2;c_{*}(m)\right)\eta\,, 
\,\widetilde{\bf L}^{\omega_2}_{\ell,q}\theta\biggr\rangle_{N\backslash G} = \cr 
 &={1\over 4}\sum_{m=1}^M {S_{{\frak a},{\frak b}}\!\left(\omega_1,\omega_2;c_{*}(m)\right)\over 
\left| c_{*}(m)\right|^2}
\left\langle 
\widetilde{\bf L}^{\omega_2}_{\ell,q}\kappa\!\left(\omega_1,\omega_2;c_{*}(m)\right)\eta\,,  
\,\widetilde{\bf L}^{\omega_2}_{\ell,q}\theta\right\rangle_{N\backslash G} = \cr 
 &={1\over 4\pi i}  
\sum_{m=1}^M {S_{{\frak a},{\frak b}}\!\left(\omega_1,\omega_2;c_{*}(m)\right)\over 
\left| c_{*}(m)\right|^2}
\sum_{p=-\ell}^{\ell}\ \int\limits_{(0)} 
\bigl(\kappa\!\left(\omega_1,\omega_2;c_{*}(m)\bigr)\eta\right)\!(\nu,p)
\,\overline{\theta(\nu,p)}\,\lambda^{*}_{\ell}(\nu,p)\left( p^2-\nu^2\right) 
{\rm d}\nu = \cr 
 &=\sum_{m=1}^M {S_{{\frak a},{\frak b}}\!\left(\omega_1,\omega_2;c_{*}(m)\right)\over 
\left| c_{*}(m)\right|^2}\,\left( {\bf B}h_{\ell}\right)
\!\left( {2\pi\sqrt{\omega_1\omega_2}\over c_{*}(m)}\right) ,}$$
where $h_{\ell}(\nu,p)$, $\lambda^{*}_{\ell}(\nu,p)$ and the ${\bf B}$-transform 
are defined as in (6.6.3), (6.6.4) and (1.9.3)-(1.9.6);  
it therefore follows by Equation~(6.6.37) that  
$$\sum_{m=1}^{\infty} {S_{{\frak a},{\frak b}}\!\left(\omega_1,\omega_2;c_{*}(m)\right)\over 
\left| c_{*}(m)\right|^2}\,\left( {\bf B}h_{\ell}\right)
\!\left( {2\pi\sqrt{\omega_1\omega_2}\over c_{*}(m)}\right) 
={\left[\Gamma_{\frak a} : \Gamma_{\frak a}'\right]\over 4\pi^2} 
\int_{N\backslash G} F(g)\,\overline{\bigl(\widetilde{\bf L}^{\omega_2}_{\ell,q}\theta\bigr) 
(g)}\,{\rm d}\dot{g}\;.\eqno(6.6.38)$$
\par 
The right-hand side of Equation~(6.6.38) is 
independent of the choice of $c_{*}$. 
Moreover, given any permutation $\Lambda$ on ${\Bbb N}$, we may 
substitute $c_{*}\circ\Lambda$ for $c_{*}$ in the above: following any such  
substitution the function $c_{*}$ will still be  
a one-to-one function with domain ${\Bbb N}$ and range 
${}^{\frak a}{\cal C}^{\frak b}$, and so the result (6.6.38) will remain valid. 
Consequently, in light of  
Riemann's theorem (Theorem~8.33 of [1]) on conditionally convergent series, 
it may be deduced that the series on the left-hand side of Equation~(6.6.38) is 
absolutely convergent.  The equation (6.6.38) therefore has the equivalent 
formulation 
$${\left[\Gamma_{\frak a} : \Gamma_{\frak a}'\right]\over 4\pi^2} 
\,\int_{N\backslash G} F(g)\,\overline{\bigl(\widetilde{\bf L}^{\omega_2}_{\ell,q}\theta\bigr) 
(g)}\,{\rm d}\dot{g} 
=\sum_{c\in{}^{\frak a}{\cal C}^{\frak b}}  
{S_{{\frak a},{\frak b}}\!\left(\omega_1,\omega_2;c\right)\over 
\left| c\right|^2}\,\left( {\bf B}h_{\ell}\right)
\!\left( {2\pi\sqrt{\omega_1\omega_2}\over c}\right) \;.\eqno(6.6.39)$$
By (6.6.31), (6.6.32) and (6.6.39) we obtain the final equality in (6.6.2), 
which is the desired geometric description of 
$\langle\phi_1,\phi_2\rangle_{\Gamma\backslash G}\ $ $\square$  

\bigskip

\noindent{\bf Part~II of the proof of Proposition~6.6.1: 
the spectral description of 
$\langle\phi_1 , \phi_2\rangle_{\Gamma\backslash G}$.}\quad  
The first equality of (6.6.2) remains to be proved. In the above we obtained 
the geometric description of $\langle\phi_1,\phi_2\rangle_{\Gamma\backslash G}$ 
by applying Lemma~6.6.2 directly to 
$\langle\phi_2,\phi_1\rangle_{\Gamma\backslash G}
=\overline{\langle\phi_1,\phi_2\rangle_{\Gamma\backslash G}}\,$ (see, in particular,  
(6.6.31)). In obtaining the spectral description 
of $\langle\phi_1,\phi_2\rangle_{\Gamma\backslash G}$ we shall 
apply Theorem~A for $f_1=\phi_1$, $f_2=\phi_2$, and, rather than   
applying Lemma~6.6.2 directly to $\langle\phi_2,\phi_1\rangle_{\Gamma\backslash G}\,$ 
(or even to $\langle\phi_1,\phi_2\rangle_{\Gamma\backslash G}$), we shall  
instead apply it to each of the terms occurring, when 
$f_1=\phi_1$ and $f_2=\phi_2$, on the right-hand side of the Parseval 
identity~(1.8.8). These applications of Lemma~6.6.2 are our first concern in 
what follows (the application of the Parseval identity is to be discussed later). 
Accordingly, we suppose now that 
$$\phi\in\left\{ 1\,,\,T_V\varphi_{\ell,q}\left(\nu_V,p_V\right)\,,
\,E^{\frak c}_{\ell,q}\left( it_{*},p_{*}\right)\right\} ,$$
where $V$ is any one of the irreducible `cuspidal' subspaces of 
$L^2(\Gamma\backslash G)$ occurring as a factor in the direct sum in (1.7.4)  
$\,$($\nu_V,p_V$ being the associated spectral parameters), while 
${\frak c}$ is a cusp contained in some given set of representatives 
${\frak C}(\Gamma)\,$  of the $\Gamma$-equivalence classes 
of cusps, the Eisenstein series 
$E^{\frak c}_{\ell,q}(\nu,p)$ is as defined in 
Subsection~1.8, and one has  
$t_{*}\in{\Bbb R}$ and $p_{*}\in[-\ell,\ell]\cap(1/2)[\Gamma_{\frak c} : 
\Gamma_{\frak c}']{\Bbb Z}$. 
\par 
In applying Lemma~6.6.2 to $\langle\phi_1,\phi\rangle_{\Gamma\backslash G} 
=\langle P^{{\frak a},*}\widetilde{\bf L}^{\omega_1}_{\ell,q}\eta , 
\phi\rangle_{\Gamma\backslash G}$ we are led to consider the term 
$F^{\frak a}_{\omega_1}\phi$, from the Fourier expansion 
of $\phi$ at the cusp ${\frak a}$. With regard to the case $\phi =1$, 
one finds by (1.4.2)-(1.4.3) that 
$$F^{\frak a}_{\omega_1}1=0\;.\eqno(6.6.40)$$
Considering next the case $\phi =T_V\varphi_{\ell,q}(\nu_V,p_V)$ 
we observe that, in the discussion around (1.7.10)-(1.7.13), it was 
implicitly found that 
$F^{\frak a}_{\omega_1}T_V\varphi_{\ell,q}(\nu_V,p_V) 
=c_V^{\frak a}(\omega_1) {\bf J}_{\omega_1}\varphi_{\ell,q}(\nu_V,p_V)$; 
in terms of the modified Fourier coefficients defined in (1.7.15), this 
result becomes: 
$$F^{\frak a}_{\omega_1}T_V\varphi_{\ell,q}(\nu_V,p_V)  
=C_V^{\frak a}\left(\omega_1;\nu_V,p_V\right) 
\left(\pi\left|\omega_1\right|\right)^{-\nu_V}
\!\left(\omega_1 /\left|\omega_1\right|\right)^{p_V} 
\!{\bf J}_{\omega_1}\varphi_{\ell,q}(\nu_V,p_V)\;.\eqno(6.6.41)$$ 
Similarly, by what is noted in the discussion around (1.8.4)-(1.8.7) 
it may be  inferred that, for all $(\nu,g)\in{\Bbb C}\times G$ such that 
${\rm Re}(\nu)\geq 0$, one has 
$F^{\frak a}_{\omega_1} E^{\frak c}_{\ell,q}(\nu,p) 
=[\Gamma_{\frak c} : \Gamma_{\frak c}']^{-1} D_{\frak c}^{\frak a}(\omega_1;\nu,p) 
{\bf J}_{\omega_1}\varphi_{\ell,q}(\nu,p)$; in particular, one has  
$$F^{\frak a}_{\omega_1} E^{\frak c}_{\ell,q}\left( it_{*},p_{*}\right) 
=\left[\Gamma_{\frak c} : \Gamma_{\frak c}'\right]^{-1} 
\left(\pi\left|\omega_1\right|\right)^{-it_{*}} 
\left(\omega_1 /\left|\omega_1\right|\right)^{p_{*}} 
B_{\frak c}^{\frak a}\left(\omega_1;it_{*},p_{*}\right) 
{\bf J}_{\omega_1}\varphi_{\ell,q}\left( it_{*},p_{*}\right)\;,\eqno(6.6.42)$$ 
where the `modified Fourier coefficient' 
$B_{\frak c}^{\frak a}\left(\omega;\nu,p\right)$ is defined as 
Equation~(1.8.9) would indicate. 
\par 
By a calculation somewhat similar to that which yields the bound (6.6.30), 
it may be deduced from the estimates (6.5.15) of Lemma~6.5.3 (and from (1.5.16)) 
that, for $(\nu,p)\in{\Bbb C}\times\{-\ell,1-\ell,\ldots ,\ell\}$ and 
$\beta\in(1,\infty)$, one has: 
$$L^{1+\beta}(N\backslash G)\ni\bigl|{\bf J}_{\omega_1}\varphi_{\ell,q}(\nu,p)\bigr|\qquad\ 
\hbox{if $\ |{\rm Re}(\nu)|<1-2/(1+\beta)$.}\eqno(6.6.43)$$
By points noted in Subsection~1.7 (see in particular the paragraphs containing 
(1.7.5)-(1.7.8) and (1.7.14)), we have either 
$(\nu_V,p_V)\in (i{\Bbb R})\times\{-\ell,1-\ell,\ldots,\ell\}$, or else 
$p_V=0$ and $0\neq\nu_V\in(-1,1)$. In the former `principal series' case, it follows 
by (6.6.41) and (6.6.43) that one has 
$L^{1+\beta}(N\backslash G)\ni 
|F^{\frak a}_{\omega_1}T_V\varphi_{\ell,q}(\nu_V,p_V)|$ for all $\beta\in(1,\infty)$; 
in the latter `complementary series' case it follows by the same results 
that, for all real $\beta >(1+|\nu_V|)/(1-|\nu_V|)$, one has 
$L^{1+\beta}(N\backslash G)\ni 
|F^{\frak a}_{\omega_1}T_V\varphi_{\ell,q}(\nu_V,p_V)|$. By (6.6.42) and (6.6.43), one has also 
$L^{1+\beta}(N\backslash G)\ni 
|F^{\frak a}_{\omega_1} E^{\frak c}_{\ell,q}( it_{*},p_{*})|$ for 
all $\beta\in(1,\infty)$. Given the 
identity (6.6.40), and given the content of the 
last few observations (subsequent to (6.6.43)), 
it is certainly the case that 
$$L^{1+\beta}(N\backslash G)\ni 
\bigl|F^{\frak a}_{\omega_1}\phi\bigr|\qquad\ 
\hbox{for some $\,\beta=\beta(\phi)\in (0,\infty)$.}\eqno(6.6.44)$$ 
\par 
By (6.6.44), the function $\phi$ will satisfy all the relevant hypotheses of 
Lemma~6.6.2 if, for some $\alpha\in(0,\infty)$, it is contained in 
$C^0(G)\cap L^{1+\alpha}(\Gamma\backslash G)$. If $\phi =1$, then 
it is trivially the case that one has $\phi\in C^0(G)\cap L^{\infty}(\Gamma\backslash G)$. 
If $\phi =T_V\varphi_{\ell,q}(\nu_V,p_V)$ then by (1.7.3), (1.7.7), (1.7.8), 
(1.7.10) and the definitions (1.4.4)-(1.4.7) one has 
$\phi\in C^{\infty}(G)\cap{}^{0}L^{2}(\Gamma\backslash G)\subset 
C^0(G)\cap L^2(\Gamma\backslash G)$: in fact, by (1.7.10) and the bound 
(1.4.13) on the growth of any cusp form $f\in A^0_{\Gamma}(\Upsilon_{\nu,p};\ell,q)$, 
it follows that the $T_V\varphi_{\ell,q}(\nu_V,p_V)$ is a bounded function on $G$, so that 
one has also $\phi=T_V\varphi_{\ell,q}(\nu_V,p_V)\in L^{\infty}(\Gamma\backslash G)$. 
If $\phi=E^{\frak c}_{\ell,q}(it_{*},p_{*})$ then, 
by (1.3.2) and the `cusp-sector estimate' (6.2.20) of Lemma~6.2.8, 
it follows that there exists some $r_0=r_0(\Gamma,\ell,t_{*})\in[1,\infty)$  
such that, 
for all ${\frak d}\in{\Bbb Q}(i)\cup\{\infty\}$, and all 
$g\in G\,$ such that $\rho(g)\geq r_0$, 
one has  
$\phi(g_{\frak d}g)=O_{\Gamma,\ell,t_{*}}(\rho(g))\,$  
(it being assumed here that the scaling matrix $g_{\frak d}\in G$ is 
such that (1.1.16) and (1.1.20)-(1.1.21) hold when ${\frak d}$ is 
substituted for ${\frak c}$),  
and so in this case one finds (by computations similar to those seen 
in the proof of Corollary~6.2.10) that 
the space $L^{1+\alpha}(\Gamma\backslash G)$ 
contains $\phi$ whenever $\alpha$ satisfies $0<\alpha <1$. 
Moreover (as is asserted in Subsection~1.8) one has 
$E^{\frak c}_{\ell,q}(it_{*},p_{*})\in C^{\infty}(\Gamma\backslash G)$.  
This follows, by arguments similar to those employed in the proof of Lemma~6.5.15,  
from the analytic continuations of the 
summands occurring on the right-hand side of Equation~(1.8.4): note, in particular, 
that a bound such as the estimate~(11.49) of [32] enables one to establish that
the Fourier expansion of $E^{\frak c}_{\ell,q}(it_{*},p_{*})$ at the cusp 
$\infty$ is uniformly convergent on any compact subset of $G$. 
Given (6.5.7), and given what has so far been ascertained in the present paragraph, 
we may add to (6.6.44) the conclusion that 
$$\phi\in C^0(G)\cap L^{3/2}(\Gamma\backslash G)\;.\eqno(6.6.45)$$ 
\par 
By (6.6.44) and (6.6.45), Lemma~6.6.2 
applies when (as we assume) 
$\phi\in\{ 1 , T_V\varphi_{\ell,q}(\nu_V,p_V) , 
E^{\frak c}_{\ell,q}(it_{*},p_{*})\}$, and so, bearing in mind (6.6.40)-(6.6.42), 
it follows by the case $\omega =\omega_1$ of Equation~(6.6.5) that  
$$\eqalign{ 
 &\left[\Gamma_{\frak a} : \Gamma_{\frak a}'\right]  
\langle\phi_1 , \phi\rangle_{\Gamma\backslash G} = \cr 
 &\quad = \left[\Gamma_{\frak a} :\Gamma_{\frak a}'\right] 
\langle P^{{\frak a},*}\widetilde{\bf L}^{\omega_1}_{\ell,q}\eta , 
\phi\rangle_{\Gamma\backslash G}  = \cr 
 &\quad =\langle\widetilde{\bf L}^{\omega_1}_{\ell,q}\eta , 
F^{\frak a}_{\omega_1}\phi\rangle_{N\backslash G} ={\ \atop\matrix{\ }} \cr 
 &\quad =\cases{0 &if $\phi =1$; \cr 
\quad & \quad \cr 
\overline{C_V^{\frak a}\left(\omega_1;\nu_V,p_V\right)}
\ \bigl\langle\widetilde{\bf L}^{\omega_1}_{\ell,q}\eta \,,  
\,\left|\pi\omega_1\right|^{-\nu_V}
\!\left(\omega_1 /\left|\omega_1\right|\right)^{p_V} 
\!{\bf J}_{\omega_1}\varphi_{\ell,q}(\nu_V,p_V)\bigr\rangle_{N\backslash G} 
 &if $\phi =T_V\varphi_{\ell,q}(\nu_v,p_V)$; \cr 
\quad & \quad \cr 
\displaystyle{1\over\left[\Gamma_{\frak c} : \Gamma_{\frak c}'\right]} 
\,\overline{B_{\frak c}^{\frak a}\left(\omega_1;it_{*},p_{*}\right)}  
\ \bigl\langle\widetilde{\bf L}^{\omega_1}_{\ell,q}\eta \,,    
\,\left|\pi\omega_1\right|^{-it_{*}} 
\left(\omega_1 /\left|\omega_1\right|\right)^{p_{*}} 
{\bf J}_{\omega_1}\varphi_{\ell,q}\left( it_{*},p_{*}\right)\bigr\rangle_{N\backslash G} 
 &if $\phi =E^{\frak c}_{\ell,q}\left( it_{*},p_{*}\right)$.}
}$$
This result, when expressed in terms of the Lebedev transform operator 
${\bf L}^{\omega_1}_{\ell,q}\,$ (defined in (6.4.2)) becomes: 
$$\eqalign{ 
 &\left[\Gamma_{\frak a} : \Gamma_{\frak a}'\right]  
\langle\phi_1 , \phi\rangle_{\Gamma\backslash G} = {\ \atop\matrix{\ }}\cr 
 &\quad =\cases{ 
0 &if $\phi =1$; \cr 
\quad & \quad \cr 
\displaystyle\overline{C_V^{\frak a}\left(\omega_1;\nu_V,p_V\right)}
\ {(-i)^{p_V}\pi^2\left\|\Phi^{\ell}_{p_V,q}\right\|_K\over\Gamma\left(\ell +1+\overline{\nu_V}\right)} 
\,\bigl( {\bf L}^{\omega_1}_{\ell,q}\widetilde{\bf L}^{\omega_1}_{\ell,q}\eta\bigr)\!\left( 
-\overline{\nu_V} , p_V\right) 
 &if $\phi =T_V\varphi_{\ell,q}(\nu_v,p_V)$; \cr 
\quad & \quad \cr 
\displaystyle{1\over\left[\Gamma_{\frak c} : \Gamma_{\frak c}'\right]} 
\,\overline{B_{\frak c}^{\frak a}\left(\omega_1;it_{*},p_{*}\right)}  
\ {(-i)^{p_{*}}\pi^2\left\|\Phi^{\ell}_{p_{*},q}\right\|_K\over 
\Gamma\left(\ell +1+\overline{\left( it_{*}\right)}\right)} 
\,\bigl( {\bf L}^{\omega_1}_{\ell,q}\widetilde{\bf L}^{\omega_1}_{\ell,q}\eta\bigr)\!\left( 
-\overline{\left( it_{*}\right)} , p_{*}\right) 
 &if $\phi =E^{\frak c}_{\ell,q}\left( it_{*},p_{*}\right)$.}
}$$ 
Moreover, since $\eta\in{\cal T}^{\ell}_{\sigma}$, and since  
$\sigma\in(1,2)$, it follows by the results (6.4.7) and (6.4.10) of 
Theorem~6.4.1, and (1.6.5) and (1.6.6) (for $\ell'=\ell$, $q'=q$), that for  
$(\nu,p)\in((i{\Bbb R})\times\{-\ell,1-\ell,\ldots ,\ell\})\cup((-1,1)\times\{ 0\})$ 
one has:  
$$\eqalign{
{(-i)^{p}\pi^2\left\|\Phi^{\ell}_{p,q}\right\|_K\over\Gamma\left(\ell +1+\overline{\nu}\right)} 
\,\bigl( &{\bf L}^{\omega_1}_{\ell,q}\widetilde{\bf L}^{\omega_1}_{\ell,q}\eta\bigr)\!\left( 
-\overline{\nu} , p\right) = \cr 
 &=(-2\pi)(-i)^p\,{\sin(\pi\nu)\over(\pi\nu)} 
\,{\nu^{1+\epsilon(p)}\over\left(\nu^2 -p^2\right)} 
\,\eta\left(-\overline{\nu} , p\right)\ \times {\ \atop\matrix{\ \cr\ }}\cr 
 &\quad\ \times\cases{ 
\left\|\varphi_{\ell,q}(\nu,p)\right\|_{\rm ps}\Gamma(\ell +1+\nu) 
 &if $(\nu,p)\in (i{\Bbb R})\times\{ -\ell,1-\ell,\ldots ,\ell\}$;\cr 
\quad &\quad\cr 
\left\|\varphi_{\ell,q}(\nu,0)\right\|_{\rm cs} 
\sqrt{\Gamma(\ell +1+\nu)\Gamma(\ell +1-\nu)} 
 &if $(0,0)\neq(\nu,p)\in(-1,1)\times\{ 0\}$.}
}$$  
We therefore find that 
$${\left[\Gamma_{\frak a} : \Gamma_{\frak a}'\right]\over (-2\pi)} 
\,\langle\phi_1 , \phi\rangle_{\Gamma\backslash G} = 0\qquad\ 
\hbox{if $\,\phi =1$,}\eqno(6.6.46)$$ 
whereas if $\phi =T_V\varphi_{\ell,q}(\nu_V,p_V)$ then 
(when the norm $\|T_V\varphi_{\ell,q}(\nu_V,p_V)\|_{\Gamma\backslash G}$ is 
defined as in (1.7.14)) one has   
$$\eqalign{{\left[\Gamma_{\frak a} : \Gamma_{\frak a}'\right]\over (-2\pi)}   
\,\langle\phi_1 , \phi\rangle_{\Gamma\backslash G} 
 &=(-i)^{p_V}\left\|T_V\varphi_{\ell,q}(\nu_V,p_V)\right\|_{\Gamma\backslash G}
\,\overline{C_V^{\frak a}\left(\omega_1;\nu_V,p_V\right)} 
\ {\sin\left(\pi\nu_V\right)\over\left(\pi\nu_V\right)} 
\,{\nu_V^{1+\epsilon\left(p_V\right)}\over\left(\nu_V^2 -p_V^2\right)}\ \times{\ \atop\matrix{\ \cr\ }} \cr 
 &\quad\ \,\times\eta\left(-\overline{\nu_V}\,, p_V\right)\cdot\cases{ 
\Gamma\left(\ell +1+\nu_V\right) 
 &if $\,\nu_V^2\leq 0$,\cr 
\quad &\quad \cr 
\sqrt{\Gamma\left(\ell +1+\nu_V\right)\,\Gamma\left(\ell +1-\nu_V\right)} 
 &if $\,1>\nu_V^2>0=p_V$,}
}\eqno(6.6.47)$$
and if it is instead the case that $\phi =E^{\frak c}_{\ell,q}(it_{*},p_{*})$ then  
$$\eqalign{
 &{\left[\Gamma_{\frak a} : \Gamma_{\frak a}'\right]\over (-2\pi)}   
\,\langle\phi_1 , \phi\rangle_{\Gamma\backslash G} = \matrix{\ \cr\ \cr\ }\cr 
 &\quad\,={(-i)^{p_{*}}\left\|\varphi_{\ell,q}\left( it_{*} , p_{*}\right)\right\|_{\rm ps}\over 
\left[\Gamma_{\frak c} : \Gamma_{\frak c}'\right]} 
\,\overline{B^{\frak a}_{\frak c}\left(\omega_1;it_{*},p_{*}\right)} 
\ {\sin\left(\pi it_{*}\right)\over\left(\pi it_{*}\right)} 
\,{\left( it_{*}\right)^{1+\epsilon\left(p_{*}\right)}\over 
\bigl(\left( it_{*}\right)^2 -p_{*}^2\bigr)} 
\ \eta\bigl(-\overline{\left( it_{*}\right)}\,, p_{*}\bigr) 
\Gamma\left(\ell +1+it_{*}\right) .}\eqno(6.6.48)$$ 
\par 
By the symmetry apparent in our hypothesis (6.6.1), 
formulae corresponding to the above may be obtained for 
$(-2\pi)^{-1}[\Gamma_{\frak b} : \Gamma_{\frak b}']\langle\phi_2 , \phi\rangle_{\Gamma\backslash G} 
=(-2\pi)^{-1}[\Gamma_{\frak b} : \Gamma_{\frak b}']\langle 
P^{{\frak b},*}\widetilde{\bf L}^{\omega_1}_{\ell,q}\theta , 
\phi\rangle_{\Gamma\backslash G}$. Consequently, if we put 
$$F_{\phi_1,\phi_2}(\phi)={\left[\Gamma_{\frak a} : \Gamma_{\frak a}'\right] 
\left[\Gamma_{\frak b} : \Gamma_{\frak b}'\right]\over 4\pi^2}   
\,\langle\phi_1 , \phi\rangle_{\Gamma\backslash G} 
\,\langle\phi , \phi_2\rangle_{\Gamma\backslash G} 
=\overline{\biggl({\left[\Gamma_{\frak b} : \Gamma_{\frak b}'\right]\over (-2\pi)}
\,\langle\phi_2 , \phi\rangle_{\Gamma\backslash G}\biggr)}
\,{\left[\Gamma_{\frak a} : \Gamma_{\frak a}'\right]\over (-2\pi)}   
\,\langle\phi_1 , \phi\rangle_{\Gamma\backslash G}\;,\eqno(6.6.49)$$
and take $h_{\ell} : \{\nu\in{\Bbb C} : |{\rm Re}(\nu)|\leq\sigma\}\times{\Bbb Z}\rightarrow{\Bbb C}$ 
to be given by (6.6.3) and~(6.6.4), then, by (6.6.46), (6.6.47) and (6.6.48),  
and the corresponding formulae 
for $(-2\pi)^{-1}[\Gamma_{\frak b} : \Gamma_{\frak b}']\langle\phi_2 , 
\phi\rangle_{\Gamma\backslash G}$, it follows that 
$$F_{\phi_1,\phi_2}(\phi) 
=\cases{ 
0 &if $\,\phi =1$; \cr 
\quad &\quad \cr 
\left\|T_V\varphi_{\ell,q}(\nu_V,p_V)\right\|_{\Gamma\backslash G}^2 
\,\overline{C_V^{\frak a}\left(\omega_1;\nu_V,p_V\right)} 
\ C_V^{\frak b}\left(\omega_2;\nu_V,p_V\right)
h_{\ell}\!\left(\nu_V,p_V\right)  
 &if $\,\phi =T_V\varphi_{\ell,q}\left(\nu_V,p_V\right)$;\cr 
\quad &\quad \cr 
\displaystyle{\left\|\varphi_{\ell,q}\left( it_{*} , p_{*}\right)\right\|_{\rm ps}^2\over 
\left[\Gamma_{\frak c} : \Gamma_{\frak c}'\right]^2} 
\,\overline{B^{\frak a}_{\frak c}\left(\omega_1;it_{*},p_{*}\right)} 
\ B^{\frak b}_{\frak c}\left(\omega_2;it_{*},p_{*}\right) 
h_{\ell}\!\left( it_{*} , p_{*}\right) 
 &if $\phi =E^{\frak c}_{\ell,q}\left( it_{*} , p_{*}\right)$.}\eqno(6.6.50)$$
\bigskip\noindent 
In arriving at the results stated in (6.6.50) we have used both the fact that   
for ${\rm Re}(\nu)>-1$ one has 
$\Gamma(\ell + 1 +\overline{\nu})=\overline{\Gamma(\ell+1+\nu)}$, and 
the fact that the  
function $\alpha : {\Bbb C}\times{\Bbb Z}\rightarrow{\Bbb C}$ given by 
$$\alpha(\nu,p) 
=\ {\sin\left(\pi\nu\right)\over\left(\pi\nu\right)} 
\,{\nu^{1+\epsilon\left(p\right)}\over\left(\nu^2 -p^2\right)}
={(-1)^p\over\Gamma\left( 1+\nu+|p|\right)\Gamma\left( 1-\nu+|p|\right)} 
\,\prod_{1\leq m<|p|}\left(\nu^2 -m^2\right)\qquad\quad 
\hbox{($\nu\in{\Bbb C}$, $p\in{\Bbb Z}$)}$$
satisfies both $\alpha(\overline{\nu},p)=\overline{\alpha(\nu ,p)}$ and 
$\alpha(-\nu,p)=\alpha(\nu,p)$, for all $(\nu,p)\in{\Bbb C}\times{\Bbb Z}$.  
We have also made use of both the observation that 
if $(\nu,p)\in{\Bbb C}\times{\Bbb Z}$ 
and $\nu^2\leq 0$ then $\overline{\nu}=-\nu$, and the (complementary)  
observation that if it is instead  
the case that $p=0$ and $1>\nu^2>0$ then 
$-\overline{\nu}=-\nu$, $p=-p$,    
$\sqrt{\Gamma\left(\ell+1+\nu\right)\Gamma\left(\ell+1-\nu\right)}\in{\Bbb R}$ 
and $\eta(-\nu,-p)=\eta(\nu,p)\,$ 
(the last equality following by virtue of the condition 
(T1) stated below (6.4.3), and our hypothesis 
that $\eta\in{\cal T}^{\ell}_{\sigma}\,$). 
\par 
Observe now that, by (6.6.1), Lemma~6.5.17 and (6.5.7), the hypotheses 
of Theorem~A are satisfied when one has (there) $f_1=\phi_1$ and $f_2=\phi_2$. 
It is therefore implied by Theorem~A that, for $f_1=\phi_1$ and $f_2=\phi_2$, 
the `Parseval identity' stated in (1.8.8) is valid. By the case $f_1=\phi_1$, 
$f_2=\phi_2$ of (1.8.8), combined with (6.6.49)-(6.6.50), we obtain the 
first inequality in (6.6.2), which is the desired spectral description of 
$\langle\phi_1,\phi_2\rangle_{\Gamma\backslash G}$. 
The sums and integrals appearing on the left-hand side of the first equality 
in~(6.6.2) are simply an alternative formulation of the sums and integrals 
which (in the case $f_1=\phi_1$, $f_2=\phi_2$) appear on the right-hand side 
of Equation~(1.8.8); it is therefore a corollary of Theorem~A 
that these sums and integrals are absolutely convergent. This completes our 
proof of the preliminary sum formula, Proposition~6.6.1\quad$\blacksquare$ 

\bigskip 

\centerline{\bf \S 6.7 Completing the proof of the spectral summation formula.} 

\medskip 

In this final subsection of 
our Appendix we show that 
the preliminary sum formula (Proposition~6.6.1) implies the more general 
result asserted  in Theorem~B. Our proof of this is closely modelled on 
the `Extension Method' employed by Lokvenec-Guleska in 
[32], Subsection~11.2 and Subsection~11.3. 
\par 
It is to be assumed henceforth  
that $\omega_1$, $\omega_2$, ${\frak a}$, ${\frak b}$, $g_{\frak a}$ and $g_{\frak b}$ 
are given, with $0\neq\omega_1,\omega_2\in{\frak O}$, and with  
${\frak a},{\frak b}\in{\Bbb Q}(i)\cup\{\infty\}$ and $g_{\frak b},g_{\frak b}\in G$ 
such that (1.1.16) and (1.1.20)-(1.1.21) hold for ${\frak c}\in\{ {\frak a},{\frak b}\}$. 
We begin our implementation of   
Lokvenec-Guleska's `Extension Method' by defining  
some of the relevant terminology. 
\par 
When $\sigma\in(0,\infty)$ and $\varrho,\vartheta \in{\Bbb R}$, let ${\cal H}_0^{\sigma}(\varrho,\vartheta )$ 
be the space of functions 
$h : \{\nu\in{\Bbb C} : |{\rm Re}(\nu)|\leq\sigma\}\times{\Bbb Z}\rightarrow{\Bbb C}$ 
satisfying, for the given choice of $\sigma$, $\varrho$ and $\vartheta$, 
the conditions~(i)-(iii) of Theorem~B. 
\par 
For $\sigma\in(1/2,1)$, $\varrho\in(2,\infty)$ and $\vartheta\in(3,\infty)$, 
we define ${\cal H}_{\star}^{\sigma}(\varrho,\vartheta )$ to be equal to the space ${\cal H}_0^{\sigma}(\varrho,\vartheta )$; 
for $\sigma\in(1,2)$, $\varrho\in(2,\infty)$ and $\vartheta\in(3,\infty)$, we define  
${\cal H}_{\star}^{\sigma}(\varrho,\vartheta )$ to be the subspace of ${\cal H}_0^{\sigma}(\varrho,\vartheta )$ 
containing just those members $h$ of ${\cal H}_0^{\sigma}(\varrho,\vartheta )$ 
that (in addition to satisfying the  
conditions (i)-(iii) of Theorem~B) are such that, for each $p\in{\Bbb Z}$, the 
function $\nu\mapsto h(\nu,p)$ has a zero of order at least $2$ at  
the point $\nu=1$. 
\par 
For $\sigma\in(1/2,1)\cup(1,2)$, $\varrho,\vartheta \in(3,\infty)$, $h\in{\cal H}_{\star}^{\sigma}(\varrho,\vartheta )$, 
$\omega,\omega'\in\{\omega_1,\omega_2\}$ and ${\frak d},{\frak d}'\in\{ 
{\frak a},{\frak b}\}$, we put 
$$\eqalignno{
\chi^{{\frak d},{\frak d}'}_{\omega,\omega'}(h) 
 &={1\over 4\pi^3 i}\,\delta^{{\frak d},{\frak d}'}_{\omega,\omega'}\,
\sum_{p\in{\Bbb Z}}\ \int\limits_{(0)}^{\hbox{\ }}h(\nu,p) 
\left( p^2 -\nu^2\right){\rm d}\nu\;, &(6.7.1)\cr  
X^{{\frak d},{\frak d}'}_{\omega,\omega'}(h) 
 &=\sum_{c\in {}^{\frak d}{\cal C}^{{\frak d}'}} 
\,{S_{{\frak d},{\frak d}'}\left(\omega , \omega' ; c\right)\over |c|^2}\,
\left({\bf B}h\right)
\!\left( {2\pi\sqrt{\omega\omega'}\over c}\right) &(6.7.2)}$$
(with ${\bf B}h : {\Bbb C}^{*}\rightarrow{\Bbb C}$ as 
defined in (1.9.3)-(1.9.6)) and, 
subject to the absolute convergence of the relevant sums and integrals, 
$$\eqalign{ 
Y^{{\frak d},{\frak d}'}_{\omega,\omega'}(h) 
 &=\sum_{V}\,\overline{C_V^{\frak d}\left(\omega ;\nu_V,p_V\right)}
\,C_V^{{\frak d}'}\left(\omega' ;\nu_V,p_V\right) h\left(\nu_V , p_V\right)\,+ \cr 
 &\quad\ +\sum_{{\frak c}\in{\frak C}(\Gamma)}
{1\over 4\pi i\left[\Gamma_{\frak c} : \Gamma_{\frak c}'\right]} 
\sum_{p\in{1\over 2}\left[\Gamma_{\frak c} : 
\Gamma_{\frak c}'\right]{\Bbb Z}} 
\ \int\limits_{(0)} \overline{B_{\frak c}^{\frak d}\left(\omega;\nu,p\right)}\,
B_{\frak c}^{{\frak d}'}\left(\omega' ;\nu,p\right) h(\nu,p)\,{\rm d}\nu\;.}\eqno(6.7.3)$$ 

By the phrase `the sum formula for $Y^{{\frak d},{\frak d}'}_{\omega,\omega'}(h)$ 
is valid', we shall mean that all sums and integrals on the right-hand side 
of Equation~(6.7.3) are absolutely convergent, and that one has 
$$Y^{{\frak d},{\frak d}'}_{\omega,\omega'}(h) 
=\chi^{{\frak d},{\frak d}'}_{\omega,\omega'}(h) 
+X^{{\frak d},{\frak d}'}_{\omega,\omega'}(h)\;.\eqno(6.7.4)$$

\par 
\noindent{\bf Remark~6.7.1.}\quad 
Let $\sigma$, $\varrho$, $\vartheta$ and $h$ satisfy the hypotheses 
assumed in defining $\chi^{{\frak d},{\frak d}'}_{\omega,\omega'}(h)$, 
$X^{{\frak d},{\frak d}'}_{\omega,\omega'}(h)$ and 
$Y^{{\frak d},{\frak d}'}_{\omega,\omega'}(h)$. 
Then $\varrho >3$, $\vartheta >3$ and  
$h$, $\sigma$, $\varrho$ and $\vartheta$ are, in particular, such that the 
conditions (ii) and (iii) 
of Theorem~B are satisfied; therefore the integrals and sum 
on the right-hand side of Equation~(6.7.1) are absolutely convergent. 
It moreover follows by [32], Lemma~11.1.1 and Lemma~11.1.2, that 
the transform ${\bf B}h : {\Bbb C}^{*}\rightarrow{\Bbb C}$ is a well-defined function 
(by virtue of the integrals and sum on the right-hand 
side of Equation~(1.9.3) being absolutely convergent for all $u\in{\Bbb C}^{*}$), 
and that this transform ${\bf B}h$ satisfies  
$$\sup\,\bigl\{\,|u|^{-2\min\{\sigma ,1\}}|({\bf B}h)(u)|\,: 
\,u\in{\Bbb C}\ \,{\rm and}\ \,0<|u|\leq r_1\bigr\}<\infty\qquad\  
\hbox{for some $\,r_1\in (0,\infty)$.}\eqno(6.7.5)$$ 
Given our assumption of the hypothesis that $\sigma\in(1/2,1)\cup(1,2)$, 
it follows by (6.7.5), Lemma~6.5.9 and the result (6.1.26) of Lemma~6.1.5  
that the sum over $c\in{}^{\frak d}{\cal C}^{{\frak d}'}$ occurring in (6.7.2)  
is absolutely convergent.  Note too that 
by Remark~1.9.2 it follows that 
when $h\in{\cal H}_{\star}^{\sigma}(\varrho,\vartheta )$ and $\sigma >1/2$  it is then the case that 
all the summands 
of the first sum on the right-hand side of Equation~(6.7.3) are defined.  
The absolute convergence of all the sums and integrals occurring 
on the right-hand side of Equation~(6.7.3) will be established later on,  
in the course of our proof of Theorem~B; in the meantime it cannot be 
taken for granted.

\bigskip 

\par 
We divide our proof of Theorem~B into two principal stages. In the first of these stages  
we deduce from Proposition~6.6.1 (the Preliminary Sum Formula) the following result. 

\bigskip 

\proclaim Proposition 6.7.2 (weak sum formula). Let $\sigma\in(1,2)$, let 
$\varrho,\vartheta \in(3,\infty)$, and let $f\in{\cal H}_{\star}^{\sigma}(\varrho,\vartheta )$. Then the 
sum formula for $Y^{{\frak a},{\frak b}}_{\omega_1,\omega_2}(f)$ is valid.  

\bigskip 

Our proof of Proposition~6.7.2 is modelled on the initial steps in 
Lokvenec-Guleska's `Extension Method' of [32], 
Subsection~11.2 and Subsection~11.3; 
we approach it via nine lemmas, two of which (Lemma~6.7.9 and Lemma~6.7.11) are 
used again in the final stage of our proof of Theorem~B. 
We omit the proofs of Lemma~6.7.7 and Lemma~6.7.8, 
which differ from the proofs of Lemma~6.7.4 and Lemma~6.7.5 (respectively) only 
in that they involve the application of Lebesgue's `dominated convergence' theorem 
(Theorem~1.34 of [40]), whereas it is Lebesgue's `monotone convergence' theorem  
that is applied in the proofs of Lemma~6.7.4 and Lemma~6.7.5.

\bigskip 

\proclaim Lemma 6.7.3. Let $\sigma\in(0,1)\cup(1,2)$, $\varrho\in(2,\infty)$ and 
$\vartheta\in(3,\infty)$,  
let $f\in{\cal H}_{\star}^{\sigma}(\varrho,\vartheta )$, and let $(f_n)_{n\in{\Bbb N}}$ be a  
sequence of elements of ${\cal H}_{\star}^{\sigma}(\varrho,\vartheta )$ converging pointwise to 
$f$ and satisfying 
$$\sup\,\bigl\{ \left( 1+\left|{\rm Im}(\nu)\right|\right)^{\varrho} 
\left( 1+|p|\right)^{\vartheta} \left| f_n(\nu,p)\right|\,: 
\,n\in{\Bbb N} ,\,p\in{\Bbb Z} ,\,\nu\in{\Bbb C}\ \,{\rm and}
\ \,{\rm Re}(\nu)=\sigma\bigr\} 
<\infty\;.\eqno(6.7.6)$$ 
Then the integrals and sum defining the transform ${\bf B}f$ are 
absolutely convergent: each of the transforms in the sequence 
${\bf B}f,{\bf B}f_1,{\bf B}f_2,{\bf B}f_3\ldots $ is a complex valued function with 
domain ${\Bbb C}^{*}$. One has, moreover, 
$$\lim_{n\rightarrow\infty}\left( 
\sup\left\{ |u|^{-2\min\{\sigma,1\}}
\left|\left( {\bf B}f_n - {\bf B}f\right)\!(u)\right| \,:\, 
u\in{\Bbb C}\ \,{\rm and}\ \,0<|u|\leq r\right\}\right) 
=0\qquad\hbox{for each $\,r\in(0,\infty)$.}\eqno(6.7.7)$$ 

\medskip 

\noindent{\bf Proof.}\quad 
This lemma is a minor refinement on Lemma~11.1.3 of [32]. 
Observe firstly that, since $f$ and each function in 
the series $(f_n)_{n\in{\Bbb N}}$ lies in the space 
${\cal H}_{\star}^{\sigma}(\varrho,\vartheta )$, and since we have (6.7.6), there therefore 
exists some $C_0\in[0,\infty)$ such that  
$$\left|\left( f_n - f\right)(\nu,p)\right|\leq {C_0\over  
\left( 1+\left|{\rm Im}(\nu)\right|\right)^{\varrho} 
\left( 1+|p|\right)^{\vartheta}}\quad\ \,  
\hbox{for all 
$\,(\nu,p,n)\in{\Bbb C}\times{\Bbb Z}\times{\Bbb N}\,$ 
such that $\,|{\rm Re}(\nu)|=\sigma$}
\eqno(6.7.8)$$ 
(i.e. this follows since our hypotheses imply that the condition (iii) of Theorem~B 
is satisfied when $h=f$, and that the condition (i) of Theorem~B is satisfied 
when $h\in\{ f \}\cup\{ f_n : n\in{\Bbb N}\}$). Our hypotheses 
imply, moreover, that the conditions (ii) and (iii) 
of Theorem~B are satisfied when $h\in\{ f \}\cup\{ f_n : n\in{\Bbb N}\}$ and 
$\sigma$, $\varrho$ and $\vartheta$ are as given (with, in particular, $\varrho >2>0$, so that 
one has $(f_n-f)(\nu,p)\rightarrow 0$ as $|{\rm Im}(\nu)|\rightarrow\infty$ with 
$n$ and $p$ fixed and $\nu$ constrained to satisfy $|{\rm Re}(\nu)|\leq\sigma$). 
Therefore, by an application of the maximum principle for analytic functions, 
one may deduce from (6.7.8) that, 
for all $n\in{\Bbb N}$ and all $p\in{\Bbb Z}$, one has 
$$\eqalign{\sup\,\bigl\{\,\left|\left( f_n - f\right)(\nu,p)\right|\,:
\,\nu\in{\Bbb C}\ \,{\rm and}\ \,|{\rm Re}(\nu)|\leq\sigma\bigr\} 
 &=\sup\,\bigl\{ \,\left|\left( f_n - f\right)(\nu,p)\right|\,:
\,\nu\in{\Bbb C}\ \,{\rm and}\ \,|{\rm Re}(\nu)|=\sigma\bigr\} \leq \cr 
 &\leq C_0\left( 1+|p|\right)^{-\vartheta}\;.}$$ 
Hence, in addition to (6.7.8), one has  
$\left|\left( f_n - f\right)(0,p)\right|\leq C_0 (1+|p|)^{-\vartheta}$ for 
all $n\in{\Bbb N}$ and all $p\in{\Bbb Z}$. Therefore our hypothesis~(6.7.6), 
despite being weaker than the corresponding hypothesis in Lemma~11.1.3 of [32], 
does in fact imply all of the assumptions that are relied upon in the 
(sketchy, but valid) proof supplied, in [32], for that lemma, 
and so (by the steps indicated on Page~100 of [32]) the result~(6.7.7) follows.  
The results stated between (6.7.6) and (6.7.7) follow immediately 
from Lemma~11.1.1 of~[32]\quad$\blacksquare$ 

\bigskip   

\proclaim Lemma 6.7.4. Let $(F_n)_{n\in{\Bbb N}}$ be a sequence of 
mappings from ${\Bbb N}$ into ${\Bbb R}\cup\{\infty\}$. Suppose moreover that 
$$0\leq F_1(m)\leq F_2(m)\leq\ldots\qquad\ \hbox{for each $m\in{\Bbb N}$.}\eqno(6.7.9)$$ 
Then, for some $\lambda\in{\Bbb R}\cup\{\infty\}$, one has  
$$\sum_{m=1}^{\infty}\Bigl(\,\lim_{n\rightarrow\infty} F_n(m)\Bigr) 
=\lambda =\lim_{n\rightarrow\infty}\biggl(\,\sum_{m=1}^{\infty} F_n(m)\biggr)\;.\eqno(6.7.10)$$ 

\medskip 

\noindent{\bf Proof.}\quad 
Let ${\bf P}({\Bbb N})$ be the set of all subsets of ${\Bbb N}$. 
Define the function $\mu_{\Bbb Z} : {\bf P}({\Bbb N})\rightarrow [0,\infty]$ by 
setting $\mu_{\Bbb Z}(A)$ equal to the cardinality of $A$ whenever ${\Bbb N}\supseteq A$. 
Then $\mu_{\Bbb Z}$ is a positive measure on ${\Bbb Z}$, and for $n\in{\Bbb N}$ one has 
$\sum_{m=1}^{\infty} F_n(m)=\int_{\Bbb Z} F_n {\rm d}\mu_{\Bbb Z}\,$ (the expression on the right-hand side 
of this equation denoting the Lebesgue integral of $F_n$ over ${\Bbb Z}$, with 
respect to the measure $\mu_{\Bbb Z}$). The result (6.7.10) is therefore 
equivalent to a special case of Lebesgue's `monotone convergence' theorem, 
Theorem~1.26 of [40]\quad$\blacksquare$ 

\bigskip 

\proclaim Lemma 6.7.5. Let $({\Phi}_n)_{n\in{\Bbb N}}$ be a sequence of mappings  
from ${\Bbb R}\times{\Bbb Z}$ into ${\Bbb R}\cup\{\infty\}$. 
Suppose moreover that 
$$0\leq {\Phi}_1(t,p)\leq {\Phi}_2(t,p)\leq\ldots\qquad\ 
\hbox{for each $(t,p)\in{\Bbb R}\times{\Bbb Z}$,}\eqno(6.7.11)$$ 
and that, for all $n\in{\Bbb N}$ and all $p\in{\Bbb Z}$, the 
mapping $t\mapsto {\Phi}_n(t,p)$ is a function on ${\Bbb R}$ that 
is measurable with respect to the Lebesgue measure on ${\Bbb R}$. 
Then, for some $\lambda\in{\Bbb R}\cup\{\infty\}$, one has 
$$\sum_{p=-\infty}^{\infty}\ \int\limits_{-\infty}^{\infty}
\Bigl(\,\lim_{n\rightarrow\infty} {\Phi}_n(t,p)\Bigr)\,{\rm d}t  
=\lambda 
=\lim_{n\rightarrow\infty}\biggl(\,\sum_{p=-\infty}^{\infty}
\ \int\limits_{-\infty}^{\infty} 
{\Phi}_n(t,p)\,{\rm d}t\biggr)\;.\eqno(6.7.12)$$ 

\medskip 

\noindent{\bf Proof.}\quad  
For $n\in{\Bbb N}$, one has 
$\sum_{p=-\infty}^{\infty}
\int_{-\infty}^{\infty} {\Phi}_n(t,p)\,{\rm d}t 
=\int_{{\Bbb R}\times{\Bbb Z}}{\Phi}_n {\rm d}\mu_{{\Bbb R}\times{\Bbb Z}}$, 
where $\mu_{{\Bbb R}\times{\Bbb Z}}$ is a positive measure on 
${\Bbb R}\times{\Bbb Z}$ 
(defined so that, for all $p\in{\Bbb Z}$, all $M\in[0,\infty]$ and all 
$A\subseteq{\Bbb R}$ such that $A$ has Lebesgue measure~$M$, one has 
$\mu_{{\Bbb R}\times{\Bbb Z}}(A\times\{ p\})=M$). 
Therefore, like the lemma which preceded it, this lemma is equivalent to a special case 
of Lebesgue's `monotone convergence' theorem\quad$\blacksquare$ 

\bigskip

\proclaim Lemma 6.7.6. Let $\sigma\in(1/2,1)\cup(1,2)$, $\varrho,\vartheta \in(3,\infty)$, 
${\frak d}\in\{ {\frak a} , {\frak b}\}$, $\omega\in\{\omega_1 , \omega_2\}$ 
and $f\in{\cal H}_{\star}^{\sigma}(\varrho,\vartheta )$. Suppose moreover that 
$$E=\bigl( (i{\Bbb R})\times{\Bbb Z} \bigr)\cup\bigl( [-2/9 , 2/9 ]\times\{ 0\}\bigr)\;,\eqno(6.7.13)$$ 
and that $(f_n)_{n\in{\Bbb N}}$ is a sequence of elements of 
${\cal H}_{\star}^{\sigma}(\varrho,\vartheta )$ which satisfies the condition (6.7.6) of Lemma~6.7.3 and 
the following conditions:\smallskip 
\hskip 15mm (i)\quad for all $n\in{\Bbb N}$, the sum formula for 
$Y^{{\frak d},{\frak d}}_{\omega,\omega}(f_n)$ is valid;\hfill\smallskip 
\hskip 14.5mm (ii)\quad for all $(\nu,p)\in E$ one has 
$0\leq f_1(\nu,p)\leq f_2(\nu,p)\leq\ldots\ \,$;\smallskip 
\hskip 14mm (iii)\quad for all $(\nu,p)\in{\Bbb C}\times{\Bbb Z}$ such that 
$\,|{\rm Re}(\nu)|\leq\sigma$, one has 
$\,\lim_{n\rightarrow\infty} f_n(\nu,p)=f(\nu,p)$. \medskip 
\noindent Then  the sum formula for $Y^{{\frak d},{\frak d}}_{\omega,\omega}(f)$ 
is valid, and one has 
$$\eqalignno{
\chi^{{\frak d},{\frak d}}_{\omega,\omega}(f) 
 &=\lim_{n\rightarrow\infty}\chi^{{\frak d},{\frak d}}_{\omega,\omega}(f_n)\;, &(6.7.14) \cr 
X^{{\frak d},{\frak d}}_{\omega,\omega}(f) 
 &=\lim_{n\rightarrow\infty}X^{{\frak d},{\frak d}}_{\omega,\omega}(f_n)\;, &(6.7.15) \cr 
Y^{{\frak d},{\frak d}}_{\omega,\omega}(f) 
 &=\lim_{n\rightarrow\infty}Y^{{\frak d},{\frak d}}_{\omega,\omega}(f_n)\;. &(6.7.16)}$$ 

\medskip 

\noindent{\bf Proof.}\quad 
Given the definition of $\chi^{{\frak d},{\frak d}'}_{\omega,\omega'}(h)$ in (6.7.1), 
and given (1.9.2) (where the 
notation `$\delta^{{\frak a},{\frak b}}_{\omega_1,\omega_2}$' is defined), 
it follows by (1.1.20) and Lemma~4.2 that, 
for $h\in{\cal H}_{\star}^{\sigma}(\varrho,\vartheta )$ we have 
$$\chi^{{\frak d},{\frak d}}_{\omega,\omega}(h) 
={1\over 2\pi^3}\sum_{p\in{\Bbb Z}}\ \int\limits_{-\infty}^{\infty} 
h(it,p)\left( t^2 + p^2\right) {\rm d}t\qquad\ 
{\rm and}\qquad\ \chi^{{\frak d},{\frak d}}_{\omega,\omega}(h)\in{\Bbb C}\eqno(6.7.17)$$ 
(the integrals and sum on the right-hand side of this equation being 
absolutely convergent, as discussed in our Remark~6.7.1).  
In particular, the hypotheses of the lemma suffice to ensure that (6.7.17) holds 
for all $h\in\{ f_n : n\in{\Bbb N}\}\cup\{ f\}$. Therefore, since the 
form $t^2+p^2$ is positive definite, and since (by hypothesis)    
the conditions (ii) and (iii) of the lemma are satisfied,
it follows that (6.7.11) and all the other hypotheses of Lemma~6.7.5 are satisfied 
if $({\Phi}_n)_{n\in{\Bbb N}}=({\Phi}_n(t,p))_{n\in{\Bbb N}}=((t^2+p^2) f_n(it,p))_{n\in{\Bbb N}}$.   
Hence the result (6.7.12) of Lemma~6.7.5 implies the equality in (6.7.14).
\par 
Given the definition of $X^{{\frak d},{\frak d}'}_{\omega,\omega'}(h)$ in (6.7.2) 
and what is observed in our Remark~6.7.1, 
and given our hypotheses concerning $f$ and the sequence $(f_n)$, 
it follows that 
$$X^{{\frak d},{\frak d}}_{\omega,\omega}(h)\in{\Bbb C}\qquad\    
\hbox{for all $\ h\in\{ f_n : n\in{\Bbb N}\}\cup\{ f\}$,}\eqno(6.7.18)$$ 
and that, by virtue of the result (6.1.26) of Lemma~6.1.5,  one has, 
for all $n\in{\Bbb N}$, 
$$\eqalign{ 
 &\bigl| X^{{\frak d},{\frak d}}_{\omega,\omega}(f) 
-X^{{\frak d},{\frak d}}_{\omega,\omega}\left( f_n\right)\bigr| \leq \cr 
 &\qquad\qquad\ \leq\sum_{c\in {}^{\frak d}{\cal C}^{{\frak d}}} 
\,{\left| S_{{\frak d},{\frak d}}\left(\omega , \omega ; c\right)\right|\over |c|^2}\,
\left|\left( {\bf B}f -{\bf B}f_n\right) 
\!\left( {2\pi\omega\over c}\right)\right| \leq \cr 
 &\qquad\qquad\ \leq |2\pi\omega|^{2\min\{ 1 , \sigma\}}  
\left(\sup\left\{ {\left|\left( {\bf B}f -{\bf B}f_n\right)(u)\right|\over 
|u|^{2\min\{\sigma , 1\}}} 
\,:\, u\in{\Bbb C}\ {\rm and}\ 0<|u|\leq{2\pi |\omega|\over\left| m_{\frak d}\right|}\right\}\right)  
\sum_{c\in {}^{\frak d}{\cal C}^{{\frak d}}} 
\,{\left| S_{{\frak d},{\frak d}}\left(\omega , \omega ; c\right)\right|\over 
|c|^{2\min\{1+\sigma , 2\}}}\;.
}$$ 
Therefore, since $\sigma >1/2$, and since the sequence 
$(f_n)_{n\in{\Bbb N}}$ satisfies the relevant hypotheses, it follows by 
Lemma~6.5.9 and Lemma~6.7.3 that the equality in (6.7.15) must hold. 
\par 
Now that (6.7.14) and (6.7.15) have been proved 
(and given that it therefore follows by (6.7.17) and (6.7.18) that the relevant 
sequences, $(\chi^{{\frak d},{\frak d}}_{\omega,\omega}( f_n))_{n\in{\Bbb N}}$ and   
$(X^{{\frak d},{\frak d}}_{\omega,\omega}( f_n))_{n\in{\Bbb N}}$, 
are convergent sequences of complex numbers), it  
will suffice for the completion of the proof of the lemma that we show that the equality 
in (6.7.16) holds: for, when that is achieved, 
it may then be inferred from (6.7.14)-(6.7.16) that, since 
the condition (i) of the lemma is satisfied, one has 
$$\chi^{{\frak d},{\frak d}}_{\omega,\omega}(f) 
+X^{{\frak d},{\frak d}}_{\omega,\omega}(f) 
-Y^{{\frak d},{\frak d}}_{\omega,\omega}(f) 
=\lim_{n\rightarrow\infty}
\left( \chi^{{\frak d},{\frak d}}_{\omega,\omega}\left( f_n\right)  
+X^{{\frak d},{\frak d}}_{\omega,\omega}\left( f_n\right)  
-Y^{{\frak d},{\frak d}}_{\omega,\omega}\left( f_n\right)\right)
=\lim_{n\rightarrow\infty} 0 =0\;,\eqno(6.7.19)$$ 
so that the sum formula for 
$Y^{{\frak d},{\frak d}}_{\omega,\omega}(f)$ is valid 
(no proof of (6.7.16) being complete without it 
having been shown that 
all the sums and integrals occurring on the right-hand side of 
Equation~(6.7.3) are absolutely convergent 
when $h=f$, ${\frak d}'={\frak d}$ and $\omega'=\omega$).  
\par 
Given that the condition (i) of the lemma is satisfied, 
it follows by (6.7.14), (6.7.15), (6.7.18) and the case $h=f$ of (6.7.17) 
that we have now 
$${\Bbb C}\ni \chi^{{\frak d},{\frak d}}_{\omega,\omega}(f) 
+X^{{\frak d},{\frak d}}_{\omega,\omega}(f) 
=\lim_{n\rightarrow\infty}
\left( \chi^{{\frak d},{\frak d}}_{\omega,\omega}\left( f_n\right)  
+X^{{\frak d},{\frak d}}_{\omega,\omega}\left( f_n\right)\right)  
=\lim_{n\rightarrow\infty} 
Y^{{\frak d},{\frak d}}_{\omega,\omega}\left( f_n\right) .\eqno(6.7.20)$$  
Moreover, by the definition (6.7.3) of 
$Y^{{\frak d},{\frak d}'}_{\omega,\omega'}(h)$, 
we have also 
$$\eqalign{ 
 &Y^{{\frak d},{\frak d}}_{\omega,\omega}(h) = \cr 
 &\quad\ =\sum_{V}\,\bigl| C_V^{\frak d}\left(\omega ;\nu_V,p_V\right)\bigr|^2 
h\left(\nu_V , p_V\right) 
+\!\sum_{{\frak c}\in{\frak C}(\Gamma)}
{1\over 4\pi\left[\Gamma_{\frak c} : \Gamma_{\frak c}'\right]} 
\sum_{p\in{1\over 2}\left[\Gamma_{\frak c} : 
\Gamma_{\frak c}'\right]{\Bbb Z}} 
\ \,\int\limits_{-\infty}^{\infty} \bigl| B_{\frak c}^{\frak d}\left(\omega;it,p\right)\bigr|^2 
h(it,p)\,{\rm d}t\;,}\eqno(6.7.21)$$
for all $h\in{\cal H}_{\star}^{\sigma}(\varrho,\vartheta )$ such that the sums and integrals on 
the right-hand side of this equation are absolutely convergent. 
Since the conditions (i) and (ii) of the lemma are satisfied,  
and since each factor $| C_V^{\frak d}\left(\omega ;\nu_V,p_V\right)|$ 
and $\bigl| B_{\frak c}^{\frak d}\left(\omega;it,p\right)\bigr|$ 
occurring on the right-hand side of Equation~(6.7.21) is real and non-negative, 
it follows that 
(6.7.21) holds for all $h\in\{ f_n : n\in{\Bbb N}\}$, and that,  
given (6.7.13) and what is noted in Remark~1.9.2, one has 
$0\leq Y^{{\frak d},{\frak d}}_{\omega,\omega}(f_1)
\leq Y^{{\frak d},{\frak d}}_{\omega,\omega}(f_2)\leq\ldots\ $.  
Hence, and by (6.7.20), the complex number 
$\chi^{{\frak d},{\frak d}}_{\omega,\omega}(f) 
+X^{{\frak d},{\frak d}}_{\omega,\omega}(f)$ 
must be real and non-negative (certainly one can give a 
more direct proof that it is real), and one must have 
$$\infty >\chi^{{\frak d},{\frak d}}_{\omega,\omega}(f) 
+X^{{\frak d},{\frak d}}_{\omega,\omega}(f) 
\geq Y^{{\frak d},{\frak d}}_{\omega,\omega}\left( f_n\right)\geq 0\qquad\ 
\hbox{for all $\,n\in{\Bbb N}$.}\eqno(6.7.22)$$ 
\par 
In the first sum on the right-hand side of Equation~(6.7.21) one sums 
over only countably many cuspidal irreducible subspaces $V$ 
(i.e. just those occurring in the decomposition (1.7.4) of the space 
${}^0L^{2}(\Gamma\backslash G)$). Let the relevant 
subspaces $V$ be arranged in a sequence, $V(1) , V(2) , \ldots\ $ (say). 
Then, since (6.7.21) holds for all $h\in\{ f_n : n\in{\Bbb N}\}$, 
we have 
$$Y^{{\frak d},{\frak d}}_{\omega,\omega}\left( f_n\right) 
=\sum_{m=1}^{\infty} F_n(m) + \sum_{{\frak c}\in{\frak C}(\Gamma)}
{1\over 4\pi\left[\Gamma_{\frak c} : \Gamma_{\frak c}'\right]} 
\sum_{P\in{\Bbb Z}} 
\ \,\int\limits_{-\infty}^{\infty}\Phi_{{\frak c},n}(t,P)\,{\rm d}t\qquad\qquad 
\hbox{($n\in{\Bbb N}$),}\eqno(6.7.23)$$
where, for $n,m\in{\Bbb N}$,
$$F_n(m)=\,\bigl| c_{V(m)}^{\frak d}\left(\omega ;\nu_{V(m)},p_{V(m)}\right)\bigr|^2 
f_n\!\left(\nu_{V(m)} , p_{V(m)}\right) ,\eqno(6.7.24)$$ 
and, for $n\in{\Bbb N}$, ${\frak c}\in{\frak C}(\Gamma)$ and 
$(t,P)\in{\Bbb R}\times{\Bbb Z}$,
$$\Phi_{{\frak c},n}(t,P)
=\left| B_{\frak c}^{\frak d}\left(\omega;it,
\textstyle{1\over 2}\left[\Gamma_{\frak c} : 
\Gamma_{\frak c}'\right] P\right)\right|^2 
f_n\!\left( it,\textstyle{1\over 2}\left[\Gamma_{\frak c} : 
\Gamma_{\frak c}'\right] P\right)\;.\eqno(6.7.25)$$
\par 
By Remark~1.9.2, each term in the sequence 
$(\nu_{V(1)},p_{V(1)}) , (\nu_{V(1)},p_{V(1)}) , \ldots\ $ 
is contained in the set $E$ defined in (6.7.13), and so, by virtue of 
the hypothesis that condition (ii) of the lemma is satisfied,  
it follows that the sequence $(F_n)_{n\in{\Bbb N}}$ defined in (6.7.24) 
satisfies the condition (6.7.9) of Lemma~6.7.4. Consequently, 
by applying Lemma~6.7.4 and the definition (6.7.24), we find that 
$$\eqalignno{
\lim_{n\rightarrow\infty}\biggl(\,\sum_{m=1}^{\infty} F_n(m)\biggr) 
=\sum_{m=1}^{\infty}\Bigl(\,\lim_{n\rightarrow\infty} F_n(m)\Bigr)
 &=\sum_{m=1}^{\infty}
\,\bigl| c_{V(m)}^{\frak d}\left(\omega ;\nu_{V(m)},p_{V(m)}\right)\bigr|^2 
\Bigl(\,\lim_{n\rightarrow\infty} f_n\!\left(\nu_{V(m)} , p_{V(m)}\right)\Bigr) = \cr 
 &=\sum_{m=1}^{\infty}
\,\bigl| c_{V(m)}^{\frak d}\left(\omega ;\nu_{V(m)},p_{V(m)}\right)\bigr|^2 
f\!\left(\nu_{V(m)} , p_{V(m)}\right) = \cr 
 &=\sum_V 
\,\bigl| C_V^{\frak d}\left(\omega ;\nu_V,p_V\right)\bigr|^2 
f\!\left(\nu_V , p_V\right)
 &(6.7.26)}$$ 
(the penultimate equality following by virtue of the hypothesis that 
the condition (iii) of the lemma is satisfied). Moreover, since 
one may infer from (6.7.21), (6.7.22) and (6.7.25) that, for 
$n\in{\Bbb N}$, ${\frak c}\in{\frak C}(\Gamma)$ and $P\in{\Bbb Z}$, 
the mapping $t\mapsto\Phi_n({\frak c};t,P)$ is measurable with 
respect to the Lebesgue measure on ${\Bbb R}$, one finds (in parallel 
with the preceding) that, for each ${\frak c}\in{\frak C}(\Gamma)$, 
the sequence $(\Phi_n)_{n\in{\Bbb N}}=(\Phi_{{\frak c},n})_{n\in{\Bbb N}}$ 
satisfies the hypotheses of Lemma~6.7.5. Therefore it follows by Lemma~6.7.5 
that, for each cusp ${\frak c}\in{\frak C}(\Gamma)$, one has 
$$\eqalignno{
\lim_{n\rightarrow\infty}\biggl(\,\sum_{P\in{\Bbb Z}}
\ \,\int\limits_{-\infty}^{\infty} 
{\Phi}_{{\frak c},n}(t,P)\,{\rm d}t\biggr)
 &=\sum_{P\in{\Bbb Z}}\ \int\limits_{-\infty}^{\infty}
\Bigl(\,\lim_{n\rightarrow\infty} {\Phi}_{{\frak c},n}(t,P)\Bigr)\,{\rm d}t = \cr 
 &=\sum_{P\in{\Bbb Z}}\ \,\int\limits_{-\infty}^{\infty}
\left| B_{\frak c}^{\frak d}\left(\omega;it,
\textstyle{1\over 2}\left[\Gamma_{\frak c} : 
\Gamma_{\frak c}'\right] P\right)\right|^2 
\Bigl(\,\lim_{n\rightarrow\infty} 
f_n\!\left( it,\textstyle{1\over 2}\left[\Gamma_{\frak c} : 
\Gamma_{\frak c}'\right] P\right)
\Bigr)\,{\rm d}t = \cr 
 &=\sum_{p\in{1\over 2}\left[\Gamma_{\frak c} : 
\Gamma_{\frak c}'\right]{\Bbb Z}}\ \,\int\limits_{-\infty}^{\infty}
\left| B_{\frak c}^{\frak d}\left(\omega;it,p\right)\right|^2 
f\!\left( it,p\right)\,{\rm d}t\;. &(6.7.27)}$$ 
\par 
By Lemma~2.2 the set ${\frak C}(\Gamma)$ is finite. Therefore it 
follows by (6.7.22), (6.7.23), (6.7.26) and (6.7.27) that one has 
$$\eqalign{
\infty &> \lim_{n\rightarrow\infty} 
Y^{{\frak d},{\frak d}}_{\omega,\omega}\left( f_n\right) = \cr 
 &\qquad\qquad =\sum_{V}\,\bigl| C_V^{\frak d}\left(\omega ;\nu_V,p_V\right)\bigr|^2 
f\left(\nu_V , p_V\right) 
+\!\sum_{{\frak c}\in{\frak C}(\Gamma)}
\!{1\over 4\pi\!\left[\Gamma_{\frak c} : \Gamma_{\frak c}'\right]} 
\sum_{p\in{1\over 2}\left[\Gamma_{\frak c} : 
\Gamma_{\frak c}'\right]{\Bbb Z}} 
\ \,\int\limits_{-\infty}^{\infty}\!\bigl| B_{\frak c}^{\frak d}\left(\omega;it,p\right)\bigr|^2 
f(it,p){\rm d}t\;.}$$
This implies that  
all the sums and integrals occurring on the right-hand side of 
Equation~(6.7.21) are absolutely convergent 
when $h=f$, and so  implies also that 
$Y^{{\frak d},{\frak d}}_{\omega,\omega}(f)$  is defined, and that 
the equation (6.7.21) holds 
when $h=f$; by reformulating the above result, more concisely, as the 
statement that one has $\infty > \lim_{n\rightarrow\infty} 
Y^{{\frak d},{\frak d}}_{\omega,\omega}\left( f_n\right)  
=Y^{{\frak d},{\frak d}}_{\omega,\omega}(f)$, one completes the proof of 
(6.7.16), and so too that of the lemma 
\ $\blacksquare$ 

\bigskip   

\proclaim Lemma 6.7.7. Let $(F_n)_{n\in{\Bbb N}}$ be a sequence of 
mappings from ${\Bbb N}$ into ${\Bbb C}$ such that, for each $m\in{\Bbb N}$, 
the sequence $(F_n(m))_{n\in{\Bbb N}}$ is convergent. Suppose moreover that, for some 
mapping $D$ from ${\Bbb N}$ into $[0,\infty)$, one has both 
$$D(m)\geq\left| F_n(m)\right|\qquad\qquad\hbox{($m,n\in{\Bbb N}$).}\eqno(6.7.28)$$ 
and 
$$\sum_{m=1}^{\infty} D(m) < \infty\eqno(6.7.29)$$ 
Then
$$\sum_{m=1}^{\infty}\left|\,\lim_{n\rightarrow\infty} F_n(m)\right| < \infty\;,\eqno(6.7.30)$$ 
and, for some $\lambda\in{\Bbb C}$, both equalities in 
(6.7.10) hold simultaneously. 

\medskip 

\noindent{\bf Proof.}\quad 
See the end of the paragraph above Lemma~6.7.3\quad$\blacksquare$ 

\bigskip 

\proclaim Lemma 6.7.8. Let $({\Phi}_n)_{n\in{\Bbb N}}$ be a sequence of mappings  
from ${\Bbb R}\times{\Bbb Z}$ into ${\Bbb C}$ such that, for each 
$(t,p)\in{\Bbb R}\times{\Bbb Z}$, the sequence $(\Phi_n(t,p))_{n\in{\Bbb N}}$ 
is convergent. 
Suppose moreover that, 
for all $n\in{\Bbb N}$ and all $p\in{\Bbb Z}$, the 
mapping $t\mapsto {\Phi}_n(t,p)$ is a complex valued function on ${\Bbb R}$ that 
is measurable with respect to the Lebesgue measure on~${\Bbb R}$, 
and that, for some mapping $\Delta$ from ${\Bbb R}\times{\Bbb Z}$ into $[0,\infty)$, 
one has both   
$$\Delta(t,p)\geq\left| {\Phi}_n(t,p)\right|\qquad\qquad 
\hbox{($t\in{\Bbb R}$, $p\in{\Bbb Z}$, $n\in{\Bbb N}$)}\eqno(6.7.31)$$ 
and 
$$\sum_{p=-\infty}^{\infty}
\ \int\limits_{-\infty}^{\infty} 
\Delta(t,p)\,{\rm d}t < \infty\;.\eqno(6.7.32)$$
Then 
$$\sum_{p=-\infty}^{\infty}\ \int\limits_{-\infty}^{\infty}
\left|\,\lim_{n\rightarrow\infty} {\Phi}_n(t,p)\right|\,{\rm d}t  < \infty\eqno(6.7.33)$$ 
and, for some $\lambda\in{\Bbb C}$, both equalities in (6.7.12) hold 
simultaneously. 

\medskip 

\noindent{\bf Proof.}\quad  
See the end of the paragraph above Lemma~6.7.3\quad$\blacksquare$ 

\bigskip  

\proclaim Lemma 6.7.9. Let $\sigma\in(1/2,1)\cup(1,2)$ and $\varrho,\vartheta \in(3,\infty)$. 
Let $d\in{\cal H}_{\star}^{\sigma}(\varrho,\vartheta )$ be such that the sum formulae for 
$Y^{{\frak a},{\frak a}}_{\omega_1,\omega_1}(d)$ and 
$Y^{{\frak b},{\frak b}}_{\omega_2,\omega_2}(d)$ are valid. 
Suppose moreover that $E$ is the set defined in (6.7.13), 
that $f\in{\cal H}_{\star}^{\sigma}(\varrho,\vartheta )$, and that $(f_n)_{n\in{\Bbb N}}$ 
is a sequence of elements of ${\cal H}_{\star}^{\sigma}(\varrho,\vartheta )$ which satisfies the 
condition (6.7.6) of Lemma~6.7.3 and the following conditions:\smallskip 
\hskip 15mm (i)\quad for all $n\in{\Bbb N}$, the sum formula for 
$Y^{{\frak a},{\frak b}}_{\omega_1,\omega_2}(f_n)$ is valid;\hfill\smallskip 
\hskip 14.5mm (ii)\quad for all $(\nu,p)\in E$ and all $n\in{\Bbb N}$, one has 
$d(\nu,p)\geq |f_n(\nu,p)|$;\smallskip 
\hskip 14mm (iii)\quad for all $(\nu,p)\in{\Bbb C}\times{\Bbb Z}$ such that 
$\,|{\rm Re}(\nu)|\leq\sigma$, one has 
$\,\lim_{n\rightarrow\infty} f_n(\nu,p)=f(\nu,p)$. \medskip 
\noindent Then  the sum formula for $Y^{{\frak a},{\frak b}}_{\omega_1,\omega_2}(f)$ 
is valid, and one has 
$$\eqalignno{
\chi^{{\frak a},{\frak b}}_{\omega_1,\omega_2}(f) 
 &=\lim_{n\rightarrow\infty}\chi^{{\frak a},{\frak b}}_{\omega_1,\omega_2}(f_n)\;, &(6.7.34) \cr 
X^{{\frak a},{\frak b}}_{\omega_1,\omega_2}(f) 
 &=\lim_{n\rightarrow\infty}X^{{\frak a},{\frak b}}_{\omega_1,\omega_2}(f_n)\;, &(6.7.35) \cr 
Y^{{\frak a},{\frak b}}_{\omega_1,\omega_2}(f) 
 &=\lim_{n\rightarrow\infty}Y^{{\frak a},{\frak b}}_{\omega_1,\omega_2}(f_n)\;. &(6.7.36)}$$ 
 
\medskip 

\noindent{\bf Proof.}\quad 
Since ${\cal H}_{\star}^{\sigma}(\varrho,\vartheta )\supset\{ f_n : n\in{\Bbb N}\}\cup\{ f , d\}$, it 
follows from the definition (6.7.1) that one has 
$${\Bbb C}\ni\chi^{{\frak a},{\frak b}}_{\omega_1,\omega_2}(h) 
={1\over 4\pi^3}\,\delta^{{\frak a},{\frak b}}_{\omega_1,\omega_2} 
\sum_{p\in{\Bbb Z}}\ \int\limits_{-\infty}^{\infty} 
h(it,p)\left( t^2 + p^2\right) {\rm d}t\qquad\  
\hbox{for all $\ h\in\{ f_n : n\in{\Bbb N}\}\cup\{ f , d\}$.}\eqno(6.7.37)$$ 
In particular, since the condition (ii) of the lemma implies that $d(\nu,p)\geq 0$ 
for all $(\nu,p)\in E$, one has  
$$0\leq\sum_{p\in{\Bbb Z}}\ \int\limits_{-\infty}^{\infty} 
d(it,p)\left( t^2 + p^2\right) {\rm d}t < \infty\;.\eqno(6.7.38)$$ 
By virtue of (6.7.13), (6.7.37), (6.7.38) and the conditions (ii) and (iii) 
of the lemma, one may verify that the hypotheses of Lemma~6.7.8 are satisfied 
when, for each $n\in{\Bbb N}$, the mapping $\Phi_n : {\Bbb R}\times{\Bbb Z}\rightarrow{\Bbb C}$ 
is given by $\Phi_n(t,p)=(t^2+p^2) f_n(it,p)\,$  
(note that the conditions (6.7.31) and (6.7.32) are then satisfied if 
one specifies that $\Delta(t,p)=(t^2+p^2) d(it,p)$ for $t\in{\Bbb R}$, $p\in{\Bbb Z}$). 
Therefore, bearing in mind the condition~(iii) of the lemma, it follows by 
Lemma~6.7.8 and (6.7.37) that the equation (6.7.34) holds. 
\par 
We observe next that, since $f$ and the sequence $(f_n)_{n\in{\Bbb N}}$ satisfy 
the relevant hypotheses of Lemma~6.7.3, the result (6.7.35) may therefore be obtained 
similarly to how (6.7.15) was obtained (within the proof of Lemma~6.7.6): 
note, with regard to the relevant application of Lemma~6.5.9, that no distinction 
need be made between those cases where ${\frak a}={\frak b}$ and 
$\omega_1=\omega_2$ and those where either 
${\frak a}\neq{\frak b}$ or else $\omega_1\neq\omega_2$. 
\par 
Given that the proofs of the the results (6.7.34) and (6.7.35) have 
been adequately described, 
and given our hypothesis that the condition~(i) of the lemma is satisfied, 
it will suffice for the completion of the 
proof of the lemma that we show now that the equation (6.7.36) holds   
(our reasoning on this point is similar to that used 
in the paragraph containing (6.7.19), within our proof of Lemma~6.7.6). 
\par 
We begin the proof of (6.7.36) by observing that 
it follows from the condition (i) of the lemma and the definition (6.7.3) that, 
for all $n\in{\Bbb N}$, one has 
$$Y^{{\frak a},{\frak b}}_{\omega_1,\omega_2}\left( f_n\right) 
=\sum_{m=1}^{\infty} F_n(m) + \sum_{{\frak c}\in{\frak C}(\Gamma)}
{1\over 4\pi\left[\Gamma_{\frak c} : \Gamma_{\frak c}'\right]} 
\sum_{P\in{\Bbb Z}} 
\ \,\int\limits_{-\infty}^{\infty}\Phi_{{\frak c},n}(t,P)\,{\rm d}t\;,\eqno(6.7.39)$$
where 
$$\Phi_{{\frak c},n}(t,P)
=\,\overline{B_{\frak c}^{\frak a}\left(\omega_1;it,
\textstyle{1\over 2}\left[\Gamma_{\frak c} : 
\Gamma_{\frak c}'\right] P\right)}  
\,B_{\frak c}^{\frak b}\left(\omega_2;it,
\textstyle{1\over 2}\left[\Gamma_{\frak c} : 
\Gamma_{\frak c}'\right] P\right) 
f_n\!\left( it,\textstyle{1\over 2}\left[\Gamma_{\frak c} : 
\Gamma_{\frak c}'\right] P\right)\eqno(6.7.40)$$
and 
$$F_n(m)
=\,\overline{ c_{V(m)}^{\frak a}\left(\omega_1 ;\nu_{V(m)},p_{V(m)}\right)} 
\,c_{V(m)}^{\frak b}\left(\omega_2 ;\nu_{V(m)},p_{V(m)}\right)  
f_n\!\left(\nu_{V(m)} , p_{V(m)}\right) ,\eqno(6.7.41)$$ 
with the sequence of spaces $V(1) , V(2) , \ldots\ $ being as indicated 
between (6.7.22) and (6.7.23), within the proof of Lemma~6.7.6.
Note that, by (6.7.40), (6.7.41), 
the condition (iii) of the lemma and the hypothesis that 
$f\in{\cal H}_{\star}^{\sigma}(\varrho,\vartheta )$,   
one has, for $m\in{\Bbb N}$,  
$$\lim_{n\rightarrow\infty} F_n(m) 
=\,\overline{ c_{V(m)}^{\frak a}\!\left(\omega_1 ;\nu_{V(m)},p_{V(m)}\right)} 
\,c_{V(m)}^{\frak b}\!\left(\omega_2 ;\nu_{V(m)},p_{V(m)}\right)  
f\!\left(\nu_{V(m)} , p_{V(m)}\right)\in{\Bbb C}\eqno(6.7.42)$$ 
and, for  $(t,P)\in{\Bbb R}\times{\Bbb Z}$ and ${\frak c}\in{\frak C}(\Gamma)$, 
$$\lim_{n\rightarrow\infty}\Phi_{{\frak c},n}(t,P) 
=\,\overline{B_{\frak c}^{\frak a}\left(\omega_1;it,
\textstyle{1\over 2}\left[\Gamma_{\frak c} : 
\Gamma_{\frak c}'\right] P\right)}  
\,B_{\frak c}^{\frak b}\left(\omega_2;it,
\textstyle{1\over 2}\left[\Gamma_{\frak c} : 
\Gamma_{\frak c}'\right] P\right) 
f\!\left( it,\textstyle{1\over 2}\left[\Gamma_{\frak c} : 
\Gamma_{\frak c}'\right] P\right)\in{\Bbb C}\;.
\eqno(6.7.43)$$
\par 
Since the sum formulae for $Y^{{\frak a},{\frak a}}_{\omega_1,\omega_1}(d)$ and 
$Y^{{\frak b},{\frak b}}_{\omega_2,\omega_2}(d)$ are valid, one has also 
$${\Bbb C}\ni 
{Y^{{\frak a},{\frak a}}_{\omega_1,\omega_1}(d) 
+Y^{{\frak b},{\frak b}}_{\omega_2,\omega_2}(d)\over 2} 
=\sum_{m=1}^{\infty} D(m) + \sum_{{\frak c}\in{\frak C}(\Gamma)}
{1\over 4\pi\left[\Gamma_{\frak c} : \Gamma_{\frak c}'\right]} 
\sum_{P\in{\Bbb Z}} 
\ \,\int\limits_{-\infty}^{\infty}\Delta_{\frak c}(t,P)\,{\rm d}t\;,\eqno(6.7.44)$$
where 
$$D(m) 
={1\over 2}\left( 
\bigl| c_{V(m)}^{\frak a}\!\left(\omega_1 ;\nu_{V(m)},p_{V(m)}\right)\bigr|^2 
+\bigl| c_{V(m)}^{\frak b}\!\left(\omega_2 ;\nu_{V(m)},p_{V(m)}\right)\bigr|^2 
\right) d\left(\nu_{V(m)},p_{V(m)}\right)\eqno(6.7.45)$$ 
and 
$$\Delta_{\frak c}(t,P) 
={1\over 2}\left( 
\bigl| B_{\frak c}^{\frak a}\left(\omega_1;it,
\textstyle{1\over 2}\left[\Gamma_{\frak c} : 
\Gamma_{\frak c}'\right] P\right)\bigr|^2 
+\bigl| B_{\frak c}^{\frak b}\left(\omega_2;it,
\textstyle{1\over 2}\left[\Gamma_{\frak c} : 
\Gamma_{\frak c}'\right] P\right)\bigr|^2  
\right) d\!\left( it,\textstyle{1\over 2}\left[\Gamma_{\frak c} : 
\Gamma_{\frak c}'\right] P\right)\;.\eqno(6.7.46)$$ 
Note that, by the condition (ii) of the lemma, and the definition (6.7.13), 
all terms of the sums occurring on the right-hand side of the equality sign 
in (6.7.44) are real and non-negative; it therefore follows from (6.7.44) 
that the function $D : {\Bbb N}\rightarrow [0,\infty)$ satisfies the 
condition (6.7.29) of Lemma~6.7.7, and that one has 
$$\sum_{P=-\infty}^{\infty} 
\ \int\limits_{-\infty}^{\infty}\Delta_{\frak c}(t,P)\,{\rm d}t < \infty\qquad\  
\hbox{for ${\frak c}\in{\frak C}(\Gamma)$.}$$
Moreover, by (6.7.40), (6.7.41), (6.7.45), (6.7.46), the arithmetic-geometric mean 
inequality, Remark~1.9.2 and the condition~(ii) of the  lemma, one finds that 
$D$ and the sequence $(F_n)_{n\in{\Bbb N}}$ satisfy the condition (6.7.28) of 
Lemma~6.7.7, and that 
$$\Delta_{\frak c}(t,p)\geq\left|\Phi_{{\frak c},n}(t,P)\right|\qquad\qquad 
\hbox{(${\frak c}\in{\frak C}(\Gamma)$, $(t,P)\in{\Bbb R}\times{\Bbb Z}$, 
$n\in{\Bbb N}$).}$$ 
\par 
Given (6.7.42), (6.7.43) and the points just noted in the preceding paragraph, 
and given that (6.7.39) holds for all $n\in{\Bbb N}$, it follows by 
Lemma~6.7.7 and Lemma~6.7.8 that one has 
$$\lim_{n\rightarrow\infty}\sum_{m=1}^{\infty} F_n(m) 
=\sum_{m=1}^{\infty} 
\,\overline{ c_{V(m)}^{\frak a}\left(\omega_1 ;\nu_{V(m)},p_{V(m)}\right)} 
\,c_{V(m)}^{\frak b}\left(\omega_2 ;\nu_{V(m)},p_{V(m)}\right)  
f\!\left(\nu_{V(m)} , p_{V(m)}\right)\eqno(6.7.47)$$ 
and, for each ${\frak c}\in{\frak C}(\Gamma)$, 
$$\lim_{n\rightarrow\infty} 
\sum_{P\in{\Bbb Z}} 
\ \,\int\limits_{-\infty}^{\infty}\Phi_{{\frak c},n}(t,P)\,{\rm d}t 
={1\over i}\sum_{p\in{1\over 2}\left[\Gamma_{\frak c} : \Gamma_{\frak c}'\right]{\Bbb Z}} 
\ \,\int\limits_{(0)} 
\,\overline{B_{\frak c}^{\frak a}\left(\omega_1;\nu,p\right)}  
\,B_{\frak c}^{\frak b}\left(\omega_2;\nu,p\right) 
f(\nu,p)\,{\rm d}\nu\;.\eqno(6.7.48)$$ 
It is worth clarifying here that, by virtue of the results (6.7.30) and (6.7.33) 
of the lemmas just applied, all of the integrals and sums occurring on the 
right-hand sides of the equations~(6.7.47) and (6.7.48) are absolutely convergent. 
Therefore, since the set ${\frak C}(\Gamma)$ is finite, and since 
$\{ V(m) : m\in{\Bbb N}\}$ is the set of cuspidal irreducible subspaces 
occurring in the decomposition (1.7.4) of the space ${}^0L^2(\Gamma\backslash G)$, 
it follows that all of the integrals and sums occurring on the right-hand side of 
Equation~(6.7.3) are absolutely convergent when 
${\frak d}={\frak a}$, ${\frak d}'={\frak b}$, $\omega=\omega_1$, $\omega'=\omega_2$ 
and $h=f$. The case 
$({\frak d},{\frak d}',\omega,\omega',h)=({\frak a},{\frak b},\omega_1,\omega_2,f)$ 
of~(6.7.3) therefore defines $Y^{{\frak a},{\frak b}}_{\omega_1,\omega_2}(f)$, and 
defines it in such a way that one has 
$Y^{{\frak a},{\frak b}}_{\omega_1,\omega_2}(f)\in{\Bbb C}$. 
We observe moreover that the limits occurring on the left-hand sides of   
the equations~(6.7.47) and~(6.7.48) exist 
(i.e. they are limits of convergent sequences), and so, 
by (6.7.39) (for $n\in{\Bbb N}$), 
the case $({\frak d},{\frak d}',\omega,\omega',h)=({\frak a},{\frak b},\omega_1,\omega_2,f)$ 
of~(6.7.3) and the 
finiteness of the set ${\frak C}(\Gamma)$, it follows from 
(6.7.47) and (6.7.48) that the equation (6.7.36) holds\quad$\blacksquare$  

\bigskip 

\proclaim Lemma 6.7.10. Let $\sigma\in(1,2)$, let $\varrho,\vartheta \in(3,\infty)$,  
and let $f\in{\cal H}_{\star}^{\sigma}(\varrho,\vartheta )$. Put 
$$C_f(\varrho,\vartheta ) 
=\sup\left\{ 
\left( 1+\left|{\rm Im}(\nu)\right|\right)^{\varrho} 
\left( 1+|p|\right)^{\vartheta} |f(\nu,p)| \,:\, p\in{\Bbb Z} ,\ \nu\in{\Bbb C}\ {\rm and}\ 
\left| {\rm Re}(\nu)\right|\leq\sigma\right\}\;.\eqno(6.7.49)$$ 
Suppose moreover that $\ell\in{\Bbb N}$, and that 
$f_{\ell}$ is the mapping 
from $\{\nu\in{\Bbb C} : |{\rm Re}(\nu)|\leq\sigma\}\times{\Bbb Z}$ into 
${\Bbb C}$ given by 
$$f_{\ell}(\nu,p)=\cases{f(\nu,p)\exp\left( -\nu^4 /\ell\right) 
 &if $|p|\leq\ell$; \cr 0 &otherwise.}\eqno(6.7.50)$$ 
Then $f_{\ell}\in{\cal H}_{\star}^{\sigma}(\varrho,\vartheta )$, one has 
$$\sup\left\{ 
\left( 1+\left|{\rm Im}(\nu)\right|\right)^{\varrho} 
\left( 1+|p|\right)^{\vartheta}\left| f_{\ell}(\nu,p)\right| \,:\, p\in{\Bbb Z} , 
\ \nu\in{\Bbb C}\ {\rm and}\ 
|{\rm Re}(\nu)|\leq\sigma\right\} 
\leq\exp\!\left( 8\sigma^{4}/\ell\right) C_f(\varrho,\vartheta )<\infty\;,\eqno(6.7.51)$$ 
and, for ${\frak d},{\frak d}'\in\{ {\frak a} , {\frak b}\}$ and 
$\omega,\omega'\in\{\omega_1,\omega_2\}$, the sum formula for 
$Y^{{\frak d},{\frak d}'}_{\omega,\omega'}(f_{\ell})$ is valid. 

\medskip 

\noindent{\bf Proof.}\quad 
Since $f\in{\cal H}_{\star}^{\sigma}(\varrho,\vartheta )$, since the function $p\mapsto |p|$ is even, 
and since the function $\nu\mapsto\exp(-\nu^4 /\ell)$ is both 
even and holomorphic on ${\Bbb C}$, 
one finds,  when checking that $f_{\ell}\in{\cal H}_{\star}^{\sigma}(\varrho,\vartheta )$, 
that it is only the condition~(iii) of Theorem~B 
that requires more than cursory consideration; 
moreover, 
given the definition (6.7.50), and given that  $f\in{\cal H}_{\star}^{\sigma}(\varrho,\vartheta )$ and 
$\ell\in{\Bbb N}$, 
one need do little more than observe that   
$$\left|\exp\left( -\nu^4 /\ell\right)\right| 
\leq\exp\left( 8\sigma^4 /\ell\right)\qquad\qquad 
\hbox{($|{\rm Re}(\nu)|\leq\sigma$)}\eqno(6.7.52)$$ 
in order to verify that the condition~(iii) of Theorem~B is 
satisfied when $h=f_{\ell}$.
Since further discussion of the proof that $f_{\ell}$ lies in 
the space ${\cal H}_{\star}^{\sigma}(\varrho,\vartheta )$ is probably unnecessary, we skip it. 
\par 
We can also be brief in discussing the proof of the result (6.7.51): 
this result follows, by virtue of (6.7.52), from the definitions (6.7.49) and (6.7.50), 
and the hypothesis that $f\in{\cal H}_{\star}^{\sigma}(\varrho,\vartheta )$.  
\par 
What remains to be demonstrated is the validity of the sum formula 
for $Y^{{\frak d},{\frak d}'}_{\omega,\omega'}$, 
when ${\frak d},{\frak d}'\in\{ {\frak a} , {\frak b}\}$ and 
$\omega,\omega'\in\{\omega_1,\omega_2\}$. 
Let ${\frak d},{\frak d}'\in\{ {\frak a} , {\frak b}\}$ and 
$\omega,\omega'\in\{\omega_1,\omega_2\}$. 
Suppose also that $\lambda^{*}_{\ell}$ is the mapping from 
${\Bbb C}\times\{ -\ell,1-\ell,\ldots ,\ell\}$ into ${\Bbb C}$ defined by 
the equations~(6.6.4) of Proposition~6.6.1. 
As a first step towards verifying  
the sum formula for $Y^{{\frak d},{\frak d}'}_{\omega,\omega'}$, we observe that, when $p\in{\Bbb Z}$, 
$\nu\in{\Bbb C}$ and $|{\rm Re}(\nu)|\leq\sigma$, it follows from (6.7.50) 
that 
$$f_{\ell}(\nu,p)=\cases{ 
\lambda^{*}_{\ell}(\nu,p)\,\overline{\theta_{\ell}\!\left( -\overline{\nu}\,,p\right)} 
\,\eta_{\ell}(\nu,p) &if $\,|p|\leq\ell$, \cr 0 &otherwise,}\eqno(6.7.53)$$ 
where $\eta_{\ell}$ and $\theta_{\ell}$ are the mappings from 
$\{\nu\in{\Bbb C} : |{\rm Re}(\nu)|\leq\sigma\}\times\{ -\ell,1-\ell,\ldots ,\ell\}$ 
into ${\Bbb C}$ given by 
$$\eta_{\ell}(\nu,p)={f(\nu,p)\over\lambda^{*}_{\ell}(\nu,p)}\;,\eqno(6.7.54)$$ 
$$\theta_{\ell}(\nu,p)=\exp\left( -\nu^4 /\ell\right)\;.\eqno(6.7.55)$$ 
In cases where $\lambda^{*}_{\ell}(\nu,p)=0$, we take (6.7.54) to 
mean that $\eta_{\ell}(\nu,p)=\lim_{\delta\rightarrow 0+} 
f(\nu +i\delta,p)/\lambda^{*}_{\ell}(\nu+i\delta,p)$.
 
\par 
If the functions $\eta_{\ell}$ and $\theta_{\ell}$  
lie in the space ${\cal T}^{\ell}_{\sigma}$ defined in the paragraph containing~(6.4.3), 
then by Proposition~6.6.1 (applied with $\eta_{\ell}$, $\theta_{\ell}$, ${\frak d}$, 
${\frak d}'$, $\omega$ and $\omega'$ substituted for 
$\eta$, $\theta$, ${\frak a}$, ${\frak b}$, $\omega_1$ and $\omega_2$, respectively)  
one obtains, given (6.7.53), the case $h=f_{\ell}$ of Equation~(6.7.4) 
(the terms occurring in that equation being defined by (6.7.1)-(6.7.3)). 
Moreover, by the same application of Proposition~6.6.1, one finds that if 
$\eta_{\ell},\theta_{\ell}\in{\cal T}^{\ell}_{\sigma}$ then the 
sums and integrals involved in the definition (via (6.7.3)) of 
$Y^{{\frak d},{\frak d}'}_{\omega,\omega'}(f_{\ell})$ are all 
absolutely convergent. One may therefore conclude that 
$$\hbox{the sum formula for}\ 
Y^{{\frak d},{\frak d}'}_{\omega,\omega'}\!\left(f_{\ell}\right)\ 
\hbox{is valid if $\,\eta_{\ell},\theta_{\ell}\in{\cal T}^{\ell}_{\sigma}$.}
\eqno(6.7.56)$$ 
\par 
We show next that 
if $\eta\in\{\eta_{\ell} , \theta_{\ell}\}$ then 
the conditions (T1)-(T3) below (6.4.3) are satisfied; 
we thereby establish that the functions 
$\eta_{\ell}$ and $\theta_{\ell}$ do lie in the space 
${\cal T}^{\ell}_{\sigma}$. 
\par 
Starting with the easier case, $\eta =\theta_{\ell}$, we note firstly 
that the mapping 
$(\nu,p)\mapsto(-\nu,-p)$ is a permutation of the domain of $\theta_{\ell}$, and that 
by (6.7.55) one has 
$\theta_{\ell}(-\nu,-p)=\exp(-(-\nu)^4/ \ell)=\exp(-\nu^4 /\ell)=\theta_{\ell}(\nu,p)$, 
for all $(\nu,p)$ in that domain. Therefore the condition (T1) below (6.4.3) is 
satisfied when $\eta =\theta_{\ell}$. 
Since the complex function $\nu\mapsto\exp(-\nu^4 /\ell)$ is entire, we find, secondly, 
that the condition (T2) is satisfied when $\eta =\theta_{\ell}$. 
Thirdly, we note that, for $A>0$, $\alpha\in[-\sigma , -\sigma]$, $t\geq 0$ and 
$\nu =\alpha \pm it$, one has 
$$\eqalign{ 
{\left|\exp\left( -\nu^4 /\ell\right)\right|\over 
\left( 1+\left|{\rm Im}(\nu)\right|\right)^{-A} e^{-(\pi /2)\left|{\rm Im}(\nu)\right|}} 
 &\leq\exp\!\left(\left( A+{\pi\over 2}\right) t 
-{\left( t^4+\alpha^4 -6\alpha^2 t^2\right)\over\ell}\right) \leq \cr 
 &\leq\exp\left( 
\max_{x\in{\Bbb R}} \left( ( A+{\pi\over 2})\,\ell^{1/4} x 
+6\sigma^2\ell^{-1/2} x^2-x^4\right)\right) <\infty\;,}$$
so that the condition (T3) is satisfied when $\eta =\theta_{\ell}$. 
\par 
We consider next the conditions (T1)-(T3) in the case where $\eta$ is 
equal to the function $\eta_{\ell}$ defined in Equation~(6.7.54). 
\par 
Given the definition of $\lambda^{*}_{\ell}(\nu,p)$ in (6.6.4) 
(noting, in particular, the second equality there), and given that $\ell$ is 
a positive integer, it follows that for each $p\in\{-\ell,1-\ell,\ldots ,\ell\}$ 
the mapping 
$\nu\mapsto (\nu^2 -1)^2 /\lambda^{*}_{\ell}(\nu,p)$ is 
is a holomorphic function on 
the open strip $\{\nu\in{\Bbb C} : |{\rm Re}(\nu)|<2\}$. 
Therefore, since the hypotheses that $\sigma\in(1,2)$ and 
$f\in{\cal H}_{\star}^{\sigma}(\varrho,\vartheta )$ imply that for each $p\in{\Bbb Z}$ the mapping 
$\nu\mapsto f(\nu,p)/(\nu^2 -1)^2$ can be holomorphically continued into  
a neighbourhood of the closed strip 
$\{\nu\in{\Bbb C} : |{\rm Re}(\nu)|\leq\sigma\}$,  it follows
from (6.7.54) that the condition (T2) below (6.4.3) is satisfied when 
$\eta =\eta_{\ell}$. 
\par 
By the definition~(6.6.4) and the hypothesis that $f\in{\cal H}_{\star}^{\sigma}(\varrho,\vartheta )$, one has 
both $1/\lambda^{*}_{\ell}(-\nu,-p)=1/\lambda^{*}_{\ell}(\nu,p)$ and 
$f(-\nu,-p)=f(\nu,p)$ for all 
$(\nu,p)\in({\Bbb C}-{\Bbb Z})\times\{-\ell,1-\ell,\ldots,\ell\}$ 
such that $|{\rm Re}(\nu)\leq\sigma$. 
Therefore, given the definition~(6.7.54), and given that the condition 
(T2) below (6.4.3) is satisfied when $\eta =\eta_{\ell}$, one may deduce that 
$\eta_{\ell}(-\nu,-p)=\eta_{\ell}(\nu,p)$ for all 
$(\nu,p)\in{\Bbb C}\times{\Bbb Z}$ such that 
$|{\rm Re}(\nu)|\leq\sigma$ and $|p|\leq\ell$. It follows that 
the condition (T1) below~(6.4.3) is satisfied when $\eta =\eta_{\ell}$. 
\par 
We now have only to check that the condition (T3) below (6.4.3) is 
satisfied when $\eta =\eta_{\ell}$; in order that this 
be verified, it will suffice to show that, for all $A>0$, one has  
$$\sup\left\{ \left( 1+\left|{\rm Im}(\nu)\right|\right)^A 
e^{(\pi /2)|{\rm Im}(\nu)|} \left| \eta_{\ell}(\nu,p)\right| 
\,:\ \nu\in{\Bbb C},\ p\in{\Bbb Z},\ \left| {\rm Re}(\nu)\right|\leq\sigma   
\ {\rm and}\ |p|\leq\ell\right\} <\infty\;.\eqno(6.7.57)$$ 
\par 
Since the condition (T2) below (6.4.3) is satisfied when $\eta =\eta_{\ell}$, 
the function $(\nu,p)\mapsto|\eta_{\ell}(\nu,p)|$ is bounded 
on the compact set 
$\{\nu\in{\Bbb C} : |{\rm Re}(\nu)|, |{\rm Im}(\nu)|\leq\sigma\} 
\times\{ p\in{\Bbb Z} : |p|\leq\ell\}=T(\sigma,\ell)\,$ (say). 
Moreover, given any $A>0$, the function 
$(\nu,p)\mapsto\left( 1+\left|{\rm Im}(\nu)\right|\right)^A 
e^{(\pi /2)|{\rm Im}(\nu)|}$ is bounded on the same set, 
$T(\sigma,\ell)$. 
Therefore the condition (6.7.57) is satisfied if and only if 
it is the case that, for each $A>0$, the function 
$(\nu,p)\mapsto\left( 1+\left|{\rm Im}(\nu)\right|\right)^A 
e^{(\pi /2)|{\rm Im}(\nu)|}|\eta_{\ell}(\nu,p)|$ is bounded on the set 
$U(\sigma,\ell)$ given by 
$$U(\sigma,\ell) 
=\left\{ (\nu,p)\in{\Bbb C}\times{\Bbb Z} 
\,:\ |{\rm Re}(\nu)|\leq\sigma\leq|{\rm Im}(\nu)|\ \,{\rm and}   
\ \,|p|\leq\ell\right\}\;.\eqno(6.7.58)$$   
\par 
Since $\sigma\in(1,2)$, and since $\ell\in{\Bbb N}$, it follows by (6.6.4), (6.4.5) 
and (6.7.58) that  
$$\eqalignno{ 
{1\over\left|\lambda^{*}_{\ell}(\nu,p)\right|} 
 &={\pi |\nu|\over\left|\sin(\pi\nu)\right|} 
\,{\left|(\nu +p)(\nu -p)\right|\over |\nu|^{1+\epsilon(p)}} 
\prod_{\scriptstyle 0<m\leq\ell\atop\scriptstyle m\neq |p|} 
{1\over\left| (\nu +m)(\nu -m)\right|} \leq 
\qquad\qquad\qquad\qquad\qquad\qquad\qquad\qquad\quad\cr 
 &\leq {2\pi\left( 1+\left|{\rm Im}(\nu)\right|\right)\over\sinh\left( 
\pi\left|{\rm Im}(\nu)\right|\right)} 
\,{4\ell^2 \left( 1+\left|{\rm Im}(\nu)\right|\right)^2\over 
\left|{\rm Im}(\nu)\right|^{2\ell +2\epsilon(p)}} 
\leq {8\pi\ell^2 \left( 1+\left|{\rm Im}(\nu)\right|\right)^3\over 
(1/3)\exp\left(\pi \left|{\rm Im}(\nu)\right|\right)}\qquad\ 
\hbox{for $\,(\nu,p)\in U(\sigma,\ell)$.} &(6.7.59)}$$ 
We are also given that $f\in{\cal H}_{\star}^{\sigma}(\varrho,\vartheta )$, and that 
$(\varrho,\vartheta )\in(3,\infty)$, and so may deduce from (6.7.54), (6.7.49),  
(6.7.58), (6.7.59) and the final inequality in (6.7.51) 
that, for all real $A>0$, one has 
$$\eqalign{
 &\sup_{(\nu,p)\in U(\sigma,\ell)} 
\left( 1+\left|{\rm Im}(\nu)\right|\right)^A 
e^{(\pi /2)|{\rm Im}(\nu)|} \left| \eta_{\ell}(\nu,p)\right| \leq \cr 
 &\qquad\qquad\qquad\qquad\qquad\qquad\leq 24\pi\ell^2 C_f(\varrho,\vartheta ) 
\sup_{(\nu,p)\in U(\sigma,\ell)} 
\left( 1+\left|{\rm Im}(\nu)\right|\right)^{A+3-\varrho} 
\left( 1+|p|\right)^{-\vartheta}  
e^{-(\pi /2)|{\rm Im}(\nu)|}  \leq \cr 
 &\qquad\qquad\qquad\qquad\qquad\qquad\leq 24\pi\ell^2 C_f(\varrho,\vartheta ) 
\max_{t\geq 0}\ (1+t)^A e^{-(\pi /2)t} < \infty\;.}$$
By this result, in combination with what has been discussed in the paragraph  
containing (6.7.58), it follows that the condition (6.7.57) is satisfied. 
It has therefore been verified that the condition (T3) below (6.4.3) 
is satisfied when $\eta =\eta_{\ell}$. 
\par 
We have now found that each of the conditions (T1)-(T3) below (6.4.3) is satisfied 
if one has either $\eta =\theta_{\ell}$ or $\eta =\eta_{\ell}$. 
Therefore the functions $\theta_{\ell}$ and $\eta_{\ell}$ lie in the space 
${\cal T}^{\ell}_{\sigma}$, and so, by (6.7.56), the sum formula for 
$Y^{{\frak d},{\frak d}'}_{\omega,\omega'}(f_{\ell})$ is valid\quad$\blacksquare$ 

\bigskip 

\proclaim Lemma 6.7.11. Let $\sigma\in(1,2)$, let $\varrho,\vartheta \in(3,\infty)$, 
and let $E$ be the set defined in (6.7.13). Suppose moreover that 
$d^{\flat}$ is the mapping from $\{\nu\in{\Bbb C} : |{\rm Re}(\nu)|\leq\sigma\} 
\times{\Bbb Z}$ into ${\Bbb C}$ given by 
$$d^{\flat}(\nu,p)=\left( 1-\nu^2\right)^2 
\left( 4-\nu^2\right)^{-(\varrho +4)/2}\left( 1+|p|\right)^{-\vartheta}\;.\eqno(6.7.60)$$ 
Then $d^{\flat}\in{\cal H}_{\star}^{\sigma}(\varrho,\vartheta )$, one has 
$$d^{\flat}(\nu,p)\geq 2^{-(\varrho +5)}\left( 1+\left|{\rm Im}(\nu)\right|\right)^{-\varrho} 
\left( 1+|p|\right)^{-\vartheta}\qquad\qquad\hbox{($\,(\nu,p)\in E$),}\eqno(6.7.61)$$ 
and, for ${\frak d}\in\{ {\frak a} , {\frak b}\}$ and 
$\omega\in\{ \omega_1 , \omega_2\}$, the sum formula for 
$Y^{{\frak d},{\frak d}}_{\omega,\omega}(d^{\flat})$ is valid. 

\medskip 

\noindent{\bf Proof.}\quad 
The function $d^{\flat}$ defined in (6.7.60) is a very minor  
modification of the function `$f_{a,b}$' which is defined and used 
in the proof of Proposition~11.3.2 of [32]. 
Since one has ${\rm Re}(4-\nu^2)\geq 4-\sigma^2>0$ when $|{\rm Re}(\nu)|\leq\sigma$, 
it follows that the mapping 
$\nu\mapsto(1-\nu^2)^2 (4-\nu^2)^{-(\varrho +4)/2}$ is holomorphic on a neighbourhood of 
the strip $\{\nu\in{\Bbb C} : |{\rm Re}(\nu)|\leq\sigma\}\,$ (i.e. one has 
$(4-\nu^2)^{-(\varrho +4)/2}=\exp(-(1/2)(\varrho +4)\log(4-\nu^2))$, where 
$\log(z)$ denotes the principal branch of the logarithm function). 
Consequently the condition (ii) of Theorem~B is satisfied when 
$h=d^{\flat}$. Given (6.7.60), and given that $1<\sigma<2$, one has also 
$d^{\flat}(-\nu,-p)=d^{\flat}(\nu,p)$ and, via a short calculation, 
$| d^{\flat}(\nu,p)|\leq 2^{(\varrho +4)/2} (\sigma +1)^4 (2-\sigma)^{-(\varrho +4)} 
(1+|{\rm Im}(\nu)|)^{-\varrho} (1+|p|)^{-\vartheta}$ for all $p\in{\Bbb Z}$ and all 
$\nu\in{\Bbb C}$ such that $|{\rm Re}(\nu)|\leq\sigma$.  
Therefore the conditions (i) and (iii) of Theorem~B are satisfied when 
$h=d^{\flat}$. Since it is moreover the case that the function 
$\nu\mapsto (1-\nu^2)^2 (4-\nu^2)^{-(\varrho +4)/2}$ has a zero of order $2$ 
at $\nu =1$, we may conclude (given what has already been noted above) 
that $d^{\flat}$ lies in the space ${\cal H}_{\star}^{\sigma}(\varrho,\vartheta )$. 
\par 
When $(\nu,p)\in (i{\Bbb R})\times{\Bbb Z}$, the inequality stated in (6.7.61) follows 
by virtue of the definition (6.7.60), the hypothesis that one has $\varrho,\vartheta \in(3,\infty)$, 
and the fact that $0<4+t^2\leq 4(1+t^2)\leq 4(1+|t|)^2$ for $t\in{\Bbb R}$. 
Consequently, given the definition (6.7.13) of the set $E$, 
the proof of (6.7.61) may be completed by making use of 
the observation that
for $\nu\in[-2/9,2/9]$ and $\lambda =1-\nu^2$ one has $1\geq\lambda >3/4$, and so  
$\lambda^2 (\lambda +3)^{-(\varrho +4)/2}\geq 4^{-\varrho /2} (1+3/\lambda)^{-2}\geq 2^{-\varrho}/25 
>2^{-(\varrho +5)}$. 
\par 
Suppose now that ${\frak d}\in\{ {\frak a} , {\frak b}\}$ and 
$\omega\in\{\omega_1 , \omega_2\}$. The proof of the lemma will be 
complete if we are able to show that the sum formula for 
$Y^{{\frak d},{\frak d}}_{\omega,\omega}(d^{\flat})$ is valid. 
As a first step towards this, we let $d^{\flat}_{\ell}$ denote 
(when $\ell\in{\Bbb N}$) the mapping from 
$\{\nu\in{\Bbb C} : |{\rm Re}(\nu)|\leq\sigma\}\times{\Bbb Z}$ into 
${\Bbb C}$ given by 
$$d^{\flat}_{\ell}(\nu,p) 
=\cases{d^{\flat}(\nu,p)\exp\left( -\nu^4 /\ell\right) 
 &if $|p|\leq\ell$; \cr 0 &otherwise.}\eqno(6.7.62)$$ 
\par 
By the case $f=d^{\flat}$ of Lemma~6.7.10, it follows that 
the condition (6.7.6) of Lemma~6.7.3 is satisfied when 
$(f_n)_{n\in{\Bbb N}}$ is the sequence $(d^{\flat}_{\ell})_{\ell\in{\Bbb N}}$, 
and that, for each $\ell\in{\Bbb N}$, the function  
$d^{\flat}_{\ell}$ lies in the space ${\cal H}_{\star}^{\sigma}(\varrho,\vartheta )$
and the sum formula for $Y^{{\frak d},{\frak d}}_{\omega,\omega}(d^{\flat}_{\ell})$ 
is valid. By (6.7.61), (6.7.62) and (6.7.13), we have moreover 
$$0\leq d^{\flat}_1(\nu,p)\leq d^{\flat}_2(\nu,p)\leq\ldots \qquad\qquad 
\hbox{($\,(\nu,p)\in E$)}$$ 
and 
$$\lim_{\ell\rightarrow\infty} d^{\flat}_{\ell}(\nu,p) 
=d^{\flat}(\nu,p)\qquad\qquad 
\hbox{($p\in{\Bbb Z}$, $\nu\in{\Bbb C}$ and $|{\rm Re}(\nu)|\leq\sigma$),}$$ 
and so may conclude that, when 
$\sigma$, $\varrho$, $\vartheta$, ${\frak d}$, $\omega$ and $E$ are as we currently 
suppose, when $f=d^{\flat}$, and when 
$(f_n)_{n\in{\Bbb N}}$ is the sequence $(d^{\flat}_{\ell})_{\ell\in{\Bbb N}}$, 
it is then the case that all of the hypotheses of Lemma~6.7.6 
are satisfied. Therefore it follows by Lemma~6.7.6 
that the sum formula for $Y^{{\frak d},{\frak d}}_{\omega,\omega}(d^{\flat})$ 
is valid\quad$\blacksquare$ 

\bigskip 

\noindent{\bf The proof of Proposition~6.7.2.}\quad  
Let $\sigma$, $\varrho$, $\vartheta$ and $f$ satisfy the hypotheses of Proposition~6.7.2. 
For each $\ell\in{\Bbb N}$, define $f_{\ell}$ to be the mapping from 
$\{\nu\in{\Bbb C} : |{\rm Re}(\nu)|\leq\sigma\}\times{\Bbb Z}$ into ${\Bbb C}$ 
given by the equation (6.7.50) of Lemma~6.7.10. 
Define the set $E$ as in (6.7.13), and put 
$$d=D_f (\varrho,\vartheta ) d^{\flat}\;,\eqno(6.7.63)$$ 
where $d^{\flat}$ is the mapping from 
$\{\nu\in{\Bbb C} : |{\rm Re}(\nu)|\leq\sigma\}\times{\Bbb Z}$ into ${\Bbb C}$ 
defined in Lemma~6.7.11, while 
$$D_f (\varrho,\vartheta ) = 2^{\varrho +5} \exp\left( 8\sigma^4\right) C_f (\varrho,\vartheta )\eqno(6.7.64)$$ 
with $C_f (\varrho,\vartheta )$ being the constant defined in Lemma~6.7.10 
(so that, as follows by the result (6.7.51) of that lemma, 
one has $0\leq C_f (\varrho,\vartheta ) <\infty$). 
\par 
We show next that 
when $\sigma$, $\varrho$, $\vartheta$ and $f$ are as we currently suppose,   
and when the set $E$, the function $d$ and the sequence 
$(f_{\ell})_{\ell\in{\Bbb N}}$ are as just defined above, 
all of the hypotheses of Lemma~6.7.9 are satisfied. 
Once this is achieved the proof of Proposition~6.7.2 will be essentially complete: 
for it will then follow, by Lemma~6.7.9,  that the 
sum formula for $Y^{{\frak a},{\frak b}}_{\omega_1,\omega_2}(f)$ is valid. 
\par 
Our current hypotheses concerning $\sigma$, $\varrho$, $\vartheta$ and $f$ imply that 
those of the hypotheses of Lemma~6.7.9 that are concerned solely with 
$\sigma$, $\varrho$, $\vartheta$ and $f$ are satisfied. Likewise, our definition of the 
set $E$ is the same as that which is posited in Lemma~6.7.9. 
Consequently the only hypotheses of Lemma~6.7.9 requiring further 
consideration are those concerning either the function $d$ or the sequence 
$(f_n)_{n\in{\Bbb N}}$. 
\par
By Lemma~6.7.10 it follows that ${\cal H}_{\star}^{\sigma}(\varrho,\vartheta )\supseteq\{ f_n : n\in{\Bbb N}\}$, 
and that for each $n\in{\Bbb N}$ the sum formula for 
$Y^{{\frak a},{\frak b}}_{\omega_1,\omega_2}(f_n)$ is valid. 
Given the result (6.7.51) of Lemma~6.7.10, 
and the definition (6.7.50), it moreover follows that 
the sequence $(f_n)_{n\in{\Bbb N}}$ satisfies the condition (6.7.6) of 
Lemma~6.7.3, and that, for all $p\in{\Bbb Z}$, and all $\nu\in{\Bbb C}$ such that 
$|{\rm Re}(\nu)|\leq\sigma$, one has 
$\lim_{n\rightarrow\infty} f_n(\nu,p)=f(\nu,p)$.  
These observations enable one to conclude 
that those of the hypotheses of Lemma~6.7.9 that concern 
$(f_n)_{n\in{\Bbb N}}$, but not $d$, are satisfied. 
Therefore it now only remains to be shown that the function 
$d$ defined in (6.7.63) and (6.7.64) does satisfy the relevant 
hypotheses  of Lemma~6.7.9. 
\par 
Consider firstly the hypothesis that $d\in{\cal H}_{\star}^{\sigma}(\varrho,\vartheta )$. 
By Lemma~6.7.11 one has $d^{\flat}\in{\cal H}_{\star}^{\sigma}(\varrho,\vartheta )$; 
it therefore follows, given (6.7.63) and (6.7.64), that 
since ${\cal H}_{\star}^{\sigma}(\varrho,\vartheta )$ is a vector space over ${\Bbb C}$  
one has also  
$d\in{\cal H}_{\star}^{\sigma}(\varrho,\vartheta )$. Similarly, since Lemma~6.7.11 implies that 
the sum formulae for $Y^{{\frak a},{\frak a}}_{\omega_1,\omega_1}(d^{\flat})$ 
and $Y^{{\frak b},{\frak b}}_{\omega_2,\omega_2}(d^{\flat})$ are valid, 
one may deduce that the sum formulae 
for $Y^{{\frak a},{\frak a}}_{\omega_1,\omega_1}(d)$ 
and $Y^{{\frak b},{\frak b}}_{\omega_2,\omega_2}(d)$  
are also valid (this following in view of the intrinsic linearity, as regards 
their dependence on  
the test function $h$, of all the transforms, integrals and sums occurring in 
the definitions (6.7.1)-(6.7.3) of 
$\chi^{{\frak d},{\frak d}'}_{\omega,\omega'}(h)$, 
$X^{{\frak d},{\frak d}'}_{\omega,\omega'}(h)$ and 
$Y^{{\frak d},{\frak d}'}_{\omega,\omega'}(h)$). 
Finally, since $\sigma >1>2/9$ it follows, by (6.7.13), 
the results (6.7.51) and (6.7.61) of Lemma~6.7.10 and Lemma~6.7.11, 
and the definitions (6.7.63) and (6.7.64),    
that for all $(\nu,p)\in E$ and all $n\in{\Bbb N}$ one has 
$$|f_n(\nu,p)|\leq 
\exp(8\sigma^4) C_f(\varrho,\vartheta ) (1+|{\rm Im}(\nu)|)^{-\varrho} (1+|p|)^{-\vartheta} 
\leq 2^{\varrho +5}\exp(8\sigma^4) C_f(\varrho,\vartheta ) d^{\flat}(\nu,p)=d(\nu,p)\;.$$
\par 
Since it has been found that all of the hypotheses of Lemma~6.7.9 are 
satisfied, it therefore follows by that lemma that the sum formula 
for $Y^{{\frak a},{\frak b}}_{\omega_1,\omega_2}(f)$ is valid\quad$\blacksquare$ 

\bigskip 

\noindent{\bf The proof of Theorem~B.}\quad  
We shall deduce Theorem~B from the Weak Sum Formula (Proposition~6.7.2). 
In order to achieve this we employ the same method as that 
employed, in similar contexts, on both Page~64 of [5] and 
Pages~107-109 of [32]. 
\par 
We shall make use of the terminology introduced in the first few 
paragraphs of the current subsection (up to, and including, the 
paragraph containing Equation~(6.7.4)); 
our hypotheses concerning $\omega_1$, $\omega_2$, ${\frak a}$, ${\frak b}$, 
$g_{\frak a}$ and $g_{\frak b}$ are as stated in the second paragraph of 
this subsection. We suppose moreover that the set $E$ is as 
given by the equation~(6.7.13) of 
Lemma~6.7.6, and that $\sigma$, $\varrho$, $\vartheta$ and $h$ satisfy 
the relevant hypotheses of Theorem~B, so that 
$\sigma\in(1/2,1)$, $\varrho,\vartheta \in(3,\infty)$ and 
$h\in{\cal H}_{\star}^{\sigma}(\varrho,\vartheta ) ={\cal H}_0^{\sigma}(\varrho,\vartheta )$. 
If it follows that the sum formula for $Y^{{\frak a},{\frak b}}_{\omega_1,\omega_2}(h)$ 
is valid then Theorem~B is true. By means of 
an application of Lemma~6.7.9 we shall succeed in deducing 
the validity of the sum 
formula for $Y^{{\frak a},{\frak b}}_{\omega_1,\omega_2}(h)$, and shall thereby 
prove Theorem~B.  For the greater part of this proof we shall be concerned with 
the preliminary steps that enable this application of Lemma~6.7.9. 
\par 
Let $g$ be the mapping from $\{\nu\in{\Bbb C} : |{\rm Re}(\nu)|\leq\sigma\}\times{\Bbb Z}$ 
into ${\Bbb C}$ given by 
$$g(\nu,p)={h(\nu,p)\over\left( 1-\nu^2\right)^2\left( 4-\nu^2\right)^{-2}}\,.\eqno(6.7.65)$$ 
The function $j(\nu)=(1-\nu^2)^2 (4-\nu^2)^{-2}$  is  
meromorphic on ${\Bbb C}$ and even; given that $1/2<\sigma<1$, one has  
$${1\over\sqrt{|j(\nu)|}}=\left|{\nu -2\over\nu -1}\right|\left|{\nu +2\over\nu +1}\right| 
\leq\left( 1+{1\over |\nu -1|}\right)\left(1+{1\over |\nu +1|}\right) 
\leq 2\left( 1+{1\over 1-\sigma}\right) <\infty\qquad\qquad 
\hbox{($|{\rm Re}(\nu)|\leq\sigma$).}$$ 
Therefore, and since $h\in{\cal H}_{\star}^{\sigma}(\varrho,\vartheta )$, it follows from the definition 
(6.7.65) that one has 
$$g\in{\cal H}_{\star}^{\sigma}(\varrho,\vartheta )\;.\eqno(6.7.66)$$ 
\par 
For each $n\in{\Bbb N}$, we define $g^{\flat}_n$ to be the mapping from 
${\Bbb C}\times{\Bbb Z}$ into ${\Bbb C}$ given by 
$$g^{\flat}_n(\nu,p) 
=-i\sqrt{n\over\pi}\int\limits_{(0)} g(\xi,p) 
\exp\left( n (\nu-\xi)^2\right) {\rm d}\xi\;.\eqno(6.7.67)$$ 
As is shown by (11.28)-(11.30) of [32], it follows from (6.7.66) and (6.7.67) that, 
for $n\in{\Bbb N}$, $p\in{\Bbb Z}$ and $\nu\in{\Bbb C}$, one has  
$$\eqalignno{
\bigl| g^{\flat}_n(\nu,p)\bigr| 
 &\leq \sqrt{n\over\pi}\int\limits_{-\infty}^{\infty} \left| g(i\tau ,p)  
\exp\left( n (\nu-i\tau )^2\right)\right| {\rm d}\tau \ll_{g,\varrho,\vartheta } \cr 
 &\ll_{g,\varrho,\vartheta }\ \left( 1+|{\rm Im}(\nu)|\right)^{-\varrho}\left( 1+|p|\right)^{-\vartheta} 
\exp\!\left( n\left( {\rm Re}(\nu)\right)^2\right) . &(6.7.68)}$$
In particular, the integral on the right-hand side of Equation~(6.7.67) is absolutely convergent 
when $n\in{\Bbb N}$, $p\in{\Bbb Z}$ and $\nu\in{\Bbb C}$. 
Indeed, using only the fact that (by virtue of (6.7.66) and the hypothesis that 
$\varrho,\vartheta \in(3,\infty)$) the function $(\tau,p)\mapsto |g(i\tau,p)|$ 
is bounded on ${\Bbb R}\times{\Bbb Z}$, one can show that when 
$n\in{\Bbb N}$, $p\in{\Bbb Z}$ and $T\in(0,\infty)$ the integral 
on the right-hand side of Equation~(6.7.67) converges uniformly 
for all complex $\nu$ such that $\max\{ |{\rm Re}(\nu)| , |{\rm Im}(\nu)|\}\leq T$. 
It therefore follows, by (for example) Section~2.83 and Section~2.84 of [43], that 
$$\hbox{for $n\in{\Bbb N}$ and $p\in{\Bbb Z}$ the mapping 
$\,\nu\mapsto g^{\flat}_n(\nu,p)\,$ is holomorphic on ${\Bbb C}$.}\eqno(6.7.69)$$ 
By means of the change of variable $\xi =-\phi\,$ (say), and by the 
application of the identity $g(-\phi,p)=g(\phi,-p)\,$ (inferred from 
(6.7.66)), one finds that it moreover follows from (6.7.67) that one has 
$$g^{\flat}_n(\nu,p) 
=g^{\flat}_n(-\nu ,-p)\qquad\qquad\hbox{($n\in{\Bbb N}$, $p\in{\Bbb Z}$ 
and $\nu\in{\Bbb C}$).}\eqno(6.7.70)$$ 
\par 
We shall later have need of two 
further properties of the sequence $(g^{\flat}_n)_{n\in{\Bbb N}}$. 
The properties in question are, firstly, that 
$$\lim_{n\rightarrow\infty} \,g^{\flat}_n(\nu,p) 
=g(\nu,p)\qquad\qquad 
\hbox{($p\in{\Bbb Z}$, $\nu\in{\Bbb C}$ and $|{\rm Re}(\nu)|\leq\sigma$),}\eqno(6.7.71)$$ 
and, secondly, that 
$$\sup\left\{ 
\left( 1+|{\rm Im}(\nu)|\right)^{\varrho}\left( 1+|p|\right)^{\vartheta}\bigl| g^{\flat}_n(\nu,p)\bigr| 
\ :\ n\in{\Bbb N},\,p\in{\Bbb Z},\,\nu\in{\Bbb C}\ \,{\rm and} 
\ \,|{\rm Re}(\nu)|\leq\sigma\right\} 
<\infty\;.\eqno(6.7.72)$$ 
The proofs of both of these properties depend on the fact 
(used in both [5] and [32]) that, for $n\in{\Bbb N}$, $p\in{\Bbb Z}$ and 
$\nu\in{\Bbb C}$, one has 
$$\int\limits_{(0)} g(\xi,p) 
\exp\left( n (\nu-\xi)^2\right) {\rm d}\xi 
=\int\limits_{(\alpha)} g(\xi,p) 
\exp\left( n (\nu-\xi)^2\right) {\rm d}\xi\qquad\qquad 
\hbox{($-\sigma\leq\alpha\leq\sigma$).}\eqno(6.7.73)$$ 
By (6.7.67) and the case $\alpha ={\rm Re}(\nu)$ of (6.7.73), one obtains 
$$g^{\flat}_n(\alpha +it,p) 
=\sqrt{n\over\pi}\int\limits_{-\infty}^{\infty} 
g(\alpha +i\tau,p)\exp\left( -n (t-\tau)^2\right) {\rm d}\tau\qquad  
\hbox{($n\in{\Bbb N}$, $p\in{\Bbb Z}$, $t\in{\Bbb R}$ and $-\sigma\leq\alpha\leq\sigma$).} 
\eqno(6.7.74)$$ 
Given (6.7.66), the property (6.7.71) follows by (6.7.74) and 
the equation~1.17.6 of [38], while the property (6.7.72) follows by 
(6.7.74), (6.7.66) and the upper bound deduced in (11.29) of [32]. 
We remark that the result (6.7.73) is a corollary of the relation 
(6.7.66), which implies that 
the relevant mappings $\xi\mapsto g(\xi,p)$ are holomorphic 
on a neighbourhood of the strip 
$\{\xi\in{\Bbb C} : |{\rm Re}(\xi)|\leq\sigma\}$, and that in that strip 
one has 
$$g(\xi,p)\ll_{g,\varrho,\vartheta }\ (1+|{\rm Im}(\nu)|)^{-\varrho}(1+|p|)^{-\vartheta} 
\leq (1+|{\rm Im}(\nu)|)^{-3}\;.$$
\par 
We now define $d^{\flat}$ and $\,h^{\flat}_1 , h^{\flat}_2 , \ldots\,$ to 
be the mappings from $\{\nu\in{\Bbb C} : |{\rm Re}(\nu)|\leq 3/2\}\times{\Bbb Z}$ 
into ${\Bbb C}$ given by 
$$h^{\flat}_n(\nu,p) 
=\left( 1-\nu^2\right)^2 \left( 4 -\nu^2\right)^{-2} g^{\flat}_n(\nu,p) \qquad\qquad 
\hbox{($n\in{\Bbb N}$, $p\in{\Bbb Z}$, $\nu\in{\Bbb C}$ and $|{\rm Re}(\nu)|\leq 3/2$),} 
\eqno(6.7.75)$$ 
$$d^{\flat}(\nu,p)
=\left( 1-\nu^2\right)^2 \left( 4-\nu^2\right)^{-(\varrho +4)/2} \left( 1+|p|\right)^{-\vartheta}\qquad\qquad 
\hbox{($p\in{\Bbb Z}$, $\nu\in{\Bbb C}$ and  $|{\rm Re}(\nu)|\leq 3/2$).} 
\eqno(6.7.76)$$ 
We recall that the function $j(\nu)=(1-\nu^2)^2 (4-\nu^2)^{-2}$ is meromorphic on ${\Bbb C}$ and 
even.  One has, moreover, 
$$\sqrt{|j(\nu)|} 
={|\nu -1|\over |\nu -2|}\,{|\nu -(-1)|\over |\nu -(-2)|}\leq (1)(1)=1\qquad\qquad 
\hbox{($|{\rm Re}(\nu)|\leq 3/2$)}\eqno(6.7.77)$$ 
and $j'(1)=j(1)=0$. Consequently it follows from (6.7.68), (6.7.69), (6.7.70) and (6.7.75) 
that one has 
$$h^{\flat}_n\in{\cal H}_{\star}^{3/2}(\varrho,\vartheta )\qquad\qquad\hbox{($n\in{\Bbb N}$).}\eqno(6.7.78)$$ 
Therefore Proposition~6.7.2 implies that 
$$\hbox{for $\,n\in{\Bbb N}\,$ the sum formula for 
$\,Y^{{\frak a},{\frak b}}_{\omega_1,\omega_2}\left( h^{\flat}_n\right)\,$ is valid.}\eqno(6.7.79)$$ 
At the same time, it follows by (6.7.76) and the case $\sigma =3/2$ of Lemma~6.7.11 
that 
$$d^{\flat}\in{\cal H}_{\star}^{3/2}(\varrho,\vartheta )\;,\eqno(6.7.80)$$ 
$$d^{\flat}(\nu,p)\geq 2^{-(\varrho +5)}\left( 1+|{\rm Im}(\nu)|\right)^{-\varrho} 
\left( 1+ |p|\right)^{-\vartheta}\qquad\qquad 
\hbox{($\,(\nu,p)\in E$)}\eqno(6.7.81)$$ 
and 
$$\hbox{the sum formulae for 
$\,Y^{{\frak a},{\frak a}}_{\omega_1,\omega_1}\left( d^{\flat}\right)\,$ and 
$\,Y^{{\frak b},{\frak b}}_{\omega_2,\omega_2}\left( d^{\flat}\right)\,$ are valid.}
\eqno(6.7.82)$$ 
\par 
For $n\in{\Bbb N}$, we define $h_n$ to be the restriction of the 
mapping $h^{\flat}_n$ to the set $\{\nu\in{\Bbb C} : |{\rm Re}(\nu)|\leq\sigma\}\times{\Bbb Z}$. 
Similarly, we define $d^{\sharp}$ to be the restriction of the mapping 
$d^{\flat}$ to the set $\{\nu\in{\Bbb C} : |{\rm Re}(\nu)|\leq\sigma\}\times{\Bbb Z}$. 
Since $3/2>1>\sigma>1/2>2/9$, it follows by (6.7.78), (6.7.79), (6.7.80), (6.7.81), (6.7.82), the definition 
(6.7.13) and what is observed in our Remark~1.9.2 that 
$$d^{\sharp}\in{\cal H}_{\star}^{\sigma}(\varrho,\vartheta ) 
={\cal H}_0^{\sigma}(\varrho,\vartheta )\;,\eqno(6.7.83)$$ 
$$d^{\sharp}(\nu,p)\geq 2^{-(\varrho +5)}\left( 1+|{\rm Im}(\nu)|\right)^{-\varrho} 
\left( 1+ |p|\right)^{-\vartheta}\qquad\qquad 
\hbox{($\,(\nu,p)\in E$),}\eqno(6.7.84)$$ 
$$\hbox{the sum formulae for 
$\,Y^{{\frak a},{\frak a}}_{\omega_1,\omega_1}\left( d^{\sharp}\right)\,$ and 
$\,Y^{{\frak b},{\frak b}}_{\omega_2,\omega_2}\left( d^{\sharp}\right)\,$ are valid}
\eqno(6.7.85)$$ 
and 
$$h_n\in{\cal H}_{\star}^{\sigma}(\varrho,\vartheta )\qquad\qquad\hbox{($n\in{\Bbb N}$),}\eqno(6.7.86)$$ 
and that 
$$\hbox{for $\,n\in{\Bbb N}\,$ the sum formula for 
$\,Y^{{\frak a},{\frak b}}_{\omega_1,\omega_2}\left( h_n\right)\,$ is valid.}\eqno(6.7.87)$$ 
Given our definition of the mappings $h_1 , h_2 , \ldots \,$, and given the definitions 
in (6.7.65) and (6.7.75), it moreover follows by (6.7.71), (6.7.72) and 
(6.7.77) that one has both 
$$\lim_{n\rightarrow\infty} \,h_n(\nu,p) 
=h(\nu,p)\qquad\qquad 
\hbox{($p\in{\Bbb Z}$, $\nu\in{\Bbb C}$ and $|{\rm Re}(\nu)|\leq\sigma$)}\eqno(6.7.88)$$ 
and 
$$\sup\left\{ 
\left( 1+|{\rm Im}(\nu)|\right)^{\varrho}\left( 1+|p|\right)^{\vartheta}\bigl| h_n(\nu,p)\bigr| 
\ :\ n\in{\Bbb N},\,p\in{\Bbb Z},\,\nu\in{\Bbb C}\ \,{\rm and} 
\ \,|{\rm Re}(\nu)|\leq\sigma\right\}
<\infty\;.\eqno(6.7.89)$$ 
Note that, by (6.7.89), the condition (6.7.6) of Lemma~6.7.3 is satisfied 
when one has $f_n=h_n$ for all $n\in{\Bbb N}$. 
\par 
Just one further definition will put us in a position to apply Lemma~6.7.9. 
We define the mapping $d : \{\nu\in{\Bbb C} : |{\rm Re}(\nu)|\leq\sigma\}\times{\Bbb C} 
\rightarrow{\Bbb C}$ by setting 
$$d=2^{\varrho +5} C^{*}_h(\varrho,\vartheta )\,d^{\sharp}\;,\eqno(6.7.90)$$ 
where $C^{*}_h(\varrho,\vartheta )$ denotes the non-negative real number equal to the supremum 
on the left-hand side of the inequality~(6.7.89). It then follows, by linearity, 
from (6.7.83) and (6.7.85) that 
$$d\in{\cal H}_{\star}^{\sigma}(\varrho,\vartheta )\ \,\hbox{and the sum formulae for 
$\,Y^{{\frak a},{\frak a}}_{\omega_1,\omega_1}\left( d\right)\,$ and 
$\,Y^{{\frak b},{\frak b}}_{\omega_2,\omega_2}\left( d\right)\,$ are valid.}\eqno(6.7.91)$$ 
By (6.7.89), the definition of $C^{*}_h(\varrho,\vartheta)$, 
(6.7.84) and (6.7.90) one has, moreover, 
$$\bigl| h_n(\nu,p)\bigr| 
\leq C^{*}_h(\varrho,\vartheta )\left( 1+|{\rm Im}(\nu)|\right)^{-\varrho} 
\left( 1+|p|\right)^{-\vartheta}  
\leq d(\nu,p)\qquad\qquad  
\hbox{($n\in{\Bbb N}$ and $(\nu,p)\in E$).}\eqno(6.7.92)$$ 
\par 
By our initial hypotheses 
in this proof, and by (6.7.86), (6.7.87), (6.7.88), (6.7.89), (6.7.91) and~(6.7.92), 
it follows that the hypotheses of Lemma~6.7.9 are satisfied if one substitutes 
$h$ and the sequence $(h_n)_{n\in{\Bbb N}}$ for the function 
$f$ and sequence $(f_n)_{n\in{\Bbb N}}$ of that lemma 
(while taking $d$ to be given by (6.7.90)). Lemma~6.7.9 therefore implies 
that the sum formula for $Y^{{\frak a},{\frak b}}_{\omega_1,\omega_2}(h)$ is valid. 
This conclusion means that the results of Theorem~B are obtained; since 
it is a conclusion that has been arrived at independently of any assumptions other than 
the stated hypotheses of Theorem~B, 
our proof of that theorem is now complete\ $\blacksquare$ 
                                             
\bigskip 

\centerline{\bf Acknowledgements.}

\medskip

Between 2004 and 2006 this work (and the author) were supported through a 
fellowship associated with an EPSRC funded project, 
`The Development and Application of Mean-Value Results in Multiplicative
Number Theory' (GR/T20236/01), which was led by Glyn Harman of Royal Holloway, 
University of London: the author thanks Glyn for his advice, 
encouragement and patience. 
\par 
Some of the preliminary work on the writing-up of the 
appendix to this paper was done while the author was employed 
(between 2006 and 2010) as a 
Lecturer in the School of Mathematics of Cardiff University. 
\par 
In the course of completing this research the 
author benefitted greatly from correspondence 
with Roelof Bruggeman of Utrecht. 
The author wishes to record his 
gratitude to Roelof for his detailed comments on previous draughts of this paper, 
and for his comprehensive answers to many questions: those answers 
included an 11~page PDF file outlining a proof of Theorem~B. In completing 
the proof of Theorem~B the author was aided by 
Roelof's preliminary notes on the sum formula for $SL_2({\Bbb C})$. 
These unpublished notes (contained in a PDF file of 158~pages) 
were preparatory to work reported on in [5]; some of the 
results obtained in them are used in the proof of Lemma~6.5.5. 
\par 
The author has also received (from the beginning to the end of this work)  
much valuable advice and encouragement from 
Yoichi Motohashi of Tokyo. The author is especially grateful to Yoichi for 
drawing his attention to the preprint [36]; it turned out 
that the method which Yoichi developed there could easily be adapted to 
yield bounds for the generalised Kloosterman sums considered in Subsection~6.1.  
\par 
The author thanks his parents for their unflagging 
encouragement and helpful support.   
\bigskip

\centerline{\bf References.}

\medskip

\parskip = 3 pt 

\item{[1]} T. M. Apostol, Mathematical Analysis, second ed., 
World Student Series, Addison-Wesley, Reading MA, 1974. 

\item{[2]} S. Bochner, Lectures on Fourier Integrals 
(English translation by M. Tenenbaum and H. Pollard), 
Annals of Mathematics Studies No. 42, Princeton University Press, Princeton NJ, 1959. 

\item{[3]} R. W. Bruggeman, Fourier coefficients of cusp forms, Invent. Math. 
45 (1978) 1-18.

\item{[4]} R. W. Bruggeman and R. J. Miatello,
Estimates of Kloosterman sums for groups of real rank one,
Duke Math. J. 80 (1995) 105-137.

\item{[5]} R. W. Bruggeman and Y. Motohashi, 
Sum formula for Kloosterman sums and fourth moment of the Dedekind
zeta-function over the Gaussian number field, Functiones et Approximatio,
31 (2003) 23-92.

\item{[6]} D. Bump, Automorphic Forms and Representations, 
Cambridge Stud. Adv. Math., vol. 55, 
Cambridge University Press, New~York, 1998. 

\item{[7]} J. Cogdell, J. S. Li, I. Piatetski-Shapiro, P. Sarnak, 
Poincar\'{e} series for $SO(n,1)$, Acta Math. 167 (1991) 229-285.  

\item{[8]} M. D. Coleman, The distribution of points at which binary quadratic 
forms are prime, Proc. London Math. Soc. (3) 61 (1990) 433-456.

\item{[9]} J.-M. Deshouillers and H. Iwaniec,
Kloosterman sums and Fourier coefficients of cusp forms,
Invent. Math. 70 (1982) 219-288.

\item{[10]} W. Duke, Some problems in multidimensional analytic number theory,
Acta Arith. 52 (1989) 203-228.

\item{[11]} J. Elstrodt, F. Grunewald and J. Mennicke, 
Groups Acting on Hyperbolic Space: Harmonic Analysis and Number Theory, 
Springer Monographs in Mathematics, 
Springer-Verlag, Berlin Heidelberg, 1998.

\item{[12]} R. Goodman and N. R. Wallach, 
Whittaker vectors and conical vectors, 
J. Funct. Anal. 39 (1980) 199-279.

\item{[13]} J. H. Graf, Ueber die addition und subtraction der argumente bei 
Bessel'schen functionen nebst einer anwendung, Math. Annalen, 
43 (1893) 136-144. 

\item{[14]} K. B. Gundlach, Ueber die darstellung der ganzen spitzenformen 
zu den idealstufen der Hilbertschen modulgruppe und die abschaetzung ihrer 
Fourierkoeffizienten, Acta Math. J. 
92 (1954) 309-345. 

\item{[15]} G. H. Hardy and E. M. Wright, An Introduction to the Theory of Numbers, 
fifth ed., Oxford University Press, Oxford, 1979 
(reprinted with corrections in 1983 and 1984).

\item{[16]} Harish-Chandra, Automorphic Forms on Semisimple Lie Groups  
(notes taken by J. G. M. Mars), Lecture Notes in Mathematics, vol. 62, 
Springer-Verlag, Berlin New~York, 1968. 

\item{[17]} E. Hecke, Eine neue art von zetafunktionen und ihre beziehungen
zur verteilung der primzahlen, II, Math. Zeitschrift 6 (1920) 11-51.

\item{[18]} M. N. Huxley, Area, Lattice Points and Exponential Sums, 
London Math. Soc. Monogr., New Ser. vol.~13, 
Oxford Science Publications, Clarendon Press, Oxford, 1996.

\item{[19]} A. Ivi\'{c}, The Riemann Zeta-Function: Theory and Applications,
Dover Publications, Inc., Mineola, New York, 2003.

\item{[20]} H. Iwaniec, Fourier coefficients of cusp forms and 
the Riemann zeta-function, 
S\'{e}m. Th. Nombres, Univ. Bordeaux 18 (1979/80) 1-36. 

\item{[21]} H. Iwaniec, Topics in Classical Automorphic Forms, 
Grad. Stud. Math., vol. 17, Amer. Math. Soc., Providence RI, 1997.

\item{[22]} H. Iwaniec, Spectral Methods of Automorphic Forms, second ed.,  
Grad. Stud. Math., vol. 53, Amer. Math. Soc. and  
Revista Matem\'{a}tica Iberoamericana, 2002.

\item{[23]} M. Jutila, On spectral large sieve inequalities, 
Functiones et Approximatio, 28 (2000) 7-18.

\item{[24]} M. Jutila and Y. Motohashi, 
Uniform bound for Hecke $L$-functions, 
Acta Math. 195 (2005) 61-115.

\item{[25]} H. H. Kim, On local $L$-functions and normalized intertwining operators, 
Canad. J. Math. 57 (3) (2005) 535-597. 

\item{[26]} H. H. Kim and F. Shahidi, Cuspidality of symmetric powers with applications, 
Duke Math. J. 112 (2002), 177-197. 

\item{[27]} A. W. Knapp, Representation Theory of Semisimple Groups, 
an Overview Based on Examples, 
Princeton Landmarks in Math. ed., 
Princeton University Press, Princeton Oxford, 2001.

\item{[28]} N. V. Kuznetsov, The Petersson hypothesis for forms of 
weight zero and the Linnik hypothesis, 
Preprint No. 02, Khabarovsk Complex Res. Inst., Acad. Sci. USSR, 
East Siberian Branch, Khabarovsk, 1977 (in Russian).

\item{[29]} N. V. Kuznetsov, Petersson's conjecture for cusp forms of weight zero 
and Linnik's conjecture. Sums of Kloosterman sums, Math. USSR Sb. 39 (3) (1981), 299-342.

\item{[30]} S. Lang, Real Analysis, second edition,
Addison-Wesley, Reading MA, 1983.

\item{[31]} R. P. Langlands, On the Functional Equations Satisfied by Eisenstein Series,  
Lecture Notes in Math., vol. 544, Springer-Verlag, Berlin, 1976.

\item{[32]} H. Lokvenec-Guleska, Sum Formula for SL${}_2$ over
Imaginary Quadratic Number Fields, PhD Thesis, University of Utrecht, 2004 
(includes summaries in Dutch and Macedonian).

\item{[33]} R. Miatello and N. R. Wallach, Automorphic forms constructed from Whittaker vectors, 
J. Funct. Anal. 86 (1989) 411-487. 

\item{[34]} R. Miatello and N. R. Wallach, Kuznetsov formulas for real rank one groups,  
J. Funct. Anal. 93 (1990) 171-206. 

\item{[35]} Y. Motohashi, Spectral Theory of the Riemann Zeta-Function, 
Cambridge Tracts in Math. 127, Cambridge University Press, 1997.

\item{[36]} Y. Motohashi, The Riemann zeta-function and Hecke congruence 
subgroups~II., arXiv:0709.2590v2 [math.NT], 2008. 

\item{[37]} W. M\"{u}ller, Weyl's law for the cuspidal spectrum of $SL_n$, 
Ann. of Math. 165 (2007) 275-333. 

\item{[38]} F. W. J. Olver, D. W. Lozier, R. F. Boisvert, C. W. Clark, 
NIST Handbook of Mathematical Functions, Cambridge University Press, 
New~York, 2010. 

\item{[39]} S. Rhagavan and J. Sengupta, On Fourier coefficients 
of Maass cusp forms in $3$-dimensional hyperbolic space, 
Proc. Indian Acad. Sci. (Math. Sci.) 104 (1)  
(1994) 77-92. 

\item{[40]} W. Rudin, Real and Complex Analysis, second ed.,  
Tata McGraw-Hill Publishing Co. Ltd., New York, New Delhi, 1974.   

\item{[41]} A. Selberg, On the estimation of Fourier coefficients of 
modular forms, Proc. Symp. Pure Math. 8 (1965) 1-15. 

\item{[42]} M. Sugiura, 
Unitary Representations and Harmonic Analysis: an Introduction, 
second ed., North-Holland Mathematical Library, vol. 44, 
North-Holland, Amsterdam New York Tokyo, 1990. 

\item{[43]} E. C. Titchmarsh, The Theory of Functions, second ed., 
Oxford University Press, Oxford, 1939. 

\item{[44]} G. N. Watson, A Treatise on the Theory of Bessel Functions, 
first ed., Cambridge University Press, 1922. 

\item{[45]} N. Watt, Kloosterman sums and a mean value for Dirichlet polynomials, 
J. Number Theory 53 (1995) 179-210. 

\item{[46]} N. Watt, Weighted spectral large-sieve inequalities 
for Hecke congruence subgroups of ${\rm SL}(2,{\Bbb Z}[i])$, submitted preprint. 

\item{[47]} N. Watt, Weighted fourth moments of Hecke zeta functions 
with groessencharakters,  preprint in preparation 
(incomplete and not yet submitted). 

\item{[48]} E. T. Whittaker and G. N. Watson, 
A Course of Modern Analysis, fourth ed., Cambridge University Press, 
New~York, 1927.

\bye